\documentclass[a4paper,11pt]{amsart}
\usepackage{amssymb,amsmath,amsthm,stmaryrd,bm,enumitem,hyperref,mathrsfs, xcolor} %graphicx 

\def\res{\hbox{ {\vrule height .3cm}{\leaders\hrule\hskip.3cm}}\hskip5.0\mu}

\newcommand\beqn{\begin{equation}}
\newcommand\eeqn{\end{equation}}
\newcommand\beqny{\begin{eqnarray}}
\newcommand\eeqny{\end{eqnarray}}
\newcommand\beqnyn{\begin{eqnarray*}}
\newcommand\eeqnyn{\end{eqnarray*}}

\usepackage{amsmath,amsfonts,amssymb,amsthm}
\newtheorem{theorem}{Theorem}[section]
\newtheorem{lemma}[theorem]{Lemma}

\newtheorem{corollary}[theorem]{Corollary}
\newtheorem{definition}[theorem]{Definition}
\newtheorem{remark}[theorem]{Remark} 
\newtheorem{hypothesis}[theorem]{Hypothesis}

\newtheorem{theoremX}{Theorem}

\newtheorem{theoremA}{Theorem}

\newtheorem{theoremB}{Theorem}

\newtheorem{theoremC}{Theorem}

\setlist[itemize]{labelindent=10mm,leftmargin=*}
\setlist[enumerate]{labelindent=6mm,leftmargin=*}

\allowdisplaybreaks[3]
\numberwithin{equation}{section}

\newcommand{\op}[1]{\operatorname{\text{\rm #1}}}

\title[branch points of planar frequency $\neq 2$]{Analysis of singularities of area-minimizing currents, part~III: branch points of planar frequency $ \neq 2$, higher order asymptotics, and the local topology} 

\author{Brian Krummel \& Neshan Wickramasekera}

\begin{document}
\begin{abstract}
This paper is the third in a five-part series developing a  new framework for 
analysing $n$-dimensional area-minimizing rectifiable currents, initially in ${\mathbb R}^{n+m}$ ($m \geq 2$) and extending the results to general Riemannian ambient spaces in Part~V. 
The present paper focuses on the asymptotic and topological nature of the current near the set ${\mathcal B}$ of rapid-decay branch points 
defined in Part~I. 
The key preliminary analysis done in Part~I based on the intrinsic planar frequency function introduced therein, together with the asymptotic analysis away from ${\mathcal B}$ done in Part~II, implies that the set ${\mathcal B}$ contains ${\mathcal H}^{n-2}$-almost all branch points.

The overarching series develops a  
unified
approach to several questions concerning the local structure of the current, and provides as a corollary  structural results for the singular set. Our framework diverges fundamentally from Almgren's pioneering 1983 work in its geometric philosophy and yields his optimal Hausdorff dimension upper bound $n-2$ for the singular set more efficiently. 
In Almgren's theory (as well as in more recent work of several authors adapting the Almgren framework), a technically demanding foundational component  is the iterative construction of auxiliary center manifolds, at the outset, corresponding to each branch point. By contrast, a feature of our framework---central to the structural analysis of the current we provide---is that this machinery is entirely bypassed for branch points of planar frequency $\neq 2$. 
This is achieved by replacing the center manifold with two geometric principles (established in Part~I) which arise directly from the structure of the current: a rapid-planar-decay/weak-non-planar-approximation dichotomy, and the existence of non-trivial tangent functions (blow-ups) relative to the tangent planes at every point in ${\mathcal B}$. These are both consequences of the monotonicity  properties of planar frequency. 

Using these tools, in this paper we establish three main results  for ${\mathcal B}^{(\neq 2)}$, the set of points in ${\mathcal B}$ of planar frequency $\neq 2$: 
\begin{itemize}
\item[1.] \emph{Higher-order asymptotics}: At 
${\mathcal H}^{n-2}$-a.e.\ point in ${\mathcal B}^{(\neq 2)}$, the current admits an expansion of finite order greater than $1 + \epsilon$, with precise decay estimates for the remainder term, where  $\epsilon >0$ is a universal constant.

\item[2.] \emph{Structural decomposition}: The set ${\mathcal B}^{(\neq 2)}$ locally decomposes into finitely many pairwise disjoint, locally $(n-2)$-rectifiable sets (of locally finite measure).

\item[3.] \emph{Topological control}: Near any branch point in ${\mathcal B}^{(\neq 2)}$ satisfying a specific planar-frequency criterion, the current is homeomorphic to an $n$-dimensional disk and admits  a 
$C^{1, \mu}$ parameterization. (Classical complex algebraic examples demonstrate that when this frequency criterion fails, the current need not be locally homeomorphic to an $n$-disk).
\end{itemize}

 In Part~IV these results are extended to planar frequency 2 branch points, capitalising on the fact that the center manifold is decisively more effective and geometrically canonical with respect to such points: specifically, about each planar-frequency-2 branch point $Z$ there is a center manifold (essentially unique for ${\mathcal H}^{n-2}$ a.e.\ such point) over which the current can be quantitatively approximated by the graph of a multi-valued normal map whose critical nodal set contains $Z$ and all nearby planar frequency 2 points of the same density.

Our results are a generalization to arbitrary dimension of the classical structure theory for two-dimensional area minimizing currents established by White, Chang, and Micallef–White. 
Overall, our framework brings the singularity theory of higher-codimension area minimizers closer in spirit to the codimension 1 setting and to classical PDE theory, where uniform a priori estimates serve as the primary organizing principle. 
This PDE theoretic aspect extends to the role of the center manifold within the program: the center manifold is invoked only in a canonical setting conforming to these principles, where it passes through all relevant branch points (namely, planar frequency 2 points) and remains well-centered at each, facilitating locally uniform decay estimates for the height of the current relative to the center manifold.
\end{abstract}
\maketitle 

\tableofcontents
%{
%\titleformat{\section}[runin]{\normalfont\bfseries}{\thesection}{1em}{}
%\paragraph{\sc Extended Abstract}
%\hfill

\section{Introduction}  
This is the third paper in a five-part series that includes \cite{KrumWica}, \cite{KrumWicb},  \cite{KrumWicd} and \cite{KrumWice}. In this series, we develop a new framework for studying 
singularities of $n$-dimensional area-minimizing rectifiable currents $T$ in a Riemannian manifold of dimension $n + m$, $m \geq 2$. In the present article (as in the first two papers \cite{KrumWica},  \cite{KrumWicb} and  the fourth part \cite{KrumWicd}) we shall assume that the ambient manifold is  
an open subset of $\mathbb{R}^{n+m}$. The additional technical arguments necessary to extend the main results to general Riemannian ambient spaces are  described in Part~V (\cite{KrumWice}).

Let ${\rm sing} \, T$ denote the set of interior singular points of $T$, i.e.\ points of ${\rm spt} \, T \setminus {\rm spt} \, \partial \, T$ near which ${\rm spt} \, T$ is not an $n$-dimensional embedded smooth submanifold. Unlike in the case of codimension $m=1$, when 
$m \geq 2$ the current $T$ may admit \emph{branch points}, i.e.\ points in ${\rm sing} \, T$ where one tangent cone to $T$ is supported on a plane. This difference introduces significant additional complexity to the 
analysis of singularities in higher codimension. 

\subsection{Historical context} 
Almgren's monumental classical work dating back to 1983 established the sharp Hausdorff dimension upper bound $n-2$ for 
${\rm sing} \, T.$  Several decades after Almgren had completed this work, his unpublished typewritten monograph was posthumously transcribed by Taylor and Scheffer into a 900+ page book (\cite{Almgren}).
 Subsequently, De~Lellis and Spadaro published the proof in a modular form (\cite{DeLSpa1}, \cite{DeLSpa2}, \cite{DeLSpa3}),  closely adhering to the original conceptual architecture of the argument while incorporating  technical and expository streamlining for parts of it, including clean decoupling of a certain delicate construction---the center manifold---from the rest of the proof. 
 Their work brought significant clarity and conciseness to the arguments in \cite{Almgren}, making them more accessible and thereby generating a renewed and broader interest in this fundamental problem.

 A technically demanding foundational aspect of this classical framework  is an iterative use of the auxiliary center manifolds for each branch point. Heuristically speaking, 
 these center manifolds are an approximation to the average of the ``sheets'' of the current at scales where the current remains  sufficiently weakly close to a plane.  They are utilized to devise a blow-up mechanism for $T$ that ensures ``sheet separation'' in the blow-up limit; about any given branch point and along a suitable sequence of scales tending to zero, the normal height of $T$ relative to the (changing) center manifolds is blown up. 
 The well-centering property of the center manifolds forces the blow-up limit---a multi-valued function---to be average free, and a monotonicity property of the celebrated Almgren frequency function associated with the normal height guarantees that the blow-up limit is non-zero. Thus the sheets of $T$ do not all collapse onto one-another in the blow-up limit, countering the inherent degeneracy of $T$ at branch points and thereby enabling extraction of information about the branch set from the blow-up.

 While some form of  ``modding out the average''  is arguably necessary for studying the branch-set dimension and related questions, the high technical complexity of  the classical theory is closely tied to the fact that it addresses the singular-set dimension question  without simultaneously developing quantitative information at branch points. 
 Under this objective, one must work with very little information on the current's behavior near branch points. Consequently, all branch points must be treated on an equal footing and center manifold sequences employed  uniformly across all branch points. In such generality, none of the center manifolds corresponding to a given branch point can be guaranteed to contain that point or any other branch point. 
 
 The ideal setting with regard to center manifolds on the other hand  
 would be to have, for a given  branch point $Z$: (a) a \emph{single} center manifold about $Z$ that contains $Z$ and all other nearby branch points we wish to study; (b) a single multi-valued normal map over this center manifold which serves as a good approximate graphical parameterization  of the current (in the sense that the graph of the normal map agrees with most of the current) so that the accuracy of the approximation improves when scaling about any point in the relevant branch set; and (c) that all relevant branch points are contained in the critical nodal set of this normal map. In such a setting, the study of the branch points of interest would amount to a question of PDE theoretic flavour: namely, study of the critical nodal set of a map satisfying variational properties.

Our framework reduces the role of the center manifold to this ideal case. We achieve this by shifting the focus, at the outset,  to the structural information  encoded in the current's decay properties at branch points, establishing \emph{geometric decay estimates as the organizing principle} of the analysis in place of auxiliary constructions necessitated by direct sheet-separation at all branch points. At a broad conceptual level, this is a refinement of the guiding principle in codimension 1 regularity 
theory---either for area-minimizers as in \cite{DG} where reduction to multiplicity 1 is possible, or more generally and closer to the present setting, for stable varifolds in arbitrary multiplicity as in \cite{Wic14} (see also \cite{MW}). 
This shift from the classical Almgren framework is necessitated by our new structural and local topological results for the current. Besides these new results, this viewpoint yields a substantial simplification of the classical framework for the dimension bound of the singular set.

\subsection{The present work} 
Concretely, our framework directly links locally uniform decay estimates at branch points first to the size of the singular set, and then to its rectifiability properties and to the current's local structure. Through these decay estimates, our approach restricts the need for center manifolds to a specific type of branch point: namely, those where the current satisfies a locally uniform estimate giving decay towards a unique plane at a quadratic rate in the scale. 
Rather strikingly, none  of our results concerning other singularities (i.e.\ other branch points, or non-branch point singularities)---in fact no result in any of the first three parts  \cite{KrumWica, KrumWicb} or the present paper---requires center manifolds of any kind.

Branch points of quadratic decay on the other hand constitute the canonical setting for Almgren's center manifold. For each such point $Z$, only a single center manifold is necessary, contrasting sharply with the sequence required for a general branch point. Furthermore, this center manifold  satisfies crucial additional properties:  
it contains $Z$ along  with \emph{all} nearby quadratic-decay branch points of the same density; these branch points are all contained in the critical (nodal) set of a single multi-valued ``normal map'' defined on the center manifold; the graph of this normal map provides a precise approximation to the current, with an accuracy that improves quantifiably upon approaching any quadratic-decay branch point; and additionally, for ${\mathcal H}^{n-2}$ a.e.\ such branch point $Z$, the center manifold is uniquely determined by $Z$ up to a lower order term (see Theorem~\ref{thm:B} below).

To handle  all other branch points (in fact all other singularities), our approach  develops an alternative, direct intrinsic method based on the planar frequency function introduced in \cite{KrumWica}. 
While this eliminates the technical burden of having to rely on auxiliary constructions for slow-decay branch points, our primary motivation for reducing this reliance lies elsewhere. Indeed, this structural feature of our framework is foundational 
to obtaining higher-order asymptotics at typical (i.e., 
$\mathcal{H}^{n-2}$~a.e.)
branch points (item~(iv) in Section~\ref{summary-intro}) and criteria for local topological conclusions (item~(v) in Section~\ref{summary-intro}). These two results are the main outcomes of the present article and \cite{KrumWicd}. 

To illustrate the point made above regarding when to use a center manifold (and when not to), consider 
$T_{1} = \{(z, w) \, : \, (w - z^{2})^{2} = z^{5}\} \subset {\mathbb C}^{2}$. As a real surface this is locally area-minimizing in ${\mathbb R}^{4}$, and writing it as $w = z^{2} \pm z^{5/2},$ we see that at the branch point $(0, 0),$ this decays to its tangent plane $\{w=0\}$ quadratically. 
The average of the two sheets is the graph of $w = z^{2},$ which can be taken as the center manifold. Since $z^{5/2}$ is lower order than $z^{2}$, the tangent function (blow-up) of $T_{1}$ at $(0, 0)$ relative to the tangent plane will be a multiplicity 2 copy of $w = z^{2}$, which does not retain the branching behaviour. Blowing up relative to the center manifold however will be $w = \pm z^{5/2}$, which retains branching. So in this case we must use a center manifold. 

Now consider $T_{2} = \{(z, w) \, : \, (w - z^{2})^{2} = z^{3}\} \subset {\mathbb C}^{2},$ which is the same as $w = z^{2} \pm z^{3/2},$ with the average again given by $w= z^{2}$. The classical approach will insist blowing up relative to a sequence of center manifolds approximating $w = z^{2}.$ (In this simple example the average $w = z^{2}$ can itself be taken as the center manifold.) However, since at $(0, 0)$ the decay rate of $T_{2}$ to the tangent plane is $O(|z|^{3/2})$ (less than quadratic), 
blowing up relative to the tangent plane suffices to produce the branched surface $w = \pm z^{3/2}$ in the limit. Our work establishes that this is a general phenomenon for any area minimizing current, and the tool we need to carry out the tangent-plane blowing up in a useful way  is the planar frequency function. 

\subsection{Summary of the overall main results} \label{summary-intro}
The framework outlined above yields the following new results, described somewhat informally:
\begin{itemize}
\item [(i)]  (\emph{Uniqueness of tangent cones}). $T$ has a unique tangent cone at ${\mathcal H}^{n-2}$ a.e.\ singular point; in particular, $T$ has a unique tangent plane $P_{Z}$ 
at ${\mathcal H}^{n-2}$ a.e.\ branch point $Z$  and $T$ satisfies, at ${\mathcal H}^{n-2}$ a.e.\ branch point $Z$ with density $\Theta(T, Z) = q$ for a given 
(positive integer) value $q$, a locally uniform estimate giving decay towards $Z + P_{Z}$ at a rate $o(|X - Z|^{1 + \alpha})$ where 
$\alpha = \alpha(n, m, q) \in (0, 1)$  (\cite[Theorem~1.1]{KrumWicb}, reproduced as Theorem~\ref{thm:X} below);  
\item [(ii)] (\emph{Structure of the singular set}). ${\rm sing} \, T$ is locally the union of a finite number of pairwise disjoint, locally compact, locally $(n-2)$ rectifiable sets
(Theorem~\ref{thm:C} below); 
\item [(iii)]  (\emph{Branching order}). For ${\mathcal H}^{n-2}$ a.e.\ branch point $Z$ of $T$, there is a (finite) rational number ${\mathcal O}_{T}(Z),$ the \emph{branching order} of $T$ at $Z$, determined uniquely by $Z$ and having the form ${\mathcal O}_{T}(Z) = \ell_{Z}/q_{Z}$  for relatively prime positive integers $\ell_{Z}$, $q_{Z}$ with $\ell_{Z} > q_{Z}$ and $q_{Z} \leq q = \Theta \, (T,Z),$ such that the rate at which the points of $T$ come together on approach to $Z$ is \emph{precisely} $O(|X - Z|^{{\mathcal O}_{T}(Z)}).$ (See Section~\ref{branching-order-intro} for the precise meaning of this statement, noting that this says much more than the standard doubling conditions implied by monotonicity of frequency functions). (Theorem~\ref{thm:B} below).
\item [(iv)] (\emph{Uniqueness of tangent functions}). At ${\mathcal H}^{n-2}$ a.e.\ branch point $Z$, $T$ has a unique, non-zero, homogeneous, average-free, locally Dirichlet energy minimizing $q$-valued tangent function (blow-up) $\varphi^{(Z)} \, : \, {P_{Z}} \to {\mathcal A}_{q}(P_{Z}^{\perp})$, where $q = \Theta(T, Z)$; the degree of homogeneity of $\varphi^{(Z)}$ is 
${\mathcal O}_{T}(Z)$, and $\varphi^{(Z)}$ is cylindrical, i.e.\ invariant under translation along an $(n-2)$-dimensional linear subspace of $P_{Z}$ (Theorem~\ref{thm:B} below).  
\item [(v)] (\emph{Topology near branch points}). Near a density $q$ branch point $Z$ where $T$ has a unique tangent plane $P_{Z}$ and where ${\mathcal O}_{T}(Z) = \ell_{Z}/q$ for some  $\ell_{Z}$ ($> q$) coprime to $q$,  we have that ${\rm spt} \, T$ is topologically an $n$-dimensional disk parameterized by a 
$C^{1, \mu}$ function (for some fixed $\mu  = \mu(n, m, q) \in  (0, 1)$) provided also that $T$ satisfies the additional condition that 
$Z$ is ``$(n-2)$-strongly non-isolated,''  i.e.\ near $Z$ the set 
$\{Z^{\prime} \, : \, \Theta \, (T, Z^{\prime}) = q \; \mbox{and } \; {\mathcal O}_{T}(Z^{\prime}) \geq {\mathcal O}_{T}(Z)\}$ projects fully onto an open ball of an $(n-2)$-dimensional subspace. In particular, this conclusion holds if, in place of the assumption that $\ell_{Z}, q$ are co-prime, we assume that $T$ is irreducible near $Z$ and $q$ is a prime, and keep all other assumptions as above. (Theorem~\ref{thm:A} below). 
\end{itemize}

\noindent
\begin{remark}\label{simult-work}
{\rm Contemporaneously with and independently of our work, De~Lellis and Skorobogatova initiated a program  taking the full extent of Almgren's 
iterative-center-manifolds framework as a starting point (\cite{DeLSko1, DeLSko2}), which, combined with their subsequent joint work with Minter (\cite{DeLMinSko}), establishes two of the main results (i)--(v)  described above: namely, uniqueness of tangent cones at ${\mathcal H}^{n-2}$~a.e.\ point of $T$ (item~(i)) and countable $(n-2)$-rectifiability of ${\rm sing} \, T$ (included in item~(ii)). We refer the reader to Section~\ref{comparison} below, and to the more comprehensive discussion in \cite[Section~2]{KrumWica}, for a comparison of this parallel work with ours.} 
\end{remark}

The above results (i)--(v) are counterparts or generalizations of several celebrated earlier works for area minimizers $T$ of either codimension $m=1$ or dimension $n=2$. If $m=1$, $T$ has no branch points nor singular tangent cones that are invariant under translation along a subspace of dimension $\geq n-6,$ and hence ${\rm sing} \, T$ has Hausdorff dimension $\leq n-7$ (as established by the foundational work of De~Giorgi, Federer, Fleming, Almgen and Simons in the decade 1960-1970). In this case,  by the combined work of Simon (\cite{Sim95}) and Naber--Valtorta (\cite{NV15}), we in fact have that (i) and (ii) hold with $(n-7)$ in place of $(n-2)$ and with the additional conclusion that ${\rm sing} \, T$ has locally finite ${\mathcal H}^{n-7}$ measure. Our work does not subsume these results, and in fact the techniques used in \cite{Sim95}, \cite{NV15} rely on the absence of higher multiplicity and are of a fundamentally different nature.

In case the dimension $n=2$, the combined work of White (\cite{Whi83}), Chang (\cite{Cha88}) and Micallef--White (\cite{MicWhi95}) had established all of (i)--(v) and more.  In this case, in fact, the singular set consists of isolated points (\cite{Cha88}); moreover, in the topological conclusion (v) the requirement that $Z$ be $(n-2)$-strongly non-isolated is, of course, automatically trivially satisfied. 
See Section \ref{2d-case-intro} for a more detailed discussion on the two-dimensional theory. 
In general dimensions, without the condition that $Z$ is $(n-2)$-strongly non-isolated, $T$ need not be topologically an $n$-dimensional disk near $Z$. This is seen by examples of area minimizers provided by complex algebraic varieties in ${\mathbb C}^{n}$, such as
$V = \{(x, y, w) \, : \, w^{2} = x^{3}(x^{2} - y^{3})\} \subset {\mathbb C}^{3}.$ For this example, the (real) 4-dimensional multiplicity 1 current $T_{V}$ induced by $V$ is locally area minimizing in ${\mathbb R}^{6}$ (as is the case with all holomorphic subvarieties), and has ${\rm sing} \, T_{V} = \{(0, a, 0) \, : \, a \in {\mathbb C}\}$ with every singularity a density 2 branch point. Since ${\mathcal O}_{T_{V}}(0, a, 0) = 3/2$ if $a \neq 0$ and ${\mathcal O}_{T_{V}}(0, 0, 0) = 5/2$, the branch point $(0, 0, 0)$ fails to be 
$(n-2)$-strongly non-isolated for $n=4.$  The variety $V$ has non-trivial fundamental group in every punctured neighborhood of $(0, 0, 0),$ and hence ${\rm spt} \, T_{V}$ near $(0, 0, 0)$ is not homeomorphic to a 4-dimensional disk. Note that in this example, the full singular set of $T_{V}$ is a $2$-dimensional subspace consisting of density 2 branch points. This indicates that the requirement that $Z$  be $(n-2)$-strongly non-isolated is optimal in that it cannot be relaxed to simply requiring that the  density $q$ branch set near $Z$ has full projection onto an $(n-2)$-dimensional ball. (See Section~\ref{topology-intro}).

We refer the reader to the introduction to \cite{KrumWica} for a comprehensive description of the methods in our overall program, as well as a detailed discussion of related works including the work of Almgren (\cite{Almgren}), and more details on the ways in which our approach contrasts, conceptually and technically, with \cite{Almgren}. To elaborate a little more on the aspect mentioned above regarding the use of center manifolds, the primary conceptual difference and the reason for the gains in simplicity in our approach is, as mentioned above, that it establishes uniform a priori decay estimates 
at ${\mathcal H}^{n-2}$ a.e.\ branch point of $T$  essentially at the same time as bounding the size (and analysing finer properties) of ${\rm sing} \, T$. This can be viewed, in principle, as aligning the higher codimension theory  closer with the codimension 1 regularity theory. For codimension 1 area minimizers of dimension~$n$, De~Giorgi's uniform decay estimate (\cite{DG}), applicable at points with planar tangent cones (which can be assumed to have multiplicity 1 without loss of generality), is the first step of the proof of Hausdorff $n-7$ dimensionality of the singular set. This basic approach  extends (albeit with a significantly more elaborate mechanism to implement it, on account of higher multiplicity) to the more general setting of stationary codimension 1 integral varifolds with stable regular part and no classical singularities (\cite{Wic14}). In that context again establishing a uniform decay estimate applicable at points with a planar tangent cone (now possibly with constant integer multiplicity $\geq 2$) is the  first step towards the same $n-7$ bound for the Hausdorff dimension of the singular set. In the codimension 1 setting of \cite{Wic14}, it is ultimately shown with the help of this uniform estimate and the earlier Schoen--Simon theory \cite{SS} (see also the more recent work of Bellettini (\cite{Bel})) that points with planar tangent cones are in fact regular points, ruling out branch points. In various other codimension 1 settings where classical singularities and hence branch points are a priori allowed, including the setting of mod $p$ minimizing hypersurfaces for even $p$, progress has very recently been achieved---starting with \cite{MW}---via this ``decay-estimates-first'' principle  (see \cite{DHMSS}, \cite{KMW}, \cite{MS26}, \cite{MX}).

Establishing a similar a priori decay estimate at branch points of area minimizers of higher codimension---as done in our program---has required an approach that is 
still different from the above codimension~1 settings. This is in part because the corresponding ``linear theory'' is, of necessity, weaker in higher codimension. Among the new ideas we rely on for this purpose are an intrinsic frequency function for $T$---the \emph{planar frequency function} introduced in \cite{KrumWica}---and certain novel ways of  using this frequency function. The key properties of the planar frequency function are its approximate monotonicity subject to decay of $T$ and that it takes values $\leq 1$ on stationary cones. Our first use of the planar frequency function is to establish  that branch points $Z$
where $T$ \emph{does not} decay towards a tangent plane rapidly (i.e.\ at a faster rate than $|X - Z|^{1 + \alpha}$ for some fixed small $\alpha \in (0, 1)$), if any such exist, behave \emph{in a certain quantitative manner} like singularities where no tangent cone is planar (\cite[Theorem~1.1]{KrumWica}). 
This property enables further asymptotic analysis of $T$ about ${\mathcal H}^{n-2}$ a.e.\ point where rapid decay of $T$ to a plane fails (\cite{KrumWicb}), which in turn leads to several key results: uniqueness of tangent cones to $T$ at ${\mathcal H}^{n-2}$ a.e.\ point; locally uniform decay estimates at ${\mathcal H}^{n-2}$ a.e.\ branch point of a given density; that slow-decay branch points, if any exist, form an ${\mathcal H}^{n-2}$-null set; and that the set of all non-branch-point singularities of a given density is locally $(n-2)$-rectifiable (with locally finite ${\mathcal H}^{n-2}$ measure). We can summarize these  results (from~\cite{KrumWica} and~\cite{KrumWicb}) as follows:

\begin{theoremX}\label{thm:X} {\rm (\cite[Theorem~1.1]{KrumWicb})}
Let $T$ be an $n$-dimensional locally area-minimizing rectifiable current in an open set $U \subset \mathbb{R}^{n+m}$, and let $\op{sing}T$ denote the interior singular set of $T$.  
\begin{enumerate}[itemsep=2mm,topsep=0mm]
	\item[{\rm (a)}]  For $\mathcal{H}^{n-2}$-a.e.~point $Z \in {\rm sing} \, T$, the current $T$ has a unique tangent cone $\mathbf{C}_Z$ at $Z$ of the form $\mathbf{C}_Z = \sum_{j=1}^p q_j \llbracket P_j \rrbracket$ where $p,q_1,\ldots,q_p$ are integers $\geq 1$ and $P_1,\ldots,P_p$ are distinct $n$-dimensional planes (all depending on $Z$).  Either $p = 1$ (i.e.~$\mathbf{C}_Z$ is supported on a single plane) or there is an $(n-2)$-dimensional subspace $L$ with $P_i \cap P_j = L$ for every $i \neq j$. 
	\item[{\rm (b)}]  We can express $\op{sing} T = \mathcal{S} \cup \mathcal{B}$ where:
\begin{enumerate}[itemsep=2mm,topsep=2mm]
	\item[{\rm (i)}]  $\mathcal{S} \cap \mathcal{B} = \emptyset$
	\item[{\rm (ii)}]  $\mathcal{S}$ is countably $(n-2)$-rectifiable, ${\mathcal S}_{q} = {\mathcal S} \cap \{\Theta(T, Z) = q\}$ has locally finite $(n-2)$-dimensional Hausdorff measure for each $q$, and for ${\mathcal H}^{n-2}$ a.e.\ point $Z \in {\mathcal S}$, the (unique) tangent cone 
	${\mathbf C}_{Z}$ is of the type as in {\rm (a)} above with $p \geq 2$;
	\item[{\rm (iii)}]  for each $Z \in \mathcal{B}$, the current $T$ has a unique tangent cone $\mathbf{C}_Z$ at $Z$ which is supported on an $n$-dimensional plane $P_Z$;
	\item[{\rm (iv)}]  for every compact set $K \subset U$ there exist constants $\mu_K \in (0,1)$ and $C_K \in (0,\infty)$ {\rm (}~which depend only on $n$, $m$ and $\sup \, \{\Theta(T, Z) \, : \, Z \in K\}$ {\rm )} such that the following holds true: for every $Z_0 \in \mathcal{B} \cap K$ there exists $\sigma_{Z_0} > 0$ such that for every $Z \in \mathcal{B} \cap \mathbf{B}_{\sigma_{Z_0}}(Z_0)$ with $\Theta(T,Z) = \Theta(T,Z_0)$ we have that   
		\begin{align*} 
			&\hspace{20mm} \op{dist}_{\mathcal H}^2(P_Z \cap \mathbf{B}_1(0), P_{Z_0} \cap \mathbf{B}_1(0)) 
			\\&\hspace{35mm} 
            \leq C_K \bigg(\frac{|Z-Z_0|}{\sigma_{Z_0}}\bigg)^{2\mu_K} 
				\sigma_{Z_0}^{-n-2} \int_{\mathbf{B}_{\sigma_{Z_0}}(Z_0)} \op{dist}^2(X,Z_0+P_{Z_0}) \,d\|T\|(X)  
		\end{align*}
		and 
		\begin{align*}
			&\hspace{5mm} \sigma^{-n-2} \int_{\mathbf{B}_{\sigma}(Z)} \op{dist}^2(X,Z+P_Z) \,d\|T\|(X) 
			\\&\hspace{20mm} \leq C_K \Big(\frac{\sigma}{\rho}\Big)^{2\mu_K} \rho^{-n-2} \int_{\mathbf{B}_{\rho}(Z)} \op{dist}^2(X,Z+P_Z) \,d\|T\|(X)  
		\end{align*}
		for all $\sigma,\rho$ with $0 < \sigma \leq \rho \leq \sigma_{Z_0}$.
\end{enumerate}
\end{enumerate}
\end{theoremX}

The decay estimates thus obtained and a second key use of the planar frequency monotonicity formula form the basis of the present paper, the contents of which include an elementary, standard argument for bounding the Hausdorff dimension of the set ${\mathcal B}^{(<2)}$ of branch points with sub-quadratic decay, followed by a detailed asymptotic analysis of $T$ about ${\mathcal H}^{n-2}$ a.e.\ point in ${\mathcal B}^{(< 2)}$. This analysis 
establishes, for ${\mathcal H}^{n-2}$ a.e.\ point in ${\mathcal B}^{(< 2)},$ the existence of a unique non-zero tangent function (blow-up) of $T$ relative to the (unique) tangent plane at that point (as described in items (iii) \& (iv) at the beginning of this section),  and as corollaries of that, local structural decomposition of ${\mathcal B}^{(<2)}$ as a finite union of disjoint, locally $(n-2)$-rectifiable pieces   (as in item (ii)), as well as topological conclusions for $T$ on approach to certain points in this set (as in item (v)).  These are all established without the need for center manifolds.
  
In fact  the arguments in the present article (for Euclidean ambient spaces) apply to the larger set ${\mathcal B}^{(\neq 2)}$ consisting of all branch points excluding those with \emph{precisely} quadratic decay, whence the main conclusions including the dimension bound for ${\mathcal B}^{(\neq 2)}$, uniqueness of tangent functions to $T$ relative to tangent planes at ${\mathcal H}^{n-2}$ a.e.\ point in ${\mathcal B}^{(\neq 2)}$, and the rectifiability conclusions for ${\mathcal B}^{(\neq 2)}$ are all obtained avoiding center manifolds. Additionally, with regard to item (v), the local topological $n$-disk conclusion near a branch point also follows without the need for center manifolds provided we take, everywhere in the statement in (v), the planar frequency ${\mathcal N}_{T, {\rm Pl}}(\cdot)$ in place of the branching order ${\mathcal O}_{T}(\cdot).$

Once these results are established, what remains is to analyze the dimension and structure of the set 
${\mathcal B}^{(=2)}$ of branch points of quadratic decay, and to extend the asymptotic  conclusions for $T$ near points in ${\mathcal B}^{(=2)}$ and the topological $n$-disk conclusion for $T$ near a point $Z$ as in item (v) and  having integer planar frequency $\geq 2.$  This is done in  \cite{KrumWicd} where we finally use Almgren's center manifold,  now armed with the crucial advantage of being able to capitalise on the fact that when there is quadratic or faster decay for $T$ at a density $q$ branch point $Z$, there is a \emph{single,} essentially unique center manifold that contains $Z$ and all nearby density $q$ branch points with quadratic or faster decay. The work  
in \cite{KrumWicd} in fact analyses the larger set ${\mathcal B}^{(\geq 2)}$ consisting of branch points of quadratic or faster decay; this requires no more effort than for just the set ${\mathcal B}^{(=2)},$ and moreover, proceeding this way is needed when the ambient space is a general Riemannian manifold 
(see \cite{KrumWicd}).   

For the rest of this introduction we shall focus on the results in the present article and in \cite{KrumWicd}, along with background directly relevant to them.

\subsection{The branching order ${\mathcal O}_{T}(Z)$} \label{branching-order-intro} 
We shall call the number ${\mathcal O}_{T}(Z)$ as in items (iii) \& (iv) of Section~\ref{summary-intro} above the \emph{branching order} of $T$ at the branch point $Z$. By our Theorem~\ref{thm:B} below, the branching order at $Z$ exists and is uniquely determined by $Z$ for ${\mathcal H}^{n-2}$ a.e.\ branch point $Z$, and it is the planar frequency of $T$ at $Z$ (cf.\ \cite[Lemma~5.1]{KrumWica}) if the planar frequency at $Z$ is $< 2$, or else, it is the ``Almgren frequency'' of $T$ at $Z$ relative to (an essentially unique)  center manifold for $T$ that contains $Z$ and all nearby branch points of the same density and quadratic or faster decay to the tangent plane. In the latter case, by \cite{Almgren}, there is indeed such a center manifold ${\mathcal M}_{Z},$ which is the graph of a $C^{3, \alpha}$ function $\psi_{Z}$ defined on a neighborhood of $Z$ in $Z + P_{Z},$ with the property that ${\mathcal M}_{Z}$ contains $Z$ and all nearby branch points $Z^{\prime}$ with the same density as $Z$ and with $T$ decaying at $Z^{\prime}$ quadratically or faster towards the unique tangent plane $Z^{\prime} + P_{Z^{\prime}}.$ Moreover, 
by Theorem~\ref{thm:B}(v), $\psi_{Z}$ is uniquely determined by $Z$ up to an additive term of order ${\rm o}(|X - Z|^{{\mathcal O}_{T}(Z)}).$ By a key result established in \cite{Almgren}, there exists a 
$q$-valued, Lipschitz normal map $N_{Z} \, : \, {\mathcal M}_{Z} \to {\mathcal A}_{q}({\mathcal M}_{Z}^{\perp})$ whose graph provides a good approximation for $T$ at \emph{all} sufficiently small scales about $Z$, and the usual (Almgren) frequency function associated with $N_{Z}$ is approximately 
monotone (non-decreasing) as a function of the scale $\sigma \in (0, \rho_{Z})$ for some $\rho_{Z} >0$. By ``Almgren frequency of $T$ relative to a center manifold about $Z$'' we mean the limiting value of the frequency function of $N_{Z}$ as $\sigma \to 0^{+}.$ Thus when the planar frequency at $Z$ is 
$\geq 2$, the branching order ${\mathcal O}_{T}(Z)$ is the limiting value of the frequency function of $N_{Z}$ as $\sigma \to 0^{+}.$   (Theorem~\ref{thm:B} below).  In view of the uniqueness of $\psi_{Z}$ up to 
perturbations of order ${\rm o}(|X - Z|^{{\mathcal O}_{T}(Z)}),$ it is clear that ${\mathcal O}_{T}(Z)$ indeed does not depend on the center manifold. 

In item (iii) of Section~\ref{summary-intro}, by ``the rate at which points of $T$ come together on approach to $Z$'' we mean either the rate at which the average $L^{2}$ height of $T$ at scale $\sigma >0$ about 
$Z$ relative to the tangent plane $Z + P_{Z}$ decays towards zero as $\sigma \to 0^{+}$, or the rate at which the average $L^{2}$ norm of $N_{Z}$ at scale $\sigma$ about $Z$ 
decays to zero as $\sigma \to 0^{+}$, depending on whether the planar frequency at $Z$ is $<2$ or not. Note that in either case, the approximate monotonicity of the relevant frequency function  (i.e.\ the planar frequency function of $T$ relative to $P_{Z}$ or the Almgren frequency function of $N_{Z}$) alone is not enough to establish (as asserted in (iii)) that this rate is precisely $O(\sigma^{{\mathcal O}_{T}(Z)})$.  While it is standard by the doubling properties implied by frequency function monotonicity that the decay rate is between 
$O(\sigma^{{\mathcal O}_{T}(Z)+ \epsilon(\sigma)})$~and~$O(\sigma^{{\mathcal O}_{T}(Z)})$ where $\epsilon(\sigma) \to 0$ as 
$\sigma \to 0^{+}$, the fact that this rate is precisely equal to $O(\sigma^{{\mathcal O}_{T}(Z)})$ is a substantially more subtle conclusion. In the present context, we obtain this more precise conclusion as a consequence of the assertion in item (iv) of Section~\ref{summary-intro}, i.e.\ the existence of a unique, non-zero, homogeneous tangent function with degree of homogeneity 
${\mathcal O}_{T}(Z)$. This is established in the present paper for the case ${\mathcal O}_{T}(Z) < 2$ (i.e.\ when the planar frequency $< 2$) and in 
\cite{KrumWicd} for the case ${\mathcal O}_{T}(Z) \geq 2$.

\subsection{The two-dimensional case: work of Chang and Micallef--White} \label{2d-case-intro}
If the dimension $n=2$, by the combined work of White (\cite{Whi83}), Chang (\cite{Cha88}) and Micallef--White (\cite{MicWhi95}), a strikingly complete asymptotic description of $T$ valid near any singularity has long been known.  In fact when $n=2$, it is easily seen a priori that every tangent cone at a singular point is supported on either a single plane or two or more planes intersecting at a point. Capitalising on this fact and via techniques special to two dimensions, it was shown in \cite{Whi83} that tangent cones are unique and in \cite{Cha88} that singularities are isolated. Moreover, if $T$ is locally irreducible near a singular point $Z$ (i.e.\ $T$ cannot be decomposed as a sum of two non-zero currents without boundary and with distinct supports, or more precisely, there is no $\rho>0$ such that $T \res {\mathbf B}_{\rho}(Z) = T_{1} + T_{2}$ with $T_{i}$ a non-zero $n$-dimensional rectifiable current with 
$\partial \, T_{i} \res {\mathbf B}_{\rho}(Z) = 0$ for $i=1, 2$ and with $Z \in {\rm spt} \, T_{1} \cap {\rm spt} \, T_{2}$), then $Z$ is a branch point with a unique tangent plane of some multiplicity $q \geq 2$, and near that point, $T$ is a parameterized disk (\cite{Cha88}, \cite{MicWhi95}). More precisely, assuming (without loss of generality, in view of irreducibility) that the multiplicity of $T$ on its regular set is $1$ near $Z$  and that $Z = 0$, the combined works \cite{Cha88} and \cite{MicWhi95} imply that after mapping by an ambient rotation, ${\rm spt} \, T$ near $Z$ can be 
 parameterized by a map of the form $$z \mapsto (z^{q}, f(z)), \;\; z \in B_{1}(0) \subset {\mathbb C} \approx {\mathbb R}^{2} \times \{0\},$$ where $f$ is a $C^{2}$ function taking values in  ${\mathbb R}^{m}$ with $D^{j}f(z) = O(|z|^{q+1-j})$  for $j=0, 1, 2;$ moreover, for any fixed $q^{\rm th}$ root of unity $\nu \neq 1$,  
 $$f(\nu z) - f(z) = {\rm Re} (az^{p}) + e(z)$$ where 
 \begin{itemize}
 \item $e$ (depending on $\nu$) satisfies $D^{j}e(z) = o(|z|^{p + 1-j})$  for $j=0, 1$; 
 \item $a \in {\mathbb C}^{m} \setminus \{0\}$ is a constant (depending on $\nu$) 
 satisfying $a \cdot a  = 0$ 
 and 
 \item $p$ is an integer (independent of $\nu$) with $p> q$ and $p$ not divisible by $q$ 
 \end{itemize}
 (\cite[Theorems~A~\&~B]{MicWhi95}). 
 In the notation of items (iii) and (iv) in Section~\ref{summary-intro}, this says that
 $${\mathcal O}_{T}(Z) = p/q$$ 
and that the unique ${\mathcal A}_{q}({\mathbb R}^{m})$ valued tangent function $\varphi_{Z}$ of $T$ at $Z$ is given by 
 $$\varphi_{Z}(z) = \sum_{j=0}^{q - 1}\llbracket {\rm Re} \, (a  \, \nu^{j} z^{p/q}) \rrbracket.$$ 
 This remarkable result says that the current near $Z$, modulo its average height relative to the tangent plane at $Z$, is a $C^{1}$ perturbation of its (unique, non-zero) tangent function $\varphi_{Z}$, thus providing a complete picture of the behaviour of a 2-dimensional area minimizing current near a branch point. In other words, when $n=2$, we have long known that the above results (i)--(v) (in Subsection~\ref{summary-intro}) hold (in fact with (ii) in the stronger form that singualarities are isolated). Note in particular that the 
 condition in (v) that the branch point $Z$ is ``$(n-2)$-strongly non isolated'' (see Theorem~\ref{thm:A} below for the meaning of this) holds automatically at every branch point when $n=2$.

\subsection{General dimensions: toplogical disk property near branch points}  \label{topology-intro}
For area minimizers of dimension $n \geq 3$, to what extent a Chang--Micallef--White type topological description near a branch points can be given has remained a largely open question. Simple examples show that without further hypotheses, 
 the direct analogue of the Chang--Micallef--White conclusion (giving that $T$ near $Z$ is a $C^{1}$ branched disk) need not hold;  this is the case even when ${\mathcal O}_{T}(Z) = p/q$ for $p$ relatively prime to $q = \Theta(T, Z).$ 
This is demonstrated by  multiplicity 1 currents associated with certain 
 algebraic subvarieties of ${\mathbb C}^{n}$, such as the 4 real dimensional current $T_{V_{1}}$ associated with   
 $$V_{1} = \{(x, y,w) \, : \, w^{2} = x(x^{2} - y^{4})\} \subset {\mathbb C}^{3} \approx {\mathbb R}^{6}.$$ Note that $T_{V_{1}}$ has an isolated singularity at $Z_{0} = (0, 0, 0)$ which is a branch point where 
 the tangent plane is $2|\{w = 0\}|$ and ${\mathcal O}_{T_{V_{1}}}(Z_{0}) = 3/2;$ the (unique) tangent function of $T_{V_{1}}$ at $Z_{0}$ is cylindrical, and is (up to a non-zero constant) given by $\varphi_{Z_{0}}(x, y) = \llbracket {\rm Re} \, (x^{3/2}) \rrbracket + \llbracket -{\rm Re} \, (x^{3/2}) \rrbracket.$ An example with similar behaviour near a non-isolated branch point is provided by the current $T_{V_{2}}$ associated with 
 $$V_{2} = \{(x, y, w) \, : \, w^{2} = x^{3}(x^{2} - y^{3}) \} \subset {\mathbb C}^{3}.$$ 
This has ${\rm sing} \, T_{V_{2}} = \{(0, a, 0) \, : \, a \in {\mathbb C}\}$, and each $Z \in {\rm sing} \, T_{V_{2}}$ is a density 2 branch point with unique tangent plane 
$Z + 2|\{w=0\}|$; moreover, the unique tangent function of $T_{V_{2}}$ at $Z$ is cylindrical, and is (up to a non-zero constant) given by $\varphi_{Z_{0}}(x, y) = \llbracket {\rm Re} \, (x^{5/2}) \rrbracket + \llbracket -{\rm Re} \, (x^{5/2}) \rrbracket$ for $Z_{0} = (0, 0, 0)$ and  $\varphi_{Z}(x, y) = \llbracket {\rm Re} \, (ia^{3/2}x^{3/2}) \rrbracket + \llbracket -{\rm Re} \, (ia^{3/2}x^{3/2}) \rrbracket$ for $Z  = (0, a, 0)$, $a \neq 0$. In particular, 
${\mathcal O}_{T_{V_{2}}}(0) = 5/2$ and ${\mathcal O}_{T_{V_{2}}}(Z) = 3/2$ if $Z = (0, a, 0)$, $a \neq 0$. As can be verified by direct computation, both these varieties have non-trivial fundamental group in every punctured neighborhood of $(0, 0, 0)$, and thus neither of them is homeomorphic to a 4-disk in any neighborhood of the origin.

 A main result in the present paper and in \cite{KrumWicd} is that a direct analogue of the Chang and Micallef--White results holds in arbitrary dimension at any of the ${\mathcal H}^{n-2}$ a.e.\  density $q$ branch point $Z$ where ${\mathcal O}_{T}(Z)$ exists, provided ${\mathcal O}_{T}(Z) = p/q$ for $p$ relatively prime to $q$ and the point $Z$ is ``$(n-2)$-strongly non-isolated''  in the sense of the following definition:

  \begin{definition}[$(n-2)$-strong non-isolatedness of a branch point]\label{strongly-non-isolated} Let $T$ and ${\mathcal B}$ be as in Theorem~\ref{thm:X}. Fix an integer $q \geq 2$ and let $Z_{0} \in {\mathcal B} \cap \{Z \, : \, \Theta(T, Z) = q\}.$ For $\rho > 0,$ we say that $Z_{0}$ is \emph{$(n-2)$-strongly non-isolated in ${\mathbf B}_{\rho}(Z_{0})$} if there exists an $(n-2)$-dimensional linear subspace $L$ of ${\mathbb R}^{n+m}$ such that 
\begin{equation}\label{topology hyp2}
 L \cap \mathbf{B}_{\rho/2}(0) 
		\subseteq \pi_L(\{ Z - Z_0 : Z \in \mathcal{B} \cap \mathbf{B}_{\rho}(Z_0), \;\; \Theta \, (T, Z) = q \;\; \text{ and } {\mathcal O}_{T}(Z) \geq  {\mathcal O}_{T}(Z_{0})\}) , 
\end{equation}
where $\pi_L : \mathbb{R}^{n+m} \rightarrow L$ is the orthogonal projection map onto $L.$ 

We say that $Z_{0} \in {\mathcal B} \cap \{Z \, : \, \Theta \, (T, Z) = q\}$ is \emph{$(n-2)$-strongly non-isolated} if there is a sequence $\rho_{j} \to 0^{+}$ such that $Z_{0}$ is 
$(n-2)$-strongly non-isolated in 
${\mathbf B}_{\rho_{j}}(Z_{0})$ for each $j$.
\end{definition}

\begin{theoremA}[Theorem~\ref{topology thm} below and {\cite[Theorem~1.8]{KrumWicd}}]\label{thm:A} 
Let $T$ be an $n$-dimensional locally area-minimizing rectifiable current in an open subset of ${\mathbb R}^{n+m}$  and let
$Z_{0}$ be a branch point of $T$ where $T$ has a unique tangent plane $P_{Z_{0}}$ of some multiplicity $q$ and where ${\mathcal O}_{T}(Z_{0})$ exists (noting that these conditions hold for ${\mathcal H}^{n-2}$ a.e.\ branch point). If
\begin{enumerate}
\item [{\rm (a)}] 
${\mathcal O}_{T}(Z_{0}) = p/q$ for some integer $p$ co-prime to $q$ and 
\item[{\rm (b)}] $Z_{0}$ is $(n-2)$-strongly non-isolated (cf.\ Definition~\ref{strongly-non-isolated}) 
\end{enumerate}
then ${\rm spt} \, T$ near $Z_{0}$ is homeomorphic to an $n$-dimensional disk and admits a $C^{1, \mu}$ parameterization for some fixed $\mu = \mu(n, m, q) \in (0, 1)$.
 
More precisely (and more generally), the following holds: Suppose $Z_{0}$ is a branch point of $T$ where $T$ has a unique tangent plane $P_{Z_{0}}$ of multiplicity $q$ and ${\mathcal O}_{T}(Z_{0})$ exists, and suppose that $Z_{0}$ 
satisfies hypothesis {\rm (a)}. There exists $\rho_{0} \in (0, 1)$ (depending on $Z_{0}$) such that if additionally, in place of hypothesis {\rm (b)}, we have that: 
\begin{enumerate} 
\item [{\rm (b$^{\prime}$)}] $Z_{0}$ is \emph{$(n-2)$-strongly non-isolated in ${\mathbf B}_{\rho}(Z_{0})$} for some $\rho \in (0, \rho_{0}]$ (cf.\ Definition~\ref{strongly-non-isolated}),  
\end{enumerate}
then there exist a number $\mu = \mu(n, m, q) \in (0, 1)$; a rotation $\Gamma$ of $\mathbb{R}^{n+m}$ with $\Gamma(P_{Z_0}) = \mathbb{R}^n \times \{0\}$; functions $g \in C^{1,\mu/2}(B^{n-2}_8(0);\mathbb{R}^2)$, $h \in C^{1,\mu/2}(B^{n-2}_8(0);\mathbb{R}^m)$, $\psi \in C^{1,\mu}(B^2_1(0) \times B^{n-2}_1(0);\mathbb{R}^m)$ and $f \in C^{1,\mu/2}(B^2_1(0) \times B^{n-2}_1(0);\mathbb{R}^m)$ such that: 
\begin{enumerate}[itemsep=2mm,topsep=0mm]
	\item[{\rm (i)}]  $\op{sing} \, T \cap {\mathbf B}_{\rho/8}(Z_{0}) = \{Z \in {\mathcal B} \, : \, \Theta \, (T, Z) = q, \;\; {\mathcal O}_{T}(Z) = p/q\} \cap {\mathbf B}_{\rho/8}(Z_{0})$ and 
    \begin{equation*}
		\Gamma(\eta_{Z_0,\rho/64}(\op{sing} T)) \cap \mathbf{C}_8(0) 
		= \{ (g(y), y, h(y)) : y \in B^{n-2}_8(0) \} \cap \mathbf{C}_8(0) ; 
	\end{equation*}
	\item[{\rm (ii)}]  the support of $T$ near $Z_{0}$ is parameterized by
    \begin{align*}
		&\op{spt} \Gamma_{\#} \eta_{Z_0,\rho/64\#} T \cap \mathcal{D} \times \mathbb{R}^m 
		\\=\,& \{ (g(y) + z^q, y, f(z,y)) : z \in B^2_1(0),\, y \in B^{n-2}_1(0) \} , 
	\end{align*}
	where we identify $\mathbb{R}^2 \cong \mathbb{C}$ and $\mathcal{D} = \{ (g(y)+x,y) : x \in B^2_1(0),\, y \in B^{n-2}_1(0) \}$; the mapping $(z,y) \mapsto (g(y)+z^q,y,f(z,y))$ is injective, and hence $\op{spt} \Gamma_{\#} \eta_{Z_0,\rho/64\#} T \cap \mathcal{D} \times \mathbb{R}^m$ is homeomorphic to an $n$-disk; 

	\item[{\rm (iii)}]  there exists a unique tangent function $\varphi$ to $\Gamma_{\#} \eta_{Z_0,\rho/64\#} T$ at the origin which is given by $\varphi(x_1,x_2,y) = \op{Re}(c \,(x_1+ix_2)^{p/q})$ for each $(x_1,x_2,y) \in \mathbb{R}^n$ and for some $c \in \mathbb{C}^m \setminus \{0\}$ with $c \cdot c = 0$ and 
	\begin{equation*}
		f(z,y) = \psi(z^q, y) + \op{Re}(cz^p) + e(z,y) 
	\end{equation*} 
	for $z \in B^2_1(0)$ and $y \in B^{n-2}_1(0)$, where 
	\begin{gather*}
		|e(z,y)| \leq C \,|(z^q,y)|^{p/q+\mu/(2nq)} , \quad 
		|D_z e(z,y)| \leq C \,|(z^q,y)|^{p/q-1+\mu/(2nq)} |z|^{q-1} , \\
		|D_y e(z,y)| \leq C \,|(z^q,y)|^{p/q-1+\mu/(2nq)} , 
	\end{gather*}
	for some constant  $C = C(n,m,q,p) \in (0,\infty)$; 
	
	\item[{\rm (iv)}]  for any fixed $q^{\rm th}$ root of unity $\nu \neq 1$
	\begin{equation*}
		f(z,y) - f(\nu z,y) = \op{Re}(az^p) + \widetilde{e}(z,y) 
	\end{equation*} 
	for $z \in B^2_1(0)$ and $y \in B^{n-2}_1(0)$, where $a \in \mathbb{C}^m \setminus \{0\}$ with $a \cdot a = 0$ and 
	\begin{gather*}
		|\widetilde{e}(z,y)| \leq C \,|(z^q,y)|^{p/q+\mu/(2nq)} , \quad 
		|D_z \widetilde{e}(z,y)| \leq C \,|(z^q,y)|^{p/q-1+\mu/(2nq)} |z|^{q-1} , \\
		|D_y \widetilde{e}(z,y)| \leq C \,|(z^q,y)|^{p/q-1+\mu/(2nq)} , 
	\end{gather*}
	where $C = C(n,m,q,p) \in (0,\infty)$ is a constant.
\end{enumerate}
\end{theoremA}

\begin{remark}\label{hyp-a}
{\rm (i) If ${\mathcal O}_{T}(Z) = 1 + q^{-1}$ (where $q = \Theta \, (T, Z)$), then hypothesis (a) of the theorem is automatically satisfied and hypothesis (b) takes the weaker form that \emph{near $Z$ the set of all density $q$ branch points $Z^{\prime}$  projects fully onto an open ball of an $(n-2)$-dimensional affine subspace}. Example $T_{V_{2}}$ discussed above shows that in general, even when hypothesis (a) holds, hypothesis (b) cannot be replaced by this weaker statement. Example $T_{V_{1}}$ shows that this weaker assumption cannot be dropped when ${\mathcal O}_{T}(Z) = 1 + q^{-1}$.

\noindent
(ii) Hypothesis (b) (without hypothesis (a)) a priori guarantees that one of the (non-zero, $q$-valued, average free, homogeneous of degree ${\mathcal O}_{T}(Z)$) tangent functions $\varphi_{Z}$ to $T$ at $Z$ is \emph{cylindrical} (i.e.\ invariant under translation along an $(n-2)$-dimensional subspace). 
Thus given hypothesis (b) and given such a tangent function $\varphi_{Z}$, hypothesis (a) is equivalent to the fact that  $\varphi_{Z}$ occurs with multiplicity 1 and is irreducible, i.e.\ that 
in the notation of Theorem~\ref{thm:B}{\rm (i)} below, $N_{Z} = 1$ and $m_{1}^{(Z)} = 1$. Of course then it follows that if both hypotheses (a) and (b) hold, $T$ is irreducible and has multiplicity 1 near $Z$.

\noindent 
(iii) Assume, conversely, that $T$ is irreducible and has multiplicity 1 near $Z$ (noting that there is no loss of generality in doing so). 
Whenever hypothesis~(b) holds, it follows from the analysis here and in \cite{KrumWicd} that $\varphi_{Z}$ is unique and is irreducible, i.e.\ that $N_{Z} = 1.$ If additionally $m_{1}^{(Z)} = 1$ then the conclusions of Theorem~\ref{thm:A} hold. What can be said if 
$m_{1}^{(Z)} \geq 2$ remains an important question to which we will return in future work.}
\end{remark}

\subsection{Uniqueness of tangent functions and the structure of the singular set} \label{uniqueness+structure-intro} 
Theorem~\ref{thm:A} follows directly from Theorem~\ref{topology thm} below and its counterpart \cite[Theorem~1.8]{KrumWicd}, which provide sufficient conditions under which an area minimizing current lying weakly close to a multiplicity $q$ plane in a cylinder (over the plane) is, in a smaller cylinder, represented by a multiplicity 1, $q$-valued $C^{1, \mu}$ graph over the plane. These conditions include: (I) that the current is close, in an appropriate sense, to a homogeneous, cylindrical $q$-valued Dirichlet energy minimizer with degree of homogeneity $p/q$ for $p, q$ relatively prime, and (II) that the ``full projection property'' of hypothesis~(b) of Theorem~\ref{thm:A} holds for the appropriate part of the branch set (cf.\ assumption \eqref{topology hyp2} of Theorem~\ref{topology thm}).  

Theorem~\ref{topology thm} in turn is obtained by applying iteratively the main excess-decay lemma of 
the present paper, Lemma~\ref{main excess decay lemma1}; and similarly, \cite[Theorem~1.8]{KrumWicd} results from iteratively applying 
\cite[Lemma~9.1]{KrumWicd} which is the counterpart, in \cite{KrumWicd}, to Lemma~\ref{main excess decay lemma1} here.  These two decay  results---Lemma~\ref{main excess decay lemma1} below and \cite[Lemma~9.1]{KrumWicd}---are on the other hand quite general and do not require either of the conditions (I) or (II) just mentioned, and consequently they have other important consequences. 
Specifically, in the absence of conditions (I) and (II), these decay results---in fact their respective strengthened versions Theorem~\ref{main excess decay thm} below and \cite[Theroem~9.3]{KrumWicd}---can still be iteratively applied to deduce: uniqueness of tangent functions at 
${\mathcal H}^{n-2}$ a.e.\ branch point; essential uniqueness of the center manifold at a branch point $Z$ with quadratic decay towards the tangent plane; the existence of a well-defined rate of decay of $T$ at a branch point $Z$ towards the center manifold, i.e.\ the branching order ${\mathcal O}_{T}(Z)$, determined uniquely by $Z$ (and in particular independently of the center manifold); and structure results for the branch set and consequently (by combining with the  structure 
results for the set of non-branch-point singularities established in \cite{KrumWicb}) structure results for the entire singular set ${\rm sing} \, T$. 
These  results are given in Theorem~\ref{thm:B} and Theorem~\ref{thm:C} below. 

\begin{theoremB}[Theorem~\ref{unique tangent function thm} below \& {\cite[Theorem~1.7]{KrumWicd}}]\label{thm:B} 
Let $T$ be an $n$-dimensional locally area-minimizing rectifiable current of an open subset $U \subset {\mathbb R}^{n+m}$ with 
$(\partial \, T) \res U = 0.$ 
For $\mathcal{H}^{n-2}$-a.e.~branch point $Z$ of $T$:
\begin{enumerate}[itemsep=3mm, topsep=0mm]
	\item[{\rm (a)}]  the density $q = \Theta(T,Z)$ is an integer $\geq 2$ and $T$ has a unique tangent cone at $Z$ which is supported on an $n$-dimensional plane $P_Z$;
	
	\item[{\rm (b)}]  there exists a radius $\rho_Z > 0$, branching order $\mathcal{O}_T(Z) \in (1,\infty)$, a $C^3$ single-valued function $\psi^{(Z)} : P_Z \rightarrow P_Z^{\perp}$, and a (unique) non-zero, homogeneous degree $\mathcal{O}_T(Z)$, locally Dirichlet energy minimizing $q$-valued function $\varphi^{(Z)} : P_Z \rightarrow \mathcal{A}_q(P_Z^{\perp})$ such that:
	\begin{enumerate}[itemsep=3mm,topsep=0mm]
		\item[{\rm (i)}] $\mathcal{O}_T(Z) = \ell_Z/q_Z$ for some relatively prime positive integers $\ell_Z,q_Z$ for which $q_Z \leq q$;

		\item[{\rm (ii)}] $\varphi^{(Z)}$ has the form 
		\begin{equation*}
			\varphi^{(Z)} = \sum_{j=1}^{N_Z} m^{(Z)}_j \varphi^{(Z)}_j 
		\end{equation*}
		where $N_Z,m^{(Z)}_1,\ldots,m^{(Z)}_{N_Z}$ are positive integers and $\varphi^{(Z)}_j : P_Z \rightarrow \mathcal{A}_{q^{(Z)}_j}(P_Z^{\perp})$ are distinct functions such that $\sum_{j=1}^{N_Z} m^{(Z)}_j q^{(Z)}_j = q$ and either $q^{(Z)}_j = 1$ and $\varphi^{(Z)}_j(X) = 0$ on $P_Z$ or $q^{(Z)}_j = q_Z$ and $\varphi^{(Z)}_j(X) = \op{Re}(c^{(Z)}_j (x_1+ix_2)^{\ell_Z/q_Z})$ on $P_Z$ for some $c^{(Z)}_j \in \mathbb{C}^m \setminus \{0\}$ (where we use coordinates $X = (x_1,x_2,\ldots,x_n)$ with respect to an orthonormal basis for $P_Z$); 
	
		\item[{\rm (iii)}]  $\varphi^{(Z)}$ is the unique tangent function of $T$ at $Z$ which is constructed as a blow-up of $T$ relative to the tangent plane $P_Z$ if $\mathcal{O}_T(Z) < 2$ and is constructed as a blow-up of $T$ relative to an Almgren  center manifold if $\mathcal{O}_T(Z) \geq 2$;

		\item[{\rm (iv)}]  $T$ is asymptotic to $Z+\op{spt}\op{graph}(\psi^{(Z)}+\varphi^{(Z)})$ in the sense that 
		\begin{align*}
			\hspace{23mm}&\frac{1}{\omega_n \rho^n} \int_{\mathbf{C}_{\rho}(Z,P_Z)} 
				\op{dist}^2(X, Z+\op{spt}\op{graph}(\psi^{(Z)}+\varphi^{(Z)})) \,d\|T\|(X) 
			\\&+ \sum_{j=1}^{N_Z} \frac{m^{(Z)}_j}{\omega_n \rho^n} 
				\int_{B_{\rho}(0) \setminus B^2_{\rho/16}(0) \times \mathbb{R}^{n-2}} \sum_{l=1}^{q^{(Z)}_j} 
				\op{dist}^2(Z + (X,\psi^{(Z)}(X)+\varphi^{(Z)}_{j,l}(X)), \op{spt} T) \,d\mathcal{H}^n(X) 
			\\&\leq C_Z \rho^{2\ell_Z/q_Z + 2\mu_Z} 
		\end{align*}
		for all $\rho \in (0,\rho_Z]$, where we let $\varphi^{(Z)}_j(X) = \sum_{l=1}^{q^{(Z)}_j} \llbracket \varphi^{(Z)}_{j,l}(X) \rrbracket$ and $(\psi^{(Z)} + \varphi^{(Z)})(X) = \sum_{j=1}^{N_Z} \sum_{l=1}^{q^{(Z)}_j} \llbracket \psi(X) + \varphi^{(Z)}_{j,l}(X) \rrbracket$ for each $X \in P_Z$ and where $\mu_Z \in (0,1)$ and $C_Z \in (0,\infty)$ are constants.

		\item[{\rm (v)}]  $\psi^{(Z)}$ is unique up to an additive term of order $o(\rho^{n+ \mathcal{O}_T(Z)})$ in $L^{2}(B_{\rho}(0))$; that is, if statement~{\rm (iv)} also holds true with a $C^3$ single-valued function $\widetilde{\psi}^{(Z)} : P_Z \rightarrow P_Z^{\perp}$ in place of $\psi^{(Z)}$, then 
		\begin{equation*}
			\rho^{-n/2}\|\widetilde{\psi}^{(Z)} - \psi^{(Z)}\|_{L^{2}(B_{\rho}(0))} = o(\rho^{\mathcal{O}_T(Z)}) \text{ as } \rho \rightarrow 0
		\end{equation*}
	\end{enumerate}
\end{enumerate}
\end{theoremB}

\begin{theoremC}[Theorem~\ref{B rectifiable thm} below, {\cite[Theorem~1.1]{KrumWicb}} \& {\cite[Theorem~1.6]{KrumWicd}}]\label{thm:C}
Let $T$ be an $n$-dimensional locally area-minimizing rectifiable current of an open subset
$U \subset {\mathbb R}^{n+m}$ with $(\partial \, T) \res U = 0$.  For every compact subset $K \subset U$, ${\rm sing} \, T \cap K$ is the union of a finite number of pairwise disjoint sets each of which is locally compact and locally $(n-2)$-rectifiable (and has in particular locally finite $(n-2)$-dimensional Hausdorff measure).
\end{theoremC} 

\subsection{Statements of the main theorems in the present article}  \label{main-present-paper-intro} Let $T$ be an $n$-dimensional locally area-minimizing rectifiable current of an open subset $U$ of $\mathbb{R}^{n+m}$.  Let $\mathcal{B}$ be as in Theorem~\ref{thm:X}.  Notice that for each point $Z \in \mathcal{B}$, the density $\Theta(T, Z)$ of $T$ at $Z$ is an integer.  Consequently, for each integer $q \geq 2$, we can consider the set $\mathcal{B} \cap \{ Z : \Theta(T,Z) = q \}$.   In light of Theorem~\ref{thm:X}(b)(iv) (and~\cite[Lemma~3.5]{KrumWica}), after restricting $T$ to a suitable neighborhood of a point of $\mathcal{B} \cap \{ Z : \Theta(T,Z) = q \}$ and then rotating and rescaling, we may assume that for some (fixed) constant $\mu \in (0,1)$ and a suitably small constant $\varepsilon_0 \in (0,1)$ the following hold true:

\begin{hypothesis}\label{introduction hyp}{\rm 
$T$ is an $n$-dimensional locally area-minimizing rectifiable current of $\mathbf{C}_1(0) = B_1(0) \times \mathbb{R}^m$ such that 
\begin{equation}\label{introduction hyp2}
	(\partial T) \llcorner \mathbf{C}_1(0) = 0, \quad \sup_{X \in \op{spt} T} \op{dist}(X,P_0) < \varepsilon_0, \quad 
		\pi_{\#} T = q \llbracket B_1(0) \rrbracket , 
\end{equation}
where $P_0 = \mathbb{R}^n \times \{0\}$ and $\pi : \mathbb{R}^{n+m} \rightarrow P_0$ is the orthogonal projection map onto $P_0$.  
}\end{hypothesis} 

\begin{hypothesis}\label{introduction hyp-1} 
For ${\mathcal H}^{n-2}$ a.e.\ branch point $Z \in \op{sing} T \cap \mathbf{C}_{1/2}(0)$ with $\Theta(T,Z) = q,$ $T$ has a unique tangent plane 
$P_Z$ at $Z$ such that 
\begin{equation}\label{introduction hyp4} 
    \sup_{\rho \in (0,1/8]} \rho^{-n-2-2\mu} \int_{\mathbf{B}_{2\rho}(Z)} \op{dist}^2(X,Z+P_Z) \,d\|T\|(X) < \infty .  
\end{equation}
\end{hypothesis} 

\begin{definition}\label{introduction B defn} 
Denote by $\mathcal{B}_q$ the set of all branch points $Z$ as in Hypothesis~\ref{introduction hyp-1}. 
\end{definition} 

In Definition~\ref{introduction B defn}, one can take $\mu = \mu_K$ as in Theorem~\ref{thm:X}(b)(iv) for a suitable choice of compact set $K \subset U$, and thereby guarantee that $\mathcal{H}^{n-2}$-a.e.~branch point $Z \in \mathbf{C}_{1/2}(0)$ of $T$ is an element of $\mathcal{B}_q$.  The assumption that \eqref{introduction hyp4} holds true for all $Z \in \mathcal{B}_q$ and for a fixed constant $\mu \in (0,1)$ (independent of $Z$) guarantees that the tangent plane $P_Z$ of $T$ at $Z$ is close to $P_0 = \mathbb{R}^n \times \{0\}$ and that the planar frequency function of $T$ at $Z$ is approximately monotone non-decreasing (see Lemma~\ref{uniform branch set lemma} and Theorem~\ref{frequency monotonicity thm}).  Consequently, for each $Z \in \mathcal{B}_q$ the planar frequency $\mathcal{N}_{T,{\rm Pl}}(Z) = \lim_{\rho\downarrow 0} N_{T,P_Z,Z}(\rho)$ exists.  Moreover, there is at least one tangent function $\varphi : P_Z \rightarrow \mathcal{A}_q(P_Z^{\perp})$ of $T$ at $Z$ constructed as a blow-up of $\eta_{Z,\rho_{\nu}\#} T$ relative to the tangent plane $P_Z$ for some sequence of radii $(\rho_{\nu})$ tending to zero.  Each such tangent function $\varphi$ is a non-zero, homogeneous degree $\mathcal{N}_{T,{\rm Pl}}(Z)$, locally Dirichlet energy minimizing $q$-valued function on $P_Z$.  Let 
\begin{equation*}
	\mathcal{B}^{(\neq 2)}_q = \{ X \in \mathcal{B}_q : \mathcal{N}_{T,{\rm Pl}}(Z) \neq 2 \} . 
\end{equation*}
By a standard argument based on the existence of tangent functions and upper semi-continuity of planar frequency, we derive a stratification result, Theorem~\ref{B prime stratification thm}, and in particular we obtain the optimal Hausdorff dimension bound for $\mathcal{B}^{(\neq 2)}_q$.  Using this together with~\cite[Theorem~1.1]{KrumWica} and~\cite[Theorem~1.5]{KrumWicd}, 
we recover the sharp Hausdorff dimension bound of the singular set of a locally area-minimizing rectifiable current with considerable simplification of Almgren's argument in~\cite{Almgren} (see~\cite[Corollary~1.3]{KrumWica} for further details).

\begin{theorem}\label{B prime dimension thm} 
Let $q \geq 2$ be an integer and $\mu \in (0,1)$.  There exists $\varepsilon_0 = \varepsilon_0(n,m,q,\mu) \in (0,1)$ such that if $T$ is as in Hypothesis~\ref{introduction hyp}, $\mathcal{B}_q$ is as in Definition~\ref{introduction B defn}, and $\mathcal{B}^{(\neq 2)}_q = \{ X \in \mathcal{B}_q : \mathcal{N}_{T,{\rm Pl}}(Z) \neq 2 \}$, then $\mathcal{B}^{(\neq 2)}_q$ has Hausdorff dimension at most $n-2$.
\end{theorem}

Our main results concerning the uniqueness of tangent functions at $\mathcal{H}^{n-2}$-a.e.~$Z \in \mathcal{B}^{(\neq 2)}_q$ and rectifiability of $\mathcal{B}^{(\neq 2)}_q$ are as follows.

\begin{theorem}\label{unique tangent function thm}
Let $q \geq 2$ be an integer.  There exists $\varepsilon_0 = \varepsilon_0(n,m,q) \in (0,1)$ such that the following holds true.  Let $T$ be as in Hypothesis~\ref{introduction hyp}.  
Then for $\mathcal{H}^{n-2}$-a.e.~branch point $Z \in \mathbf{C}_{1/2}(0)$ of $T$ with $\mathcal{N}_{T,{\rm Pl}}(Z) \neq 2,$ there exists a number $\rho_Z > 0$ and a (unique) homogeneous  Dirichlet energy minimizing $q$-valued function $\varphi^{(Z)} : P_Z \rightarrow \mathcal{A}_q(P_Z^{\perp})$ such that:
\begin{itemize}[itemsep=3mm,topsep=0mm]
	\item[{\rm (i)}] $\mathcal{N}_{T,{\rm Pl}}(Z) = \ell_Z/q_Z$ for some relatively prime positive integers $\ell_Z,q_Z$ for which $q_Z \leq q$;

	\item[{\rm (ii)}] $\varphi^{(Z)}$ has the form 
	\begin{equation*}
		\varphi^{(Z)} = \sum_{j=1}^{N_Z} m^{(Z)}_j \varphi^{(Z)}_j 
	\end{equation*}
	where $N_Z,m^{(Z)}_1,\ldots,m^{(Z)}_{N_Z}$ are positive integers and $\varphi^{(Z)}_j : \mathbb{R}^n \rightarrow \mathcal{A}_{q^{(Z)}_j}(\mathbb{R}^m)$ are distinct functions such that $\sum_{j=1}^{N_Z} m^{(Z)}_j q^{(Z)}_j = q$ and either $q^{(Z)}_j = 1$ and $\varphi^{(Z)}_j(X) = 0$ on $P_Z$ or $q^{(Z)}_j = q_Z$ and $\varphi^{(Z)}_j(X) = \op{Re}(c^{(Z)}_j (x_1+ix_2)^{\ell_Z/q_Z})$ on $P_Z$ for some $c^{(Z)}_j \in \mathbb{C}^m \setminus \{0\}$ (where we use coordinates $X = (x_1,x_2,\ldots,x_n)$ with respect to an orthonormal basis for $P_Z$); 
	
	\item[{\rm (iii)}]  $\varphi^{(Z)}$ is the unique tangent function of $T$ at $Z$ as constructed in Subsection~\ref{sec:B tanfun sec} (as a blow-up of $T$ relative to the tangent plane $P_Z$ in $\mathbf{C}_{\sigma}(0,P_Z)$ for all $\sigma > 0$); 
		
	\item[{\rm (iv)}]  $T$ is asymptotic to $Z+\op{spt}\op{graph}\varphi^{(Z)}$ in the sense that 
	\begin{align*}
		&\hspace{10mm} \frac{1}{\omega_{n}\rho^{n}}\int_{\mathbf{C}_{\rho}(Z,P_Z)} \op{dist}^2(X, Z+\op{spt}\op{graph}\varphi^{(Z)}) \,d\|T\|(X) 
		\\&\hspace{20mm} + \sum_{j=1}^{N_Z} \frac{m^{(Z)}_j}{\omega_{n} \rho^{n}} \int_{B_{\rho}(0) \setminus B^2_{\rho/16}(0) \times \mathbb{R}^{n-2}} 
			\sum_{l=1}^{q^{(Z)}_j} \op{dist}^2(Z + (x,y,\varphi^{(Z)}_{j,l}(x,y)), \op{spt} T) \,d\mathcal{H}^n(X) 
		\\&\hspace{30mm} \leq C_Z \rho^{2\ell_Z/q_Z + 2\mu_Z} 
	\end{align*}
	for all $\rho \in (0,\rho_Z]$, where $\varphi^{(Z)}_j(X) = \sum_{l=1}^{q^{(Z)}_j} \llbracket \varphi^{(Z)}_{j,l}(X) \rrbracket$ for each $X \in P_Z$ and where $\mu_Z \in (0,1)$ and $C_Z \in (0,\infty)$ are constants.
\end{itemize}
\end{theorem}

\begin{remark}{\rm 
It follows from the work of Micallef and White~\cite[Theorem~3.2]{MicWhi95} that the constants $c^{(Z)}_j$ in the conclusion of Theorem~\ref{unique tangent function thm} satisfy $c^{(Z)}_j \cdot c^{(Z)}_j = 0$.
}\end{remark}

\begin{theorem}\label{B rectifiable thm}
Let $q \geq 2$ be an integer and $\mu \in (0,1)$.  There exists $\varepsilon_0 = \varepsilon_0(n,m,q,\mu) \in (0,1)$ such that the following holds true.  Let $T$ be as in Hypothesis~\ref{introduction hyp}, $\mathcal{B}_q$ be as in Definition~\ref{introduction B defn} and $\mathcal{B}^{(< 2)}_q = \{ X \in \mathcal{B}_q : \mathcal{N}_{T,{\rm Pl}}(Z) < 2 \}$.  
For each compact set $K \subset \mathbf{C}_{1/2}(0)$ 
\begin{equation}\label{B rectifiable concl}
	\mathcal{B}^{(< 2)}_q \cap K = \bigcup_{j=1}^N L_j 
\end{equation}
for some positive integer $N$ and for some pairwise disjoint, locally compact, locally $(n-2)$-rectifiable sets $L_1,\ldots,L_N$ (and in particular each set $L_j$ has locally finite $\mathcal{H}^{n-2}$-measure in the sense that for each $X \in L_j$ there exists $\delta > 0$ such that $\mathcal{H}^{n-2}(L_j \cap \mathbf{B}_{\delta}(X)) < \infty$).
\end{theorem}

\begin{remark}{\rm 
Let $\mathcal{B}_q^{(=2)} = \{ Z \in \mathcal{B}_q : \mathcal{N}_{T,{\rm Pl}}(Z) = 2 \}$ and $\mathcal{B}_q^{(>2)} = \{ Z \in \mathcal{B}_q : \mathcal{N}_{T,{\rm Pl}}(Z) > 2 \}$.  By arguing as in the proof of Theorem~\ref{B rectifiable thm}, for each compact set $K \subset \mathbf{C}_{1/2}(0) \setminus \mathcal{B}_q^{(=2)}$ we can express $\mathcal{B}_q^{(> 2)} \cap K$ as a finite union of pairwise disjoint, locally compact, locally $(n-2)$-rectifiable sets (as in \eqref{B rectifiable concl} with $\mathcal{B}_q^{(>2)}$ in place of $\mathcal{B}_q^{(<2)}$).  The assumption that $K$ does not intersect $\mathcal{B}_q^{(=2)}$ is required since $\mathcal{B}_q^{(> 2)}$ is a locally compact subset of $\mathbf{C}_{1/2}(0) \setminus \mathcal{B}_q^{(=2)}$, but $\mathcal{B}_q^{(> 2)}$ is not necessarily a locally compact subset of $\mathbf{C}_{1/2}(0)$ since a sequence of points in $\mathcal{B}_q^{(> 2)}$ might converge to a point in $\mathcal{B}_q^{(=2)}$.
}\end{remark}

In addition to providing $\mathcal{H}^{n-2}$-a.e.~uniqueness of tangent functions, Theorem~\ref{unique tangent function thm} gives us important information about the topology and geometric structure of $T$ near certain branch points.

\begin{theorem}\label{topology thm}  
Let $\alpha = p/q$ where $p,q$ are relatively prime positive integers such that $p > q \geq 2$.  For each $\lambda \in (0,1/4)$ there exists $\varepsilon_0 = \varepsilon_0(n,m,q,\alpha,\lambda) \in (0,1)$ such that the following holds true.  Let $T$ be as in Hypothesis~\ref{introduction hyp} and $\mathcal{B}_q$ be as in Definition~\ref{introduction B defn}.  Suppose that $0 \in \mathcal{B}_q$ with $T$ tangent to $P_0 = \mathbb{R}^n \times \{0\}$ at the origin and 
\begin{equation}\label{topology hyp1}
	\mathcal{N}_{T,{\rm Pl}}(0) \geq \alpha, \quad N_{T,P_0,0}(1/2) < \alpha+\varepsilon_0 , 
\end{equation}
where $N_{T,P_0,0}$ is the planar frequency function of $T$ relative to the plane $P_{0}$ and with base point the origin (as in Definition~\ref{H D N defn}).  Suppose that the origin is $(n-2)$-strongly non-isolated in $B_{1/4}(0)$ in the sense that there exists an $(n-2)$-dimensional linear subspace $L \subset \mathbb{R}^{n+m}$ such that 
\begin{equation}\label{topology hyp2}
	L \cap \mathbf{B}_{\lambda}(0) \subset \pi_L(\{ X \in \mathcal{B}_q \cap \mathbf{B}_{1/4}(0) : \mathcal{N}_{T,{\rm Pl}}(X) \geq \alpha \}) 
\end{equation}
where $\pi_L : \mathbb{R}^{n+m} \rightarrow L$ is the orthogonal projection map onto $L$ (cf.~\cite[Definition~1.7]{KrumWica}).  Then near the origin, the support of $T$ is homeomorphic to an $n$-dimensional disk and admits a $C^{1, \mu}$ parameterization for some fixed $\mu = \mu(n,m,q) \in (0, 1/q]$. 

More precisely, after an orthogonal change of coordinates of $\mathbb{R}^{n+m}$ which fixes vectors in $\{0\} \times \mathbb{R}^m$, there exists $g \in C^{1,\mu/2}(B^{n-2}_{\lambda/4}(0),\mathbb{R}^2)$, $h \in C^{1,\mu/2}(B^{n-2}_{\lambda/4}(0),\mathbb{R}^m)$, and $f \in C^{1,\mu/2}(B^2_{(\lambda/32)^{1/q}}(0) \times B^{n-2}_{\lambda/32}(0), \mathbb{R}^m)$ such that:
\begin{enumerate}[itemsep=2mm,topsep=0mm]
	\item[{\rm (a)}]  $\op{sing} T \cap \mathbf{C}_{\lambda/4}(0) = \{ X \in \mathcal{B}_q : \mathcal{N}_{T,{\rm Pl}}(Z) = \alpha \} \cap \mathbf{C}_{\lambda/4}(0)$ and 
	\begin{equation*}%\label{topology concl1}
		\op{sing} T \cap \mathbf{C}_{\lambda/4}(0) = \{ (g(y),y,h(y)) : y \in B^{n-2}_{\lambda/4}(0) \} \cap \mathbf{C}_{\lambda/4}(0) ;
	\end{equation*}

	\item[{\rm (b)}]  the support of $T$ is parameterized by 
	\begin{equation*}%\label{topology concl2}
		\op{spt} T \cap \mathcal{D} \times \mathbb{R}^m 
		= \{ (g(y) + z^q, y,  f(z,y)) : z \in B^2_{(\lambda/32)^{1/q}}(0),\, y \in B^{n-2}_{\lambda/32}(0) \} , 
	\end{equation*}
	where we identify $\mathbb{R}^2 \cong \mathbb{C}$ and $\mathcal{D} = \{ (g(y)+x,y) : x \in B^2_{\lambda/32}(0),\, y \in B^{n-2}_{\lambda/32}(0) \}$; the mapping $(z,y) \mapsto (g(y)+z^q,y,f(z,y))$ is injective, and hence $\op{spt} T \cap \mathcal{D} \times \mathbb{R}^m$ is homeomorphic to an $n$-disk;
	
	\item[{\rm (c)}]  there exists a unique tangent function $\varphi^{(0)}$ of $T$ at the origin which is given by $\varphi^{(0)}(x_1,x_2,y) = \op{Re}(c \, (x_1+ ix_2)^{p/q})$ for each $(x_1,x_2,y) \in \mathbb{R}^n$ and for some $c \in \mathbb{C}^m \setminus \{0\}$ with $c \cdot c = 0$ and 
	\begin{equation}\label{topology concl3}
		f(z,y) = \op{Re}(c \, z^p) + e(z,y) 
	\end{equation}
	for each $z \in B^2_{(\lambda/32)^{1/q}}(0)$ and $y \in B^{n-2}_{\lambda/32}(0)$, where 
	\begin{gather}\label{topology concl4}
		|e(z,y)| \leq C \widehat{E} \,|(z^q,y)|^{p/q+\mu/(8nq)} , \quad 
		|D_z e(z,y)| \leq C \widehat{E} \,|(z^q,y)|^{p/q-1+\mu/(8nq)} |z|^{q-1}, \\ 
		|D_y e(z,y)| \leq C \widehat{E} \,|(z^q,y)|^{p/q-1+\mu/(8nq)} , \nonumber 
	\end{gather}
	where $\widehat{E} = \widehat{E}(T,P_0,\mathbf{C}_1(0))$ and $C = C(n,m,q,\alpha) \in (0,\infty)$ is a constant; 

	\item[{\rm (d)}]  for any fixed $q^{\rm th}$ root of unity $\nu \neq 1$  
	\begin{equation}\label{topology concl5}
		f(z,y) - f(\nu z,y) = \op{Re}(a\, z^p) + \widehat{e}(z,y) 
	\end{equation}
	for each $z \in B^2_{(\lambda/32)^{1/q}}(0)$ and $y \in B^{n-2}_{\lambda/32}(0)$, where $a \in \mathbb{C}^m \setminus \{0\}$ with $a \cdot a = 0$ and 
	\begin{gather}\label{topology concl6}
		|\widehat{e}(z,y)| \leq C \widehat{E} \,|(z^q,y)|^{p/q+\mu/(8nq)} , \quad 
		|D_z \widehat{e}(z,y)| \leq C \widehat{E} \,|(z^q,y)|^{p/q-1+\mu/(8nq)} |z|^{q-1}, \\ 
		|D_y \widehat{e}(z,y)| \leq C \widehat{E} \,|(z^q,y)|^{p/q-1+\mu/(8nq)} , \nonumber 
	\end{gather}
	where $C = C(n,m,q,\alpha) \in (0,\infty)$ is a constant.
\end{enumerate}
\end{theorem}

\begin{remark}\label{topology rmk}
{\rm 
It follows from the proof of Theorem~\ref{topology thm} that for each $Z \in \op{sing} T \cap \mathbf{C}_{\lambda/32}(0)$ there exists a rotation $Q_Z$ of $\mathbb{R}^{n+m}$ and functions $g_Z  \in C^{1,\mu/2}(B^{n-2}_{\lambda/4}(0),\mathbb{R}^2)$, $h_Z \in C^{1,\mu/2}(B^{n-2}_{\lambda/4}(0),\mathbb{R}^m)$, and $f_Z \in C^{1,\mu/2}(B^2_{(\lambda/32)^{1/q}}(0) \times B^{n-2}_{\lambda/32}(0), \mathbb{R}^m)$ such that conclusion~(a)--(d) of Theorem~\ref{topology thm} holds true with $Q_{Z\#} \eta_{Z,1\#} T,g_Z,h_Z,f_Z$ in place of $T,g,h,f$.  In particular, $T$ has a unique tangent function $\varphi^{(Z)}$ at $Z$ such that $Q_Z \varphi^{(Z)}(Q_Z^{-1} X) = \op{Re}(c_Z\, (x_1+ix_2)^{p/q})$ for each $X = (x_1,x_2,y) \in \mathbb{R}^n$, where $c_Z \in \mathbb{C}^m \setminus \{0\}$ with $c_Z \cdot c_Z = 0$.  For each $z \in B^2_{(\lambda/32)^{1/q}}(0)$ and $y \in B^{n-2}_{\lambda/32}(0)$
\begin{equation}\label{topology rmk eqn1}
	f_Z(z,y) = \op{Re}(c_Z\, z^p) + e_Z(z,y) , 
\end{equation}
where 
\begin{gather*}
	|e_Z(z,y)| \leq C \widehat{E} \,|(z^q,y)|^{p/q+\mu/(8nq)} , \quad 
	|D_z e_Z(z,y)| \leq C \widehat{E} \,|(z^q,y)|^{p/q-1+\mu/(8nq)} |z|^{q-1}, \\ 
	|D_y e_Z(z,y)| \leq C \widehat{E} \,|(z^q,y)|^{p/q-1+\mu/(8nq)} , \nonumber 
\end{gather*}
where $C = C(n,m,q,\alpha) \in (0,\infty)$ is a constant.  
To say more about the asymptotic behavior of the current on approach to branch points, we would need to know higher regularity of $\op{sing}\, T$, and of $Q_Z$ and $c_Z$ as functions of $Z$.  
In the special case $q = 2$ and $\alpha = 3/2$ (and we assume \ref{topology hyp2}),~\cite{Krum14} showed that the singular set of $T$ is an $(n-2)$-dimensional locally real analytic submanifold near the origin; otherwise, little is known.  
By a modification of the proof of Theorem~\ref{topology thm}, if near the origin we knew that $\op{sing} T$ is a $C^{2,\mu}$ $(n-2)$-dimensional submanifold and $Q_Z$ is a $C^{1,\mu}$-function of $Z \in \op{sing} T$, then we could parameterize the support of $T$ by 
\begin{equation*}
	(g(y),y,h(y)) + Q_{(g(y),y,h(y))}(z^q,0,\widetilde{f}(z,y))
\end{equation*}
for some function $\widetilde{f} \in C^{1,\mu/2}(B_{\sigma}(0),\mathbb{R}^m)$, where $\sigma > 0$ is suitably small (and for $g,h$ as in conclusion~(a) above), and for each $(z,y) \in \mathbb{C} \times \mathbb{R}^{n-2}$ near the origin 
\begin{equation}\label{topology rmk eqn2}
	\widetilde{f}(z,y) = \op{Re}(c_{(g(y),y,h(y))} \, z^p) + \widetilde{e}(z,y) 
\end{equation}
and 
\begin{equation*}
	|\widetilde{e}(z,y)| + |z| \,|D_z \widetilde{e}(z,y)| + |z|^q |D_y \widetilde{e}(z,y)| \leq C |z|^{p/q+\mu/8} E_{\varsigma/2}(T,\varphi^{(0)},\mathbf{C}_1(0)) , 
\end{equation*}
where $C = C(n,m,q,\alpha) \in (0,\infty)$ is a constant.
}\end{remark}

The proofs of Theorem~\ref{unique tangent function thm}, Theorem~\ref{B rectifiable thm} and Theorem~\ref{topology thm} involve a blow-up method based in part on the techniques originally developed in the seminal work of Simon (\cite{Sim93}) for multiplicity 1 classes, and extended to higher multiplicity settings in \cite{Wic14}, \cite{KrumWic1} and \cite{KrumWic2}. The novelty here is that we develop the blow-up method using the monotonicity formula for the planar frequency function (Theorem~\ref{frequency monotonicity thm}).  
The blueprint for our approach is provided by our previous work~\cite{KrumWic1} and~\cite{KrumWic2}, in which we established that Dirichlet energy minimizing multi-valued functions and two-valued minimal graphs are asymptotic to unique homogeneous harmonic multi-valued functions at branch point singularities.  
To briefly highlight some key new aspects in the present setting, let $T$ be an $n$-dimensional locally area-minimizing rectifiable current of $\mathbf{C}_1(0)$ which is weakly close to the plane $\mathbb{R}^n \times \{0\}$ with constant multiplicity $q$ and is significantly closer in $L^2$ to the graph of a homogeneous, cylindrical harmonic $q$-valued function $\varphi : \mathbb{R}^n \rightarrow \mathcal{A}_q(\mathbb{R}^m)$ (as in Hypothesis~$(\star)$ and Hypothesis~$(\star\star)$ in Subsection~\ref{sec:graph_rep_results}). 
As in~\cite{KrumWic1}, among other things we must control the tilt of a tangent plane at a branch points of $T$ relative to the fixed plane $\mathbb{R}^n \times \{0\}$.  In~\cite{KrumWic1}, $T$ was the area-stationary graph of a $C^{1,1/2}$ two-valued function $u$ over $\mathbb{R}^n$, and we controlled the tilt of the tangent planes by expressing $T$ as the graph of a two-valued function $\widetilde{u}$ over the tangent plane at the branch point and using an estimate---namely, \cite[Lemma~3.3]{KrumWic1}---comparing the average-free parts of $u$ and $\widetilde{u}$.  In Theorem~\ref{nonconcentration thm} of this paper, we proceed differently and bound the tilt of the tangent plane alongside the distance of ``good frequency'' branch point to the spine of 
$\varphi$. This latter distance estimate is one of the ``standard'' a priori estimates now appearing in numerous works, and was first proved in \cite{Sim93} in the  
multiplicity 1, non-degenerate setting and extended---in its more refined form appearing here---to a degenerate, arbitrary multiplicity setting in \cite{Wic14}.

Unlike in~\cite{KrumWic1} and~\cite{KrumWic2}, $T$ need not be the graph of a  $q$-valued function over $\mathbb{R}^n \times \{0\}$, and hence appropriate sheet-separation properties away from the spine of $\varphi$ are not readily available. 
We here establish (in Theorem~\ref{relharm separation3 thm}) a ``separation-into-lower-multiplicity-currents'' property via a novel uniform $L^2-L^{\infty}$ estimate for the distance of $T$ to the union of quantitatively pairwise disjoint graphs of single-valued (nearly) harmonic functions. This together with Almgren's strong Lipschitz approximation theorem is then used to establish an approximate graphical representation of $T$ relative to the graph of $\varphi$ (Theorem~\ref{graphical_thm}), which facilitates the basic blow-up procedure described in Section~\ref{sec:blowup subsec} for (a sequence of currents) $T$ relative to the graph of (a sequence of functions) $\varphi$. This procedure is central to the main decay estimate of the present work (Theorem~\ref{main excess decay thm}).

\subsection{The work of De~Lellis--Minter--Skorobogatova: a brief comparison}\label{comparison} 
In contemporaneous, independent  work \cite{DeLSko1}, \cite{DeLSko2}, \cite{DeLMinSko}
(available on the arXiv at arXiv:2304.11552, arXiv:2304.11555, arXiv:2304.11553 shortly after the initial posting of \cite{KrumWica} and \cite{KrumWicb} on the arXiv at arXiv:2304.10653 and arXiv:2304.10272), De~Lellis, Minter and Skorobogatova  also established two of the main results (i)--(v) in Section~\ref{summary-intro}, namely:
\begin{itemize}
\item[(a)] uniqueness of tangent cones  at ${\mathcal H}^{n-2}$ a.e.\ point (item (i), Section~\ref{summary-intro});  
\item[(b)] countable $(n-2)$-rectifiability of the singular set (part of item (ii), Section~\ref{summary-intro}).
\end{itemize}
Their approach is fundamentally different to ours: it is based on Almgren's framework (as presented in \cite{DeLSpa1, DeLSpa2, DeLSpa3}), taking the iterative construction of center manifolds for each branch point as the foundation. In particular, their proof of a.e.\ uniqueness of tangent cones requires center manifolds, and ours does not. 

Besides this, there is another fundamental difference between their program and ours, which is in how rectifiability properties of the singular set are established. Fixing an integer $q \geq 2$, our (stronger) rectifiability conclusions for the set of density $q$ singularities ${\rm sing}_{q} \, T$ are a corollary of an overarching asymptotic analysis of $T$ leading to ${\mathcal H}^{n-2}$ a.e.\ uniqueness of its tangent objects: tangent cones at singular points and, additionally,  tangent functions at branch points. The De~Lellis--Minter--Skorobogatova program on the other hand focuses on two different aspects of $T$ 
for different parts of ${\rm sing}_{q}\, T,$ and employs entirely different techniques for these different parts:

\begin{itemize}[itemsep=3mm,topsep=1mm]
\item[(i)] for one part $\Sigma \subset {\rm sing}_{q} \, T$, which contains the set of all non branch point singularities in ${\rm sing}_{q} \, T,$ 
their method---developed in \cite{DeLMinSko}---is similar to ours for the corresponding part ${\mathcal S}$, and is based on proving ${\mathcal H}^{n-2}$ a.e.\ uniqueness of tangent cones. This proceeds via an excess decay lemma in \cite{DeLMinSko} which has a counterpart in our work \cite{KrumWicb}. This decay lemma is the common indispensable ingredient in both programs, and in both settings, the proof of this lemma is based on \cite{Sim93} and \cite{Wic14}.\\

There is however an important difference between the two programs with regard to where and how this excess-decay lemma is applied: in their program, $\Sigma$ is first decomposed into a countable union, and a key hypothesis of the decay lemma is  verified for each piece in this decomposition. In our argument, the lemma is applied to the whole set ${\mathcal S}$  (which contains $\Sigma$, with the difference ${\mathcal S} \setminus \Sigma$ seen to be  ${\mathcal H}^{n-2}$-null a posteriori); we use the planar frequency monotonicity to show first that the hypotheses of the decay lemma are  satisfied \emph{locally uniformly} along ${\mathcal S}$. This leads  to a stronger structural conclusion: local-finiteness of the $(n-2)$-dimensional Hausforff measure of ${\mathcal S},$ in addition to contable rectifiability. This feeds into our final local structural decomposition result for ${\rm sing} \, T$ as in \ref{summary-intro}(ii), which is of course considerably stronger than the countable $(n-2)$-rectifiability of ${\rm sing} \, T$ obtained in the De~Lellis--Minter--Skorobogatova program.

\item[(ii)] the remaining part ${\rm sing}_{q} \, T \setminus \Sigma$ (which consists of branch points) is seen, subsequent to the analysis 
of $\Sigma$, to contain ${\mathcal H}^{n-2}$ a.e.\ branch point of ${\rm sing}_{q} \, T$. The analysis of  this part is carried out in \cite{DeLSko1, DeLSko2} and the methods used in this analysis are entirely different from the methods we employ for the analysis of ${\rm sing}_{q} \, T \setminus {\mathcal S}$. The method in \cite{DeLSko1, DeLSko2} is based on obtaining local upper Minkowski content bounds (though not uniform, as this set is first decomposed into a countable union, and upper Mikowski content bounds are obtained separately for each of the sets in this union), whereas our analysis of ${\rm sing}_{q} \, T \setminus {\mathcal S}$ is based on establishing uniqueness of tangent functions. Moreover, in their argument, constructing a center manifold sequence for each of the points in ${\rm sing}_{q} \, T  \setminus \Sigma$ is the technical foundation of the analysis. Our asymptotic analysis of ${\rm sing}_{q} \, T \setminus {\mathcal S}$ invokes the center manifold only in its canonical case of quadratic-decay branch points, where a single center manifold satisfying simplifying properties can be constructed. For all other branch points, instead of center manifold sequences, our framework uses the tangent plane;   intrinsic estimates  derived from the monotonicity formula for the planar frequency function provide the basis of our analysis.
\end{itemize}

To elaborate a little more on (ii), our focusing on uniqueness of the tangent functions at ${\mathcal H}^{n-2}$ a.e.\ point of ${\rm sing}_{q} \, T \setminus {\mathcal S}$ is motivated by the fact that the  tangent functions at ${\mathcal H}^{n-2}$ a.e.\ point in ${\rm sing}_{q} \, T \setminus {\mathcal S}$  have a simple, explicit form: they are translation invariant along an $(n-2)$-dimensional subspace and hence determined by multi-valued homogeneous locally Dirichlet energy minimizing functions on ${\mathbb R}^{2}$, which can be computed explicitly. Uniqueness results for tangent objects typically reveal significant information about the objects being studied; when an explicit classification of the tangent objects at singularities is a priori readily available, this provides an even more compelling reason to study their 
uniqueness as one would expect such  uniqueness to then reveal local structural/topological rigidity in addition to singular-set regularity. 
This is indeed the case with regard to branch points of area-minimizing currents, as established in the preset paper and in \cite{KrumWicd} (cf.\ Theorem~\ref{thm:A} above). 

We refer the reader to \cite[Section~2]{KrumWica} for a more comprehensive comparison between the two programs.

\subsection{Organization of the paper}  In Section~\ref{sec:preliminaries}, we provide notation for the paper and discuss some basic facts about locally area-minimizing rectifiable currents and Dirichlet energy minimizing multi-valued functions.  In Section~\ref{sec:planar freq sec}, we recall the monotonicity formula for planar frequency functions and its basic consequences from~\cite{KrumWica}.  Then we use the planar frequency function to construct tangent functions and establish the stratification result, Theorem~\ref{B prime stratification thm}, for $\mathcal{B}^{(\neq 2)}_q$.  As an immediate consequence we obtain the optimal Hausdorff dimension bound for $\mathcal{B}^{(\neq 2)}_q$ (Theorem~\ref{B prime dimension thm}).  

The rest of the paper is dedicated to the proofs of $\mathcal{H}^{n-2}$-a.e.~uniqueness of tangent functions at points of $\mathcal{B}^{(\neq 2)}_q$ (Theorem~\ref{unique tangent function thm}), local rectifiability of $\mathcal{B}^{(\neq 2)}_q$ (Theorem~\ref{B rectifiable thm}) and the topological disk property near certain branch points (Theorem~\ref{topology thm}).  In Section~\ref{sec:relharm}, we adapt arguments from Section~3 of~\cite{KrumWicb} to prove an $L^2-L^{\infty}$-estimate for the distance of a locally area-minimizing rectifiable current to a union of pairwise disjoint graphs of single-valued harmonic functions (Theorem~\ref{relharm separation3 thm}).  In Section~\ref{sec:graph rep}, we obtain graphical representation results (Theorem~\ref{graphical_thm}) and a blow-up procedure for locally area-minimizing rectifiable currents relative to graphs of homogeneous harmonic multi-valued functions.  In Section~\ref{sec:key_est_sec} and Section~\ref{sec:nonconcentration}, we derive a priori estimates for a locally area-minimizing rectifiable current $T$ which is $L^2$-close to the graph of a homogeneous harmonic multi-valued function $\varphi$, including a homogeneity estimate (Theorem~\ref{keyest thm}) which plays a fundamental role in our analysis and an estimate (Theorem~\ref{nonconcentration thm}) showing that the $L^2$-distance of $T$ to the graph of $\varphi$ cannot concentrate near the spine of $\varphi$.  Section~\ref{sec:fine blowups} and Section~\ref{sec:fine blowups decay} concern the analysis of fine blow-ups of locally area-minimizing rectifiable currents relative to graphs of homogeneous harmonic multi-valued functions, including the classification of homogeneous blow-ups (Theorem~\ref{classify in L thm}) and the asymptotic behavior of blow-ups (Theorem~\ref{blowup decay thm}).  In Sections~\ref{sec:fine excess decay}, we prove excess decay results for the $L^2$-distance of a locally area-minimizing rectifiable current to the graph of homogeneous harmonic multi-valued functions (Theorem~\ref{main excess decay thm}).  Finally, in Section~\ref{sec:proof of main results} we use the main excess decay result to complete the proofs of the main results of the paper, Theorems~\ref{unique tangent function thm}, \ref{B rectifiable thm} and \ref{topology thm}.

\section{Preliminaries}\label{sec:preliminaries} 

Throughout the paper, we shall adopt the same notation as~\cite{KrumWica} and~\cite{KrumWicb} unless otherwise stated.  See~\cite[Section~3]{KrumWica} for a discussion of general notation, as well as an overview of locally area-minimizing rectifiable currents and Dirichlet energy minimizing multi-valued functions.

We let $n \geq 1$ and $m \geq 2$ be fixed integers.  We often let $X = (x_1,x_2,\ldots,x_{n+m})$ denote a general point in $\mathbb{R}^{n+m}$.  $e_1,e_2,\ldots,e_{n+m}$ denotes the standard basis of $\mathbb{R}^{n+m}$.  Given an $n$-dimensional plane $P$ in $\mathbb{R}^{n+m}$, we let $\pi_P : \mathbb{R}^{n+m} \rightarrow P$ denote the orthogonal projection map onto $P$ and we let $\pi_{P^{\perp}} : \mathbb{R}^{n+m} \rightarrow P^{\perp}$ denote the orthogonal projection map onto the orthogonal complement $P^{\perp}$.  For each $Z \in \mathbb{R}^{n+m}$ and $\rho > 0$, we let  
\begin{align*}
	\mathbf{B}_{\rho}(Z) =\,& \{ X \in \mathbb{R}^{n+m} : |X - Z| < \rho \}, \\
	B_{\rho}(Z,P) =\,& \{ Z + X : X \in P \text{ with } |X| < \rho \}, \\
	\mathbf{C}_{\rho}(Z,P) =\,& \{ Z + X + Y : X \in P \text{ with } |X| < \rho \text{ and } Y \in P^{\perp} \} .
\end{align*}
We say that $P$ is an oriented $n$-dimensional plane if we equip $P$ with a simple $n$-vector $\vec P = \tau_1 \wedge \tau_2 \wedge\cdots\wedge \tau_n$, called the orientation of $P$, where $\tau_1, \tau_2, \ldots, \tau_n$ is an orthonormal basis for $P$.  We let $\llbracket P \rrbracket$ denote integral current with support $P$, multiplicity one, and orientation $\vec P$.

We will often let $P_0 = \mathbb{R}^n \times \{0\}$, and we equip $P_0$ with the orientation $\vec P_0 = e_1 \wedge e_2 \wedge\cdots\wedge e_n$.  In this case, we will identify $P_0 \cong \mathbb{R}^n$ and $P_0^{\perp} \cong \mathbb{R}^m$.  We let $\pi$ denote the orthogonal projection map onto $P_0$ and $\pi^{\perp}$ denote the orthogonal projection map onto $P_0^{\perp}$.  For each $z \in \mathbb{R}^n$ and $\rho > 0$, we will let $B_{\rho}(z) = \{ x \in \mathbb{R}^n : |x-z| < \rho \}$ and $\mathbf{C}_{\rho}(z) = B_{\rho}(z) \times \mathbb{R}^m$.

Let $Z \in \mathbb{R}^{n+m}$, $\rho > 0$, and $P$ be an $n$-dimensional plane in $\mathbb{R}^{n+m}$.  For each $n$-dimensional rectifiable current $T$ of $\mathbf{B}_{\rho}(Z)$, we define 
\begin{equation*}%\label{coarse height excess ball defn}
	\widehat{E}(T,P,\mathbf{B}_{\rho}(Z)) = \left( \frac{1}{\omega_n \rho^{n+2}} \int_{\mathbf{B}_{\rho}(Z)} \op{dist}^2(X,Z+P) \,d\|T\|(X) \right)^{1/2} .
\end{equation*}
Similarly, for each $n$-dimensional rectifiable current $T$ of $\mathbf{C}_{\rho}(Z,P)$ with  $\sup_{X \in \op{spt} T} \op{dist}(X,P) < \infty$, we define 
\begin{equation}\label{coarse height excess defn}
	\widehat{E}(T,P,\mathbf{C}_{\rho}(Z,P)) = \left( \frac{1}{\omega_n \rho^{n+2}} \int_{\mathbf{C}_{\rho}(Z,P)} \op{dist}^2(X,Z+P) \,d\|T\|(X) \right)^{1/2} .
\end{equation}
Next equip $P$ with an orientation $n$-vector $\vec P$.  For each $n$-dimensional rectifiable current $T$ of $\mathbf{C}_{\rho}(Z,P)$ with $(\partial T) \llcorner \mathbf{C}_{\rho}(Z,P) = 0$ and $\sup_{X \in \op{spt} T} \op{dist}(X,P) < \infty$, we let 
\begin{equation}\label{oriented tilt defn}
	\widehat{\mathcal E}(T,P,\mathbf{C}_{\rho}(Z,P))^2 = \frac{\|T\|(\mathbf{C}_{\rho}(Z,P))}{\omega_n\rho^n} 
		- \frac{\|\pi_{P\#} T\|(\mathbf{C}_{\rho}(Z,P))}{\omega_n\rho^n} . 
\end{equation}
In the case that $\pi_{\#} T = q \llbracket P \rrbracket$ for some integer $q \geq 0$, by~\cite[5.3.1]{Fed69} 
\begin{equation*}
	\widehat{\mathcal E}(T,P,\mathbf{C}_{\rho}(Z,P))^2 = \frac{1}{2\omega_n\rho^n} \int_{\mathbf{C}_{\rho}(Z,P)} |\vec T - \vec P|^2 \,d\|T\|(X) ,
\end{equation*}
where $\vec T$ is the orientation $n$-vector of $T$, and by \eqref{oriented tilt defn} 
\begin{equation}\label{oriented tilt and mass}
	\|T\|(\mathbf{C}_{\rho}(Z,P)) = (q + \widehat{\mathcal E}(T,P,\mathbf{C}_{\rho}(Z,P))^2) \,\omega_n \rho^n .
\end{equation}
Assuming $T$ is locally area-minimizing in $\mathbf{C}_{\rho}(Z,P)$, $\widehat{\mathcal E}(T,P,\mathbf{C}_{\rho}(Z,P)) = 0$ if and only if $T$ is a sum of $n$-dimensional planes which are parallel to $P$ and have orientation $\vec P$.

\begin{remark}{\rm 
The above notation differs from that of~\cite{KrumWica} and~\cite{KrumWicb}, which for instance denote the $L^2$-distance of $T$ to $Z+P$ in $\mathbf{C}_{\rho}(Z,P)$ by $E(T,P,\mathbf{C}_{\rho}(Z,P))$, rather than by $\widehat{E}(T,P,\mathbf{C}_{\rho}(Z,P))$ as in \eqref{coarse height excess defn}.  This is in order to contrast the coarse excess $\widehat{E}(T,P,\mathbf{C}_{\rho}(Z,P))$ of $T$ relative to a plane $P$ from the fine excess $E_{\varsigma,\varepsilon}(T,\varphi,\mathbf{C}_{\rho}(Z,P))$ of $T$ from a multi-valued harmonic function $\varphi : P \rightarrow \mathcal{A}_q(P^{\perp})$.
}\end{remark} 

Finally, we recall the following Lipschitz approximation theorem for locally area-minimizing rectifiable currents which was introduced by Almgren in~\cite{Almgren}.

\begin{theorem}[Almgren's strong approximation theorem]\label{lip approx thm}
For each $\gamma \in (0,1)$ there exists $\varepsilon = \varepsilon(n,m,q,\gamma) > 0$ such that the following holds true.  Let $x_0 \in \mathbb{R}^n$ and $\rho > 0$.  Let $T$ be an $n$-dimensional locally area-minimizing rectifiable current of $\mathbf{C}_{\rho}(x_0)$ such that 
\begin{equation*}
	(\partial T) \llcorner \mathbf{C}_{\rho}(x_0) = 0, \quad \sup_{X \in \op{spt} T} \op{dist}(X,P_0) < \infty, \quad 
	\pi_{\#} T = q \llbracket B_{\rho}(x_0) \rrbracket, \quad \mathcal{E} = \mathcal{E}(T,P_0,\mathbf{C}_{\rho}(x_0)) < \varepsilon .
\end{equation*}
Then there exists a Lipschitz $q$-valued function $u : B_{\gamma\rho}(x_0) \rightarrow \mathcal{A}_q(\mathbb{R}^m)$ and a closed set $K \subseteq B_{\gamma\rho}(x_0)$ such that 
\begin{gather*}%\label{lip approx concl}
	\op{Lip} u \leq C\mathcal{E}^{\varsigma}, \quad T \llcorner (K \times \mathbb{R}^m) = (\op{graph} u) \llcorner (K \times \mathbb{R}^m), \\
	\mathcal{L}^n(B_{\gamma\rho}(x_0) \setminus K) + \|T\|((B_{\gamma\rho}(x_0) \setminus K) \times \mathbb{R}^m) \leq C \mathcal{E}^{2+\varsigma} \rho^n ,
\end{gather*}
where $C= C(n,m,q,\gamma) \in (0,\infty)$ and $\varsigma = \varsigma(n,m,q) \in (0,1)$ are constants.
\end{theorem}

\begin{proof}
See~\cite[Corollary~3.29]{Almgren} or~\cite[Theorem~2.4]{DeLSpa1}.
\end{proof}

\section{Planar frequency function and Hausdorff dimension of the set of branch points of planar frequency $\neq 2$}\label{sec:planar freq sec}

Let $q \geq 2$ be an integer.  Let $T$ be as in Hypothesis~\ref{introduction hyp} and let $\mathcal{B}_q$ be as in Definition~\ref{introduction B defn}.  In this section, we want to apply the monotonicity formula for planar frequency~\cite[Theorem~4.4]{KrumWica} to study the structure of the set $\mathcal{B}^{(\neq 2)}_q$ of all branch points $Z \in \mathcal{B}_q$ at which $T$ has planar frequency $\neq 2$.  In particular, we establish that at each point $Z \in \mathcal{B}^{(\neq 2)}_q$ there exists at least one non-zero, homogeneous, Dirichlet energy minimizing function $\varphi : P_Z \rightarrow \mathcal{A}_q(P_Z^{\perp})$ constructed as a blow-up of the current $T$ at $Z$ relative to the tangent plane $P_Z$.  Moreover, using the well-known stratification theory of Almgren in~\cite{Almgren}, we prove Theorem~\ref{B prime dimension thm} of the Introduction, which asserts the optimal Hausdorff dimension bound for the set $\mathcal{B}^{(\neq 2)}$.

\subsection{Monotonicity formula for planar frequency functions at points of $\mathcal{B}_q$}\label{sec:planar freq subsec}  First, let's recall from~\cite{KrumWica} the definition of the planar frequency function as well as the monotonicity formula for planar frequency functions and its immediate consequences.

\begin{definition}[{\cite[Definition~4.1]{KrumWica}}]\label{H D N defn}
Let $Z \in \mathbb{R}^{n+m}$, $\rho_0 > 0$, and $P$ be an $n$-dimensional plane in $\mathbb{R}^{n+m}$.  Let $T$ be an $n$-dimensional locally area-minimizing rectifiable current of $\mathbf{C}_{\rho_0}(Z,P)$ such that 
\begin{equation*} 
	(\partial T) \llcorner \mathbf{C}_{\rho_0}(Z,P) = 0, \quad \sup_{X \in \op{spt} T \cap \mathbf{C}_{\rho_0}(Z,P)} \op{dist}(X,Z+P) < \infty . 
\end{equation*}
Let $\phi : [0,\infty) \rightarrow \mathbb{R}$ be the Lipschitz function defined by 
\begin{equation}\label{planar freq fun phi defn}
	\phi(s) = \begin{cases} 
		1 &\text{if } 0 \leq s < 1/2 \\
		2-2s &\text{if } 1/2 \leq s < 1 \\
		0 &\text{if } 1 \leq s < \infty .
	\end{cases}
\end{equation}
Recall from~\cite[Remark~4.2]{KrumWica} that we adopt the convention that $\phi'(1/2) = \phi'(1) = 0$.  We let 
\begin{align*}%\label{H defn} \label{D defn} 
	H_{T,P,Z}(\rho) =\,& -\rho^{1-n} \int \op{dist}^2(X,Z+P) \,|\nabla^S R|^2 \,\frac{1}{R} \,\phi'(R/\rho) \,d|T|(X,S), \\
	D_{T,P,Z}(\rho) =\,& \frac{1}{2} \rho^{2-n} \int \|\pi_S - \pi_P\|^2 \,\phi(R/\rho) \,d|T|(X,S) 
\end{align*}
for each $\rho \in (0,\rho_0]$, where $|T|$ denotes the $n$-dimensional integral varifold associated with $T$, $R(X) = |\pi_P(X-Z)|$, $\nabla^S$ denotes the gradient with respect to the plane $S$, and $\|\pi_S - \pi_P\|$ denotes the Frobenius norm of $\pi_S - \pi_P$.  The \emph{planar frequency function} $N_{T,P,Z}$ of $T$ at $Z$ relative to the tangent plane $P$ is given by 
\begin{equation*}%\label{N defn} 
	N_{T,P,Z}(\rho) = \frac{D_{T,P,Z}(\rho)}{H_{T,P,Z}(\rho)} 
\end{equation*}
for each $\rho \in (0,\rho_0]$ with $H_{T,P,Z}(\rho) > 0$.  (Note that by~\cite[Remark~4.6]{KrumWica}, $H_{T,P,Z}(\rho) = 0$ if and only if $\op{spt} T \cap \mathbf{C}_{\rho}(Z,P) \subset P$.)  Assuming $Z$ is a branch point of $T$ and $P_Z$ is the unique $n$-dimensional tangent plane to $T$ at $Z$, we define the \emph{planar frequency} of $T$ at $Z$ by 
\begin{equation*}%\label{planar freq defn} 
	\mathcal{N}_{T,{\rm Pl}}(Z) = \lim_{\rho\downarrow 0} N_{T,P_Z,Z}(\rho) 
\end{equation*}
whenever the limit exists.
\end{definition}

\begin{theorem}\label{frequency monotonicity thm}
Let $q \geq 2$ be an integer.  There exists $\delta = \delta(n,m,q) \in (0,1)$ and $\eta_0 = \eta_0(n,m,q) \in (0,1)$ such that the following holds true.  Let $Z \in \mathbb{R}^{n+m}$, $\rho_0 > 0$, and $P_Z$ be an $n$-dimensional plane in $\mathbb{R}^{n+m}$.  Let $T$ be an $n$-dimensional locally area-minimizing rectifiable current of $\mathbf{C}_{7\rho_0/4}(Z,P_Z)$ such that for some constants $\mu \in (0,1)$ and $\eta \in (0,\eta_0]$ 
\begin{gather}\label{freq mono hyp} 
	(\partial T) \llcorner \mathbf{C}_{7\rho_0/4}(Z,P_Z) = 0, \quad \Theta(T,Z) = q, \quad 
		\|T\|(\mathbf{C}_{7\rho_0/4}(Z,P_Z)) \leq (q+\delta) \,\omega_n (7\rho_0/4)^n , \\
	\widehat{E}(T,P_Z,\mathbf{C}_{7\rho/4}(Z,P_Z)) \leq \eta \Big(\frac{\rho}{\rho_0}\Big)^{\mu} \text{ for all } \rho \in (0,\rho_0] . \nonumber 
\end{gather}
Assume that $\op{spt} T \cap \mathbf{C}_{\rho}(Z,P_Z) \not\subset P_Z$ for some $\rho \in (0,\rho_0]$.  Then: 
\begin{enumerate}[itemsep=3mm,topsep=0mm]
	\item[(a)]  $\op{spt} T \cap \mathbf{C}_{\rho}(Z,P_Z) \not\subset P_Z$ for all $\rho \in (0,\rho_0]$;
	\item[(b)]  for all $0 < \sigma \leq \rho \leq \rho_0$
	\begin{equation}\label{freq mono concl1} 
		N_{T,P_Z,Z}(\sigma) \leq e^{C\eta^{\varsigma} (\rho/\rho_0)^{\mu\varsigma}} N_{T,P_Z,Z}(\rho) , 
	\end{equation}
	where $\varsigma = \varsigma(n,m,q) \in (0,1)$ is as in Theorem~\ref{lip approx thm} and $C = C(n,m,q,\mu) \in (0,\infty)$ is a constant;	
	\item[(c)]  the limit $\mathcal{N}_{T,{\rm Pl}}(Z) = \lim_{\rho\downarrow 0} N_{T,P_Z,Z}(\rho)$ exists and $\mathcal{N}_{T,{\rm Pl}}(Z) \geq 1+\mu$;
	\item[(d)]  for all $0 < \sigma \leq \rho \leq \rho_0$ 
	\begin{align*}
		e^{-C\eta^{2\varsigma} (\rho/\rho_0)^{2\mu\varsigma}/(2\mu\varsigma)} 
			&\Big(\frac{\sigma}{\rho}\Big)^{2e^{C\eta^{\varsigma} (\rho/\rho_0)^{\mu\varsigma}} N_{T,P_Z,Z}(\rho)} H_{T,P_Z,Z}(\rho) 
		\leq H_{T,P_Z,Z}(\sigma) 
		\\&\leq e^{C (\mathcal{N}_{T,{\rm Pl}}(Z)+1) \,\eta^{2\varsigma} (\rho/\rho_0)^{2\mu\varsigma}/(2\mu\varsigma)} 
			\Big(\frac{\sigma}{\rho}\Big)^{2\mathcal{N}_{T,{\rm Pl}}(Z)} H_{T,P_Z,Z}(\rho) , \nonumber 
	\end{align*}
	where $C = C(n,m,q) \in (0,\infty)$ are constants;
	\item[(e)]  for all $0 < \sigma \leq \rho \leq \rho_0$ 
	\begin{align*}
		&\frac{1}{\sqrt{2}} e^{-C\eta^{2\varsigma} (\rho/\rho_0)^{2\mu\varsigma}/(4\mu\varsigma)} 
			\Big(\frac{\sigma}{\rho}\Big)^{e^{C\eta^{\varsigma} (\rho/\rho_0)^{\mu\varsigma}} N_{T,P_Z,Z}(\rho)-1} \widehat{E}(T,P_Z,\mathbf{C}_{\rho}(Z,P_Z)) 
		\\&\hspace{10mm}\leq \widehat{E}(T,P_Z,\mathbf{C}_{\sigma}(Z,P_Z)) \nonumber
		\\&\hspace{20mm}\leq e^{C (\mathcal{N}_{T,{\rm Pl}}(Z)+2) \,\eta^{2\varsigma} (\rho/\rho_0)^{2\mu\varsigma}/(4\mu\varsigma)} 
			\Big(\frac{\sigma}{\rho}\Big)^{\mathcal{N}_{T,{\rm Pl}}(Z)-1} \widehat{E}(T,P_Z,\mathbf{C}_{\rho}(Z,P_Z)) , \nonumber 
	\end{align*}
	where $C = C(n,m,q) \in (0,\infty)$ are constants.
\end{enumerate}
\end{theorem}

\begin{proof}
By~\cite[Theorem~4.4]{KrumWica}, for every open interval $(a,b) \subset (0,\rho_0]$ such that $H_{T,P_Z,Z}(\rho) > 0$ for all $\rho \in (a,b)$, \eqref{freq mono concl1} holds true for all $a < \sigma < \rho < b$.  It follows using~\cite[Corollary~5.3(i)]{KrumWica} that $H_{T,P_Z,Z}(\rho) > 0$ for all $\rho \in (0,\rho_0]$.  Hence by~\cite[Remark~4.6]{KrumWica}, conclusion~(a) holds true.  Moreover, by~\cite[Theorem~4.4]{KrumWica}, conclusion~(b) holds true and the limit $\mathcal{N}_{T,{\rm Pl}}(Z) = \lim_{\rho\downarrow 0} N_{T,P_Z,Z}(\rho)$ exists.  By~\cite[Corollary~5.3]{KrumWica}, $\mathcal{N}_{T,{\rm Pl}}(Z_0) \geq 1+\mu$ (completing the proof of conclusion~(c)) and conclusions (d) and (e) hold true.
\end{proof}

\begin{corollary}\label{planar frequency cor}
Let $q \geq 2$ be an integer.  Let $Z \in \mathbb{R}^{n+m}$, $\sigma_0 > 0$, and $P_Z$ be an $n$-dimensional plane in $\mathbb{R}^{n+m}$.  Let $T$ be an $n$-dimensional locally area-minimizing rectifiable current of $\mathbf{C}_{7\sigma_0/4}(Z,P_Z)$ such that for some $\mu \in (0,1)$ 
\begin{equation}\label{freq cor hyp} 
	(\partial T) \llcorner \mathbf{C}_{7\sigma_0/4}(Z,P_Z) = 0, \quad \Theta(T,Z) = q , \quad
	\sup_{\rho \in (0,\sigma_0]} \rho^{-\mu} \widehat{E}(T,P_Z,\mathbf{C}_{7\rho/4}(Z,P_Z)) < \infty .  
\end{equation}
Then the planar frequency $\mathcal{N}_{T,{\rm Pl}}(Z) = \lim_{\rho\downarrow 0} N_{T,P_Z,Z}(\rho)$ exists and $\mathcal{N}_{T,{\rm Pl}}(Z) \geq 1+\mu$.  
\end{corollary}

\begin{proof}
Let $\delta = \delta(n,m,q) \in (0,1)$ and $\eta_0 = \eta_0(n,m,q) \in (0,1)$ be as in Theorem~\ref{frequency monotonicity thm}.  By \eqref{freq cor hyp}, for each $\eta \in (0,\eta_0]$ there exists $\rho_0 \in (0,\sigma_0/2]$ such that 
\begin{equation*}
	\|T\|(\mathbf{B}_{2\rho_0}(Z)) \leq (q+\delta/2) \,\omega_n (2\rho_0)^n, \quad 
	\widehat{E}(T,P_Z,\mathbf{C}_{7\rho/4}(Z,P_Z)) \leq \eta \Big(\frac{\rho}{\rho_0}\Big)^{\mu} \text{ for all } \rho \in (0,2\rho_0] .
\end{equation*}
By~\cite[Lemma~3.5]{KrumWica}, 
\begin{equation}\label{freq cor eqn1} 
	\sup_{X \in \op{spt} T \cap \mathbf{C}_{2\rho_0}(Z,P_Z)} \op{dist}(X,Z+P_Z) 
		\leq C\rho_0 \widehat{E}(T,P_Z,\mathbf{C}_{7\rho_0/2}(Z,P_Z)) \leq 2C\eta\rho_0 ,
\end{equation}
where $C = C(n,m) \in (0,\infty)$ is a constant.  Thus provided $\eta$ is small enough that $C\eta \leq 1/8$ and $(q+\delta/2)(7/4+C\eta)^n \leq (q+\delta)(7/4)^n$ (where $C$ as in \eqref{freq cor eqn1}), by the monotonicity formula for area 
\begin{align*}
	\|T\|(\mathbf{C}_{7\rho_0/4}(Z,P_Z)) \leq\,& \|T\|(\mathbf{B}_{7\rho_0/4 + 2C\eta\rho_0}(Z)) \\ 
	\leq\,& (q+\delta/2) \,\omega_n (7\rho_0/4 + 2C\eta\rho_0)^n \leq (q+\delta) \,\omega_n (7\rho_0/4)^n 
\end{align*}
(where again $C$ as in \eqref{freq cor eqn1}).  Therefore, \eqref{freq mono hyp} holds true.  Hence $N_{T,P_Z,Z}(\rho)$ is well-defined for all $\rho \in (0,\rho_0]$, and we can apply Theorem~\ref{frequency monotonicity thm}(c) to deduce that $\mathcal{N}_{T,{\rm Pl}}(Z)$ exists and $\mathcal{N}_{T,{\rm Pl}}(Z) \geq 1+\mu$.
\end{proof}

\subsection{Tangent functions}\label{sec:B tanfun sec}  Let $q \geq 2$ be an integer.  Let $Z \in \mathbb{R}^{n+m}$, $\sigma_0 > 0$, and $P_Z$ be an $n$-dimensional plane in $\mathbb{R}^{n+m}$.  Let $T$ be an $n$-dimensional locally area-minimizing rectifiable current of $\mathbf{C}_{7\sigma_0/4}(Z,P_Z)$ such that \eqref{freq cor hyp} true for some constant $\mu \in (0,1)$.   Let $\delta = \delta(n,m,q) \in (0,1)$ and $\eta_0 = \eta_0(n,m,q) \in (0,1)$ be as in Theorem~\ref{frequency monotonicity thm}.  Arguing as in the proof of Corollary~\ref{planar frequency cor}, for each $\eta \in (0,\eta_0]$ there exists $\rho_0 \in (0,\sigma_0]$ such that \eqref{freq mono hyp} holds true.  By~\cite[Lemma~5.1]{KrumWica}, upon equipping $P_Z$ with the orientation induced by $T$, $q \llbracket P_Z \rrbracket$ is the unique tangent cone to $T$ at $Z$ and 
\begin{equation*}
	\pi_{P_Z\#} (T \llcorner \mathbf{C}_{\rho_0}(Z,P_Z)) = q \llbracket B_{\rho_0}(Z,P_Z) \rrbracket .
\end{equation*}
Let $(\rho_{\nu})_{\nu=1}^{\infty}$ be a sequence of positive radii with $\rho_{\nu}\rightarrow 0^+$ as $\nu\rightarrow\infty$.  Set $\widehat{E}_{\nu} = \widehat{E}(T,P_Z,\mathbf{C}_{\rho_{\nu}}(Z,P_Z))$ for each positive integer $\nu$.  By Theorem~\ref{frequency monotonicity thm}(e) 
\begin{equation*}%\label{B tanfun eqn1}
	\widehat{E}(T,P_Z,\mathbf{C}_{4\sigma\rho_{\nu}}(Z,P_Z)) \leq C (4\sigma)^{e^{C\eta^{\varsigma}} N_{T,P_Z,Z}(\rho_0) - 1} \widehat{E}_{\nu}
\end{equation*}
for each $\sigma \in [1,\frac{\rho_0}{4\rho_{\nu}})$, where $C = C(n,m,q,\mu) \in (0,\infty)$ is a constant.  Thus following the blow-up procedure in~\cite[Subsection~6.1]{KrumWica}, for each $\sigma \in [1,\infty)$ we can construct a blow-up $\varphi_{\sigma}$ of $\eta_{Z,\rho_{\nu}\#} T$ relative to the tangent plane $P_Z$ and excess $\widehat{E}_{\nu}$ in $\mathbf{C}_{\sigma}(0,P_Z)$.  In other words, for each $\sigma \in [1,\infty)$ and each sufficiently large $\nu$, there exists a Lipschitz function $u_{\sigma,\nu} : B_{\sigma}(0,P_Z) \rightarrow \mathcal{A}_q(P_Z^{\perp})$ and a closed set $K_{\sigma,\nu} \subset B_{\sigma}(0,P_Z)$ such that 
\begin{align} 
	\label{B tanfun eqn2} &(\eta_{Z,\rho_{\nu}\#} T) \llcorner \pi_{P_Z}^{-1}(K_{\sigma,\nu}) 
		= (\op{graph} u_{\sigma,\nu}) \llcorner \pi_{P_Z}^{-1}(K_{\sigma,\nu}) , \\
	\label{B tanfun eqn3} &\mathcal{L}^n(B_{\sigma}(0,P_Z) \setminus K_{\sigma,\nu}) 
		+ \|\eta_{Z,\rho_{\nu}\#} T\|(\mathbf{C}_{\sigma}(0,P_Z) \setminus \pi_{P_Z}^{-1}(K_{\sigma,\nu})) 
		\\&\hspace{10mm} \leq C \widehat{E}(T,P_Z,\mathbf{C}_{4\sigma\rho_{\nu}}(Z,P_Z))^{2+2\varsigma} \sigma^n \nonumber 
		\leq C (4\sigma)^{n + (2\varsigma+2) (e^{C\eta^{\varsigma}} N_{T,P_Z,Z}(\rho_0) - 1)} \widehat{E}_{\nu}^{2+2\varsigma} , \nonumber \\
	\label{B tanfun eqn4} &\sup_{B_{\sigma}(0)} |u_{\sigma,\nu}| 
		\leq C \sigma \,\widehat{E}(T,P_Z,\mathbf{C}_{4\sigma\rho_{\nu}}(Z,P_Z)) 
		\leq C (4\sigma)^{e^{C\eta^{\varsigma}} N_{T,P_Z,Z}(\rho_0)} \widehat{E}_{\nu} \\
	\label{B tanfun eqn5} &\sup_{B_{\sigma}(0)} |Du_{\sigma,\nu}| 
		\leq C \widehat{E}(T,P_Z,\mathbf{C}_{4\sigma\rho_{\nu}}(Z,P_Z))^{2\varsigma} 
		\leq C (4\sigma)^{2\varsigma (e^{C\eta^{\varsigma}} N_{T,P_Z,Z}(\rho_0) - 1)} \widehat{E}_{\nu}^{2\varsigma} , 
\end{align}
where $\varsigma = \varsigma(n,m,q) \in (0,1)$ is as in Theorem~\ref{lip approx thm} and $C = C(n,m,q,\mu) \in (0,\infty)$ are constants.  It follows from \eqref{B tanfun eqn2}--\eqref{B tanfun eqn5} and~\cite[Lemma~3.7]{KrumWica} that 
\begin{equation*}
	\frac{1}{\sigma^2} \sup_{B_{\sigma}(0)} |u_{\sigma,\nu}|^2 + \frac{1}{\sigma^n} \int_{B_{\sigma}(0)} |Du_{\sigma,\nu}|^2
	\leq C  (4\sigma)^{2e^{C\eta^{\varsigma}} N_{T,P_Z,Z}(\rho_0) - 2} \widehat{E}_{\nu}^2 
\end{equation*}
for all sufficiently large $\nu$, where $C = C(n,m,q,\mu) \in (0,\infty)$ is a constant (see (6.8) in~\cite{KrumWica}).  Hence by~\cite[Proposition~2.11]{DeLSpaDirMin} and~\cite[Theorem~5.2]{DeLSpa1} (also see~\cite[Theorem~2.19]{Almgren}), for each $\sigma \in [1,\infty)$ we can find a subsequence $(\rho_{\nu_k})_{k=1}^{\infty}$ of $(\rho_{\nu})$ (depending on $\sigma$) and a locally Dirichlet energy minimizing $q$-valued function $\varphi_{\sigma} \in W^{1,2}(B_{\sigma}(0,P_Z),\mathcal{A}_q(P_Z^{\perp}))$ such that 
\begin{equation}\label{B tanfun eqn6}
	u_{\sigma,\nu_k}/\widehat{E}_{\nu_k} \rightarrow \varphi_{\sigma}
\end{equation} 
pointwise $\mathcal{L}^n$-a.e.~on $B_{\sigma}(0,P_Z)$ and strongly in $L^2(B_{\sigma}(0,P_Z),\mathcal{A}_q(P_Z^{\perp}))$ and 
\begin{equation}\label{B tanfun eqn7}
	|Du_{\sigma,\nu_k}|/\widehat{E}_{\nu_k} \rightarrow |D\varphi_{\sigma}|
\end{equation}
in $L^2(B_{\sigma}(0,P_Z))$ as $k\rightarrow\infty$.  

Notice that due to the choice of Lipschitz approximations $u_{\sigma,\nu}$, $\varphi_{\sigma}$ might depend on $\sigma$.  In Theorem~\ref{B tanfun thm}, we will show that after passing to a subsequence, there exists $\varphi \in W^{1,2}_{\rm loc}(P_Z,\mathcal{A}_qP_Z^{\perp}))$ (independent of $\sigma$) such that \eqref{B tanfun eqn6} and \eqref{B tanfun eqn7} hold true with $\varphi_{\sigma} = \varphi$ for all $\sigma \in [1,\infty)$.  We call such a function $\varphi$ a \emph{tangent function} of $T$ at $Z$.  Recall from~\cite[Subsection~3.5]{KrumWica} that the single-valued average $\varphi_a : P_Z \rightarrow P_Z^{\perp}$ of $\varphi$ is defined by $\varphi_a(X) = \frac{1}{q} \sum_{l=1}^q \varphi_l(X)$ for each $X \in P_Z$, where $\varphi(X) = \sum_{l=1}^q \llbracket \varphi_l(X) \rrbracket$.  Recall from~\cite[Subsection~3.6]{KrumWica} that we let 
\begin{equation}\label{blowup freq defn1}
    H_{\varphi,\xi}(\rho) = -\rho^{1-n} \int_{P_Z} |\varphi|^2 \,\frac{1}{R} \,\phi'(R/\rho), \quad 
    D_{\varphi,\xi}(\rho) = \rho^{2-n} \int_{P_Z} |D\varphi|^2 \,\phi(R/\rho) 
\end{equation}
for each $\rho > 0$, where $R(X) = |X-\xi|$ for each $X \in P_Z$ and $\phi$ is as in \eqref{planar freq fun phi defn}.  (More generally, $\phi : [0,\infty) \rightarrow \mathbb{R}$ can be any non-increasing Lipschitz function such that $\phi = 1$ on $[0,1/2]$ and $\phi = 0$ on $[1,\infty)$.  Letting $\phi$ increase to the characteristic function on the interval $[0,1)$ recovers the classical frequency function as defined in~\cite{Almgren}.)  The frequency function of $\varphi$ at $\xi \in P_Z$ is given by 
\begin{equation}\label{blowup freq defn2}
    N_{\varphi,\xi}(\rho) = \frac{D_{\varphi,\xi}(\rho)}{H_{\varphi,\xi}(\rho)}
\end{equation}
for each $\rho > 0$.  The frequency of $\varphi$ at $\xi$ is given by 
\begin{equation}\label{blowup freq defn3}
    \mathcal{N}_{\varphi}(\xi) = \lim_{\rho \downarrow 0} N_{\varphi,\xi}(\rho) .
\end{equation}
(More generally, if $P$ is an $n$-dimensional linear plane in $\mathbb{R}^{n+m}$ and $\varphi : \Omega \rightarrow \mathcal{A}_q(P^{\perp})$ is a locally Dirichlet energy minimizing $q$-valued function on an open subset $\Omega \subseteq P$, then for each $\xi \in \Omega$ and $\rho \in (0,\op{dist}(\xi,\partial\Omega))$ we can define $H_{\varphi,\xi}(\rho)$, $D_{\varphi,\xi}(\rho)$, $N_{\varphi,\xi}(\rho)$, and $\mathcal{N}_{\varphi}(\xi)$ as in \eqref{blowup freq defn1}, \eqref{blowup freq defn2}, and \eqref{blowup freq defn3} with $P$ in place of $P_Z$.)

\begin{theorem}\label{B tanfun thm}
Let $q \geq 2$ be an integer. 
Let $Z \in \mathbb{R}^{n+m}$, $\sigma_0 > 0$, and let $P_Z$ be an $n$-dimensional plane in $\mathbb{R}^{n+m}$.  Let $T$ be an $n$-dimensional locally area-minimizing rectifiable current of $\mathbf{C}_{7\sigma_0/4}(Z,P_Z)$ such that \eqref{freq cor hyp} holds true for some $\mu \in (0,1)$.  Equip $P_Z$ with the orientation induced by $T$ so that $q \llbracket P_Z \rrbracket$ is the unique tangent cone to $T$ at $Z$.
\begin{enumerate}[itemsep=3mm,topsep=0mm]
	\item[(i)]  For each sequence $(\rho_{\nu})_{\nu=1}^{\infty}$ of positive numbers such that $\rho_{\nu} \rightarrow 0$ as $\nu \rightarrow \infty$, there exists a subsequence $(\rho_{\nu_k})_{k=1}^{\infty}$ of $(\rho_{\nu})$ and a tangent function $\varphi \in W^{1,2}_{\rm loc}(P_Z,\mathcal{A}_q(P_Z^{\perp}))$ (independent of $\sigma$) such that $\varphi$ is the blow-up of $\eta_{Z,\rho_{\nu_k}\#} T$ relative to the tangent plane $P_Z$ and excess $\widehat{E}(T,P_Z,\mathbf{C}_{\rho_{\nu_k}}(Z,P_Z))$ in $\mathbf{C}_{\sigma}(0,P_Z)$ for every $\sigma \in [1,\infty)$. 
	\item[(ii)]  Each tangent function $\varphi$ of $T$ at $Z$ is non-zero, homogeneous degree $\mathcal{N}_{T,{\rm Pl}}(Z)$ and locally Dirichlet energy minimizing in $P_Z$.
\end{enumerate}
\end{theorem}

\begin{proof}
Let $\delta = \delta(n,m,q) \in (0,1)$ and $\eta_0 = \eta_0(n,m,q) \in (0,1)$ be as in Theorem~\ref{frequency monotonicity thm}.  It follows from \eqref{freq cor hyp} that for each $\eta \in (0,\eta_0]$ there exists $\rho_0 \in (0,\sigma_0]$ such that \eqref{freq mono hyp} holds true.  Let $(\rho_{\nu})$ be a sequence of positive numbers such that $\rho_{\nu} \rightarrow 0$ as $\nu \rightarrow \infty$.  Recall that for each $\sigma \in [1,\infty)$ and each sufficiently large $\nu$, there exists a Lipschitz function $u_{\sigma,\nu} : B_{\sigma}(0,P_Z) \rightarrow \mathcal{A}_q(P_Z^{\perp})$ and a closed set $K_{\sigma,\nu} \subset B_{\sigma}(0,P_Z)$ such that \eqref{B tanfun eqn2}--\eqref{B tanfun eqn5} holds true, and that for each $\sigma \in [1,\infty)$ there exists a subsequence $(\rho_{\nu_k})$ of $(\rho_{\nu})$ (depending on $\sigma$) and a locally Dirichlet energy minimizing $q$-valued function $\varphi_{\sigma} \in W^{1,2}(B_{\sigma}(0,P_Z),\mathcal{A}_q(P_Z^{\perp}))$ such that \eqref{B tanfun eqn6} and \eqref{B tanfun eqn7} hold true.

To see~(i), note that by \eqref{B tanfun eqn2} and \eqref{B tanfun eqn3} for each $1 \leq \sigma < \tau < \infty$ 
\begin{align}\label{B tanfun eqn8}
	\lim_{\nu \rightarrow \infty} \mathcal{L}^n(B_{\sigma}(0,P_Z) \cap \{u_{\sigma,\nu} \neq u_{\tau,\nu}\}) 
	\leq\,& \lim_{\nu \rightarrow \infty} \mathcal{L}^n(B_{\sigma}(0,P_Z) \setminus (K_{\sigma,\nu} \cap K_{\tau,\nu})) 
	\\ \leq\,& \lim_{\nu \rightarrow \infty} C (4\tau)^{n + (2\gamma+2) (e^{C\eta^{\gamma}} N_{T,P_Z,Z}(\rho_0) - 1)} \widehat{E}_{\nu}^{2+2\gamma}
	= 0, \nonumber 
\end{align}
where $C = C(n,m,q,\mu) \in (0,\infty)$ is a constant.  By a standard diagonalization argument, there exists a subsequence $(\rho_{\nu_k})$ of $(\rho_{\nu})$ such that for every $\sigma \in \{1,2,3,\ldots\}$ there exists a locally Dirichlet energy minimizing $q$-valued functions $\varphi_{\sigma} \in W^{1,2}_{\rm loc}(B_{\sigma}(0,P_Z),\mathcal{A}_q(P_Z^{\perp}))$ such that \eqref{B tanfun eqn6} and \eqref{B tanfun eqn7} hold true.  By \eqref{B tanfun eqn6} and \eqref{B tanfun eqn8}, $\varphi_{\sigma} = \varphi_{\tau}$ on $B_{\sigma}(0,P_Z)$ for all positive integers $\sigma < \tau$.  Thus we may define $\varphi \in W^{1,2}_{\rm loc}(\mathbb{R}^n,\mathcal{A}_q(\mathbb{R}^m))$ by $\varphi = \varphi_{\sigma}$ on $B_{\sigma}(0,P_Z)$ for all $\sigma \in \{1,2,3,\ldots\}$.  That way, for every $\sigma \in \{1,2,3,\ldots\}$, \eqref{B tanfun eqn6} and \eqref{B tanfun eqn7} hold true with $\varphi_{\sigma} = \varphi$.  For general $\sigma \in [1,\infty)$, choose an integer $\tau > \sigma$ so that using \eqref{B tanfun eqn8} and $u_{\tau,\nu_k}/\widehat{E}_{\nu_k} \rightarrow \varphi$ pointwise $\mathcal{L}^n$-a.e.~on $B_{\tau}(0,P_Z)$, we deduce that \eqref{B tanfun eqn6} and \eqref{B tanfun eqn7} hold true with $\varphi_{\sigma} = \varphi$.  

To see (ii), observe that by the area formula, \eqref{B tanfun eqn2}--\eqref{B tanfun eqn5}, and \eqref{B tanfun eqn6} 
\begin{align*}
	\frac{1}{\omega_n} \int_{B_1(0,P_Z)} |\varphi|^2 
	=\,& \lim_{k \rightarrow \infty} \frac{1}{\omega_n \widehat{E}_{\nu_k}^2} \int_{B_1(0,P_Z)} |u_{1,\nu_k}|^2 
	\\=\,& \lim_{k \rightarrow \infty} \frac{1}{\omega_n \rho_{\nu_k}^{n+2} \widehat{E}_{\nu_k}^2} 
		\int_{\mathbf{C}_{\rho_{\nu_k}}(Z,P_Z)} \op{dist}^2(X,P_Z) \,d\|T\|(X) = 1 
\end{align*}
and in particular $\varphi$ is not identically zero.  Since $\Theta(\eta_{Z,\rho_{\nu_k}\#} T, 0) = \Theta(T,Z) = q$, by~\cite[Lemma~6.4]{KrumWica} we have that $\varphi_a(0) = 0$.  Hence by~\cite[Lemma~6.5]{KrumWica}
\begin{equation*}
	N_{\varphi,0}(\sigma) 
	= \lim_{k \rightarrow \infty} N_{\eta_{Z,\rho_{\nu_k}\#} T,P_Z,0}(\sigma) 
	= \lim_{k \rightarrow \infty} N_{T,P_Z,Z}(\sigma\rho_{\nu_k})
	= \mathcal{N}_{T,{\rm Pl}}(Z) 
\end{equation*}
for each $\sigma \in (0,\infty)$.  Therefore, $\varphi$ is homogeneous degree $\mathcal{N}_{T,{\rm Pl}}(Z)$. 
\end{proof}

\begin{remark}{\rm 
Let $(\rho_{\nu})$ be a sequence of radii with $\rho_{\nu} \rightarrow 0^+$ as $\nu \rightarrow \infty$.  Let $(\rho_{\nu_k})$ be a subsequence of $(\rho_{\nu})$ and $\varphi \in W^{1,2}_{\rm loc}(P_Z,\mathcal{A}_q(P_Z^{\perp}))$ be a tangent function of $T$ at $Z$ such that for every $\sigma \in [1,\infty)$ we have that \eqref{B tanfun eqn6} and \eqref{B tanfun eqn7} hold true with $\varphi_{\sigma} = \varphi$ (as in Theorem~\ref{B tanfun thm}(i)).  Further assume (after passing to a subsequence) that $\widehat{E}_{\nu_k} < 2^{-k}$ for all $k$.  Recalling \eqref{B tanfun eqn3}, we can find a sequence of integers $\sigma_k$ for which $\sigma_k \rightarrow +\infty$ slowly enough that $\widetilde{u}_k = u_{\sigma_{\nu_k},\nu_k}$ and $\widetilde{K}_k = K_{\sigma_{\nu_k},\nu_k}$ satisfy 
\begin{equation*}%\label{B tanfun eqn9}
	\widetilde{u}_k/\widehat{E}_{\nu_k} \rightarrow \varphi 
\end{equation*} 
pointwise $\mathcal{L}^n$-a.e.~on $\mathbb{R}^n$ and strongly in $L^2(B_{\sigma}(0,P_Z),\mathcal{A}_q(P_Z^{\perp}))$ for each $\sigma \in [1,\infty)$, 
\begin{equation*}%\label{B tanfun eqn10}
	|D\widetilde{u}_k|/\widehat{E}_{\nu_k} \rightarrow |D\varphi|
\end{equation*}
strongly in $L^2(B_{\sigma}(0,P_Z))$ for each $\sigma \in [1,\infty)$, and 
\begin{equation*}%\label{B tanfun eqn11}
	\mathcal{L}^n(B_{\sigma_k}(0,P_Z) \setminus \widetilde{K}_k) < 2^{-k} .
\end{equation*}
}\end{remark}

\subsection{Stratification and Hausdorff dimension of $\mathcal{B}_q$}\label{sec:stratify sec}  Let $T$ be as in Hypothesis~\ref{introduction hyp} and let $\mathcal{B}_q$ be as in Definition~\ref{introduction B defn} from the Introduction.  Now we want to use the monotonicity formula for planar frequency functions and tangent functions to study the structure of $\mathcal{B}_q$.  First we will show via a modification~\cite[Theorem~8.1 and Theorem~8.2]{KrumWica} that there is a uniform choice $\rho_0 = 1/8$, $\delta = C\varepsilon^2$, and $\eta = C\varepsilon$ such that \eqref{freq mono hyp} holds true for every $Z \in \mathcal{B}_q$, where $C = C(n,m,q,\mu) \in (0,\infty)$ is a constant (which in particular depends on $\mu$ as in \eqref{introduction hyp4}).

\begin{lemma}\label{uniform branch set lemma} 
Let $q \geq 2$ be an integer and let $\mu \in (0,1)$.  There exists $\varepsilon_0 = \varepsilon_0(n,m,q,\mu) \in (0,1)$ such that if $T$ is as in Hypothesis~\ref{introduction hyp} and $\mathcal{B}_q$ is as in Definition~\ref{introduction B defn}, then for all $Z \in \mathcal{B}_q$ 
\begin{gather}
	\label{uniform branch set concl1} \op{spt} T \cap \overline{\mathbf{C}_{1/4}(Z,P_Z)} \subset \mathbf{C}_1(0), \quad 
		(\partial T) \llcorner \mathbf{C}_{1/4}(Z,P_Z) = 0 , \\
	\label{uniform branch set concl2} (\pi_{P_Z\#} T) \llcorner \mathbf{C}_{1/4}(Z,P_Z) = q \llbracket B_{1/4}(Z,P_Z) \rrbracket , \\
	\label{uniform branch set concl3} |\vec P_Z - \vec P_0| \leq C\widehat{E}(T,P_0,\mathbf{C}_1(0)) , \\
	\label{uniform branch set concl4} \|T\|(\mathbf{C}_{7/32}(Z,P_Z)) \leq (q + C\widehat{E}(T,P_0,\mathbf{C}_1(0))^2) \,\omega_n (7/32)^n , \\
	\label{uniform branch set concl5} \widehat{E}(T,P_Z,\mathbf{C}_{7\rho/4}(Z,P_Z)) \leq C\rho^{\mu} \widehat{E}(T,P_0,\mathbf{C}_1(0)) 
		\text{ for all } \rho \in (0,1/8] , 
\end{gather}
where $P_0 = \mathbb{R}^n \times \{0\}$, $P_Z$ is equipped with the orientation $\vec P_Z$ induced by $T$, and $C = C(n,m,q,\mu) \in (0,\infty)$ is a constant.
\end{lemma}

\begin{proof} 
We begin by arguing as in the proof of~\cite[Theorem~8.1]{KrumWica}.  Let $\widetilde{\mu} = \widetilde{\mu}(n,m,q,\mu) \in (0,\mu]$ and $\widetilde{\eta} = \widetilde{\eta}(n,m,q,\mu) \in (0,1)$ to be later determined.  By \eqref{introduction hyp2} and~\cite[Lemma~3.7]{KrumWica} 
\begin{equation}\label{uniform branch set eqn1}
	\widehat{\mathcal E}(T,P_0,\mathbf{C}_{1/4}(Z,P_0)) \leq C \widehat{E}(T,P_0,\mathbf{C}_1(0)) \leq C \varepsilon_0 , 
\end{equation}
where $C = C(n,m) \in (0,\infty)$ is a constant.  Thus noting that $\pi_{\#} T = q \llbracket B_1(0) \rrbracket$ and using the monotonicity formula
\begin{equation}\label{uniform branch set eqn2}
	\frac{\|T\|(\mathbf{B}_{\rho}(Z))}{\omega_n \rho^n} \leq \frac{\|T\|(\mathbf{C}_{1/4}(Z,P_0))}{\omega_n (1/4)^n} < q + C\varepsilon_0^2 
\end{equation}
for all $\rho \in (0,1/4]$, where $C = C(n,m) \in (0,\infty)$ is a constant.  For each integer $i \geq 0$, let $\rho_i = 2^{-i-5}$ and let $P_i$ be an $n$-dimensional linear plane such that 
\begin{equation*}
	\widehat{E}(T,P_i,\mathbf{B}_{\rho_i}(Z)) = \inf_{P' \in \mathcal{P}} \widehat{E}(T,P',\mathbf{B}_{\rho_i}(Z)) ,
\end{equation*}
where $\mathcal{P}$ denotes the set of all $n$-dimensional linear planes in $\mathbb{R}^{n+m}$.  Observe that for each integer $i \geq 0$, one of the following three possibilities holds true: 
\begin{enumerate}[itemsep=3mm,topsep=0mm]
	\item[(a)]  $\widehat{E}(T,P_i,\mathbf{B}_{4\rho_i}(Z)) \geq \widetilde{\eta}$; 
	\item[(b)]  $\widehat{E}(T,P_i,\mathbf{B}_{4\rho_i}(Z)) < \widetilde{\eta}$ and $\widehat{E}(T,P_i,\mathbf{B}_{2\rho_i}(Z)) \geq 2^{\widetilde{\mu}} \widehat{E}(T,P_{i-1},\mathbf{B}_{2\rho_{i-1}}(Z))$; 
	\item[(c)]  $\widehat{E}(T,P_i,\mathbf{B}_{4\rho_i}(Z)) < \widetilde{\eta}$ and $\widehat{E}(T,P_i,\mathbf{B}_{2\rho_i}(Z)) < 2^{\widetilde{\mu}} \widehat{E}(T,P_{i-1},\mathbf{B}_{2\rho_{i-1}}(Z))$.
\end{enumerate}

We want to show that (c) holds true for all $i \geq 0$.  Suppose to the contrary that either (a) or (b) holds true for some integer $i \geq 0$.  We claim that: \textit{there exist an integer $i_* \geq 0$ such that either (a) or (b) holds true when $i = i_*$ and}
\begin{equation}\label{uniform branch set eqn3}
	\widehat{E}(T,P_Z,\mathbf{B}_{2\rho}(Z)) \leq C \widetilde{\eta} \Big(\frac{\rho}{\rho_{i_*}}\Big)^{\widetilde{\mu}} 
\end{equation}
\textit{for all $\rho \in (0,\rho_{i_*}]$, where $C = C(n,m,q,\widetilde{\mu}) \in (0,\infty)$ is a constant.}  If either (a) or (b) holds true for infinitely many integers $i \geq 0$, then by \eqref{introduction hyp4} it suffices to choose $i_*$ so that either (a) or (b) holds true when $i = i_*$ and $i_*$ is large enough that 
\begin{equation*}
	\sup_{\rho \in (0,1/32]} \rho^{-\mu} \widehat{E}(T,P_Z,\mathbf{B}_{2\rho}(Z)) < \widetilde{\eta} \rho_{i_*}^{-\widetilde{\mu}} .
\end{equation*}
If instead (a) or (b) holds true for only finitely many integers $i \geq 0$, choose $i_*$ to be the largest such integer.  Then (c) holds true for every integer $i > i_*$.  Arguing as we did to obtain (8.19) and (8.22) of~\cite{KrumWica} 
\begin{equation}\label{uniform branch set eqn4}
	\|\pi_{P_i} - \pi_{P_j}\| \leq C \widetilde{\eta} \bigg(\frac{\rho_i}{\rho_{i_*}}\bigg)^{\widetilde{\mu}} 
\end{equation}
for all integers $i_* \leq i \leq j$ and 
\begin{equation}\label{uniform branch set eqn5}
	\widehat{E}(T,P_i,\mathbf{B}_{2\rho}(Z)) \leq C \widetilde{\eta} \bigg(\frac{\rho}{\rho_{i_*}}\bigg)^{\widetilde{\mu}} 
\end{equation}
for all integers $i \geq i_*$ and all $\rho \in [\rho_{i+1},2\rho_{i_*}]$, where $C = C(n,q,\widetilde{\mu}) \in (0,\infty)$ is a constant.  By \eqref{uniform branch set eqn4}, $(P_i)$ is a Cauchy sequence of $n$-dimensional planes in $\mathbb{R}^{n+m}$ and thus $(P_i)$ must converge to some $n$-dimensional plane.  In fact, by \eqref{introduction hyp4} and \eqref{uniform branch set eqn5}, $(P_i)$ converges to the unique tangent plane $P_Z$ of $T$ at $Z$.  It follows from \eqref{uniform branch set eqn5} that \eqref{uniform branch set eqn3} holds true.

Provided $\widetilde{\eta}$ is sufficiently small, by~\cite[Lemma~3.5]{KrumWica} and \eqref{uniform branch set eqn3} 
\begin{equation}\label{uniform branch set eqn6}
	\sup_{X \in \op{spt} T \cap \mathbf{B}_{15\rho/8}(Z)} \op{dist}(X,Z+P_Z) \leq C \widetilde{\eta} \rho < \rho/8 
\end{equation}
for all $\rho \in (0,\rho_{i_*}]$, where $C = C(n,m,q,\widetilde{\mu}) \in (0,\infty)$ is a constant.  Hence $\widetilde{T} = T \llcorner \mathbf{B}_{15\rho_{i_*}/8}(Z)$ is a locally area-minimizing rectifiable current of $\mathbf{C}_{7\rho_{i_*}/4}(Z,P_Z)$ with $(\partial\widetilde{T}) \llcorner \mathbf{C}_{7\rho_{i_*}/4}(Z,P_Z) = 0$.  Moreover, for each $\rho \in (0,\rho_{i_*}]$ and each $X \in \op{spt} \widetilde{T} \cap \mathbf{C}_{7\rho/4}(Z,P_Z)$ we have that $\op{dist}(X,Z+P_Z) \leq \frac{1}{15} |X-Z|$, from which it follows that $\op{dist}(X,Z+P_Z) \leq \frac{1}{14} |\pi_{P_Z}(X-Z)| \leq \rho/8$ so that $|X-Z| \leq 15\rho/8$.  That is, 
\begin{equation}\label{uniform branch set eqn7}
	\op{spt} \widetilde{T} \cap \mathbf{C}_{7\rho/4}(Z,P_Z) \subset \mathbf{B}_{15\rho/8}(Z)
\end{equation}
for all $\rho \in (0,\rho_{i_*}]$.  By \eqref{uniform branch set eqn2} and \eqref{uniform branch set eqn6} (with $\rho = \rho_{I_*}$)
\begin{equation}\label{uniform branch set eqn8}
	\|\widetilde{T}\|(\mathbf{C}_{7\rho_{i_*}/4}(Z,P_Z)) \leq \|T\|(\mathbf{B}_{7\rho_{i_*}/4 + C\widetilde{\eta}\rho_{i_*}}(Z)) 
		\leq (q+C\varepsilon_0^2) \,\omega_n (7\rho_{i_*}/4 + C\widetilde{\eta}\rho_{i_*})^n 
\end{equation}
where $C = C(n,m,q,\widetilde{\mu}) \in (0,\infty)$ are constants.  Provided $\widetilde{\eta}$ and $\widetilde{\mu}$ are sufficiently small, since (a) or (b) holds true for $i = i_*$ and we have \eqref{uniform branch set eqn2} and \eqref{uniform branch set eqn3} (both with $\rho = \rho_{i_*}$), we can apply~\cite[Lemma~7.3 and Lemma~7.7]{KrumWica} (with $\mu/2$ in place of $\nu$) to obtain 
\begin{equation*}
	N_{\widetilde{T},P_Z,Z}(\rho) < 1+\mu/4
\end{equation*}
for all $\rho \in [\rho_{i_*+1},\rho_{i_*}]$.  Hence provided $\widetilde{\eta}$ is sufficiently small, by \eqref{uniform branch set eqn6} (with $\rho = \rho_{i_*}$), \eqref{uniform branch set eqn8}, \eqref{uniform branch set eqn7}, and \eqref{uniform branch set eqn3} we can apply~\cite[Theorem~4.4]{KrumWica} to obtain 
\begin{equation*}
	N_{\widetilde{T},P_Z,Z}(\rho) \leq e^{C\widetilde{\eta}\varsigma} N_{T,P_Z,Z}(\rho_{i_*}) < 1+\mu/2
\end{equation*}
for all $\rho \in (0,\rho_{i_*}]$, where $C = C(n,m,q,\widetilde{\mu}) \in (0,\infty)$ is a constant.  Thus by \eqref{uniform branch set eqn7} and~\cite[Corollary~5.3(ii)]{KrumWica} 
\begin{equation*}
	\widehat{E}(T,P_Z,\mathbf{B}_{2\rho}(Z)) 
	\geq (7/8)^{n/2+1} \widehat{E}(\widetilde{T},P_Z,\mathbf{C}_{7\rho/4}(Z,P_Z))
	\geq c \bigg(\frac{\rho}{\rho_{i_*}}\bigg)^{3\mu/4} \widehat{E}(T,P_Z,\mathbf{B}_{\rho_{i_*}}(Z))
\end{equation*}
for all $\rho \in (0,\rho_{i_*}]$ and for some constant $c = c(n,m,q,\widetilde{\mu},\mu) > 0$, contradicting \eqref{introduction hyp4}.

Now fix a suitably small constant $\widetilde{\mu} = \widetilde{\mu}(n,m,q,\mu) \in (0,\mu]$.  We have shown that (c) holds true for all $i \geq 0$.  Arguing as we did to obtain (8.19) and (8.22) of~\cite{KrumWica}
\begin{equation*}
	\|\pi_{P_i} - \pi_{P_j}\| \leq C \widetilde{\eta} \rho_i^{\widetilde{\mu}} 
\end{equation*}
for all integers $0 \leq i \leq j$ and 
\begin{equation*}
	\widehat{E}(T,P_i,\mathbf{B}_{2\rho}(Z)) \leq C \widetilde{\eta} \rho^{\widetilde{\mu}} 
\end{equation*}
for all integers $i \geq 0$ and all $\rho \in [\rho_{i+1},1/16]$, where $C = C(n,q,\widetilde{\mu}) \in (0,\infty)$ is a constant.  It follows that $(P_i)$ converges to the unique tangent plane $P_Z$ of $T$ at $Z$ and 
\begin{equation}\label{uniform branch set eqn9}
	\widehat{E}(T,P_Z,\mathbf{B}_{2\rho}(Z)) \leq C \widetilde{\eta} \rho^{\widetilde{\mu}} 
\end{equation}
for all $\rho \in (0,1/16]$, where $C = C(n,m,q,\mu) \in (0,\infty)$ is a constant.  By the triangle inequality,~\cite[Lemma~3.4]{KrumWica}, \eqref{uniform branch set eqn9}, and \eqref{introduction hyp2}
\begin{align}\label{uniform branch set eqn10}
	\|\pi_{P_Z} - \pi\|^2 
	\leq\,& \frac{2}{\omega_n (1/16)^n} \int_{G_n(\mathbf{B}_{1/16}(Z))} \|\pi_S - \pi_{P_Z}\|^2 \,d|T|(X,S)  
		\\&+ \frac{2}{\omega_n (1/16)^n} \int_{G_n(\mathbf{B}_{1/16}(Z))} \|\pi_S - \pi\|^2 \,d|T|(X,S) \nonumber 
	\\ \leq\,& 2C \widehat{E}(T,P_Z,\mathbf{B}_{1/8}(Z)) + 2C \widehat{E}(T,P_0,\mathbf{B}_{1/8}(Z)) 
	\leq C \widetilde{\eta} + C \varepsilon_0 , \nonumber
\end{align}
where $C = C(n,m,q,\mu) \in (0,\infty)$ are constants.  By \eqref{introduction hyp2} and \eqref{uniform branch set eqn10}  
\begin{equation}\label{uniform branch set eqn11}
	\op{dist}(X,Z+P_Z) \leq \op{dist}(X,P_0) + \op{dist}(Z,P_0) + \|\pi_{P_Z} - \pi\| \,|X-Z| \leq C \widetilde{\eta} + C \epsilon_0 
\end{equation}
for all $X \in \op{spt} T \cap \mathbf{C}_{1/4}(Z,P_Z)$, where $C = C(n,m,q,\mu) \in (0,\infty)$ is a constant.  \eqref{uniform branch set eqn11} together with $(\partial T) \llcorner \mathbf{C}_1(0) = 0$ gives us \eqref{uniform branch set concl1}.  By \eqref{uniform branch set eqn11} and \eqref{uniform branch set eqn2} 
\begin{equation}\label{uniform branch set eqn12}
	\|T\|(\mathbf{C}_{7/32}(Z,P_Z)) \leq \|T\|(\mathbf{B}_{7/32 + C\widetilde{\eta} + C\varepsilon_0}(Z)) 
		\leq (q+C\varepsilon_0^2) \,\omega_n (7/32 + C\widetilde{\eta} + C\varepsilon_0)^n 
\end{equation}
where $C = C(n,m,q,\mu) \in (0,\infty)$ are constants.  By \eqref{uniform branch set eqn9} and~\cite[Lemma~3.5]{KrumWica}
\begin{equation*}
	\sup_{X \in \op{spt} T \cap \mathbf{B}_{15\rho/8}(Z)} \op{dist}(X,Z+P_Z) \leq C \widetilde{\eta} \rho^{1+\widetilde{\mu}} < \rho/8 
\end{equation*}
for all $\rho \in (0,1/16]$, where $C = C(n,m,q,\mu) \in (0,\infty)$ is a constant.  It follows arguing as we did to obtain \eqref{uniform branch set eqn7} above that $\op{spt} T \cap \mathbf{C}_{7\rho/4}(Z,P_Z) \subset \mathbf{B}_{15\rho/8}(Z)$ for all $\rho \in (0,1/16]$.  Thus by \eqref{uniform branch set eqn9} 
\begin{equation}\label{uniform branch set eqn13}
	\widehat{E}(T,P_Z,\mathbf{C}_{7\rho/4}(Z,P_Z)) \leq C \widetilde{\eta} \rho^{\widetilde{\mu}} 
\end{equation}
for all $\rho \in (0,1/16]$, where $C = C(n,m,q,\mu) \in (0,\infty)$ is a constant.  By \eqref{uniform branch set concl1} and \eqref{uniform branch set eqn13}, we can apply Corollary~\ref{planar frequency cor} to deduce that $\mathcal{N}_{T,{\rm Pl}}(Z)$ exists.  Moreover, by \eqref{introduction hyp4} and Corollary~\ref{planar frequency cor}, $\mathcal{N}_{T,{\rm Pl}}(Z) \geq 1+\mu$.  Equip $P_Z$ with the orientation induced by $T$.  By \eqref{uniform branch set concl1}, \eqref{uniform branch set eqn12}, and \eqref{uniform branch set eqn13}, we can apply~\cite[Lemma~5.1]{KrumWica} to deduce that $q \llbracket P_Z \rrbracket$ is the unique tangent cone to $T$ at $Z$ and \eqref{uniform branch set concl2} holds true.

We claim that 
\begin{equation}\label{uniform branch set eqn14}
	\widehat{E}(T,P_Z,\mathbf{C}_{7\rho/4}(Z,P_Z)) \leq C\rho^{\mu} \widehat{E}(T,P_Z,\mathbf{C}_{7/32}(Z,P_Z)) 
\end{equation}
for all $\rho \in (0,1/8]$. where $C = C(n,m,q,\mu) \in (0,\infty)$ is a constant.  By \eqref{uniform branch set concl1}, $\Theta(T,Z) = q$, \eqref{uniform branch set eqn12}, and \eqref{uniform branch set eqn9}, assuming that $\widetilde{\eta}$ is sufficiently small, we can apply Theorem~\ref{frequency monotonicity thm}(e) together with $\mathcal{N}_{T,{\rm Pl}}(Z) \geq 1+\mu$ to obtain  
\begin{align*}
	\widehat{E}(T,P_Z,\mathbf{C}_{7\rho/4}(Z,P_Z)) 
	\leq\,& e^{C(\mathcal{N}_{T,{\rm Pl}}(Z)+2) \,\widetilde{\eta}^{2\varsigma}} (\rho/8)^{\mathcal{N}_{T,{\rm Pl}}(Z)-1} 
		\widehat{E}(T,P_Z,\mathbf{C}_{7/32}(Z,P_Z)) 
	\\ \leq\,& e^{4C\widetilde{\eta}^{4\varsigma}} (e^{2C\widetilde{\eta}^{2\varsigma}} \theta/8)^{\mathcal{N}_{T,{\rm Pl}}(Z)-1-\mu} (\rho/8)^{\mu} 
		\widehat{E}(T,P_Z,\mathbf{C}_{7/32}(Z,P_Z))  
	\\ \leq\,& C\rho^{\mu} \widehat{E}(T,P_Z,\mathbf{C}_{7/32}(Z,P_Z))  
\end{align*}
for all $\rho \in (0,\theta]$, where for $\theta = \theta(n,m,q,\mu) \in (0,1/8]$ is a suitably small constant and $C = C(n,m,q,\mu) \in (0,\infty)$ are constants.  If instead $\rho \in (\theta,1/8]$, then \eqref{uniform branch set eqn14} obviously holds true with $C = 8^{-n/2-1} \theta^{-n/2-1-\mu}$.  

To show that \eqref{uniform branch set concl3}--\eqref{uniform branch set concl5} hold true, notice that by the triangle inequality,~\cite[Lemma~3.7]{KrumWica}, and \eqref{uniform branch set eqn14} 
\begin{align}\label{uniform branch set eqn15}
	|\vec P_Z - \vec P_0|^2 
	\leq\,& 2 \widehat{\mathcal E}(T,P_Z,\mathbf{C}_{7\rho/8}(Z,P_Z))^2 + 2 \rho^{-n} \widehat{\mathcal E}(T,P_0,\mathbf{C}_{7/8}(0))^2 %\nonumber 
	\\ \leq\,& C \widehat{E}(T,P_Z,\mathbf{C}_{7\rho/4}(Z,P_Z))^2 + C \rho^{-n} \widehat{E}(T,P_0,\mathbf{C}_1(0))^2 \nonumber 
	\\ \leq\,& C\rho^{2\mu} \widehat{E}(T,P_Z,\mathbf{C}_{7/32}(Z,P_Z))^2 + C \rho^{-n} \widehat{E}(T,P_0,\mathbf{C}_1(0))^2 \nonumber 
\end{align}
for all $\rho \in (0,1/8]$, where $C = C(n,m,q,\mu) \in (0,\infty)$ are constants.  By the triangle inequality, \eqref{uniform branch set eqn12}, and \eqref{uniform branch set eqn15} 
\begin{align*}
	\widehat{E}(T,P_Z,\mathbf{C}_{7/32}(Z,P_Z)) 
	\leq\,& \widehat{E}(T,P_0,\mathbf{C}_1(0)) + C |\vec P_Z - \vec P_0| 
	\\ \leq\,& C\rho^{\mu} \widehat{E}(T,P_Z,\mathbf{C}_{7/32}(Z,P_Z)) + C \rho^{-n/2} \widehat{E}(T,P_0,\mathbf{C}_1(0)) 
\end{align*}
for all $\rho \in (0,1/8]$, where $C = C(n,m,q,\mu) \in (0,\infty)$ are constants.  Choosing $\rho = \rho(n,m,q,\mu) \in (0,1/8]$ small enough that $C\rho^{\mu} < 1/2$ 
\begin{equation}\label{uniform branch set eqn16}
	\widehat{E}(T,P_Z,\mathbf{C}_{7/32}(Z,P_Z)) \leq C \widehat{E}(T,P_0,\mathbf{C}_1(0)) , 
\end{equation}
where $C = C(n,m,q,\mu) \in (0,\infty)$ is a constant.  By combining \eqref{uniform branch set eqn15} (with $\rho = 1/8$) and \eqref{uniform branch set eqn16}, we obtain \eqref{uniform branch set concl3}.  By \eqref{uniform branch set eqn11}, $\op{spt} T \cap \mathbf{C}_{7/32}(Z,P_Z) \subset \mathbf{C}_{1/4}(Z,P_0)$.  Thus by the triangle inequality, \eqref{uniform branch set eqn12}, \eqref{uniform branch set eqn1}, and \eqref{uniform branch set concl3}
\begin{equation*}
	\widehat{\mathcal E}(T,P_Z,\mathbf{C}_{7/32}(Z,P_Z)) 
	\leq \widehat{\mathcal E}(T,P_0,\mathbf{C}_{1/4}(Z,P_0)) + C |\vec P_Z - \vec P_0| 
	\leq C \widehat{E}(T,P_0,\mathbf{C}_1(0)) , 
\end{equation*}
where $C = C(n,m,q,\mu) \in (0,\infty)$ are constants.  Thus by \eqref{uniform branch set concl2} and \eqref{oriented tilt and mass}, \eqref{uniform branch set concl4} holds true.  By combining \eqref{uniform branch set eqn14} and \eqref{uniform branch set eqn16}, we obtain \eqref{uniform branch set concl5}.  
\end{proof}

Now let $q \geq 2$ be an integer and let $\mu \in (0,1)$.  For $\varepsilon_0 = \varepsilon_0(n,m,q,\mu) \in (0,1)$ suitably small, let $T$ be as in Hypothesis~\ref{introduction hyp} and $\mathcal{B}_q$ is as in Definition~\ref{introduction B defn}.  By Lemma~\ref{uniform branch set lemma} and Corollary~\ref{planar frequency cor}, for each $Z \in \mathcal{B}_q$ we have that $\mathcal{N}_{T,{\rm Pl}}(Z)$ exists and $\mathcal{N}_{T,{\rm Pl}}(Z) \geq 1+\mu$.  By Theorem~\ref{B tanfun thm}, for each $Z \in \mathcal{B}_q$ there exists at least one tangent function $\varphi$ of $T$ at $Z$.  Recall from the Introduction that $\mathcal{B}^{(\neq 2)}_q = \mathcal{B}^{(\neq 2)}_q(T) = \{ X \in \mathcal{B}_q : \mathcal{N}_{T,{\rm Pl}}(Z) \neq 2 \}$.  Note that by Lemma~\ref{uniform branch set lemma}, $\mathcal{B}_q$ is precisely the set of all branch points $Z$ of $T$ for which \eqref{uniform branch set concl5} holds true for some $n$-dimensional plane $P_Z$, and consequently $\mathcal{B}_q$ is a relatively closed subset of $\mathbf{C}_{1/2}(0)$.  For each $j \in \{0,1,2,\ldots,n-2\}$, we define $S^{(j)}_q = S^{(j)}_q(T)$ to be the set of all points $Z \in \mathcal{B}^{(\neq 2)}_q$ such that $\dim\op{spine}\varphi \leq j$ for every tangent function $\varphi$ of $T$ at $Z$, where $\op{spine}\varphi = \{ X \in P_Z : \mathcal{N}_{\varphi}(X) = \mathcal{N}_{\varphi}(0) \}$ (as in~\cite[Subsection~3.6]{KrumWica}).  Observe that 
\begin{equation*}
	\mathcal{B}^{(\neq 2)}_q = S^{(n-2)}_q \supseteq S^{(n-3)}_q \supseteq\cdots\supseteq S^{(1)}_q \supseteq S^{(0)}_q . 
\end{equation*}

\begin{theorem}\label{B prime stratification thm}
Let $q \geq 2$ be an integer and $\mu \in (0,1)$.  There exists $\varepsilon_0 = \varepsilon_0(n,m,q,\mu) \in (0,1)$ such that if $T$ is as in Hypothesis~\ref{introduction hyp} and $\mathcal{B}_q$ is as in Definition~\ref{introduction B defn}, then: 
\begin{enumerate}[itemsep=3mm,topsep=0mm]
	\item[{\rm (i)}]  for each  $j \in \{0,1,2,\ldots,n-2\}$, $S^{(j)}_{q}$ has Hausdorff dimension at most $j$;
	\item[{\rm (ii)}]  for each $\alpha > 0$, $S^{(0)}_{q} \cap \{ Z : \mathcal{N}_{T,{\rm Pl}}(Z) = \alpha \}$ is discrete;
	\item[{\rm (iii)}]  $\mathcal{B}^{(\neq 2)}_q = \{ Z \in \mathcal{B}_q : \mathcal{N}_{T,{\rm Pl}}(Z) \neq 2 \}$ has Hausdorff dimension at most $n-2$.
\end{enumerate}
In particular, Theorem~\ref{B prime dimension thm} from the Introduction holds true.
\end{theorem}

\begin{proof}
Let $Z \in \mathcal{B}^{(\neq 2)}_q$.  Set $\alpha = \mathcal{N}_{T,{\rm Pl}}(Z)$ and let $\varphi$ be any tangent function of $T$ at $Z$.  For each $\xi \in P_Z$ and $X \in P_Z$, let $\ell_{\xi}(X) = \varphi_a(\xi) + D\varphi_a(\xi) \cdot (X-\xi)$ and $(\varphi - \ell_{\xi})(X) = \sum_{l=1}^q \llbracket \varphi_l(X) - \ell_{\xi}(X) \rrbracket$, where $\varphi(X) = \sum_{l=1}^q \llbracket \varphi_l(X) \rrbracket$ and $\varphi_a(X) = \frac{1}{q} \sum_{l=1}^q \varphi_l(X)$.  Let 
\begin{equation*}
	L_{\varphi} = \{ \xi \in P_Z : \mathcal{N}_{\varphi-\ell_{\xi}}(\xi) = \alpha \} ,
\end{equation*}
where $\mathcal{N}_{\varphi-\ell_{\xi}}(\xi)$ is the frequency of $\varphi-\ell_{\xi}$ at $\xi$ (as in~\cite[Subsection~3.6]{KrumWica}).  We want to show that $\varphi_a(\xi) = 0$ and $D\varphi_a(\xi) = 0$ so that $L_{\varphi} = \op{spine}\varphi$.  Recall that $\varphi$ is homogeneous degree $\alpha$.  If $\alpha$ is not an integer, then $\varphi_a$ is a homogeneous degree $\alpha$ single-valued harmonic function and thus $\varphi_a \equiv 0$ on $P_Z$.  If instead $\alpha$ is an integer $\geq 3$, then $\varphi_a(\xi+X)  - \ell_{\xi}(\xi+X)$ is a homogeneous degree $\alpha$ harmonic polynomial and thus 
\begin{equation}\label{B prime stratification eqn1}
	\varphi_a(\xi+X)  - \ell_{\xi}(\xi+X) = \varphi_a(X)
\end{equation}
for all $X \in P_Z$.  Notice that if $\psi : P_Z \rightarrow P_Z^{\perp}$ is a homogeneous polynomial of degree $\alpha \geq 1$ and $D\psi(\xi) = 0$, then 
\begin{equation}\label{B prime stratification eqn2}
	\psi(\xi) = \psi(\xi) - \psi(0) = \int_0^1 D\psi(t\xi) \cdot \xi \,dt = \int_0^1 t^{\alpha-1} D\psi(\xi) \cdot \xi \,dt = 0 .
\end{equation}
Differentiating both sides of \eqref{B prime stratification eqn1} twice and using $\alpha \geq 3$, we obtain $D^2 \varphi_a(\xi) = D^2 \varphi_a(0) = 0$.  Thus by applying \eqref{B prime stratification eqn2} with $\psi = D\varphi_a$ and $\alpha-1$ in place of $\alpha$, we deduce that $D\varphi_a(\xi) = 0$.  By applying \eqref{B prime stratification eqn2} again with $\psi = \varphi_a$, we deduce that $\varphi_a(\xi) = 0$.

To prove the conclusions of the theorem, we adapt the well-known argument~\cite[Theorem~2.26 and Corollary~2.27]{Almgren}.  We will focus on the proof of conclusion~(i), as the proof of conclusion~(ii) involves adapting the usual argument in a similar manner and conclusion~(iii) follows directly from conclusion~(i).  To see conclusion~(i), we need to show that for each $\delta > 0$ and $Z \in S^{(j)}_{q}$ there exists $\varepsilon > 0$ (depending on $T,Z$) such that for every $\rho \in (0,\varepsilon]$ there exists a $j$-dimensional linear subspace $L \subset P_Z$ (depending on $T,Z,\rho$) such that 
\begin{equation}\label{B prime stratification eqn3}
	\{ X \in \mathcal{B}^{(\neq 2)}_q \cap \mathbf{B}_{\rho}(Z) : \mathcal{N}_{T,{\rm Pl}}(X) \geq \mathcal{N}_{T,{\rm Pl}}(Z) - \varepsilon \} 
	\subset \{ X \in \mathbf{B}_{\rho}(Z) : \op{dist}(X,Z+L) < \delta\rho \} .  
\end{equation}
To see this, fix $\delta > 0$ and $Z \in S^{(j)}_q$ and set $\alpha = \mathcal{N}_{T,{\rm Pl}}(Z)$.  Suppose to the contrary that for $\nu = 1,2,3,\ldots$ there exists numbers $\varepsilon_{\nu}$ and $\rho_{\nu}$ such that $0 < \rho_{\nu} \leq \varepsilon_{\nu}$, $\varepsilon_{\nu} \rightarrow 0$, and \eqref{B prime stratification eqn3} does not hold true for any $j$-dimensional linear subspace $L \subset P_Z$.  After passing to a subsequence, let $\varphi \in W^{1,2}_{\rm loc}(P_Z,\mathcal{A}_q(P_Z^{\perp}))$ be the blow-up of $\eta_{Z,\rho_{\nu}\#} T$ relative to $P_Z$ and excess $\widehat{E}(T,P_Z,\mathbf{C}_{\rho_{\nu}}(Z,P_Z))$ in $\mathbf{C}_{\sigma}(0,P_Z)$ for every $\sigma \in [1,\infty)$.  Let $L \subset P_Z$ be a $j$-dimensional linear subspace with $\op{spine}\varphi \subseteq L$.  Since \eqref{B prime stratification eqn3} does not hold true, for each positive integer $\nu$ there exists $X_{\nu} \in \mathcal{B}^{(\neq 2)}_q \cap \mathbf{B}_{\rho_{\nu}}(Z)$ such that $\mathcal{N}_{T,{\rm Pl}}(X_{\nu}) \geq \alpha - \varepsilon_{\nu}$ but $\op{dist}(X_{\nu},Z+L) \geq \delta\rho_{\nu}$.  After passing to a subsequence, let $\pi_{P_Z}(X_{\nu} - Z)/\rho_{\nu} \rightarrow \xi$ in $P_Z$.  By Lemma~\ref{uniform branch set lemma} and~\cite[Corollary~6.6]{KrumWica}, $\mathcal{N}_{\varphi-\ell_{\xi}}(\xi) \geq \alpha$ so that $\xi \in L_{\varphi}$.  But by the discussion above, $L_{\varphi} = \op{spine} \varphi$.  Hence $\xi \in \op{spine}\varphi \subseteq L$ and $\op{dist}(\xi,L) \geq \delta$, giving us a contradiction.  To complete the proof of conclusion~(i), we apply a well-known covering argument.
\end{proof}

\begin{remark}\label{B prime stratification rmk}{\rm 
Let $\mathcal{B}^{(=2)}_q = \{ Z \in \mathcal{B}_q : \mathcal{N}_{T,{\rm Pl}}(Z) = 2 \}$.  Suppose that $Z \in \mathcal{B}^{(=2)}_q$.  Then \eqref{B prime stratification eqn1} holds true for all $\xi \in P_Z$.  It follows that if $\varphi(X) = q \llbracket \varphi_a(X) \rrbracket$ for all $X \in P_Z$, then $L_{\varphi} = P_Z$ is an $n$-dimensional linear subspace.  Hence arguing as in the proof of Theorem~\ref{B prime stratification thm} would only give us that $\mathcal{B}^{(=2)}_q$ has Hausdorff dimension at most $n$.  For this reason, it is necessary to consider the center manifold of $T$ in order to prove the optimal Hausdorff dimension and rectifiability of $\mathcal{B}^{(=2)}_q$.
}\end{remark}

\begin{remark}\label{dist tangent plane rmk}{\rm 
When applying the semi-continuity of frequency~\cite[Corollary~6.6]{KrumWica}, we had to account for how the tangent plane $P_X$ to $T$ at a branch point $X$ differs from the tangent plane $P_Z$ to $T$ at $Z$.  We can readily show that for each $\rho \in (0,(1-|Z|)/4]$ and each $X \in \mathcal{B}_q \cap \mathbf{B}_{\rho}(Z)$, the distance between $P_X$ and $P_Z$ is bounded by the coarse excess as 
$$\op{dist}_{\mathcal H}(P_X \cap \mathbf{B}_1(0),P_Z \cap \mathbf{B}_1(0)) \leq C \widehat{E}(T,P_Z,\mathbf{C}_{4\rho}(Z,P_Z)) ,$$ 
where $C = C(n,m,q,\alpha) \in (0,\infty)$ is a constant. We in fact show that under certain assumptions stronger estimates hold true.  Let $\alpha > 1$ with $\alpha \neq 2$.  The proof of Theorem~\ref{B prime stratification thm} shows that for $\rho \in (0,(1-|Z|)/4]$ and $X \in \mathcal{B}_q \cap \mathbf{B}_{\rho}(Z)$ with $\mathcal{N}_{T,{\rm Pl}}(X) \geq \mathcal{N}_{T,{\rm Pl}}(Z) - \varepsilon(\rho)$ where $\varepsilon(\rho) \rightarrow 0^+$ as $\rho \rightarrow 0^+$
$$\op{dist}_{\mathcal H}(P_X \cap \mathbf{B}_1(0),P_Z \cap \mathbf{B}_1(0)) = o(\widehat{E}(T,P_Z,\mathbf{C}_{4\rho}(Z,P_Z))) \text{ as } \rho\rightarrow 0^+ $$
so that the tangent planes $P_X$ converge to $P_Z$ under coarse blow-up limits.  Moreover, it follows from Theorem~\ref{nonconcentration thm} below that if $\rho_{\nu} \rightarrow 0^+$ such that the blow-up $\varphi$ of $(\eta_{Z,\rho_{\nu}\#} T)$ relative to $P_Z$ has spine dimension $n-2$ and if $X_{\nu} \in \mathcal{B}_q \cap \mathbf{B}_{\rho_{\nu}}(Z)$ with $\mathcal{N}_{T,{\rm Pl}}(X_{\nu}) \geq \alpha$, then 
$$\op{dist}_{\mathcal H}(P_{X_{\nu}} \cap \mathbf{B}_1(0),P_Z \cap \mathbf{B}_1(0)) \leq C E_{\varsigma/2}(T,\varphi,\mathbf{C}_{4\rho_{\nu}}(Z,P_Z))$$ 
for all sufficiently large $\nu$, where $E_{\varsigma/2}(T,P_Z,\mathbf{C}_{4\rho_{\nu}}(Z,P_Z))$ is the fine excess defined roughly as the $L^2$-distance of $T$ to the graph of $\varphi$ in $\mathbf{C}_{4\rho_{\nu}}(Z,P_Z)$ (see \eqref{graphical E defn} below) and $\varsigma = \varsigma(n,m,q) \in (0,1)$ and $C = C(n,m,q,\alpha) \in (0,\infty)$ are constants.
}\end{remark}

An immediate consequence of Theorem~\ref{B prime stratification thm} is that for $\mathcal{H}^{n-2}$-a.e.~$Z \in \mathcal{B}^{(\neq 2)}_q$, at least one tangent function $\varphi$ of $T$ at $Z$ has $\dim \op{spine}\varphi = n-2$.  In this case, $\varphi$ must be homogeneous degree $\alpha$, where $\alpha$ is a rational number of the form $\alpha = \ell_0/q_0$ for some relatively prime positive integers $\ell_0,q_0$ with $q_0 \leq q$.  In particular, $\alpha \geq 1+1/q$.  Moreover, after an orthogonal change of coordinates, we may assume that $P_Z = \mathbb{R}^n \times \{0\}$ and $\op{spine}\varphi = \{0\}^2 \times \mathbb{R}^{n-2}$ so that 
\begin{equation*}
	\varphi = \sum_{j=1}^J m_j \varphi_j
\end{equation*}
where $J,m_1,\ldots,m_J$ are positive integers and $\varphi_j$ are distinct multi-valued functions such that either $\varphi_j : \mathbb{R}^n \rightarrow \mathbb{R}^m$ is given by $\varphi_j(x) = 0$ for all $x \in \mathbb{R}^n$ or $\varphi_j : \mathbb{R}^n \rightarrow \mathcal{A}_{q_0}(\mathbb{R}^m)$ is given by 
\begin{equation}\label{homog dirmin eqn1}
	\varphi_j(x_1,x_2,x_3,\ldots,x_n) = \op{Re}(c_j (x_1+ix_2)^{\alpha})
\end{equation}
for some $c_j \in \mathbb{C}^m \setminus \{0\}$.  By~\cite[Theorem~3.2]{MicWhi95}, whenever $\varphi_j$ is given by \eqref{homog dirmin eqn1}, $c_j \cdot c_j = 0$.

\subsection{Locally uniform weak approximation of area-minimizing currents by tangent functions}  The remainder of the paper will be focused on proving Theorem~\ref{unique tangent function thm} and Theorem~\ref{B rectifiable thm} from the Introduction.  We close this section by proving weakly locally uniform approximation of locally area-minimizing rectifiable currents by homogeneous Dirichlet energy minimizing functions at branch points.  This plays an important role in applying the blow-up method of~\cite{Sim93} to prove Theorem~\ref{unique tangent function thm} and Theorem~\ref{B rectifiable thm}, and is analogous to~\cite[Lemma~2.4]{Sim93}.

Let $q \geq 2$ be an integer.  Let $Z \in \mathbb{R}^{n+m}$, $\rho > 0$, and $P$ be an $n$-dimensional plane in $\mathbb{R}^{n+m}$.  Let $T$ be an $n$-dimensional rectifiable current of $\mathbf{C}_{\rho}(Z,P)$.  Let $\alpha > 1$ with $\alpha \neq 2$ and let $\varphi \in W^{1,2}_{\rm loc}(P,\mathcal{A}_q(P^{\perp}))$ be a non-zero, homogeneous degree $\alpha$, locally Dirichlet energy minimizing $q$-valued function.  We define 
\begin{align}\label{Q defn} 
	&Q(T,\varphi,\mathbf{C}_{\rho}(Z,P))^2 
		= \frac{1}{\omega_n \rho^{n+2}} \int_{\mathbf{C}_{\rho}(Z,P)} {\rm dist}^2(X, Z+\op{spt}\op{graph}\varphi) \,d\|T\|(X) 
		\\&\hspace{15mm} + \frac{1}{\omega_n \rho^{n+2}} \int_{B_{\rho/2}(0,P) \cap \{r > \rho/16\}} \sum_{l=1}^q {\rm dist}^2(Z+X+\varphi_l(X), \op{spt} T) 
		\,d\mathcal{H}^n(X) , \nonumber 
\end{align}
where $r(X) = \op{dist}(X,\op{spine}\varphi)$ and $\varphi(X) = \sum_{l=1}^q \llbracket \varphi_l(X) \rrbracket$ for each $X \in P$. 

\begin{remark}\label{T and varphi scaling rmk}{\rm 
Let $Z \in \mathbb{R}^{n+m}$, $\rho > 0$, and $P$ be an $n$-dimensional plane in $\mathbb{R}^{n+m}$.  Let $\lambda > 0$, $T$ be an $n$-dimensional locally area-minimizing rectifiable current of $\mathbf{C}_{\lambda\rho}(Z,P)$, and $\varphi \in W^{1,2}_{\rm loc}(P,\mathcal{A}_q(P^{\perp}))$ be a non-zero, homogeneous degree $\alpha$, locally Dirichlet energy minimizing $q$-valued function.  Recall that $\eta_{Z,\lambda}(X) = (X-Z)/\lambda$ for all $X \in \mathbb{R}^{n+m}$.  Upon rescaling by $\eta_{Z,\lambda}$, we have that $\eta_{Z,\lambda\#} T$ is a locally area-minimizing rectifiable current of $\mathbf{C}_{\rho}(0,P)$.  Applying the same rescaling to the graph of $\varphi$, we have that $\eta_{0,\lambda}(\op{graph}\varphi) = \op{graph}\varphi_{\lambda}$ where $\varphi_{\lambda}(X) = \lambda^{-1} \varphi(\lambda X) = \lambda^{\alpha-1} \varphi(X)$ for all $X \in P$.  Moreover, $Q(T,\varphi,\mathbf{C}_{\lambda\rho}(Z,P)) = Q(\eta_{Z,\lambda\#} T, \varphi_{\lambda}, \mathbf{C}_{\rho}(0,P))$.
}\end{remark}

\begin{lemma}\label{lemma2_4} 
Let $q \geq 2$ be an integer, $\mu \in (0,1)$, and $\kappa \in (1+1/q,\infty)$.  There are functions $\delta : (0,1) \rightarrow (0,1)$ and $R : (0,1) \rightarrow (64,\infty)$ depending only on $n,m,q,\mu,\kappa$ such that the following holds true.  Let $\alpha \in [1+1/q,\kappa]$ with $|\alpha - 2| \geq 1/q$ and $\varepsilon \in (0,1)$.  Let $T$ is an $n$-dimensional locally area-minimizing rectifiable current of $\mathbf{C}_{R(\varepsilon)}(0)$ such that 
\begin{gather}
	\label{lemma2_4 hyp1} (\partial T) \llcorner \mathbf{C}_{R(\varepsilon)}(0) = 0, \quad \sup_{X \in \op{spt} T} \op{dist}(X,P_0) \leq \delta(\varepsilon), \quad 
		\Theta(T,0) = q, \quad \pi_{P_0\#} T = q \llbracket B_{R(\varepsilon)}(0) \rrbracket, \\
	\label{lemma2_4 hyp3} \sup_{\rho \in (0,R(\varepsilon)/2]} \rho^{-\mu} \widehat{E}(T,P_0,\mathbf{B}_{2\rho}(0)) < \infty , \\
	\label{lemma2_4 hyp4} N_{T,P_0,0}(R(\varepsilon)/2) < \alpha + \delta(\varepsilon) , 
\end{gather}
where $P_0 = \mathbb{R}^n \times \{0\}$.  Let $Z \in \mathbf{C}_1(0)$ be a branch point of $T$ at which $\Theta(T,Z) = q$, there is a unique tangent plane $P_Z$ of $T$ at $Z$, and 
\begin{gather}
	\label{lemma2_4 hyp6} \sup_{\rho \in (0,R(\varepsilon)/8]} \rho^{-\mu} \widehat{E}(T,P_Z,\mathbf{B}_{2\rho}(Z)) < \infty , \\
	\label{lemma2_4 hyp7} \mathcal{N}_{T,{\rm Pl}}(Z) \geq \alpha . 
\end{gather}
Then the following holds true: 
\begin{enumerate}[itemsep=3mm,topsep=0mm]
	\item[(i)]  $|N_{T,P_Z,Z}(\rho) - \alpha| < \varepsilon^2$ for all $\rho \in (0,R(\varepsilon)/8]$;
	\item[(ii)]  for every $\rho \in (0,1]$ there exists a non-zero locally Dirichlet energy minimizing $q$-valued function $\varphi \in W^{1,2}_{\rm loc}(P_Z,\mathcal{A}_q(P_Z^{\perp}))$ (depending on $\rho$) such that $\varphi$ is homogeneous degree $\mathcal{N}_{\varphi}(0)$, $|\mathcal{N}_{\varphi}(0) - \alpha| < \varepsilon^2$, and 
	\begin{equation*}%\label{lemma2_4 concl1}
		Q(T,\varphi,\mathbf{C}_{\rho}(Z,P_Z)) < \varepsilon \,\widehat{E}(T,P_Z,\mathbf{C}_{\rho}(Z,P_Z)) ; 
	\end{equation*}
	\item[(iii)]  for every $\rho \in (0,1]$, either there exists a non-zero locally Dirichlet energy minimizing $q$-valued function $\varphi \in W^{1,2}_{\rm loc}(P_Z,\mathcal{A}_q(P_Z^{\perp}))$ (depending on $\rho$) such that $\varphi$ is homogeneous degree $\mathcal{N}_{\varphi}(0)$, $|\mathcal{N}_{\varphi}(0) - \alpha| < \varepsilon^2$, $\dim \op{spine}\varphi = n-2$, and 
	\begin{equation}\label{lemma2_4 concl2}
		Q(T,\varphi,\mathbf{C}_{\rho}(Z,P_Z)) < \varepsilon \,\widehat{E}(T,P_Z,\mathbf{C}_{\rho}(Z,P_Z)) 
	\end{equation}
	or there exists an $(n-3)$-dimensional linear subspace $L$ of $P_Z$ (depending on $\rho$) such that 
	\begin{equation}\label{lemma2_4 concl3}
		\{ X \in \mathcal{B}_q \cap \overline{\mathbf{C}_{\rho}(Z)} : \mathcal{N}_{T,{\rm Pl}}(X) \geq \alpha \} 
			\subset \{ X \in \mathbb{R}^{n+m} : \op{dist}(X,Z+L) < \varepsilon\rho \} , 
	\end{equation}
	where $\mathcal{B}_q$ denotes the set of all branch points $X \in \mathbf{C}_2(0)$ of $T$ at which $\Theta(T,X) = q$, $P_X$ denotes the unique tangent plane to $T$ at $X$, and $\sup_{\rho \in (0,R(\varepsilon)/8]} \rho^{-\mu} \widehat{E}(T,P_X,\mathbf{B}_{2\rho}(X)) < \infty$.
\end{enumerate}
\end{lemma}

\begin{proof}
By \eqref{lemma2_4 hyp1} and~\cite[Lemma~3.7]{KrumWica}, $\widehat{\mathcal E}(T,P_0,\mathbf{C}_{7R(\varepsilon)/8}(0)) \leq C(n,m) \,\delta(\varepsilon)$ and thus by \eqref{oriented tilt and mass} 
\begin{equation}\label{lemma2_4 eqn1}
	\|T\|(\mathbf{C}_{7R(\varepsilon)/8}(0)) \leq (q+C\delta(\varepsilon)^2) \,\omega_n (7R(\varepsilon)/8)^n , 
\end{equation}
where $C = C(n,m) \in (0,\infty)$ is a constant.  Equip $P_Z$ with the orientation induced by $T$.  Assuming $\delta(\varepsilon)$ is sufficiently small, by \eqref{lemma2_4 hyp1}, \eqref{lemma2_4 hyp3}, %\eqref{lemma2_4 hyp5}, 
and \eqref{lemma2_4 hyp6}, we can apply Lemma~\ref{uniform branch set lemma} with $\eta_{0,R(\varepsilon)\#} T$ in place of $T$ to deduce that, assuming $\delta(\varepsilon)$ is sufficiently small, 
\begin{equation}\label{lemma2_4 eqn2} 
	\widehat{E}(T,P_0,\mathbf{C}_{7\rho/4}(0)) \leq C \delta(\varepsilon) \,\bigg(\frac{\rho}{R(\varepsilon)}\bigg)^{\mu} 
		\text{ for all } \rho \in (0,R(\varepsilon)/2] 
\end{equation}
and 
\begin{gather}
	\label{lemma2_4 eqn3} \op{spt} T \cap \overline{\mathbf{C}_{R(\varepsilon)/4}(Z,P_Z)} \subset \mathbf{C}_{R(\varepsilon)}(0), \quad 
		(\partial T) \llcorner \mathbf{C}_{R(\varepsilon)/4}(Z,P_Z) = 0 , \\
	\label{lemma2_4 eqn4} \pi_{P_Z\#} (T \llcorner \mathbf{C}_{R(\varepsilon)/4}(Z,P_Z)) = q \llbracket B_{R(\varepsilon)/4}(Z,P_Z) \rrbracket , \\
	\label{lemma2_4 eqn5} \|T\|(\mathbf{C}_{7R(\varepsilon)/32}(Z,P_Z)) \leq (q+C\delta(\varepsilon)^2) \,\omega_n (7R(\varepsilon)/32)^n , \\
	\label{lemma2_4 eqn6} \widehat{E}(T,P_Z,\mathbf{C}_{7\rho/4}(Z,P_Z)) \leq C \delta(\varepsilon) \,\bigg(\frac{\rho}{R(\varepsilon)}\bigg)^{\mu} 
		\text{ for all } \rho \in (0,R(\varepsilon)/8] , 
\end{gather}
where $C = C(n,m,q,\mu) \in (0,\infty)$ is a constant.

\noindent\textit{Proof of Conclusion~(i).}  By $\Theta(T,Z) = q$, \eqref{lemma2_4 eqn3}, \eqref{lemma2_4 eqn5}, and \eqref{lemma2_4 eqn6}, we can apply Theorem~\ref{frequency monotonicity thm}(b) (with $\rho_0 = R(\varepsilon)/8$) together with \eqref{lemma2_4 hyp7} to obtain 
\begin{equation*}
	\alpha e^{-C\delta(\varepsilon)^{\varsigma}}  \leq e^{-C\delta(\varepsilon)^{\varsigma}} \mathcal{N}_{T,{\rm Pl}}(Z) \leq N_{T,P_Z,Z}(\rho) 
		\leq e^{C\delta(\varepsilon)^{\varsigma}} N_{T,P_Z,Z}(R(\varepsilon)/8)
\end{equation*}
for all $\rho \in (0,R(\varepsilon)/8]$, where $\varsigma = \varsigma(n,m,q) \in (0,1)$ is as in Theorem~\ref{lip approx thm} and $C = C(n,m,q,\mu) \in (0,\infty)$ is a constant.  Thus, assuming $\delta(\varepsilon)$ is suitably small, it suffices to show that $N_{T,P_Z,Z}(R(\varepsilon)/8) < \alpha + \varepsilon^2/4$.

Fix $\varepsilon \in (0,1)$.  Suppose to the contrary that for $\nu = 1,2,3,\ldots$ there exists $\delta_{\nu} \rightarrow 0^+$, $R_{\nu} \rightarrow \infty$, $\alpha_{\nu} \in [1+1/q,\kappa]$ with $|\alpha_{\nu} - 2| \geq 1/q$, an $n$-dimensional locally area-minimizing rectifiable current $T_{\nu}$ of $\mathbf{C}_{R_{\nu}}(0)$, $Z_{\nu} \in \mathbf{C}_1(0)$ with $\Theta(T_{\nu},Z_{\nu}) = q$, and an $n$-dimensional plane $P_{\nu}$ such that \eqref{lemma2_4 hyp1}--\eqref{lemma2_4 hyp7} hold true with $\delta_{\nu},R_{\nu},\alpha_{\nu},T_{\nu},Z_{\nu},P_{\nu}$ in place of $\delta(\varepsilon),R(\varepsilon),\alpha,T,Z,P_Z$ but 
\begin{equation}\label{lemma2_4 eqn7}
	N_{T_{\nu},P_{\nu},Z_{\nu}}(R_{\nu}/8) \geq \alpha_{\nu} + \varepsilon^2/4 . 
\end{equation}
After passing to a subsequence, let $\alpha_{\nu} \rightarrow \alpha$ in $\mathbb{R}$ as $\nu\rightarrow\infty$.  By \eqref{lemma2_4 hyp1} and \eqref{lemma2_4 eqn1}, we can follow the blow-up procedure of~\cite[Subsection~6.1]{KrumWica} to construct a blow-up $w \in W^{1,2}_{\rm loc}(B_{7/8}(0),\mathcal{A}_q(\mathbb{R}^m))$ of $\eta_{0,R_{\nu}\#} T_{\nu}$ with respect to the plane $P_0$ and excess $\widehat{E}_{\nu} = \widehat{E}(T_{\nu},P_0,\mathbf{C}_{7R_{\nu}/8}(0))$.  By \eqref{lemma2_4 hyp4} and Theorem~\ref{frequency monotonicity thm}(e), $\widehat{E}(T_{\nu},P_0,\mathbf{C}_{R_{\nu}/2}(0)) \geq c \widehat{E}_{\nu}$ for some constant $c = c(n,m,q,\mu,K) > 0$, and thus $\|w\|_{L^2(B_{1/2}(0))} \geq c \omega_n^{1/2} (1/2)^{n/2+1}$ so that $w$ is non-zero.  Notice that $Z_{\nu}/R_{\nu} \rightarrow 0$ as $\nu \rightarrow \infty$.  Let $w_a(x) = \frac{1}{q} \sum_{l=1}^q w_l(x)$ for each $x \in B_{7/4}(0)$, where $w(x) = \sum_{l=1}^q \llbracket w_l(x) \rrbracket$.  Since $\Theta(T_{\nu},0) = q$, by~\cite[Lemma~6.4]{KrumWica} we have that $w_a(0) = 0$.  By the proof of~\cite[Corollary~6.6]{KrumWica}, since $P_0$ is the tangent plane to $T_{\nu}$ at $0$ and satisfies \eqref{lemma2_4 eqn2}, we have that $Dw_a(0) = 0$.  Also by the proof of~\cite[Corollary~6.6]{KrumWica}, since $Z_{\nu}/R_{\nu} \rightarrow 0$, $P_{\nu}$ is the tangent plane to $T_{\nu}$ at $Z_{\nu}$, and $P_{\nu}$ satisfies \eqref{lemma2_4 eqn6}, we can express $P_{\nu} = \{ (x,A_{\nu} x) : x \in \mathbb{R}^n \}$ where $A_{\nu}$ is an $m \times n$ matrix such that $A_{\nu}/\widehat{E}_{\nu} \rightarrow Dw_a(0) = 0$.  Thus by~\cite[Lemma~6.5]{KrumWica}, Theorem~\ref{frequency monotonicity thm}(b), \eqref{lemma2_4 hyp4}, and \eqref{lemma2_4 eqn7}
\begin{align*}
	N_{w,0}(1/8) =\,& \lim_{\nu\rightarrow\infty} N_{T_{\nu},P_0,0}(R_{\nu}/8) \leq \lim_{\nu\rightarrow\infty} N_{T_{\nu},P_0,0}(R_{\nu}/2) \leq \alpha , \\
	N_{w,0}(1/8) =\,& \lim_{\nu\rightarrow\infty} N_{T_{\nu},P_{\nu},Z_{\nu}}(R_{\nu}/8) \geq \alpha + \varepsilon^2/4 , 
\end{align*}
where $N_{w,0}$ is the frequency function of $w$ at the origin (as in~\cite[Subsection~3.6]{KrumWica}), giving us a contradiction.

\noindent\textit{Proof of Conclusion~(ii).}  Fix $\varepsilon \in (0,1)$.  For $\nu = 1,2,3,\ldots$ let $\delta_{\nu} \rightarrow 0^+$, $R_{\nu} \rightarrow \infty$, $\alpha_{\nu} \in [1+1/q,\kappa]$ with $|\alpha_{\nu} - 2| \geq 1/q$, $T_{\nu}$ be an $n$-dimensional locally area-minimizing rectifiable current of $\mathbf{C}_{R_{\nu}}(0)$, $Z_{\nu} \in \mathbf{C}_1(0)$ with $\Theta(T_{\nu},Z_{\nu}) = q$, $P_{\nu}$ be an $n$-dimensional plane, and $\rho_{\nu} \in (0,1]$ such that \eqref{lemma2_4 hyp1}--\eqref{lemma2_4 hyp7} hold true with $\delta_{\nu},R_{\nu},\alpha_{\nu},T_{\nu},Z_{\nu},P_{\nu}$ in place of $\delta(\varepsilon),R(\varepsilon),\alpha,T,Z,P_Z$ and 
\begin{equation}\label{lemma2_4 eqn8}
	\lim_{\nu\rightarrow\infty} \sup_{\sigma \in (0,R_{\nu}/8]} |N_{T_{\nu},P_{\nu},Z_{\nu}}(\sigma) - \alpha_{\nu}| = 0.
\end{equation}
In light of the arbitrary choice of sequences, it suffices to show that for infinitely many $\nu$, conclusion~(ii) holds true with $\alpha_{\nu},T_{\nu},Z_{\nu},P_{\nu},\rho_{\nu}$ in place of $\alpha,T,Z,P_Z,\rho$.

After passing to a subsequence, let $\alpha_{\nu} \rightarrow \alpha$ in $\mathbb{R}$ as $\nu\rightarrow\infty$ and note that $\alpha > 1$ and $\alpha \neq 2$.  For each $\nu$, let $\Gamma_{\nu}$ be a rotation of $\mathbb{R}^{n+m}$ such that $\Gamma_{\nu}(P_{\nu}) = P_0$.  Set $\widetilde{T}_{\nu} = \Gamma_{\nu\#} \eta_{Z_{\nu},\rho_{\nu}\#} T_{\nu}$ and set $\widetilde{E}_{\nu} = \widehat{E}(T_{\nu},P_{\nu},\mathbf{C}_{\rho_{\nu}}(Z_{\nu},P_{\nu}))$.  By \eqref{lemma2_4 eqn6} (with $\rho = \rho_{\nu}$), $\widetilde{E}_{\nu} \rightarrow 0$ as $\nu \rightarrow \infty$.  By~\cite[Lemma~3.5]{KrumWica} and \eqref{lemma2_4 eqn6} (with $\rho = 8\rho_{\nu}$) 
\begin{equation*}
	\sup_{X \in \op{spt} T_{\nu} \cap \mathbf{C}_{6\rho_{\nu}}(Z_{\nu},P_{\nu})} \op{dist}(X,Z_{\nu}+P_{\nu}) \leq C \delta_{\nu} \rho_{\nu} , 
\end{equation*}
where $C = C(n,m,q,\mu) \in (0,\infty)$ is a constant.  Hence by the monotonicity formula and \eqref{lemma2_4 eqn5} 
\begin{equation}\label{lemma2_4 eqn9}
	\|T_{\nu}\|(\mathbf{C}_{4\rho_{\nu}}(Z_{\nu},P_{\nu})) 
	\leq \|T_{\nu}\|(\mathbf{B}_{4\rho_{\nu} + C \delta_{\nu} \rho_{\nu}}(Z_{\nu})) 
	\leq (q+C\delta_{\nu}^2) \,\omega_n (4\rho_{\nu} + C \delta_{\nu} \rho_{\nu})^n
\end{equation}
where $C = C(n,m,q,\mu) \in (0,\infty)$ are constants.  Noting that by \eqref{lemma2_4 eqn8} we have $N_{T_{\nu},P_{\nu},Z_{\nu}}(4\rho_{\nu}) \leq \kappa+1$ for all large $\nu$ and using $\Theta(T_{\nu},Z_{\nu}) = q$, \eqref{lemma2_4 eqn3}, \eqref{lemma2_4 eqn5}, and \eqref{lemma2_4 eqn6}, we can apply Theorem~\ref{frequency monotonicity thm}(e) to obtain 
\begin{equation}\label{lemma2_4 eqn10}
	\widehat{E}(T_{\nu},P_{\nu},\mathbf{C}_{4\rho_{\nu}}(Z_{\nu},P_{\nu})) \leq C \widetilde{E}_{\nu}
\end{equation}
for some constant $C = C(n,m,q,\mu,\kappa) \in (0,\infty)$.  By $\Theta(T_{\nu},Z_{\nu}) = q$, \eqref{lemma2_4 eqn3}, \eqref{lemma2_4 eqn4}, \eqref{lemma2_4 eqn9}, \eqref{lemma2_4 eqn10}, and $\widetilde{E}_{\nu} \rightarrow 0$ as $\nu\rightarrow\infty$, we can follow the blow-up procedure of~\cite[Subsection~6.1]{KrumWica} to construct a blow-up $\widetilde{\varphi} \in W^{1,2}(B_1(0),\mathcal{A}_q(\mathbb{R}^m))$ of $\widetilde{T}_{\nu}$ with respect to the plane $P_0$ and excess $\widetilde{E}_{\nu}$.  That is, there exists Lipschitz functions $\widetilde{u}^{(\nu)} : B_1(0) \rightarrow \mathcal{A}_q(\mathbb{R}^m)$ and sets $\widetilde{K}_{\nu} \subset B_1(0)$ such that 
\begin{gather}\label{lemma2_4 eqn11}
	\widetilde{T}_{\nu} \llcorner (\widetilde{K}_{\nu} \times \mathbb{R}^m) 
		= (\op{graph} \widetilde{u}^{(\nu)}) \llcorner (\widetilde{K}_{\nu} \times \mathbb{R}^m) , \\
	\mathcal{L}^n(B_1(0) \setminus \widetilde{K}_{\nu}) + \|\widetilde{T}_{\nu}\|((B_1(0) \setminus \widetilde{K}_{\nu}) \times \mathbb{R}^m) 
		\leq C \widetilde{E}_{\nu}^{2+\varsigma} , \nonumber \\
	\sup_{B_1(0)} |\widetilde{u}^{(\nu)}| \leq C \widetilde{E}_{\nu}, \quad 
		\sup_{B_1(0)} |D\widetilde{u}^{(\nu)}| \leq C \widetilde{E}_{\nu}^{2+\varsigma} , \nonumber 
\end{gather}
where $\varsigma = \varsigma(n,m,q) \in (0,1)$ is as in Theorem~\ref{lip approx thm} and $C = C(n,m,q,\mu,\kappa) \in (0,\infty)$ is a constant.  It follows from \eqref{lemma2_4 eqn11} and~\cite[Lemma~3.7]{KrumWica} that 
\begin{equation*}
	\sup_{B_1(0)} |\widetilde{u}^{(\nu)}|^2 + \int_{B_1(0)} |D\widetilde{u}^{(\nu)}|^2 \leq C \widetilde{E}_{\nu}^2 
\end{equation*}
for all sufficiently large $\nu$, where $C = C(n,m,q,\mu,\kappa) \in (0,\infty)$ is a constant (see (6.8) in~\cite{KrumWica}).  Hence by~\cite[Proposition~2.11]{DeLSpaDirMin} and~\cite[Theorem~5.2]{DeLSpa1}, after passing to a subsequence we find a $q$-valued function $\widetilde{\varphi} \in W^{1,2}(B_1(0),\mathcal{A}_q(\mathbb{R}^m))$ such that 
\begin{equation}\label{lemma2_4 eqn12}
	\widetilde{u}^{(\nu)}/\widetilde{E}_{\nu} \rightarrow \widetilde{\varphi}
\end{equation}
pointwise $\mathcal{L}^n$-a.e.~on $B_1(0)$ and strongly in $L^2(B_1(0),\mathcal{A}_q(\mathbb{R}^m))$ as $\nu\rightarrow\infty$ and $|D\widetilde{u}^{(\nu)}|/\widetilde{E}_{\nu} \rightarrow |D\widetilde{\varphi}|$ strongly in $L^2(B_1(0))$ as $\nu\rightarrow\infty$.  Let $\widetilde{\varphi}_a(x) = \frac{1}{q} \sum_{l=1}^q \widetilde{\varphi}_l(x)$ for each $x \in B_1(0)$, where $\widetilde{\varphi}(x) = \sum_{l=1}^q \llbracket \widetilde{\varphi}_l(x) \rrbracket$.  By \eqref{lemma2_4 eqn11} and \eqref{lemma2_4 eqn12}
\begin{equation*}
	\int_{B_1(0)} |\widetilde{\varphi}|^2 
	= \lim_{\nu\rightarrow\infty} \frac{1}{\widetilde{E}_{\nu}^2} \int_{B_1(0)} |\widetilde{u}^{(\nu)}|^2 
	= \lim_{\nu\rightarrow\infty} \frac{1}{\widetilde{E}_{\nu}^2 \rho_{\nu}^{n+2}} \int_{\mathbf{C}_{\rho_{\nu}}(Z_{\nu},P_{\nu})} 
		\op{dist}^2(X,Z_{\nu}+P_{\nu}) \,d\|T_{\nu}\|(X) 
	= 1 
\end{equation*}
and thus $\widetilde{\varphi}$ is non-zero.  Since $\Theta(\widetilde{T}_{\nu},0) = \Theta(T_{\nu},Z_{\nu}) = q$, by~\cite[Lemma~6.4]{KrumWica} we have that $\widetilde{\varphi}_a(0) = 0$.  By~\cite[Lemma~6.5]{KrumWica} and \eqref{lemma2_4 eqn8} 
\begin{equation*}
	N_{\widetilde{\varphi},0}(\sigma) = \lim_{\nu\rightarrow\infty} N_{T_{\nu},P_{\nu},Z_{\nu}}(\rho_{\nu} \sigma) = \alpha 
\end{equation*}
for all $\sigma \in (0,\infty)$.  Therefore, $\widetilde{\varphi}$ is homogeneous degree $\alpha$.  Extend $\widetilde{\varphi}$ to a homogeneous degree $\alpha$ function on $\mathbb{R}^n$.

For each positive integer $\nu$, define $\varphi^{(\nu)} \in W^{1,2}_{\rm loc}(P_{\nu},\mathcal{A}_q(P_{\nu}^{\perp}))$ by $\varphi^{(\nu)}(X) = \widetilde{E}_{\nu} \rho_{\nu}^{1-\alpha} \,\Gamma_{\nu}^{-1}\widetilde{\varphi}(\Gamma_{\nu} X)$ for each $X \in P_{\nu}$ (recalling that we identify $P_0 \cong \mathbb{R}^n$ and $P_0^{\perp} \cong \mathbb{R}^m$).  Clearly $\varphi^{(\nu)}$ is homogeneous degree $\alpha$ and $|\mathcal{N}_{\varphi^{(\nu)}}(0) - \alpha_{\nu}| = |\alpha - \alpha_{\nu}| < \varepsilon^2$ for all sufficiently large $\nu$.  By \eqref{lemma2_4 eqn11} and \eqref{lemma2_4 eqn12} 
\begin{align}\label{lemma2_4 eqn13}
	&\lim_{\nu\rightarrow\infty} \frac{1}{\widetilde{E}_{\nu}^2 \rho_{\nu}^{n+2}} \int_{\mathbf{C}_{\rho_{\nu}}(Z_{\nu},P_{\nu})} 
		\op{dist}^2(X,Z_{\nu}+\op{spt}\op{graph}\varphi^{(\nu)}) \,d\|T_{\nu}\|(X) 
	\\ =\,& \lim_{\nu\rightarrow\infty} \frac{1}{\widetilde{E}_{\nu}^2} \int_{\mathbf{C}_1(0)} 
		\op{dist}^2(X, \op{spt}\op{graph}(\widetilde{E}_{\nu} \widetilde{\varphi})) \,d\|\widetilde{T}_{\nu}\|(X) \nonumber 
	\\ =\,& \lim_{\nu\rightarrow\infty} \frac{1}{\widetilde{E}_{\nu}^2} \int_{B_1(0)} \sum_{l=1}^q 
		\op{dist}^2(x + \widetilde{u}^{(\nu)}_l(x), \op{spt}\op{graph}(\widetilde{E}_{\nu} \widetilde{\varphi})) \,d\mathcal{L}^n(x) \nonumber 
	\\ \leq\,& \lim_{\nu\rightarrow\infty} \frac{1}{\widetilde{E}_{\nu}^2} \int_{B_1(0)} \mathcal{G}(\widetilde{u}^{(\nu)}(x), 
		\widetilde{E}_{\nu} \widetilde{\varphi}(x))^2 \,d\mathcal{L}^n(x) = 0, \nonumber
\end{align}
where $\widetilde{u}^{(\nu)}(x) = \sum_{l=1}^q \llbracket \widetilde{u}^{(\nu)}_l(x) \rrbracket$ for each $x \in B_1(0)$.  By similar reasoning, 
\begin{align}\label{lemma2_4 eqn14}
	&\lim_{\nu\rightarrow\infty} \frac{1}{\widetilde{E}_{\nu}^2 \rho_{\nu}^{n+2}} \int_{B_{\rho_{\nu}}(0,P_{\nu})} \sum_{l=1}^q 
		\op{dist}^2(X+Z_{\nu}+\varphi^{(\nu)}_l(X), \op{spt} T_{\nu}) \,d\mathcal{H}^n(X) 
	\\ \leq\,& \lim_{\nu\rightarrow\infty} \frac{1}{\widetilde{E}_{\nu}^2} \int_{B_1(0)} \mathcal{G}(\widetilde{u}^{(\nu)}(x), 
		\widetilde{E}_{\nu} \widetilde{\varphi}(x))^2 \,d\mathcal{L}^n(x) 
	= 0, \nonumber 
\end{align}
where $\varphi^{(\nu)}(x) = \sum_{l=1}^q \llbracket \varphi^{(\nu)}_l(x) \rrbracket$ for each $x \in B_1(0)$, completing the proof that conclusion~(ii) holds true with $T_{\nu},Z_{\nu},P_{\nu},\rho_{\nu},\varphi^{(\nu)}$ in place of $T,Z,P_Z,\rho,\varphi$ for all sufficiently large $\nu$.

\noindent\textit{Proof of Conclusion~(iii).}  We argue exactly as we did for conclusion~(ii), except we observe the following.  If $\op{dim}\op{spine}\widetilde{\varphi} = n-2$, then arguing as we did in \eqref{lemma2_4 eqn13} and \eqref{lemma2_4 eqn14}, we have that \eqref{lemma2_4 concl2} holds true with $T_{\nu},Z_{\nu},P_{\nu},\rho_{\nu},\varphi^{(\nu)}$ in place of $T,Z,P_Z,\rho,\varphi$ for all sufficiently large $\nu$.  If instead $\op{dim}\op{spine}\widetilde{\varphi} \leq n-3$, let $L_{\nu}$ be an $(n-3)$-dimensional linear subspace such that $\op{spine}\varphi^{(\nu)} \subseteq L_{\nu}$.  Since $\widetilde{\varphi}$ is homogeneous degree $\alpha$ where $\alpha > 1$ and $\alpha \neq 2$, we can argue as in the proof of Theorem~\ref{B prime stratification thm} that 
\begin{equation*}
	\op{spine}\widetilde{\varphi} = \{ \xi \in \mathbb{R}^n : \mathcal{N}_{\widetilde{\varphi}-\ell_{\xi}}(\xi) = \alpha \} ,
\end{equation*}
where $\ell_{\xi}(x) = \widetilde{\varphi}_a(\xi) + D\widetilde{\varphi}_a(\xi) \cdot (x-\xi)$ and $(\widetilde{\varphi}-\ell_{\xi})(x) = \sum_{l=1}^q \llbracket \widetilde{\varphi}_l(x) - \ell_{\xi}(x) \rrbracket$ for each $x,\xi \in \mathbb{R}^n$.  It follows using~\cite[Corollary~6.6]{KrumWica} that \eqref{lemma2_4 concl3} holds true with $T_{\nu},Z_{\nu},P_{\nu},\rho_{\nu},L_{\nu}$ in place of $T,Z,P_Z,\rho,L$ for all sufficiently large $\nu$.
\end{proof}

\section{Bounding the distance of an area-minimizing current to a union of non-intersecting single-valued nearly harmonic graphs}\label{sec:relharm}  

In this section, we derive local estimates for locally area-minimizing rectifiable currents relative to graphs of multi-valued harmonic functions, which will then be used in Section~\ref{sec:graph rep} to obtain a graphical representation of locally area-minimizing rectifiable currents relative to graphs of homogeneous multi-valued harmonic functions.  Let $T$ be an $n$-dimensional locally area-minimizing rectifiable current of $\mathbf{C}_1(0)$ which is weakly close to the plane $P_0 = \mathbb{R}^n \times \{0\}$ with multiplicity $q$.  Let $\varphi : B_1(0) \rightarrow \mathcal{A}_s(\mathbb{R}^m)$ be a multi-valued function given by $\varphi(x) = \sum_{i=1}^s \llbracket \varphi_i(x) \rrbracket$, where $s \geq 1$ is an integer and $\varphi_i : B_1(0) \rightarrow \mathbb{R}^m$ are single-valued harmonic functions such that either $\varphi_i \equiv \varphi_j$ on $B_1(0)$ or $\op{graph} \varphi_i \cap \op{graph} \varphi_j = \emptyset$ for all $i \neq j$.  That is, the support of the graph of $\varphi$ in $\mathbf{C}_1(0)$ is a union of pairwise disjoint graphs of single-valued harmonic functions.  (In order to apply our results in~\cite{KrumWicd}, we will in fact consider the more general case where each $\varphi_i$ is a $C^2$ single-valued functions, not necessarily harmonic.)  In this section, we extend the results of~\cite[Section~3]{KrumWicb} to $L^{\infty}$-estimates for the area-minimizer $T$ relative to the graph of $\varphi$.  In particular, in Theorem~\ref{relharm separation3 thm}, we show that if the components of $\varphi$ are well-separated as in \eqref{relharm varphi wellsep hyp}, then the $L^{\infty}$-distance of $T$ to the graph of $\varphi$ is bounded in terms of the $L^2$-distance of $T$ to the graph of $\varphi$.

One possible approach to proving Theorem~\ref{relharm separation3 thm} is to locally approximating the graph of $\varphi$ by its tangent planes at a point and apply~\cite[Theorem~3.11]{KrumWicb}.  However, the distance of the graph of $\varphi$ to its tangent plane is bounded in terms of $\|D^2 \varphi\|_{C^0}$, and we want to Theorem~\ref{relharm separation3 thm} to apply in the situation that $\|D^2 \varphi\|_{C^0}$ is much larger than the $L^2$-distance of $T$ to the graph of $\varphi$ or the minimum (non-zero) distance between the components $\varphi_i$ of $\varphi$.

Our approach in this section is to follow the strategy of the proof of~\cite[Theorem~3.11]{KrumWicb}.  This mostly goes through with obvious modifications; in particular, instead of considering the distance and tilt excess of $T$ relative to a union of planes as in~\cite{KrumWicb}, we now consider the distance and tilt excess of $T$ relative to a sum of harmonic functions $\varphi(x) = \sum_{i=1}^s \llbracket \varphi_i(x) \rrbracket$ as above (also see Definition~\ref{parallel harmonic defn} and Hypothesis~\ref{relharm varphi defn}).  Notable changes include a harmonic approximation result (Theorem~\ref{relharm harmonic approx lemma}), a simplification of the inductive argument, and an $L^{\infty}$-estimate for $T$ relative to a single-valued harmonic function $\varphi$ (Lemma~\ref{relharm separation0 lemma}).  In contrast with~\cite[Lemma~3.18]{KrumWicb}, which uses the fact that coordinate function $x_i$ are harmonic on the minimal submanifold $T$, Lemma~\ref{relharm separation0 lemma} does not follow immediately from~\cite[Theorem~7.5(6)]{Allard} since $\varphi$ is harmonic function on $\mathbb{R}^n$ but not on $T$.  Rather, we derive Lemma~\ref{relharm separation0 lemma} directly from Theorem~\ref{relharm separation1 thm}, in which we bound the $L^{\infty}$-distance of $T$ to the graph of $\varphi$ in terms of the $L^2$ tilt excess of $T$ relative to the graph of $\varphi$.

\subsection{Energy estimate for an area-minimizer relative to a single-valued nearly flat $C^{2}$ function}  Throughout this section, we let $X = (x,y)$ denote a point in $\mathbb{R}^{n+m}$, where $x \in \mathbb{R}^n$ and $y \in \mathbb{R}^m$.  We have the following energy estimate for a locally area-minimizing rectifiable current $T$ relative to the graph of a $C^2$ single-valued function $\varphi$.

\begin{lemma}[Energy estimate]\label{relharm energy lemma} 
Let $q$ be a positive integer and $0 < \gamma < \overline{\gamma} < 1$.  There exists $\varepsilon = \varepsilon(n,m,q,\overline{\gamma}) \in (0,1)$ such that the following holds true.  Let $x_0 \in \mathbb{R}^n$, and $\rho > 0$.  Let $T$ be an $n$-dimensional locally area-minimizing rectifiable current of $\mathbf{C}_{\rho}(x_0)$ such that 
\begin{equation*}
	(\partial T) \llcorner \mathbf{C}_{\rho}(x_0) = 0, \quad \sup_{X \in \op{spt} T} \op{dist}(X,P_0) < \infty, \quad
	\pi_{\#} T = q \llbracket B_{\rho}(x_0) \rrbracket, \quad \widehat{\mathcal E}(T,P_0,\mathbf{C}_{\rho}(x_0)) < \varepsilon .  
\end{equation*}
Let $u : B_{\overline{\gamma}\rho}(x_0) \rightarrow \mathcal{A}_q(\mathbb{R}^m)$ be a Lipschitz approximation of $T$ and $K \subseteq B_{\overline{\gamma}\rho}(x_0)$ be a closed set as in Theorem~\ref{lip approx thm} with $\overline{\gamma}$ in place of $\gamma$.  Let $\varphi \in C^2(B_{\rho}(x_0), \mathbb{R}^m)$ be a single-valued function with $\|D\varphi\|_{C^0(B_{\rho}(x_0))} < \varepsilon^2$.  Then 
\begin{align}\label{relharm energy concl}
	&\frac{1}{\rho^n} \int_{B_{\gamma\rho}(x_0)} \mathcal{G}(Du,q\llbracket D\varphi \rrbracket)^2 
		+ \frac{1}{\rho^n} \int_{G_n(\mathbf{C}_{\gamma\rho}(x_0))} \sum_{\kappa=1}^m |\pi_S(-D\varphi^{\kappa}(x),e_{\kappa})|^2 \,d|T|(X,S) 
	\\ \leq\,& \frac{C}{\rho^{n+2}} \int_{\mathbf{C}_{\overline{\gamma}\rho}(x_0)} |y - \varphi(x)|^2 \,d\|T\|(X) 
		+ C\rho^2 \|\Delta\varphi\|_{C^0(B_{\rho}(x_0))}^2 + C \widehat{\mathcal E}(T,P_0,\mathbf{C}_{\rho}(x_0))^{2+\varsigma} , \nonumber
\end{align}
where $|T|$ is the $n$-dimensional integral varifold associated with $T$, $e_1,\ldots,e_m$ is the standard basis of $\mathbb{R}^m$, $\varsigma = \varsigma(n,m,q) \in (0,1)$ is as in Theorem~\ref{lip approx thm}, and $C = C(n,m,q,\gamma,\overline{\gamma}) \in (0,\infty)$ is a constant.
\end{lemma}

\begin{proof}
Without loss of generality assume $x_0 = 0$ and $\rho = 1$.  By translating, we may assume that $\varphi(0) = 0$ so that $\|\varphi\|_{C^1(B_1(0))} < \varepsilon$.  Since $\widehat{\mathcal E}(T,P_0,\mathbf{C}_{\rho}(x_0)) < \varepsilon$ and $\|D\varphi\|_{C^0(B_1(0))} < \varepsilon$, the left-hand side of \eqref{relharm energy concl} is $\leq C(n,m,q,\gamma,\overline{\gamma})\,\varepsilon^2$.  Hence we may assume that 
\begin{equation*}
    \int_{\mathbf{C}_{\overline{\gamma}}(0)} |y - \varphi(x)|^2 \,d\|T\|(X) < \varepsilon^2 
\end{equation*}
as otherwise \eqref{relharm energy concl} must hold true.  Hence using $\|\varphi\|_{C^1(B_1(0))} < \varepsilon$ and~\cite[Lemma~3.5]{KrumWica} 
\begin{equation}\label{relharm energy eqn1}
    \sup_{X \in \op{spt}T \cap \mathbf{C}_{(\gamma+\overline{\gamma})/2}(0)} \op{dist}^2(X,P_0) \leq C \int_{\mathbf{C}_{\overline{\gamma}}(0)} |y|^2 \,d\|T\|(X) < C \varepsilon^2 ,
\end{equation}
where $C = C(n,m,q,\gamma,\overline{\gamma}) \in (0,\infty)$ is a constant.  It follows that $|u| \leq C\varepsilon$ on $K$, which since $\op{Lip} u \leq C\varepsilon^{\varsigma}$ guarantees that $|u| \leq C\varepsilon^{\varsigma}$ on $B_{\overline{\gamma}}(0)$, where $C = C(n,m,q,\gamma,\overline{\gamma}) \in (0,\infty)$ are constants.

Recall that since $T$ is a locally area-minimizing rectifiable current of $\mathbf{C}_1(0)$ with $(\partial T) \llcorner \mathbf{C}_1(0) = 0$, the integral varifold $|T|$ associated with $T$ is area-stationary in the sense that 
\begin{equation}\label{first variation}
	\int_{G_n(\mathbf{C}_1(0))} \op{div}_S \zeta(X) \,d|T|(X,S) = 0 
\end{equation}
for all $\zeta \in C^1_c(\mathbf{C}_1(0),\mathbb{R}^{n+m})$, where $\op{div}_S$ denotes the divergence with respect to the plane $S$ (see~\cite[Definition~16.3 and Lemma~33.2]{SimonGMT}).  Let $\eta \in C^1_c(B_1(0))$ such that $0 \leq \eta \leq 1$ on $B_1(0)$, $\eta = 1$ on $B_{\gamma}(0)$, $\eta = 0$ on $B_1(0) \setminus B_{(\gamma+\overline{\gamma})/2}(0)$, and $|D\eta| \leq \frac{4}{\overline{\gamma}-\gamma}$.  Extend $\eta(x)$ to a function $\eta(x,y)$ of $x \in B_1(0)$ and $y \in \mathbb{R}^m$ which is independent of $y$.  Note that for $\kappa = 1,2,\ldots,m$, the gradient of $y^{\kappa} - \varphi^{\kappa}(x)$ is $(-D\varphi^{\kappa}(x),e_{\kappa})$ at each $(x,y) \in \mathbb{R}^{n+m}$.  In \eqref{first variation}, setting $\zeta(x,y) = \eta(x)^2 \,(0,y - \varphi(x))$ for each $(x,y) \in B_1(0) \times \mathbb{R}^m$ near the support of $T$ gives us 
\begin{align*}
	0 =\,& \int_{G_n(\mathbf{C}_1(0))} \sum_{\kappa=1}^m \eta(x)^2 \,\pi_S (0,e_{\kappa}) \cdot \pi_S (-D\varphi^{\kappa}(x),e_{\kappa}) \,d|T|(X,S) 
		\\&+ 2 \int_{G_n(\mathbf{C}_1(0))} \sum_{\kappa=1}^m \eta(x) \,(y^{\kappa} - \varphi^{\kappa}(x)) \,\nabla^S \eta(x,y) \cdot \pi_S (0,e_{\kappa}) 
		\,d|T|(X,S) \nonumber 
\end{align*}
where $\nabla^S = \pi_S \circ \nabla$ is the tangential gradient with respect to $S$.  Using Theorem~\ref{lip approx thm}, the area formula, $\|\varphi\|_{C^1(B_1(0))} < \varepsilon$, and \eqref{relharm energy eqn1} 
\begin{equation}\label{relharm energy eqn2}
	\int_{B_{\overline{\gamma}}(0)} \sum_{l=1}^q \eta^2 Du_l \cdot (Du_l - D\varphi) 
	\leq -2 \int_{B_{\overline{\gamma}}(0)} \sum_{l=1}^q \sum_{i=1}^n \eta \,D_i \eta \,(u_l - \varphi) \cdot D_i u_l 
		+ C \widehat{\mathcal E}(T,P_0,\mathbf{C}_1(0))^{2+\varsigma} , 
\end{equation}
where $u(x) = \sum_{l=1}^q \llbracket u_l(x) \rrbracket$ and $Du(x) = \sum_{l=1}^q \llbracket Du_l(x) \rrbracket$ at $\mathcal{L}^n$-a.e.~$x \in B_{\overline{\gamma}}(0)$ following the convention from~\cite[Subsection~3.5]{KrumWica} and where $C = C(n,m,q,\gamma,\overline{\gamma}) \in (0,\infty)$ is a constant.  On the other hand, by the divergence theorem 
\begin{equation}\label{relharm energy eqn3}
	\int_{B_{\overline{\gamma}}(0)} \eta^2 D\varphi \cdot (Du_a - D\varphi) = -2 \int_{B_{\overline{\gamma}}(0)} \sum_{i=1}^n \eta \,D_i \eta \,(u_a - \varphi) \cdot D_i \varphi 
		- \int_{B_{\overline{\gamma}}(0)} \eta^2 \,(u_a - \varphi) \cdot \Delta\varphi ,
\end{equation}
where $u_a(x) = \frac{1}{q} \sum_{l=1}^q u_l(x)$ for each $x \in B_{\overline{\gamma}}(0)$.  Subtracting \eqref{relharm energy eqn2} and \eqref{relharm energy eqn3} gives us 
\begin{align*}
	\int_{B_{\overline{\gamma}}(0)} \eta^2 \mathcal{G}(Du,q\llbracket D\varphi \rrbracket)^2 
	\leq\,& -2 \int_{B_{\overline{\gamma}}(0)} \sum_{l=1}^q \sum_{i=1}^n \eta \,D_i \eta \,(u_l - \varphi) \cdot (D_i u_l - D_i \varphi)
		\\&- q \int_{B_{\overline{\gamma}}(0)} \eta^2 \,(u_a - \varphi) \cdot \Delta\varphi + C \widehat{\mathcal E}(T,P_0,\mathbf{C}_1(0))^{2+\varsigma} , 
\end{align*}
Using Cauchy's inequality 
\begin{equation}\label{relharm energy eqn4}
	\int_{B_{\overline{\gamma}}(0)} \eta^2 \,\mathcal{G}(Du,q\llbracket D\varphi \rrbracket)^2 
	\leq C \int_{B_{\overline{\gamma}}(0)} \mathcal{G}(u,q\llbracket \varphi \rrbracket)^2 + C \|\Delta\varphi\|_{C^0(B_1(0))}^2 
		+ C \widehat{\mathcal E}(T,P_0,\mathbf{C}_1(0))^{2+\varsigma} ,  
\end{equation}
where $C = C(n,m,q,\gamma,\overline{\gamma}) \in (0,\infty)$ is a constant.  

Recalling Rademacher's theorem~\cite[Theorem~1.13]{DeLSpaDirMin}, let $x \in K$ be a point at which $K$ has density one and $u$ is strongly differentiable.  At $(x,u_l(x))$, since the gradient of $u$ is small and $(-Du^{\kappa}_l(x),e_{\kappa})$ is orthogonal to $\op{graph} u$ 
\begin{equation*}
	|\pi_S(-D\varphi^{\kappa}(x),e_{\kappa})| 
	= |\pi_S((-D\varphi^{\kappa}(x),e_{\kappa}) - (-Du^{\kappa}_l(x),e_{\kappa}))| 
	= |\pi_S(Du^{\kappa}_l(x) - D\varphi^{\kappa}(x), 0)|  
\end{equation*}
and thus 
\begin{equation}\label{relharm energy eqn5}
	\tfrac{1}{2} |Du^{\kappa}_l(x) - D\varphi^{\kappa}(x)| \leq |\pi_S(-D\varphi^{\kappa}(x),e_{\kappa})| \leq |Du^{\kappa}_l(x) - D\varphi^{\kappa}(x)| ,
\end{equation}
where $S$ is the tangent plane to $T$ at $(x,u_l(x))$ and $\varphi^{\kappa}(x)$ and $u^{\kappa}(x)$ are the $\kappa$-th coordinates of $\varphi(x)$ and $u_l(x)$ respectively.  Using Theorem~\ref{lip approx thm} and the area formula, \eqref{relharm energy eqn4} gives us \eqref{relharm energy concl}
\end{proof}

\subsection{Harmonic approximation of an area-minimizing current relative to a single-valued nearly harmonic graph}  We define the following excess $\mathcal{E}_{\varsigma,\varepsilon}(T,\varphi,\mathbf{C}_{\rho}(x_0))$ of $T$ relative to the graph of $C^2$ single-valued $\varphi$.

\begin{definition}\label{relharm excess defn1} 
Let $q$ be a positive integer, $\varsigma,\varepsilon \in (0,1]$, $x_0 \in \mathbb{R}^n$, and $\rho > 0$.  Let $T$ be an $n$-dimensional rectifiable current of $\mathbf{C}_{\rho}(x_0)$ and let $\varphi \in C^2(B_{\rho}(x_0), \mathbb{R}^m)$.  We define 
\begin{align}\label{relharm excess defn1 eqn}  
	\mathcal{E}_{\varsigma,\varepsilon}(T,\varphi,\mathbf{C}_{\rho}(x_0))^2 =\,& 
		\frac{1}{\omega_n \rho^n} \int_{G_n(\mathbf{C}_{\rho}(x_0))} \sum_{\kappa=1}^m |\pi_S(-D\varphi^{\kappa}(x),e_{\kappa})|^2 \,d|T|(X,S) 
		\\&+ \frac{1}{\varepsilon^2} \rho^2 \,\|\Delta\varphi\|_{C^0(B_{\rho}(x_0))}^2 
		+ \|D\varphi\|_{C^0(B_{\rho}(x_0))} \,\widehat{\mathcal E}(T,P_0,\mathbf{C}_{\rho}(x_0))^2 \nonumber
		\\&+ \widehat{\mathcal E}(T,P_0,\mathbf{C}_{\rho}(x_0))^{2+\varsigma} . \nonumber
\end{align}
where $|T|$ is the $n$-dimensional varifold associated with $T$ and $\widehat{\mathcal E}(T,P_0,\mathbf{C}_{\rho}(x_0))$ is as in \eqref{oriented tilt defn}.  When $\varepsilon = 1$, we let $\mathcal{E}_{\varsigma}(T,\varphi,\mathbf{C}_{\rho}(x_0)) = \mathcal{E}_{\varsigma,1}(T,\varphi,\mathbf{C}_{\rho}(x_0))$.
\end{definition}

The first term on the right-hand side of \eqref{relharm excess defn1 eqn} quantifies the tilt excess of $T$ relative to the graph of $\varphi$ (cf.~the proof of Lemma~\ref{relharm energy lemma}).  The second term $\varepsilon^{-1} \rho \|\Delta\varphi\|_{C^0(B_{\rho}(x_0))}$ vanishes if $\varphi$ is harmonic and otherwise guarantees that $\Delta\varphi$ is negligible in the sense that $\rho\|\Delta\varphi\|_{C^0(B_{\rho}(x_0))} \leq \varepsilon \,\mathcal{E}_{\varsigma,\varepsilon}(T,\varphi,\mathbf{C}_{\rho}(x_0))$.  The remaining terms in \eqref{relharm excess defn1 eqn} are error terms which arise when using Almgren's strong Lipschitz approximation theorem and in the harmonic approximation theorem Lemma~\ref{relharm harmonic approx lemma} below.

\begin{remark}\label{relharm excess rmk1}{\rm
(1)  Assuming that $T \neq 0$ in $\mathbf{C}_{\rho}(x_0)$, $\mathcal{E}_{\varsigma,\varepsilon}(T,\varphi,\mathbf{C}_{\rho}(x_0)) = 0$ if and only if $\varphi$ is a constant function and $T$ is a sum of $n$-dimensional planes which are parallel to $P_0$ and have orientation $\vec P_0$.

\noindent (2)  Given $\varepsilon \in (0,1)$, if $\widehat{\mathcal E}(T,P_0,\mathbf{C}_{\rho}(x_0)) < \varepsilon$ and $\|D\varphi\|_{C^0(B_{\rho}(x_0))} + \rho \|\Delta\varphi\|_{C^0(B_{\rho}(x_0))} < \varepsilon^2$, then $\mathcal{E}_{\varsigma,\varepsilon}(T,\varphi,\mathbf{C}_{\rho}(x_0)) < C(n,m,q) \,\varepsilon$.  

\noindent (3)  Conversely, given $\varepsilon_0, \varepsilon > 0$, if $\mathcal{E}_{\varsigma,\varepsilon}(T,\varphi,\mathbf{C}_{\rho}(x_0)) < \varepsilon_0$, then $\widehat{\mathcal E}(T,P_0,\mathbf{C}_{\rho}(x_0)) < \varepsilon_0^{\frac{2}{2+\varsigma}}$.
}\end{remark}

\begin{definition}\label{parallel harmonic defn} 
Let $\varphi : B_{\rho}(x_0) \rightarrow \mathbb{R}^m$ be a single-valued function and let $w : B_{\rho}(x_0) \rightarrow \mathcal{A}_q(\mathbb{R}^m)$ be a $q$-valued function.  We define $(\varphi+w) : B_{\rho}(x_0) \rightarrow \mathcal{A}_q(\mathbb{R}^m)$
\begin{equation*}%\label{varphi plus w defn}
	(\varphi+w)(x) = \sum_{l=1}^q \llbracket \varphi(x) + w_l(x) \rrbracket
\end{equation*}
for each $x \in B_{\rho}(x_0)$, where $w(x) = \sum_{l=1}^q \llbracket w_l(x) \rrbracket$.  In particular, given $a = \sum_{l=1}^q \llbracket a_l \rrbracket \in \mathcal{A}_q(\mathbb{R}^m)$, we define $(\varphi+a) : B_{\rho}(x_0) \rightarrow \mathcal{A}_q(\mathbb{R}^m)$
\begin{equation}\label{varphi plus a defn}
	(\varphi+a)(x) = \sum_{l=1}^q \llbracket \varphi(x) + a_l \rrbracket
\end{equation}
for each $x \in B_{\rho}(x_0)$.
\end{definition}

\begin{definition}\label{a sep width defn} 
Let $a = \sum_{i=1}^N q_i \llbracket a_i \rrbracket \in \mathcal{A}_q(\mathbb{R}^m)$, where $N,q_1,\ldots,q_N$ are positive integers and $a_1,\ldots,a_N \in \mathbb{R}^m$ are distinct.  If $N \geq 2$, we define the separation $\op{sep} a$ and width $\op{width} a$ by 
\begin{equation*}
	\op{sep} a = \min_{i \neq j} |a_i - a_j|, \quad \op{width} a = \max_{i \neq j} |a_i - a_j| .
\end{equation*}
If instead $N = 1$, we define $\op{sep} a = \infty$ and $\op{width} a = 0$.
\end{definition}

We have the following harmonic approximation theorem for an area-minimizing rectifiable current $T$ relative to the graph of $C^2$ single-valued $\varphi$.

\begin{lemma}[Harmonic approximation lemma]\label{relharm harmonic approx lemma} 
Let $q$ be a positive integer.  For every $\eta > 0$ and $\gamma \in (0,1)$ there exists $\varepsilon = \varepsilon(n,m,q,\gamma,\eta) \in (0,1)$ such that the following holds true.  Let $x_0 \in \mathbb{R}^n$, and $\rho > 0$.  Let $T$ be an $n$-dimensional locally area-minimizing rectifiable current of $\mathbf{C}_{\rho}(x_0)$ such that 
\begin{align}\label{relharm harmonic approx hyp}
	&(\partial T) \llcorner \mathbf{C}_{\rho}(x_0) = 0, \quad\quad \sup_{X \in \op{spt} T} \op{dist}(X,P_0) < \infty, 
	\\& \pi_{\#} T = q \llbracket B_{\rho}(x_0) \rrbracket, \quad\quad \widehat{\mathcal E}(T,P_0,\mathbf{C}_{\rho}(x_0)) < \varepsilon . \nonumber 
\end{align}
Let $u : B_{(1+\gamma)\rho/2}(x_0) \rightarrow \mathcal{A}_q(\mathbb{R}^m)$ be a Lipschitz approximation of $T$ as in Theorem~\ref{lip approx thm} with $(1+\gamma)/2$ in place of $\gamma$.  Let $\varphi \in C^2(B_{\rho}(x_0), \mathbb{R}^m)$ be a single-valued function with $\|D\varphi\|_{C^0(B_{\rho}(x_0))} + \rho \|\Delta\varphi\|_{C^0(B_{\rho}(x_0))} < \varepsilon^2$.  Then there exists a Dirichlet energy minimizing $q$-valued function $w : B_{\gamma\rho}(x_0) \rightarrow \mathcal{A}_q(\mathbb{R}^m)$ such that 
\begin{equation}\label{relharm harmonic approx concl}
	\frac{1}{\rho^{n+2}} \int_{B_{\gamma\rho}(x_0)} \mathcal{G}(u,\varphi+w)^2 + \frac{1}{\rho^n} \int_{B_{\gamma\rho}(x_0)} (|D(u-\varphi)| - |Dw|)^2 
	\leq \eta \,\mathcal{E}_{\varsigma/2,\varepsilon}(T,\varphi,\mathbf{C}_{\rho}(x_0))^2 ,  
\end{equation}
where $\varsigma = \varsigma(n,m,q) \in (0,1)$ is as in Theorem~\ref{lip approx thm}.
\end{lemma}

\begin{proof}
Assume that $x_0 = 0$ and $\rho = 1$.  If $\mathcal{E}_{\varsigma/2,\varepsilon}(T,\varphi,\mathbf{C}_1(0)) = 0$, then $\widehat{\mathcal E}(T,P_0,\mathbf{C}_1(0)) = 0$ and thus the set $K \subseteq B_{(1+\gamma)/2}(0)$ as in Theorem~\ref{lip approx thm} (with $(1+\gamma)/2$ in place of $\gamma$) satisfies $K = B_{(1+\gamma)/2}(0)$.  It follows from Remark~\ref{relharm excess rmk1}(1) that $u$ and $\varphi$ are constant functions and thus \eqref{relharm harmonic approx concl} holds true with $w$ as a constant $q$-valued function.  Hence we may assume that $\mathcal{E}_{\varsigma/2,\varepsilon}(T,\varphi,\mathbf{C}_1(0)) > 0$.

For $\nu = 1,2,3,\ldots$ let $\varepsilon_{\nu} \rightarrow 0^+$ and $T^{(\nu)}$ be an $n$-dimensional locally area-minimizing rectifiable current of $\mathbf{C}_1(0)$ such that \eqref{relharm harmonic approx hyp} holds true with $\varepsilon_{\nu}$ and $T^{(\nu)}$ in place of $\varepsilon$ and $T$.  Let $\varphi^{(\nu)} \in C^2(B_1(0), \mathbb{R}^m)$ be a single-valued function such that $\|D\varphi^{(\nu)}\|_{C^1(B_1(0))} + \|\Delta\varphi^{(\nu)}\|_{C^0(B_1(0))} < \varepsilon_{\nu}^2$.  Let $\widehat{\mathcal E}_{\nu} = \widehat{\mathcal E}(T^{(\nu)},P_0,\mathbf{C}_1(0))$ and ${\mathcal E}_{\nu} = \mathcal{E}_{\varsigma/2,\varepsilon_{\nu}}(T^{(\nu)}, \varphi^{(\nu)}, \mathbf{C}_1(0)) > 0$.  Let $u^{(\nu)} : B_{(1+\gamma)/2}(0) \rightarrow \mathcal{A}_q(\mathbb{R}^m)$ be a Lipschitz approximation of $T^{(\nu)}$ and $K_{\nu} \subseteq B_{(1+\gamma)/2}(0)$ be a closed set as in Theorem~\ref{lip approx thm} (with $T^{(\nu)}, u^{(\nu)}, K_{\nu}, (1+\gamma)/2$ in place of $T, u, K, \gamma$) so that 
\begin{gather}\label{relharm harmonic approx eqn1}
	T^{(\nu)} \llcorner (K_{\nu} \times \mathbb{R}^m) = (\op{graph} u^{(\nu)}) \llcorner (K_{\nu} \times \mathbb{R}^m) , \quad 
	\op{Lip} u^{(\nu)} \leq \widehat{\mathcal E}_{\nu}^{\varsigma} , \\
	\mathcal{L}^n(B_{(1+\gamma)/2}(0) \setminus K_{\nu}) + \|T^{(\nu)}\|((B_{(1+\gamma)/2}(0) \setminus K_{\nu}) \times \mathbb{R}^m) 
		\leq C \widehat{\mathcal E}_{\nu}^{2+\varsigma}, \nonumber 
\end{gather}
where $C = C(n,m,q,\gamma) \in (0,\infty)$ is a constant.  It suffices to show that for infinitely many positive integers $\nu$ there exists a locally Dirichlet energy minimizing $q$-valued function $w$ (depending on $\nu$) such that \eqref{relharm harmonic approx concl} holds true with $T^{(\nu)},\varphi^{(\nu)},u^{(\nu)}$ in place of $T,\varphi,u$.

By Remark~\ref{relharm excess rmk1}(2), $\mathcal{E}_{\nu} \rightarrow 0$ as $\nu\rightarrow\infty$.  Let 
\begin{equation*}
	g^{(\nu)}(x) = \sum_{l=1}^q \left\llbracket \frac{u^{(\nu)}_l(x) - \varphi^{(\nu)}(x)}{{\mathcal E}_{\nu}} \right\rrbracket
\end{equation*}
for each $x \in B_{(1+\gamma)/2}(0)$, where $u^{(\nu)}(x) = \sum_{l=1}^q \llbracket u^{(\nu)}_l(x) \rrbracket$.  By \eqref{relharm energy eqn5} and \eqref{relharm harmonic approx eqn1}, 
\begin{equation}\label{relharm harmonic approx eqn2}
	\int_{B_{(1+\gamma)/2}(0)} |Dg^{(\nu)}|^2 
	= \frac{1}{{\mathcal E}_{\nu}^2} \int_{B_{(1+\gamma)/2}(0)} \mathcal{G}(Du^{(\nu)}, q \llbracket D\varphi^{(\nu)} \rrbracket)^2 \leq C 
\end{equation}
where $C = C(n,m,q,\gamma) \in (0,\infty)$ is a constant.  Hence by~\cite[Proposition~4.3]{DeLSpa1}, after passing to a subsequence there exists 
\begin{enumerate}[itemsep=2mm,topsep=0mm]
	\item[(i)]  positive integers $J$ and $q_1,q_2,\ldots,q_J$ with $\sum_{j=1}^J q_j = q$; 
	\item[(ii)]  $y^{(\nu)}_1,y^{(\nu)}_2,\ldots,y^{(\nu)}_J \in \mathbb{R}^m$ with $\lim_{\nu\rightarrow\infty} |y^{(\nu)}_i - y^{(\nu)}_j| = \infty$ for all $i \neq j$; 
	\item[(iii)]  $q_j$-valued functions $w_j \in W^{1,2}(B_{(1+\gamma)/2}(0),\mathcal{A}_{q_j}(\mathbb{R}^m))$ for each $j \in \{1,2,\ldots,J\}$; 
	\item[(iv)]  a sequence of translating sheets $(h^{(\nu)})$ in $W^{1,2}(B_{(1+\gamma)/2}(0),\mathcal{A}_q(\mathbb{R}^m))$ such that $h^{(\nu)}(x) = \sum_{j=1}^J \sum_{l=1}^{q_j} \llbracket y^{(\nu)}_j + w_{j,l}(x) \rrbracket$ for each $x \in B_{(1+\gamma)/2}(0)$, where $w_j(x) = \sum_{l=1}^{q_j} \llbracket w_{j,l}(x) \rrbracket$;
\end{enumerate}
such that 
\begin{gather}
	\label{relharm harmonic approx eqn3a} \lim_{\nu\rightarrow\infty} \int_{B_{(1+\gamma)/2}(0)} \mathcal{G}(g^{(\nu)}, h^{(\nu)})^2 = 0 , \\
	\label{relharm harmonic approx eqn3b} \sum_{j=1}^J \int_{B_{\gamma}(0)} |Dw_j|^2 \leq \liminf_{\nu\rightarrow\infty} \int_{B_{\gamma}(0)} |Dg^{(\nu)}|^2 , \\
	\label{relharm harmonic approx eqn3c} \limsup_{\nu\rightarrow\infty} \int_{B_{\gamma}(0)} (|Dg^{(\nu)}| - |Dh^{(\nu)}|)^2 
		\leq \limsup_{\nu\rightarrow\infty} \int_{B_{\gamma}(0)} |Dg^{(\nu)}|^2 - \sum_{j=1}^J \int_{B_{\gamma}(0)} |Dw_j|^2 .
\end{gather}
(Note that $|Dh^{(\nu)}|^2 = \sum_{j=1}^J |Dw_j|^2$ on $B_{(1+\gamma)/2}(0)$ for all $\nu$.)

We claim that $w_j$ and $h^{(\nu)}$ are locally Dirichlet energy minimizing in $B_{\gamma}(0)$ for all sufficiently large $\nu$.  If will suffice to show that each $w_j$ is locally Dirichlet energy minimizing $q_j$-valued function in $B_{\gamma}(0)$.  To see this, suppose that each $w_j$ is locally Dirichlet energy minimizing in $B_{\gamma}(0)$.  Note that each $w_j$ is continuous in $B_{\gamma}(0)$.  Let $\rho \in (0,\gamma)$ and let $\widehat{h}^{(\nu)} \in W^{1,2}(B_{\rho}(0),\mathcal{A}_q(\mathbb{R}^m))$ be a Dirichlet energy minimizing $q$-valued function such that $\widehat{h}^{(\nu)} = h^{(\nu)}$ on $\partial B_{\rho}(0)$.  By the maximum principle~\cite[Proposition~3.5]{DeLSpaDirMin}, $\widehat{h}^{(\nu)} = \sum_{j=1}^J \llbracket y^{(\nu)}_j + \widehat{w}_{j,l}(x) \rrbracket$ on $\overline{B_{\rho}(0)}$ where $\widehat{w}_j(x) = \sum_{l=1}^{q_j} \llbracket \widehat{w}_{j,l}(x) \rrbracket$ are Dirichlet energy minimizing $q_j$-valued functions on $B_{\rho}(0)$ such that $\widehat{w}_j = w_j$ on $\partial B_{\rho}(0)$.  Since each $w_j$ is locally Dirichlet energy minimizing in $B_{\gamma}(0)$, $\int_{B_{\rho}(0)} |Dw_j|^2 = \int_{B_{\rho}(0)} |D\widehat{w}_j|^2$ for all $j$ and thus $\int_{B_{\rho}(0)} |Dh^{(\nu)}|^2 = \int_{B_{\rho}(0)} |D\widehat{h}^{(\nu)}|^2$.  Hence $h^{(\nu)}$ must be Dirichlet energy minimizing in $B_{\rho}(0)$.  

Let $j_0 \in \{1,2,\ldots,J\}$.  To see that $w_{j_0}$ is locally Dirichlet energy minimizing in $B_{\gamma}(0)$, suppose to the contrary that there exists $\widetilde{w}_{j_0} \in W^{1,2}(B_{\gamma}(0);\mathcal{A}_{q_{j_0}}(\mathbb{R}^m))$ such that 
\begin{gather*}
	\widetilde{w}_{j_0} = w_{j_0} \text{ on } B_{(1+\gamma)/2}(0) \setminus B_{\gamma}(0), \nonumber \\ 
	\int_{B_{\gamma}(0)} |D\widetilde{w}_{j_0}|^2 \leq \int_{B_{\gamma}(0)} |Dw_{j_0}|^2 - \delta . %\label{relharm harmonic approx eqn4}
\end{gather*}
Set $\widetilde{w}_j = w_j$ for each $j \in \{1,\ldots,J\}$ with $j \neq j_0$.  Set $\overline{h}^{(\nu)}(x) = \sum_{j=1}^J \sum_{l=1}^{q_j} \llbracket y^{(\nu)}_j + \widetilde{w}_{j,l}(x) \rrbracket$ for each $x \in B_{(1+\gamma)/2}(0)$, where $\widetilde{w}_j(x) = \sum_{l=1}^{q_j} \llbracket \widetilde{w}_{j,l}(x) \rrbracket$, so that $\overline{h}^{(\nu)} = h^{(\nu)}$ on $B_{(1+\gamma)/2}(0) \setminus B_{\gamma}(0)$ and $\int_{B_{\gamma}(0)} |D\overline{h}^{(\nu)}|^2 \leq \int_{B_{\gamma}(0)} |Dh^{(\nu)}|^2 - \delta$.  Hence by \eqref{relharm harmonic approx eqn2}, \eqref{relharm harmonic approx eqn3a}, and \eqref{relharm harmonic approx eqn3b}, we can apply~\cite[Proposition~4.4]{DeLSpa1} to interpolate between $g^{(\nu)}$ and $\overline{h}^{(\nu)}$ and, after passing to a subsequence, find a Lipschitz $q$-valued function $\widetilde{h}^{(\nu)} : B_{(1+\gamma)/2}(0) \rightarrow \mathcal{A}_q(\mathbb{R}^m)$ such that 
\begin{gather} 
	\label{relharm harmonic approx eqn6a} \widetilde{h}^{(\nu)} = g^{(\nu)} \text{ on } B_{(1+\gamma)/2}(0) \setminus B_{(1+3\gamma)/4}(0), \\
	\label{relharm harmonic approx eqn6b} \op{Lip} \widetilde{h}^{(\nu)} \leq C^* \big( \op{Lip} g^{(\nu)} + 1 \big) , \\
	\label{relharm harmonic approx eqn6c} \int_{B_{(1+3\gamma)/4}(0)} |D\widetilde{h}^{(\nu)}|^2 
		\leq \int_{B_{(1+3\gamma)/4}(0)} |Dg^{(\nu)}|^2 - \frac{3\delta}{4} , 
\end{gather}
where $C^* \in (0,\infty)$ is a constant depending on $\delta$ and the sequences $(h^{(\nu)})$ and $(g^{(\nu)})$ but independent of $\nu$.  Let $g^{(\nu)}_a(x) = \frac{1}{q} \sum_{l=1}^q g^{(\nu)}_l(x)$ for each $x \in B_{(1+\gamma)/2}(0)$, where $g^{(\nu)}(x) = \sum_{l=1}^q \llbracket g^{(\nu)}_l(x) \rrbracket$.  Similarly, let $\widetilde{h}^{(\nu)}_a(x) = \frac{1}{q} \sum_{l=1}^q \widetilde{h}^{(\nu)}_l(x)$ for each $x \in B_{(1+\gamma)/2}(0)$, where $\widetilde{h}^{(\nu)}(x) = \sum_{l=1}^q \llbracket \widetilde{h}^{(\nu)}_l(x) \rrbracket$.  By the Poincar\'e inequality, \eqref{relharm harmonic approx eqn6a}, \eqref{relharm harmonic approx eqn2}, and \eqref{relharm harmonic approx eqn6c}, and  
\begin{align}\label{relharm harmonic approx eqn6d}
	\int_{B_{(1+3\gamma)/4}(0)} |g^{(\nu)}_a - \widetilde{h}^{(\nu)}_a|^2 
	\leq\,& C \int_{B_{(1+3\gamma)/4}(0)} |Dg^{(\nu)}_a - D\widetilde{h}^{(\nu)}_a|^2
	\\ \leq\,& 2C \int_{B_{(1+3\gamma)/4}(0)} (|Dg^{(\nu)}|^2 + |D\widetilde{h}^{(\nu)}|^2)
	\leq C, \nonumber 
\end{align}
where $C = C(n,m,q,\gamma) \in (0,\infty)$ are constants.  Set $\widetilde{u}^{(\nu)}(x) = \sum_{l=1}^{q} \llbracket \varphi^{(\nu)}(x) + {\mathcal E}_{\nu} \widetilde{h}^{(\nu)}_l(x) \rrbracket$ for each $x \in B_{(1+\gamma)/2}(0)$, where $\widetilde{h}^{(\nu)}(x) = \sum_{l=1}^q \llbracket \widetilde{h}^{(\nu)}_l(x) \rrbracket$.  By \eqref{relharm harmonic approx eqn6a} 
\begin{equation}\label{relharm harmonic approx eqn7a} 
	\widetilde{u}^{(\nu)} = u^{(\nu)} \text{ on } B_{(1+\gamma)/2}(0) \setminus B_{(1+3\gamma)/4}(0) .
\end{equation}
By noting that $\op{Lip} g^{(\nu)} \leq \mathcal{E}_{\nu}^{-1} (\sqrt{q} \op{Lip} \varphi^{(\nu)} + \op{Lip} u^{(\nu)})$ and using \eqref{relharm harmonic approx eqn6b} and $\op{Lip} u^{(\nu)} \leq \widehat{\mathcal E}_{\nu}^{\varsigma}$ (by \eqref{relharm harmonic approx eqn1}) 
\begin{equation}\label{relharm harmonic approx eqn7b} 
	\op{Lip} \widetilde{u}^{(\nu)} \leq C^* \|D\varphi^{(\nu)}\|_{C^0(B_1(0))} + C^* \widehat{\mathcal E}_{\nu}^{\varsigma} + C^* {\mathcal E}_{\nu} , 
\end{equation}
where $C^* \in (0,\infty)$ is a constant independent of $\nu$.  By \eqref{relharm harmonic approx eqn6a}, $\|\Delta\varphi\|_{C^0(B_1(0))} \leq \varepsilon_{\nu} \mathcal{E}_{\nu}$, and \eqref{relharm harmonic approx eqn6d} 
\begin{equation*}
	\int_{B_{(1+3\gamma)/4}(0)} D\varphi^{(\nu)} \cdot (D\widetilde{h}^{(\nu)}_a - Dg^{(\nu)}_a) 
	= -\int_{B_{(1+3\gamma)/4}(0)} \Delta\varphi^{(\nu)} \cdot (\widetilde{h}^{(\nu)}_a - g^{(\nu)}_a) 
	\leq C \varepsilon_{\nu} \mathcal{E}_{\nu}, 
\end{equation*}
where $C = C(n,\gamma) \in (0,\infty)$ is a constant.  Hence it follows from \eqref{relharm harmonic approx eqn6c} that 
\begin{align}\label{relharm harmonic approx eqn7c} 
	&\int_{B_{(1+3\gamma)/4}(0)} |D\widetilde{u}^{(\nu)}|^2 - \int_{B_{(1+3\gamma)/4}(0)} |Du^{(\nu)}|^2 
	\\=\,& 2q {\mathcal E}_{\nu} \int_{B_{(1+3\gamma)/4}(0)} D\varphi^{(\nu)} \cdot (D\widetilde{h}^{(\nu)}_a - Dg^{(\nu)}_a) 
		+ {\mathcal E}_{\nu}^2 \int_{B_{(1+3\gamma)/4}(0)} (|D\widetilde{h}^{(\nu)}|^2 - |Dg^{(\nu)}|^2) \nonumber 
	\\ \leq\,& -\frac{\delta}{2} {\mathcal E}_{\nu}^2 . \nonumber
\end{align}
By \eqref{relharm harmonic approx eqn7a}, \eqref{relharm harmonic approx eqn7c}, and \eqref{relharm harmonic approx eqn1} 
\begin{equation}\label{relharm harmonic approx eqn8} 
	\int_{B_{(1+\gamma)/2}(0)} |D\widetilde{u}^{(\nu)}|^2 \leq \int_{B_{(1+\gamma)/2}(0)} |Du^{(\nu)}|^2 \leq 4 \widehat{\mathcal E}_{\nu}^2 
\end{equation}
for all sufficiently large $\nu$.

Now we argue as in~\cite[Theorem~5.2]{DeLSpa1}.  In particular, using the coarea formula, Fatou's lemma, and \eqref{relharm harmonic approx eqn1}, there exists $\overline{\rho} \in \big(\frac{1+3\gamma}{4},\frac{1+\gamma}{2}\big)$ such that the slice $\langle T^{(\nu)} - \op{graph} u^{(\nu)}, r, \overline{\rho} \rangle$ exists and is an $(n-1)$-dimensional rectifiable current with mass $\leq C(n,m,q,\gamma) \,\widehat{\mathcal E}_{\nu}^{2+\varsigma}$, where $r(x,y) = |x|$.  By the isoperimetric inequality, there exists an $n$-dimensional area-minimizing rectifiable current $R^{(\nu)}$ of $\mathbb{R}^{n+m}$ such that 
\begin{equation}\label{relharm harmonic approx eqn9} 
	\partial R^{(\nu)} = \langle T^{(\nu)} - \op{graph} u^{(\nu)}, r, \overline{\rho} \rangle , \quad\quad 
	\mathbf{M}(R^{(\nu)}) \leq C(n,m,q,\gamma) \,\widehat{\mathcal E}_{\nu}^{\frac{n(2+\varsigma)}{n-1}} .
\end{equation}
Let $S^{(\nu)}$ be a competitor for $T^{(\nu)}$ defined by  
\begin{equation*}
	S^{(\nu)} = \op{graph} \widetilde{u}^{(\nu)} \llcorner \mathbf{C}_{\overline{\rho}}(0) + R^{(\nu)} 
		+ T^{(\nu)} \llcorner \mathbf{C}_1(0) \setminus \mathbf{C}_{\overline{\rho}}(0) . 
\end{equation*}
Using \eqref{relharm harmonic approx eqn7a} and \eqref{relharm harmonic approx eqn9}, it is readily verified that $\partial S^{(\nu)} = \partial T^{(\nu)}$ and $S^{(\nu)} \llcorner \mathbf{C}_1(0) \setminus \mathbf{C}_{(1+\gamma)/2}(0) = T^{(\nu)} \llcorner \mathbf{C}_1(0) \setminus \mathbf{C}_{(1+\gamma)/2}(0)$.  By \eqref{relharm harmonic approx eqn1}, \eqref{relharm harmonic approx eqn7a}, \eqref{relharm harmonic approx eqn7b}, \eqref{relharm harmonic approx eqn7c}, \eqref{relharm harmonic approx eqn8}, and \eqref{relharm harmonic approx eqn9} 
\begin{align*}
	\mathbf{M}(S^{(\nu)}) - \mathbf{M}(T^{(\nu)}) 
	\leq\,& \mathbf{M}(\op{graph} \widetilde{u}^{(\nu)} \llcorner \mathbf{C}_{\overline{\rho}}(0)) - \mathbf{M}(T^{(\nu)} \llcorner \mathbf{C}_{\overline{\rho}}(0)) 
		+ \mathbf{M}(R^{(\nu)}) 
	\\ \leq\,& \frac{1}{2} \int_{B_{\overline{\rho}}(0)} |D\widetilde{u}^{(\nu)}|^2 - \frac{1}{2} \int_{B_{\overline{\rho}}(0)} |Du^{(\nu)}|^2 
		\\&+ C \|D\varphi^{(\nu)}\|_{C^0(B_1(0))}^2 \widehat{\mathcal E}_{\nu}^2 + C {\mathcal E}_{\nu}^2 \widehat{\mathcal E}_{\nu}^2 
		+ C \widehat{\mathcal E}_{\nu}^{2+\varsigma}
	\\ \leq\,& -\frac{\delta}{8} {\mathcal E}_{\nu}^2 < 0
\end{align*}
for all sufficiently large $\nu$ and some constant $C \in (0,\infty)$ independent of $\nu$.  
This contradicts $T^{(\nu)}$ being locally area-minimizing in $\mathbf{C}_1(0)$.

Next we claim that 
\begin{equation}\label{relharm harmonic approx eqn10} 
	\lim_{\nu\rightarrow\infty} \int_{B_{\gamma}(0)} |Dg^{(\nu)}|^2 = \sum_{j=1}^J \int_{B_{\gamma}(0)} |Dw_j|^2 .
\end{equation}
In light of \eqref{relharm harmonic approx eqn3b}, it suffices to assume to the contrary that after passing to a subsequence there exists $\delta > 0$ such that 
\begin{equation}\label{relharm harmonic approx eqn11} 
	\sum_{j=1}^J \int_{B_{\gamma}(0)} |Dw_j|^2 \leq \int_{B_{\gamma}(0)} |Dg^{(\nu)}|^2 - \delta . 
\end{equation}
By \eqref{relharm harmonic approx eqn2}, \eqref{relharm harmonic approx eqn3a}, and \eqref{relharm harmonic approx eqn11}, we can apply~\cite[Proposition~4.4]{DeLSpa1} to interpolate between $g^{(\nu)}$ and $h^{(\nu)}$ and find a Lipschitz $q$-valued function $\widetilde{h}^{(\nu)} : B_{(1+\gamma)/2}(0) \rightarrow \mathcal{A}_q(\mathbb{R}^m)$ such that \eqref{relharm harmonic approx eqn6a}--\eqref{relharm harmonic approx eqn6c} hold true.  Repeating the above argument, we again derive a contradiction with $T^{(\nu)}$ being locally area-minimizing in $\mathbf{C}_1(0)$.  Therefore, \eqref{relharm harmonic approx eqn10} holds true.  By \eqref{relharm harmonic approx eqn3a}, \eqref{relharm harmonic approx eqn3c}, and \eqref{relharm harmonic approx eqn10}, \eqref{relharm harmonic approx concl} holds true with $T^{(\nu)},\varphi^{(\nu)},u^{(\nu)},\mathcal{E}_{\nu} h^{(\nu)}$ in place fo $T,\varphi,u,w$ for infinitely many $\nu$.
\end{proof}

\subsection{$L^{\infty}$-distance estimate in terms of tilt excess relative to a nearly harmonic graph}  
Let $T$ be an $n$-dimensional locally area-minimizing rectifiable current of $\mathbf{C}_1(0)$ and $\varphi \in C^2(B_1(0),\mathbb{R}^m)$.  In this subsection, we want to prove Theorem~\ref{relharm separation1 thm}, in which we bound the $L^{\infty}$-distance of $T$ to some union of graphs parallel $\op{graph} \varphi$ in terms of the tilt excess $\mathcal{E}_{\varsigma/2,\varepsilon}(T,\varphi,\mathbf{C}_1(0))$ of $T$ relative to the graph of $\varphi$.

\begin{lemma}[Controlled growth of tilt excess]\label{relharm tilt growth lemma}
For every positive integer $q$ and every $\theta \in (0,1/16)$ there exists $\varepsilon = \varepsilon(n,m,q,\theta) \in (0,1)$ such that the following holds true.  Let $x_0 \in \mathbb{R}^n$ and $\rho > 0$.  Let $T$ be an $n$-dimensional locally area-minimizing rectifiable current of $\mathbf{C}_{\rho}(x_0)$ such that 
\begin{equation*}
	(\partial T) \llcorner \mathbf{C}_{\rho}(x_0) = 0, \quad \sup_{X \in \op{spt} T} \op{dist}(X,P_0) < \infty, \quad
	\pi_{\#} T = q \llbracket B_{\rho}(x_0) \rrbracket, \quad \widehat{\mathcal E}(T,P_0,\mathbf{C}_{\rho}(x_0)) < \varepsilon .  
\end{equation*}
Let $\varphi \in C^2(B_{\rho}(x_0), \mathbb{R}^m)$ be a single-valued function with $\|D\varphi\|_{C^0(B_{\rho}(x_0))} + \rho \|\Delta\varphi\|_{C^0(B_{\rho}(x_0))} < \varepsilon^2$.  Then 
\begin{equation}\label{relharm tilt growth concl}
	\mathcal{E}_{\varsigma/2,\varepsilon}(T,\varphi,\mathbf{C}_{\theta\rho}(x_0)) 
	\leq C \theta^{\mu-1} \mathcal{E}_{\varsigma/2,\varepsilon}(T,\varphi,\mathbf{C}_{\rho}(x_0)) ,  
\end{equation}
where $\varsigma = \varsigma(n,m,q) \in (0,1)$ is a constant (which is small enough that Theorem~\ref{lip approx thm} holds true) and $\mu = \mu(n,m,q) \in (0,1)$ and $C = C(n,m,q) \in (0,\infty)$ are constants (independent of $\theta$).
\end{lemma}

\begin{proof}
Assume $x_0 = 0$ and $\rho = 1$.  Let $u : B_{5/8}(0) \rightarrow \mathcal{A}_q(\mathbb{R}^m)$ be a Lipschitz approximation of $T$ as in Theorem~\ref{lip approx thm} with $\gamma = 5/8$.  Let $w : B_{1/4}(0) \rightarrow \mathcal{A}_q(\mathbb{R}^m)$ be as in Lemma~\ref{relharm harmonic approx lemma} with $\gamma = 1/4$ and with $\eta = \eta(n,m,q,\theta) \in (0,1)$ to be later determined.  By~\cite[Theorem~2.13]{Almgren}, there exists a constant $\mu = \mu(n,m,q) \in (0,1)$ such that 
\begin{equation}\label{relharm tilt growth eqn1}
	\theta^{2-n} \int_{B_{\theta}(0)} |Dw|^2 \leq 4^{n-2+2\mu} \theta^{2\mu} \int_{B_{1/4}(0)} |Dw|^2 .
\end{equation}
Hence by the estimates for $u$ and $w$ from Lemma~\ref{relharm harmonic approx lemma}
\begin{align*}
	\frac{1}{\theta^n} \int_{B_{\theta}(0)} |D(u-\varphi)|^2 \leq 4^{n-1+2\mu} \theta^{2\mu-2} \int_{B_{1/4}(0)} |D(u-\varphi)|^2 + C_{\theta} \eta \mathcal{E}_{\varsigma/2,\varepsilon}(T,\varphi,\mathbf{C}_1(0)) , 
\end{align*}
where $\varsigma = \varsigma(n,m,q) \in (0,1)$ is as in Theorem~\ref{lip approx thm} and $C_{\theta} = C_{\theta}(n,m,q,\theta) \in (0,\infty)$ is a constant.  By the area formula, the estimates for $u$ from Theorem~\ref{lip approx thm}, and \eqref{relharm energy eqn5} 
\begin{align*}
	&\frac{1}{\theta^n} \int_{G_n(\mathbf{C}_{\theta}(0))} \sum_{\kappa=1}^m |\pi_S(-D\varphi^{\kappa}(x),e_{\kappa})|^2 \,d|T|(X,S)
	\\ \leq\,& 2 \cdot 4^{n-1+2\mu} \theta^{2\mu-2} \int_{G_n(\mathbf{C}_1(0))} \sum_{\kappa=1}^m |\pi_S(-D\varphi^{\kappa}(x),e_{\kappa})|^2 \,d|T|(X,S)  
		\\&+ C_{\theta} (\eta \,\mathcal{E}_{\varsigma/2,\varepsilon}(T,\varphi,\mathbf{C}_1(0)) 
		+ \widehat{\mathcal E}(T,P_0,B_1(0))^{2+\varsigma}) ,  
\end{align*}
where $C_{\theta} = C_{\theta}(n,m,q,\theta) \in (0,\infty)$ is a constant.
Provided $\eta = \eta(n,m,q,\theta)$ and $\varepsilon = \varepsilon(n,m,q,\theta)$ are sufficiently small 
\begin{equation}\label{relharm tilt growth eqn2}
	\frac{1}{\omega_n \theta^n} \int_{G_n(\mathbf{C}_{\theta}(0))} \sum_{\kappa=1}^m |\pi_S(-D\varphi^{\kappa}(x),e_{\kappa})|^2 \,d|T|(X,S)
	\leq C \theta^{\mu-1} \mathcal{E}_{\varsigma/2,\varepsilon}(T,\varphi,\mathbf{C}_1(0)) , 
\end{equation}
where $C = C(n,m,q) \in (0,\infty)$ is a constant (independent of $\theta$).  By~\cite[Lemma~3.15]{KrumWicb}, 
\begin{equation}\label{relharm tilt growth eqn3}
	\widehat{\mathcal E}(T,P_0,\mathbf{C}_{\theta}(0)) \leq C \theta^{\mu-1} \widehat{\mathcal E}(T,P_0,\mathbf{C}_1(0)) , 
\end{equation}
where $\mu = \mu(n,m,q) \in (0,1)$ is as in \eqref{relharm tilt growth eqn1} and $C = C(n,m,q) \in (0,\infty)$ is a constant.  After replacing $\varsigma = \varsigma(n,m,q)$ with a constant small enough that $\varsigma < \frac{2\mu}{1-\mu}$, and thus $(2+\varsigma/2)(1-\mu) < 2(1-\mu/2)$, 
\begin{align}\label{relharm tilt growth eqn4}
	\widehat{\mathcal E}(T,P_0,\mathbf{C}_{\theta}(0))^{2+\varsigma/2} 
	\leq\,& C \theta^{(2+\varsigma/2)(\mu-1)} \widehat{\mathcal E}(T,P_0,\mathbf{C}_1(0))^{2+\varsigma/2} 
	\\ \leq\,& C \theta^{2(\mu/2-1)} \widehat{\mathcal E}(T,P_0,\mathbf{C}_1(0))^{2+\varsigma/2} . \nonumber 
\end{align}
Combining \eqref{relharm tilt growth eqn2}, \eqref{relharm tilt growth eqn3}, and \eqref{relharm tilt growth eqn4} gives us \eqref{relharm tilt growth concl} with $\mu/2$ in place of $\mu$.
\end{proof}

\begin{lemma}[Poincar\'e-type inequality]\label{relharm poincare lemma}
Let $q$ be a positive integer and $\gamma \in (0,1)$.  There exists $\varepsilon = \varepsilon(n,m,q,\gamma) \in (0,1)$ such that the following holds true.  Let $x_0 \in \mathbb{R}^n$ and $\rho > 0$.  Let $T$ be an $n$-dimensional locally area-minimizing rectifiable current of $\mathbf{C}_{\rho}(x_0)$ such that 
\begin{equation*}
	(\partial T) \llcorner \mathbf{C}_{\rho}(x_0) = 0, \quad \sup_{X \in \op{spt} T} \op{dist}(X,P_0) < \infty, \quad
	\pi_{\#} T = q \llbracket B_{\rho}(x_0) \rrbracket, \quad \widehat{\mathcal E}(T,P_0,\mathbf{C}_{\rho}(x_0)) < \varepsilon .  
\end{equation*}
Let $u : B_{(1+\gamma)\rho/2}(x_0) \rightarrow \mathcal{A}_q(\mathbb{R}^m)$ be a Lipschitz approximation of $T$ as in Theorem~\ref{lip approx thm} with $(1+\gamma)/2$ in place of $\gamma$.  Let $\varphi \in C^2(B_{\rho}(x_0), \mathbb{R}^m)$ be a single-valued function.  Then there exists $a \in \mathcal{A}_q(\mathbb{R}^m)$ such that 
\begin{align}\label{relharm poincare concl}
	&\frac{1}{\rho^{n+2}} \int_{\mathbf{C}_{\gamma\rho}(x_0)} \op{dist}^2(X, \op{spt} {\rm graph}(\varphi+a)) \,d\|T\|(X) 
	\\&\hspace{15mm} + \frac{1}{\rho^{n+2}} \int_{B_{\gamma\rho}(x_0)} \mathcal{G}(u,\varphi+a)^2 
	\leq C \mathcal{E}_{\varsigma}(T,\varphi,\mathbf{C}_{\rho}(x_0))^2 , \nonumber 
\end{align}
where $\varphi+a$ is as in \eqref{varphi plus a defn}, $\varsigma = \varsigma(n,m,q) \in (0,1)$ is as in Theorem~\ref{lip approx thm}, and $C = C(n,m,q,\gamma) \in (0,\infty)$ is a constant.
\end{lemma}

\begin{proof}
Assume $x_0 = 0$ and $\rho = 1$.  Let $u$ and $K$ be as in Theorem~\ref{lip approx thm} with $(1+\gamma)/2$ in place of $\gamma$.  We apply the Poincar\'e inequality~\cite[Proposition~2.12]{DeLSpaDirMin} to $u - \varphi$ to find $a \in \mathcal{A}_q(\mathbb{R}^m)$ such that 
\begin{equation}\label{relharm poincare eqn1}
	\int_{B_{(1+\gamma)/2}(0)} \mathcal{G}(u,\varphi+a)^2 \leq C \int_{B_{(1+\gamma)/2}(0)} \mathcal{G}(Du,q \llbracket D\varphi \rrbracket)^2 , 
\end{equation}
where $C = C(n,m,q,\gamma) \in (0,\infty)$ is a constant.  It follows from the area formula, Theorem~\ref{lip approx thm}, \eqref{relharm energy eqn5}, and \eqref{relharm poincare eqn1} that 
\begin{gather}
	\label{relharm poincare eqn2} \int_{B_{(1+\gamma)/2}(0)} \mathcal{G}(u,\varphi+a)^2 \leq C \mathcal{E}_{\varsigma}(T,\varphi,\mathbf{C}_{\rho}(x_0))^2 , \\
	\label{relharm poincare eqn3} \int_{K \times \mathbb{R}^m} \op{dist}^2(X, \op{spt} \op{graph}(\varphi+a)) \,d\|T\|(X) 
		\leq C \mathcal{E}_{\varsigma}(T,\varphi,\mathbf{C}_{\rho}(x_0))^2 
\end{gather}
for some constant $C = C(n,m,q,\gamma) \in (0,\infty)$.  Notice that for each $Z \in \op{spt} T \cap \mathbf{C}_{\gamma}(0)$
\begin{equation}\label{relharm poincare eqn4}
    \op{dist}(Z, \op{spt} \op{graph}(\varphi+a)) \leq 1-\gamma .
\end{equation}
since otherwise by \eqref{relharm poincare eqn3} and $\|T\|((B_{(1+\gamma)/2}(0) \setminus K) \times \mathbb{R}^m) \leq C \widehat{\mathcal E}(T,P_0,\mathbf{C}_1(0))^{2+\varsigma}$ (by Theorem~\ref{lip approx thm}) 
\begin{equation*}
    \omega_n \bigg(\frac{1-\gamma}{2}\bigg)^{n+2} \leq \int_{\mathbf{B}_{(1-\gamma)/2}(Z)} \op{dist}^2(X,\op{spt}\op{graph}(\varphi+a)) \,d\|T\|(X) \leq C \mathcal{E}_{\varsigma}(T,\varphi,\mathbf{C}_{\rho}(x_0))^2 
\end{equation*}
for some constant $C = C(n,m,q,\gamma) \in (0,\infty)$, which provided $\varepsilon = \varepsilon(n,m,q,\gamma)$ is sufficiently small gives us a contradiction. Thus again using $\|T\|((B_{(1+\gamma)/2}(0) \setminus K) \times \mathbb{R}^m) \leq C \widehat{\mathcal E}(T,P_0,\mathbf{C}_1(0))^{2+\varsigma}$
\begin{equation}\label{relharm poincare eqn5}
	 \int_{(B_{\gamma}(0) \setminus K) \times \mathbb{R}^m} \op{dist}^2(X, \op{spt} \op{graph}(\varphi+a)) \,d\|T\|(X) 
	 	\leq C \widehat{\mathcal E}(T,P_0,\mathbf{C}_1(0))^{2+\varsigma}  
\end{equation}
for some constant $C = C(n,m,q,\gamma) \in (0,\infty)$.  Combining \eqref{relharm poincare eqn2}, \eqref{relharm poincare eqn3}, and \eqref{relharm poincare eqn5} gives us \eqref{relharm poincare concl}.
\end{proof}

\begin{theorem}\label{relharm separation1 thm} 
Let $q$ be a positive integer and $\gamma \in (0,1)$.  There exists $\varepsilon = \varepsilon(n,m,q,\gamma) \in (0,1)$ such that the following holds true.  Let $T$ be an $n$-dimensional locally area-minimizing rectifiable current of $\mathbf{C}_1(0)$ such that 
\begin{equation*}
	(\partial T) \llcorner \mathbf{C}_1(0) = 0, \quad \sup_{X \in \op{spt} T} \op{dist}(X,P_0) < \infty, \quad \pi_{\#} T = q \llbracket B_1(0) \rrbracket, \quad 
	\widehat{\mathcal E}(T,P_0,\mathbf{C}_1(0)) < \varepsilon .  
\end{equation*}
Let $\varphi \in C^2(B_1(0), \mathbb{R}^m)$ be a single-valued function with $\|D\varphi\|_{C^0(B_1(0))} + \|\Delta\varphi\|_{C^0(B_1(0))} < \varepsilon^2$.  Then there exists $a \in \mathcal{A}_q(\mathbb{R}^m)$ such that 
\begin{equation}\label{relharm separation1 concl}
	\sup_{X \in \op{spt} T \cap \mathbf{C}_{\gamma}(0)} \op{dist}(X, \op{spt} {\rm graph}(\varphi+a)) 
	\leq C \mathcal{E}_{\varsigma/2}(T,\varphi,\mathbf{C}_1(0)) , 
\end{equation}
where $\varphi+a$ is as in \eqref{varphi plus a defn}, $\varsigma = \varsigma(n,m,q) \in (0,1)$ is as in Lemma~\ref{relharm tilt growth lemma}, and $C = C(n,m,q,\gamma) \in (0,\infty)$ is a constant.
\end{theorem}

\begin{proof}
In light of Remark~\ref{relharm excess rmk1}(1), we may assume that $\mathcal{E}_{\varsigma/2,\varepsilon}(T,\varphi,\mathbf{C}_1(0)) > 0$ as otherwise \eqref{relharm separation1 concl} obviously holds true.  Choose $\theta = \theta(n,m,q) \in (0,1/16)$ so that $C\theta^{\mu/2} < 1$, where $\mu$ and $C$ are as in Lemma~\ref{relharm tilt growth lemma}.  Let $u : B_{(3+\gamma)/4}(0) \rightarrow \mathcal{A}_q(\mathbb{R}^m)$ and $K \subset B_{(3+\gamma)/4}(0)$ be as in Theorem~\ref{lip approx thm} with $x_0 = 0$, $\rho = 1$, and $(3+\gamma)/4$ in place of $\gamma$.  By Lemma~\ref{relharm poincare lemma} (with $(1+\gamma)/2$ in place of $\gamma$) there exists $a \in \mathcal{A}_q(\mathbb{R}^m)$ such that 
\begin{align}\label{relharm separation1 eqn1}
	&\int_{\mathbf{C}_{(1+\gamma)/2}(0)} \op{dist}^2(X, \op{spt} {\rm graph}(\varphi+a)) \,d\|T\|(X) 
	\\&\hspace{15mm} + \int_{B_{(1+\gamma)/2}(0)} \mathcal{G}(u,\varphi+a)^2 
	\leq C_1 \mathcal{E}_{\varsigma/2,\varepsilon}(T,\varphi,\mathbf{C}_1(0))^2 , \nonumber 
\end{align}
where $C_1 = C_1(n,m,q,\gamma) \in (0,\infty)$ is a constant.  

Let $\varepsilon_0 = \varepsilon_0(n,m,q,\gamma) \in (0,1)$ be a suitably small constant.  By Remark~\ref{relharm excess rmk1}(2), provided $\varepsilon$ is sufficiently small we have that $\mathcal{E}_{\varsigma/2,\varepsilon}(T, \varphi, \mathbf{C}_1(0)) < \varepsilon_0$.  Let $i_0 \geq 1$ be the smallest positive integer such that $\theta^{i_0 (\mu/2-1)} \mathcal{E}_{\varsigma/2,\varepsilon}(T, \varphi, \mathbf{C}_1(0)) \geq \varepsilon_0$.  Fix any $\xi \in B_{\gamma}(0)$.  By inductively applying Lemma~\ref{relharm tilt growth lemma}
\begin{equation}\label{relharm separation1 eqn2}
	\mathcal{E}_{\varsigma/2,\varepsilon}(T, \varphi, \mathbf{C}_{\theta^i (1-\gamma)/2}(\xi)) 
	\leq \bigg(\frac{2}{1-\gamma}\bigg)^{n(1+\varsigma/2)} \theta^{i(\mu/2-1)} \mathcal{E}_{\varsigma/2,\varepsilon}(T, \varphi, \mathbf{C}_1(0)) 
\end{equation}
for all $i \in \{0,1,2,\ldots,i_0\}$.  In particular, \eqref{relharm separation1 eqn2} with $i = 0$ is readily verified using Definition~\ref{relharm excess defn1}.  If \eqref{relharm separation1 eqn2} holds true for some $i \in \{0,1,\ldots,i_0-1\}$, then $\mathcal{E}_{\varsigma/2,\varepsilon}(T,\varphi,\mathbf{C}_{\theta^i (1-\gamma)/2}(\xi)) \leq \big(\frac{2}{1-\gamma}\big)^{n(1+\varsigma/2)} \varepsilon_0$.  Thus by Remark~\ref{relharm excess rmk1}(3), provided $\varepsilon_0$ is sufficiently small, we can apply Lemma~\ref{relharm tilt growth lemma} to obtain \eqref{relharm separation1 eqn2} with $i+1$ in place of $i$.  

For each $i \in \{1,2,\ldots,i_0\}$, let $u_i : B_{(3+\gamma)\theta^i/8}(\xi) \rightarrow \mathcal{A}_q(\mathbb{R}^m)$ and $K \subset B_{(3+\gamma)\theta^i/8}(\xi)$ be as in Theorem~\ref{lip approx thm} with $\xi$, $\theta^i(1-\gamma)/2$, $3/4$, $u_i$, and $K_i$ in place of $x_0$, $\rho$, $\gamma$, $u$, and $K$.  By Remark~\ref{relharm excess rmk1}(3) we can use Lemma~\ref{relharm poincare lemma} (with $\gamma = 1/2$) and \eqref{relharm separation1 eqn2} to find $a_i \in \mathcal{A}_q(\mathbb{R}^m)$ such that 
\begin{align}\label{relharm separation1 eqn3}
	&\frac{1}{\theta^{i(n+2)}} \int_{\mathbf{C}_{\theta^i(1-\gamma)/4}(\xi)} \op{dist}^2(X, \op{spt} {\rm graph}(\varphi+a_i)) \,d\|T\|(X) 
	\\&\hspace{15mm} + \frac{1}{\theta^{i(n+2)}} \int_{B_{\theta^i(1-\gamma)/4}(\xi)} \mathcal{G}(u_i,\varphi+a_i)^2 
	\leq C \theta^{2i (\mu/2-1)} \mathcal{E}_{\varsigma/2,\varepsilon}(T,\varphi,\mathbf{C}_1(0))^2 , \nonumber 
\end{align}
where $C = C(n,m,q,\gamma) \in (0,\infty)$ is a constant.  Let $u_0 = u$, $K_0 = K$, and $a_0 = a$ be as chosen above so that by \eqref{relharm separation1 eqn1}, \eqref{relharm separation1 eqn3} holds true when $i = 0$.  Noting that $u_i = u_{i+1}$ on $B_{(3+\gamma)\theta^{i+1}/8}(\xi) \cap K_i \cap K_{i+1}$ and $\mathcal{L}^n(B_{\theta^{i+1}(1+\gamma)/4}(\xi) \setminus (K_i \cup K_{i+1})) \leq C(n,m,q,\gamma) \, \theta^{in} \widehat{\mathcal E}(T,P_0,\mathbf{C}_{\theta^{i+1}(1+\gamma)/2}(\xi))^{2+\varsigma}$, by the triangle inequality and \eqref{relharm separation1 eqn3} 
\begin{equation*}
	\mathcal{G}(a_i,a_{i+1}) \leq C \theta^{i\mu/2} \mathcal{E}_{\varsigma/2,\varepsilon}(T,\varphi,\mathbf{C}_1(0)) 
\end{equation*}
for all $i \in \{0,1,2,\ldots,i_0-1\}$, where $C = C(n,m,q,\gamma) \in (0,\infty)$ is a constant.  Hence by the triangle inequality 
\begin{equation}\label{relharm separation1 eqn4}
	\mathcal{G}(a_i,a) \leq \sum_{j=0}^{i-1} \mathcal{G}(a_j,a_{j+1}) \leq C \mathcal{E}_{\varsigma/2,\varepsilon}(T,\varphi,\mathbf{C}_1(0)) 
\end{equation}
for all $i \in \{0,1,2,\ldots,i_0\}$, where $C = C(n,m,q,\gamma) \in (0,\infty)$ is a constant.  By \eqref{relharm separation1 eqn3} and~\cite[Lemma~3.6]{KrumWica} with $K = \op{spt} {\rm graph}(\varphi+a_i)$ 
\begin{equation}\label{relharm separation1 eqn5}
	\sup_{X \in \op{spt} T \cap \mathbf{C}_{\theta^i(1-\gamma)/8}(\xi)} \op{dist}(X, \op{spt}\op{graph}(\varphi+a_i))
	\leq C \theta^{i \frac{n+\mu}{n+2}} \mathcal{E}_{\varsigma/2,\varepsilon}(T,\varphi,\mathbf{C}_1(0))^{\frac{2}{n+2}}
\end{equation}
for all $i \in \{0,1,2,\ldots,i_0\}$, where $C = C(n,m,q,\gamma) \in (0,\infty)$ is a constant.  

Now notice that since $\theta^{i_0 (\mu/2-1)} \mathcal{E}_{\varsigma/2}(T, \varphi, \mathbf{C}_1(0)) \geq \varepsilon_0$, 
\begin{equation*}
	\theta^{i_0} \leq \left(\frac{\mathcal{E}_{\varsigma/2}(T, \varphi, \mathbf{C}_1(0))}{\varepsilon_0}\right)^{\frac{1}{1-\mu/2}} . 
\end{equation*}
Hence by \eqref{relharm separation1 eqn5} with $i = i_0$ 
\begin{align}\label{relharm separation1 eqn6}
	\sup_{X \in \op{spt} T \cap \mathbf{C}_{\theta^{i_0} (1-\gamma)/8}(\xi)} \op{dist}(X, \op{spt}\op{graph}(\varphi+a_{i_0}))
	\leq\,& C \mathcal{E}_{\varsigma/2,\varepsilon}(T,\varphi,\mathbf{C}_1(0))^{\frac{1}{1-\mu/2} \cdot \frac{n+\mu}{n+2} + \frac{2}{n+2}} 
	\\=\,& C \mathcal{E}_{\varsigma/2,\varepsilon}(T,\varphi,\mathbf{C}_1(0))^{\frac{1}{1-\mu/2}} \nonumber 
\end{align}
where $C = C(n,m,q,\gamma) \in (0,\infty)$ is a constant.  Combining \eqref{relharm separation1 eqn4} with $i = i_0$ and \eqref{relharm separation1 eqn6} 
\begin{equation*}
	\sup_{X \in \op{spt} T \cap \mathbf{C}_{\theta^{i_0} (1-\gamma)/8}(\xi)} \op{dist}(X, \op{spt}\op{graph}(\varphi+a))
	\leq C \mathcal{E}_{\varsigma/2,\varepsilon}(T,\varphi,\mathbf{C}_1(0))
\end{equation*}
for some constant $C = C(n,m,q,\gamma) \in (0,\infty)$.  Since $\xi \in \mathbf{C}_{\gamma}(0)$ is arbitrary, this proves that 
\begin{equation}\label{relharm separation1 eqn7}
	\sup_{X \in \op{spt} T \cap \mathbf{C}_{\gamma}(0)} \op{dist}(X, \op{spt} {\rm graph}(\varphi+a)) 
	\leq C_0 \mathcal{E}_{\varsigma/2,\varepsilon}(T,\varphi,\mathbf{C}_1(0)) 
\end{equation}
for some constant $C_0 = C_0(n,m,q,\gamma) \in (0,\infty)$.  Since $\varepsilon$ depends on $n,m,q,\gamma$, \eqref{relharm separation1 eqn7} gives us \eqref{relharm separation1 concl} with $C = C_0/\varepsilon$.
\end{proof}

\begin{remark}\label{relharm separation1 rmk}{\rm 
Assume that $\mathcal{E}_{\varsigma/2,\varepsilon}(T,\varphi,\mathbf{C}_1(0)) > 0$.  Recall from the proof of Theorem~\ref{relharm separation1 rmk} that we chose $a \in \mathcal{A}_q(\mathbb{R}^m)$ so that \eqref{relharm separation1 eqn1} holds true, and we proved that \eqref{relharm separation1 eqn7} holds true.  Let  
\begin{equation*}
	\{ X \in \mathbf{C}_{\gamma}(0) : \op{dist}(X, \op{spt}\op{graph}(\varphi+a)) < 2C_0 \mathcal{E}_{\varsigma/2,\varepsilon}(T,\varphi,C_1(0)) \} 
	= \bigcup_{i=1}^N U_i
\end{equation*}
where $\{U_i\}$ is a collection of mutually disjoint, connected, open subsets of $\mathbf{C}_{\gamma}(0)$ and $C_0$ is as in \eqref{relharm separation1 eqn7}.  Then 
\begin{equation*}
	T \llcorner \mathbf{C}_{\gamma}(0) = \sum_{i=1}^N T_i  \quad\text{where}\quad T_i = T \llcorner U_i 
\end{equation*}
and $T_i$ are locally area-minimizing rectifiable currents of $\mathbf{C}_{\gamma}(0)$ with $(\partial T_i) \llcorner \mathbf{C}_{\gamma}(0) = 0$.  By the constancy theorem 
\begin{equation*}
	\pi_{\#} T_i = q_i \llbracket B_{\gamma}(0) \rrbracket 
\end{equation*}
for some integers $q_i$ with $\sum_{i=1}^N q_i = q$.  Assuming $\varepsilon_0$ is sufficiently small, by~\cite[Lemma~3.5]{KrumWicb} we have that $q_i \geq 0$ for all $i \in \{1,2,\ldots,N\}$.  By~\cite[Lemma~3.4]{KrumWicb}, $q_i = 0$ if and only if $\op{spt} T_i \cap \mathbf{C}_{\gamma/2}(0) = \emptyset$.  We claim that if \eqref{relharm separation1 eqn1} and \eqref{relharm separation1 eqn7} hold true for $C_0 = C_0(n,m,q,\gamma)$ sufficiently large, then $q_i > 0$ for all $i \in \{1,2,\ldots,N\}$.  Observe that by \eqref{relharm separation1 eqn1} and Theorem~\ref{lip approx thm}, there exists a set $\Omega \subset K \cap B_{\gamma/2}(0)$ such that $\mathcal{L}^n(\Omega) \geq \frac{1}{2} \omega_n (\gamma/2)^n$ and $\mathcal{G}(u(x),\varphi(x)+a) \leq C_2 \mathcal{E}_{\varsigma/2,\varepsilon}(T,\varphi,\mathbf{C}_1(0))$ for all $x \in \Omega$, where $C_2 = \sqrt{\frac{4C_1}{\omega_n(\gamma/2)^n}}$ (for $C_1$ as in \eqref{relharm separation1 eqn1}).  Take any $i \in \{1,2,\ldots,N\}$ and $x \in \Omega$.  Let $a = \sum_{j=1}^q \llbracket a_j \rrbracket$ and $u(x) = \sum_{j=1}^q \llbracket u_j(x) \rrbracket$, where $a_j, u_j(x) \in \mathbb{R}^m$.  There exists $k(i) \in \{1,2,\ldots,q\}$ such that $\op{graph}(\varphi+a_{k(i)}) \cap \mathbf{C}_{\gamma}(0) \subset U_i$, and there exists $j(i) \in \{1,2,\ldots,q\}$ such that $X = (x,u_{j(i)}(x)) \in \op{spt} T$ and   
\begin{align*}
	\op{dist}(X,\op{graph}(\varphi+a_{k(i)})) \leq\,& |u_{j(i)}(x) - (\varphi(x) + a_{k(i)})| \leq \mathcal{G}(u(x),\varphi(x)+a) 
		\\ \leq\,& C_2 \mathcal{E}_{\varsigma/2,\varepsilon}(T,\varphi,\mathbf{C}_1(0)) .
\end{align*}
Hence provided we take $C_0$ in \eqref{relharm separation1 eqn7} large enough that $C_2 < 2C_0$, $X \in \op{spt} T_i \cap \mathbf{C}_{\gamma/2}(0)$.  Therefore, $q_i > 0$ for all $i \in \{1,2,\ldots,N\}$.
}\end{remark}

\subsection{$L^{\infty}$-distance estimate in terms $L^2$-distance relative to parallel nearly harmonic graphs}  We define the following excess $E_{\varsigma,\varepsilon}(T,\varphi,\mathbf{C}_{\rho}(x_0))$ of $T$ relative to the graph of $\varphi(x) = \sum_{l=1}^s \llbracket \varphi_l(x) \rrbracket$, where $\varphi_l \in C^2(B_{\rho}(x_0), \mathbb{R}^m)$ are single-valued functions.

\begin{definition}\label{relharm excess defn2}
Let $\varsigma,\varepsilon \in (0,1]$, $x_0 \in \mathbb{R}^n$, and $\rho > 0$.  Let $T$ be an $n$-dimensional rectifiable current $T$ of $\mathbf{C}_{\rho}(x_0)$.  Let $\varphi : B_{\rho}(x_0) \rightarrow \mathcal{A}_s(\mathbb{R}^m)$ such that $\varphi(x) = \sum_{l=1}^s \llbracket \varphi_l(x) \rrbracket$ where $\varphi_l \in C^2(B_{\rho}(x_0), \mathbb{R}^m)$ are single-valued functions.  
(Here we include the case $s = 1$ where $\varphi : B_{\rho}(x_0) \rightarrow \mathbb{R}^m$ is a single-valued function.)  We define 
\begin{align}\label{relharm excess defn2 eqn}
	E_{\varsigma,\varepsilon}(T,\varphi,\mathbf{C}_{\rho}(x_0))^2 
		=\,& \frac{1}{\omega_n \rho^{n+2}} \int_{\mathbf{C}_{\rho}(x_0)} \op{dist}^2(X, \op{spt} \op{graph} \varphi) \,d\|T\|(X) 
		\\&+ \sum_{l=1}^s \|D\varphi_l\|_{C^0(B_{\rho}(x_0))} \,\widehat{\mathcal E}(T,P_0,\mathbf{C}_{\rho}(x_0))^2  \nonumber
		\\&+ \frac{1}{\varepsilon^2} \rho^2 \sum_{l=1}^s \|\Delta\varphi_l\|_{C^0(B_{\rho}(x_0))}^2  
		+ \widehat{\mathcal E}(T,P_0,\mathbf{C}_{\rho}(x_0))^{2+\varsigma} , \nonumber 
\end{align}
where $\widehat{\mathcal E}(T,P_0,\mathbf{C}_{\rho}(x_0))$ is as in \eqref{oriented tilt defn}.  When $\varepsilon = 1$, we let $E_{\varsigma}(T,\varphi,\mathbf{C}_{\rho}(x_0)) = E_{\varsigma,1}(T,\varphi,\mathbf{C}_{\rho}(x_0))$.
\end{definition}

The first term on the right-hand side of \eqref{relharm excess defn2 eqn} quantifies the $L^2$-distance of $T$ relative to the graph of $\varphi$.  As was the case in Definition~\ref{relharm excess defn1}, the remaining terms on the right-hand side of \eqref{relharm excess defn2 eqn} control $\Delta\varphi$ and error terms arising from Almgren's strong Lipschitz approximation theorem and Lemma~\ref{relharm harmonic approx lemma}.

\begin{remark}\label{relharm excess rmk2}{\rm
(1)  $E_{\varsigma}(T,\varphi,\mathbf{C}_{\rho}(x_0)) = 0$ if and only if $\varphi$ is harmonic in $B_{\rho}(x_0)$ and $T$ is a sum of $n$-dimensional planes which are parallel to $P_0$, have orientation $\vec P_0$, and satisfy $\op{spt} T \subseteq \op{spt} \op{graph} \varphi$.  (Note that $E_{\varsigma}(T,\varphi,\mathbf{C}_{\rho}(x_0)) = 0$ does not guarantee that $\op{spt} T = \op{spt} \op{graph} \varphi$ and thus $\varphi$ need not be constant.)

\noindent (2)  Given $\varepsilon > 0$, if $E_{\varsigma}(T,\varphi,\mathbf{C}_{\rho}(x_0)) < \varepsilon$, then $\widehat{\mathcal E}(T,P_0,\mathbf{C}_{\rho}(x_0)) < \varepsilon^{\frac{2}{2+\varsigma}}$.
}\end{remark}

Now we prove the following estimate on the $L^{\infty}$-distance of a locally area-minimizing rectifiable current $T$ to the graph of exactly one single-valued $C^2$-function $\varphi$.

\begin{lemma}\label{relharm separation0 lemma}
Let $q$ be a positive integer and $\gamma \in (0,1)$.  There exists $\varepsilon = \varepsilon(n,m,q,\gamma) \in (0,1)$ such that the following holds true.  Let $T$ be an $n$-dimensional locally area-minimizing rectifiable current of $\mathbf{C}_1(0)$ such that 
\begin{equation*}%\label{relharm separation1 hyp}
	(\partial T) \llcorner \mathbf{C}_1(0) = 0, \quad \sup_{X \in \op{spt} T} \op{dist}(X,P_0) < \infty, \quad \pi_{\#} T = q \llbracket B_1(0) \rrbracket, \quad 
	\widehat{\mathcal E}(T,P_0,\mathbf{C}_1(0)) < \varepsilon . 
\end{equation*}
Let $\varphi \in C^2(B_1(0),\mathbb{R}^m)$ be a single-valued function such that $\|D\varphi\|_{C^0(B_1(0))} + \|\Delta\varphi\|_{C^0(B_1(0))} < \varepsilon$.  Then 
\begin{equation}\label{relharm separation0 concl}
	\sup_{X \in \op{spt} T \cap \mathbf{C}_{\gamma}(0)} \op{dist}(X, \op{graph} \varphi) \leq C E_{\varsigma/2}(T,\varphi,\mathbf{C}_1(0)) , 
\end{equation}
where $\varsigma = \varsigma(n,m,q) \in (0,1)$ is as in Lemma~\ref{relharm tilt growth lemma} and $C = C(n,m,q,\gamma) \in (0,\infty)$ is a constant.
\end{lemma}

\begin{proof}
In light of Remark~\ref{relharm excess rmk2}(1), we may assume that $E_{\varsigma/2}(T,\varphi,\mathbf{C}_1(0)) > 0$ as otherwise \eqref{relharm separation0 concl} trivially holds true.  By Lemma~\ref{relharm energy lemma} and Theorem~\ref{relharm separation1 thm}, there exists $a = \sum_{l=1}^q \llbracket a_l \rrbracket \in \mathcal{A}_q(\mathbb{R}^m)$ such that 
\begin{equation}\label{relharm separation0 eqn1}
	\sup_{(x,y) \in \op{spt} T \cap \mathbf{C}_{\gamma}(0)} \op{dist}(X, \op{spt} \op{graph}(\varphi+a)) < C_0 E_{\varsigma/2}(T,\varphi,\mathbf{C}_1(0)) , 
\end{equation}
where $C_0 = C_0(n,m,q,\gamma) \in (0,\infty)$ is a constant.  Let 
\begin{equation*}
	\{ X \in \mathbf{C}_{\gamma}(0) : \op{dist}(X, \op{spt} \op{graph}(\varphi+a)) < C_0 E_{\varsigma/2}(T,\varphi,\mathbf{C}_1(0)) \} 
	= \bigcup_{i=1}^N U_i 
\end{equation*}
where $\{U_i\}$ is a collection of mutually disjoint, connected, open subsets of $\mathbf{C}_{\gamma}(0)$ (and $C_0$ is as in \eqref{relharm separation0 eqn1}).  Then by \eqref{relharm separation0 eqn1}, 
\begin{equation*}
	T \llcorner \mathbf{C}_{\gamma}(0) = \sum_{i=1}^N T_i \quad\text{where}\quad T_i = T \llcorner U_i 
\end{equation*}
and $T_i$ are locally area-minimizing rectifiable currents of $\mathbf{C}_{\gamma}(0)$ with $(\partial T_i) \llcorner \mathbf{C}_{\gamma}(0) = 0$.  By Remark~\ref{relharm separation1 rmk}, we may assume that $\op{spt} T_i \cap \mathbf{C}_{\gamma/2}(0) \neq \emptyset$ for all $i \in \{1,2,\ldots,N\}$.  In particular, by the monotonicity formula $\|T_i\|(\mathbf{C}_{\gamma}(0)) \geq \omega_n (\gamma/2)^n$ for all $i \in \{1,2,\ldots,N\}$.  For each $l \in \{1,2,\ldots,q\}$ let $i(l) \in \{1,2,\ldots,N\}$ such that $\op{graph}(\varphi+a_l) \cap \mathbf{C}_{\gamma}(0) \subset U_{i(l)}$.  It follows from \eqref{relharm separation0 eqn1} that 
\begin{equation}\label{relharm separation0 eqn2}
	|a_l| \leq |y - \varphi(x)| + |y - (\varphi(x) + a_l)| 
	\leq 2 \op{dist}(X,\op{graph}\varphi) + 4qC_0 E_{\varsigma/2}(T,\varphi,\mathbf{C}_1(0))  
\end{equation}
for all $X = (x,y) \in \op{spt} T_{i(l)}$ (where $C_0$ is as in \eqref{relharm separation0 eqn1}).  Integrating \eqref{relharm separation0 eqn2} over $X \in \op{spt} T_{i(l)}$ and recalling that $\|T_{i(l)}\|(\mathbf{C}_{\gamma}(0)) \geq \omega_n (\gamma/2)^n$, we obtain 
\begin{equation}\label{relharm separation0 eqn3}
	|a_l| \leq C E_{\varsigma/2}(T,\varphi,\mathbf{C}_1(0))
\end{equation}
for all $l \in \{1,2,\ldots,q\}$, where $C = C(n,m,q,\gamma) \in (0,\infty)$ is a constant.  By the triangle inequality, \eqref{relharm separation0 eqn1}, and \eqref{relharm separation0 eqn3}, we deduce that \eqref{relharm separation0 concl} holds true.
\end{proof}

For the remainder of the subsection, we shall focus on proving Lemma~\ref{relharm separation2 lemma}, in which we bound the $L^{\infty}$-distance of a locally area-minimizing rectifiable current $T$ to the union of two or more graphs which are parallel to the graph of a single-valued $C^2$-function $\varphi$.

\begin{lemma}\label{relharm separation2 lemma}
Let $q \geq 1$ and $s \geq 2$ be integers.  For all $\gamma \in (0,1)$ and $\Lambda \in (0,\infty)$ there exists $\varepsilon = \varepsilon(n,m,q,s,\gamma,\Lambda) \in (0,1)$ and $\beta = \beta(n,m,q,s,\gamma,\Lambda) \in (0,1)$ such that the following holds true.  Let $T$ be an $n$-dimensional locally area-minimizing rectifiable current of $\mathbf{C}_1(0)$ such that 
\begin{equation}\label{relharm separation2 hyp1}
	(\partial T) \llcorner \mathbf{C}_1(0) = 0, \quad \sup_{X \in \op{spt} T} \op{dist}(X,P_0) < \infty, \quad 
	\pi_{\#} T = q \llbracket B_1(0) \rrbracket, \quad \widehat{\mathcal E}(T,P_0,\mathbf{C}_1(0)) < \varepsilon . 
\end{equation}
Let $\varphi \in C^2(B_1(0),\mathbb{R}^m)$ be a single-valued function such that $\|D\varphi\|_{C^0(B_1(0))} + \|\Delta\varphi\|_{C^0(B_1(0))} < \varepsilon^2$ and let $a \in \mathcal{A}_s(\mathbb{R}^m)$.  Suppose that 
\begin{gather}
	\label{relharm separation2 hyp2} \op{width} a \leq \Lambda \op{sep} a, \\
	\label{relharm separation2 hyp3} E_{\varsigma/2,\varepsilon}(T, \varphi+a, \mathbf{C}_1(0)) \leq \beta \op{sep} a ,
\end{gather}
where $\varsigma = \varsigma(n,m,q) \in (0,1)$ is as in Lemma~\ref{relharm tilt growth lemma},   Then 
\begin{equation}\label{relharm separation2 concl}
	\sup_{X \in \op{spt} T \cap \mathbf{C}_{\gamma}(0)} \op{dist}(X, \op{spt}\op{graph}(\varphi+a)) \leq C E_{\varsigma/2}(T,\varphi+a,\mathbf{C}_1(0)) , 
\end{equation}
where $C = C(n,m,q,\gamma) \in (0,\infty)$ is a constant.
\end{lemma}

To prove Lemma~\ref{relharm separation2 lemma}, we first need the following lemma.

\begin{lemma}\label{relharm wellsep2 lemma}
Let $q \geq 1$ and $s \geq 2$ be integers.  For all $\gamma \in (0,1)$, $\Lambda \in (0,\infty)$, and $\eta \in (0,1)$ there exists $\varepsilon = \varepsilon(n,m,q,s,\gamma,\Lambda,\eta) \in (0,1)$ and $\beta = \beta(n,m,q,s,\gamma,\Lambda,\eta) \in (0,1)$ such that the following holds true.  Let $T$ be an $n$-dimensional locally area-minimizing rectifiable current of $\mathbf{C}_1(0)$.  Let $\varphi \in C^2(B_1(0),\mathbb{R}^m)$ be a single-valued function with $\|D\varphi\|_{C^0(B_1(0))} + \|\Delta\varphi\|_{C^0(B_1(0))} < \varepsilon^2$ and let $a \in \mathcal{A}_s(\mathbb{R}^m)$.  Suppose that \eqref{relharm separation2 hyp1}, \eqref{relharm separation2 hyp2}, and \eqref{relharm separation2 hyp3} hold true, where $\varsigma = \varsigma(n,m,q) \in (0,1)$ is as in Lemma~\ref{relharm tilt growth lemma},   Then 
\begin{equation*}%\label{relharm wellsep2 concl}
	\mathcal{E}_{\varsigma/2,\varepsilon}(T, \varphi, \mathbf{C}_{\gamma}(0)) \leq \eta \op{sep} a .
\end{equation*}
\end{lemma}

\begin{proof}  
Throughout the proof, let $C \in (0,\infty)$ denote constants depending only on $n,m,q,\gamma,\Lambda$.  Fix $\eta \in (0,1)$ and integers $q \geq 1$ and $s \geq 2$.  Suppose to the contrary that for $\nu = 1,2,3,\ldots$ there are $\varepsilon_{\nu} \rightarrow 0^+$, $\beta_{\nu} \rightarrow 0^+$, $n$-dimensional locally area-minimizing rectifiable currents $T^{(\nu)}$, functions $\varphi^{(\nu)} \in C^2(B_1(0),\mathbb{R}^m)$, and $a^{(\nu)} = \sum_{i=1}^s \llbracket a^{(\nu)}_i \rrbracket \in \mathcal{A}_s(\mathbb{R}^m)$ such that $\|D\varphi^{(\nu)}\|_{C^0(B_1(0))} + \|\Delta\varphi^{(\nu)}\|_{C^0(B_1(0))} < \varepsilon_{\nu}^2$ and \eqref{relharm separation2 hyp1}, \eqref{relharm separation2 hyp2}, and \eqref{relharm separation2 hyp3} hold true with $\varepsilon_{\nu},\beta_{\nu},T^{(\nu)},\varphi^{(\nu)},a^{(\nu)}$ in place of $\varepsilon,\beta,T,\varphi,a$ but 
\begin{equation}\label{relharm wellsep2 eqn1}
	\mathcal{E}_{\varsigma/2,\varepsilon_{\nu}}(T^{(\nu)}, \varphi^{(\nu)}, \mathbf{C}_{\gamma}(0)) > \eta \op{sep} a^{(\nu)} .
\end{equation}
By Remark~\ref{relharm excess rmk1}(2), $\mathcal{E}_{\varsigma/2,\varepsilon_{\nu}}(T^{(\nu)}, \varphi^{(\nu)}, B_1(0)) \rightarrow 0$ as $\nu\rightarrow\infty$.  Hence by \eqref{relharm wellsep2 eqn1}, $\op{sep} a^{(\nu)} \rightarrow 0$ as $\nu\rightarrow\infty$.  By translating, assume that $\varphi^{(\nu)}(0) = 0$ and $a^{(\nu)}_1 = 0$.  By \eqref{relharm separation2 hyp2} 
\begin{equation}\label{relharm wellsep2 eqn2}
	|a^{(\nu)}_i| \leq \Lambda \op{sep} a^{(\nu)}
\end{equation}
for all $i \in \{1,2,\ldots,s\}$.  By the triangle inequality, \eqref{relharm separation2 hyp3}, and \eqref{relharm wellsep2 eqn2} 
\begin{equation}\label{relharm wellsep2 eqn3}
	E_{\varsigma/2,\varepsilon_{\nu}}(T^{(\nu)}, \varphi^{(\nu)}, \mathbf{C}_1(0)) \leq C \op{sep} a^{(\nu)} .
\end{equation}
By $\varphi^{(\nu)}(0) = 0$, $\|D\varphi^{(\nu)}\|_{C^0(B_1(0))} < \varepsilon_{\nu}^2$, and \eqref{relharm wellsep2 eqn3}, $\widehat{E}(T^{(\nu)},P_0,\mathbf{C}_1(0)) \leq C \op{sep} a^{(\nu)}$ and thus by~\cite[Lemma~3.5]{KrumWica}, $\op{dist}(X,P_0) \leq C \op{sep} a^{(\nu)}$ for all $X \in \op{spt}T \cap \mathbf{C}_{(1+\gamma)/2}(0)$.  By Lemma~\ref{relharm energy lemma} and \eqref{relharm wellsep2 eqn3} 
\begin{equation}\label{relharm wellsep2 eqn4}
	{\mathcal E}_{\varsigma/2,\varepsilon_{\nu}}(T^{(\nu)},\varphi^{(\nu)},\mathbf{C}_{(3+\gamma)/4}(0)) \leq C \op{sep} a^{(\nu)} .
\end{equation}
Let $u^{(\nu)} : B_{(5+3\gamma)/8}(0) \rightarrow \mathcal{A}_q(\mathbb{R}^m)$ be a Lipschitz approximation of $T^{(\nu)}$ as in Theorem~\ref{lip approx thm} with $0, \frac{3+\gamma}{4}, \frac{5+3\gamma}{6+2\gamma}, T^{(\nu)}, u^{(\nu)}$ in place of $x_0, \rho, \gamma, T, u$.  By truncating $u^{(\nu)}$, we may assume that $|u^{(\nu)}| \leq C\op{sep}a^{(\nu)}$ on $B_{(5+3\gamma)/8}(0)$.  By Lemma~\ref{relharm harmonic approx lemma} and \eqref{relharm wellsep2 eqn4}, there exists locally Dirichlet energy minimizing $q$-valued functions $w^{(\nu)} : B_{(1+\gamma)/2}(0) \rightarrow \mathcal{A}_q(\mathbb{R}^m)$ such that 
\begin{gather}
    \label{relharm wellsep2 eqn5} \lim_{\nu\rightarrow\infty} \frac{1}{(\op{sep}a^{(\nu)})^2} \int_{B_{(1+\gamma)/2}(0)} \mathcal{G}(u^{(\nu)},\varphi^{(\nu)}+w^{(\nu)})^2 = 0, \\ 
    \label{relharm wellsep2 eqn6} \lim_{\nu\rightarrow\infty} \frac{1}{(\op{sep}a^{(\nu)})^2} \int_{B_{(1+\gamma)/2}(0)} (|D(u^{(\nu)}-\varphi^{(\nu)})| - |Dw^{(\nu)}|)^2 = 0 .
\end{gather}
By the area formula, Theorem~\ref{lip approx thm}, \eqref{relharm wellsep2 eqn3}, and \eqref{relharm separation2 hyp3} 
\begin{align}
	\label{relharm wellsep2 eqn7} &\int_{B_{(1+\gamma)/2}(0)} \mathcal{G}(u^{(\nu)}, q \llbracket \varphi^{(\nu)} \rrbracket)^2 
		\leq C E_{\varsigma/2,\varepsilon_{\nu}}(T^{(\nu)}, \varphi^{(\nu)}, \mathbf{C}_1(0))^2 \leq C (\op{sep} a^{(\nu)})^2, \\
	\label{relharm wellsep2 eqn8} &\int_{B_{(1+\gamma)/2}(0)} \sum_{l=1}^q \op{dist}^2(u^{(\nu)}_l(x), \op{spt}(\varphi^{(\nu)}(x) + a^{(\nu)})) \,dx
		\\&\hspace{30mm} \leq C E_{\varsigma/2,\varepsilon_{\nu}}(T^{(\nu)}, \varphi^{(\nu)} + a^{(\nu)}, \mathbf{C}_1(0))^2 
		\leq C \beta_{\nu}^2 (\op{sep} a^{(\nu)})^2 , \nonumber
		\end{align}
where $u^{(\nu)}(x) = \sum_{l=1}^q \llbracket u^{(\nu)}_l(x) \rrbracket$ for each $x \in B_{(1+\gamma)/2}(0)$.  By \eqref{relharm wellsep2 eqn2}, after passing to a subsequence $a^{(\nu)}/\op{sep} a^{(\nu)} \rightarrow \widetilde{a}$ in $\mathcal{A}_q(\mathbb{R}^m)$ as $\nu\rightarrow\infty$.  By \eqref{relharm wellsep2 eqn7} and \eqref{relharm wellsep2 eqn5}, $\|w^{(\nu)}\|_{L^2(B_{(1+\gamma)/2}(0))} \leq C \op{sep} a^{(\nu)}$.  Hence after passing to a subsequence there exists a Dirichlet energy minimizing $q$-valued function $\widetilde{w} : B_{\gamma}(0) \rightarrow \mathcal{A}_q(\mathbb{R}^m)$ such that $w^{(\nu)}/\op{sep} a^{(\nu)} \rightarrow \widetilde{w}$ uniformly on $B_{\gamma}(0)$ and $|Dw^{(\nu)}|/\op{sep} a^{(\nu)} \rightarrow |D\widetilde{w}|$ in $L^2(B_{\gamma}(0))$ as $\nu\rightarrow\infty$.  It follows using \eqref{relharm wellsep2 eqn8} and \eqref{relharm wellsep2 eqn5} that $\op{spt} \widetilde{w}(x) \subseteq \op{spt} \widetilde{a}$ for all $x \in B_{\gamma}(0)$.  Hence by the continuity of Dirichlet energy minimizing $q$-valued functions, $\widetilde{w}$ is a constant $q$-valued function.  Therefore, by the coarea formula, Theorem~\ref{lip approx thm}, and \eqref{relharm wellsep2 eqn6} 
\begin{equation*}
	\lim_{\nu\rightarrow\infty} \frac{\mathcal{E}_{\varsigma/2,\varepsilon_{\nu}}(T^{(\nu)}, \varphi^{(\nu)}, \mathbf{C}_{\gamma}(0))}{(\op{sep} a^{(\nu)})^2} 
	= \lim_{\nu\rightarrow\infty} \frac{1}{\omega_n \gamma^n (\op{sep} a^{(\nu)})^2} \int_{B_{\gamma}(0)} |Dw^{(\nu)}|^2 =  \frac{1}{\omega_n \gamma^n} \int_{B_{\gamma}(0)} |D\widetilde{w}|^2 = 0, 
\end{equation*}
contradicting \eqref{relharm wellsep2 eqn1}.
\end{proof}

\begin{proof}[Proof of Lemma~\ref{relharm separation2 lemma}]
Let $\eta = \eta(n,m,q,s,\gamma) \in (0,1)$ to be later determined.  By Lemma~\ref{relharm wellsep2 lemma}, provided $\varepsilon = \varepsilon(n,m,q,s,\gamma,\eta,\Lambda)$ and $\beta = \beta(n,m,q,s,\gamma,\eta,\Lambda)$ are sufficiently small 
\begin{equation*}
	\mathcal{E}_{\varsigma/2,\varepsilon}(T,\varphi,\mathbf{C}_{(3+5\gamma)/8}(0)) < \eta \op{sep} a . 
\end{equation*}
Thus by Theorem~\ref{relharm separation1 thm} there exists $b = \sum_{l=1}^q \llbracket b_l \rrbracket \in \mathcal{A}_q(\mathbb{R}^m)$ such that 
\begin{equation}\label{relharm separation2 eqn1}
	\sup_{X \in \op{spt} T \cap \mathbf{C}_{(1+3\gamma)/4}(0)} \op{dist}(X, \op{spt}\op{graph}(\varphi+b)) < C_1 \eta \op{sep} a ,
\end{equation}
where $C_1 = C_1(n,m,q,s,\gamma) \in (0,\infty)$ is a constant.  Let   
\begin{equation*}
	\{ X \in \mathbf{C}_{(1+3\gamma)/4}(0) : \op{dist}(X,\op{spt}\op{graph}(\varphi+b)) < C_1 \eta \op{sep} a \} = \bigcup_{i=1}^{\widehat{N}} \widehat{U}_i 
\end{equation*}
for some collection of mutually disjoint connected open subsets $\{\widehat{U}_i\}$ of $\mathbf{C}_{(1+3\gamma)/4}(0)$ (where $C_1$ is as in \eqref{relharm separation2 eqn1}).  Then we can express 
\begin{equation*}
	T \llcorner \mathbf{C}_{(1+3\gamma)/4}(0) = \sum_{i=1}^N \widehat{T}_i \quad\text{where}\quad \widehat{T}_i = T \llcorner \widehat{U}_i 
\end{equation*}
and $\widehat{T}_i$ are locally area-minimizing rectifiable currents of $\mathbf{C}_{(1+3\gamma)/4}(0)$ with $(\partial\widehat{T}_i) \llcorner \mathbf{C}_{(1+3\gamma)/4}(0) = 0$.  By the constancy theorem and Remark~\ref{relharm separation1 rmk}
\begin{equation*}
	\pi_{\#} \widehat{T}_i = \widehat{q}_i \llbracket B_{(1+3\gamma)/4}(0) \rrbracket 
\end{equation*}
for some integers $\widehat{q}_i > 0$ with $\sum_{i=1}^{\widehat{N}} \widehat{q}_i = q$.  For each $l \in \{1,2,\ldots,q\}$, there exists $j(l) \in \{1,2,\ldots,\widehat{N}\}$ such that $\op{graph}(\varphi+b_l) \cap \mathbf{C}_{(1+3\gamma)/4}(0) \subseteq \widehat{U}_{j(l)}$.  Thus each point $X = (x,y) \in \op{spt} \widehat{T}_{j(l)}$ is a uniform distance $\leq 2q C_1 \eta \op{sep} a$ from $\op{graph}(\varphi+b_l)$.  Hence by the triangle inequality 
\begin{align}\label{relharm separation2 eqn2}
	\op{dist}(b_l, \op{spt} a) 
	=\,& \op{dist}(\varphi(x) + b_l, \op{spt}(\varphi(x)+a)) 
	\leq \op{dist}(y, \op{spt}(\varphi(x)+a)) + |y - \varphi(x) - b_l|
	\\ \leq\,& 2 \op{dist}(X, \op{spt}\op{graph}(\varphi+a)) + 4q C_1 \eta \op{sep} a \nonumber 
\end{align}
for all $X = (x,y) \in \op{spt} \widehat{T}_{j(l)}$.  Integrating \eqref{relharm separation2 eqn2} over $X \in \op{spt} \widehat{T}_{j(l)}$ and then using \eqref{relharm separation2 hyp2}
\begin{align}\label{relharm separation2 eqn3}
	\op{dist}(b_l, \op{spt} a) \leq\,& C E_{\varsigma/2,\varepsilon}(\widehat{T}_{j(l)},\varphi+a,\mathbf{C}_{(1+3\gamma)/4}(0)) + C \eta \op{sep} a
	\\ \leq\,& C E_{\varsigma/2,\varepsilon}(T,\varphi+a,\mathbf{C}_1(0)) + C \eta \op{sep} a
	\leq C (\beta+\eta) \op{sep} a , \nonumber 
\end{align}
where $C = C(n,m,q,s,\gamma) \in (0,\infty)$ are constants.  Thus by \eqref{relharm separation2 eqn1} and \eqref{relharm separation2 eqn3}
\begin{equation}\label{relharm separation2 eqn4}
	\sup_{X \in \op{spt} T \cap \mathbf{C}_{(1+3\gamma)/4}(0)} \op{dist}(X, \op{spt}\op{graph}(\varphi+a)) \leq C (\beta+\eta) \op{sep} a ,
\end{equation}
where $C = C(n,m,q,s,\gamma) \in (0,\infty)$ is a constant.  Choosing $\eta$ and $\beta$ small enough that $C (\beta+\eta) < 1/3$ (where $C$ is as in \eqref{relharm separation2 eqn4}), we deduce that 
\begin{equation*}
	T = \sum_{l=1}^s T_l \quad\text{where}\quad 
	T_l = T \llcorner \Big\{ X \in \mathbf{C}_{(1+3\gamma)/4}(0) : \op{dist}(X,\op{graph}(\varphi+a_l)) < \tfrac{1}{3} \op{sep} a \Big\} 
\end{equation*}
and $T_l$ are locally area-minimizing rectifiable currents of $\mathbf{C}_{(1+3\gamma)/4}(0)$ with $(\partial T_l) \llcorner \mathbf{C}_{(1+3\gamma)/4}(0) = 0$.  Applying Lemma~\ref{relharm separation0 lemma} with $\eta_{0,(1+3\gamma)/4\#} T_l$ and $\big(\frac{1+3\gamma}{4}\big)^{-1} \big( \varphi\big(\frac{1+3\gamma}{4}x\big) + a_l \big)$ in place of $T$ and $\varphi$, we conclude that \eqref{relharm separation2 concl} holds true.
\end{proof}

\subsection{$L^{\infty}$-distance estimate in terms $L^2$-distance relative to disjoint, not-necessarily parallel nearly harmonic graphs}\label{sec:relharm sep3}  Given integers $1 \leq p \leq s$, a suitable small constant $\varepsilon \in (0,1)$, and a constant $\kappa \in [0,\infty)$, we shall assume the following:

\begin{hypothesis}\label{relharm varphi defn}{\rm 
Let $1 \leq p \leq s$ be integers.  $\varphi : B_1(0) \rightarrow \mathcal{A}_s(\mathbb{R}^m)$ is given by 
\begin{equation}\label{relharm varphi defn eqn}
	\varphi(x) = \sum_{i=1}^p s_i \llbracket \varphi_i(x) \rrbracket
\end{equation} 
for each $x \in B_1(0)$, where $s_1,\ldots,s_p$ are positive integers such that $\sum_{i=1}^p s_i = s$ and $\varphi_i \in C^2(B_1(0),\mathbb{R}^m)$ are distinct single-valued functions such that $\|D\varphi_i\|_{C^0(B_1(0))} + \|\Delta\varphi_i\|_{C^0(B_1(0))} < \varepsilon$ for all $i \in \{1,2,\ldots,p\}$ and  
\begin{equation}\label{relharm varphi wellsep hyp}
	\sup_{B_1(0)} |D\varphi_i - D\varphi_j| \leq \kappa \inf_{B_1(0)} |\varphi_i - \varphi_j| 
\end{equation}
for all $i,j \in \{1,2,\ldots,p\}$.
}\end{hypothesis}

\begin{definition}\label{relharm minsep width} 
Let $\varphi : B_1(0) \rightarrow \mathcal{A}_s(\mathbb{R}^m)$ be given by \eqref{relharm varphi defn eqn} where $p,s_1,\ldots,s_p$ are positive integers and $\varphi_i \in C^2(B_1(0),\mathbb{R}^m)$ are distinct single-valued functions.  When $p > 1$, we define 
\begin{align*}
	\op{minsep} \varphi =\,& \min_{1 \leq i < j \leq p} \,\inf_{x \in B_1(0)} |\varphi_i(x) - \varphi_j(x)|, \\
	\op{width} \varphi =\,& \max_{1 \leq i < j \leq p} \,\sup_{x \in B_1(0)} |\varphi_i(x) - \varphi_j(x)|.
\end{align*}
When $p = 1$, we define $\op{minsep} \varphi = \infty$ and $\op{width} \varphi = 0$. 
\end{definition}

Let $p > 1$ and $\varphi$ be as in Hypothesis~\ref{relharm varphi defn}.  For all $\rho \in (0,1]$ and all $i \neq j$, by \eqref{relharm varphi wellsep hyp}
\begin{align}\label{relharm dist to sep}
	\sup_{B_{\rho}(0)} |\varphi_i - \varphi_j| 
	\leq\,& \inf_{B_{\rho}(0)} |\varphi_i - \varphi_j| + 2\rho \sup_{B_{\rho}(0)} |D\varphi_i - D\varphi_j| 
	\\ \leq\,& (1 + 2\kappa\rho) \inf_{B_{\rho}(0)} |\varphi_i - \varphi_j| . \nonumber 
\end{align}

\begin{lemma}\label{relharm wellsep3 lemma} 
Let $q \geq 1$ and $1 < p \leq s$ be integers.  For all $\gamma \in (0,1)$ and $\kappa,\lambda \in (0,\infty)$ there exists $\varepsilon = \varepsilon(n,m,q,s,\gamma,\kappa,\lambda) \in (0,1)$, $\eta = \eta(n,m,q,s,\gamma,\kappa,\lambda) \in (0,1)$, and $\beta = \beta(n,m,q,s,\gamma,\kappa,\lambda) \in (0,1)$ such that the following holds true.  Let $T$ be an $n$-dimensional locally area-minimizing rectifiable current of $\mathbf{C}_1(0)$ such that 
\begin{equation*}
	(\partial T) \llcorner \mathbf{C}_1(0) = 0, \quad \sup_{X \in \op{spt} T} \op{dist}(X,P_0) < \infty, \quad \pi_{\#} T = q \llbracket B_1(0) \rrbracket, \quad 
	\widehat{\mathcal E}(T,P_0,\mathbf{C}_1(0)) < \varepsilon .  
\end{equation*}
Let $\varphi : B_1(0) \rightarrow \mathcal{A}_s(\mathbb{R}^m)$ be as in Hypothesis~\ref{relharm varphi defn}.  Suppose that 
\begin{gather}
	\label{relharm wellsep3 hyp1} \op{width} \varphi \leq \lambda \op{minsep} \varphi , \\
	\label{relharm wellsep3 hyp2} E_{\varsigma/2}(T, \varphi, \mathbf{C}_1(0)) \leq \beta \op{minsep} \varphi ,
\end{gather}
where $\varsigma = \varsigma(n,m,q) \in (0,1)$ is as in Lemma~\ref{relharm tilt growth lemma}.  Then 
\begin{equation}\label{relharm wellsep3 concl}
	\sup_{X \in \op{spt} T \cap \mathbf{C}_{\gamma}(0)} \op{dist}(X, \op{spt} \op{graph} \varphi) \leq C E_{\varsigma/2}(T,\varphi,\mathbf{C}_1(0)) , 
\end{equation}
where $C = C(n,m,q,\gamma) \in (0,\infty)$ is a constant.
\end{lemma}

\begin{proof}
Let $\sigma = \sigma(n,m,q,\gamma,\kappa,\lambda) \in (0,(1-\gamma)/32)$ be a constant to be later determined.  Let $\{B_{\sigma}(x_k) : k = 1,2,\ldots,K\}$ be a collection of balls such that $x_k \in B_{(5+3\gamma)/8}(0)$ for each $k$, $B_{(5+3\gamma)/8}(0) \subseteq \bigcup_{k=1}^K B_{\sigma}(x_k)$, and $K \leq C(n,\gamma)\,\sigma^{-n}$.  For each $k \in \{1,2,\ldots,K\}$, define $\widehat{\varphi}_k : B_{2\sigma}(x_k) \rightarrow \mathcal{A}_s(\mathbb{R}^m)$ by 
\begin{equation*}
	\widehat{\varphi}_k(x) = \sum_{i=1}^s \llbracket \widehat{\varphi}_{k,i}(x) \rrbracket \quad\text{where}\quad 
	\widehat{\varphi}_{k,i}(x) = \varphi_1(x) + \varphi_i(x_k) - \varphi_1(x_k)  
\end{equation*}
for each $x \in B_{2\sigma}(x_k)$.  Notice that $\op{sep} \widehat{\varphi}_k(x)$ is a constant function of $x \in B_{2\sigma}(x_k)$, whose value we can denote by $\op{sep} \widehat{\varphi}_k$.  (In fact, $\op{sep} \widehat{\varphi}_k(x) = \op{sep} \varphi(x_k)$ for all $x \in B_{2\sigma}(x_k)$.)  Since $\widehat{\varphi}_{k,i}(x_k) = \varphi_i(x_k)$, we have that $\op{minsep} \varphi \leq \op{sep} \widehat{\varphi}_k$.  Moreover, by \eqref{relharm wellsep3 hyp1}
\begin{equation*}
	\op{width} \widehat{\varphi}_k \leq \op{width} \varphi \leq \lambda \op{minsep} \varphi \leq \lambda \op{sep}\widehat{\varphi}_k
\end{equation*}
By \eqref{relharm varphi wellsep hyp} and \eqref{relharm wellsep3 hyp1} 
\begin{align}\label{relharm wellsep3 eqn1}
	\sup_{B_{2\sigma}(x_k)} |\widehat{\varphi}_{k,i} - \varphi_i| 
	\leq\,& 2\sigma \sup_{B_{2\sigma}(x_k)} |D\widehat{\varphi}_{k,i} - D\varphi_i| 
	= 2\sigma \sup_{B_{2\sigma}(x_k)} |D\varphi_1 - D\varphi_i| 
	\\ \leq\,& 2\kappa\sigma \inf_{B_1(0)} |\varphi_1 - \varphi_i| 
	\leq 2\kappa\sigma \op{width} \varphi
	\leq 2\kappa\lambda\sigma \op{minsep} \varphi . \nonumber 
\end{align}
By the triangle inequality, \eqref{relharm wellsep3 hyp2}, and \eqref{relharm wellsep3 eqn1} 
\begin{align*}
	E_{\varsigma/2}(T,\widehat{\varphi}_k,\mathbf{C}_{2\sigma}(x_k)) 
	\leq\,& E_{\varsigma/2}(T,\varphi,\mathbf{C}_{2\sigma}(x_k)) + \frac{C(n,q)}{\sigma} \sup_{B_{2\sigma}(x_k)} \mathcal{G}(\widehat{\varphi}_k,\varphi) 
	\\ \leq\,& C \left( \frac{\beta}{\sigma^{n/2}} + \sigma \right) \frac{1}{\sigma} \op{minsep} \varphi 
    \leq C \left( \frac{\beta}{\sigma^{n/2}} + \sigma \right) \frac{1}{\sigma} \op{sep}\widehat{\varphi}_k \nonumber 
\end{align*}
for all $k \in \{1,2,\ldots,K\}$, where $C = C(n,q,\kappa,\lambda) \in (0,\infty)$ is a constant.  Assuming that $\sigma = \sigma(n,m,q,s,\gamma,\lambda)$ and $\beta = \beta(n,m,q,s,\gamma,\lambda,\sigma)$ are sufficiently small, by Lemma~\ref{relharm separation2 lemma} with $\eta_{(x_k,0),2\sigma\#} T$ and $(2\sigma)^{-1} \varphi(x_9 + 2\sigma x)$ in place of $T$ and $\varphi$
\begin{equation}\label{relharm wellsep3 eqn2}
	\sup_{X \in \op{spt} T \cap \mathbf{C}_{\sigma}(x_k)} \op{dist}(X, \op{spt} \op{graph} \widehat{\varphi}_k) 
	\leq C \left( \frac{\beta}{\sigma^{n/2}} + \sigma \right) \op{sep}\widehat{\varphi}_k < \frac{2}{9} \op{sep} \widehat{\varphi}_k , 
\end{equation}
where $C = C(n,m,q,\kappa,\lambda) \in (0,\infty)$ is a constant.  In the last step of \eqref{relharm wellsep3 eqn2}, we assumed that $\sigma$ is small enough that $C\sigma < 1/9$ and $\beta$ is  small enough that $C \beta \sigma^{-n/2} < 1/9$.  By \eqref{relharm wellsep3 eqn2}, 
\begin{equation*}
	T \llcorner \mathbf{C}_{\sigma}(x_k) = \sum_{i=1}^s T_{k,i} \quad\text{where}\quad 
	T_{k,i} = T \llcorner \Big\{ X \in \mathbf{C}_{\sigma}(x_k) : \op{dist}(X, \op{graph} \widehat{\varphi}_{k,i}) < \tfrac{2}{9} \op{sep} \widehat{\varphi}_k \Big\} 
\end{equation*}
and $T_{k,i}$ are locally area-minimizing rectifiable currents of $\mathbf{C}_{\sigma}(x_k)$ with $(\partial T_{k,i}) \llcorner \mathbf{C}_{\sigma}(x_k) = 0$.  Assuming that $2\kappa\lambda\sigma < 1/27$, by the triangle inequality, \eqref{relharm wellsep3 eqn1}, and \eqref{relharm wellsep3 eqn2}  
\begin{equation*}
	\sup_{X \in \op{spt} T_{k,i}} \op{dist}(X, \op{graph} \varphi_i) < \tfrac{1}{3} \op{minsep} \varphi .
\end{equation*}
That is, $\op{dist}(X, \op{spt} \op{graph} \varphi) < \frac{1}{3} \op{minsep} \varphi$ for all $X \in \op{spt} T \cap \mathbf{C}_{(5+3\gamma)/8}(0)$.  Hence 
\begin{equation*}
	T \llcorner \mathbf{C}_{(5+3\gamma)/8}(0) = \sum_{i=1}^s T_i \quad\text{where}\quad 
	T_i = T \llcorner \Big\{ X \in \mathbf{C}_{(5+3\gamma)/8}(0) : \op{dist}(X, \op{graph} \varphi_i) < \tfrac{1}{3} \op{minsep} \varphi \Big\} 
\end{equation*}
and $T_i$ are locally area-minimizing rectifiable currents of $\mathbf{C}_{(5+3\gamma)/8}(0)$ with $(\partial T_i) \llcorner \mathbf{C}_{(5+3\gamma)/8}(0) = 0$.  (In fact, $T_i \llcorner (\mathbf{C}_{(5+3\gamma)/8}(0) \cap \mathbf{C}_{\sigma}(x_k)) = T_{k,i} \llcorner (\mathbf{C}_{(5+3\gamma)/8}(0) \cap \mathbf{C}_{\sigma}(x_k))$ for all $i$ and $k$.)  Applying Lemma~\ref{relharm separation0 lemma} with $\eta_{0,(5+3\gamma)/8\#} T_i$ and $\big(\frac{5+3\gamma}{8}\big)^{-1} \varphi_i\big(\frac{5+3\gamma}{8}x\big)$ in place of $T$ and $\varphi$ gives us \eqref{relharm wellsep3 concl}.
\end{proof}

\begin{theorem}\label{relharm separation3 thm}
Let $q \geq 1$ and $1 \leq p \leq s$ be integers.  For all $\gamma \in (0,1)$ and $\kappa \in [0,\infty)$ there exists $\varepsilon = \varepsilon(n,m,q,s,\gamma,\kappa) \in (0,1)$ such that the following holds true.  Let $T$ be an $n$-dimensional locally area-minimizing rectifiable current of $\mathbf{C}_1(0)$ such that 
\begin{equation}\label{relharm separation3 hyp}
	(\partial T) \llcorner \mathbf{C}_1(0) = 0, \quad \sup_{X \in \op{spt} T} \op{dist}(X,P_0) < \infty, \quad \pi_{\#} T = q \llbracket B_1(0) \rrbracket, \quad 
	\widehat{\mathcal E}(T,P_0,\mathbf{C}_1(0)) < \varepsilon .  
\end{equation}
Let $\varphi : B_1(0) \rightarrow \mathcal{A}_s(\mathbb{R}^m)$ be as in Hypothesis~\ref{relharm varphi defn}.  Then 
\begin{equation}\label{relharm separation3 concl}
	\sup_{X \in \op{spt} T \cap \mathbf{C}_{\gamma}(0)} \op{dist}(X, \op{spt} \op{graph} \varphi) \leq C E_{\varsigma/2}(T,\varphi,\mathbf{C}_1(0)) , 
\end{equation}
where $\varsigma = \varsigma(n,m,q) \in (0,1)$ is as in Lemma~\ref{relharm tilt growth lemma} and $C = C(n,m,q,s,\gamma) \in (0,\infty)$ are constants.
\end{theorem}

\begin{proof}
We proceed by double induction on $q$ and $s$.    
The base case $s = 1$ follows from Lemma~\ref{relharm separation0 lemma}.  Suppose that $q_0 \geq 1$ and $s_0 > 1$ are integers such that 
\begin{enumerate}[topsep=0mm]
	\item[(H1)]  Theorem~\ref{relharm separation3 thm} holds true if $q = q_0$ and $s \in \{1,2,\ldots,s_0-1\}$
\end{enumerate}
and either (i) $q_0 = 1$ or (ii) $q_0 > 1$ and 
\begin{enumerate}[topsep=0mm]
	\item[(H2)]  Theorem~\ref{relharm separation3 thm} holds true if $q \in \{1,2,\ldots,q_0-1\}$ and $s \in \{1,2,\ldots,s_0\}$.
\end{enumerate}
Let $T$ be an $n$-dimensional locally area-minimizing rectifiable current of $\mathbf{C}_1(0)$ such that \eqref{relharm separation3 hyp} holds true with $q = q_0$ and $\varphi : B_1(0) \rightarrow \mathcal{A}_{s_0}(\mathbb{R}^m)$ be as in Hypothesis~\ref{relharm varphi defn} with $s = s_0$.  Without loss of generality assume that $\varphi = \sum_{l=1}^{s_0} \llbracket \varphi_l \rrbracket$ for distinct single-valued functions $\varphi_l \in C^2(B_1(0),\mathbb{R}^m)$.  Further assume that 
\begin{equation*}
	\op{minsep} \varphi = \inf_{B_1(0)} |\varphi_{s_0-1} - \varphi_{s_0}|
\end{equation*}
and set $\widetilde{\varphi} = \sum_{l=1}^{s_0-1} \llbracket \varphi_l \rrbracket$.  Note that by \eqref{relharm dist to sep}
\begin{equation}\label{relharm separation3 eqn1}
	\sup_{X \in \op{spt}\op{graph}\varphi} \op{dist}(X,\op{spt}\op{graph}\widetilde{\varphi}) \leq (1+2\kappa) \op{minsep} \varphi .
\end{equation}
Let $\widetilde{H} = E_{\varsigma/2}(T,\widetilde{\varphi},\mathbf{C}_1(0))$.  Assume that $\widetilde{H} > 0$, as otherwise by Remark~\ref{relharm excess rmk2}(1) we have that $\op{spt} T \subseteq \op{spt}\op{graph}\widetilde{\varphi} \subseteq \op{spt}\op{graph}\varphi$.  By (H1), we can apply Theorem~\ref{relharm separation3 thm} to obtain 
\begin{equation}\label{relharm separation3 eqn2}
	\sup_{X \in \op{spt} T \cap \mathbf{C}_{(3+\gamma)/4}(0)} \op{dist}(X,\op{spt}\op{graph}\widetilde{\varphi}) < C_0 \widetilde{H} 
\end{equation}
for some constant $C_0 = C_0(n,m,q_0,s_0,\gamma) \in (0,\infty)$.  Let $\beta = \beta(n,m,q_0,s_0,\gamma) \in (0,1)$ to be later determined.  If $\op{minsep}\varphi \leq \beta^{-1} \,E_{\varsigma/2,\varepsilon}(T,\varphi,\mathbf{C}_1(0))$, then by the triangle inequality and \eqref{relharm separation3 eqn1} 
\begin{equation*}
    \widetilde{H} = E_{\varsigma/2}(T,\widetilde{\varphi},\mathbf{C}_1(0)) \leq E_{\varsigma/2,\varepsilon}(T,\varphi,\mathbf{C}_1(0)) + C \op{minsep}\varphi \leq C \beta^{-1} E_{\varsigma/2,\varepsilon}(T,\varphi,\mathbf{C}_1(0)) ,
\end{equation*}
where $C = C(n,q,\kappa) \in (0,\infty)$ are constants.  Thus by \eqref{relharm separation3 eqn2}, \eqref{relharm separation3 concl} holds true for some constant $C = C(n,m,q_0,s_0,\gamma,\kappa)$.  Hence we may assume that 
\begin{equation}\label{relharm separation3 eqn3}
	E_{\varsigma/2,\varepsilon}(T,\varphi,\mathbf{C}_1(0)) \leq \beta \op{minsep} \varphi .
\end{equation}
Moreover, letting 
\begin{equation*}
	\{ X \in \mathbf{C}_{(3+\gamma)/4}(0) : \op{dist}(X,\op{spt}\op{graph}\varphi) < C_0 \widetilde{H} \} = \bigcup_{i=1}^{\widetilde{N}} \widetilde{U}_i 
\end{equation*}
for some collection of mutually disjoint connected open subsets $\{\widetilde{U}_i\}$ of $\mathbf{C}_{(3+\gamma)/4}(0)$ (where $C_0$ is as in \eqref{relharm separation3 eqn2}), by \eqref{relharm separation3 eqn2} we can express 
\begin{equation*}
	T = \sum_{i=1}^{\widetilde{N}} \widetilde{T}_i \quad\text{where}\quad \widetilde{T}_i = T \llcorner \widetilde{U}_i 
\end{equation*}
and $\widetilde{T}_i$ are locally area-minimizing rectifiable currents of $\mathbf{C}_{(3+\gamma)/4}(0)$ with $(\partial\widetilde{T}_i) \llcorner \mathbf{C}_{(3+\gamma)/4}(0) = 0$.  By the constancy theorem and~\cite[Lemma~3.5]{KrumWicb}, $\pi_{\#} \widetilde{T}_i = \widetilde{q}_i \llbracket B_{(3+\gamma)/4}(0) \rrbracket$ for some integers $\widetilde{q}_i \geq 0$ with $\sum_{i=1}^{\widetilde{N}} \widetilde{q}_i = q$.  If $\widetilde{q}_i = 0$ for some $i \in \{1,2,\ldots,\widetilde{N}\}$, then after removing from $\varphi$ those components $\varphi_l$ whose graph is contained in $\widetilde{U}_i$, by (H1) we can apply Theorem~\ref{relharm separation3 thm} (with $\eta_{0,(3+\gamma)/4\#} T_i$ in place of $T$) to deduce that \eqref{relharm separation3 concl} holds true.  Otherwise, if $\widetilde{N} > 1$ so that $\widetilde{q}_i < q_0$ for all $i \in \{1,2,\ldots,\widetilde{N}\}$, then by (H2) we can apply Theorem~\ref{relharm separation3 thm} (with $\eta_{0,(3+\gamma)/4\#} T_i$ in place of $T$) to deduce that \eqref{relharm separation3 concl} holds true.  

It remains to consider the case where $\widetilde{N} = 1$ and \eqref{relharm separation3 eqn3} holds true.  By the triangle inequality and \eqref{relharm separation3 eqn3} 
\begin{equation}\label{relharm separation3 eqn4}
    \widetilde{H} = E_{\varsigma/2}(T,\widetilde{\varphi},\mathbf{C}_1(0)) \leq E_{\varsigma/2,\varepsilon}(T,\varphi,\mathbf{C}_1(0)) + C \op{minsep}\varphi \leq C \op{minsep}\varphi ,
\end{equation}
where $C = C(n,q,\kappa) \in (0,\infty)$ are constants.  Since $\widetilde{N} = 1$, assuming $s_0 > 2$ we can relabel the components $\varphi_l$ of $\widetilde{\varphi}$ so that for each $l \in \{2,3,\ldots,s_0-1\}$ there exists $j(l) \in \{1,2,\ldots,l-1\}$ such that 
\begin{equation*}
	\inf_{B_{(3+\gamma)/4}(0)} |\varphi_l - \varphi_{j(l)}| \leq 2C_0 \widetilde{H} .
\end{equation*}
By \eqref{relharm dist to sep} 
\begin{equation*}
	\sup_{B_{(3+\gamma)/4}(0)} |\varphi_l - \varphi_{j(l)}| \leq 2 (1+2\kappa) \,C_0 \widetilde{H} .
\end{equation*}
It follows that 
\begin{equation}\label{relharm separation3 eqn5}
	\op{width}\widetilde{\varphi} \leq 2q_0 (1+2\kappa) \,C_0 \widetilde{H} .
\end{equation}
(Note that in the case $s_0 = 2$, $\widetilde{\varphi}$ has exactly one component and thus $\op{width}\widetilde{\varphi} = 0$.)  By \eqref{relharm separation3 eqn1}, \eqref{relharm separation3 eqn5}, and \eqref{relharm separation3 eqn4}

\begin{equation}\label{relharm separation3 eqn6}
	\op{width} \varphi \leq C \op{minsep} \varphi
\end{equation}
for some constant $C = C(n,m,q_0,s_0,\gamma,\kappa) \in (0,\infty)$.  In light of \eqref{relharm separation3 eqn3} and \eqref{relharm separation3 eqn6}, we can apply Lemma~\ref{relharm wellsep3 lemma} to deduce that provided $\beta$ is sufficiently small, \eqref{relharm separation3 concl} holds true.
\end{proof}

\section{A graphical representation and a blow-up procedure}\label{sec:graph rep}

Let $\varphi : \mathbb{R}^n \rightarrow \mathcal{A}_q(\mathbb{R}^m)$ be a non-zero, homogeneous degree $\alpha$, harmonic function which is translation invariant along $\{0\} \times \mathbb{R}^{n-2}$ (as in Definition~\ref{Phi-p_defn} below).  Let $T$ be an $n$-dimensional locally area-minimizing rectifiable current of $\mathbf{C}_1(0)$ such that $(\partial T) \llcorner \mathbf{C}_1(0) = 0$ and $T$ is weakly close to $\mathbb{R}^n \times \{0\}$ with integer multiplicity $q$.  In the case $\varphi$ has two or more components, further assume that $T$ is significantly closer in $L^2$ to $\op{graph} \varphi$ than to $\op{graph} \varphi'$ for any other homogeneous harmonic function $\varphi'$ with fewer non-zero components than $\varphi$ (as in \eqref{graphical hyp2}).  In particular, by setting $\varphi' = 0$ we assume that $T$ is significantly closer in $L^2$ to $\op{graph}\varphi$ than to the plane $\mathbb{R}^n \times \{0\}$.  We want to express $T$, away from the singular axis $\{0\} \times \mathbb{R}^{n-2}$, as the graph of an appropriate function over the graph of $\varphi$.  As in~\cite{KrumWicb}, there is no regularity theory applicable to $T$ that would guarantee $T$ is the graph of a multi-valued Lipschitz function away from the axis $\{0\} \times \mathbb{R}^{n-2}$.  Accordingly, in Theorem~\ref{graphical_thm} (and in Lemma~\ref{graphical_lemma0} and Lemma~\ref{graphical_lemma}) we use the height estimate, Theorem~\ref{relharm separation3 thm}, to show that away from $\{0\} \times \mathbb{R}^{n-2}$, $T$ separates as the sum of locally area-minimizing rectifiable currents $T_i$ such that $T_i$ is uniformly close to the graph of a component $\varphi_i$ of $\varphi$, and we use Almgren's Strong Lipschitz Approximation (Theorem~\ref{lip approx thm}) to approximate $T_i$ by the graph of a Lipschitz multi-valued function $v_i$ over a domain in the graph of $\varphi_i$.  As a consequence of Theorem~\ref{graphical_thm}, we can construct ``blow-ups'' of a sequence $(T^{(\nu)})$ of locally area-minimizing rectifiable currents relative to a sequence of homogeneous, harmonic $q$-valued functions $(\varphi^{(\nu)})$ (as in Definition~\ref{Phi-p_defn}).  These blow-ups play a key role in the proof of Lemma~\ref{graphical_lemma} (from which Theorem~\ref{graphical_thm} follows) and the proof of our main excess decay lemma Theorem~\ref{main excess decay thm} for locally area-minimizing rectifiable currents.

\subsection{Notation and terminology}\label{sec:graphical_rep_notation}  Throughout the remainder of the paper, we let $X = (x,y,z)$ denote points in $\mathbb{R}^{n+m}$, where $x \in \mathbb{R}^2$, $y \in \mathbb{R}^{n-2}$, and $z \in \mathbb{R}^m$.  

We define the following set $\Phi_{\alpha,q,p}$ of homogeneous, cylindrical, harmonic $q$-valued functions, which is based on the classification of homogeneous locally Dirichlet energy minimizing $q$-valued functions on $\mathbb{R}^n$ with spine dimension $n-2$ (see the end of Subsection~\ref{sec:stratify sec}) and the class $\Phi_{\varepsilon,p}(\varphi^{(0)})$ in~\cite[Definition~5.2]{KrumWic2}.  In particular, $\Phi_{\alpha,q,p}$ is the set of all $q$-valued functions $\varphi : \mathbb{R}^n \rightarrow \mathcal{A}_q(\mathbb{R}^m)$ which are homogeneous degree $\alpha$, translation invariant along $\{0\} \times \mathbb{R}^{n-2}$, and has exactly $p$ distinct non-zero components.

\begin{definition}\label{Phi-p_defn}
Let $q \geq 2$ and $p \geq 0$ be integers.  Let $\alpha = \ell_0/q_0$ for relatively prime positive integers $\ell_0$ and $q_0$ with $q_0 \leq q$.  We define $\Phi_{\alpha,q,p}$ to be the set of functions $\varphi : \mathbb{R}^n \rightarrow \mathcal{A}_q(\mathbb{R}^m)$ of the form 
\begin{equation*}
	\varphi = \sum_{j=1}^J m_j \varphi_j 
\end{equation*}
where 
\begin{enumerate}[itemsep=3mm,topsep=0mm]
	\item[{\rm (i)}]  either $\varphi_j : \mathbb{R}^n \rightarrow \mathbb{R}^m$ is the single-valued zero function $\varphi_j(\,\cdot\,) = 0$ or $\varphi_j : \mathbb{R}^n \rightarrow \mathcal{A}_{q_0}(\mathbb{R}^m)$ is the $q_0$-valued function given by 
\begin{equation*}%\label{varphi_form eqn2}
	\varphi_j(x,y) = \op{Re}(c_j (x_1+ix_2)^{\alpha})
\end{equation*}
for all $x = (x_1,x_2) \in \mathbb{R}^2$ and $y \in \mathbb{R}^{n-2}$, where $c_j \in \mathbb{C}^m \setminus \{0\}$ and $i = \sqrt{-1}$;

	\item[{\rm (ii)}]  $\varphi_j \neq \varphi_{j'}$ whenever $j \neq j'$;
	
	\item[{\rm (iii)}]  $J$ and $m_j$ are positive integers such that $J = p+1$ if $\varphi_{j_1}(\,\cdot\,) = 0$ for some $j_1 \in \{1,2,\ldots,J\}$ and $J = p$ otherwise.
\end{enumerate}
\end{definition}

\begin{remark}{\rm (1)  Throughout we let $q_j = 1$ if $\varphi_j$ identically zero and $q_j = q_0$ otherwise.  That way, each component $\varphi_j$ is a $q_j$-valued function and the multiplicities $m_j$ satisfy $\sum_{j=1}^J q_j m_j = q$.

\noindent
(2)  We have that $q_0 \sum_{j=1}^p m_j = q$ in the case that $\varphi_j \not\equiv 0$ for all $j$ and $m_{j_1} + q_0 \sum_{j \neq j_1} m_j = q$ in the case that $\varphi_{j_1} \equiv 0$ for some (unique) $j_1$.  In particular, $q_0 p \leq q$ and thus $1 \leq p \leq \lfloor q/q_0 \rfloor$.

\noindent 
(3)  We do not assume that the functions $\varphi \in \Phi_{\alpha,q, p}$ are locally Dirichlet energy minimizing nor do we assume that they are average-free.  In particular, $\varphi$ is automatically average-free whenever $\alpha$ is not an integer (i.e.~$q_0 \geq 2$), but $\varphi$ need not be average-free when $\alpha$ is an integer (i.e.~$q_0 = 1$).
} \end{remark}

\begin{definition}\label{Phi_defn}
Let $q$ be a positive integer.  Let $\alpha = \ell_0/q_0$ for relatively prime positive integers $\ell_0$ and $q_0$ with $q_0 \leq q$.  We let $\Phi_{\alpha,q} = \bigcup_{p=0}^{\lfloor q/q_0 \rfloor} \Phi_{\alpha,q,p}$.
\end{definition}

For $p \in \{0,1,\ldots,\lfloor q/q_0 \rfloor\}$, let $\varphi \in \Phi_{\alpha,q,p}$ be as in Definition~\ref{Phi-p_defn}.  For each $j \in \{1,2,\ldots,J\}$, the component $\varphi_j$ is a $q_j$-valued function with integer multiplicity $m_j$ and the graph of the component $\varphi_j$ is a smoothly immersed submanifold of $(\mathbb{R}^n \setminus \{0\} \times \mathbb{R}^{n-2}) \times {\mathbb R}^{m}$.  Thus corresponding to any ball $B \subset \mathbb{R}^n \setminus \{0\} \times \mathbb{R}^{n-2}$, there are single-valued harmonic functions $\varphi_{j,l} : B \rightarrow \mathbb{R}^m$ such that 
\begin{equation}\label{varphi_localized}
	\varphi_{j}(x,y) = \sum_{l=1}^{q_j} \llbracket \varphi_{j,l}(x,y) \rrbracket 
\end{equation}
for each $(x,y) \in B$, where we recall that $q_j = 1$ if $\varphi_j \equiv 0$ and $q_j = q_0$ otherwise.

\begin{definition}\label{v_class defn}
Let $\Omega \subset \mathbb{R}^n \setminus \{0\} \times {\mathbb R}^{n-2}$ be a connected open set.  For $p \in \{0,1,\ldots,\lfloor q/q_0 \rfloor\}$, let $\varphi \in \Phi_{\alpha,q,p}$ be as in Definition~\ref{Phi-p_defn}.  We let 
\begin{equation*}
	\mathcal{F}_{\varphi,\Omega} = \prod_{j=1}^J C^{0,1}(\op{graph} \varphi_j |_{\Omega}; \mathcal{A}_{m_j}(\mathbb{R}^m))
\end{equation*}
so that an element $v = (v_j) \in \mathcal{F}_{\varphi,\Omega}$ consists of Lipschitz functions $v_j : \op{graph} \varphi_j |_{\Omega} \rightarrow \mathcal{A}_{m_j}(\mathbb{R}^m)$.
\end{definition}

Let $v = (v_j) \in \mathcal{F}_{\varphi,\Omega}$ (as in Definition~\ref{v_class defn}).  Let $B \subset \Omega$ be an open ball and let $\varphi_{j,l}$ is as in \eqref{varphi_localized}.  For each $j \in \{1,2,\ldots,J\}$ and $l \in \{1, \ldots, q_j\}$ we define $v_{j,l} : B \to {\mathcal A}_{m_{j, k}}({\mathbb R}^{m})$ by 
\begin{equation}\label{v_localized}
	v_{j,l}(x,y) = v_j(x,y, \varphi_{j,l}(x,y)) = \sum_{h=1}^{m_j} \llbracket v_{j,l,h}(x,y) \rrbracket
\end{equation}
for each $(x,y) \in B$, where $v_{j,l,h}(x,y) \in {\mathbb R}^m$ for each $h \in \{1,2,\ldots,m_j\}$.  For each $j \in \{1,2,\ldots,J\}$, we define a $q_j m_j$-valued function $u_j : \Omega \rightarrow \mathcal{A}_{q_j m_j}(\mathbb{R}^m)$ by 
\begin{equation} \label{u defn}
	u_j(x,y) = \sum_{l=1}^{q_j} \sum_{h=1}^{m_j} \llbracket u_{j,l,h}(x,y) \rrbracket \quad\text{where}\quad u_{j,l,h}(x,y) = \varphi_{j,l}(x,y) + v_{j,l,h}(x,y)
\end{equation}
for every ball $B \subset \Omega$ and every $(x,y) \in B$.  Hence $u = \sum_{j=1}^J u_j$ is a $q$-valued function associated with $v = (v_j)$.

\begin{definition}\label{v comp minimizing defn}
Let $\Omega$, $\varphi$, and $v = (v_j) \in \mathcal{F}_{\varphi,\Omega}$ be as in Definition~\ref{Phi-p_defn}.  We say that $v$ is \emph{component-wise Dirichlet energy minimizing in $\Omega$} if for every ball $B \subset\subset \Omega$, each function $v_{j,l} : B \to {\mathcal A}_{m_j}({\mathbb R}^m)$ (as in \eqref{v_localized}) is in $W^{1,2}(B;\mathcal{A}_{m_j}(\mathbb{R}^m))$ and is locally Dirichlet energy minimizing in $B$.
\end{definition}

Let $\Omega \subset \mathbb{R}^n \setminus \{0\} \times {\mathbb R}^{n-2}$ be a connected open set which is rotationally symmetric about $\{0\} \times {\mathbb R}^{n-2}$, i.e.~if $x,\widetilde{x} \in \mathbb{R}^2$ and $y \in \mathbb{R}^{n-2}$ such that $(x,y) \in \Omega$ and $|x| = |\widetilde{x}|$ then $(\widetilde{x},y) \in \Omega$.  
Let $p \in \{0,1,\ldots,\lfloor q/q_0 \rfloor\}$ and let $\varphi \in \Phi_{\alpha,q,p}$ be as in Definition~\ref{Phi-p_defn}.  For $\varepsilon = \varepsilon(n,m,q,\Omega) \in (0,\infty)$ suitably small, assume that $\|\varphi\|_{L^{\infty}(\Omega)} < \varepsilon$.  Let $T$ be an $n$-dimensional locally area-minimizing rectifiable current of $\Omega \times \mathbb{R}^m$ and such that $(\partial T) \llcorner \Omega \times \mathbb{R}^m = 0$ and 
\begin{equation*}
	\sup_{X \in \op{spt} T} \op{dist}(X,\op{spt} \op{graph} \varphi) < \tfrac{1}{3} \inf_{(x,y) \in \Omega} \op{sep} \varphi(x,y) ,
\end{equation*}
where the separation $\op{sep} \varphi(x,y)$ is as in Definition~\ref{a sep width defn}.  Then there exists (unique) $n$-dimensional locally area-minimizing rectifiable currents $T_j$ of $\Omega \times \mathbb{R}^m$ such that 
\begin{gather*}%\label{T_separated eqn3}
	T \llcorner (\Omega \times \mathbb{R}^m) = \sum_{j=1}^J T_j , \quad\quad
	(\partial T_j) \llcorner (\Omega \times \mathbb{R}^m) = 0, \\
	\sup_{X \in \op{spt} T_j} \op{dist}(X, \op{graph} \varphi_j) < \tfrac{1}{3} \inf_{(x,y) \in \Omega} \op{sep} \varphi(x,y) . \nonumber 
\end{gather*}
Moreover, corresponding to any ball $B \subset \Omega$ there are $n$-dimensional locally area-minimizing rectifiable currents $T_{j,l}$ of $B \times \mathbb{R}^m$ such that 
\begin{gather}\label{T_separated eqn4}
	T_j \llcorner (B \times \mathbb{R}^m) = \sum_{l=1}^{q_j} T_{j,l} , \quad\quad
	(\partial T_{j,l}) \llcorner (B \times \mathbb{R}^m) = 0, \\
	\sup_{X \in \op{spt} T_{j,l}} \op{dist}(X, \op{graph} \varphi_{j,l}) < \tfrac{1}{3} \inf_{(x,y) \in \Omega} \op{sep} \varphi(x,y) , \nonumber 
\end{gather}
where $\varphi_{j,l}$ is as in \eqref{varphi_localized}.  By the constancy theorem, for each $j \in \{1,\ldots,J\}$ there exists an integer $\widehat{m}_j$ such that 
\begin{equation*}%\label{T_separated eqn5}
	\pi_{\#} T_{j,l} = \widehat{m}_j \llbracket B \rrbracket 
\end{equation*}
for every ball $B \subset \Omega$ and $l \in \{1,\ldots,q_j\}$.  Note that since $\Omega$ is connected and rotationally symmetric about $\{0\} \times {\mathbb R}^{n-2}$, $\op{graph} \varphi_j |_{\Omega}$ is connected and thus it follows that $\widehat{m}_j$ is independent of the ball $B$ and index $l$.  In particular,
\begin{equation*}
	\pi_{\#} T_j = \widehat{m}_j q_j \llbracket \Omega \rrbracket .
\end{equation*}

\subsection{Graphical representation of area-minimizing currents with small excess}\label{sec:graph_rep_results}  For each $\varsigma,\varepsilon \in (0,1]$, $\rho > 0$, $n$-dimensional rectifiable current $T$ of $\mathbf{C}_{\rho}(0)$, and $\varphi \in \Phi_{\alpha,q}$ 
\begin{align}
	\label{graphical E defn} &E_{\varsigma,\varepsilon}(T,\varphi,\mathbf{C}_{\rho}(0))^2 = \frac{1}{\omega_n \rho^{n+2}} \int_{\mathbf{C}_{\rho}(0)} 
		{\rm dist}^2(X, \op{spt} \op{graph} \varphi) \,d\|T\|(X) \\&\hspace{15mm}+ \varepsilon^{-2} \widehat{E}(T,P_0,\mathbf{C}_{\rho}(0))^{2+\varsigma}, \nonumber\\
	\label{graphical Q defn} &Q_{\varsigma,\varepsilon}(T,\varphi,\mathbf{C}_{\rho}(0))^2 = \frac{1}{\omega_n \rho^{n+2}} 
		\int_{\mathbf{C}_{\rho}(0)} {\rm dist}^2(X, \op{spt} \op{graph} \varphi) \,d\|T\|(X) \\&\hspace{15mm} + \frac{1}{\omega_n \rho^{n+2}} \sum_{j=1}^J m_j
		\int_{B_{\rho/2}(0) \cap \{r > \rho/16\}} \sum_{l=1}^{q_j} {\rm dist}^2((x,y,\varphi_{j,l}(x,y)), \op{spt} T) \,dx \,dy \nonumber 
		\\&\hspace{15mm}+ \varepsilon^{-2} \widehat{E}(T,P_0,\mathbf{C}_{\rho}(0))^{2+\varsigma} , \nonumber 
\end{align}
where $r = r(x,y) = |x|$ for each $(x,y) \in \mathbb{R}^n$, $\varphi_{j,l}$ are locally defined by \eqref{varphi_localized}, and $\widehat{E}(T,P_0,\mathbf{C}_{\rho}(0))$ is as in \eqref{coarse height excess defn} with $P = P_0$.  When $\varepsilon = 1$, we let $E_{\varsigma}(T,\varphi,\mathbf{C}_{\rho}(0)) = E_{\varsigma,1}(T,\varphi,\mathbf{C}_{\rho}(0))$ and $Q_{\varsigma}(T,\varphi,\mathbf{C}_{\rho}(0)) = Q_{\varsigma,1}(T,\varphi,\mathbf{C}_{\rho}(0))$.

Let $\lambda,\rho > 0$, $T$ be a $n$-dimensional rectifiable current of $\mathbf{C}_{\lambda\rho}(0)$, and $\varphi \in \Phi_{\alpha,q}$.  Arguing as in Remark~\ref{T and varphi scaling rmk}, $\eta_{0,\lambda}(\op{graph}\varphi) = \op{graph}\varphi_{\lambda}$ where $\varphi_{\lambda}(x,y) = \lambda^{-1} \varphi(\lambda x,\lambda y) = \lambda^{\alpha-1} \varphi(x,y)$ for each $(x,y) \in \mathbb{R}^n$.  Moreover, $E_{\varsigma}(T,\varphi,\mathbf{C}_{\lambda\rho}(0)) = E_{\varsigma}(\eta_{0,\lambda\#} T, \varphi_{\lambda}, \mathbf{C}_{\rho}(0))$ and $Q_{\varsigma}(T,\varphi,\mathbf{C}_{\lambda\rho}(0)) = Q_{\varsigma}(\eta_{0,\lambda\#} T, \varphi_{\lambda}, \mathbf{C}_{\rho}(0))$.

\begin{theorem}\label{graphical_thm}
Let $q \geq 2$ be an integer.  Let $\alpha = \ell_0/q_0 \neq 2$ where $\ell_0$ and $q_0$ are relatively prime positive integers with $q_0 \leq q$.  Let $0 < \tau < \gamma < 1$ and $p \in \{1,2,\ldots, \lfloor q/q_0 \rfloor \}$.  Let $\varsigma = \varsigma(n,m,q) \in (0,1)$ be as in Theorem~\ref{relharm separation3 thm}.  There exists $\varepsilon_0 = \varepsilon_0(n,m,q,p,\alpha,\gamma,\tau) \in (0,1)$ and $\beta_0 = \beta_0(n,m,q,p,\alpha,\gamma,\tau) \in (0,1)$ such that if $\varphi \in \Phi_{\alpha,q,p}$ and $T$ is an $n$-dimensional rectifiable current of $\mathbf{C}_1(0)$ such that 
\begin{gather}\label{graphical hyp1}
	(\partial T) \llcorner \mathbf{C}_1(0) = 0, \quad\sup_{X \in \op{spt} T} \op{dist}(X,P_0) < \infty, \quad \|T\|(\mathbf{C}_1(0)) \leq (q+1/2) \,\omega_n , \\ 
	\pi_{\#} T = q \llbracket B_1(0) \rrbracket, \quad \widehat{E}(T,P_0,\mathbf{C}_1(0)) < \varepsilon_0 . \nonumber 
\end{gather}
and 
\begin{equation}\label{graphical hyp2} 
	Q_{\varsigma/2}(T, \varphi, \mathbf{C}_1(0)) \leq \beta_0 \inf_{\varphi' \in \bigcup_{p'=0}^{p-1} \Phi_{\alpha,q,p'}} Q_{\varsigma/2}(T, \varphi', \mathbf{C}_1(0)) , 
\end{equation}
then: 
\begin{enumerate}[itemsep=3mm,topsep=0mm]
	\item[(A)]  
	\begin{equation}\label{graphical conclA1} 
		\inf_{x \in \mathbb{S}^1} \op{sep} \varphi(x,0) \geq c \inf_{\varphi' \in \bigcup_{p'=0}^{p-1} \Phi_{\alpha,q,p'}} Q_{\varsigma/2}(T, \varphi', \mathbf{C}_1(0)) , 
	\end{equation}
	where the separation $\op{sep} \varphi(x,0)$ is as in Definition~\ref{a sep width defn} and $c = c(n,m,q,p,\alpha) > 0$ is a constant, and 
	\begin{equation}\label{graphical conclA2}
		\|\varphi_{j,l} - \varphi_{j',l'}\|_{L^2(B_{1/2}(x_0,0))} \leq C \inf_{B_{1/2}(x_0,0)} |\varphi_{j,l} - \varphi_{j',l'}|
	\end{equation}
	for all $x_0 \in \mathbb{S}^1$, $j,j' \in \{1,\ldots,J\}$, $l \in \{1,\ldots,q_j\}$, and $l' \in \{1,\ldots,q_{j'}\}$, where $\varphi_{j,l}$ and $\varphi_{j',l'}$ are as in \eqref{varphi_localized} with $B = B_{1/2}(x_0,0)$ and where $C = C(n,m,q,p,\alpha) \in (0,\infty)$ is a constant;

	\item[(B)]  up to changing the multiplicities $m_j$ of $\varphi$, there exist $n$-dimensional locally area-minimizing rectifiable currents $T_j$ of $\mathbf{C}_{(1+\gamma)/2}(0) \cap \{r > \tau/2\}$ such that 
	\begin{gather*}%\label{graphical conclB}
		T \llcorner \mathbf{C}_{(1+\gamma)/2}(0) \cap \{r > \tau/2\} = \sum_{j=1}^J T_j , \\
		(\partial T_j) \llcorner \mathbf{C}_{(1+\gamma)/2}(0) \cap \{r > \tau/2\} = 0, \nonumber \\
		\pi_{\#} T_j = m_j q_j \llbracket B_{(1+\gamma)/2}(0) \cap \{r > \tau/2\} \rrbracket, \nonumber \\
		\sup_{X \in \op{spt} T_j} \op{dist}(X, \op{spt} \op{graph} \varphi_j) \leq C_{\sigma} E_{\varsigma/2}(T, \varphi, \mathbf{C}_1(0)) 
			\text{ for all } \sigma \in [\tau/2,1/2], \nonumber 
	\end{gather*}
	 where (as a slight abuse of notation) $r = r(x,y,z) = |x|$ for each $(x,y,z) \in \mathbb{R}^{n+m}$ and $C_{\sigma} = C_{\sigma}(n,m,q,p,\alpha,\gamma,\sigma) \in (0,\infty)$ is a constant;

	\item[(C)]  for each $j \in \{1,2,\ldots,J\}$ there exists a Lipschitz $m_j$-valued function $v_j : \op{graph} \varphi_j |_{B_{\gamma}(0) \cap \{r > \tau\}} \rightarrow \mathcal{A}_{m_j}(\mathbb{R}^m)$ and a closed set $K_j \subseteq B_{\gamma}(0)$ such that 
	\begin{gather*}
    	T_j \llcorner (K_j \times \mathbb{R}^m) = (\op{graph} u_j) \llcorner (K_j \times \mathbb{R}^m), \\
		\mathcal{H}^n(B_{\gamma}(0) \cap \{r > \sigma\} \setminus K_j) + \|T_j\|((B_{\gamma}(0) \cap \{r > \sigma\} \setminus K_j) \times \mathbb{R}^m) 
			\leq C_{\sigma} \widehat{E}^{2+\varsigma}, \\
		\sup_{B_{\gamma}(0) \cap \{r > \sigma\}} |u_j| \leq C_{\sigma} \widehat{E} , \quad 
		\sup_{B_{\gamma}(0) \cap \{r > \sigma\}} |Du_j| \leq C_{\sigma} \widehat{E}^{\varsigma} , \\
		\sup_{B_{\gamma}(0) \cap \{r > \sigma\}} \sum_{j=1}^J \sum_{l=1}^{q_j} |v_{j,l}|^2 
			+ \int_{B_{\gamma}(0) \cap \{r > \sigma\}} \sum_{j=1}^J \sum_{l=1}^{q_j} |Dv_{j,l}|^2 
			\leq C_{\sigma} E_{\varsigma/2}(T,\varphi,\mathbf{C}_1(0))^2 , 
	\end{gather*}
	for all $\sigma \in [\tau,\gamma)$, where $u_j$ is as in \eqref{u defn}, $v_{j,l}$ is as in \eqref{v_localized}, $\widehat{E} = \widehat{E}(T,P_0,\mathbf{C}_1(0))$, and $C_{\sigma} = C_{\sigma}(n,m,q,p,\alpha,\gamma,\sigma) \in (0,\infty)$ is a constant. 
\end{enumerate}
\end{theorem}

\begin{remark}{\rm 
In the special case that $q_0 = 1$ and $\varphi$ has exactly one component, we recall from Definition~\ref{a sep width defn} that $\op{sep} \varphi(x,0) = \infty$ for all $x \in \mathbb{S}^1$ and thus conclusion~(a) of Theorem~\ref{graphical_thm} trivially holds true.
}\end{remark}

\begin{remark}\label{changing multiplicities rmk}{\rm 
Suppose that $\varepsilon_0,\beta_0 \in (0,1)$, $p \in \{1,2,\ldots,\lfloor q/q_0 \rfloor \}$, $\varphi \in \Phi_{\alpha,q,p}$, and $T$ be an $n$-dimensional rectifiable current of $\mathbf{C}_1(0)$ such that \eqref{graphical hyp1} and \eqref{graphical hyp2} hold true.  If $\widehat{\varphi} = \sum_{j=1}^J \widehat{m}_j \varphi_j \in \Phi_{\alpha,q,p}$ where $\widehat{m}_j$ are positive integers such that $\sum_{j=1}^J q_j \widehat{m}_j = q$, then \eqref{graphical hyp1} and \eqref{graphical hyp2} continue to be hold true with $\widehat{\varphi}$, $\sqrt{q} \beta_0$ in place of $\varphi$, $\beta_0$ respectively.
}\end{remark}

In Lemma~\ref{graphical_lemma0}, Lemma~\ref{graphical_lemma}, and subsequently we shall use the following notation: 
For every $\gamma \in (0,1)$, $\zeta \in \mathbb{R}^{n-2}$, $\rho > 0$, and $\kappa \in (0,1]$, we let 
\begin{equation}\label{annuli_defn}
	A_{\rho,\kappa}(\zeta) = \left\{ (x,y) \in \mathbb{R}^2 \times \mathbb{R}^{n-2} : (|x| - \rho)^2 + |y - \zeta|^2 < \kappa^2 (1-\gamma)^2 \rho^2/64 \right\} 
\end{equation}
and $\mathbf{A}_{\rho,\kappa}(\zeta) = A_{\rho,\kappa}(\zeta) \times \mathbb{R}^m$.  For each $\varsigma \in (0,1]$, $n$-dimensional rectifiable current $T$ of $\mathbf{A}_{\rho,\kappa}(\zeta)$, and $\varphi \in \Phi_{\alpha,q,p}$ we define 
\begin{align*}
	\widehat{E}(T,P_0,\mathbf{A}_{\rho,\kappa}(\zeta))^2 
		=\,& \frac{1}{\rho^{n+2}} \int_{\mathbf{A}_{\rho,\kappa}(\zeta)} {\rm dist}^2(X,P_0) \,d\|T\|(X) , \\
	E_{\varsigma}(T,\varphi,\mathbf{A}_{\rho,\kappa}(\zeta))^2 =\,& \frac{1}{\rho^{n+2}} \int_{\mathbf{A}_{\rho,\kappa}(\zeta)} 
		{\rm dist}^2(X, \op{spt} \op{graph} \varphi) \,d\|T\|(X) + \widehat{E}(T,P_0,\mathbf{A}_{\rho,\kappa}(\zeta))^{2+\varsigma} , \\
	Q_{\varsigma}(T,\varphi,\mathbf{A}_{\rho,\kappa}(\zeta))^2 
		=\,& \frac{1}{\rho^{n+2}} \int_{\mathbf{A}_{\rho,\kappa}(\zeta)} {\rm dist}^2(X, \op{spt} \op{graph} \varphi) \,d\|T\|(X) \\& 
		+ \frac{1}{\rho^{n+2}} \sum_{j=1}^J m_j \int_{A_{\rho,\kappa/2}(\zeta)} \sum_{l=1}^{q_j} {\rm dist}^2((x,y,\varphi_{j,l}(x,y)), \op{spt} T) \,dx \,dy 
		\\&+ \widehat{E}(T,P_0,\mathbf{A}_{\rho,\kappa}(\zeta))^{2+\varsigma}, 
\end{align*}
where $\varphi_{j,l}$ are locally defined by \eqref{varphi_localized}.

\begin{lemma}\label{graphical_lemma0}
Let $q \geq 2$ be an integer.  Let $\alpha = \ell_0/q_0 \neq 2$ where $\ell_0$ and $q_0$ are relatively prime positive integers with $q_0 \leq q$.  Let $\varsigma = \varsigma(n,m,q) \in (0,1)$ be as in Theorem~\ref{relharm separation3 thm}.  Let $\gamma,\kappa \in (0,1)$, $\mu \in (0,\infty)$, and $p \in \{0,1,2,\ldots, \lfloor q/q_0 \rfloor \}$.  There exists $\overline{\varepsilon} = \overline{\varepsilon}(n,m,q,p,\alpha,\gamma,\kappa,\mu) \in (0,1)$ and $\overline{\beta} = \overline{\beta}(n,m,q,p,\alpha,\gamma,\kappa,\mu) \in (0,1)$ such that the following holds true.  Suppose that $\varphi \in \Phi_{\alpha,q,p}$ as in Definition~\ref{Phi-p_defn} and $T$ is an $n$-dimensional locally area-minimizing rectifiable current of $\mathbf{A}_{1,1}(0)$ such that $\|\varphi\|_{L^2(B_1(0))} < \overline{\varepsilon}$ and 
\begin{align}\label{graphical annuli hyp1} 
	&(\partial T) \llcorner \mathbf{A}_{1,1}(0) = 0, \quad\quad \sup_{X \in \op{spt} T} \op{dist}(X,P_0) < \infty, \\ 
	&\pi_{\#} T = q \llbracket A_{1,1}(0) \rrbracket, \quad\quad \widehat{E}(T,P_0,\mathbf{A}_{1,1}(0)) < \overline{\varepsilon} ,  \nonumber 
\end{align}
and 
\begin{equation}\label{graphical annuli hyp2} 
	\|\varphi_{j,l} - \varphi_{j',l'}\|_{L^2(B_{1/2}(x_0,0))} \leq \mu \inf_{B_{1/2}(x_0,0)} |\varphi_{j,l} - \varphi_{j',l'}|
\end{equation}
for all $x_0 \in \mathbb{S}^1$, $j,j' \in \{1,\ldots,J\}$, $l \in \{1,\ldots,q_j\}$, and $l' \in \{1,\ldots,q_{j'}\}$, where $\varphi_{j,l}$ and $\varphi_{j',l'}$ are as in \eqref{varphi_localized} with $B = B_{1/2}(x_0,0)$.  Then:
\begin{enumerate}[itemsep=3mm,topsep=0mm]
	\item[(A)]  there exists integers $\widehat{m}_j \geq 0$ with $\sum_{j=1}^J q_j \widehat{m}_j = q$ and $n$-dimensional locally area-minimizing rectifiable current $T_j$ of $\mathbf{A}_{1,(1+\kappa)/2}(0)$ such that 
	\begin{gather}
		\label{graphical annuli conclB1} T \llcorner \mathbf{A}_{1,(1+\kappa)/2}(0) = \sum_{j=1}^J T_j , \quad 
		(\partial T_j) \llcorner \mathbf{A}_{1,(1+\kappa)/2}(0) = 0, \\
		\label{graphical annuli conclB3} \pi_{\#} T_j = \widehat{m}_j q_j \llbracket A_{1,(1+\kappa)/2}(0) \rrbracket, \\
		\label{graphical annuli conclB2} \sup_{X \in \op{spt} T_j} \op{dist}(X, \op{spt} \op{graph} \varphi_j) \leq C E , 
	\end{gather}
	where $C = C(n,m,q,p,\alpha,\gamma,\kappa,\mu) \in (0,\infty)$ is a constant and 
	\begin{equation}\label{graphical annuli E defn}
		E = E_{\varsigma/2}(T,\varphi,\mathbf{A}_{1,1}(0)) + \|D\varphi\|_{C^0(B_1(0))}^{1/2} \,\widehat{E}(T,P_0,\mathbf{A}_{1,1}(0)) ;
	\end{equation}
	
	\item[(B)]  for each $j \in \{1,2,\ldots,J\}$ with $\widehat{m}_j > 0$ there exists a Lipschitz $\widehat{m}_j$-valued function $v_j : \op{graph} \varphi_j |_{A_{1,\kappa}(0)} \rightarrow \mathcal{A}_{\widehat{m}_j}(\mathbb{R}^m)$ and a closed set $K_j \subseteq A_{1,\kappa}(0)$ such that 
	\begin{gather}\label{graphical annuli conclC1}
		T_j \llcorner (K_j \times \mathbb{R}^m) = (\op{graph} u_j) \llcorner (K_j \times \mathbb{R}^m) , \\
		\mathcal{H}^n(A_{1,\kappa}(0) \setminus K_j) + \|T_j\|((A_{1,\kappa}(0) \setminus K_j) \times \mathbb{R}^m) 
			\leq C \widehat{E}(T,P_0,\mathbf{A}_{1,1}(0))^{2+\varsigma}, \nonumber\\
		\sup_{A_{1,\kappa}(0)} |u_j| \leq C \widehat{E}(T,P_0,\mathbf{A}_{1,1}(0)) , \quad 
		\sup_{A_{1,\kappa}(0)} |Du_j| \leq C \widehat{E}(T,P_0,\mathbf{A}_{1,1}(0))^{\varsigma} , \nonumber
	\end{gather}
	and 
	\begin{equation}\label{graphical annuli conclC2}
		\sup_{A_{1,\kappa}(0)} \sum_{j=1}^J \sum_{l=1}^{q_j} |v_{j,l}|^2  + \int_{A_{1,\kappa}(0)} \sum_{j=1}^J \sum_{l=1}^{q_j} |Dv_{j,l}|^2  \leq C E^2 , 
	\end{equation}
	where $u_j$ is as in \eqref{u defn} (with $\widehat{m}_j$ in place of $m_j$), $v_{j,l}$ is as in \eqref{v_localized}, and $C = C(n,m,q,p,\alpha,\gamma,\kappa) \in (0,\infty)$ is a constant. 
\end{enumerate}
\end{lemma}

\begin{lemma}\label{graphical_lemma} 
Let $q \geq 2$ be an integer.  Let $\alpha = \ell_0/q_0 \neq 2$ where $\ell_0$ and $q_0$ are relatively prime positive integers with $q_0 \leq q$.  Let $\varsigma = \varsigma(n,m,q) \in (0,1)$ be as in Theorem~\ref{relharm separation3 thm}.  Let $\gamma,\kappa,\delta \in (0,1)$ and $p \in \{0,1,2,\ldots, \lfloor q/q_0 \rfloor \}$.  There exists $\overline{\varepsilon} = \overline{\varepsilon}(n,m,q,p,\alpha,\gamma,\kappa) \in (0,1)$ and $\overline{\beta} = \overline{\beta}(n,m,q,p,\alpha,\gamma,\kappa) \in (0,1)$ such that the following holds true.  Suppose that $\varphi \in \Phi_{\alpha,q,p}$ as in Definition~\ref{Phi-p_defn} and $T$ is an $n$-dimensional locally area mininimizing rectifiable current of $\mathbf{A}_{1,1}(0)$ such that 
\begin{gather}\label{graphical annuli hyp3} 
	(\partial T) \llcorner \mathbf{A}_{1,1}(0) = 0, \quad \sup_{X \in \op{spt} T} \op{dist}(X,P_0) < \infty, \quad 
		\|T\|(\mathbf{A}_{1,1}(0)) \leq (q+1/2) \,\mathcal{L}^n(A_{1,1}(0)) , \\ 
	\pi_{\#} T = q \llbracket A_{1,1}(0) \rrbracket, \quad \widehat{E}(T,P_0,\mathbf{A}_{1,1}(0)) < \overline{\varepsilon} , \nonumber 
\end{gather}
and either 
\begin{enumerate}[itemsep=3mm,topsep=0mm]
	\item[(i)]  $p = 0$ or 
	\item[(ii)]  $p > 0$ and 
	\begin{equation}\label{graphical annuli hyp4} 
		Q_{\varsigma/2}(T,\varphi,\mathbf{A}_{1,1}(0)) < \overline{\beta} \inf_{\varphi' \in \bigcup_{p'=0}^{p-1} \Phi_{\alpha,q,p'}} 
			Q_{\varsigma/2}(T, \varphi', \mathbf{A}_{1,1}(0)) .
	\end{equation}
\end{enumerate}
Then:
\begin{enumerate}[itemsep=3mm,topsep=0mm]
	\item[(A)]  when $p > 0$ 
	\begin{equation}\label{graphical annuli conclA1} 
		\inf_{x \in \mathbb{S}^1} \op{sep} \varphi(x,0) \geq c \inf_{\varphi' \in \bigcup_{p'=0}^{p-1} \Phi_{\alpha,q,p'}} Q_{\varsigma/2}(T, \varphi', \mathbf{A}_{1,1}(0)) ,
	\end{equation}
	where the separation $\op{sep} \varphi(x,0)$ is as in Definition~\ref{a sep width defn} and $c = c(n,m,q,p,\alpha,\gamma) > 0$ is a constant, and 
	\begin{equation}\label{graphical annuli conclA2} 
		\|\varphi_{j,l} - \varphi_{j',l'}\|_{L^2(B_{1/2}(x_0,0))} \leq C \inf_{B_{1/2}(x_0,0)} |\varphi_{j,l} - \varphi_{j',l'}|
	\end{equation}
	for all $x_0 \in \mathbb{S}^1$, $j,j' \in \{1,\ldots,p\}$, $l \in \{1,\ldots,q_j\}$, and $l' \in \{1,\ldots,q_{j'}\}$, where $\varphi_{j,l}$ and $\varphi_{j',l'}$ are as in \eqref{varphi_localized} with $B = B_{1/2}(x_0,0)$ and where $C = C(n,m,q,p,\alpha,\gamma) \in (0,\infty)$ is a constant;
	
	\item[(B)]  up to changing the values of the multiplicities $m_j$ of $\varphi$, there exists $n$-dimensional locally area-minimizing rectifiable currents $T_j$ of $\mathbf{A}_{1,(1+\kappa)/2}(0)$ such that \eqref{graphical annuli conclB1} holds true, 
	\begin{gather}
    	\pi_{\#} T_j = m_j q_j \llbracket A_{1,(1+\kappa)/2}(0) \rrbracket, \nonumber\\
		\label{graphical annuli conclB4} \sup_{X \in \op{spt} T_j} \op{dist}(X, \op{spt} \op{graph} \varphi_j) \leq C E_{\varsigma/2}(T,\varphi,\mathbf{A}_{1,1}(0))
	\end{gather}
	 for some constant $C = C(n,m,q,p,\alpha,\gamma,\kappa) \in (0,\infty)$;

	\item[(C)]  for each $j \in \{1,2,\ldots,J\}$ there exists a Lipschitz $m_j$-valued function $v_j : \op{graph} \varphi_j |_{A_{1,\kappa}(0)} \rightarrow \mathcal{A}_{m_j}(\mathbb{R}^m)$ and a closed set $K_j \subseteq A_{1,\kappa}(0)$ satisfying \eqref{graphical annuli conclC1} 
	and 
	\begin{equation}\label{graphical annuli conclC3}
		\sup_{A_{1,\kappa}(0)} \sum_{j=1}^J \sum_{l=1}^{q_j} |v_{j,l}|^2  + \int_{A_{1,\kappa}(0)} \sum_{j=1}^J \sum_{l=1}^{q_j} |Dv_{j,l}|^2 
		\leq C E_{\varsigma/2}(T,\varphi,\mathbf{A}_{1,1}(0))^2 , 
	\end{equation}
	where $u_j$ is as in \eqref{u defn}, $v_j$ is as in \eqref{v_localized}, and $C = C(n,m,q,p,\alpha,\gamma,\kappa) \in (0,\infty)$ is a constant.
\end{enumerate}
\end{lemma}

Before proceeding, we observe the following consequences of \eqref{graphical hyp1} and \eqref{graphical hyp2}.  Note that similar results hold true with obvious modifications in the case that $p \in \{1,2,\ldots,\lfloor q/q_0 \rfloor\}$, $\varphi \in \Phi_{\alpha,q,p}$, and $T$ is a current of $\mathbf{A}_{1,1}(0)$ such that \eqref{graphical annuli hyp3} and \eqref{graphical annuli hyp4} hold true.

\begin{remark}\label{graphical rmk}{\rm 
(1)  \textit{There exists $\varepsilon_0 = \varepsilon_0(n,m) \in (0,\infty)$ such that if $p \in \{1,2,\ldots,\lfloor q/q_0 \rfloor\}$, $\varphi \in \Phi_{\alpha,q,p}$, and $T$ is an $n$-dimensional rectifiable current of $\mathbf{C}_1(0)$ such that \eqref{graphical hyp1} holds true, then 
\begin{equation}\label{graphical rmk eqn1}
	Q_{\varsigma/2}(T,q \llbracket 0 \rrbracket,\mathbf{C}_1(0)) \leq C \widehat{E}(T,P_0,\mathbf{C}_1(0)) ,
\end{equation}
where $C = C(n,m,q) \in (0,\infty)$ is a constant.  Consequently, if additionally \eqref{graphical hyp2} holds true for some $\beta_0 \in (0,1)$, then 
\begin{equation*}
	Q_{\varsigma/2}(T,\varphi,\mathbf{C}_1(0)) \leq C \beta_0 \widehat{E}(T,P_0,\mathbf{C}_1(0))
\end{equation*}
for some constant $C = C(n,m,q) \in (0,\infty)$.}

To see \eqref{graphical rmk eqn1}, observe that $q \llbracket 0 \rrbracket \in \Phi_{\alpha,q,0}$ and the graph of $q \llbracket 0 \rrbracket$ is the plane $P_0 = \mathbb{R}^n \times \{0\}$ with multiplicity $q$.  By~\cite[Lemma~3.5]{KrumWica} 
\begin{equation}\label{graphical rmk eqn2}
	\sup_{X \in \op{spt} T \cap \mathbf{C}_{7/8}(0)} \op{dist}(X,P_0) \leq C_0 \widehat{E}(T,P_0,\mathbf{C}_1(0)) ,
\end{equation}
where $C_0 = C_0(n,m) \in (0,\infty)$ is a constant.  By $\pi_{\#} T = q \llbracket B_1(0) \rrbracket$ and \eqref{graphical rmk eqn2}, for each $(x,y) \in B_{1/2}(0)$ there exists $z \in \mathbb{R}^m$ such that $(x,y,z) \in \op{spt} T$ and $|z| \leq C_0 \widehat{E}(T,P_0,\mathbf{C}_1(0))$ (where $C_0$ is as in \eqref{graphical rmk eqn2}).  Hence using $\|T\|(\mathbf{C}_1(0)) \leq (q+1/2) \,\omega_n$
\begin{align*}
	Q_{\varsigma/2}(T, q \llbracket 0 \rrbracket, \mathbf{C}_1(0))^2 
	\leq\,& E_{\varsigma/2}(T, q \llbracket 0 \rrbracket, \mathbf{C}_1(0))^2 + \frac{q}{\omega_n} \int_{B_{1/2}(0)} \op{dist}^2((x,y,0), \op{spt} T) \,dx \,dy
	\\ \leq\,& C \widehat{E}(T,P_0, \mathbf{C}_1(0))^2 , \nonumber 
\end{align*}
where $C = C(n,m,q) \in (0,\infty)$ is a constant.

\noindent (2)  \textit{There exist $\varepsilon_0 = \varepsilon_0(n,m,q,p,\alpha) \in (0,1)$ and $\beta_0 = \beta_0(q) \in (0,1)$ such that if $p \in \{1,2,\ldots,\lfloor q/q_0 \rfloor\}$, $\varphi \in \Phi_{\alpha,q,p}$, and $T$ is an $n$-dimensional rectifiable current of $\mathbf{C}_1(0)$ such that \eqref{graphical hyp1} holds true and }
\begin{equation}\label{graphical rmk eqn3}
	Q_{\varsigma/2}(T,\varphi,\mathbf{C}_1(0)) \leq C \beta_0 \widehat{E}(T,P_0,\mathbf{C}_1(0)) ,
\end{equation}
\textit{then }
\begin{equation}\label{graphical rmk eqn4}
	\frac{1}{C} \widehat{E}(T,P_0,\mathbf{C}_1(0)) \leq \|\varphi\|_{L^2(B_1(0))} \leq C \widehat{E}(T,P_0,\mathbf{C}_1(0)) , 
\end{equation}
\textit{where $C = C(n,m,q,\alpha) \in (1,\infty)$ is a constant.}

To see the first inequality in \eqref{graphical rmk eqn4}, by the triangle inequality, \eqref{graphical rmk eqn3}, and $\|T\|(\mathbf{C}_1(0)) \leq (q+1/2) \,\omega_n$ 
\begin{align*}
	\widehat{E}(T,P_0,\mathbf{C}_1(0)) 
	\leq\,& E_{\varsigma/2}(T,\varphi,\mathbf{C}_1(0)) + (q+1/2)^{1/2} \|\varphi\|_{L^{\infty}(B_1(0))}
	\\ \leq\,& \beta_0 \widehat{E}(T,P_0,\mathbf{C}_1(0)) + (q+1/2)^{1/2} \|\varphi\|_{L^{\infty}(B_1(0))} . 
\end{align*}
Taking $\beta_0 = 1/2$, it follows that the first inequality in \eqref{graphical rmk eqn4} holds true.  

To see the second inequality in \eqref{graphical rmk eqn4}, let $j \in \{1,2,\ldots,J\}$ and $\varphi_j(x,y) = \op{Re}(c_j (x_1+ix_2)^{\alpha})$ where $c_j \in \mathbb{C}^m \setminus \{0\}$.  Assume that the function $|{\rm Re}(c_j e^{i\alpha\theta})|$ of $\theta \in \mathbb{R}$ attains its maximum value of $|{\rm Re}(c_j)|$ at $\theta = 0$.  Then ${\rm Re}(c_j)$ and ${\rm Im}(c_j)$ are orthogonal in $\mathbb{R}^m$ and $|{\rm Im}(c_j)| \leq |{\rm Re}(c_j)|$.  Let $\varphi_{j,l}$ be as in \eqref{varphi_localized} with $B = B_{1/4}(e_1/2)$ (where $e_1 = (1,0,\ldots,0)$) and let $\varphi_{j,1}(re^{i\theta},y) = {\rm Re}(c_j r^{\alpha} e^{i\alpha\theta})$ for each $(re^{i\theta},y) \in B_{1/4}(e_1/2)$.  Let $W = \big\{ (re^{i\theta},y) \in B_{1/4}(e_1/2) : |\theta| < \frac{\pi}{6\alpha} \big\}$.  For each $(re^{i\theta},y) \in W$  
\begin{align}\label{graphical rmk eqn5}
	|\varphi_{j,1}(re^{i\theta},y)| =\,& |{\rm Re}(c_j r^{\alpha} e^{i\alpha\theta})| 
	\geq |{\rm Re}(c_j)| \,r^{\alpha} \cos(\alpha\theta) - |{\rm Im}(c_j)| \,r^{\alpha} |\sin(\alpha\theta)| 
	\\ \geq\,& \tfrac{\sqrt{3} - 1}{2} \,4^{-\alpha} |{\rm Re}(c_j)| . \nonumber 
\end{align}
If $(x,y) \in W$ such that $\op{dist}((x,y,\varphi_{j,1}(x,y)), \op{spt} T) < 1/8$, then the closest point to $(x,y,\varphi_{j,1}(x,y))$ on $\op{spt} T$ lies in $\mathbf{C}_{7/8}(0)$.  Thus by the triangle inequality, \eqref{graphical rmk eqn2}, and \eqref{graphical rmk eqn5} 
\begin{align}\label{graphical rmk eqn6}
	\op{dist}((x,y,\varphi_{j,1}(x,y)), \op{spt} T) 
	\geq\,& \op{dist}((x,y,\varphi_{j,1}(x,y)), P_0) - \sup_{Y \in \op{spt} T \cap \mathbf{C}_{7/8}(0)} \op{dist}(Y,P_0) 
	\\ =\,& |\varphi_{j,1}(re^{i\theta},y)| - \sup_{Y \in \op{spt} T \cap \mathbf{C}_{7/8}(0)} \op{dist}(Y,P_0) \nonumber
	\\ \geq\,& \tfrac{\sqrt{3} - 1}{2} \,4^{-\alpha} |{\rm Re}(c_j)| - C_0 \widehat{E}(T,P_0,\mathbf{C}_1(0)) \nonumber 
\end{align}
(where $C_0$ is as in \eqref{graphical rmk eqn2}).  Integrating \eqref{graphical rmk eqn6} over $(x,y) \in W$ and using \eqref{graphical rmk eqn3} that 
\begin{align*}
	&\min\left\{ \tfrac{1}{8}, \tfrac{\sqrt{3} - 1}{2} \,4^{-\alpha} |{\rm Re}(c_j)| - C_0 \widehat{E}(T,P_0,\mathbf{B}_1(0)) \right\} 
	\\&\hspace{15mm} \leq C \left( \int_{W} \op{dist}^2((x,y,\varphi_{j,1}(x,y)), \op{spt} T) \,dx \,dy \right)^{1/2}
	\leq C \beta_0 \widehat{E}(T,P_0,\mathbf{C}_1(0)) , 
\end{align*}
where $C_0$ is as in \eqref{graphical rmk eqn2} and $C = C(n,m,q,\alpha) \in (0,\infty)$ are constants.  Assuming $\varepsilon_0$ is sufficiently small, since $\widehat{E}(T,P_0,\mathbf{C}_1(0)) < \varepsilon_0$ we must have that $\frac{\sqrt{3} - 1}{2} \,4^{-\alpha} \op{Re}(c_j) - C_0 \widehat{E}(T,P_0,\mathbf{C}_1(0)) \leq \tfrac{1}{8}$ and thus $|c_j| \leq \sqrt{2} \,|{\rm Re}(c_j)| \leq C \widehat{E}(T,P_0,\mathbf{C}_1(0))$, where $C = C(n,m,q,\alpha) \in (0,\infty)$ is a constant, thereby proving the second inequality in \eqref{graphical rmk eqn4}.

\noindent (3)  \textit{There exist $\varepsilon_0 = \varepsilon_0(n,m,q,p,\alpha) \in (0,1)$ and $\beta_0 = \beta_0(q) \in (0,1)$ such that if $p \in \{1,2,\ldots,\lfloor q/q_0 \rfloor\}$, $\varphi \in \Phi_{\alpha,q,p}$, and $T$ is an $n$-dimensional rectifiable current of $\mathbf{C}_1(0)$ such that \eqref{graphical hyp1} and \eqref{graphical hyp2} hold true, then 
\begin{equation}\label{graphical rmk eqn7}
	\left( \int_{B_1(0)} \mathcal{G}(\varphi_j,\varphi_{j'})^2 \right)^{1/2} 
		\geq c \inf_{\varphi' \in \bigcup_{p'=1}^{p-1} \Phi_{\alpha,q,p'}} Q_{\varsigma/2}(T, \varphi',\mathbf{C}_1(0)) 
\end{equation}
whenever $j \neq j'$ and $\varphi_j, \varphi_{j'}$ are both non-zero and 
\begin{equation}\label{graphical rmk eqn8}
	\left( \int_{B_1(0)} |\varphi_j|^2 \right)^{1/2} \geq c \inf_{\varphi' \in \bigcup_{p'=1}^{p-1} \Phi_{\alpha,q,p'}} Q_{\varsigma/2}(T, \varphi',\mathbf{C}_1(0))  
\end{equation}
whenever $\varphi_j$ is non-zero, where $c = c(n,m,q,\alpha) > 0$ is a constant.}  

Notice that \eqref{graphical rmk eqn7} and \eqref{graphical rmk eqn8} are stronger conclusions than \eqref{graphical conclA1}, as it is still possible that the graph of $\varphi$ has self-intersections (or is close to self-intersecting) along an $(n-1)$-dimensional half-plane.  To see \eqref{graphical rmk eqn7}, let $\varphi' \in \Phi_{\alpha,q,p-1}$ be given by $\varphi' = m_j \varphi_{j'} + \sum_{k \neq j} m_k \varphi_k$ (i.e.~we collapse the components $\varphi_j,\varphi_{j'}$ to the same component $\varphi_{j'}$).  By the triangle inequality and \eqref{graphical hyp2}
\begin{align*}
	Q_{\varsigma/2}(T,\varphi',\mathbf{C}_1(0))^2 
	&\leq 2 Q_{\varsigma/2}(T,\varphi,\mathbf{C}_1(0))^2 + C \sup_{X \in \op{spt} \op{graph} \varphi} \op{dist}^2(X,\op{spt} \op{graph} \varphi')
	\\&\leq 2 \beta_0^2 Q_{\varsigma/2}(T,\varphi',\mathbf{C}_1(0))^2 + C \sup_{B_1(0)} \mathcal{G}(\varphi_j,\varphi_{j'})^2 , 
\end{align*}
where $C = C(n,q) \in (0,\infty)$ are constants.  Taking $\beta_0 \leq \frac{1}{2\sqrt{q}}$ and using~\cite[Lemma~A.5]{KrumWic2} gives us \eqref{graphical rmk eqn7}.  \eqref{graphical rmk eqn8} follows by a similar argument in which we set $\varphi' = m_j \llbracket 0 \rrbracket + \sum_{k \neq j} m_k \varphi_k$ (i.e.~collapse the component $\varphi_j$ to zero). 
}\end{remark}

We shall give the proofs of Theorem~\ref{graphical_thm}, Lemma~\ref{graphical_lemma0}, and Lemma~\ref{graphical_lemma} at the end of this section.  The proof of Lemma~\ref{graphical_lemma} will proceed by induction on $p$ and will make use of a blow-up procedure that inductively uses the lemma itself.  We shall next discuss this blow-up procedure, which will also play an essential role in the proof of the excess decay lemma, Theorem~\ref{main excess decay thm}.

\subsection{Blow-ups of area-minimizing currents relative to cylindrical harmonic functions}\label{sec:blowup subsec}  We follow a blow-up procedure based on that of~\cite[Subsection~6.3]{KrumWic2}.  Let $\Omega = B_1(0)$ or $\Omega = A_{1,1}(0)$.  Let $p \in \{1,2,\ldots,\lfloor q/q_0 \rfloor\}$.  Suppose that $\varepsilon_{\nu} \rightarrow 0^+$ as $\nu \rightarrow \infty$ and that $0 < \beta_{\nu} \leq \overline{\beta}$ if $\Omega = A_{1,1}(0)$, where $\overline{\beta}$ is as in Lemma~\ref{graphical_lemma}, or $\beta_{\nu} \rightarrow 0^+$ as $\nu \rightarrow \infty$ if $\Omega = B_1(0)$.  For $\nu = 1,2,3,\ldots$, let $\varphi^{(\nu)} \in \Phi_{\alpha,q,p}$ and $T^{(\nu)}$ be an $n$-dimensional locally area-minimizing rectifiable current such that 
\begin{gather*}
	(\partial T^{(\nu)}) \llcorner (\Omega \times \mathbb{R}^m) = 0, \quad\quad \sup_{X \in \op{spt} T^{(\nu)}} \op{dist}(X,P_0) < \infty, \quad\quad
		\|T^{(\nu)}\|(\Omega \times \mathbb{R}^m) \leq (q+1/2) \,\mathcal{L}^n(\Omega), \\
	\pi_{\#} T^{(\nu)} = q \llbracket \Omega \rrbracket, \quad\quad \widehat{E}(T^{(\nu)},P_0,\Omega\times\mathbb{R}^m) < \varepsilon_{\nu} , \\
	Q_{\varsigma/2}(T^{(\nu)}, \varphi^{(\nu)}, \Omega\times\mathbb{R}^m) 
		\leq \beta_{\nu} \inf_{\varphi' \in \bigcup_{p'=0}^{p-1} \Phi_{\alpha,q,p'}} Q_{\varsigma/2}(T^{(\nu)}, \varphi', \Omega) ,
\end{gather*}
where $\varsigma = \varsigma(n,m,q) \in (0,1)$ is as in Theorem~\ref{relharm separation3 thm}.  Express 
\begin{equation}\label{blowup eqn1}
	\varphi^{(\nu)} = \sum_{j=1}^J m_j \varphi^{(\nu)}_j 
\end{equation}
on $\mathbb{R}^n$, where $\varphi^{(\nu)}_j$ denote distinct $q_j$-valued components of $\varphi^{(\nu)}$ with multiplicity $m_j$ (as in Definition~\ref{Phi-p_defn}).  After passing to a subsequence of $(\nu)$, we may assume that $J$, $q_j$, and $m_j$ are independent of $\nu$.  For each $\nu$, let $\widehat{E}_{\nu} = \widehat{E}(T^{(\nu)},P_0,\Omega\times\mathbb{R}^m)$ and $E_{\nu} = E_{\varsigma/2}(T^{(\nu)},\varphi^{(\nu)},\Omega\times\mathbb{R}^m)$.  (Note that more generally, we can take $E_{\nu}$ to be any positive real number such that $E_{\varsigma/2}(T^{(\nu)},\varphi^{(\nu)},\Omega\times\mathbb{R}^m) \leq C E_{\nu}$ for a constant $C \in (0,\infty)$.)  By applying Theorem~\ref{graphical_thm} if $\Omega = B_1(0)$ and Lemma~\ref{graphical_lemma} if $\Omega = A_{1,1}(0)$, up to changing the values of the multiplicities $m_j$ and passing to a subsequence of $(\nu)$, there exists:
\begin{enumerate}[itemsep=3mm,topsep=0mm]
	\item[(i)]  open sets $\Omega_{\nu} \subset\subset \Omega \setminus \{0\} \times \mathbb{R}^{n-2}$ such that 
	\begin{equation*}
		\Omega_{\nu} \subseteq \Omega_{\nu+1} \text{ for all $\nu$,} \quad\quad 
		\Omega \setminus \{0\} \times \mathbb{R}^{n-2} = \bigcup_{\nu=1}^{\infty} \Omega_{\nu};
	\end{equation*}
	\item[(ii)]  $n$-dimensional locally area-minimizing rectifiable currents $T^{(\nu)}_j$ in $\Omega_{\nu}\times\mathbb{R}^m$ such that 
	\begin{gather*}
		T^{(\nu)} \llcorner (\Omega_{\nu} \times \mathbb{R}^m) = \sum_{j=1}^J T^{(\nu)}_j , \quad
		(\partial T^{(\nu)}_j) \llcorner (\Omega_{\nu} \times \mathbb{R}^m) = 0, \quad 
		\pi_{\#} T^{(\nu)}_j = m_j q_j \llbracket \Omega_{\nu} \rrbracket , \\
		\sup_{X \in \op{spt} T^{(\nu)}_j \cap (\Omega' \times \mathbb{R}^m)} \op{dist}(X, \op{spt} \op{graph} \varphi^{(\nu)}_j) \leq C E_{\nu} \nonumber 
	\end{gather*}
	whenever $\Omega' \subset\subset \Omega_{\nu}$, 
	where $C = C(n,m,q,\alpha,\Omega',\Omega) \in (0,\infty)$ is a constant independent of $\nu$; 
	
	\item[(iii)]  	Lipschitz $m_j$-valued functions $v^{(\nu)}_j : \op{graph} \varphi^{(\nu)}_j |_{\Omega_{\nu}} \rightarrow \mathcal{A}_{m_j}(\mathbb{R}^m)$ and closed sets $K^{(\nu)}_j \subseteq \Omega_{\nu}$ such that 
	\begin{gather*}
    %\label{blowup eqn2}
		T_j \llcorner (K^{(\nu)}_j \times \mathbb{R}^m) = (\op{graph} u^{(\nu)}_j) \llcorner (K^{(\nu)}_j \times \mathbb{R}^m) , \\
		\mathcal{H}^n(\Omega' \setminus K^{(\nu)}_j) + \|T^{(\nu)}_j\|((\Omega' \setminus K^{(\nu)}_j) \times \mathbb{R}^m) 
			\leq C \widehat{E}_{\nu}^{2+\varsigma}, \nonumber\\
		\sup_{\Omega'} |u^{(\nu)}_j| \leq C \widehat{E}_{\nu}, \quad \sup_{\Omega'} |Du^{(\nu)}_j| \leq C \widehat{E}_{\nu}^{\varsigma} \nonumber 
	\end{gather*}
	and 
	\begin{equation}\label{blowup eqn3}
		\sup_{\Omega'} \sum_{j=1}^J \sum_{l=1}^{q_j} |v^{(\nu)}_{j,l}|^2 + \int_{\Omega'} \sum_{j=1}^J \sum_{l=1}^{q_j} |Dv^{(\nu)}_{j,l}|^2 \leq C E_{\nu}^2
	\end{equation}
	whenever $\Omega' \subset\subset \Omega_{\nu}$, where $v^{(\nu)}_{j,l}$ and $u^{(\nu)}_j$ are as in \eqref{v_localized} and \eqref{u defn} with $\varphi^{(\nu)},u^{(\nu)}_j,v_{j,l}^{(\nu)}$ in place of $\varphi,u_j,v_j$ and $C = C(n,m,q,\alpha,\Omega',\Omega) \in (0,\infty)$ is a constant independent of $\nu$.
\end{enumerate}
By passing to a subsequence of $(\nu)$, we may assume that for each $j \in \{1,\ldots,J\}$ one of the following three possibilities holds:
\begin{enumerate}[itemsep=3mm,topsep=0mm]
	\item[(a)]  $\varphi^{(\nu)}_j$ is not identically zero for all $\nu = 1,2,3,\ldots$; 
	\item[(b)]  $\varphi^{(\nu)}_j \equiv 0$ for all $\nu = 1,2,3,\ldots$.
\end{enumerate}
We shall construct cylindrical functions $\varphi^{(0)}_j$ and $\varphi^{(\infty)}_j$ and functions (blow-ups) 
$$w_j : \op{graph} \varphi^{(\infty)}_j |_{\Omega \setminus \{0\} \times \mathbb{R}^{n-2}} \rightarrow \mathcal{A}_{m_j}(\mathbb{R}^m)$$
by considering these two cases as follows:

\noindent\textbf{Case~(a):}  Let $(x,y) = (re^{i\theta},y)$ denote cylindrical coordinates on $\mathbb{R}^n$, where $r > 0$, $\theta \in \mathbb{R}$, and $y \in \mathbb{R}^{n-2}$.  Let $\varphi^{(\nu)}_j(re^{i\theta},y) = \op{Re}(c^{(\nu)}_j r^{\alpha} e^{i\alpha\theta})$ for $c^{(\nu)}_j \in \mathbb{C}^m \setminus \{0\}$.  Recalling Remark~\ref{graphical rmk}(2), $|c^{(\nu)}_j| \leq C(n,m,q,\alpha,\Omega) \,\widehat{E}_{\nu}$.  Hence after passing to a subsequence of $(\nu)$ let 
\begin{equation*}%\label{blowup eqn4}
	c^{(0)}_j = \lim_{\nu\rightarrow\infty} \frac{c^{(\nu)}_j}{\widehat{E}_{\nu}} , \quad\quad 
	 c^{(\infty)}_j = \lim_{\nu\rightarrow\infty} \frac{c^{(\nu)}_j}{|c^{(\nu)}_j|} . 
\end{equation*}
(Note that $c^{(0)}_j$ may or may not equal zero, whereas $c^{(\infty)}_j$ is a unit vector.)  Define $\varphi^{(0)}_j : \mathbb{R}^n \rightarrow \mathcal{A}_{q_0}(\mathbb{R}^m)$ and $\varphi^{(\infty)}_j : \mathbb{R}^n \rightarrow \mathcal{A}_{q_0}(\mathbb{R}^m)$ by 
\begin{equation*}%\label{blowup eqn5}
	\varphi^{(0)}_j(x_1,x_2,y) = \op{Re}(c^{(0)}_j (x_1+ix_2)^{\alpha}) , \quad 
	\varphi^{(\infty)}_j(x_1,x_2,y) = \op{Re}(c^{(\infty)}_j (x_1+ix_2)^{\alpha}) 
\end{equation*}
for each $(x_1,x_2,y) \in \mathbb{R}^n$.  Observe that $\op{Re}(c^{(\nu)}_{j,k} r^{\alpha} e^{i\alpha\theta})$, $\op{Re}(c^{(0)}_j r^{\alpha} e^{i\alpha\theta})$, and $\op{Re}(c^{(\infty)}_j r^{\alpha} e^{i\alpha\theta})$ are well-defined single-valued functions of $r > 0$, $\theta \in \mathbb{R}$, and $y \in \mathbb{R}^{n-2}$ which are $2\pi q_0$-periodic as functions of $\theta$.  In particular, for each $r > 0$, $\theta \in \mathbb{R}$, and $y \in \mathbb{R}^{n-2}$, $(re^{i\theta},y,\op{Re}(c^{(\nu)}_j r^{\alpha} e^{i\alpha\theta}))$ and $(re^{i\theta},y,\op{Re}(c^{(\infty)}_j r^{\alpha} e^{i\alpha\theta}))$ are well-defined points on $\op{graph} \varphi^{(\nu)}$ and $\op{graph} \varphi^{(\infty)}$ respectively.  Define $w^{(\nu)}_j : \op{graph} \varphi^{(\infty)}_j |_{\Omega_{\nu}} \rightarrow \mathcal{A}_{m_j}(\mathbb{R}^m)$ by 
\begin{equation*}
	w^{(\nu)}_j(re^{i\theta},y,\op{Re}(c^{(\infty)}_j r^{\alpha} e^{i\alpha\theta})) 
	= \frac{v^{(\nu)}_j(re^{i\theta},y,\op{Re}(c^{(\nu)}_j r^{\alpha} e^{i\alpha\theta}))}{E_{\nu}}
\end{equation*}
for every $r > 0$, $\theta \in \mathbb{R}$, and $y \in \mathbb{R}^{n-2}$ such that $(re^{i\theta},y) \in \Omega_{\nu}$.  By \eqref{blowup eqn3} and the Rellich compactness lemma for multi-valued Sobolev functions~\cite[Proposition~2.11]{DeLSpaDirMin}, after passing to a subsequence of $(\nu)$ there exists a function $w_j : \op{graph} \varphi^{(\infty)}_j |_{\Omega \setminus \{0\} \times \mathbb{R}^{n-2}} \rightarrow \mathcal{A}_{m_j}(\mathbb{R}^m)$ such that $w^{(\nu)}_j \rightarrow w_j$ in $L^2(\op{graph} \varphi^{(\infty)}_j |_{\Omega'}, \mathcal{A}_{m_j}(\mathbb{R}^m))$ for each $\Omega' \subset\subset \Omega\setminus\{0\}\times\mathbb{R}^{n-2}$.  Arguing as in Lemma~\ref{relharm harmonic approx lemma}, $w_j$ is component-wise Dirichlet energy minimizing in $\Omega\setminus\{0\}\times\mathbb{R}^{n-2}$ (as in Definition~\ref{v comp minimizing defn}) and $|Dw^{(\nu)}_j| \rightarrow |Dw_j|$ in $L^2(\op{graph} \varphi^{(\infty)}_j |_{\Omega'})$ for each $\Omega' \subset\subset \Omega\setminus\{0\}\times\mathbb{R}^{n-2}$.

\noindent\textbf{Case~(b):}  Define $\varphi^{(0)} : \mathbb{R}^n \rightarrow \mathbb{R}^m$ and $\varphi^{(\infty)} : \mathbb{R}^n \rightarrow \mathbb{R}^m$ by $\varphi^{(0)}_j(x,y) = \varphi^{(\infty)}_j(x,y) = 0$ for all $(x,y) \in \mathbb{R}^n$.  Define $w^{(\nu)}_j : \op{graph} \varphi^{(\nu)}_j |_{\Omega_{\nu}} \rightarrow \mathcal{A}_{m_j}(\mathbb{R}^m)$ by 
\begin{equation*}
	w^{(\nu)}_j(x,y,0) = \frac{v^{(\nu)}_j(x,y,0)}{E_{\nu}} 
\end{equation*}
for all $(x,y) \in \Omega_{\nu}$.  By \eqref{blowup eqn3} and the Rellich compactness lemma for multi-valued Sobolev functions~\cite[Proposition~2.11]{DeLSpaDirMin}, after passing to a subsequence of $(\nu)$ there exists a function $w_j : \op{graph} \varphi^{(\infty)}_j |_{\Omega \setminus \{0\} \times \mathbb{R}^{n-2}} \rightarrow \mathcal{A}_{m_j}(\mathbb{R}^m)$ such that $w^{(\nu)}_j \rightarrow w_j$ in $L^2(\op{graph} \varphi^{(\infty)}_j |_{\Omega'}, \mathcal{A}_{m_j}(\mathbb{R}^m))$ for each $\Omega' \subset\subset \Omega\setminus\{0\}\times\mathbb{R}^{n-2}$.  Arguing as in Lemma~\ref{relharm harmonic approx lemma}, $w_j$ is component-wise Dirichlet energy minimizing in $\Omega\setminus\{0\}\times\mathbb{R}^{n-2}$ and $|Dw^{(\nu)}_j| \rightarrow |Dw_j|$ in $L^2(\op{graph} \varphi^{(\infty)}_j |_{\Omega'})$ for each $\Omega' \subset\subset \Omega\setminus\{0\}\times\mathbb{R}^{n-2}$.

We say that $w = (w_j)$ is a \emph{blow-up of $T^{(\nu)}$ relative to $\varphi^{(\nu)}$ by excess $E_{\nu}$}.

\begin{remark}{\rm 
Suppose instead that $\Omega = \mathbb{R}^n$ and for $\varepsilon_{\nu} \rightarrow 0^+$, $\beta_{\nu} \rightarrow 0^+$, $p \in \{1,2,\ldots,\lfloor q/q_0 \rfloor\}$, $\varphi^{(\nu)} \in \Phi_{\alpha,q,p}$, and $u^{(\nu)} \in \Phi_{\alpha,q}$ such that 
\begin{gather}
	\|\varphi^{(\nu)}\|_{L^2(B_1(0))} < \varepsilon_{\nu}, \quad \|u^{(\nu)}\|_{L^2(B_1(0))} < \varepsilon_{\nu}, \nonumber \\  
	\label{blowup eqn6} \int_{B_1(0)} \mathcal{G}(u^{(\nu)}, \varphi^{(\nu)})^2 
		\leq \beta_{\nu} \inf_{\varphi' \in \bigcup_{p'=0}^{p-1} \Phi_{\alpha,q,p'}} \int_{B_1(0)} \mathcal{G}(u^{(\nu)}, \varphi')^2 .
\end{gather}
Express $\varphi^{(\nu)}$ as in \eqref{blowup eqn1}.  After passing to a subsequence of $(\nu)$, assume that $J$, $q_j$, and $m_j$ are independent of $\nu$.  By \eqref{blowup eqn6} and~\cite[Lemma~A5 and Remark~6.6(b)]{KrumWic2}, each component of $u^{(\nu)}$ is uniformly $C E_{\nu}$-close to a unique component of $\varphi^{(\nu)}$ in $B_1(0)$, where $E_{\nu} = \left( \int_{B_1(0)} \mathcal{G}(u^{(\nu)}, \varphi^{(\nu)})^2 \right)^{1/2}$ and $C = C(n,m,q,\alpha) \in (0,\infty)$ is a constant.  In particular, after passing to a subsequence of $(\nu)$, we may express 
\begin{equation*}
	u^{(\nu)} = \sum_{j=1}^J \sum_{k=1}^{s_j} m_{j,k} u^{(\nu)}_{j,k}
\end{equation*}
where $u^{(\nu)}_{j,k}$ are distinct $q_{j,k}$-valued components of $u^{(\nu)}$ (as in Definition~\ref{Phi-p_defn}) such that either 
\begin{enumerate}[itemsep=3mm,topsep=0mm]
	\item[(i)]  $\varphi^{(\nu)}_j$ and $u^{(\nu)}_{j,k}$ are both non-zero and $\sup_{B_1(0)} \mathcal{G}(u^{(\nu)}_{j,k},\varphi^{(\nu)}_j) \leq C(n,m,q,\alpha) \,E_{\nu}$, 
	\item[(ii)]  $\varphi^{(\nu)}_j$ and $u^{(\nu)}_{j,k}$ are both identically zero, or 
	\item[(iii)]  $\varphi^{(\nu)}_j$ is identically zero, $u^{(\nu)}_{j,k}$ is non-zero, and $\sup_{B_1(0)} |u^{(\nu)}_{j,k}| \leq C(n,m,q,\alpha) \,E_{\nu}$ 
\end{enumerate}
and $s_j$ and $m_{j,k}$ are positive integers such that $\sum_{j=1}^J \sum_{k=1}^{s_j} q_{j,k} m_{j,k} = q$.  After passing to a subsequence of $(\nu)$, we may assume that $s_j$, $m_{j,k}$, and $q_{j,k}$ are independent of $\nu$.  Set $u^{(\nu)}_j = \sum_{k=1}^{s_j} m_{j,k} u^{(\nu)}_{j,k}$.  Set $\widehat{m}_j = \sum_{k=1}^{s_j} m_{j,k}$ if $\varphi^{(\nu)}_j$ is non-zero and $\widehat{m}_j = \sum_{k=1}^{s_j} q_{j,k} m_{j,k}$ if $\varphi^{(\nu)}_j \equiv 0$.  Since $u^{(\nu)}_j$ is uniformly $C E_{\nu}$-close to $\varphi^{(\nu)}_j$, there exists unique continuous function $v^{(\nu)}_j : \op{graph} \varphi^{(\nu)}_j \rightarrow \mathcal{A}_{\widehat{m}_j}(\mathbb{R}^m)$ such that $|v^{(\nu)}_j| \leq C E_{\nu}$ on $\op{graph} \varphi^{(\nu)}_j |_{B_1(0)}$ and \eqref{u defn} holds true with $\varphi^{(\nu)}_j, v^{(\nu)}_j, u^{(\nu)}_j$ in place of $\varphi_j, v_j, u_j$.  By following the procedure above, we produce a blow-up $w = (w_j)$ of $u^{(\nu)}$ relative to $\varphi^{(\nu)}$ by the excess $E_{\nu}$ in $\mathbb{R}^n$ (cf.~\cite[Remark~6.7]{KrumWic2}).
}\end{remark}

\subsection{Proofs of the graphical representation results}

\begin{proof}[Proof of Lemma~\ref{graphical_lemma0}]
Notice that if $E_{\varsigma}(T,\varphi,\mathbf{A}_{1,1}(0)) = 0$, then $\op{spt} T$ is a union of planes parallel to $P_0$ and $\op{spt} T \subseteq \op{spt} \op{graph} \varphi$.  Thus $T = q \llbracket P_0 \rrbracket \llcorner \mathbf{A}_{1,1}(0)$ and there exists $j_1 \in \{1,\ldots,J\}$ such that $\varphi_{j_1} \equiv 0$ so that the conclusion of the lemma holds true with $T_{j_1} = T \llcorner \mathbf{A}_{1,(1+\kappa)/2}(0)$, $v_{j_1} = \widehat{m}_{j_1} \llbracket 0 \rrbracket$ for some integer $\widehat{m}_{j_1}$, $K_{j_1} = A_{1,\kappa}(0)$, and $T_j = 0$ for $j \neq j_1$.  Assume that $E_{\varsigma}(T,\varphi,\mathbf{A}_{1,1}(0)) > 0$.  
Since $\varphi$ is homogeneous degee $\alpha$ and harmonic, by standard elliptic estimates \eqref{graphical annuli hyp2} implies that 
\begin{equation}\label{graphical0 eqn1} 
    \sup_{B_{(1-\kappa)(1-\gamma)/32}(x_0,y_0)} |D\varphi_{j,l} - D\varphi_{j',l'}| \leq C\mu \inf_{B_{(1-\kappa)(1-\gamma)/32}(x_0,y_0)} |\varphi_{j,l} - \varphi_{j',l'}| 
\end{equation}
$j,j' \in \{1,\ldots,J\}$, $l \in \{1,\ldots,q_j\}$, and $l' \in \{1,\ldots,q_{j'}\}$, where $\varphi_{j,l}$ and $\varphi_{j',l'}$ are as in \eqref{varphi_localized} with $B = B_{(1-\kappa)(1-\gamma)/32}(x_0,y_0)$ and $C = C(n,m,\kappa,\gamma) \in (0,\infty)$ is a constant.  
By $\|\varphi\|_{L^2(B_1(0))} < \overline{\varepsilon}$, \eqref{graphical annuli hyp1}, and \eqref{graphical0 eqn1}, we can apply Theorem~\ref{relharm separation3 thm} to obtain 
\begin{equation}\label{graphical0 eqn2} 
	\sup_{X \in \op{spt} T \cap \mathbf{C}_{(1-\kappa)(1-\gamma)/32}(x_0,y_0)} \op{dist}(X,\op{spt} \op{graph} \varphi) < CE 
\end{equation}
for all $(x_0,y_0) \in \mathbf{A}_{1,(1+\kappa)/2}(0)$, where $E$ is as in \eqref{graphical annuli E defn} and $C = C(n,m,q,\gamma,\kappa,\mu) \in (0,\infty)$ is a constant.  That is, 
\begin{equation}\label{graphical0 eqn3} 
	\sup_{X \in \op{spt} T \cap \mathbf{A}_{1,(1+\kappa)/2}(0)} \op{dist}(X,\op{spt} \op{graph} \varphi) < C E
\end{equation}
(where $C$ is as in \eqref{graphical0 eqn2}).  Express  
\begin{equation*}
	\{ X \in \mathbf{A}_{1,(1+\kappa)/2}(0) : \op{dist}(X,\op{spt} \op{graph} \varphi) < C E \} = \bigcup_{k=1}^N U_k 
\end{equation*}
for some finite collection $\{U_k\}$ of pairwise disjoint connected open subsets of $\mathbf{A}_{1,(1+\kappa)/2}(0)$ (where $C$ is as in \eqref{graphical0 eqn2}).  For each $k \in \{1,2,\ldots,N\}$, there exists $j(k) \in \{1,2,\ldots,J\}$ such that $\op{graph} \varphi_{j(k)} \subset U_k$.  Set $T_{j(k)} = T \llcorner U_k$ for all $k \in \{1,2,\ldots,N\}$ and set $T_j = 0$ of $j \not\in \{j(1),j(2),\ldots,j(N)\}$.  One readily verifies that \eqref{graphical annuli conclB1} and \eqref{graphical annuli conclB2} hold true.  By the constancy theorem,~\cite[Lemma~3.5]{KrumWicb}, \eqref{graphical annuli conclB1}, and $\pi_{\#} T = q \llbracket A_{1,1}(0) \rrbracket$, there exist integers $\widehat{m}_j \geq 0$ such that $\sum_{j=1}^J q_j \widehat{m}_j = q$ and \eqref{graphical annuli conclB3} holds true.  This completes the proof of Conclusion~(A).  To see Conclusion~(B), take any $(x_0,y_0) \in A_{1,\kappa}(0)$ and let $T_{j,l}$ be as in \eqref{T_separated eqn4} where $B = B_{(1-\kappa)(1-\gamma)/32}(x_0,y_0)$.  Let $u_{j,l} : B_{(1-\kappa)(1-\gamma)/128}(x_0,y_0) \rightarrow \mathcal{A}_{\widehat{m}_j}(\mathbb{R}^m)$ be a Lipschitz function $K_{j,l} \subseteq B_{(1-\kappa)(1-\gamma)/128}(x_0,y_0)$ be a closed set as in Theorem~\ref{lip approx thm} with $(x_0,y_0), (1-\kappa)(1-\gamma)/32, 1/4, T_{j,l}, u_{j,l}, K_{j,l}$ in place of $x_0,\rho,\gamma,T,u,K$.  By replacing $K_{j,l}$ with $K_j = \bigcap_{l'=1}^{q_j} K_{j,l'}$, assume that $K_{j,l} = K_j$ independent of $l$.  Define $v_{j,l} : \op{graph} \varphi_{j,l} |_{B_{(1-\kappa)(1-\gamma)/64}(x_0,y_0)} \rightarrow \mathcal{A}_{\widehat{m}_j}(\mathbb{R}^m)$ by 
\begin{equation*}
	v_{j,l}(x,y,\varphi_{j,l}(x,y)) = \sum_{h=1}^{\widehat{m}_j} \llbracket u_{j,l,h}(x,y) - \varphi_{j,l}(x,y) \rrbracket 
\end{equation*}
for each $(x,y) \in B_{(1-\kappa)(1-\gamma)/128}(x_0,y_0)$, where $u_{j,l}(x,y) = \sum_{h=1}^{\widehat{m}_j} \llbracket u_{j,l,h}(x,y) \rrbracket$.  Using a partition of unity argument (like in the proof of~\cite[Theorem~4.4]{KrumWicb}), we can construct $v_j : \op{graph} \varphi_j |_{A_{1,\kappa}(0)} \rightarrow \mathcal{A}_{\widehat{m}_j}(\mathbb{R}^m)$ and a closed set $K \subseteq A_{1,\kappa}(0)$ satisfying \eqref{graphical annuli conclC1}, where $u_j$ is as in \eqref{u defn}.  By \eqref{graphical0 eqn3} and truncation, we can guarantee the supremum estimate for $v_j$ in \eqref{graphical annuli conclC2}.  By Lemma~\ref{relharm energy lemma} and \eqref{graphical annuli conclC1}, we obtain the gradient $L^2$-estimate for $v_j$ in \eqref{graphical annuli conclC2}.
\end{proof}

The proof of Lemma~\ref{graphical_lemma} proceeds by induction on $p$, first proving the base case $p = 0$ and then assuming for $p_0 \in \{1,2,\ldots,\lfloor q/q_0 \rfloor\}$ that 
\begin{enumerate}[topsep=0mm]
	\item[(H3)]  Lemma~\ref{graphical_lemma} holds true for all $p \in \{0,1,2,\ldots,p_0-1\}$.
\end{enumerate}

\begin{remark}\label{tilde varphi rmk}{\rm  
(1)  Let $\widetilde{\beta} \in (0,1)$ be an arbitrary constant.  Suppose that $\varsigma \in (0,1]$, $1 \leq p_0 \leq q$, $\varphi \in \Phi_{\alpha,q,p_0}$, and $T$ is an $n$-dimensional rectifiable current of $\mathbf{A}_{1,1}(0)$.  Arguing as in~\cite[Remark~4.8(1)]{KrumWicb} we can choose an integer $\widetilde{p} \in \{0,1,2,\ldots,p_0-1\}$ and $\widetilde{\varphi} \in \Phi_{\alpha,q,\widetilde{p}}$ such that 
\begin{equation}\label{tilde varphi rmk eqn1}
	Q_{\varsigma/2}(T,\widetilde{\varphi},\mathbf{A}_{1,1}(0)) 
	\leq 2^{p_0} \widetilde{\beta}^{1-p_0} \inf_{\varphi' \in \bigcup_{p'=0}^{p_0-1} \Phi_{\alpha,q,p'}} Q_{\varsigma/2}(T,\varphi',\mathbf{A}_{1,1}(0)) 
\end{equation}
and either $\widetilde{p} = 0$ or $\widetilde{p} > 0$ and 
\begin{equation}\label{tilde varphi rmk eqn2} 
	Q_{\varsigma/2}(T,\widetilde{\varphi},\mathbf{A}_{1,1}(0)) 
	\leq \widetilde{\beta} \inf_{\varphi' \in \bigcup_{p'=0}^{\widetilde{p}-1} \Phi_{\alpha,q,p'}} Q_{\varsigma/2}(T,\varphi',\mathbf{A}_{1,1}(0)) .
\end{equation}

\noindent (2)  Suppose that $p_0 \in \{1,2,\ldots,\lfloor q/q_0 \rfloor\}$ such that (H3) holds true and that $\varphi \in \Phi_{\alpha,q,p_0}$ and $T$ satisfies the hypotheses of Lemma~\ref{graphical_lemma}.  Let $\varsigma = \varsigma(n,m,q) \in (0,1)$ be as in Theorem~\ref{relharm separation3 thm} and let $\widetilde{\beta} = \widetilde{\beta}(n,m,q,p_0,\alpha,\gamma) \in (0,1)$ be a suitably small constant.  Let $\widetilde{p} \in \{0,1,2,\ldots,p_0-1\}$ and $\widetilde{\varphi} \in \Phi_{\alpha,q,\widetilde{p}}$ such that \eqref{tilde varphi rmk eqn1} and \eqref{tilde varphi rmk eqn2} holds true.  Let 
\begin{equation}\label{tilde varphi rmk eqn3}
	\widetilde{\varphi} = \sum_{j=1}^{\widetilde{J}} \widetilde{m}_j \widetilde{\varphi}_j
\end{equation}
with distinct $\widetilde{q}_j$-valued components $\widetilde{\varphi}_j$ with multiplicity $\widetilde{m}_j$ (as in Definition~\ref{Phi-p_defn}).  Here $\widetilde{q}_j = 1$ if $\widetilde{\varphi}_j \equiv 0$ and $\widetilde{q}_j = q_0$ otherwise.  By (H3), \eqref{graphical annuli hyp3}, and \eqref{tilde varphi rmk eqn2}, we can apply Lemma~\ref{graphical_lemma} with $\widetilde{\varphi}$ in place of $\varphi$ to deduce the following.  By Lemma~\ref{graphical_lemma}(A), if $\widetilde{p} > 0$ then 
\begin{equation}\label{tilde varphi rmk eqn4} 
	\inf_{x \in \mathbb{S}^1} \op{sep} \widetilde{\varphi}(x,0) \geq c \inf_{\varphi' \in \bigcup_{p'=0}^{\widetilde{p}-1} \Phi_{\alpha,q,p'}} 
	Q_{\varsigma/2}(T, \varphi', \mathbf{A}_{1,1}(0)) , 
\end{equation}
where $c = c(n,m,q,p_0,\alpha,\gamma) > 0$ is a constant, and 
\begin{equation}\label{tilde varphi rmk eqn5} 
	\|\widetilde{\varphi}_{j,l} - \widetilde{\varphi}_{j',l'}\|_{L^2(B_{1/2}(x_0,0))} 
	\leq C \inf_{B_{1/2}(x_0,0)} |\widetilde{\varphi}_{j,l} - \widetilde{\varphi}_{j',l'}|
\end{equation}
for all $x_0 \in \mathbb{S}^1$, $j,j' \in \{1,\ldots,\widetilde{J}\}$, $l \in \{1,\ldots,\widetilde{q}_j\}$, and $l' \in \{1,\ldots,\widetilde{q}_{j'}\}$, where $\widetilde{\varphi}_j = \sum_{l=1}^{\widetilde{q}_j} \llbracket \widetilde{\varphi}_{j,l} \rrbracket$ in $B = B_{1/2}(x_0,0)$ for some single-valued harmonic functions $\widetilde{\varphi}_{j,l} : B_{1/2}(x_0,0) \rightarrow \mathbb{R}^m$ (as in \eqref{varphi_localized}) and $C = C(n,m,q,p_0,\alpha,\gamma) \in (0,\infty)$ is a constant.  By Lemma~\ref{graphical_lemma}(B), up to changing the values of $\widetilde{m}_j$, there exists $n$-dimensional locally area-minimizing rectifiable current $\widetilde{T}_j$ in $\mathbf{A}_{1,(7+\kappa)/8}(0)$ such that 
\begin{gather}\label{tilde varphi rmk eqn6}
	T \llcorner \mathbf{A}_{1,(7+\kappa)/8}(0) = \sum_{j=1}^{\widetilde{J}} \widetilde{T}_j , \quad
	(\partial \widetilde{T}_j) \llcorner \mathbf{A}_{1,(7+\kappa)/8}(0) = 0, \quad
	\pi_{\#} \widetilde{T}_j = \widetilde{m}_j q_j \llbracket \mathbf{A}_{1,(7+\kappa)/8}(0) \rrbracket , \\	
	\sup_{X \in \op{spt} \widetilde{T}_j} \op{dist}(X, \op{spt} \op{graph} \widetilde{\varphi}_j) \leq C E_{\varsigma/2}(T,\widetilde{\varphi},\mathbf{A}_{1,1}(0)) , 
		\nonumber
\end{gather}
where $C = C(n,m,q,p_0,\alpha,\gamma,\kappa) \in (0,\infty)$ is a constant.  By Lemma~\ref{graphical_lemma}(C), for each $j \in \{1,2,\ldots,\widetilde{J}\}$ there exists a Lipschitz $\widetilde{m}_j$-valued function $\widetilde{v}_j : \op{graph} \widetilde{\varphi}_j |_{A_{1,(3+\kappa)/4}(0)} \rightarrow \mathcal{A}_{\widetilde{m}_j}(\mathbb{R}^m)$ and a closed set $\widetilde{K}_j \subseteq A_{1,(3+\kappa)/4}(0)$ such that 
\begin{gather}\label{tilde varphi rmk eqn7}
	\widetilde{T}_j \llcorner (\widetilde{K}_j \times \mathbb{R}^m) = (\op{graph} \widetilde{u}_j) \llcorner (\widetilde{K}_j \times \mathbb{R}^m) , \\
	\mathcal{H}^n(A_{1,(3+\kappa)/4}(0) \setminus \widetilde{K}_j) 
		+ \|\widetilde{T}_j\|((A_{1,(3+\kappa)/4}(0) \setminus \widetilde{K}_j) \times \mathbb{R}^m) 
		\leq C \widehat{E}(T,P_0,\mathbf{A}_{1,1}(0))^{2+\varsigma}, \nonumber\\
	\sup_{A_{1,(3+\kappa)/4}(0)} |\widetilde{u}_j| \leq C \widehat{E}(T,P_0,\mathbf{A}_{1,1}(0)) , \quad 
	\sup_{A_{1,(3+\kappa)/4}(0)} |D\widetilde{u}_j| \leq C \widehat{E}(T,P_0,\mathbf{A}_{1,1}(0))^{\varsigma} , \nonumber\\
	\sup_{A_{1,(3+\kappa)/4}(0)} \sum_{j=1}^{\widetilde{J}} \sum_{l=1}^{\widetilde{q}_j} |\widetilde{v}_{j,l}|^2 
		+ \int_{A_{1,(3+\kappa)/4}(0)} \sum_{j=1}^{\widetilde{J}} \sum_{l=1}^{\widetilde{q}_j} |D\widetilde{v}_{j,l}|^2 
		\leq C E_{\varsigma/2}(T,\widetilde{\varphi},\mathbf{A}_{1,1}(0))^2 , \nonumber
\end{gather}
where $\widetilde{v}_{j,l}$ and $\widetilde{u}_j$ are as in \eqref{v_localized} and \eqref{u defn} with $\widetilde{\varphi},\widetilde{v}_j,\widetilde{u}_j$ in place of $\varphi,v_j,u_j$ and $C = C(n,m,q,p_0,\alpha,\gamma,\kappa) \in (0,\infty)$ is a constant.

\noindent (3)   
Let $\varphi = \sum_{j=1}^J m_j \varphi_j$ with $q_j$-valued components $\varphi_j$ with positive integer multiplicity $m_j$. By the triangle inequality, \eqref{graphical annuli hyp4}, \eqref{tilde varphi rmk eqn6}, and \eqref{graphical hyp2}
\begin{align*}
	&\sum_{j=1}^J m_j \int_{\mathbf{A}_{1,1/2}(0)} \sum_{l=1}^{q_j} \op{dist}^2((x,y,\varphi_{j,l}(x,y)), \op{spt}\op{graph}\widetilde{\varphi}) \,dx \,dy 
	\\ \leq\,& 2 \sum_{j=1}^J m_j \int_{\mathbf{A}_{1,1/2}(0)} \sum_{l=1}^{q_j} \op{dist}^2((x,y,\varphi_{j,l}(x,y)), \op{spt} T) \,dx \,dy 
		\\&+ C \sup_{X \in \op{spt} T \cap \mathbf{A}_{1/3,4}(0)} \op{dist}^2(X, \op{spt}\op{graph}\widetilde{\varphi}) 
	\\ \leq\,& 
		C Q_{\varsigma/2}(T, \widetilde{\varphi}, \mathbf{A}_{1,1}(0)) , 
\end{align*}
where $\varphi_{j,l}$ are as in \eqref{varphi_localized} 
 and $C = C(n,m,q,p_0,\alpha,\gamma) \in (0,\infty)$ are constants.  On the other hand, by the triangle inequality, \eqref{tilde varphi rmk eqn7}, and \eqref{graphical hyp2} 
\begin{align*}
	&\sum_{j=1}^{\widetilde{J}} \widetilde{m}_j \int_{A_{1,1/2}(0)} \sum_{l=1}^{\widetilde{q}_j} 
		\op{dist}^2((x,y,\widetilde{\varphi}_{j,l}(x,y)), \op{spt}\op{graph}\varphi) \,dx \,dy 
	\\ \leq\,& 2 \int_{A_{1,1/2}(0)} \sum_{j=1}^{\widetilde{J}} \sum_{l=1}^{\widetilde{q}_j} \sum_{h=1}^{\widetilde{m}_j} 
		\op{dist}^2((x,y,\widetilde{u}_{j,l,h}(x,y)), \op{spt}\op{graph}\varphi) \,dx \,dy 
		+ C \sup_{A_{1,1/2}(0)}  \sum_{j=1}^{\widetilde{J}} \sum_{l=1}^{\widetilde{q}_j} |\widetilde{v}_{j,l}|^2
	\\ \leq\,& 
		C Q_{\varsigma/2}(T, \widetilde{\varphi}, \mathbf{A}_{1,1}(0)) , 
\end{align*}
where $C = C(n,m,q,p_0,\alpha,\gamma) \in (0,\infty)$ is a constant.  Noting that $\varphi,\widetilde{\varphi} \in \Phi_{\alpha,q}$, it follows that 
\begin{equation}\label{tilde varphi rmk eqn8}
	\op{dist}_{\mathcal H}(\op{spt} \op{graph} \varphi \cap \mathbf{C}_1(0), \op{spt} \op{graph} \widetilde{\varphi} \cap \mathbf{C}_1(0)) 
	\leq C Q_{\varsigma/2}(T,\widetilde{\varphi},\mathbf{A}_{1,1}(0)) ,
\end{equation}
where $C = C(n,m,q,p_0,\alpha,\gamma) \in (0,\infty)$ is a constant.

\noindent (4)   By \eqref{tilde varphi rmk eqn2} and \eqref{tilde varphi rmk eqn4}, $Q_{\varsigma/2}(T,\widetilde{\varphi},\mathbf{A}_{1,1}(0)) \leq C \widetilde{\beta} \inf_{x \in \mathbb{S}^1} \op{sep} \widetilde{\varphi}(x,0)$ for some constant $C = C(n,m,q,p_0,\alpha,\gamma) \in (0,\infty)$.  Thus assuming that $\widetilde{\beta}$ is sufficiently small, by \eqref{tilde varphi rmk eqn8} 
\begin{equation}\label{tilde varphi rmk eqn9}
	\varphi = \sum_{j=1}^{\widetilde{J}} \sum_{k=1}^{s_j} m_{j,k} \varphi_{j,k}  
\end{equation}
where $\varphi_{j,k}$ are distinct $q_{j,k}$-valued components of $\varphi$ (as in Definition~\ref{Phi-p_defn}) such that either 
\begin{enumerate}[itemsep=3mm,topsep=0mm]
	\item[(i)]  $\widetilde{\varphi}_j$ and $\varphi_{j,k}$ are both non-zero and 
	\begin{equation}\label{tilde varphi rmk eqn10}
		\sup_{B_1(0)} \mathcal{G}(\varphi_{j,k}, \widetilde{\varphi}_j) \leq C Q_{\varsigma/2}(T,\widetilde{\varphi},\mathbf{A}_{1,1}(0)) 
	\end{equation}
	for some constant $C = C(n,m,q,p_0,\alpha,\gamma) \in (0,\infty)$, 
	\item[(ii)]  $\widetilde{\varphi}_j$ and $\varphi_{j,k}$ are both identically zero, or 
	\item[(iii)]  $\widetilde{\varphi}_j$ is identically zero, $\varphi_{j,k}$ is non-zero, and 
	\begin{equation}\label{tilde varphi rmk eqn11}
		\sup_{B_1(0)} |\varphi_{j,k}| \leq C Q_{\varsigma/2}(T,\widetilde{\varphi},\mathbf{A}_{1,1}(0)) 
	\end{equation} 
	for some constant $C = C(n,m,q,p_0,\alpha,\gamma) \in (0,\infty)$ 
\end{enumerate}
and $s_j$ and $m_{j,k}$ are positive integers such that $\sum_{j=1}^{\widetilde{J}} \sum_{k=1}^{s_j} m_{j,k} q_{j,k} = q$.  By $p_0 > \widetilde{p}$, either there exists $j \in \{1,\ldots,\widetilde{J}\}$ and $k \in \{1,\ldots,s_j\}$ such that $\widetilde{\varphi}_j \equiv 0$ and $\varphi_{j,k}$ is a non-zero component of $\varphi$ which is close to $\widetilde{\varphi}_j$ and thus by \eqref{tilde varphi rmk eqn11}
\begin{equation*}%\label{tilde varphi rmk eqn12}
	\left( \int_{B_1(0)} |\varphi_{j,k}|^2 \right)^{1/2} \leq C Q_{\varsigma/2}(T, \widetilde{\varphi}, \mathbf{A}_{1,1}(0)) 
\end{equation*}
or there exists $j \in \{1,\ldots,\widetilde{J}\}$ and $k \neq k'$ such that $\widetilde{\varphi}_j$ is non-zero and $\varphi_{j,k}$ and $\varphi_{j,k'}$ are distinct non-zero components of $\varphi$ which are close to $\widetilde{\varphi}_j$ and thus by \eqref{tilde varphi rmk eqn10}
\begin{equation*}%\label{tilde varphi rmk eqn13}
	\left( \int_{B_1(0)} \mathcal{G}(\varphi_{j,k},\varphi_{j,k'})^2 \right)^{1/2} \leq C Q_{\varsigma/2}(T, \widetilde{\varphi}, \mathbf{A}_{1,1}(0)) .
\end{equation*}
}\end{remark}

\begin{proof}[Proof of Lemma~\ref{graphical_lemma}]
We proceed by induction on $p$.
The base case $p = 0$, in which $\varphi = q \llbracket 0 \rrbracket$, is a special case of~\cite[Lemma~4.6]{KrumWicb}.  
Suppose that $p_0 \in \{1,2,\ldots,\lfloor q/q_0 \rfloor\}$ such that (H3) holds true.

\noindent\textit{Proof of Lemma~\ref{graphical_lemma}(A).}  It suffices to suppose that either $q_0 = 1$ and $\varphi$ has two or more components or $q_0 \geq 2$, as otherwise \eqref{graphical annuli conclA1} and \eqref{graphical annuli conclA2} trivially hold true.  To see \eqref{graphical annuli conclA1}, suppose to the contrary that for $\nu = 1,2,3,\ldots$ there are $\varepsilon_{\nu} \rightarrow 0^+$, $\beta_{\nu} \rightarrow 0^+$, $\varphi^{(\nu)} \in \Phi_{\alpha,q,p_0}$, an $n$-dimensional locally area mininimizing rectifiable current $T^{(\nu)}$ of $\mathbf{A}_{1,1}(0)$ such that \eqref{graphical annuli hyp3} and \eqref{graphical annuli hyp4} hold true with $\varepsilon_{\nu}, \beta_{\nu}, \varphi^{(\nu)}, T^{(\nu)}$ in place of $\overline{\varepsilon}, \overline{\beta}, \varphi, T$ but 
\begin{equation}\label{graphical annuli eqn1}
	\inf_{x \in \mathbb{S}^1} \op{sep} \varphi^{(\nu)}(x,0) 
	\leq \frac{1}{\nu} \inf_{\varphi' \in \bigcup_{p'=1}^{p_0-1} \Phi_{\alpha,q,p'}} Q_{\varsigma/2}(T^{(\nu)},\varphi',\mathbf{A}_{1,1}(0)) .
\end{equation}
Let $\widetilde{\beta} = \widetilde{\beta}(n,m,q,p_0,\alpha,\gamma) \in (0,1)$ be a suitably small constant.  Let $\widetilde{p} \in \{0,1,2,\ldots,p_0-1\}$ and $\widetilde{\varphi}^{(\nu)} \in \Phi_{\alpha,q,\widetilde{p}}$ such that \eqref{tilde varphi rmk eqn1} and \eqref{tilde varphi rmk eqn2} hold true with $T^{(\nu)}$ and $\widetilde{\varphi}^{(\nu)}$ in place of $T$ and $\widetilde{\varphi}$.  Let 
\begin{equation*}
	\widetilde{\varphi}^{(\nu)} = \sum_{j=1}^{\widetilde{J}} \widetilde{m}_j \widetilde{\varphi}^{(\nu)}_j  
\end{equation*}
where $\widetilde{\varphi}^{(\nu)}_j$ are distinct $\widetilde{q}_j$-valued components of $\widetilde{\varphi}^{(\nu)}$ with positive integer multiplicity $\widetilde{m}_j$ (as in \eqref{tilde varphi rmk eqn3}).  After passing to a subsequence, we may assume that $\widetilde{J}$, $\widetilde{p}$, and $\widetilde{m}_j$ are independent of $\nu$, and that for each $j$ either $\widetilde{q}_j = 1$ and $\widetilde{\varphi}^{(\nu)}_j \equiv 0$ for all $\nu$ or $\widetilde{q}_j = q_0$ and $\widetilde{\varphi}^{(\nu)}_j$ is non-zero for all $\nu$.  Let $\widetilde{Q}_{\nu} = Q_{\varsigma/2}(T^{(\nu)},\widetilde{\varphi}^{(\nu)},\mathbf{A}_{1,1}(0))$.  By Remark~\ref{tilde varphi rmk}(4), 
\begin{equation*}
	\varphi^{(\nu)} = \sum_{j=1}^{\widetilde{J}} \sum_{k=1}^{s_j} m_{j,k} \varphi^{(\nu)}_{j,k}  
\end{equation*}
where $\varphi^{(\nu)}_{j,k}$ are distinct $q_{j,k}$-valued components of $\varphi^{(\nu)}$ which are uniformly $C \widetilde{Q}_{\nu}$-close to $\widetilde{\varphi}^{(\nu)}_j$ for some constant $C = C(n,m,q,p_0,\alpha,\gamma) \in (0,\infty)$ and where $s_j$ and $m_{j,k}$ are positive integers such that $\sum_{j=1}^{\widetilde{J}} \sum_{k=1}^{s_j} m_{j,k} q_{j,k} = q$ (as in~\eqref{tilde varphi rmk eqn9}).  After passing to a subsequence, assume that $s_j$, $m_{j,k}$, and $q_{j,k}$ are independent of $\nu$.  Set $\varphi^{(\nu)}_j = \sum_{k=1}^{s_j} m_{j,k} \varphi^{(\nu)}_{j,k}$ for each $j \in \{1,2,\ldots,\widetilde{J}\}$.  

Let $w = (w_j)$ be the blow-up of $T^{(\nu)}$ relative to $\widetilde{\varphi}^{(\nu)}$ by the excess $\widetilde{Q}_{\nu}$, and let $\psi = (\psi_j)$ be the blow-up of $\varphi^{(\nu)}$ relative to $\widetilde{\varphi}^{(\nu)}$ by the excess $\widetilde{Q}_{\nu}$ (as in Subsection~\ref{sec:blowup subsec}).  If $\widetilde{\varphi}^{(\nu)}_j \equiv 0$, let $\widetilde{\varphi}^{(\infty)}_j \equiv 0$.  If instead $\widetilde{\varphi}^{(\nu)}_j(x,y) = \op{Re}(\widetilde{c}^{(\nu)}_j (x_1+ix_2)^{\alpha})$ for some $\widetilde{c}^{(\nu)}_j \in \mathbb{C}^m \setminus \{0\}$, let $\widetilde{\varphi}^{(\nu)}_j/|\widetilde{c}^{(\nu)}_j| \rightarrow \widetilde{\varphi}^{(\infty)}_j$ uniformly on compact subsets of $\mathbb{R}^n$ as $\nu\rightarrow\infty$.  That way, $w_j$ is a function on $\op{graph} \widetilde{\varphi}^{(\infty)}_j |_{A_{1,1}(0)}$ and $\psi_j$ are functions on $\op{graph} \widetilde{\varphi}^{(\infty)}_j$ (as in Subsection~\ref{sec:blowup subsec}).  Let $\widetilde{v}^{(\nu)}_j : \op{graph} \widetilde{\varphi}^{(\nu)}_j |_{A_{1,7/8}(0)} \rightarrow \mathcal{A}_{\widetilde{m}_j}(\mathbb{R}^m)$, $\widetilde{u}^{(\nu)}_j : A_{1,7/8}(0) \rightarrow \mathcal{A}_{\widetilde{q}_j \widetilde{m}_j}(\mathbb{R}^m)$, and $\widetilde{K}^{(\nu)}_j \subseteq A_{1,7/8}(0)$ be a closed set such that \eqref{u defn} holds true with $\widetilde{\varphi}^{(\nu)}_j, \widetilde{v}^{(\nu)}_j, \widetilde{u}^{(\nu)}_j$ in place of $\varphi_j, v_j, u_j$ and \eqref{tilde varphi rmk eqn7} holds true with $1/2,\widetilde{\varphi}^{(\nu)}, T^{(\nu)}, \widetilde{v}^{(\nu)}_j, \widetilde{u}^{(\nu)}_j, \widetilde{K}^{(\nu)}_j$ in place of $\kappa,\widetilde{\varphi}, T, \widetilde{v}_j, \widetilde{u}_j, \widetilde{K}_j$.  It follows from \eqref{graphical annuli hyp4} using the area formula and \eqref{tilde varphi rmk eqn7} that 
\begin{align}\label{graphical annuli eqn2}
	\int_{A_{1,1/2}(0)} \sum_{j=1}^{\widetilde{J}} \sum_{l=1}^{\widetilde{q}_j} \sum_{h=1}^{\widetilde{m}_j} 
		\op{dist}^2(\widetilde{u}^{(\nu)}_{j,l,h}(x,y), \op{spt} \varphi^{(\nu)}_j(x,y)) \,dx \,dy \hspace{20mm}&\\
	\leq 4 \beta_{\nu}^2 \widetilde{Q}_{\nu}^2 + C \widehat{E}(T^{(\nu)},P_0,\mathbf{A}_{1,1}(0))^{2+\varsigma} & \nonumber 
\end{align}
where $\widetilde{u}^{(\nu)}_j(x,y) = \sum_{l=1}^{\widetilde{q}_j} \sum_{h=1}^{\widetilde{m}_j} \llbracket \widetilde{u}^{(\nu)}_{j,l,h}(x,y) \rrbracket$ for each $(x,y) \in A_{1,1/2}(0)$ and $C = C(n,m,q,p_0,\alpha,\gamma) \in (0,\infty)$ is a constant.  Dividing both sides of \eqref{graphical annuli eqn2} by $\widetilde{Q}_{\nu}^2$ and letting $\nu\rightarrow\infty$, we deduce that $\op{spt} w_j(X) \subseteq \op{spt}\psi_j(X)$ for all $X \in \op{graph} \widetilde{\varphi}^{(\infty)}_j |_{A_{1,1/2}(0)}$ and all $j \in \{1,2,\ldots,\widetilde{J}\}$.  By similar reasoning, $\op{spt} \psi_j(X) \subseteq \op{spt}w_j(X)$ for all $X \in \op{graph} \widetilde{\varphi}^{(\infty)}_j |_{A_{1,1/2}(0)}$ and all $j \in \{1,2,\ldots,\widetilde{J}\}$.  Therefore, $\op{spt} w_j(X) = \op{spt}\psi_j(X)$ for all $X \in \op{graph} \widetilde{\varphi}^{(\infty)}_j |_{A_{1,1/2}(0)}$ and all $j \in \{1,2,\ldots,\widetilde{J}\}$.

By \eqref{graphical annuli eqn1}, there exists $x_{\nu} \in \mathbb{S}^1$ such that $\op{sep} \varphi^{(\nu)}(x_{\nu},0) \leq \frac{1}{\nu} \,\widetilde{Q}_{\nu}$.  After passing to a subsequence, let $x_{\nu} \rightarrow x_0$ in $\mathbb{S}^1$ so that 
\begin{equation}\label{graphical annuli eqn3}
	\inf_{x \in \mathbb{S}^1 \cap B_{1/2}(x_0,0)} \op{sep} \varphi^{(\nu)}(x,0) \leq \frac{1}{\nu} \,\widetilde{Q}_{\nu} .
\end{equation}
Express 
\begin{equation*}
	\widetilde{\varphi}^{(\nu)}_j(x,y) = \sum_{l=1}^{\widetilde{q}_j} \llbracket \widetilde{\varphi}^{(\nu)}_{j,l}(x,y) \rrbracket , \quad 
	\varphi^{(\nu)}_{j,k}(x,y) = \sum_{l=1}^{q_{j,k}} \llbracket \varphi^{(\nu)}_{j,k,l}(x,y) \rrbracket
\end{equation*}
for each $(x,y) \in B_{1/2}(x_0,0)$, where $\widetilde{\varphi}^{(\nu)}_{j,l} : B_{1/2}(x_0,0) \rightarrow \mathbb{R}^m$ and $\varphi^{(\nu)}_{j,k,l} : B_{1/2}(x_0,0) \rightarrow \mathbb{R}^m$ are single-valued harmonic functions (as in \eqref{varphi_localized}).  Recall that by \eqref{tilde varphi rmk eqn2} and \eqref{tilde varphi rmk eqn4}, $\widetilde{Q}_{\nu} \leq C \widetilde{\beta} \inf_{x \in \mathbb{S}^1} \op{sep} \widetilde{\varphi}^{(\nu)}(x,0)$ for some constant $C = C(n,m,q,p_0,\alpha,\gamma) \in (0,\infty)$.  Thus by \eqref{tilde varphi rmk eqn10}, we may assume that if $\widetilde{\varphi}^{(\nu)}_j$ is non-zero for all $\nu$, then 
\begin{equation}\label{graphical annuli eqn4}
	\sup_{B_{1/2}(x_0,0)} |\varphi^{(\nu)}_{j,k,l} - \widetilde{\varphi}^{(\nu)}_{j,l}| \leq C \widetilde{Q}_{\nu}
\end{equation}
for all $k \in \{1,\ldots,s_j\}$ and $l \in \{1,\ldots,q_0\}$, where $C = C(n,m,q,p_0,\alpha,\gamma) \in (0,\infty)$ is a constant.  Hence assuming $\widetilde{\beta} = \widetilde{\beta}(n,m,q,p_0,\alpha,\gamma) \in (0,1)$ is sufficiently small, by $\widetilde{Q}_{\nu} \leq C \widetilde{\beta} \inf_{x \in \mathbb{S}^1} \op{sep} \widetilde{\varphi}^{(\nu)}(x,0)$, \eqref{graphical annuli eqn4}, and \eqref{tilde varphi rmk eqn11}
\begin{equation*}
	\inf_{B_{1/2}(x_0,0)} |\varphi^{(\nu)}_{j,k,l} - \varphi^{(\nu)}_{j',k',l'}| 
	\geq \tfrac{1}{2} \inf_{B_{1/2}(x_0,0)} |\widetilde{\varphi}^{(\nu)}_{j,l} - \widetilde{\varphi}^{(\nu)}_{j',l'}| 
	\geq c \widetilde{\beta}^{-1} \widetilde{Q}_{\nu}
\end{equation*}
whenever $(j,l) \neq (j',l')$, where $c = c(n,m,q,p_0,\alpha,\gamma) > 0$ is a constant.  Thus by \eqref{graphical annuli eqn3}, there exists $j_1 \in \{1,\ldots,\widetilde{J}\}$, $k \neq k'$, and $l_1 \in \{1,2,\ldots,\widetilde{q}_{j_1}\}$ such that 
\begin{equation}\label{graphical annuli eqn5}
	\inf_{x \in \mathbb{S}^1 \cap B_{1/2}(x_0)} |\varphi^{(\nu)}_{j_1,k,l_1}(x,0) - \varphi^{(\nu)}_{j_1,k',l_1}(x,0)| \leq \frac{1}{\nu} \,\widetilde{Q}_{\nu} .
\end{equation}
If $\widetilde{\varphi}^{(\nu)}_{j_1} \equiv 0$ note that $l_1 = 1$ and let $\widetilde{\varphi}^{(\infty)}_{j_1,1} \equiv 0$.  If instead $\widetilde{\varphi}^{(\nu)}_{j_1}(x,y) = \op{Re}(\widetilde{c}^{(\nu)}_{j_1} (x_1+ix_2)^{\alpha})$ for some $\widetilde{c}^{(\nu)}_{j_1} \in \mathbb{C}^m \setminus \{0\}$, let $\widetilde{\varphi}^{(\nu)}_{j_1,l_1}/|\widetilde{c}^{(\nu)}_{j_1}| \rightarrow \widetilde{\varphi}^{(\infty)}_{j_1,l_1}$ uniformly on compact subsets of $\mathbb{R}^n$ as $\nu\rightarrow\infty$ (as in Subsection~\ref{sec:blowup subsec}).  Upon blowing-up, \eqref{graphical annuli eqn5} gives us that the graph of $\psi_{j_1}(x,y,\widetilde{\varphi}^{(\infty)}_{j_1,l_1}(x,y))$ intersects itself along an $(n-1)$-dimensional half-plane.   Note that by Remark~\ref{graphical rmk}(3) (with $\mathbf{A}_{1,1}(0)$ in place of $\mathbf{C}_1(0)$), $\|\varphi^{(\nu)}_{j_1}\|_{L^2(B _1(0))} \geq c \widetilde{Q}_{\nu}$ if $\widetilde{\varphi}^{(\nu)}_{j_1} \equiv 0$ and $\|\mathcal{G}(\varphi^{(\nu)}_{j_1},\widetilde{\varphi}^{(\nu)}_{j_1})\|_{L^2(B _1(0))} \geq c \widetilde{Q}_{\nu}$ if $\widetilde{\varphi}^{(\nu)}_{j_1}$ is non-zero.  Thus $\psi_{j_1}$ is non-zero.  Since $\op{spt} w_{j_1}(X) = \op{spt}\psi_{j_1}(X)$ for all $X \in \op{graph} \widetilde{\varphi}^{(\infty)}_{j_1} |_{A_{1,1/2}(0)}$, $w_{j_1}(x,y,\widetilde{\varphi}^{(\infty)}_{j_1,l_1}(x,y))$ is a non-zero, homogeneous degree $\alpha$, Dirichlet energy minimizing multi-valued function on $B_{1/2}(x_0,0) \cap A_{1,1/2}(0)$ which is translation invariant along $\{0\} \times \mathbb{R}^{n-2}$ and whose singular set contains an $(n-1)$-dimensional half-plane, giving us a contradiction.  Therefore, \eqref{graphical annuli conclA1} holds true.

For the rest of the proof of Lemma~\ref{graphical_lemma}, let $\varepsilon,\beta \in (0,1)$ be suitably small constants depending only on $n,m,q,p_0,\alpha,\gamma,\kappa$ and let $\varphi \in \Phi_{\alpha,q,p_0}$ and $T$ be an $n$-dimensional locally area mininimizing rectifiable current of $\mathbf{A}_{1,1}(0)$ such that \eqref{graphical annuli hyp3} and \eqref{graphical annuli hyp4} hold true.  To see \eqref{graphical annuli conclA2}, let $\widetilde{\beta} = \widetilde{\beta}(n,m,q,p_0,\alpha,\gamma) \in (0,1)$ be a suitably small constant.  Let $\widetilde{p} \in \{0,1,2,\ldots,p_0-1\}$ and $\widetilde{\varphi} \in \Phi_{\alpha,q,\widetilde{p}}$ such that \eqref{tilde varphi rmk eqn1} and \eqref{tilde varphi rmk eqn2} hold true.  Let $\widetilde{\varphi}$ be as in \eqref{tilde varphi rmk eqn3} and $\varphi$ be as in \eqref{tilde varphi rmk eqn9}.  Let $x_0 \in \mathbb{S}^1$ and express 
\begin{equation*}
	\widetilde{\varphi}_j(x,y) = \sum_{l=1}^{q_j} \llbracket \widetilde{\varphi}_{j,l}(x,y) \rrbracket , \quad 
	\varphi_{j,k}(x,y) = \sum_{l=1}^{q_j} \llbracket \varphi_{j,k,l}(x,y) \rrbracket
\end{equation*}
for each $(x,y) \in B_{1/2}(x_0,0)$, where $\widetilde{\varphi}_{j,l} : B_{1/2}(x_0,0) \rightarrow \mathbb{R}^m$ and $\varphi_{j,k,l} : B_{1/2}(x_0,0) \rightarrow \mathbb{R}^m$ are single-valued harmonic functions (as in \eqref{varphi_localized}).  By \eqref{tilde varphi rmk eqn10}, we may assume that if $\widetilde{\varphi}_j$ is non-zero then 
\begin{equation}\label{graphical annuli eqn6}
	\|\varphi_{j,k,l} - \widetilde{\varphi}_{j,l}\|_{L^{\infty}(B_{1/2}(x_0,0))} \leq C Q_{\varsigma/2}(T, \widetilde{\varphi}, \mathbf{A}_{1,1}(0)) 
\end{equation}
for all $k \in \{1,\ldots,s_j\}$ and $l \in \{1,\ldots,\widetilde{q}_j\}$, where $C = C(n,m,q,p_0,\alpha,\gamma) \in (0,\infty)$ is a constant.  In the case $(j,l) = (j',l')$, by \eqref{graphical annuli eqn6}, \eqref{tilde varphi rmk eqn11}, \eqref{tilde varphi rmk eqn1}, and \eqref{graphical annuli conclA1}
\begin{align*}
	\|\varphi_{j,k,l} - \varphi_{j,k',l}\|_{L^2(B_{1/2}(x_0,0))} 
	\leq\,& C Q_{\varsigma/2}(T, \widetilde{\varphi}, \mathbf{A}_{1,1}(0)) 
	\leq C \inf_{x \in \mathbb{S}^1} \op{sep} \varphi(x,0)
	\\ \leq\,& C \inf_{x \in B_{1/2}(x_0,0)} |\varphi_{j,k,l} - \varphi_{j,k',l}|
\end{align*}
whenever $k \neq k'$.  Moreover, by \eqref{graphical annuli eqn6}, \eqref{tilde varphi rmk eqn2}, and \eqref{tilde varphi rmk eqn4} we have $\|\varphi_{j,k,l} - \widetilde{\varphi}_{j,l}\|_{L^{\infty}(B_{1/2}(x_0,0))} \leq C\widetilde{\beta} \inf_{x \in \mathbb{S}^1} \op{sep} \widetilde{\varphi}(x,0)$ if $\varphi_j$ is non-zero, and similarly using \eqref{tilde varphi rmk eqn11} we have $\|\varphi_{j,k,l}\|_{L^{\infty}(B_{1/2}(x_0,0))} \leq C\widetilde{\beta} \inf_{x \in \mathbb{S}^1} \op{sep} \widetilde{\varphi}(x,0)$ if $\varphi_j \equiv 0$.  This together with \eqref{tilde varphi rmk eqn5} gives us 
\begin{align*}
	\|\varphi_{j,k,l} - \varphi_{j',k',l'}\|_{L^2(B_{1/2}(x_0,0))} 
	\leq\,& 2 \|\widetilde{\varphi}_{j,l} - \widetilde{\varphi}_{j',l'}\|_{L^2(B_{1/2}(x_0,0))} 
	\\ \leq\,& 2C \inf_{x \in B_{1/2}(x_0,0)} |\widetilde{\varphi}_{j,l} - \widetilde{\varphi}_{j',l'}|
	\\ \leq\,& 4C \inf_{x \in B_{1/2}(x_0,0)} |\varphi_{j,k,l} - \varphi_{j',k',l'}|
\end{align*}
whenever $(j,l) \neq (j',l')$.  

\noindent\textit{Proof of Lemma~\ref{graphical_lemma}(B)(C).}  Note that by Remark~\ref{graphical rmk}(2) (with $\mathbf{A}_{1,1}(0)$ in place of $\mathbf{C}_1(0)$), $\|\varphi\|_{L^2(B_1(0))} \leq C \widehat{E}(T,P_0,\mathbf{A}_{1,1}(0))$ for some constant $C = C(n,m,q,\alpha,\gamma) \in (0,\infty)$.  
By \eqref{graphical annuli hyp4} and \eqref{graphical annuli conclA1}, $\inf_{x \in \mathbb{S}^1} \op{sep}\varphi(x,0) \geq c \beta^{-1} Q_{\varsigma/2}(T,\varphi,\mathbf{A}_{1,1}(0))$ for some constant $c = c(n,m,q,p,\alpha,\gamma) \in (0,\infty)$.  
Thus by \eqref{graphical annuli hyp3} and \eqref{graphical annuli conclA2}, we can apply Lemma~\ref{graphical_lemma0} to deduce that there exists integers $\widehat{m}_j \geq 0$ with $\sum_{j=1}^J q_j \widehat{m}_j = q$ and (unique) $n$-dimensional locally area-minimizing rectifiable currents $T_j$ of $\mathbf{A}_{1,(1+\kappa)/2}(0)$ such that \eqref{graphical annuli conclB1}, \eqref{graphical annuli conclB3}, and \eqref{graphical annuli conclB4} hold true.  Moreover, whenever $\widehat{m}_j > 0$ there exists a Lipschitz $\widehat{m}_j$-valued functions $v_j : \op{graph} \varphi_j |_{A_{1,\kappa}(0)} \rightarrow \mathcal{A}_{\widehat{m}_j}(\mathbb{R}^m)$ and a closed set $K_j \subset A_{1,\kappa}(0)$ such that \eqref{graphical annuli conclC1} and \eqref{graphical annuli conclC3} hold true.

Suppose that $\widehat{m}_j = 0$ for some $j \in \{1,2,\ldots,J\}$.  Then by~\cite[Lemma~3.4]{KrumWicb}, $\op{spt} T_j \cap \mathbf{A}_{1,(3+5\kappa)/8}(0) = \emptyset$.  We claim that for each $X \in \op{graph} \varphi_j |_{A_{1,1/4}(0)}$ 
\begin{equation}\label{graphical annuli eqn7}
	\op{dist}(X, \op{spt} T) \geq c \inf_{\varphi' \in \bigcup_{p'=1}^{p_0-1} \Phi_{\alpha,q,p'}} Q_{\varsigma/2}(T,\varphi',\mathbf{A}_{1,1}(0)) , 
\end{equation}
for some constant $c = c(n,m,q,p_0,\alpha,\gamma,\kappa) > 0$.  To see \eqref{graphical annuli eqn7}, we may assume that $\op{dist}(X, \op{spt} T) < (1-\gamma)/64$ so that the closest point $Y$ to $X$ on $\op{spt} T$ is in $\mathbf{A}_{1,3/8}(0)$.  Since $\op{spt} T_j \cap \mathbf{A}_{1,3/8}(0) = \emptyset$, the closest point to $Y$ on the support of $\op{graph}\varphi$ must lie on $\op{graph}\varphi_{j'}$ for some $j' \neq j$.  By \eqref{graphical annuli conclB2} and \eqref{graphical annuli hyp4}, the distance of $Y$ to $\op{spt}\op{graph}\varphi_{j'}$ is $\leq C\overline{\beta} \inf_{\varphi' \in \bigcup_{p'=1}^{p_0-1} \Phi_{\alpha,q,p'}} Q_{\varsigma/2}(T,\varphi',\mathbf{A}_{1,1}(0))$ for some constant $C = C(n,m,q,p_0,\alpha,\gamma,\kappa) \in (0,\infty)$.  However, by \eqref{graphical annuli conclA1}, the distance of $X$ to $\op{graph} \varphi_{j'} |_{A_{1,1/4}(0)}$ is $\geq c \inf_{\varphi' \in \bigcup_{p'=1}^{p_0-1} \Phi_{\alpha,q,p'}} Q_{\varsigma/2}(T,\varphi',\mathbf{A}_{1,1}(0))$ for some constant $c = c(n,m,q,p_0,\alpha,\gamma,\kappa) > 0$, and thus \eqref{graphical annuli eqn7} must hold true.  
Integrating \eqref{graphical annuli eqn7} over $X \in \op{graph} \varphi_j |_{A_{1,1/4}(0)}$, we derive a contradiction with \eqref{graphical annuli hyp4}.  Therefore, $\widehat{m}_j > 0$ for all $j \in \{1,2,\ldots,J\}$.  Up to changing the multiplicities $m_j$, we may take $\widehat{m}_j = m_j$ for all $j$.
\end{proof}

\begin{proof}[Proof of Theorem~\ref{graphical_thm}] 
Let $0 < \tau < \gamma < 1$.  Let $\varepsilon_0, \beta_0 \in (0,1)$, $p \in \{1,2,\ldots,\lfloor q/q_0 \rfloor \}$, $\varphi \in \Phi_{\alpha,q,p}$, and $T$ be an $n$-dimensional rectifiable current of $\mathbf{C}_1(0)$ such that \eqref{graphical hyp1} and \eqref{graphical hyp2} hold true.  Let $\rho \geq \tau/4$ and $\zeta \in \mathbb{R}^{n-2}$ such that $\rho^2 + |\zeta|^2 \leq (3+\gamma)^2/16$.  It follows immediately from \eqref{graphical hyp1} that 
\begin{gather}\label{graphical eqn1} 
	(\partial T) \llcorner \mathbf{A}_{\rho,1}(\zeta) = 0 , \quad 
	\sup_{X \in \op{spt} T} \op{dist}(X,P_0) < \infty, \quad 
	\pi_{\#} (T \llcorner \mathbf{A}_{\rho,1}(\zeta)) = q \llbracket A_{\rho,1}(\zeta) \rrbracket , \\ 
	\widehat{E}(T,P_0,\mathbf{A}_{\rho,1}(\zeta)) \leq \omega_n^{1/2} \rho^{-(n+2)/2} \widehat{E}(T,P_0,\mathbf{C}_1(0)) 
		< \omega_n^{1/2} (\tau/4)^{-(n+2)/2} \varepsilon_0 . \nonumber 
\end{gather}
By similar reasoning using~\cite[Lemma~3.5 and Lemma~3.7]{KrumWica}, $\rho^{-n} \int_{\mathbf{A}_{\rho,1}(\zeta)} |\vec T - \vec P_0|^2 \,d\|T\| \leq C\tau^{-(n+2)/2} \widehat{E}(T,P_0,\mathbf{C}_1(0)) < C \tau^{-(n+2)/2} \varepsilon_0$ and thus 
\begin{equation}\label{graphical eqn2} 
	\|T\|(\mathbf{A}_{\rho,1}(\zeta)) \leq (q+1/2) \,\mathcal{L}^n(A_{\rho,1}(0)) .
\end{equation}

Let $\widetilde{\beta} = \beta(n,m,q,p_0,\alpha) \in (0,1)$ to be later determined.  By Remark~\ref{tilde varphi rmk}(1)(4), either there exists a non-zero component $\varphi_j$ of $\varphi$ such that 
\begin{equation}\label{graphical eqn3} 
	\left(\int_{B_1(0)} |\varphi_j|^2\right)^{1/2} \leq C \inf_{\varphi' \in \bigcup_{p'=0}^{p_0-1} \Phi_{\alpha,q,p'}} Q_{\varsigma/2}(T,\varphi',\mathbf{A}_{1/4,1}(0)) 
\end{equation}
or there exist two distinct non-zero components $\varphi_j$ and $\varphi_{j'}$ of $\varphi$ such that 
\begin{equation}\label{graphical eqn4} 
	\left(\int_{B_1(0)} \mathcal{G}(\varphi_j,\varphi_{j'})^2\right)^{1/2} 
		\leq C \inf_{\varphi' \in \bigcup_{p'=0}^{p_0-1} \Phi_{\alpha,q,p'}} Q_{\varsigma/2}(T,\varphi',\mathbf{A}_{1/4,1}(0)) ,
\end{equation}
where $C = C(n,m,q,\alpha,\gamma) \in (0,\infty)$ is a constant.  Let's assume that \eqref{graphical eqn4} holds true, as the case where \eqref{graphical eqn3} holds true is similar.  By $A_{1/4,1}(0) \subset B_1(0)$, $A_{1/4,1/2}(0) \subset B_{1/2}(0) \cap \{r > 1/16\}$, the homogeneity of $\varphi$, \eqref{graphical hyp2}, and \eqref{graphical rmk eqn7}
\begin{equation}\label{graphical eqn5} 
	Q_{\varsigma/2}(T,\varphi,\mathbf{A}_{1/4,1}(0)) \leq Q_{\varsigma/2}(T,\varphi,\mathbf{C}_1(0)) 
	\leq C\beta_0 \left(\int_{B_1(0)} \mathcal{G}(\varphi_j,\varphi_{j'})^2\right)^{1/2} , 
\end{equation}
where $C = C(n,m,q,\alpha,\gamma) \in (0,\infty)$ is a constant.  Combining \eqref{graphical eqn4} and \eqref{graphical eqn5}
\begin{equation}\label{graphical eqn6} 
	Q_{\varsigma/2}(T,\varphi,\mathbf{A}_{1/4,1}(0)) 
	\leq C\beta_0 \inf_{\varphi' \in \bigcup_{p'=0}^{p_0-1} \Phi_{\alpha,q,p'}} Q_{\varsigma/2}(T,\varphi',\mathbf{A}_{1/4,1}(0)) 
\end{equation}
for some constant $C = C(n,m,q,\alpha,\gamma) \in (0,\infty)$.  By \eqref{graphical eqn4} and \eqref{graphical rmk eqn7}
\begin{equation}\label{graphical eqn7} 
	\inf_{\varphi' \in \bigcup_{p'=0}^{p_0-1} \Phi_{\alpha,q,p'}} Q_{\varsigma/2}(T,\varphi',\mathbf{C}_1(0)) 
		\leq C \inf_{\varphi' \in \bigcup_{p'=0}^{p_0-1} \Phi_{\alpha,q,p'}} Q_{\varsigma/2}(T,\varphi',\mathbf{A}_{1/4,1}(0)) ,
\end{equation}
where $C = C(n,m,q,\alpha,\gamma) \in (0,\infty)$ is a constant.  

Now in light of \eqref{graphical eqn1}, \eqref{graphical eqn2}, \eqref{graphical eqn6}, and \eqref{graphical eqn7}, we can apply Lemma~\ref{graphical_lemma}(A) with $\eta_{0,1/4\#} T$ and $4^{1-\alpha} \varphi$ in place of $T$ and $\varphi$ to deduce that Theorem~\ref{graphical_thm}(A) holds true.  Moreover, by the Vitali covering lemma there is a finite collection of $\rho_k \geq \tau/4$ and $\zeta_k \in \mathbb{R}^{n-2}$ such that $\rho_k^2 +|\zeta_k|^2 \leq (3+\gamma)^2/16$, $\rho_{k+1} \leq \rho_k$, $\{A_{\rho_k,1/8}(\rho_k,\zeta_k)\}$ covers $B_{(3+\gamma)/4}(0) \cap \{r > \tau/4\}$, and $\{A_{\rho_k,1/40}(\rho_k,\zeta_k)\}$ is pairwise disjoint.  By \eqref{graphical eqn1}, \eqref{graphical conclA2}, and Remark~\ref{graphical rmk}(2), we can apply Lemma~\ref{graphical_lemma0} with $\eta_{(0,\zeta_k),\rho_k\#} T$ and $\rho_k^{\alpha-1} \varphi$ in place of $T$ and $\varphi$ to deduce that there exists integers $\widehat{m}_{k,j} \geq 0$ with $\sum_{j=1}^J \widehat{m}_{k,j} q_j = q$ and $n$-dimensional locally area-minimizing rectifiable currents $T_{k,j}$ of $\mathbf{A}_{\rho_k,3/4}(\zeta_k)$ such that 
\begin{gather}\label{graphical eqn8}
	T \llcorner \mathbf{A}_{\rho_k,3/4}(\zeta_k) = \sum_{j=1}^J T_{k,j} , \quad 
	(\partial T_{k,j}) \llcorner \mathbf{A}_{\rho_k,3/4}(\zeta_k) = 0, \quad 
	\pi_{\#} T_{k,j} = \widehat{m}_{k,j} q_j \llbracket \mathbf{A}_{\rho_k,3/4}(\zeta_k) \rrbracket , \\	
	\sup_{X \in \op{spt} T_{k,j}} \op{dist}(X, \op{spt} \op{graph} \varphi_j) \leq C \rho_k^{-(n+2)(4+\varsigma)/8} E_{\varsigma}(T,\varphi,\mathbf{C}_1(0)) , \nonumber
\end{gather}
and when $\widehat{m}_j > 0$ there exists $\widehat{m}_j$-valued functions $v_{k,j} : \op{graph} \varphi_j |_{A_{\rho_k,\kappa}(\zeta_k)} \rightarrow \mathcal{A}_{\widehat{m}_{k,j}}(\mathbb{R}^m)$ and closed sets $K_{k,j} \subseteq A_{\rho_k,\kappa}(\zeta_k)$ such that 
\begin{gather*}%\label{graphical eqn9}
	T_{k,j} \llcorner (K_{k,j} \times \mathbb{R}^m) = (\op{graph} u_{k,j}) \llcorner (K_{k,j} \times \mathbb{R}^m) , \\
	\mathcal{H}^n(A_{\rho_k,1/2}(\zeta_k) \setminus K_{k,j}) + \|T_j\|((A_{\rho_k,1/2}(\zeta_k) \setminus K_{k,j}) \times \mathbb{R}^m) \leq C \rho_k^{-2-\varsigma(n+2)/2} \widehat{E}(T,P_0,\mathbf{C}_1(0))^{2+\varsigma}, \nonumber\\
	\sup_{A_{\rho_k,1/2}(\zeta_k)} |u_{k,j}| \leq C \rho_k^{-n/2} \widehat{E}(T,P_0,\mathbf{C}_1(0)), \quad 
	\sup_{A_{\rho_k,1/2}(\zeta_k)} |Du_{k,j}| \leq C \rho_k^{-\varsigma (n+2)/2} \widehat{E}(T,P_0,\mathbf{C}_1(0))^{\varsigma} , \nonumber \\
	\frac{1}{\rho_k^2} \sup_{A_{\rho_k,1/2}(\zeta_k)} \sum_{j=1}^J \sum_{l=1}^{q_j} |v_{k,j,l}|^2 
		+ \frac{1}{\rho_k^n} \int_{A_{\rho_k,1/2}(\zeta_k)} \sum_{j=1}^J \sum_{l=1}^{q_j} |Dv_{k,j,l}|^2 
		\leq C \rho_k^{-(n+2)(2+\varsigma)/2} E_{\varsigma/2}(T,\varphi,\mathbf{C}_1(0))^2 , \nonumber
\end{gather*}
where $v_{k,j,l}$ and $u_{k,j}$ are as in \eqref{v_localized} and \eqref{u defn} with $v_{k,j}, u_{k,j}$ in place of $v_j, u_j$ and $C = C(n,m,q,\alpha,\gamma) \in (0,\infty)$ are constants.  By \eqref{graphical eqn8}, if $\mathbf{A}_{\rho_k,3/4}(\zeta_k) \cap \mathbf{A}_{\rho_{k'},3/4}(\zeta_{k'}) \neq \emptyset$ then $T_{k,j} = T_{k',j}$ on $\mathbf{A}_{\rho_k,3/4}(\zeta_k) \cap \mathbf{A}_{\rho_{k'},3/4}(\zeta_{k'})$ and $\widehat{m}_{k,j} = \widehat{m}_{k',j}$.  Hence for each $j \in \{1,2,\ldots,J\}$, $\widehat{m}_j = \widehat{m}_{k,j}$ independent of $k$ and there exists an $n$-dimensional locally area-minimizing rectifiable current $T_j$ of $\mathbf{B}_{(1+\gamma)/2}(0) \cap \{r > \tau/2\}$ such that $T_j = T_{k,j}$ on $\mathbf{A}_{\rho_k,3/4}(\zeta_k) \cap \mathbf{B}_{(1+\gamma)/2}(0) \cap \{r > \tau/2\}$.  In particular, $T_j$ satisfies Theorem~\ref{graphical_thm}(B).  
By Lemma~\ref{graphical_lemma}(B) with $\eta_{0,1/4\#} T$ and $4^{1-\alpha} \varphi$ in place of $T$ and $\varphi$, $\widehat{m}_j > 0$ for all $j$ and so up to changing the multiplicities $m_j$ we may take $\widehat{m}_j = m_j$ for all $j$.  Using a partition of unity argument (as in the proof of~\cite[Theorem~4.4]{KrumWicb}), we can glue together the multi-valued functions $v_{k,j}$ and the sets $K_{k,j}$ to form $m_j$-valued functions $v_j : \op{graph} \varphi_j |_{B_{\gamma}(0) \cap \{r > \tau\}} \rightarrow \mathcal{A}_{m_j}(\mathbb{R}^m)$ and closed sets $K_k \subseteq B_{\gamma}(0) \cap \{r > \tau\}$ satisfying Theorem~\ref{graphical_thm}(C).
\end{proof}

\begin{remark}\label{graphical triangle rmk}{\rm 
(1)  \textit{Given $\rho \in (0,1)$ there exists $\varepsilon_0 = \varepsilon_0(n,m,q,\alpha,\rho) \in (0,1)$ and $\beta_0 = \beta_0(n,m,q,\alpha,\rho) \in (0,1)$ such that the following holds true.  Let $p \in \{1,2,\ldots,\lfloor q/q_0 \rfloor \}$, $\varphi \in \Phi_{\alpha,q,p}$, and $T$ be an $n$-dimensional rectifiable current of $\mathbf{C}_1(0)$ such that \eqref{graphical hyp1} and \eqref{graphical hyp2} hold true, and let $\varphi^* \in \Phi_{\alpha,q}$.  Let $\varphi = \sum_{j=1}^J m_j \varphi_j$ with $q_j$-valued components $\varphi_j$ with positive integer multiplicity $m_j$, and let $\varphi^* = \sum_{j=1}^{J^*} m^*_j \varphi^*_j$ with $q^*_j$-valued components $\varphi^*_j$ with positive integer multiplicity $m^*_j$ (as in Definition~\ref{Phi-p_defn}).  Then 
\begin{align}\label{graphical triangle rmk eqn1}
	&\op{dist}_{\mathcal H}(\op{spt} \op{graph} \varphi^* \cap \mathbf{C}_1(0), \op{spt} \op{graph} \varphi \cap \mathbf{C}_1(0)) 
	\\&\hspace{15mm} \leq C Q_{\varsigma/2}(T,\varphi^*,\mathbf{C}_{\rho}(0)) + C E_{\varsigma/2}(T,\varphi,\mathbf{C}_1(0)) , \nonumber 
\end{align}
where $C = C(n,m,q,\alpha,\rho) \in (0,\infty)$ is a constant.}

To see this, by the triangle inequality and Theorem~\ref{graphical_thm}(B) 
\begin{align}\label{graphical triangle rmk eqn2}
	&\frac{2}{\rho^{n+2}} \sum_{j=1}^J m_j \int_{B_{\rho/4}(0) \cap \{r > \rho/8\}} \sum_{l=1}^{q_j} 
		\op{dist}^2((x,y,\varphi^*_{j,l}(x,y)), \op{spt}\op{graph}\varphi) \,dx \,dy 
	\\ \leq\,& \frac{2}{\rho^{n+2}} \sum_{j=1}^J m_j \int_{B_{\rho/4}(0) \cap \{r > \rho/8\}} \sum_{l=1}^{q_j} 
		\op{dist}^2((x,y,\varphi^*_{j,l}(x,y)), \op{spt} T) \,dx \,dy \nonumber 
		\\&+ C \sup_{X \in \op{spt} T \cap B_{\rho/2}(0) \cap \{r > \rho/16\}} \op{dist}^2(X, \op{spt}\op{graph}\varphi) \nonumber
	\\ \leq\,& 2 Q_{\varsigma/2}(T, \varphi^*, \mathbf{C}_{\rho}(0))^2 + C E_{\varsigma/2}(T,\varphi,\mathbf{C}_1(0))^2, \nonumber
\end{align}
where $\varphi^*_j(x,y) = \sum_{l=1}^{q^*_j} \llbracket \varphi^*_{j,l}(x,y) \rrbracket$ for each $(x,y) \in \mathbb{R}^n$ (as in \eqref{varphi_localized}) and $C = C(n,m,q,\alpha,\rho) \in (0,\infty)$ are constants.  Let $K_j$, $u_j$, and $v_j$ be as in Theorem~\ref{graphical_thm}(C) with $\gamma = \rho/2$ and $\tau = \rho/16$.  By the triangle inequality and Theorem~\ref{graphical_thm}(C) 
\begin{align}\label{graphical triangle rmk eqn3}
	&\frac{2}{\rho^{n+2}} \sum_{j=1}^J m_j \int_{B_{\rho/4}(0) \cap \{r > \rho/8\}} 
		\sum_{l=1}^{q_j} \op{dist}^2((x,y,\varphi_{j,l}(x,y)), \op{spt}\op{graph}\varphi^*) \,dx \,dy 
	\\ \leq\,& \frac{2}{\rho^{n+2}} \int_{B_{\rho/2}(0) \cap \{r > \rho/16\}} \sum_{j=1}^J \sum_{l=1}^{q_j} \sum_{h=1}^{m_j} 
		\op{dist}^2((x,y,u_{j,l,h}(x,y)), \op{spt}\op{graph}\varphi^*) \,dx \,dy \nonumber 
		\\&+ \frac{2}{\rho^{n+2}} \int_{B_{\rho/2}(0) \cap \{r > \rho/16\}} \sum_{j=1}^J \sum_{l=1}^{q_j}  |v_{j,l}(x,y)|^2 \,dx \,dy \nonumber
	\\ \leq\,& 2 Q_{\varsigma/2}(T, \varphi^*, \mathbf{C}_{\rho}(0))^2 + 8 E_{\varsigma/2}(T, \varphi, \mathbf{C}_1(0))^2 
		+ C \widehat{E}(T, P_0, \mathbf{C}_1(0))^{2+\varsigma} , \nonumber
\end{align}
where $u_{j,l,h}$ and $v_{j,l}$ are as in \eqref{u defn} and \eqref{v_localized} and $C = C(n,m,q,\alpha,\rho) \in (0,\infty)$ is a constant.  Combining \eqref{graphical triangle rmk eqn2} and \eqref{graphical triangle rmk eqn3} and noting that $\varphi,\varphi^* \in \Phi_{\alpha,q}$ gives us \eqref{graphical triangle rmk eqn1}.

\noindent (2)  By arguing as in \eqref{graphical triangle rmk eqn3} with $\varphi^* = q \llbracket 0 \rrbracket$ and assuming for $\varepsilon_0$ small that $C  \widehat{E}(T, P_0, \mathbf{C}_1(0))^{\varsigma} < 1$, we obtain:}  Given $\rho \in (0,1)$ there exists $\varepsilon_0 = \varepsilon_0(n,m,q,\alpha,\rho) \in (0,1)$ and $\beta_0 = \beta_0(n,m,q,\alpha,\rho) \in (0,1)$ such that if $p \in \{1,2,\ldots,\lfloor q/q_0 \rfloor \}$, $\varphi \in \Phi_{\alpha,q,p}$, and $T$ be an $n$-dimensional rectifiable current of $\mathbf{C}_1(0)$ such that \eqref{graphical hyp1} and \eqref{graphical hyp2} hold true, then 
\begin{equation*}
	\rho^{\alpha-1} \|\varphi\|_{L^{\infty}(B_1(0))} \leq C \widehat{E}(T,P_0,\mathbf{C}_{\rho}(0)) + C E_{\varsigma/2}(T,\varphi,\mathbf{C}_1(0)) , 
\end{equation*}
where $C = C(n,m,q,\alpha) \in (0,\infty)$ is a constant (independent of $\rho$).
\end{remark}

\section{A priori estimates: Part~I}\label{sec:key_est_sec}

Let $\varphi : \mathbb{R}^n \rightarrow \mathcal{A}_q(\mathbb{R}^m)$ be a non-zero, homogeneous degree $\alpha$, harmonic function which is translation invariant along $\{0\} \times \mathbb{R}^{n-2}$ (as in Definition~\ref{Phi-p_defn}).  Let $T$ be an $n$-dimensional locally area-minimizing rectifiable current of $\mathbf{C}_1(0)$ such that $(\partial T) \llcorner \mathbf{C}_1(0) = 0$ and $T$ is weakly close to $\mathbb{R}^n \times \{0\}$ with integer multiplicity $q$.  In this section and Section~\ref{sec:nonconcentration}, we establish key integral estimates for $T$ which are inspired by the results of~\cite{Sim93} and \cite{Wic14}.  These estimates will play a fundamental role in the proof of the main excess decay estimate, Theorem~\ref{main excess decay thm}.

The main results of this section are Theorem~\ref{keyest thm} and Theorem~\ref{radial noncon thm}.  Theorem~\ref{keyest thm} will be used to obtain Theorem~\ref{radial noncon thm} as well as further estimates in Section~\ref{sec:nonconcentration}.  Additionally, Theorem~\ref{keyest thm} together with the estimates in Section~\ref{sec:nonconcentration} directly imply estimates for blow-ups (as constructed in Subsection~\ref{sec:blowup subsec}) which form the basis for the analysis of blow-ups in Section~\ref{sec:fine blowups}.  Theorem~\ref{radial noncon thm} as well as the estimates in the subsequent section (Theorem~\ref{nonconcentration thm} and Corollary~\ref{nonconcentration cor}) rule out the concentration of $L^2$-excess of $T$ relative to $\varphi$ near branch points of planar frequency $\geq \alpha$.  This in turn guarantees convergence of blow-up sequences in $L^2_{\rm loc}(B_{1/2}(0))$, which is essential to obtaining excess improvement as in option~(B) of Theorem~\ref{main excess decay thm}.

The proof of Theorem~\ref{keyest thm} is based on the monotonicity formula for the planar frequency function~\cite[Theorem~4.4]{KrumWica} and adaptations of arguments from~\cite[Theorem~8.2]{KrumWic1} and~\cite[Theorem~7.2]{KrumWic2}, as well as arguments inspired by~\cite[Theorem~10.1]{Wic14} to obtain estimates for $T$ over the (annular) cylinders $A_{\rho,1}(\zeta) \times \mathbb{R}^m$.  The proof of Theorem~\ref{radial noncon thm} follows a standard argument based on~\cite[Lemma~3.4]{Sim93} (also see~\cite[Corollary~8.3]{KrumWic1} and~\cite[Lemma~8.1]{KrumWic2}).

Throughout this section and Sections~\ref{sec:nonconcentration}--\ref{sec:fine excess decay}, we shall assume the first or both of the following hypotheses for appropriate choices of small constants $\varepsilon_0 \in (0,1)$ and $\beta_0 \in (0,1)$:

\noindent\textbf{Hypothesis~$(\star)$.}  $q \geq 2$ is an integer.  $\alpha = \ell_0/q_0 \neq 2$ where $\ell_0$ and $q_0$ are relatively prime positive integers with $q_0 \leq q$.  $p \in \{1,2,\ldots,\lfloor q/q_0 \rfloor\}$, $\varphi \in \Phi_{\alpha,q,p}$, and $T$ is an $n$-dimensional locally area-minimizing rectifiable current of $\mathbf{C}_1(0)$ such that 
\begin{gather}\label{main hyp eqn1} 
	(\partial T) \llcorner \mathbf{C}_1(0) = 0, \quad\quad \sup_{X \in \op{spt} T} \op{dist}(X,P_0) < \varepsilon_0, \\ 
	\|T\|(\mathbf{C}_1(0)) \leq (q+1/2) \,\omega_n , \quad\quad \pi_{\#} T = q \llbracket B_1(0) \rrbracket . \nonumber 
\end{gather}

\noindent\textbf{Hypothesis~$(\star\star)$.}  $q \geq 2$ is an integer.  $\alpha = \ell_0/q_0 \neq 2$ where $\ell_0$ and $q_0$ are relatively prime positive integers with $q_0 \leq q$.  $p \in \{1,\ldots,\lfloor q/q_0 \rfloor\}$, $\varphi \in \Phi_{\alpha,q,p}$, and $T$ is an $n$-dimensional  locally area-minimizing rectifiable current of $\mathbf{C}_1(0)$ such that 
\begin{equation}\label{main hyp eqn2} 
	Q_{\varsigma/2}(T, \varphi, \mathbf{C}_1(0)) \leq \beta_0 \inf_{\varphi' \in \bigcup_{p'=0}^{p-1} \Phi_{\alpha,q,p'}} Q_{\varsigma/2}(T, \varphi', \mathbf{C}_1(0)) , 
\end{equation}
where $\varsigma = \varsigma(n,m,q) \in (0,1)$ is as in Theorem~\ref{relharm separation3 thm} (i.e.~$\varsigma$ is as in Lemma~\ref{relharm tilt growth lemma}).

\begin{remark}\label{planar freq rmk}{\rm 
Let $\varepsilon_0 = \varepsilon_0(n,m,q) \in (0,1)$ be suitably small and let $T$ and $\alpha$ be as in Hypothesis~$(\star)$.  Let $\mu \in (0,1/q]$ and $\mathcal{B}_q$ be as in Definition~\ref{introduction B defn}.  In Sections~\ref{sec:key_est_sec}--\ref{sec:proof of main results}, we will be interested in points $Z \in \mathcal{B}_q$ such that $\mathcal{N}_{T,{\rm Pl}}(Z) \geq \alpha$.  By Lemma~\ref{uniform branch set lemma} and Theorem~\ref{frequency monotonicity thm}(e), such points $Z$ can be characterized as the set of all branch points $Z$ of $T$ such that $\Theta(T,Z) = q$, $T$ has a unique tangent plane $P_Z$ at $Z$, and 
\begin{equation*} 
	\sup_{\rho \in (0,1/8]} \rho^{-n-2\alpha} \int_{\mathbf{B}_{2\rho}(Z)} \op{dist}^2(X,Z+P_Z) \,d\|T\|(X) < \infty .  
\end{equation*}
In particular, the set $\{ Z \in \mathcal{B}_q : \mathcal{N}_{T,{\rm Pl}}(Z) \geq \alpha \}$ is independent of the choice of parameter $\mu$ as in Hypothesis~\ref{introduction hyp-1} and Definition~\ref{introduction B defn}.  By noting that $\alpha \geq 1+1/q$ and applying Lemma~\ref{uniform branch set lemma} (with $\mu = 1/q$), \eqref{uniform branch set concl1}--\eqref{uniform branch set concl5} hold true for some constant $C = C(n,m,q) \in (0,\infty)$.
}\end{remark}

Now recall that given $\varepsilon_0 \in (0,1)$ suitably small, $T$ as in Hypothesis~$(\star)$, and $\mathcal{B}_q$ as in Definition~\ref{introduction B defn}, for each $\rho \in (0,1]$ we define $H_{T,P_0,0}(\rho)$ and $D_{T,P_0,0}(\rho)$ as in Definition~\ref{H D N defn} (with $Z = 0$ and $P = P_0 = \mathbb{R}^n \times \{0\}$).  For each $Z \in \mathcal{B}_q$ and $\rho \in (0,1/2]$, we let $P_Z$ denote the unique tangent plane of $T$ at $Z$ and define $H_{T,P_Z,Z}(\rho)$, $D_{T,P_Z,Z}(\rho)$, $N_{T,P_Z,Z}(\rho)$, and $\mathcal{N}_{T,{\rm Pl}}(Z)$ as in Definition~\ref{H D N defn}.  Recall from Subsection~\ref{sec:graphical_rep_notation} that we represent a point $X \in \mathbb{R}^{n+m}$ by $X = (x,y,z)$ for $x \in \mathbb{R}^2$, $y \in \mathbb{R}^{n-2}$, and $z = (z^1,z^2,\ldots,z^m) \in \mathbb{R}^m$. 

\begin{lemma}\label{weiss lemma}
For each $\gamma \in (0,1)$ there exists $\varepsilon_0 = \varepsilon_0(n,m,q,\gamma) \in (0,1)$ such that the following holds true.  Let $T$ be an $n$-dimensional locally area-minimizing rectifiable current of $\mathbf{C}_1(0)$ satisfying \eqref{main hyp eqn1}.  Let $\mathcal{B}_q$ be as in Definition~\ref{introduction B defn} (with $\mu = 1/q$) and assume $0 \in \mathcal{B}_q$.  Then for each $\alpha \in \mathbb{R}$ 
\begin{align}\label{weiss concl}
	&\frac{d}{d\rho} \Big( \rho^{-2\alpha} (D_{T,P_0,0}(\rho) - \alpha \,H_{T,P_0,0}(\rho)) \Big) 
	\\ \geq\,& -2\rho^{-n-2\alpha} \int_{G_n(\mathbf{C}_{\rho}(0)) \cap \{\nabla^S R \neq 0\}} 
		\sum_{\kappa=1}^m R^{2\alpha+1} \left|\nabla^S R \cdot \nabla^S \bigg(\frac{z^{\kappa}}{R^{\alpha}}\bigg)\right|^2 
		\frac{\phi'(R/\rho)}{|\nabla^S R|^2} \,d|T|(X,S) 
		\nonumber\\&- C \rho^{-2\alpha+1} \widehat{E}(T,P_0,\mathbf{C}_{\rho/\gamma}(0))^{2+\varsigma} \nonumber
\end{align}
for all $\rho \in (0,\gamma]$, where 
$R = R(X) = |\pi(X)|$ for each $X \in \mathbb{R}^{n+m}$, $\varsigma = \varsigma(n,m,q) \in (0,1)$ is as in Theorem~\ref{lip approx thm}, and $C = C(n,m,q,\gamma) \in (0,\infty)$ is a constant.
\end{lemma}

\begin{proof}
For each $\rho \in (0,\gamma]$, let $\widehat{E}_{\rho} = \widehat{E}(T,P_0,\mathbf{C}_{\rho/\gamma}(0))$.  Note that since $0 \in \mathcal{B}_q$, by Lemma~\ref{uniform branch set lemma} we have that $\widehat{E}_{\rho} \leq C\varepsilon_0 \rho^{1/q} \leq C\varepsilon_0$, where $C = C(n,m,q,\gamma) \in (0,\infty)$ is a constant.  By~\cite[Lemma~3.5 and Lemma~3.7]{KrumWica} and Theorem~\ref{lip approx thm}, for each $\rho \in (0,\gamma]$ there exists a Lipschitz function $u_{\rho} : B_{\rho}(0) \rightarrow \mathcal{A}_q(\mathbb{R}^m)$ and a closed set $K_{\rho} \subset B_{\rho}(0)$ such that 
\begin{gather}\label{weiss eqn1}
	T \llcorner (K_{\rho} \times \mathbb{R}^m) = (\op{graph} u_{\rho}) \llcorner (K_{\rho} \times \mathbb{R}^m), \quad 
	\sup_{B_{\rho}(0)} |u_{\rho}| \leq C \widehat{E}_{\rho} \,\rho, \quad \op{Lip} u_{\rho} \leq C \widehat{E}_{\rho}^{\varsigma}, \\
	\mathcal{L}^n(B_{\rho}(0) \setminus K_{\rho}) + \|T\|((B_{\rho}(0) \setminus K_{\rho}) \times \mathbb{R}^m) 
		\leq C \widehat{E}_{\rho}^{2+\varsigma} \,\rho^n , 
	\nonumber
\end{gather}
where $C = C(n,m,q,\gamma) \in (0,\infty)$ is a constant.  By the area formula, \eqref{weiss eqn1}, and~\cite[Lemma~3.5]{KrumWica}
\begin{align}\label{weiss eqn2}
	&\frac{1}{\rho^{n+2}} \int_{G_n(\mathbf{C}_{\rho}(0))} \op{dist}^2(X,P_0) \,\|\pi_S - \pi\|^2 \,d|T|(X,S) 
	\leq \frac{1}{\rho^{n+2}} \int_{B_{\rho}(0)} |u_{\rho}|^2 |Du_{\rho}|^2 + C \widehat{E}_{\rho}^{4+\varsigma} 
	\\ \leq\,& \frac{C\widehat{E}_{\rho}^{2\varsigma}}{\rho^{n+2}} \int_{B_{\rho}(0)} |u_{\rho}|^2 + C \widehat{E}_{\rho}^{4+\varsigma} 
	\leq \frac{C\widehat{E}_{\rho}^{2\varsigma}}{\rho^{n+2}} \int_{\mathbf{C}_{\rho}(0)} \op{dist}^2(X,P_0) \,d\|T\|(X) 
		+ C \widehat{E}_{\rho}^{4+\varsigma} 
	\leq C \widehat{E}_{\rho}^{2+\varsigma} \nonumber 
\end{align}
and similarly using~\cite[Lemma~3.7]{KrumWica}
\begin{align}\label{weiss eqn3}
	&\frac{1}{\rho^n} \int_{\mathbf{C}_{\rho}(0)} |\vec T - \vec P_0|^4 \,d\|T\|(X) 
	\leq \frac{C\widehat{E}_{\rho}^{2\varsigma}}{\rho^n} \int_{\mathbf{C}_{\rho}(0)} |\vec T - \vec P_0|^2 \,d\|T\|(X) 
		+ C \widehat{E}_{\rho}^{2+\varsigma}  
	\\ \leq\,& \frac{C\widehat{E}_{\rho}^{2\varsigma}}{\rho^{n+2}} \int_{\mathbf{C}_{\rho/\gamma}(0)} \op{dist}^2(X,P_0) \,d\|T\|(X) 
		+ C \widehat{E}_{\rho}^{2+\varsigma}  
	\leq C \widehat{E}_{\rho}^{2+\varsigma} \nonumber 
\end{align}
for all $\rho \in (0,\gamma]$, where $C = C(n,m,q,\gamma) \in (0,\infty)$ are constants.  Hence by~\cite[Lemma~4.5]{KrumWica} together with (4.61), and (4.62) in~\cite{KrumWica} 
\begin{align*}
	H'_{T,P_0,0}(\rho) \leq\,& -2\rho^{-n} \int \pi^{\perp}(X) \cdot \pi^{\perp}(\nabla^S R) \,\phi'(R/\rho) \,d|T|(X,S) 
		+ C\rho \,\widehat{E}_{\rho}^{2+\varsigma} , \\
	D_{T,P_0,0}(\rho) =\,& -\rho^{1-n} \int \pi^{\perp}(X) \cdot \pi^{\perp}(\nabla^S R) \,\phi'(R/\rho) \,d|T|(X,S) , \\
	D'_{T,P_0,0}(\rho) \geq\,& -2\rho^{-n} \int_{G_n(\mathbf{C}_{\rho}(0)) \cap \{\nabla^S R \neq 0\}} \frac{|\pi^{\perp}(\nabla^S R)|^2}{|\nabla^S R|^2} 
		\,R \,\phi'(R/\rho) \,d|T|(X,S) - C\rho \,\widehat{E}_{\rho}^{2+\varsigma} 
\end{align*}
for $\mathcal{L}^1$-a.e.~$\rho \in (0,\gamma]$, where $C = C(n,m,q,\gamma) \in (0,\infty)$ is a constant.  Hence 
\begin{align*}
	&\frac{d}{d\rho} \Big( \rho^{-2\alpha} (D_{T,P_0,0}(\rho) - \alpha \,H_{T,P_0,0}(\rho)) \Big) 
	\\=\,& \rho^{-2\alpha} (D'_{T,P_0,0}(\rho) - \alpha \,H'_{T,P_0,0}(\rho)) - 2\alpha\rho^{-2\alpha-1} (D_{T,P_0,0}(\rho) - \alpha \,H_{T,P_0,0}(\rho))  
	\\ \geq\,& -2\rho^{-n-2\alpha} \int_{G_n(\mathbf{C}_{\rho}(0)) \cap \{\nabla^S R \neq 0\}} 
		\left| \frac{R \,\pi^{\perp}(\nabla^S R)}{|\nabla^S R|^2} - \alpha \,\pi^{\perp}(X) \right|^2 |\nabla^S R|^2 \,\frac{1}{R} \,\phi'(R/\rho) \,d|T|(X,S) 
		\\&- C \rho^{-2\alpha+1} \widehat{E}_{\rho}^{2+\varsigma} 
\end{align*}
for $\mathcal{L}^1$-a.e.~$\rho \in (0,\gamma]$, where $C = C(n,m,q,\gamma) \in (0,\infty)$ is a constant.  To complete the proof, we note that 
\begin{equation*}
	\sum_{\kappa=1}^m \nabla^S R \cdot \nabla^S \bigg(\frac{z^{\kappa}}{R^{\alpha}}\bigg) \,e_{m+\kappa} 
		= \frac{R \,\pi^{\perp}(\nabla^S R) - \alpha \,\pi^{\perp}(X) \,|\nabla^S R|^2}{R^{\alpha+1}}
\end{equation*}
where $e_1,e_2,\ldots,e_{n+m}$ is the standard basis on $\mathbb{R}^{n+m}$.
\end{proof}

\begin{theorem}\label{keyest thm}
Let $\alpha = \ell_0/q_0 \neq 2$ where $\ell_0$ and $q_0$ are relatively prime positive integers with $q_0 \leq q$.  Let $\varsigma = \varsigma(n,m,q) \in (0,1)$ is as in Theorem~\ref{relharm separation3 thm}.  For each $\gamma \in (0,1)$ there exists $\varepsilon_0 = \varepsilon_0(n,m,q,\alpha,\gamma) \in (0,1)$ and $\beta_0 = \beta_0(n,m,q,\alpha,\gamma) \in (0,1)$ such that the following holds true.  Let $\varphi$ and $T$ satisfy Hypothesis~$(\star)$ and Hypothesis~$(\star\star)$.  Let $\mathcal{B}_q$ be as in Definition~\ref{introduction B defn} (with $\mu = 1/q$).  If $0 \in \mathcal{B}_q$ such that $T$ is tangent to $P_0 = \mathbb{R}^n \times \{0\}$ at the origin and $\mathcal{N}_{T,{\rm Pl}}(0) \geq \alpha$, then 
\begin{align}\label{keyest concl}
	&\int_{G_n(\mathbf{C}_{\gamma}(0)) \cap \{\nabla^S R \neq 0\}} 
		\sum_{\kappa=1}^m \frac{R^{2-n}}{|\nabla^S R|^2} \left|\nabla^S R \cdot \nabla^S \bigg(\frac{z^{\kappa}}{R^{\alpha}}\bigg)\right|^2 \,d|T|(X,S) 
			\\&+ \int_{G_n(\mathbf{C}_{\gamma}(0))} \sum_{j=1}^{n-2} |\pi_S(e_{2+j})|^2 \,d|T|(X,S)
	\leq C E_{\varsigma/2}(T,\varphi,\mathbf{C}_1(0))^2 , \nonumber 
\end{align}
where $e_1,e_2,\ldots,e_{n+m}$ is the standard basis for $\mathbb{R}^{n+m}$and $C = C(n,m,q,\alpha,\gamma) \in (0,\infty)$ is a constant.
\end{theorem}

\begin{proof}
Without loss of generality, assume that $\gamma \in [1/2,1)$.  Set $\overline{\gamma} = \frac{2\gamma}{1+\gamma}$ and note that $1/2 \leq \gamma < \overline{\gamma} < 1$.  By Theorem~\ref{frequency monotonicity thm}(c)(d), either $\mathcal{N}_{T,{\rm Pl}}(0) = \alpha$ and $\limsup_{\rho\rightarrow 0^+} \rho^{-\alpha} H_{T,P_0,0}(\rho) < \infty$ or $\mathcal{N}_{T,{\rm Pl}}(0) > \alpha$ and $\lim_{\rho\rightarrow 0^+} \rho^{-\alpha} H_{T,P_0,0}(\rho) = 0$.  Hence 
\begin{equation}\label{keyest eqn1}
	\lim_{\rho\rightarrow 0^+} \rho^{-2\alpha} \,(D_{T,P_0,0}(\rho) - \alpha \,H_{T,P_0,0}(\rho)) 
	= \lim_{\rho\rightarrow 0^+} \rho^{-2\alpha} H_{T,P_0,0}(\rho) \,(N_{T,P_0,0}(\rho) - \alpha) = 0.
\end{equation}
Recalling Remark~\ref{planar freq rmk} and applying Theorem~\ref{frequency monotonicity thm}(e) 
\begin{equation}\label{keyest eqn2}
	\widehat{E}(T,P_0,\mathbf{C}_{\rho}(0)) \leq C \rho^{\alpha-1} \widehat{E}(T,P_0,\mathbf{C}_1(0))
\end{equation}
for all $\rho \in (0,1]$, where $C = C(n,m,q,\alpha) \in (0,\infty)$ is a constant.  Hence integrating \eqref{weiss concl} and using \eqref{keyest eqn1} and \eqref{keyest eqn2}, 
\begin{align}\label{keyest eqn3}
	&\rho^{-2\alpha} \,(D_{T,P,0}(\rho) - \alpha \,H_{T,P,0}(\rho))
	\\ \geq\,& -2 \int_{0}^{\rho} \tau^{-n-2\alpha} \int \sum_{\kappa=1}^m R^{2\alpha+1} \left|\nabla^S R 
		\cdot \nabla^S \bigg(\frac{z^{\kappa}}{R^{\alpha}}\bigg)\right|^2 \frac{\phi'(R/\tau)}{|\nabla^S R|^2} \,d|T|(X,S) \,d\tau 
		\nonumber\\&- C \int_{0}^{\rho} \tau^{\varsigma (\alpha-1) - 1} \widehat{E}(T,P_0,\mathbf{C}_1(0))^{2+\varsigma} \,d\tau \nonumber 
	\\ \geq\,& -2 \int_{0}^{\rho} \tau^{-n-2\alpha} \int_{G_n(\mathbf{C}_{\overline{\gamma}\rho}(0))} \sum_{\kappa=1}^m R^{2\alpha+1} \left|\nabla^S R 
		\cdot \nabla^S \bigg(\frac{z^{\kappa}}{R^{\alpha}}\bigg)\right|^2 \frac{\phi'(R/\tau)}{|\nabla^S R|^2} \,d|T|(X,S) \,d\tau 
		\nonumber\\&- C \rho^{\varsigma (\alpha-1)} \widehat{E}(T,P_0,\mathbf{C}_1(0))^{2+\varsigma} \,d\tau \nonumber 
\end{align}
for all $\rho \in \big(0,\frac{3+\gamma}{4}\big]$, where $C = C(n,m,q,\alpha,\gamma) \in (0,\infty)$ are constants.  Recalling that $1/2 \leq \overline{\gamma} < 1$, if $\rho \in \big(0,\frac{3+\gamma}{4}\big]$ and $X \in \mathbf{C}_{\overline{\gamma}\rho}(0)$, then $R < R/\overline{\gamma} \leq \min\{\rho,2R\}$ and thus 
\begin{equation}\label{keyest eqn4}
	-2 \int_{0}^{\rho} \tau^{-n-2\alpha} \,\phi'(R/\tau) \,d\tau 
	= 4 \int_{R}^{\min\{\rho,2R\}} \tau^{-n-2\alpha} \,d\tau 
	\geq c R^{1-n-2\alpha} , 
\end{equation}
where $c = \frac{4}{n+2\alpha-1} (1 - \overline{\gamma}^{1-n-2\alpha}) > 0$.  Hence using Fubini's theorem and \eqref{keyest eqn4}, \eqref{keyest eqn3} gives us 
\begin{align}\label{keyest eqn5}
	\rho^{-2\alpha} \,(D_{T,P,0}(\rho) - \alpha \,H_{T,P,0}(\rho))
	\geq\,& c \int_{G_n(\mathbf{C}_{\overline{\gamma}\rho}(0))} \sum_{\kappa=1}^m \frac{R^{2-n}}{|\nabla^S R|^2} 
		\left|\nabla^S R \cdot \nabla^S \bigg(\frac{z^{\kappa}}{R^{\alpha}}\bigg)\right|^2 \,d|T|(X,S) 
		\\&- C \rho^{\varsigma (\alpha-1)} \widehat{E}(T,P_0,\mathbf{C}_1(0))^{2+\varsigma} \nonumber 
\end{align}
for all $\rho \in \big(0,\frac{3+\gamma}{4}\big]$, where $c = c(n,m,q,\alpha,\gamma) > 0$ and $C = C(n,m,q,\alpha,\gamma) \in (0,\infty)$ are constants.  On the other hand, by the definitions of $D_{T,P_0,0}$ and $H_{T,P_0,0}$ 
\begin{align}\label{keyest eqn6}
	&D_{T,P_0,0}(\rho) - \alpha \,H_{T,P_0,0}(\rho)
	\\=\,& \rho^{2-n} \int \left( \frac{1}{2} \,\|\pi_S - \pi\|^2 \,\phi(r/\rho) 
		+ \alpha \,|z|^2 \,|\nabla^S R|^2 \,\frac{1}{\rho R} \,\phi'(R/\rho) \right) d|T|(X,S) \nonumber 
\end{align}
for all $\rho \in \big(0,\frac{3+\gamma}{4}\big]$.  Thus by combining \eqref{keyest eqn5} and \eqref{keyest eqn6}
\begin{align}\label{keyest eqn7}
	&c \int_{G_n(\mathbf{C}_{\overline{\gamma}\rho}(0))} \sum_{\kappa=1}^m \frac{R^{2-n}}{|\nabla^S R|^2} 
		\left|\nabla^S R \cdot \nabla^S \bigg(\frac{z^{\kappa}}{R^{\alpha}}\bigg)\right|^2 \,d|T|(X,S) 
		\\ \leq\,& \rho^{2-n-2\alpha} \int \left( \frac{1}{2} \,\|\pi_S - \pi\|^2 \,\phi(R/\rho) + \alpha \,|z|^2 \,|\nabla^S R|^2 \,\frac{1}{\rho R} \,\phi'(R/\rho) \right) 
		d|T|(X,S) \nonumber\\&+ C \rho^{\varsigma (\alpha-1)} \widehat{E}(T,P_0,\mathbf{C}_1(0))^{2+\varsigma} \nonumber 
\end{align}
for all $\rho \in \big(0,\frac{3+\gamma}{4}\big]$, where $c = c(n,m,q,\alpha,\gamma) > 0$ and $C = C(n,m,q,\alpha,\gamma) \in (0,\infty)$ are constants.  Let $\chi : (0,\infty) \rightarrow [0,\infty)$ be a smooth function such that $\op{spt} \chi = \big[\frac{1+\gamma}{2},\frac{3+\gamma}{4}\big]$ and $\int_0^{\infty} \chi(\rho) \,d\rho = 1$.  Define $\psi : [0,\infty) \rightarrow \mathbb{R}$ by 
\begin{equation*}
	 \psi(s) = \int_0^{\infty} \chi(\rho) \,\phi(s/\rho) \,d\rho 
\end{equation*}
for all $s \in [0,\infty)$.  Using change of variables $\sigma = s/\rho$, $\psi(s) = -\int_0^{\infty} \chi(s/\sigma) \,\phi(\sigma) \,s \,\sigma^{-2} \,d\sigma$ and it follows that $\psi \in C^{\infty}([0,\infty))$.  Since $\phi \in C^{0,1}([0,\infty))$, $\psi'(s) = \int_0^{\infty} \chi(\rho) \,\frac{1}{\rho} \,\phi'(s/\rho) \,d\rho$ for all $s \in [0,\infty)$.  Since $\phi(s) = 1$ if $0 \leq s \leq 1/2$ and $\op{spt} \phi = \big[0,\frac{3+\gamma}{4}\big]$, $\psi(s) = 1$ if $0 \leq s \leq \frac{1+\gamma}{4}$ and $\op{spt} \psi = \big[0,\frac{3+\gamma}{4}\big]$.  Note that $\overline{\gamma} \rho = \frac{2\gamma}{1+\gamma} \,\rho \geq \gamma$ for all $\rho \in \big[\frac{1+\gamma}{2}, \frac{3+\gamma}{4}\big]$.  Hence multiplying both sides of \eqref{keyest eqn7} by $\rho^{n-2+2\alpha} \,\chi(\rho)$ and integrating over $\rho \in \big[\frac{1+\gamma}{2}, \frac{3+\gamma}{4}\big]$, 
\begin{align}\label{keyest eqn8}
	&c \int_{G_n(\mathbf{C}_{\gamma}(0))} \sum_{\kappa=1}^m \frac{R^{2-n}}{|\nabla^S R|^2} 
		\left|\nabla^S R \cdot \nabla^S \bigg(\frac{z^{\kappa}}{R^{\alpha}}\bigg)\right|^2 \,d|T|(X,S) 
		\\ \leq\,& \int \left( \frac{1}{2} \,\|\pi_S - \pi\|^2 \,\psi(R) + \alpha \,|z|^2 \,|\nabla^S R|^2 \,\frac{1}{R} \,\psi'(R) \right) d|T|(X,S) 
		+ C \widehat{E}(T,P_0,\mathbf{C}_1(0))^{2+\varsigma} , \nonumber 
\end{align}
where $c = c(n,m,q,\alpha,\gamma) > 0$ and $C = C(n,m,q,\alpha,\gamma) \in (0,\infty)$ are constants.

Setting $\zeta = \psi(R) \,(x,0,\alpha z)$ in \eqref{first variation} and noting that 
\begin{align*}
	\op{div}_S (x,0,z) =\,& n - \op{div}_S (0,y,0) = n - \sum_{j=1}^{n-2} |\pi_S(e_{2+j})|^2 = 2 + \sum_{j=1}^{n-2} |\pi_{S^{\perp}}(e_{2+j})|^2 , \\
	\op{div}_S (0,0,z) =\,& \sum_{j=1}^m |\pi_S(e_{n+j})|^2 = \frac{1}{2} \,\|\pi_S - \pi\|^2 , \\
	\nabla^S R \cdot (x,0,\alpha z) =\,& \frac{1}{R} \,\pi_S(x,y,0) \cdot (x,0,\alpha z) = \frac{1}{R} \big( |x|^2 - \pi_{S^{\perp}}(x,y,0) \cdot (x,0,\alpha z) \big) 
\end{align*}
gives us 
\begin{align}\label{keyest eqn9}
	0 =\,& \int \left( 2 + \frac{\alpha-1}{2} \,\|\pi_S - \pi\|^2 + \sum_{j=1}^{n-2} |\pi_{S^{\perp}}(e_{2+j})|^2 \right) \psi(R) \,d|T|(X,S) 
		\\&+ \int \left( |x|^2 - \pi_{S^{\perp}}(x,y,0) \cdot (x,0,\alpha z) \right) \frac{1}{R} \,\psi'(R) \,d|T|(X,S) . \nonumber 
\end{align}
Replacing $T$ with $q \llbracket P \rrbracket$ in \eqref{keyest eqn9} and then using $\pi_{\#} T = q \llbracket P \rrbracket$ and the area formula 
\begin{equation*}
	0 = q \int_P \left( 2 \psi(R) + \frac{|x|^2}{R} \,\psi'(R) \right) d\mathcal{H}^n  
	= \int \left( 2 \psi(R) + \frac{|x|^2}{R} \,\psi'(R) \right) J\pi_P \,d\|T\|(X) , 
\end{equation*}
where $J\pi_P : \op{spt} T \rightarrow \mathbb{R}$ is the Jacobian on $T$ defined by $\pi_{P\#} \vec T(X) = J\pi_P(X) \,\vec P$ for $\|T\|$-a.e.~$X \in \mathbf{C}_1(0)$.  By Taylor's theorem 
\begin{equation*}
	\left| J\pi_P - 1 + \frac{1}{4} \,\|\pi_S - \pi_P\|^2 \right| \leq C |\vec T - \vec P|^4 
\end{equation*}
for some constant $C = C(n,m) \in (0,\infty)$.  Thus 
\begin{equation}\label{keyest eqn10}
	0 \leq \int \left( 2 \psi(R) + \frac{|x|^2}{R} \,\psi'(R) \right) \left( 1 - \frac{1}{4} \,\|\pi_S - \pi_P\|^2 + C |\vec T - \vec P|^4 \right) d|T|(X,S) , 
\end{equation}
where $C = C(n,m) \in (0,\infty)$ is a constant.  Arguing as in \eqref{weiss eqn3}, $\int_{\mathbf{C}_{(3+\gamma)/4}(0)} |\vec T - \vec P_0|^4 \,d\|T\|(X) \leq C(n,m,q,\gamma) \,\widehat{E}(T,P_0,\mathbf{C}_1(0))^{2+\varsigma}$.  Thus subtracting \eqref{keyest eqn9} and \eqref{keyest eqn10} 
\begin{align}\label{keyest eqn11}
	0 \geq\,& \int \left( \frac{\alpha}{2} \,\|\pi_S - \pi\|^2 + \sum_{j=1}^{n-2} |\pi_{S^{\perp}}(e_{2+j})|^2 \right) \psi(R) \,d|T|(X,S)
		\\&+ \int \left( \frac{1}{4} |x|^2 \|\pi_S - \pi\|^2 - \pi_{S^{\perp}} (x,y,0) \cdot (x,0,\alpha z) \right) \frac{1}{R} \,\psi'(R) \,d|T|(X,S) \nonumber 
		\\&- C \widehat{E}(T,P_0,\mathbf{C}_1(0))^{2+\varsigma} , \nonumber
\end{align}
where $C = C(n,m,q,\gamma) \in (0,\infty)$ is a constant.  Combining \eqref{keyest eqn8} and \eqref{keyest eqn11} 
\begin{align*}%\label{keyest eqn12}
	&c \int_{G_n(\mathbf{C}_{\gamma}(0))} \sum_{\kappa=1}^m \frac{R^{2-n}}{|\nabla^S R|^2} 
		\left|\nabla^S R \cdot \nabla^S \bigg(\frac{z^{\kappa}}{R^{\alpha}}\bigg)\right|^2 \,d|T|(X,S) 
		+ \int \sum_{j=1}^{n-2} |\pi_{S^{\perp}}(e_{2+j})|^2 \,\psi(R) \,d|T|(X,S) 
	\\ \leq\,& \int \left( -\frac{1}{4} \,|x|^2 \,\|\pi_S - \pi\|^2 + \pi_{S^{\perp}} (x,y,0) \cdot (x,0,\alpha z) + \alpha^2 |z|^2 |\nabla^S R|^2 \right) 
		\frac{1}{R} \,\psi'(R) \,d|T|(X,S) \nonumber \\&+ C \widehat{E}(T,P_0,\mathbf{C}_1(0))^{2+\varsigma} , \nonumber
\end{align*}
where $c = (n,m,q,\alpha,\gamma) > 0$ and $C = C(n,m,q,\alpha,\gamma) \in (0,\infty)$ are constants.  
Using the fact that $|\psi'(s)|^2/\psi(s) \leq 2 \|\psi''\|_{L^{\infty}([0,\infty))} \leq C(\gamma)$ for all $s \in \big[0,\frac{3+\gamma}{4}\big)$
and Cauchy's inequality 
\begin{equation*}
	\pi_{S^{\perp}}(0,y,0) \cdot (x,0,\alpha z) \,\psi'(R) 
	\leq \frac{1}{2} \sum_{j=1}^{n-2} |\pi_{S^{\perp}}(e_{2+j})|^2 \,\psi(R) + \frac{1}{2} |\pi_{S^{\perp}}(x,0,\alpha z)|^2 \,\frac{|\psi'(R)|^2}{\psi(R)}
\end{equation*}
for $|T|$-a.e.~$(X,S) \in G_n(\mathbf{C}_{(3+\gamma)/4}(0))$, we obtain 
\begin{align}\label{keyest eqn13} 
	&c \int_{G_n(\mathbf{C}_{\gamma}(0))} \sum_{\kappa=1}^m \frac{R^{2-n}}{|\nabla^S R|^2} 
		\left|\nabla^S R \cdot \nabla^S \bigg(\frac{z^{\kappa}}{R^{\alpha}}\bigg)\right|^2 \,d|T|(X,S) 
	\\&+ \frac{1}{2} \int \sum_{j=1}^{n-2} |\pi_{S^{\perp}}(e_{2+j})|^2 \,\psi(R) \,d|T|(X,S) \nonumber
	\\ \leq\,& \int \left( -\frac{1}{4} \,|x|^2 \,\|\pi_S - \pi\|^2 + \pi_{S^{\perp}} (x,0,0) \cdot (x,0,\alpha z) + \alpha^2 |z|^2 |\nabla^S R|^2 \right) 
		\frac{1}{R} \,\psi'(R) \,d|T|(X,S) \nonumber \\&+ \frac{1}{2} \int_{G_n(\mathbf{C}_{(3+\gamma)/4}(0))} |\pi_{S^{\perp}}(x,0,\alpha z)|^2 \,d|T|(X,S)	
		+ C \widehat{E}(T,P_0,\mathbf{C}_1(0))^{2+\varsigma} , \nonumber
\end{align}
where $c = (n,m,q,\alpha,\gamma) > 0$ and $C = C(n,m,q,\alpha,\gamma) \in (0,\infty)$ are constants.  

Recall that $A_{\rho,\kappa}(\zeta) = \{ (r\omega,y) : (r,y) \in B^{n-1}_{\kappa (1-\gamma) \rho/8}(\rho,\zeta), \,\omega \in \mathbb{S}^1 \}$ for each $\kappa \in (0,1]$, $\rho > 0$, and $\zeta \in \mathbb{R}^{n-2}$.  Note that $A_{\rho,\kappa}(\zeta) \cap A_{\rho',\kappa}(\zeta') \neq \emptyset$ $\iff$ $|(\rho,\zeta) - (\rho',\zeta')| < \frac{1}{8} \kappa (1-\gamma)(\rho+\rho')$, in which case $\frac{8-\kappa(1-\gamma)}{8+\kappa(1-\gamma)} \,\rho < \rho' < \frac{8+\kappa(1-\gamma)}{8-\kappa(1-\gamma)} \,\rho$.  By applying the Vitali covering lemma, there exists a countable collection $\mathcal{I}$ of $(\rho,\zeta)$ with $\rho > 0$, $\zeta \in \mathbb{R}^{n-2}$, and $\rho^2 + |\zeta|^2 < (3+\gamma)^2/16$ such that $\{A_{\rho,1/20}(\zeta)\}_{(\rho,\zeta) \in \mathcal{I}}$ is a collection of pairwise disjoint annuli and $\{A_{\rho,1/4}(\zeta)\}_{(\rho,\zeta) \in \mathcal{I}}$ covers $B_{(3+\gamma)/4}(0)$.  Observe that if $A_{\rho,1}(\zeta) \cap A_{\rho',1}(\zeta') \neq \emptyset$ then $|(\rho,\zeta) - (\rho',\zeta')| < \frac{1}{8} (1-\gamma)(\rho+\rho')$ and thus $\frac{7+\gamma}{9-\gamma} \,\rho < \rho' < \frac{9-\gamma}{7+\gamma} \,\rho$.  Hence 
\begin{equation*}
	B^{n-1}_{\frac{(1-\gamma)(7+\gamma) \rho}{160(9-\gamma)}}(\rho',\zeta') 
	\subset B^{n-1}_{\frac{(1-\gamma)\rho'}{160}}(\rho',\zeta') 
	\subset B^{n-1}_{ \frac{(1-\gamma)\rho}{160} \left( 20 + \frac{21 (9-\gamma)}{7+\gamma} \right)}(\rho,\zeta)
\end{equation*}
which since $\{B^{n-1}_{(1-\gamma)\rho'/160}(\rho',\zeta')\}_{(\rho',\zeta') \in \mathcal{I}}$ is a pairwise disjoint collection of ball implies that for each $(\rho,\zeta) \in \mathcal{I}$ 
\begin{equation*}
	\#\{ (\rho',\zeta') \in \mathcal{I} : A_{\rho,1}(\zeta) \cap A_{\rho',1}(\zeta') \neq \emptyset \} \leq C(n,\gamma) .
\end{equation*}
Thus there exists an integer $N \leq C(n,\gamma)$ and pairwise disjoint sets $\mathcal{I}_1, \mathcal{I}_2, \ldots, \mathcal{I}_N \subset \mathcal{I}$ such that $\mathcal{I} = \bigcup_{j=1}^N \mathcal{I}_j$ and $\{A_{\rho,1}(\zeta)\}_{(\rho,\zeta) \in \mathcal{I}_j}$ is a collection of pairwise disjoint annuli for each $j \in \{1,2,\ldots,N\}$.  Let $\{\eta_{(\rho,\zeta)}\}_{(\rho,\zeta) \in \mathcal{I}}$ is a smooth partition of unity subordinate to $\{A_{\rho,1/2}(\zeta)\}_{(\rho,\zeta) \in \mathcal{I}}$ such that 
\begin{gather*}%\label{keyest eqn14}
	\eta_{(\rho,\zeta)}(x,y) = \eta_{(\rho,\zeta)}(\widetilde{x},y) \text{ whenever } |x| = |\widetilde{x}|, \\
	\op{spt} \eta_{(\rho,\zeta)} \subset A_{\rho,1/2}(\zeta), \quad\quad 
	\sum_{(\rho,\zeta) \in \mathcal{I}} \eta_{(\rho,\zeta)} = 1 \text{ on } B_{(3+\gamma)/4}(0). \nonumber
\end{gather*}
We claim that for each $(\rho,\zeta) \in \mathcal{I}$ 
\begin{align}\label{keyest eqn15} 
	&\int \left( -\frac{1}{4} |x|^2 \|\pi_S - \pi\|^2 + \pi_{S^{\perp}} (x,0,0) \cdot (x,0,\alpha z) + \alpha^2 |z|^2 |\nabla^S R|^2 \right) 
		\frac{1}{R} \psi'(R) \,\eta_{(\rho,\zeta)} \,d|T|(X,S) 
	\\&\leq C \int_{\mathbf{A}_{\rho,1}(\zeta)} \op{dist}^2(X,\op{spt} \op{graph} \varphi) \,d\|T\|(X) 
		+ C \widehat{E}^{\varsigma/2} \int_{\mathbf{A}_{\rho,1}(\zeta)} \op{dist}^2(X,P_0) \,d\|T\|(X) \nonumber
\end{align}
and 
\begin{align}\label{keyest eqn16} 
	&\int_{G_n(\mathbf{C}_{(3+\gamma)/4}(0))} |\pi_{S^{\perp}}(x,0,\alpha z)|^2 \,\eta_{(\rho,\zeta)} \,d|T|(X,S) 
	\\&\leq C \int_{\mathbf{A}_{\rho,1}(\zeta)} \op{dist}^2(X,\op{spt} \op{graph} \varphi) \,d\|T\|(X) 
		+ C \widehat{E}^{\varsigma/2} \int_{\mathbf{A}_{\rho,1}(\zeta)} \op{dist}^2(X,P_0) \,d\|T\|(X) , \nonumber
\end{align}
where $\widehat{E} = \widehat{E}(T,P_0,\mathbf{C}_1(0))$ and $C = C(n,m,q,\alpha,\gamma) \in (0,\infty)$ is a constant.  Then by summing \eqref{keyest eqn15} and \eqref{keyest eqn16} over $(\rho,\zeta) \in \mathcal{I} = \bigcup_{j=1}^N \mathcal{I}_j$ and keeping in mind that $\psi(R) = 0$ on $\mathbf{C}_1(0) \setminus \mathbf{C}_{(3+\gamma)/4}(0)$ 
\begin{align}\label{keyest eqn17} 
	&\int \left( -\frac{1}{4} \,|x|^2 \,\|\pi_S - \pi\|^2 + \pi_{S^{\perp}} (x,0,0) \cdot (x,0,\alpha z) + \alpha^2 |z|^2 |\nabla^S R|^2 \right) 
		\frac{1}{R} \,\psi'(R) \,d|T|(X,S) 
	\\&\leq C E_{\varsigma/2}(T,P_0,\mathbf{C}_1(0)) , \nonumber
\end{align}
and 
\begin{align}\label{keyest eqn18} 
	\int_{G_n(\mathbf{C}_{(3+\gamma)/4}(0))} |\pi_{S^{\perp}}(x,0,\alpha z)|^2 \,d|T|(X,S)	 \leq C E_{\varsigma/2}(T,P_0,\mathbf{C}_1(0)) , 
\end{align}
where $C = C(n,m,q,\alpha,\gamma) \in (0,\infty)$ is a constant.  Combining \eqref{keyest eqn13}, \eqref{keyest eqn17}, and \eqref{keyest eqn18} will then prove \eqref{keyest concl}. 

To see \eqref{keyest eqn15} and \eqref{keyest eqn16}, fix $(\rho,\zeta) \in \mathcal{I}$.  Let 
\begin{equation*}
	\Xi = \{ (x,y,z) \in \mathbf{C}_1(0) : |z| \leq 2 \|\varphi\|_{L^{\infty}(B_1(0))} |x|^{\alpha} \} .
\end{equation*}
Cover $A_{\rho,1}(\zeta)$ with a finite collection of ball $\{B_{3(1-\gamma)\rho/16}(x_k,\zeta)\}_{k=1}^N$, where $x_k \in \partial B^2_{\rho}(0)$ and $N$ is an integer with $N \leq C(n,\gamma)$.  Since $\alpha > 1$ and $\|\varphi\|_{L^2(B_1(0))} \leq C\varepsilon_0$, $|z| \leq C\varepsilon_0 |x| \text{ for all } (x,y,z) \in \Xi$.  Thus assuming $\varepsilon_0$ is sufficiently small, $\{\mathbf{B}_{3(1-\gamma)\rho/16}(x_k,\zeta,0)\}_{k=1}^N$ covers $\op{spt} T \cap \mathbf{A}_{\rho,1}(\zeta) \cap \Xi$.  By the monotonicity formula and $\|T\|(\mathbf{C}_1(0)) \leq (q+1/2) \,\omega_n$ 
\begin{align}\label{keyest eqn19} 
	\|T\|(\mathbf{A}_{\rho,1}(\zeta) \cap \Xi) 
	\leq\,& \sum_{k=1}^N \|T\|(\mathbf{B}_{(1-\gamma)\rho/4}(x_k,\zeta,0)) 
	\leq \sum_{k=1}^N \rho^n \|T\|(\mathbf{B}_{(1-\gamma)/4}(x_k,\zeta,0)) 
	\\ \leq\,& N \rho^n \|T\|(\mathbf{C}_1(0)) 
	\leq C \rho^n , \nonumber 
\end{align}
where $C = C(n,q,\gamma) \in (0,\infty)$ is a constant.  By the triangle inequality and \eqref{keyest eqn19} 
\begin{align}\label{keyest eqn20}
	\int_{\mathbf{A}_{\rho,1}(\zeta) \cap \Xi} \op{dist}^2(X,P_0) \,d\|T\|(X) 
	\leq\,& 2 \int_{\mathbf{A}_{\rho,1}(\zeta) \cap \Xi} \op{dist}^2(X, \op{spt} \op{graph} \varphi) \,d\|T\|(X) 
		\\&+ C \rho^{n+2\alpha} \,\|\varphi\|_{L^{\infty}(B_1(0))}^2 , \nonumber 
\end{align}
where $C = C(n,q,\alpha,\gamma) \in (0,\infty)$ is a constant.  On the other hand, by the triangle inequality for each $(x,y,z) \in \mathbf{C}_1(0) \setminus \Xi$, $|z| \leq |z - \varphi(x,y)| + |\varphi(x,y)| \leq |z - \varphi(x,y)| + \frac{1}{2} |z|$ and thus $|z| \leq 2 |z - \varphi(x,y)|$.  Hence 
\begin{align}\label{keyest eqn21}
	\int_{\mathbf{A}_{\rho,1}(\zeta) \setminus \Xi} \op{dist}^2(X,P_0) \,d\|T\|(X) 
	\leq 16 \int_{\mathbf{A}_{\rho,1}(\zeta) \setminus \Xi} \op{dist}^2(X, \op{spt} \op{graph} \varphi) \,d\|T\|(X) .
\end{align}
Adding \eqref{keyest eqn20} and \eqref{keyest eqn21} 
\begin{align}\label{keyest eqn22}
	\int_{\mathbf{A}_{\rho,1}(\zeta)} \op{dist}^2(X,P_0) \,d\|T\|(X) 
	\leq\,& 16 \int_{\mathbf{A}_{\rho,1}(\zeta)} \op{dist}^2(X, \op{spt} \op{graph} \varphi) \,d\|T\|(X) 
		\\&+ C \rho^{n+2\alpha} \,\|\varphi\|_{L^{\infty}(B_1(0))}^2 , \nonumber 
\end{align}
where $C = C(n,q,\alpha,\gamma) \in (0,\infty)$ is a constant.  

Suppose that  
\begin{equation}\label{keyest eqn23}
	\int_{\mathbf{A}_{\rho,1}(\zeta)} \op{dist}^2(X, \op{spt} \op{graph} \varphi) \,d\|T\|(X) \leq \rho^{n+2\alpha} \|\varphi\|_{L^{\infty}(B_1(0))}^2 . 
\end{equation}
(We later consider the case in which \eqref{keyest eqn29} holds true below.)  Let 
\begin{align*}%\label{keyest eqn24}
	E^2 =\,& \frac{1}{\rho^{n+2}} \int_{\mathbf{A}_{\rho,1}(\zeta)} \op{dist}^2(X, \op{spt} \op{graph} \varphi) \,d\|T\|(X) 
		\\&+ \widehat{E}(T,P_0,\mathbf{C}_1(0))^{\varsigma/2} \cdot \frac{1}{\rho^{n+2}} \int_{\mathbf{A}_{\rho,1}(\zeta)} \op{dist}^2(X,P_0) \,d\|T\|(X) . \nonumber 
\end{align*}
(for simplicity omitting the dependence of $E$ on $T$, $\varphi$, $\rho$, and $\zeta$).  By Remark~\ref{graphical rmk}(2), $\|\varphi\|_{L^{\infty}(B_1(0))} \leq C \widehat{E}(T,P_0,\mathbf{C}_1(0))$, where $C = C(n,m,q,\alpha,\gamma) \in (0,\infty)$ is a constant.  Thus by $\alpha > 1$, \eqref{keyest eqn22}, and \eqref{keyest eqn23}, $\widehat{E}(T,P_0,\mathbf{A}_{\rho,1}(\zeta)) \leq C \rho^{\alpha-1} \,\|\varphi\|_{L^{\infty}(B_1(0))} \leq C \widehat{E}(T,P_0,\mathbf{B}_1(0))$, where $C = C(n,m,q,\alpha,\gamma) \in (0,\infty)$ is a constant.  Hence $\widehat{E}(T,P_0,\mathbf{A}_{\rho,1}(\zeta)) \leq C \varepsilon_0$ and 
\begin{equation*}
	E_{\varsigma/2}(T,\varphi,\mathbf{A}_{\rho,1}(\zeta)) 
		+ \|D\varphi\|_{L^{\infty}(B_{\rho}(0,\zeta))}^{1/2} \,\widehat{E}(T,P_0,\mathbf{A}_{\rho,1}(\zeta)) \leq C E 
\end{equation*}
for some constant $C = C(n,m,q,\alpha,\gamma) \in (0,\infty)$.  By Hypothesis~$(\star)$ and Hypothesis~$(\star\star)$, we can apply Theorem~\ref{graphical_thm}(A) to obtain \eqref{graphical conclA2}.  By Hypothesis~$(\star)$, $(\partial T) \llcorner \mathbf{A}_{\rho,1}(\zeta) = 0$, $T \llcorner \mathbf{A}_{\rho,1}(\zeta)$ has finite height relative to $P_0$, and $\pi_{\#} (T \llcorner \mathbf{A}_{\rho,1}(\zeta)) = q \llbracket A_{\rho,1}(\zeta) \rrbracket$.  By this together with $\widehat{E}(T,P_0,\mathbf{A}_{\rho,1}(\zeta)) \leq C \varepsilon_0$ and \eqref{graphical conclA2}, we can apply Lemma~\ref{graphical_lemma0} with $\eta_{(0,\zeta),\rho\#} T$ and $\rho^{\alpha-1} \varphi$ in place of $T$ and $\varphi$ to show the following.  By Lemma~\ref{graphical_lemma0}(A), there are integers $\widehat{m}_j \geq 0$ with $\sum_{j=1}^J \widehat{m}_j q_j \leq q$ and $n$-dimensional locally area-minimizing rectifiable currents $T_j$ of $\mathbf{A}_{\rho,3/4}(\zeta)$ ($\widehat{m}_j$ and $T_j$ possibly zero) such that \begin{gather}\label{keyest eqn25} 
	T \llcorner \mathbf{A}_{\rho,3/4}(\zeta) = \sum_{j=1}^J T_j , \quad
	(\partial T_j) \llcorner \mathbf{A}_{\rho,3/4}(\zeta) = 0, \\ 
	\pi_{\#} T_j = \widehat{m}_j q_j \llbracket \mathbf{A}_{\rho,3/4}(\zeta) \rrbracket , \quad 
	\sup_{X \in \op{spt} T_j} \op{dist}(X, \op{spt} \op{graph} \varphi_j) \leq C E \rho, \nonumber
\end{gather} 
where $C = C(n,m,q,\alpha,\gamma,\kappa) \in (0,\infty)$ is a constant.  By Lemma~\ref{graphical_lemma0}(B), for each $j \in \{1,2,\ldots,J\}$ with $\widehat{m}_j > 0$ there exists a Lipschitz $\widehat{m}_j$-valued function $v_j : \op{graph} \varphi_j |_{A_{\rho,1/2}(\zeta)} \rightarrow \mathcal{A}_{\widehat{m}_j}(\mathbb{R}^m)$ and a closed set $K_j \subseteq A_{\rho,1/2}(\zeta)$ such that 
\begin{gather}\label{keyest eqn26}
	T_j \llcorner (K_j \times \mathbb{R}^m) = (\op{graph} u_j) \llcorner (K_j \times \mathbb{R}^m) , \\
	\mathcal{H}^n(A_{\rho,1/2}(\zeta) \setminus K_j) + \|T_j\|((A_{\rho,1/2}(\zeta) \setminus K_j) \times \mathbb{R}^m) 
		\leq C\rho^n \,\widehat{E}(T,P_0,\mathbf{A}_{\rho,1}(\zeta))^{2+\varsigma}, \nonumber\\
	\sup_{A_{1,1/2}(0)} |u_j| \leq C \rho \,\widehat{E}(T,P_0,\mathbf{A}_{\rho,1}(\zeta)) , \quad 
	\sup_{A_{1,1/2}(0)} |Du_j| \leq C \widehat{E}(T,P_0,\mathbf{A}_{\rho,1}(\zeta))^{\varsigma} , \nonumber
\end{gather}
and 
\begin{equation}\label{keyest eqn27}
	\frac{1}{\rho^2} \sup_{A_{\rho,1/2}(\zeta)} \sum_{j=1}^J \sum_{l=1}^{q_j} |v_{j,l}|^2 + \int_{A_{1,1/2}(\zeta)} \sum_{j=1}^J \sum_{l=1}^{q_j} |Dv_{j,l}|^2  \leq C E^2 , 
\end{equation}
where $v_{j,l}$ and $u_j$ are as in \eqref{v_localized} and \eqref{u defn} (with $\widehat{m}_j$ in place of $m_j$) and $C = C(n,m,q,p,\alpha,\gamma,\kappa) \in (0,\infty)$ is a constant.  By the area formula, \eqref{keyest eqn25}, and \eqref{keyest eqn26} 
\begin{align*} 
	&\int \left( -\frac{1}{4} |x|^2 \,\|\pi_S - \pi\|^2 + \pi_{S^{\perp}} (x,0,0) \cdot (x,0,\alpha z) + \alpha^2 |z|^2 |\nabla^S R|^2 \right) \frac{1}{R} \psi'(R) 
		\eta_{(\rho,\zeta)} \,d|T|(X,S) 
	\\ \leq\,& \int \sum_{j=1}^J \sum_{l=1}^{q_j} \sum_{h=1}^{\widehat{m}_j} \sum_{\kappa=1}^m \left( -\frac{1}{2} |x|^2 |Du_{j,l,h}^{\kappa}|^2 
		+ x \cdot D_x u_{j,l,h}^{\kappa} (x \cdot D_x u_{j,l,h}^{\kappa} - \alpha u_{j,l,h}^{\kappa}) \right) 
		\frac{1}{R} \,\psi'(R) \,\eta_{(\rho,\zeta)} \\&+ \int \sum_{j=1}^J \sum_{l=1}^{q_j} \sum_{h=1}^{\widehat{m}_j} \alpha^2 |u_{j,l,h}|^2 
		\,\frac{1}{R} \,\psi'(R) \,\eta_{(\rho,\zeta)} + C\rho^{n+2} \,\widehat{E}(T,P_0,\mathbf{A}_{\rho,1}(\zeta))^{2+\varsigma} , 
\end{align*}
where $u_{j,l,h}$ is as in \eqref{u defn}, $D_x u_{j,l,h}(x,y)$ is the derivative of $u_{j,l,h}$ with respect to the $x$-variable, and $C = C(n,m,q,p,\alpha,\gamma,\kappa) \in (0,\infty)$ is a constant.  Recall that $\widehat{E}(T,P_0,\mathbf{A}_{\rho,1}(\zeta)) \leq C \widehat{E}(T,P_0,\mathbf{C}_1(0))$, where $C = C(n,m,q,\alpha,\gamma) \in (0,\infty)$ is a constant.  Expressing $u_{j,l,h} = \varphi_{j,l} + v_{j,l,h}$ as in \eqref{u defn} and using \eqref{keyest eqn27}, 
\begin{align}\label{keyest eqn28}
	&\int \left( -\frac{1}{4} |x|^2 \|\pi_S - \pi\|^2 + \pi_{S^{\perp}} (x,0,0) \cdot (x,0,\alpha z) + \alpha^2 |z|^2 |\nabla^S R|^2 \right) \frac{1}{R} \psi'(R) 
		\eta_{(\rho,\zeta)} \,d|T|(X,S) 
	\\ \leq\,& \sum_{j=1}^J \widehat{m}_j \int \left( -\frac{1}{2} |x|^2 |D\varphi_j|^2 + \alpha^2 |\varphi_j|^2 \right) \frac{1}{R} \,\psi'(R) \,\eta_{(\rho,\zeta)} 
		\nonumber \\&+ \sum_{j=1}^J \widehat{m}_j \int \sum_{l=1}^{q_j} \left( -|x|^2 D\varphi_{j,l}^{\kappa} \cdot Dv_{j,l;a}^{\kappa} 
		+ x \cdot D_x \varphi_{j,l}^{\kappa} (x \cdot D_x v_{j,l;a}^{\kappa} - \alpha v_{j,l;a}^{\kappa}) \right) \frac{1}{R} \,\psi'(R) \,\eta_{(\rho,\zeta)} \nonumber 
		\\&+ \sum_{j=1}^J \widehat{m}_j \int \sum_{l=1}^{q_j} 2\alpha^2 \varphi_{j,l}^{\kappa} v_{j,l;a}^{\kappa} \,\frac{1}{R} \,\psi'(R) \,\eta_{(\rho,\zeta)} 
		+ C\rho^{n+2} E^2 , \nonumber
\end{align}
where $v_{j,l;a}(x,y) = \frac{1}{\widehat{m}_j} \sum_{h=1}^{\widehat{m}_j} v_{j,l,h}(x,y)$ for each $(x,y) \in A_{\rho,1/2}(\zeta)$ and $C = C(n,m,q,p,\alpha,\gamma,\kappa) \in (0,\infty)$ is a constant.  Arguing as on pages~40--41 of~\cite{KrumWic2} using integration by parts and $\varphi_j$ being homogeneous degree $\alpha$, harmonic, and independent of $y$ 
\begin{gather*}
	\int \left( -\frac{1}{2} |x|^2 |D\varphi_j|^2 + \alpha^2 |\varphi_j|^2 \right) \frac{1}{R} \,\psi'(R) \,\eta_{(\rho,\zeta)} = 0, \\
	\int \sum_{l=1}^{q_j} \left( -|x|^2 D\varphi_{j,l}^{\kappa} \cdot Dv_{j,l;a}^{\kappa} + x \cdot D_x \varphi_{j,l}^{\kappa} (x \cdot D_x v_{j,l;a}^{\kappa} 
		- \alpha v_{j,l;a}^{\kappa}) + 2\alpha^2 \varphi_{j,l}^{\kappa} v_{j,l;a}^{\kappa} \right) \frac{1}{R} \,\psi'(R) \,\eta_{(\rho,\zeta)} = 0
\end{gather*}
and thus \eqref{keyest eqn28} gives us \eqref{keyest eqn15}.  Similarly, by the area formula, \eqref{keyest eqn25}, and \eqref{keyest eqn26} 
\begin{align*} 
	&\int_{G_n(\mathbf{C}_{(3+\gamma)/4}(0))} |\pi_{S^{\perp}}(x,0,\alpha z)|^2 \,\eta_{(\rho,\zeta)} \,d|T|(X,S) 
	\\ \leq\,& \int_{B_{(3+\gamma)/4}(0)} \sum_{j=1}^J \sum_{l=1}^{q_j} \sum_{h=1}^{\widehat{m}_j} |x \cdot D_x u_{j,l,h} - \alpha u_{j,l,h}|^2 
	+ C\rho^{n+2} \,\widehat{E}(T,P_0,\mathbf{A}_{\rho,1}(\zeta))^{2+\varsigma} , 
\end{align*}
where $C = C(n,m,q,p,\alpha,\gamma,\kappa) \in (0,\infty)$ is a constant.  Recall $\widehat{E}(T,P_0,\mathbf{A}_{\rho,1}(\zeta)) \leq C \widehat{E}(T,P_0,\mathbf{C}_1(0))$, where $C = C(n,m,q,\alpha,\gamma) \in (0,\infty)$ is a constant.  Expressing $u_{j,l,h} = \varphi_{j,l} + v_{j,l,h}$ as in \eqref{u defn} and using \eqref{keyest eqn27} gives us \eqref{keyest eqn16}. 

Now suppose instead that $(\rho,\zeta) \in \mathcal{I}$ with 
\begin{equation}\label{keyest eqn29}
	\int_{\mathbf{A}_{\rho,1}(\zeta)} \op{dist}^2(X, \op{spt} \op{graph} \varphi) \,d\|T\|(X) \geq \rho^{n+2\alpha} \|\varphi\|_{L^{\infty}(B_1(0))}^2 . 
\end{equation}
By \eqref{keyest eqn22} and \eqref{keyest eqn29} 
\begin{equation}\label{keyest eqn30}
	\int_{\mathbf{A}_{\rho,1}(\zeta)} \op{dist}^2(X,P_0) \,d\|T\|(X) 
	\leq C \int_{\mathbf{A}_{\rho,1}(\zeta)} \op{dist}^2(X, \op{spt} \op{graph} \varphi) \,d\|T\|(X) 
\end{equation}
where $C = C(n,q,\alpha,\gamma) \in (0,\infty)$ is a constant.  By~\cite[Lemma~3.4]{KrumWica} and \eqref{keyest eqn30}, 
\begin{equation}\label{keyest eqn31}
	\int_{G_n(\mathbf{A}_{\rho,1/2}(\zeta))} (|z|^2 + |x|^2 \,\|\pi_S - \pi\|^2) \,d|T|(X,S) 
	\leq C \int_{\mathbf{A}_{\rho,1}(\zeta)} \op{dist}^2(X, \op{spt} \op{graph} \varphi) \,d\|T\|(X) ,
\end{equation}
where $C = C(n,m,q,\alpha,\gamma) \in (0,\infty)$ is a constant.  Using \eqref{keyest eqn31} and noting that $|\pi_{S^{\perp}}(x,0,0)| \leq |x| \,\|\pi_S - \pi\|$ for all $(X,S) \in G_n(\mathbf{A}_{\rho,1/2}(\zeta))$, it follows that \eqref{keyest eqn15} and \eqref{keyest eqn16} hold true.
\end{proof}

\begin{theorem}\label{radial noncon thm}
Let $\alpha = \ell_0/q_0 \neq 2$ where $\ell_0$ and $q_0$ are relatively prime positive integers with $q_0 \leq q$.  Let $\varsigma = \varsigma(n,m,q) \in (0,1)$ is as in Theorem~\ref{relharm separation3 thm}.  For each $\gamma \in (0,1)$ and $\sigma \in (0,1)$ there exists $\varepsilon_0 = \varepsilon_0(n,m,q,\alpha,\gamma,\sigma) \in (0,1)$ and $\beta_0 = \beta_0(n,m,q,\alpha,\gamma,\sigma) \in (0,1)$ the following holds true.  Let $\varphi$ and $T$ satisfy Hypothesis~$(\star)$ and Hypothesis~$(\star\star)$.  Let $\mathcal{B}_q$ be as in Definition~\ref{introduction B defn} (with $\mu = 1/q$).  If $0 \in \mathcal{B}_q$ such that $T$ is tangent to $P_0 = \mathbb{R}^n \times \{0\}$ at the origin and $\mathcal{N}_{T,{\rm Pl}}(0) \geq \alpha$, then   
\begin{equation}\label{radial noncon concl}
	\int_{\mathbf{C}_{\gamma}(0)} \frac{\op{dist}^2(X, \op{spt}\op{graph}\varphi)}{R^{n+2\alpha-\sigma}} \,d\|T\|(X) 
	\leq C E_{\varsigma/2}(T,\varphi,\mathbf{C}_1(0))^2 ,
\end{equation}
where $R(x,y,z) = |(x,y)|$ and $C = C(n,m,q,\alpha,\gamma,\sigma) \in (0,\infty)$ is a constant.
\end{theorem}

\begin{proof}
  Let $\varphi = \sum_{j=1}^J m_j \varphi_j$ for some $q_j$-valued components $\varphi_j$ with positive integer multiplicity $m_j$ (as in Definition~\ref{Phi-p_defn}).  Take any open ball $B \subset\subset \mathbb{R}^n \setminus \{0\} \times \mathbb{R}^{n-2}$ and let $\varphi_j = \sum_{l=1}^{q_j} \llbracket \varphi_{j,l} \rrbracket$ on $B$ where $\varphi_{j,l}$ are single-valued harmonic functions on $B$ (as in \eqref{varphi_localized}).  Let $\widehat{\chi} : [0,\infty) \rightarrow \mathbb{R}$ be a smooth function such that $0 \leq \widehat{\chi}(s) \leq 1$, $\widehat{\chi}(s) = 0$ if $0 \leq s \leq 1/6$, $\widehat{\chi}(s) = 1$ if $s \geq 1/3$, and $0 \leq \widehat{\chi\,}'(s) \leq 12$.  For each $j \in \{1,2,\ldots,J\}$ and $l \in \{1,2,\ldots,q_j\}$ define $\chi_{j,l} : B \times \mathbb{R}^m \rightarrow \mathbb{R}$ by 
\begin{equation}\label{radial noncon eqn1}
	\chi_{j,l}(x,y,z) = \prod_{(j',l') \neq (j,l)} \widehat{\chi}\bigg(\frac{|z - \varphi_{j',l'}(x,y)|}{|\varphi_{j,l}(x,y) - \varphi_{j',l'}(x,y)|}\bigg)
\end{equation}
for all $(x,y) \in B$ and $z \in \mathbb{R}^m$.  Note that by Theorem~\ref{graphical_thm}(A), $\op{sep} \varphi(x,y) > 0$ for each $(x,y) \in B$ and thus $\chi_{j,l}$ is well-defined and smooth on $B$.  Next define $\widetilde{d} : \mathbb{R}^{n+m} \rightarrow [0,\infty)$ by 
\begin{equation}\label{radial noncon eqn2}
	\widetilde{d}(x,y,z) = \frac{\sum_{j=1}^J \sum_{l=1}^{q_j} \chi_{j,l}(x,y,z) \,|z - \varphi_{j,l}(x,y)|}{\sum_{j=1}^J \sum_{l=1}^{q_j} \chi_{j,l}(x,y,z)}
\end{equation}
for all $(x,y,z) \in \mathbb{R}^{n+m}$ with $x \neq 0$, where $\varphi_{j,l}$ and $\chi_{j,l}$ are as in \eqref{varphi_localized} and \eqref{radial noncon eqn1} with $B = B_{|x|/2}(x,y)$, and $\widetilde{d}(0,y,z) = |z|$ for all $y \in \mathbb{R}^{n-2}$ and $z \in \mathbb{R}^m$.  Note that if $(x,y,z) \in \mathbb{R}^{n+m}$ with $x \neq 0$ and if $j_0 \in \{1,2,\ldots,J\}$ and $l_0 \in \{1,2,\ldots,q_{j_0}\}$ such that 
\begin{equation}\label{radial noncon eqn3}
	|z - \varphi_{j_0,l_0}(x,y)| \leq \min_{1 \leq j \leq J} \min_{1 \leq l \leq q_{j}} |z - \varphi_{j,l}(x,y)| , 
\end{equation}
then $\chi_{j_0,l_0}(x,y,z) = 1$.  Otherwise, $|z - \varphi_{j',l'}(x,y)| \leq \frac{1}{3} |\varphi_{j,l}(x,y) - \varphi_{j',l'}(x,y)|$ and thus by the triangle inequality $|z - \varphi_{j',l'}(x,y)| \leq \frac{1}{3} (|z - \varphi_{j,l}(x,y)| + |z - \varphi_{j',l'}(x,y)|)$ so that $|z - \varphi_{j',l'}(x,y)| \leq \tfrac{1}{2} |z - \varphi_{j,l}(x,y)|$, contradicting \eqref{radial noncon eqn3}.  As a result, $\widetilde{d}$ is well-defined and continuous on $\mathbb{R}^{n+m}$ and $\widetilde{d}$ is smooth on $\mathbb{R}^{n+m} \setminus \{0\} \times \mathbb{R}^{n-2+m}$.

We claim that 
\begin{gather}
	\label{radial noncon eqn4a} \widetilde{d}(\lambda x, \lambda y, \lambda^{\alpha} z) = \lambda^{\alpha} \widetilde{d}(x,y,z) \text{ for all } \lambda > 0, \\
	\label{radial noncon eqn4b} \op{dist}(X,\op{spt}\op{graph}\varphi) \leq \widetilde{d}(x,y,z) \leq 14 \op{dist}(X,\op{spt}\op{graph}\varphi) , \\
	\label{radial noncon eqn4c} |D_x \widetilde{d}(x,y,z)| \leq C \widehat{E} \,|x|^{\alpha-1}, \quad D_y \widetilde{d}(x,y,z) = 0, \quad |D_z \widetilde{d}(x,y,z)| \leq C, 
\end{gather}
for all $X = (x,y,z) \in \mathbb{R}^{n+m} \setminus \{0\} \times \mathbb{R}^{n-2+m}$, where $\widehat{E} = \widehat{E}(T,P_0,\mathbf{C}_1(0))$ and $C = C(n,m,q,\alpha) \in (0,\infty)$ is a constant.  \eqref{radial noncon eqn4a} follows immediately from \eqref{radial noncon eqn1}, \eqref{radial noncon eqn2}, and the homogeneity of $\varphi$.  Let $B \subset\subset \mathbb{R}^n \setminus \{0\} \times \mathbb{R}^{n-2}$ be any open ball and let $\varphi_{j,l}$ and $\chi_{j,l}$ be as in \eqref{varphi_localized} and \eqref{radial noncon eqn1}.  Note that $\widetilde{d}$ is given by \eqref{radial noncon eqn2} on $B \times \mathbb{R}^m$.  The first inequality in \eqref{radial noncon eqn4b} follows from \eqref{radial noncon eqn2} and $\op{dist}(X,\op{spt}\op{graph}\varphi) \leq |z - \varphi_{j,l}(x,y)|$ for all $j \in \{1,\ldots,J\}$ and $l \in \{1,\ldots,q_j\}$.  To see the second inequality in \eqref{radial noncon eqn4b}, let $(x,y,z) \in B \times \mathbb{R}^m$ and let $j_0 \in \{1,\ldots,J\}$ and $l_0 \in \{1,\ldots,q_{j_0}\}$ such that \eqref{radial noncon eqn3} holds true.  If $j \in \{1,\ldots,J\}$ and $l \in \{1,\ldots,q_j\}$ such that $\chi_{j,l}(x,y,z) > 0$, then $|z - \varphi_{j_0,l_0}(x,y)| > \tfrac{1}{6} |\varphi_{j,l}(x,y) - \varphi_{j_0,l_0}(x,y)|$ and thus by the triangle inequality 
\begin{align*}
	|z - \varphi_{j,l}(x,y)| \leq\,& |z - \varphi_{j_0,l_0}(x,y)| + |\varphi_{j,l}(x,y) - \varphi_{j_0,l_0}(x,y)| 
		\\ \leq\,& 7 |z - \varphi_{j_0,l_0}(x,y)| \leq 14 \op{dist}(X,\op{spt}\op{graph}\varphi) .
\end{align*}
Therefore, by \eqref{radial noncon eqn2}, the second inequality in \eqref{radial noncon eqn4b} holds true.  To see \eqref{radial noncon eqn4c}, first observe that by Remark~\ref{graphical rmk}(2), $|D_x \varphi_{j,l}(x,y)| \leq C \widehat{E} \,|x|^{\alpha-1}$ for all $(x,y) \in B$.  Also, $\varphi_{j,l}(x,y)$ is independent of $y$.  If $(x,y,z) \in \op{spt} D\chi_{j,l}$, then for each $(j',l') \neq (j,l)$ we have that $|z - \varphi_{j',l'}(x,y)| \leq \frac{1}{3} |\varphi_{j,l}(x,y) - \varphi_{j',l'}(x,y)|$ and thus by the triangle inequality $|z - \varphi_{j,l}(x,y)| \leq \frac{4}{3} |\varphi_{j,l}(x,y) - \varphi_{j',l'}(x,y)|$.  It follows by differentiating \eqref{radial noncon eqn1} that 
\begin{equation*}
	|D_x \chi_{j,l}(x,y,z)| \leq \frac{C \widehat{E} \,|x|^{\alpha-1}}{|z - \varphi_{j,l}(x,y)|} , \quad 
	D_y \chi_{j,l}(x,z) = 0, \quad 
	|D_z \chi_{j,l}(x,y,z)| \leq \frac{C}{|z - \varphi_{j,l}(x,y)|} 
\end{equation*}
for all $(x,y) \in B$ and $z \in \mathbb{R}$, where $C = C(n,m,q,\alpha) \in (0,\infty)$ is a constant.  \eqref{radial noncon eqn4c} then follows by differentiating \eqref{radial noncon eqn2}.

Let $\psi : [0,1] \rightarrow \mathbb{R}$ be a smooth function such that $\psi(s) = 1$ if $0 \leq s \leq \gamma$, $\psi(s) = 0$ if $(1+\gamma)/2 \leq s \leq 1$, and $-3/(1-\gamma) \leq \psi'(s) \leq 0$.  For each $\delta \in (0,1/2]$, let $\eta_{\delta} : [0,1] \rightarrow \mathbb{R}$ be a smooth function such that $\eta_{\delta}(s) = 0$ if $0 \leq s \leq \delta/2$, $\eta_{\delta}(s) = 1$ if $\delta \leq s \leq 1$, and $0 \leq \eta'_{\delta}(s) \leq 3/\delta$.  We compute 
\begin{align*}
	\op{div}_S(R^{-n+\sigma} (x,y,0)) =\,& R^{-n+\sigma} \big( n - \tfrac{1}{2} \|\pi_S - \pi\|^2 \big) - (n-\sigma) \,R^{-n+\sigma} |\nabla^S R|^2 
	\\=\,& -\tfrac{1}{2} R^{-n+\sigma} \|\pi_S - \pi\|^2 + \sigma \,R^{-n+\sigma} + (n-\sigma) \,R^{-n+\sigma} |\nabla^{S^{\perp}} R|^2
\end{align*}
for all $(x,y,z) \in \mathbb{R}^{n+m}$ and all $n$-dimensional plane $S$ (where $R = |(x,y,0)|$).  Thus by \eqref{first variation} with 
\begin{equation*}
	\zeta(X) = \psi(R)^2 \,\eta_{\delta}(R) \,R^{-n+\sigma-2\alpha} \widetilde{d}(X)^2 (x,y,0) 
\end{equation*}
on an open neighborhood of $\op{spt} T$ 
\begin{align}\label{radial noncon eqn5}
	&\sigma \int_{G_n(\mathbf{C}_1(0))} \psi(R)^2 \,\eta_{\delta}(R) \,R^{-n+\sigma-2\alpha} \,\widetilde{d}^2 \,d|T|(X,S)
	\\=\,& \int_{G_n(\mathbf{C}_1(0))} \psi(R)^2 \,\eta_{\delta}(R) \,R^{-n+\sigma-2\alpha} \,\widetilde{d}^2 
		\left( \tfrac{1}{2} \,\|\pi_S - \pi\|^2 - (n-\sigma) \,|\nabla^{S^{\perp}} R|^2 \right) d|T|(X,S) \nonumber 
	\\&- \int_{G_n(\mathbf{C}_1(0))} \psi(R)^2 \,\eta_{\delta}(R) \,R^{1-n+\sigma} \,\nabla^S R \cdot \nabla^S \bigg(\frac{\widetilde{d}^2}{R^{2\alpha}}\bigg) 
		\,d|T|(X,S) \nonumber
	\\&- 2 \int_{G_n(\mathbf{C}_1(0))} \psi(R) \,\psi'(R) \,\eta_{\delta}(R) \,R^{1-n+\sigma-2\alpha} \,\widetilde{d}^2 \,|\nabla^S R|^2 \,d|T|(X,S) \nonumber
	\\&- \int_{G_n(\mathbf{C}_1(0))} \psi(R)^2 \,\eta'_{\delta}(R) \,R^{1-n+\sigma-2\alpha} \,\widetilde{d}^2 \,|\nabla^S R|^2 \,d|T|(X,S) , \nonumber
\end{align}
where $\nabla^{S^{\perp}}$ is the component of the gradient orthogonal to $S$.  To bound the second term on the right-hand side of \eqref{radial noncon eqn5}, notice that by \eqref{radial noncon eqn4a} for each $(x,y,z) \in \mathbb{R}^{n+m}$ 
\begin{align*}
	\nabla \bigg(\frac{\widetilde{d}(x,y,z)}{R^{\alpha}}\bigg) 
	=\,& \nabla \bigg(\widetilde{d}\bigg(\frac{x}{R}, \frac{y}{R}, \frac{z}{R^{\alpha}}\bigg)\bigg)
	\\=\,& (D_{(x,y)} \widetilde{d})\bigg(\frac{x}{R}, \frac{y}{R}, \frac{z}{R^{\alpha}}\bigg) \cdot \frac{\pi_P - \nabla R \otimes \nabla R}{R} 
		+ (D_z \widetilde{d})\bigg(\frac{x}{R}, \frac{y}{R}, \frac{z}{R^{\alpha}}\bigg) \cdot \nabla\bigg(\frac{z}{R^{\alpha}}\bigg) 
		 \\=\,& D_{(x,y)} \widetilde{d}(x,y,z) \cdot \frac{\pi_P - \nabla R \otimes \nabla R}{R^{\alpha}} 
		+ D_z \widetilde{d}(x,y,z) \cdot \nabla\bigg(\frac{z}{R^{\alpha}}\bigg) . 
\end{align*}
Since $(\pi_P - \nabla R \otimes \nabla R) \,\nabla R = 0$ on $\mathbb{R}^{n+m}$, for each $(X,S) \in G_n(\mathbf{C}_1(0))$ 
\begin{align*}
	-\nabla^S R \cdot \nabla^S \bigg(\frac{\widetilde{d}}{R^{\alpha}}\bigg)  
	=\,& \frac{D_{(x,y)} \widetilde{d}(x,y,z)}{R^{\alpha}} \cdot \big( \pi_P(\nabla^{S^{\perp}} R) - |\nabla^{S^{\perp}} R|^2 \,\nabla R \big)
		\\&- \sum_{\kappa=1}^m D_{z^{\kappa}} \widetilde{d}(x,y,z) \,\nabla^S R \cdot \nabla^S \bigg(\frac{z^{\kappa}}{R^{\alpha}}\bigg) .
\end{align*}
By \eqref{radial noncon eqn4c}, for each $(X,S) \in G_n(\mathbf{C}_1(0))$ 
\begin{equation}\label{radial noncon eqn6}
	-\nabla^S R \cdot \nabla^S \bigg(\frac{\widetilde{d}}{R^{\alpha}}\bigg)  
	\leq \frac{C \widehat{E} \,\|\pi_S - \pi\|^2}{R} + C \sum_{\kappa=1}^m \left| \nabla^S R \cdot \nabla^S \bigg(\frac{z^{\kappa}}{R^{\alpha}}\bigg) \right| ,
\end{equation}
where $C = C(n,m,q,\alpha) \in (0,\infty)$ is a constant.  Since $\eta'_{\delta}(R) \geq 0$, the last terms on the right-hand side of \eqref{radial noncon eqn5} is $\leq 0$.  
Hence using \eqref{radial noncon eqn6} and Cauchy's inequality in \eqref{radial noncon eqn5} 
\begin{align}\label{radial noncon eqn7}
	&\int_{\mathbf{C}_1(0)} \psi(R)^2 \,\eta_{\delta}(R) \,R^{-n+\sigma-2\alpha} \,\widetilde{d}^2 \,d\|T\|(X)
	\\ \leq\,& C \int_{G_n(\mathbf{C}_1(0))} \psi(R)^2 \,R^{-n+\sigma-2\alpha} \,\widetilde{d}^2 \,\|\pi_S - \pi\|^2 \,d|T|(X,S) \nonumber
	\\&+ C \int_{G_n(\mathbf{C}_{\gamma}(0)) \cap \{\nabla^S R \neq 0\}} 
		\sum_{\kappa=1}^m R^{2-n} \left| \nabla^S R \cdot \nabla^S \bigg(\frac{z^{\kappa}}{R^{\alpha}}\bigg) \right|^2 \,d|T|(X,S) \nonumber
	\\&+ C \widehat{E}^2 \int_{G_n(\mathbf{C}_1(0))} R^{-n+\sigma} \psi(R)^2 \,\|\pi_S - \pi\|^4 \,d|T|(X,S) 
	+ C \int_{\mathbf{C}_1(0)} \widetilde{d}^2 \,d\|T\|(X) , \nonumber
\end{align}
where $C = C(n,m,q,\alpha,\gamma,\sigma) \in (0,\infty)$ is a constant.  

To bound the first and third terms on the right-hand side of \eqref{radial noncon eqn7}, for each $\rho \in (0,1]$ let $\widehat{E}_{\rho} = \widehat{E}(T,P_0,\mathbf{C}_{\rho}(0))$ (and recall that $\widehat{E} = \widehat{E}(T,P_0,\mathbf{C}_1(0))$).  By $\mathcal{N}_{T,{\rm Pl}}(0) \geq \alpha$, Remark~\ref{planar freq rmk}, and Theorem~\ref{frequency monotonicity thm}(e) we have that \eqref{keyest eqn2} holds true, i.e.~$\widehat{E}_{\rho} \leq C\rho^{\alpha-1} \widehat{E}$ for all $\rho \in (0,1]$, and in particular $\widehat{E}_{\rho} < C\varepsilon_0$ for all $\rho \in (0,1]$, where $C = C(n,m,q,\alpha) \in (0,\infty)$ are constants.  By Theorem~\ref{lip approx thm} and~\cite[Lemma~3.5 and Lemma~3.7]{KrumWica}, for each $\rho \in (0,\gamma]$ there exists a Lipschitz function $u_{\rho} : B_{(1+\gamma)\rho/2}(0) \rightarrow \mathcal{A}_q(\mathbb{R}^m)$ and a closed set $K_{\rho} \subset B_{(1+\gamma)\rho/2}(0)$ such that 
\begin{gather}\label{radial noncon eqn8}
	T \llcorner (K_{\rho} \times \mathbb{R}^m) = (\op{graph} u_{\rho}) \llcorner (K_{\rho} \times \mathbb{R}^m), \quad 
	\sup_{B_{(1+\gamma)\rho/2}(0)} |u_{\rho}| \leq C \widehat{E}_{\rho} \,\rho, \quad \op{Lip} u_{\rho} \leq C \widehat{E}_{\rho}^{\varsigma}, \\
	\mathcal{L}^n(B_{(1+\gamma)\rho/2}(0) \setminus K_{\rho}) + \|T\|((B_{(1+\gamma)\rho/2}(0) \setminus K_{\rho}) \times \mathbb{R}^m) 
		\leq C \widehat{E}_{\rho}^{2+\varsigma} \,\rho^n , 
	\nonumber
\end{gather}
where $\widehat{E}_{\rho} = \widehat{E}(T,P_0,\mathbf{C}_{\rho}(0))$ and $C = C(n,m,q,\gamma) \in (0,\infty)$ is a constant.  By the area formula, \eqref{radial noncon eqn8}, and Remark~\ref{graphical rmk}(2) 
\begin{align*}
	&\frac{1}{\rho^{n+2}} \int_{G_n(\mathbf{C}_{(1+\gamma)\rho/2}(0))} \widetilde{d}^2 \,\|\pi_S - \pi\|^2 \,d|T|(X,S)
	\\ \leq\,& \frac{1}{\rho^{n+2}} \int_{G_n(\mathbf{C}_{(1+\gamma)\rho/2}(0))} (\op{dist}^2(X,P_0) + |\varphi(x,y)|^2) \,\|\pi_S - \pi\|^2 \,d|T|(X,S) \nonumber
	\\ \leq\,& \frac{1}{\rho^{n+2}} \int_{B_{(1+\gamma)\rho/2}(0)} (|u_{\rho}|^2 + |\varphi|^2) \,|Du_{\rho}|^2 
		+ C \widehat{E}_{\rho}^{2+\varsigma} \nonumber
	\\ \leq\,& C \widehat{E}_{\rho}^{2+\varsigma} \nonumber
\end{align*} 
for all $\rho \in (0,1]$, where $C = C(n,m,q,\alpha,\gamma) \in (0,\infty)$ are constants.  Thus by \eqref{keyest eqn2} 
\begin{equation}\label{radial noncon eqn9}
	\frac{1}{\rho^{n+2}} \int_{G_n(\mathbf{C}_{(1+\gamma)\rho/2}(0))} \widetilde{d}^2 \,\|\pi_S - \pi\|^2 \,d|T|(X,S) 
		\leq C \rho^{(\alpha-1) (2+\varsigma)} \widehat{E}^{2+\varsigma} 
\end{equation} 
for all $\rho \in (0,1]$, where $C = C(n,m,q,\alpha,\gamma) \in (0,\infty)$ is a constant.  Multiplying both sides of \eqref{radial noncon eqn9} by $\rho^{\sigma-2(\alpha-1)}$, setting $\rho = 2^{-k}$, and summing over $k = 0,1,2,\ldots$ 
\begin{align}\label{radial noncon eqn10}
	&\int_{G_n(\mathbf{C}_1(0))} \psi(R)^2 \,R^{-n+\sigma-2\alpha} \,\widetilde{d}^2 \,\|\pi_S - \pi\|^2 \,d|T|(X,S)
	\\ \leq\,& C \sum_{k=1}^{\infty} 2^{k (n-\sigma+2\alpha)} \int_{G_n(\mathbf{C}_{(1+\gamma) 2^{-k}}(0))} \widetilde{d}^2 \,\|\pi_S - \pi\|^2 \,d|T|(X,S) \nonumber
	\\ \leq\,& C \sum_{k=1}^{\infty} 2^{-k (\sigma + (\alpha-1) \varsigma)} \widehat{E}^{2+\varsigma} 
	\leq C \widehat{E}^{2+\varsigma} , \nonumber 
\end{align} 
where $C = C(n,m,q,\alpha,\gamma) \in (0,\infty)$ are constants.  By similar reasoning, 
\begin{equation*}
	\frac{1}{\rho^n} \int_{G_n(\mathbf{C}_{(1+\gamma)\rho/2}(0))} \|\pi_S - \pi\|^4 \,d|T|(X,S)
	\leq C \widehat{E}_{\rho}^{2+\varsigma} 
	\leq C \rho^{(\alpha-1) (2+\varsigma)} \widehat{E}^{2+\varsigma} 
\end{equation*} 
for all $\rho \in (0,1]$ and thus 
\begin{align}\label{radial noncon eqn11}
	&\int_{G_n(\mathbf{C}_1(0))} R^{-n+\sigma} \psi(R)^2 \,\|\pi_S - \pi\|^4 \,d|T|(X,S)
	\\ \leq\,& C \sum_{k=1}^{\infty} 2^{k (n-\sigma)} \int_{G_n(\mathbf{C}_{(1+\gamma) 2^{-k}}(0))} \|\pi_S - \pi\|^4 \,d|T|(X,S) \nonumber
	\\ \leq\,& C \sum_{k=1}^{\infty} 2^{-k (\sigma + (\alpha-1)(2+\varsigma))} \widehat{E}^{2+\varsigma} 
	\leq C \widehat{E}^{2+\varsigma} , \nonumber 
\end{align} 
where $C = C(n,m,q,\alpha,\gamma) \in (0,\infty)$ are constants.  Using Theorem~\ref{keyest thm}, \eqref{radial noncon eqn10}, and \eqref{radial noncon eqn11} to bound the right-hand side of \eqref{radial noncon eqn7} and letting $\delta\rightarrow 0^+$ gives us \eqref{radial noncon concl}.
\end{proof}

\section{A priori estimates: Part~II}\label{sec:nonconcentration} 

Let $\varphi : \mathbb{R}^n \rightarrow \mathcal{A}_q(\mathbb{R}^m)$ be a non-zero, homogeneous degree $\alpha$, harmonic function which is translation invariant along $\{0\} \times \mathbb{R}^{n-2}$ (as in Definition~\ref{Phi-p_defn}).  Let $T$ be an $n$-dimensional locally area-minimizing rectifiable current of $\mathbf{C}_1(0)$ such that $(\partial T) \llcorner \mathbf{C}_1(0) = 0$ and $T$ is weakly close to $\mathbb{R}^n \times \{0\}$ with integer multiplicity $q$.  In Theorem~\ref{nonconcentration thm} below, we obtain estimates on the distance of a branch point $Z$ of planar frequency $\geq \alpha$ to the singular axis $\{0\}^2 \times \mathbb{R}^{n-2} \times \{0\}^m$ as well as integral estimates which rule out the concentration of $L^2$-excess of $T$ relative to $\varphi$ near $Z$.  These estimates and their proofs are inspired by the corresponding results in~\cite{Sim93} and the additional arguments from~\cite{Wic14} for establishing such estimates in the presence of higher multiplicities.

Moreover, we introduce a new estimate (contained in Theorem~\ref{nonconcentration thm}(A)) on the distance of the tangent plane $P_Z$ of $T$ at $Z$ to a fixed plane $\mathbb{R}^n \times \{0\}$.  Such an estimate is necessary to study branch points of minimal submanifolds, where harmonic approximations and frequency functions are computed relative to the tangent plane at each branch point.  Similar considerations arose in the authors' past work of~\cite{KrumWic1}; however, in~\cite{KrumWic1}, minimal submanifolds were expressed as graphs of $C^{1,\mu}$ two-valued functions $\widehat{u}^{(Z)}$ relative to the tangent plane at each branch point $Z$ and frequency functions were computed for the average-free part of $\widehat{u}^{(Z)}$.  Consequently, the analysis involved a different strategy utilizing an estimate~\cite[Lemma~3.3]{KrumWic1} comparing the average-free parts of $\widehat{u}^{(Z)}$ at different branch points.

Let $\varepsilon_0 = \varepsilon_0(n,m,q) \in (0,1)$ and $\beta_0 = \beta_0(n,m,q) \in (0,1)$ be suitably small constants.  Let $\varphi \in \Phi_{\alpha,q}$ and $T$ satisfy Hypothesis~$(\star)$ and Hypothesis~$(\star\star)$.  Let $\mathcal{B}_q$ be as in Definition~\ref{introduction B defn} (with $\mu = 1/q$).  Let $\varphi = \sum_{j=1}^J m_j \varphi_j$ with $q_j$-valued components $\varphi_j$ with multiplicity $m_j$ (as in Definition~\ref{Phi-p_defn}).  Let $Z = (\xi,\zeta,\chi) \in \mathcal{B}_q$ such that $\mathcal{N}_{T,{\rm Pl}}(Z) \geq \alpha$ and let $P_Z$ be the tangent plane to $T$ at $Z$.  By Lemma~\ref{uniform branch set lemma} (in particular \eqref{uniform branch set concl3}), we can express $P_Z$ as 
\begin{equation}\label{noncon tangent plane}
	P_Z = \{ (x,y,A(x,y)) : (x,y) \in \mathbb{R}^n \} 
\end{equation}
where $A : \mathbb{R}^n \rightarrow \mathbb{R}^m$ is a linear transformation such that $\|A\| \leq C \op{dist}_{\mathcal H}(P_Z \cap \mathbf{B}_1(0), P_0 \cap \mathbf{B}_1(0))$, where $P_0 = \mathbb{R}^n \times \{0\}$ and $C = C(n,m) \in (0,\infty)$ is a constant.  Note that by Lemma~\ref{uniform branch set lemma} we have that $\|A\| \leq C(n,m,q) \,\widehat{E}(T,P_0,\mathbf{C}_1(0))$; we will obtain better bounds on $\|A\|$ in Theorem~\ref{nonconcentration thm} and its proof below.  It follows that there exists a rotation $\Gamma$ of $\mathbb{R}^{n+m}$ such that $\Gamma(\mathbb{R}^n \times \{0\}) = P_Z$ and 
\begin{equation}\label{noncon rotate} 
	\Gamma = I + \left[\begin{array}{rr} 0 & -A^T \\ A & 0 \end{array}\right] + \mathcal{R} , 
\end{equation}
where $I$ is the identity on $\mathbb{R}^{n+m}$, $A^T$ is the transpose of $A$, and $|\mathcal{R}| \leq C(n,m,q) \,\|A\|^2$.  In  Theorem~\ref{nonconcentration thm} and throughout this section, we define $\widehat{\varphi} : \mathbb{R}^n \rightarrow \mathcal{A}_q(\mathbb{R}^m)$ is defined by 
\begin{align}\label{noncon rotvarphi} 
	\widehat{\varphi}(x,y) = \sum_{j=1}^J \sum_{l=1}^{q_j} m_j \llbracket \chi - A(x-\xi,y-\zeta) + \varphi_{j,l}(x,y) - D_x \varphi_{j,l}(x,y) \cdot \xi \rrbracket  
\end{align} 
for each $(x,y) \in \mathbb{R}^n$, where $\varphi_{j,l}$ are locally defined by \eqref{varphi_localized}.

\begin{theorem}\label{nonconcentration thm} 
Let $\alpha = \ell_0/q_0 \neq 2$ where $\ell_0$ and $q_0$ are relatively prime positive integers with $q_0 \leq q$.  Let $\varsigma = \varsigma(n,m,q) \in (0,1)$ be as in Theorem~\ref{relharm separation3 thm}.  Let $p \in \{1,2,\ldots,\lfloor q/q_0 \rfloor\}$, $0 < \tau < \gamma < 1$, and $\sigma \in (0,1)$.  There exists $\varepsilon_0 = \varepsilon_0(n,m,q,p,\alpha,\gamma,\tau,\sigma) \in (0,1)$ and $\beta_0 = \beta_0(n,m,q,p,\alpha,\gamma,\tau,\sigma) \in (0,1)$ such that the following holds true.  Let $\varphi \in \Phi_{\alpha,q,p}$ and $T$ satisfy Hypothesis~$(\star)$ and Hypothesis~$(\star\star)$.  Let $\mathcal{B}_q$ be as in Definition~\ref{introduction B defn} (with $\mu = 1/q$).  If $Z = (\xi,\zeta,\chi) \in \mathcal{B}_q$ such that $T$ is tangent to $P_Z$ at $Z$ and $\mathcal{N}_{T,{\rm Pl}}(Z) \geq \alpha$, then 
\begin{align*}
	&(A) \quad \widehat{E} |\xi| + |\chi| + \op{dist}_{\mathcal H}(P_Z \cap \mathbf{B}_1(0), P_0 \cap \mathbf{B}_1(0)) 
		\leq C E_{\varsigma/2}(T,\varphi,\mathbf{C}_1(0)), \\ 
	&(B) \quad \int_{\mathbf{C}_1(0)} \op{dist}^2(X, Z+\Gamma(\op{spt} \op{graph} \varphi)) \,d\|T\|(X) \leq C E_{\varsigma/2}(T,\varphi,\mathbf{C}_1(0))^2 , 
\end{align*}
where $\widehat{E} = \widehat{E}(T,P_0,\mathbf{C}_1(0)$ and $C = C(n,m,q,p,\alpha) \in (0,\infty)$ is a constant.  Moreover, 
\begin{align*}
	&(C) \quad \int_{\mathbf{C}_{\gamma}(0)} \frac{\op{dist}^2(X, \op{spt} \op{graph} \varphi)}{|X-Z|^{n-\sigma}} \,d\|T\|(X)
		\leq C E_{\varsigma/2}(T,\varphi,\mathbf{C}_1(0))^2 , \\
	&(D) \quad \int_{\mathbf{C}_{\gamma}(0) \cap \{r > \tau\}} \frac{\op{dist}^2(X, \op{spt} \op{graph} \widehat{\varphi})}{|X-Z|^{n+2\alpha-\sigma}} \,d\|T\|(X) 
		\leq C E_{\varsigma/2}(T,\varphi,\mathbf{C}_1(0))^2 , 
\end{align*}
where $\widehat{\varphi}$ is as in \eqref{noncon rotvarphi} and $C = C(n,m,q,p,\alpha,\gamma,\sigma) \in (0,\infty)$ is a constant.
\end{theorem}

We will prove Theorem~\ref{nonconcentration thm} by induction on $p$, assuming that: 
\begin{enumerate}[topsep=0mm]
	\item[(H4)]  whenever $\widetilde{p} \in \{1,2,\ldots,p-1\}$, Theorem~\ref{nonconcentration thm} holds true with $\widetilde{p}$ in place of $p$.
\end{enumerate}
First we have a few preliminary lemmas, Lemma~\ref{nonconcentration3 lemma}, Lemma~\ref{nonconcentration4 lemma}, and Lemma~\ref{nonconcentration5 lemma}.

\begin{lemma}\label{nonconcentration3 lemma}
Let $\alpha = \ell_0/q_0 \neq 2$ where $\ell_0$ and $q_0$ are relatively prime positive integers with $q_0 \leq q$.  Let $\varsigma = \varsigma(n,m,q) \in (0,1)$ be as in Theorem~\ref{relharm separation3 thm}.  For all $\delta \in (0,1)$ there exists $\varepsilon_0 = \varepsilon_0(n,m,q,\alpha,\delta) \in (0,1)$ and $\beta_0 = \beta_0(n,m,q,\alpha,\delta) \in (0,1)$ such that the following holds true.  Let $p \in \{1,2,\ldots,\lfloor q/q_0 \rfloor\}$, $\varphi \in \Phi_{\alpha,q,p}$, and $T$ satisfy Hypothesis~$(\star)$ and 
\begin{equation}\label{noncon3 hyp} 
	Q_{\varsigma/2}(T,\varphi,\mathbf{C}_1(0)) \leq \beta_0 \widehat{E}(T,P_0,\mathbf{C}_1(0)) .
\end{equation}
Let $\mathcal{B}_q = \mathcal{B}_q(T)$ be as in Definition~\ref{introduction B defn} (with $\mu = 1/q$).  If $Z = (\xi,\zeta,\chi) \in \mathcal{B}_q$ such that $T$ is tangent to $P_Z$ at $Z$ and $\mathcal{N}_{T,{\rm Pl}}(Z) \geq \alpha$, then 
\begin{equation}\label{noncon3 concl} 
	\widehat{E} |\xi| + |\chi| + \op{dist}_{\mathcal H}(P_Z \cap \mathbf{B}_1(0), P_0 \cap \mathbf{B}_1(0)) \leq \delta \widehat{E} , 
\end{equation}
where $\widehat{E} = \widehat{E}(T,P_0,\mathbf{C}_1(0))$.
\end{lemma}

\begin{proof}
We proceed by contradiction.  Fix $\delta > 0$ and $p \in \{1,2,\ldots,\lfloor q/q_0 \rfloor\}$.  Suppose that for $\nu = 1,2,3,\ldots$ there exists $\varepsilon_{\nu}\rightarrow 0^+$, $\beta_{\nu}\rightarrow 0^+$, $\varphi^{(\nu)} \in \Phi_{\alpha,q,p}$, an $n$-dimensional locally area-minimizing rectifiable current $T^{(\nu)}$ of $\mathbf{C}_1(0)$, $Z_{\nu} = (\xi_{\nu},\zeta_{\nu},\chi_{\nu}) \in \mathcal{B}_q(T^{(\nu)})$, and a (unique) tangent plane $P_{\nu}$ of $T^{(\nu)}$ at $Z_{\nu}$ such that Hypothesis~$(\star)$ and \eqref{noncon3 hyp} hold true with $\varepsilon_{\nu},\beta_{\nu},\varphi^{(\nu)},T^{(\nu)}$ in place of $\varepsilon_0,\beta_0,\varphi,T$ and $\mathcal{N}_{T^{(\nu)},{\rm Pl}}(Z_{\nu}) \geq \alpha$ but 
\begin{equation}\label{noncon3 eqn1} 
	\widehat{E}_{\nu} |\xi_{\nu}| + |\chi_{\nu}| + \op{dist}_{\mathcal H}(P_{\nu} \cap \mathbf{B}_1(0), P_0 \cap \mathbf{B}_1(0)) 
	\geq \delta \widehat{E}_{\nu} , 
\end{equation}
where $\widehat{E}_{\nu} = \widehat{E}(T^{(\nu)},P_0,\mathbf{C}_1(0))$.

In view of Hypothesis~$(\star)$, we can apply the blow procedure of~\cite[Subsection~6.1]{KrumWica} to construct a blow-up $w$ of $T^{(\nu)}$ relative to $P_0$ by excess $\widehat{E}_{\nu}$.  By Remark~\ref{graphical rmk}(2), after passing to a subsequence let $\varphi^{(\nu)}/\widehat{E}_{\nu} \rightarrow \psi$ uniformly on compact subsets of $\mathbb{R}^n$ so that $\psi$ is the blow-up of $\varphi^{(\nu)}$ relative to $P_0$ by by excess $\widehat{E}_{\nu}$ and $\psi$ is non-zero, homogeneous degree $\alpha$, and translation invariant along $\{0\} \times \mathbb{R}^{n-2}$.  By \eqref{noncon3 hyp}, we can argue as we did in \eqref{graphical annuli eqn2} to deduce that $\op{spt} w(x,y) = \op{spt} \psi(x,y)$ for all $(x,y) \in B_1(0)$.  In particular, $w$ is non-zero, homogeneous degree $\alpha$, and translation invariant along $\{0\} \times \mathbb{R}^{n-2}$. 

Let $(\xi_{\nu},\zeta_{\nu}) \rightarrow (\xi,\zeta)$ in $\mathbb{R}^n$ as $\nu\rightarrow\infty$.  Let $w_a(x,y) = \frac{1}{q} \sum_{l=1}^q w_l(x,y)$ for each $(x,y) \in B_1(0)$, where $w(x,y) = \sum_{l=1}^q \llbracket w_l(x,y) \rrbracket$.  By applying the Hardt-Simon inequality~\cite[Lemma~6.4]{KrumWica}, $\chi_{\nu}/\widehat{E}_{\nu} \rightarrow w_a(\xi,\zeta)$ as $\nu\rightarrow\infty$ and $w(\xi,\zeta) = q \llbracket w_a(\xi,\zeta) \rrbracket$.  Recall that by Lemma~\ref{uniform branch set lemma}, we can express $P_{\nu} = \{ (x,y,z) : z = A_{\nu}(x,y) \}$ for some linear function $A_{\nu} : \mathbb{R}^n \rightarrow \mathbb{R}^m$ such that $\|A_{\nu}\|_{L^{\infty}(B_1(0))} \leq C(n,m,q,\alpha) \,\widehat{E}_{\nu}$ (as in \eqref{noncon tangent plane}).  By~\cite[Corollary~6.6]{KrumWica}, $A_{\nu}/\widehat{E}_{\nu} \rightarrow Dw_a(\xi)$ and $\mathcal{N}_{w - \ell}(\xi,\zeta) \geq \alpha$, where $\ell(x,y) = w_a(\xi,\zeta) + Dw_a(\xi,\zeta) \cdot (x-\xi,y-\zeta)$ for all $(x,y) \in \mathbb{R}^n$.  Since $w$ is non-zero and homogeneous degree $\alpha$ where $\alpha > 1$ and $\alpha \neq 2$, we can argue as we did in the proof of Theorem~\ref{B prime stratification thm} (with $w$ in place of $\varphi$) to obtain $w_a(\xi,\zeta) = 0$, $Dw_a(\xi,\zeta) = 0$, and $\mathcal{N}_w(\xi,\zeta) \geq \alpha$ so that $(\xi,\zeta) \in \op{spine} w$.  Since $\op{spine} w = \{0\} \times \mathbb{R}^{n-2}$, $\xi = 0$.  Therefore, $\xi_{\nu} \rightarrow 0$, $\chi_{\nu}/\widehat{E}_{\nu} \rightarrow 0$, and $A_{\nu}/\widehat{E}_{\nu} \rightarrow 0$ as $\nu\rightarrow\infty$, contradicting \eqref{noncon3 eqn1}. 
\end{proof}

\begin{lemma}\label{nonconcentration4 lemma} 
Let $\alpha = \ell_0/q_0 \neq 2$ where $\ell_0$ and $q_0$ are relatively prime positive integers with $q_0 \leq q$.  Assume that $p \in \{2,3,\ldots,\lfloor q/q_0 \rfloor\}$ such that (H4) holds true.  Let $\varsigma = \varsigma(n,m,q) \in (0,1)$ is as in Theorem~\ref{relharm separation3 thm}.  For all $\delta \in (0,1)$ there exists $\varepsilon_0 = \varepsilon_0(n,m,q,p,\alpha,\delta) \in (0,1)$, $\beta_0 = \beta_0(n,m,q,p,\alpha,\delta) \in (0,1)$, and $\gamma_0 = \gamma_0(n,m,q,p,\alpha,\delta) \in (0,1)$ such that the following holds true.  Let $\varphi \in \Phi_{\alpha,q,p}$ and $T$ satisfy Hypothesis~$(\star)$ and Hypothesis~$(\star\star)$.  Let $\mathcal{B}_q = \mathcal{B}_q(T)$ be as in Definition~\ref{introduction B defn} (with $\mu = 1/q$).  Suppose there exists $s \in \{1,2,\ldots,p-1\}$ such that 
\begin{equation}\label{noncon4 hyp} 
	\inf_{\varphi' \in \bigcup_{p'=0}^s \Phi_{\alpha,q,p'}} Q_{\varsigma/2}(T,\varphi',\mathbf{C}_1(0)) 
		\leq \gamma_0 \inf_{\varphi' \in \bigcup_{p'=0}^{s-1} \Phi_{\alpha,q,p'}} Q_{\varsigma/2}(T,\varphi',\mathbf{C}_1(0)) .
\end{equation}
If $Z = (\xi,\zeta,\chi) \in \mathcal{B}_q$ such that $T$ is tangent to $P_Z$ at $Z$ and $\mathcal{N}_{T,{\rm Pl}}(Z) \geq \alpha$, then 
\begin{equation*}%\label{noncon4 concl} 
	\widehat{E} |\xi| + |\chi| + \op{dist}_{\mathcal H}(P_Z \cap \mathbf{B}_1(0), P_0 \cap \mathbf{B}_1(0)) 
	\leq \delta \inf_{\varphi' \in \bigcup_{p'=0}^s \Phi_{\alpha,q,p'}} Q_{\varsigma/2}(T,\varphi',\mathbf{C}_1(0)) ,
\end{equation*}
where $\widehat{E} = \widehat{E}(T,P_0,\mathbf{C}_1(0))$.
\end{lemma}

\begin{proof}
We proceed by contradiction.  Fix $p \in \{2,3,\ldots,\lfloor q/q_0 \rfloor\}$ such that (H4) holds true.  Fix $\delta > 0$ and $s \in \{1,2,\ldots,p-1\}$.  Suppose that for $\nu = 1,2,3,\ldots$ there exists $\varepsilon_{\nu}\rightarrow 0^+$, $\beta_{\nu}\rightarrow 0^+$, $\gamma_{\nu}\rightarrow 0^+$, $\varphi^{(\nu)} \in \Phi_{\alpha,q,p}$, an $n$-dimensional locally area-minimizing rectifiable current $T^{(\nu)}$ of $\mathbf{C}_1(0)$, $Z_{\nu} = (\xi_{\nu},\zeta_{\nu},\chi_{\nu}) \in \mathcal{B}_q(T^{(\nu)})$, and a (unique) tangent plane $P_{\nu}$ of $T^{(\nu)}$ at $Z_{\nu}$ such that Hypothesis~$(\star)$, Hypothesis~$(\star\star)$, and \eqref{noncon4 hyp} hold true with $\varepsilon_{\nu},\beta_{\nu},\gamma_{\nu},\varphi^{(\nu)},T^{(\nu)}$ in place of $\varepsilon_0,\beta_0,\gamma_0,\varphi,T$ and $\mathcal{N}_{T^{(\nu)},{\rm Pl}}(Z_{\nu}) \geq \alpha$ but 
\begin{equation}\label{noncon4 eqn1} 
	\widehat{E}_{\nu} |\xi_{\nu}| + |\chi_{\nu}| + \op{dist}_{\mathcal H}(P_{\nu} \cap \mathbf{B}_1(0), P_0 \cap \mathbf{B}_1(0)) 
	\geq \delta \inf_{\varphi' \in \bigcup_{p'=0}^s \Phi_{\alpha,q,p'}} Q_{\varsigma/2}(T^{(\nu)},\varphi',\mathbf{C}_1(0)) , 
\end{equation}
where $\widehat{E}_{\nu} = \widehat{E}(T^{(\nu)},P_0,\mathbf{C}_1(0))$.  

Let $\widetilde{\varphi}^{(\nu)} \in \Phi_{\alpha,q,s}$ such that 
\begin{equation*} 
	Q_{\varsigma/2}(T,\widetilde{\varphi}^{(\nu)},\mathbf{C}_1(0)) \leq 2 \inf_{\varphi' \in \Phi_{\alpha,q,s}} Q_{\varsigma/2}(T,\varphi',\mathbf{C}_1(0)) .
\end{equation*}
Let $\widetilde{\varphi}^{(\nu)} = \sum_{j=1}^{\widetilde{J}} \widetilde{m}_j \widetilde{\varphi}^{(\nu)}_j$ with $\widetilde{q}_j$-valued components $\widetilde{\varphi}^{(\nu)}_j$ with multiplicity $\widetilde{m}_j$ (as in Definition~\ref{Phi-p_defn}).  After passing to a subsequence, we may assume that $\widetilde{J}, \widetilde{q}_j, \widetilde{m}_j$ are independent of $\nu$.  Set $\widetilde{Q}_{\nu} = Q_{\varsigma/2}(T^{(\nu)},\widetilde{\varphi}^{(\nu)},\mathbf{C}_1(0))$.  In view of Hypothesis~$(\star)$ and \eqref{noncon4 hyp}, following the blow procedure in Subsection~\ref{sec:blowup subsec}, we obtain a blow-up $w = (w_j)$ of $T^{(\nu)}$ relative to $\widetilde{\varphi}^{(\nu)}$ by excess $\widetilde{Q}_{\nu}$, and we obtain a blow-up $\psi = (\psi_j)$ of $\varphi^{(\nu)}$ relative to $\widetilde{\varphi}^{(\nu)}$ by excess $\widetilde{Q}_{\nu}$.  By Hypothesis~$(\star\star)$, arguing like we did in \eqref{graphical annuli eqn2} it follows that $\op{spt} w_j = \op{spt} \psi_j$ on $\op{graph} \widetilde{\varphi}^{(\infty)}_j |_{B_1(0)}$ (where $\widetilde{\varphi}^{(\infty)}_j$ is as in Subsection~\ref{sec:blowup subsec} with $\widetilde{\varphi}^{(\nu)}_j$ in place of $\varphi^{(\nu)}_j$).  In particular, $w$ is homogeneous degree $\alpha$ and translation invariant along $\{0\} \times \mathbb{R}^{n-2}$.  

Let $\widetilde{\varphi}^{(\nu)}_j/\widehat{E}_{\nu} \rightarrow \widetilde{\varphi}^{(0)}_j$ uniformly on compact subsets of $\mathbb{R}^n$ as $\nu\rightarrow\infty$.  By Remark~\ref{graphical rmk}(2), $\widetilde{\varphi}^{(0)}_j$ is not identically zero for some $j \in \{1,2,\ldots,J\}$.  By the assumption (H4) above and Theorem~\ref{nonconcentration thm}(A)
\begin{equation*}%\label{noncon4 eqn2} 
	\widehat{E}_{\nu} |\xi_{\nu}| + |\chi_{\nu}| + \op{dist}_{\mathcal H}(P_{\nu} \cap \mathbf{B}_1(0), P_0 \cap \mathbf{B}_1(0)) \leq C \widetilde{Q}_{\nu}
\end{equation*}
for some constant $C = C(n,m,q,p,\alpha) \in (0,\infty)$.  Thus after passing to a subsequence $\widehat{E}_{\nu} \xi_{\nu}/\widetilde{Q}_{\nu} \rightarrow \kappa$ in $\mathbb{R}^2$, $\zeta_{\nu} \rightarrow \zeta$ in $\mathbb{R}^{n-2}$, and $\chi_{\nu}/\widetilde{Q}_{\nu} \rightarrow \lambda$ in $\mathbb{R}^m$ as $\nu\rightarrow\infty$.  Additionally, $P^{(\nu)} = \{ (x,y,z) : z = A_{\nu}(x,y) \}$ for some linear function $A_{\nu} : \mathbb{R}^n \rightarrow \mathbb{R}^m$ such that $\|A_{\nu}\|_{L^{\infty}(B_1(0))} \leq C(n,m,q,\alpha) \,\widetilde{Q}_{\nu}$ (as in \eqref{noncon tangent plane}).  After passing to a further subsequence let $A_{\nu}/\widetilde{Q}_{\nu} \rightarrow \Lambda$ as linear transformations from $\mathbb{R}^n$ to $\mathbb{R}^m$ as $\nu\rightarrow\infty$.  By (H4) and Theorem~\ref{nonconcentration thm}(D), for every $\tau > 0$ and $\nu$ sufficiently large (depending on $\tau$) 
\begin{equation*}%\label{noncon4 eqn3} 
	\int_{B_{3/4}(0) \cap \{r > \tau\}} \sum_{j=1}^J \sum_{l=1}^{q_j} \sum_{h=1}^{m_j} 
		\frac{|v^{(\nu)}_{j,l,h}(x,y) - \chi_{\nu} - A_{\nu}(x-\xi_{\nu},y-\zeta_{\nu}) - D\widetilde{\varphi}^{(\nu)}_{j,l}(x,y) \cdot \xi_{\nu}|^2}{
		|(x,y) - (\xi_{\nu},\zeta_{\nu})|^{n+2\alpha-\sigma}} \,dx \,dy
	\leq C \widetilde{Q}_{\nu}^2 , 
\end{equation*}
where $v^{(\nu)}_j$ is as in Theorem~\ref{graphical_thm} with $T^{(\nu)},\varphi^{(\nu)}$ in place of $T,\varphi$, $v^{(\nu)}_{j,l,h}$ is as in \eqref{v_localized} with $\varphi^{(\nu)}_j, v^{(\nu)}_j$ in place of $\varphi_j,v_j$, and $C = C(n,m,q,\alpha,\sigma) \in (0,\infty)$ is a constant.  Dividing both sides 
by $\widetilde{Q}_{\nu}^2$ and letting $\nu\rightarrow\infty$ 
\begin{equation*} 
	\int_{B_{3/4}(0)} \sum_{j=1}^J \sum_{l=1}^{q_j} \sum_{h=1}^{m_j} 
		\frac{|w_{j,l,h}(x,y) - \lambda - \Lambda (x,y-\zeta) - D\widetilde{\varphi}^{(0)}_{j,l}(x,y) \cdot \kappa|^2}{
		|(x,y) - (0,\zeta)|^{n+2\alpha-\sigma}} \,dx \,dy \leq C , 
\end{equation*}
where $w_j(x,y,\varphi_{j,l}(x,y)) = \sum_{l=1}^{q_j} \llbracket w_{j,l,h}(x,y) \rrbracket$ for each $(x,y) \in B_{3/4}(0)$ and $C = C(n,m,q,\alpha,\sigma) \in (0,\infty)$ is a constant.  Recall that $w_j(x,y,\varphi_{j,l}(x,y))$ is homogeneous degree $\alpha$ and translation invariant along $\{0\} \times \mathbb{R}^{n-2}$.  Since $\alpha > 1$ and $\alpha \neq 2$, the functions $\lambda$, $\Lambda (x,y)$, and $\frac{1}{q_j} \sum_{l=1}^{q_j} D\widetilde{\varphi}^{(0)}_{j,l}(x,y) \cdot \kappa$ are mutually $L^2(\mathbb{S}^{n-1})$-orthogonal.  (In fact, $\widetilde{\varphi}^{(0)}_j$ is a homogeneous degree $\alpha$, harmonic, single-valued function if $q_j = 1$ and is an average-free $q_j$-valued function if $q_j \geq 2$.)  Hence 
\begin{equation*} 
	\sum_{j=1}^J m_j \int_{B_{1/4}(0,\zeta)} \sum_{l=1}^{q_j} 
		\frac{|\lambda|^2 + |\Lambda (x,y-\zeta)|^2 + |D\widetilde{\varphi}^{(0)}_{j,l}(x,y) \cdot \kappa|^2}{|(x,y) - (0,\zeta)|^{n+2\alpha-\sigma}} \,dx \,dy < \infty . 
\end{equation*}
Therefore, $\kappa = 0$, $\lambda = 0$, and $\Lambda = 0$.  That is, $\xi_{\nu} \rightarrow 0$, $\chi_{\nu}/\widetilde{Q}_{\nu} \rightarrow 0$, and $A_{\nu}/\widetilde{Q}_{\nu} \rightarrow 0$ as $\nu\rightarrow\infty$, contradicting \eqref{noncon4 eqn1}. 
\end{proof}

\begin{lemma}\label{nonconcentration5 lemma}
Let $\alpha = \ell_0/q_0 \neq 2$ where $\ell_0$ and $q_0$ are relatively prime positive integers with $q_0 \leq q$.  Assume that $p \in \{1,2,\ldots,\lfloor q/q_0 \rfloor\}$ such that either $p = 1$ or $p > 1$ and (H4) holds true.  Let $\varsigma = \varsigma(n,m,q) \in (0,1)$ be as in Theorem~\ref{relharm separation3 thm}.  For all $\delta \in (0,1)$ there exists $\varepsilon_0 = \varepsilon_0(n,m,q,p,\alpha,\gamma,\delta) \in (0,1)$ and $\beta_0 = \beta_0(n,m,q,p,\alpha,\gamma,\delta) \in (0,1)$ such that the following holds true.  Let $\varphi \in \Phi_{\alpha,q,p}$ and $T$ satisfy Hypothesis~$(\star)$ and Hypothesis~$(\star\star)$.  Let $\mathcal{B}_q$ be as in Definition~\ref{introduction B defn} (with $\mu = 1/q$).  If $Z = (\xi,\zeta,\chi) \in \mathcal{B}_q$ with $\mathcal{N}_{T,{\rm Pl}}(Z) \geq \alpha$, then 
\begin{equation}\label{noncon5 concl} 
	\widehat{E} |\xi| + |\chi| + \op{dist}_{\mathcal H}(P_Z \cap \mathbf{B}_1(0), P_0 \cap \mathbf{B}_1(0)) 
	\leq \delta \inf_{\varphi' \in \bigcup_{p'=0}^{p-1} \Phi_{\alpha,q,p'}} Q_{\varsigma/2}(T,\varphi',\mathbf{C}_1(0)) , 
\end{equation}
where $\widehat{E} = \widehat{E}(T,P_0,\mathbf{C}_1(0))$.  
\end{lemma}

\begin{proof}
Lemma~\ref{nonconcentration5 lemma} follows from Lemma~\ref{nonconcentration3 lemma} and Lemma~\ref{nonconcentration4 lemma} by arguing as in~\cite[Lemma~8.4]{KrumWic2} (for which Lemma~8.4 follows from Lemma~8.5).
\end{proof}

\begin{remark}\label{noncon triangle rmk}{\rm 
We make the following observations which we shall rely on in the proof of Theorem~\ref{nonconcentration thm}.  Let $p \in \{1,2,\ldots,\lfloor q/q_0 \rfloor\}$ such that either $p = 1$ or $p > 1$ and (H4) holds true.  Suppose that for small constants $\varepsilon_0,\beta_0 \in (0,1)$, let $\varphi \in \Phi_{\alpha,q,p}$, and $T$ satisfy Hypothesis~$(\star)$ and Hypothesis~$(\star\star)$.  Let $\mathcal{B}_q$ be as in Definition~\ref{introduction B defn} (with $\mu = 1/q$).   
Let $Z = (\xi,\zeta,\chi) \in \mathcal{B}_q$ such that $\mathcal{N}_{T,{\rm Pl}}(Z) \geq \alpha$ and let $P_Z$ be the tangent plane to $T$ at $Z$.  Let $\widehat{E} = \widehat{E}(T,P_0,\mathbf{C}_1(0))$.  Let $A : \mathbb{R}^n \rightarrow \mathbb{R}^m$ is a linear transformation such that $P_Z$ is given by \eqref{noncon tangent plane} and $\|A\| \leq C(n,m,q) \,\widehat{E}$.  Let $\Gamma$ be a rotation of $\mathbb{R}^{n+m}$ such that $\Gamma(\mathbb{R}^n \times \{0\}) = P_Z$ and \eqref{noncon rotate} holds true. 

\noindent (1)  \textit{For each $\gamma \in (0,1)$ there exists $\varepsilon_0 = \varepsilon_0(n,m,q,\alpha,\gamma) \in (0,1)$ and $\beta_0 = \beta_0(n,m,q,\alpha,\gamma) \in (0,1)$ such that for each $X = (x,y,z) \in \op{spt} T \cap \mathbf{C}_{\gamma}(0)$ 
\begin{equation}\label{noncon triangle1} 
	\op{dist}(X, \op{spt}\op{graph}\widehat{\varphi}) \leq 2 \op{dist}(X, Z + \Gamma(\op{spt}\op{graph}\varphi)) 
		+ C \widehat{E} \,\|A\| + C \widehat{E} \,|x|^{\alpha-2} \,|\xi|^2 , 
\end{equation}
where $\widehat{\varphi}$ is as in \eqref{noncon rotvarphi}, $\widehat{E} = \widehat{E}(T,P_0,\mathbf{C}_1(0)$, and $C = C(n,m,q,\alpha) \in (0,\infty)$ is a constant.}  

To see this, note that by~\cite[Lemma~3.5]{KrumWica}, $|z - \chi| \leq C(n,m,\gamma) \,\widehat{E}$.  Thus by \eqref{noncon rotate} 
\begin{equation}\label{noncon triangle eqn1} 
	\Gamma^{-1}(X-Z) = (x-\xi, y-\zeta, z-\chi - A(x-\xi,y-\zeta)) + \mathcal{R}(X) ,
\end{equation}
where $|\mathcal{R}(X)| \leq C(n,m,q,\gamma) \,\widehat{E} \,\|A\|$.  Let $\varphi_{j,l}$ be as in \eqref{varphi_localized} with $B = B_{|x|/2}(x,y)$ and let $\widehat{\varphi}_{j,l}(x,y) = \chi - A(x-\xi,y-\zeta) + \varphi_{j,l}(x,y) - D_x \varphi_{j,l}(x,y) \cdot \xi$.  If $|\xi| \leq |x|/2$, then by Taylor's theorem and Remark~\ref{graphical rmk}(2) 
\begin{equation}\label{noncon triangle eqn2}
	\varphi_{j,l}(x-\xi,y-\zeta) = \varphi_{j,l}(x,y) - D_x \varphi_{j,l}(x,y) \cdot \xi + \mathcal{R}_{j,l}(x,y) 
\end{equation}
for all $j \in \{1,\ldots,J\}$ and $l \in \{1,\ldots,q_j\}$, where $|\mathcal{R}_{j,l}(x,y)| \leq C(n,m,q,\alpha) \,\widehat{E} \,|x|^{\alpha-2} \,|\xi|^2$.  Choosing $j \in \{1,\ldots,J\}$ and $l \in \{1,\ldots,q_j\}$ so that $\op{dist}(X, Z+\Gamma(\op{spt}\op{graph}\varphi)) = \op{dist}(X, Z+\Gamma(\op{graph}\varphi_{j,l}))$ and using \eqref{noncon triangle eqn1} and \eqref{noncon triangle eqn2}
\begin{align*}
	\op{dist}(X, \op{spt}\op{graph}\widehat{\varphi})
	\leq\,& |z - \widehat{\varphi}_{j,l}(x,y)|
	\\ \leq\,& |z-\chi - A(x-\xi,y-\zeta) - \varphi_{j,l}(x-\xi,y-\zeta)| + C \widehat{E} \,|x|^{\alpha-2} \,|\xi|^2 \nonumber 
	\\ \leq\,& 2 \op{dist}(\Gamma^{-1}(X-Z), \op{spt}\op{graph}\varphi) + C \widehat{E} \,\|A\| + C \widehat{E} \,|x|^{\alpha-2} \,|\xi|^2 \nonumber
	\\ =\,& 2 \op{dist}(X, Z + \Gamma(\op{spt}\op{graph}\varphi)) + C \widehat{E} \,\|A\| + C \widehat{E} \,|x|^{\alpha-2} \,|\xi|^2 , \nonumber
\end{align*}
where $C = C(n,m,q,\alpha) \in (0,\infty)$ are constants.
	
\noindent (2)  \textit{For each $\gamma \in (0,1)$ and $\tau \in (0,1/2)$ there exists $\varepsilon_0 = \varepsilon_0(n,m,q,\alpha,\gamma,\tau)$ and $\beta_0 = \beta_0(n,m,q,\alpha,\gamma,\tau)$ such that the following holds true.  Let $X = (x,y,z) \in \op{spt} T \cap \mathbf{C}_{\gamma}(0)$ such that $|\xi| < |x|/2$ and $|x| \geq \tau/2$.  Let $\varphi_{j,l}$ be as in \eqref{varphi_localized} with $B = B_{|x|/2}(x,y)$ and let $j_0 \in \{1,\ldots,J\}$ and $l_0 \in \{1,\ldots,q_{j_0}\}$ such that $\op{dist}(X,\op{spt}\op{graph}\varphi) = \op{dist}(X,\op{graph}\varphi_{j_0,l_0})$.  Then 
\begin{align}\label{noncon triangle2} 
	&|\chi + A(x-\xi,y-\zeta) - D_x \varphi_{j_0,l_0}(x,y) \cdot \xi| 
	\leq 2 \op{dist}(X, \op{spt} \op{graph} \varphi) \\&\hspace{15mm} + 2 \op{dist}(X, Z + \Gamma(\op{spt} \op{graph} \varphi)) 
		+ C \widehat{E} \,\|A\| + C \widehat{E} \,|x|^{\alpha-2} \,|\xi|^2 , \nonumber 
\end{align}
where $C = C(n,m,q,\alpha) \in (0,\infty)$ is a constant.}

By the fundamental theorem of calculus and Remark~\ref{graphical rmk}(2), $|\varphi_{j,l}(x,y) - \varphi_{j,l}(x-\xi,y-\zeta)| \leq C\widehat{E} \,|x|^{\alpha-1} |\xi|$, where $C = C(n,m,q,\alpha) \in (0,\infty)$ is a constant.  Thus by \eqref{noncon5 concl} and \eqref{graphical conclA1} 
\begin{align}\label{noncon triangle eqn4}
	&|\varphi_{j,l}(x,y) - \chi - A(x-\xi,y-\zeta) - \varphi_{j,l}(x-\xi,y-\zeta)| + \widehat{E} \|A\| 
	\\ \leq\,& C (\widehat{E} |\xi| + |\chi| + \|A\|) 
	\leq C \delta \inf_{\mathbb{S}^1} \op{sep} \varphi 
	\leq C \delta |x|^{-\alpha} \op{sep} \varphi(x,y) \nonumber
	\\ \leq\,& C\delta (\tau/2)^{-\alpha} \op{sep} \varphi(x,y) \nonumber
\end{align}
for all $j \in \{1,\ldots,J\}$ and $l \in \{1,\ldots,q_j\}$, where $C = C(n,m,q,\alpha,\gamma) \in (0,\infty)$ is a constant.  Taking $\delta$ to be sufficiently small, by \eqref{noncon triangle eqn1} and \eqref{noncon triangle eqn4}, $\op{dist}(X, Z + \Gamma(\op{spt}\op{graph}\varphi)) = \op{dist}(X, Z + \Gamma(\op{graph}\varphi_{j_0,l_0}))$.  Therefore, by the triangle inequality and \eqref{noncon triangle1}
\begin{align*}%\label{noncon triangle2} 
	&|\chi + A(x-\xi,y-\zeta) - D_x \varphi_{j_0,l_0}(x,y) \cdot \xi| 
	\leq |z - \varphi_{j_0,l_0}(x,y)| + |z - \widehat{\varphi}_{j_0,l_0}(x,y)| 
	\\ \leq\,& 2 \op{dist}(X, \op{spt} \op{graph} \varphi) + 2 \op{dist}(X, Z + \Gamma(\op{spt} \op{graph} \varphi)) 
		+ C \widehat{E} \,\|A\| + C \widehat{E} \,|x|^{\alpha-2} \,|\xi|^2 , \nonumber
\end{align*}
where $\widehat{\varphi}_{j_0,l_0}(x,y) = \chi - A(x-\xi,y-\zeta) + \varphi_{j_0,l_0}(x,y) - D_x \varphi_{j_0,l_0}(x,y) \cdot \xi$ and $C = C(n,m,q,\alpha) \in (0,\infty)$ is a constant.

\noindent (3)  By the fundamental theorem of calculus and Remark~\ref{graphical rmk}(2) 
\begin{equation*}
	\mathcal{G}(\varphi(x,y),\varphi(x-\xi,y-\zeta)) \leq \int_0^1 |D_x \varphi(x-t\xi,y-t\zeta)| \,|\xi| \,dt \leq C \widehat{E} |\xi| 
\end{equation*}
for each $(x,y) \in B_1(0)$, where $C = C(n,m,q,\alpha) \in (0,\infty)$ is a constant.  Hence using $\|\Gamma - I\| \leq C(n,m) \,\|A\|$ 
\begin{equation}\label{noncon triangle3}
	\mathcal{G}((x,y,\varphi(x,y)), Z + \Gamma(x-\xi,y-\zeta,\varphi(x-\xi,y-\zeta))) \leq C (\widehat{E} |\xi| + |\chi| + \|A\|) 
\end{equation}
for each $(x,y) \in B_1(0)$, where $C = C(n,m,q,\alpha) \in (0,\infty)$ is a constant.  Here for each $(x,y) \in \mathbb{R}^n$ we let $(x,y,\varphi(x,y)) = \sum_{j=1}^J \sum_{l=1}^{q_j} m_j \llbracket (x,y,\varphi_{j,l}(x,y)) \rrbracket$ and $Z + \Gamma(x,y,\varphi(x,y)) = \sum_{j=1}^J \sum_{l=1}^{q_j} m_j \llbracket Z + \Gamma(x,y,\varphi_{j,l}(x,y)) \rrbracket$ as elements of $\mathcal{A}_q(\mathbb{R}^{n+m})$, where $\varphi(x,y) = \sum_{j=1}^J \sum_{l=1}^{q_j} m_j \llbracket \varphi_{j,l}(x,y) \rrbracket$.
}\end{remark}

\begin{remark}\label{noncon rescale rmk}
Let $p \in \{1,2,\ldots,\lfloor q/q_0 \rfloor\}$ such that either $p = 1$ or $p > 1$ and (H4) holds true.  For each $\rho \in (0,1/8]$ there exists $\delta = \delta(n,m,q,\alpha,\rho) > 0$, $\varepsilon_0 = \varepsilon_0(n,m,q,\alpha,\rho)$, and $\beta_0 = \beta_0(n,m,q,\alpha,\rho)$ such that the following holds true.  Let $\varphi \in \Phi_{\alpha,q,p}$ and $T$ satisfy Hypothesis~$(\star)$ and Hypothesis~$(\star\star)$.  Let $\mathcal{B}_q$ be as in Definition~\ref{introduction B defn} (with $\mu = 1/q$).  Let $Z = (\xi,\zeta,\chi) \in \mathcal{B}_q$ such that $\mathcal{N}_{T,{\rm Pl}}(Z) \geq \alpha$ and let $P_Z$ be the tangent plane to $T$ at $Z$.  Let $\Gamma$ be a rotation of $\mathbb{R}^{n+m}$ such that $\Gamma(\mathbb{R}^n \times \{0\}) = P_Z$ and \eqref{noncon rotate} holds true.  Then Hypothesis~$(\star)$ and Hypothesis~$(\star\star)$ hold true with $C\varepsilon_0$, $C\delta + C\beta_0$, $\Gamma^{-1}_{\#} \eta_{Z,\rho\#} T$, and $\rho^{1-\alpha} \varphi$ in place of $\varepsilon_0$, $\beta_0$, $T$, and $\varphi$, where $C = C(n,m,q,\alpha,\rho) \in (0,\infty)$ is a constant.

{\rm To see this, let $T^{(Z)} = \Gamma^{-1}_{\#} \eta_{Z,\rho\#} T$.  Recall that from Remark~\ref{planar freq rmk} that we can apply Lemma~\ref{uniform branch set lemma} with $\mu = 1/q$.  By \eqref{uniform branch set concl1} and \eqref{uniform branch set concl2}, we have that $T$ is a locally area-minimizing rectifiable current of $\mathbf{C}_{1/4}(Z,P_Z)$ with $(\partial T) \llcorner \mathbf{C}_{1/4}(Z,P_Z) = 0$ and $\pi_{P_Z\#} (T \llcorner \mathbf{C}_{1/4}(Z,P_Z)) = q \llbracket B_{1/4}(Z,P_Z) \rrbracket$.  By the triangle inequality and \eqref{main hyp eqn1}, and \eqref{uniform branch set concl3} 
\begin{equation}\label{noncon triangle eqn6}
	\sup_{X \in \op{spt} T \cap \mathbf{C}_{1/4}(Z,P_Z)} \op{dist}(X,Z+P_Z) 
	\leq 2 \sup_{X \in \op{spt} T} \op{dist}(X,P_0) + C  |\vec P_Z - \vec P_0| 
	\leq C \varepsilon_0
\end{equation}
where $C = C(n,m,q) \in (0,\infty)$ are constants.  By \eqref{noncon triangle eqn6},~\cite[Lemma~3.7]{KrumWica}, and \eqref{oriented tilt and mass}, we have that $\|T\|(\mathbf{C}_{\rho}(Z,P_Z)) \leq (q+C\varepsilon_0^2) \,\omega_n \rho^n < (q+1/2) \,\omega_n \rho^n$, where $C = C(n,m,q) \in (0,\infty)$ is a constant.  Therefore, Hypothesis~$(\star)$ holds true with $T^{(Z)}$ and $C(n,m,q) \,\rho^{-1} \varepsilon_0$ in place of $T$ and $\varepsilon_0$.
	
We want to show that 
\begin{equation}\label{noncon triangle eqn8}
	Q_{\varsigma/2}(T^{(Z)},\rho^{1-\alpha}\varphi,\mathbf{C}_1(0)) 
		\leq C (\beta_0 + \delta) \inf_{\varphi' \in \bigcup_{p'=0}^{p-1} \Phi_{\alpha,q,p'}} Q_{\varsigma/2}(T^{(Z)},\varphi',\mathbf{C}_1(0)) 
\end{equation}
for some constant $C = C(n,m,q,\alpha,\rho) \in (0,\infty)$ is a constant (as in Hypothesis~$(\star\star)$ with $C (\beta_0 + \delta)$, $T^{(Z)}$, and $\rho^{1-\alpha}\varphi$ in place of $\beta_0$, $T$, and $\varphi$).  Let $\widehat{E} = \widehat{E}(T,P_0,\mathbf{C}_1(0))$ and let $A : \mathbb{R}^n \rightarrow \mathbb{R}^m$ is a linear transformation such that $P_Z$ is given by \eqref{noncon tangent plane} and $\|A\| \leq C(n,m,q) \,\widehat{E}$.  By the triangle inequality, \eqref{noncon triangle3}, \eqref{main hyp eqn2}, and \eqref{noncon5 concl}
\begin{align}\label{noncon triangle eqn9}
	&\int_{\mathbf{C}_{\rho}(Z,P_Z)} \op{dist}^2(X, Z + \Gamma(\op{spt}\op{graph}\varphi)) \,d\|T\|(X) 
	\\ \leq\,& 2 \int_{\mathbf{C}_1(0)} \op{dist}^2(X,\op{spt}\op{graph}\varphi) \,d\|T\|(X) + C (|\chi|^2 + \widehat{E}^2 |\xi|^2 + \|A\|^2) \nonumber 
	\\ \leq\,& C (\beta_0^2 + \delta^2) \inf_{\varphi' \in \bigcup_{p'=0}^{p-1} \Phi_{\alpha,q,p'}} Q_{\varsigma/2}(T,\varphi',\mathbf{C}_1(0))^2 , \nonumber
\end{align}
where $C = C(n,m,q,\alpha) \in (0,\infty)$ are constants.  By the triangle inequality, \eqref{noncon triangle3}, Remark~\ref{graphical rmk}(2), Theorem~\ref{graphical_thm}(B), \eqref{main hyp eqn2}, and \eqref{noncon5 concl} 
\begin{align}\label{noncon triangle eqn10}
	&\int_{B_{\rho/2}(0) \cap \{r > \rho/16\}} \op{dist}^2(Z + \Gamma(x,y,\varphi_{j,l}(x,y)), \op{spt} T) \,dx \,dy 
	\\ \leq\,& 2 \int_{B_{\rho}(0,\zeta) \cap \{r > \rho/32\}} \sum_{j=1}^J \sum_{l=1}^{q_j} \op{dist}^2((x,y,\varphi_{j,l}(x,y)), \op{spt} T) \,dx \,dy 
		\nonumber\\&+ C (|\chi|^2 + \widehat{E}^2 |\xi|^2 + \|A\|^2) \nonumber
	\\ \leq\,& C E_{\varsigma/2}(T,\varphi,\mathbf{C}_1(0))^2 + C (|\chi|^2 + \widehat{E}^2 |\xi|^2 + \|A\|^2) \nonumber
	\\ \leq\,& C (\beta_0^2 + \delta^2) \inf_{\varphi' \in \bigcup_{p'=0}^{p-1} \Phi_{\alpha,q,p'}} Q_{\varsigma/2}(T,\varphi',\mathbf{C}_1(0))^2 , \nonumber
\end{align}
where $C = C(n,m,q,\alpha,\rho) \in (0,\infty)$ are constants.  Therefore, by \eqref{noncon triangle eqn9} and \eqref{noncon triangle eqn10} 
\begin{equation}\label{noncon triangle eqn11}
	Q_{\varsigma/2}(T^{(Z)},\rho^{1-\alpha}\varphi,\mathbf{C}_1(0)) 
	\leq C (\beta_0 + \delta) \inf_{\varphi' \in \bigcup_{p'=0}^{p-1} \Phi_{\alpha,q,p'}} Q_{\varsigma/2}(T,\varphi',\mathbf{C}_1(0)) 
\end{equation}
for some constant $C = C(n,m,q,\alpha,\rho) \in (0,\infty)$.  

Let $p^* \in \{0,1,\ldots,p-1\}$ and $\varphi^* \in \Phi_{\alpha,q,p^*}$ such that 
\begin{equation}\label{noncon triangle eqn12}
	Q_{\varsigma/2}(T^{(Z)},\rho^{1-\alpha} \varphi^*,\mathbf{C}_1(0)) 
		\leq 2 \inf_{\varphi' \in \bigcup_{p'=0}^{p-1} \Phi_{\alpha,q,p'}} Q_{\varsigma/2}(T^{(Z)},\varphi',\mathbf{C}_1(0)) . 
\end{equation}
Let $\varphi^* = \sum_{j=1}^{J^*} m^*_j \varphi^*_j$ with $q^*_j$-valued components $\varphi^*_j$ with positive integer multiplicity $m^*_j$ (as in Definition~\ref{Phi-p_defn}).  By \eqref{uniform branch set concl5} from Lemma~\ref{uniform branch set lemma} (with $\mu = 1/q$), $\widehat{E}(T^{(Z)},P_0,\mathbf{C}_1(0)) = \widehat{E}(T,P_Z,\mathbf{C}_{\rho}(Z,P_Z)) \leq C(n,m,q) \,\widehat{E}(T,P_0,\mathbf{C}_1(0))$.  Thus arguing as in Remark~\ref{graphical rmk}(2) 
\begin{align*}
	\min\left\{ \tfrac{\rho}{8}, c \|\varphi^*\|_{L^{\infty}(B_1(0))} - C \widehat{E}(T,P_0,\mathbf{C}_1(0)) \right\} 
	\leq C Q_{\varsigma/2}(T^{(Z)},\varphi^*,\mathbf{C}_1(0)) 
\end{align*}
where $c = c(n,m,q,\alpha,\rho) > 0$ and $C = C(n,m,q,\alpha,\rho) \in (0,\infty)$ are constants.  By \eqref{noncon triangle eqn12} and Remark~\ref{graphical rmk}(1), $Q_{\varsigma/2}(T^{(Z)},\varphi^*,\mathbf{C}_1(0)) \leq C \widehat{E}(T^{(Z)},P_0,\mathbf{C}_1(0)) \leq C \widehat{E}(T,P_0,\mathbf{C}_1(0)) < C\varepsilon_0^2$, where $C = C(n,m,q,\alpha,\rho) \in (0,\infty)$ are constants.  Thus provided $\varepsilon_0$ is sufficiently small, $\|\varphi^*\|_{L^{\infty}(B_1(0))} \leq C \widehat{E}(T,P_0,\mathbf{C}_1(0))$, where $C = C(n,m,q,\alpha,\rho) \in (0,\infty)$ is a constant.  Thus arguing as in Remark~\ref{noncon triangle rmk}(3) above, \eqref{noncon triangle3} holds true with $\varphi^*$ in place of $\varphi$.  Provided $\beta_0 \leq 1/2$, by the triangle inequality and \eqref{main hyp eqn2} 
\begin{align}\label{noncon triangle eqn13}
	&\sqrt{2q \omega_n} \op{dist}_{\mathcal H}(\op{spt}\op{graph}\varphi \cap \mathbf{C}_1(0), \op{spt}\op{graph}\varphi^* \cap \mathbf{C}_1(0)) 
	\\ \geq\,& Q_{\varsigma/2}(T,\varphi^*,\mathbf{C}_1(0)) - Q_{\varsigma/2}(T,\varphi,\mathbf{C}_1(0)) \nonumber 
	\\ \geq\,& (1-\beta_0) \,Q_{\varsigma/2}(T,\varphi^*,\mathbf{C}_1(0)) \geq \tfrac{1}{2} Q_{\varsigma/2}(T,\varphi^*,\mathbf{C}_1(0)) . \nonumber
\end{align}
On the other hand, let $K_j, u_j, v_j$ be as in Theorem~\ref{graphical_thm}(C) with $\gamma = \rho$ and $\tau = \rho/8$.  By the triangle inequality, \eqref{noncon triangle3} with $\varphi^*$ in place of $\varphi$, Theorem~\ref{graphical_thm}(C), \eqref{main hyp eqn2}, and \eqref{noncon5 concl} 
	\begin{align}\label{noncon triangle eqn14}
		&\rho^{2\alpha-2} \sum_{j=1}^J m_j \int_{B_{1/2}(0,\zeta) \cap \{r > 1/16\}} \sum_{l=1}^{q_j} \op{dist}^2((x,y,\varphi_{j,l}(x,y)), 
			\op{spt}\op{graph}\varphi^*) \,dx \,dy 
		\\ \leq\,& \frac{2}{\rho^{n+2}} \int_{B_{\rho/2}(0,\zeta) \cap \{r > \rho/16\}} \sum_{j=1}^J \sum_{l=1}^{q_j} \op{dist}^2((x,y,u_{j,l}(x,y)),
			\op{spt}\op{graph}\varphi^*) \,dx \,dy \nonumber \\&+ C \sup_{B_{\rho/2}(0,\zeta) \cap \{r > \rho/16\}} 
			\sum_{j=1}^J \sum_{l=1}^{q_j} |v_{j,l}|^2 \nonumber
		\\ \leq\,& \frac{2}{\rho^{n+2}} \int_{\mathbf{C}_{\rho/2}(0,\zeta) \cap \{r > \rho/16\}} \op{dist}^2(X,\op{spt}\op{graph}\varphi^*) \,d\|T\|(X) 
			+ C E_{\varsigma/2}(T,\varphi,\mathbf{C}_1(0))^2 \nonumber
		\\ \leq\,& \frac{4}{\rho^{n+2}} \int_{\mathbf{C}_{\rho/2}(0,\zeta) \cap \{r > \rho/16\}} \op{dist}^2(X,Z+\Gamma(\op{spt}\op{graph}\varphi^*)) \,d\|T\|(X) 
			\nonumber \\&+ C (|\chi|^2 + \widehat{E}^2 |\xi|^2 + \|A\|^2) + C E_{\varsigma/2}(T,\varphi,\mathbf{C}_1(0))^2 \nonumber
		\\ \leq\,& 4 Q_{\varsigma/2}(T^{(Z)},\varphi^*,\mathbf{C}_1(0))^2 + C (\beta_0^2 + \delta^2) \inf_{\varphi' \in \bigcup_{p'=0}^{p-1} \Phi_{\alpha,q,p'}} 
			Q_{\varsigma/2}(T,\varphi',\mathbf{C}_1(0))^2 , \nonumber
	\end{align}
	where $v_{j,l}$ and $u_j$ are as in \eqref{v_localized} and \eqref{u defn} and $C = C(n,m,q,\alpha,\rho) \in (0,\infty)$ are constants.  By the triangle inequality, \eqref{noncon triangle3} with $\varphi^*$ in place of $\varphi$, Theorem~\ref{graphical_thm}(B), \eqref{noncon triangle eqn1}, \eqref{noncon5 concl}, and \eqref{main hyp eqn2} 
\begin{align}\label{noncon triangle eqn15}
	&\rho^{2\alpha-2} \sum_{j=1}^{J^*} m^*_j \int_{B_{1/2}(0) \cap \{r > 1/16\}} \sum_{l=1}^{q^*_j} \op{dist}^2((x,y,\varphi^*_{j,l}(x,y)), 
		\op{spt}\op{graph}\varphi) \,dx \,dy 
	\\ \leq\,& \frac{2}{\rho^{n+2}} \sum_{j=1}^{J^*} m^*_j \int_{B_{\rho/2}(0) \cap \{r > \rho/16\}} \sum_{l=1}^{q^*_j} \op{dist}^2((x,y,\varphi^*_{j,l}(x,y)), 
		\op{spt} T) \,dx \,dy \nonumber \\&+ C \sup_{X \in \op{spt}\op{graph}\varphi \cap \mathbf{C}_{\rho}(0) \cap \{r > \rho/32\}} 
		\op{dist}^2(X,\op{spt} T) \nonumber 
	\\ \leq\,& \frac{4}{\rho^{n+2}} \sum_{j=1}^{J^*} m^*_j \int_{B_{\rho/2}(0) \cap \{r > \rho/16\}} \sum_{l=1}^{q^*_j} 
		\op{dist}^2(Z+\Gamma(x,y,\varphi^*_{j,l}(x,y)), \op{spt} T) \,dx \,dy \nonumber \\&+ C (|\chi|^2 + \widehat{E}^2 |\xi|^2 + \|A\|^2) 
		+ C E_{\varsigma/2}(T,\varphi,\mathbf{C}_1(0))^2 \nonumber 
	\\ \leq\,& 4 Q_{\varsigma/2}(T^{(Z)},\varphi^*,\mathbf{C}_1(0))^2 
		+ C (\beta_0^2 + \delta^2) \inf_{\varphi^* \in \bigcup_{p'=0}^{p-1} \Phi_{\alpha,q,p'}} Q_{\varsigma/2}(T,\varphi^*,\mathbf{C}_1(0))^2 , \nonumber 
\end{align}
where $\varphi^*_j(x,y) = \sum_{l=1}^{q^*_j} m^*_j \llbracket \varphi^*_{j,l}(x,y) \rrbracket$ for each $(x,y) \in \mathbb{R}^n$ (as in \eqref{varphi_localized}) and $C = C(n,m,q,\alpha,\rho) \in (0,\infty)$ are constants.  Using $\varphi, \varphi^* \in \Phi_{\alpha,q}$ and combining \eqref{noncon triangle eqn14} and \eqref{noncon triangle eqn15} 
\begin{align}\label{noncon triangle eqn16}
	&\op{dist}_{\mathcal H}(\op{spt}\op{graph}\varphi \cap \mathbf{C}_1(0), \op{spt}\op{graph}\varphi^* \cap \mathbf{C}_1(0)) 
	\\ \leq\,& C Q_{\varsigma/2}(T^{(Z)},\rho^{1-\alpha}\varphi^*,\mathbf{C}_1(0)) 
		+ C (\beta_0 + \delta) \inf_{\varphi' \in \bigcup_{p'=0}^{p-1} \Phi_{\alpha,q,p'}} Q_{\varsigma/2}(T,\varphi',\mathbf{C}_1(0)) , \nonumber
\end{align}
where $C = C(n,m,q,\alpha,\rho) \in (0,\infty)$ is a constant.  Combining \eqref{noncon triangle eqn13}, \eqref{noncon triangle eqn16}, and \eqref{noncon triangle eqn12} 
\begin{equation}\label{noncon triangle eqn17}
	\inf_{\varphi' \in \bigcup_{p'=0}^{p-1} \Phi_{\alpha,q,p'}} Q_{\varsigma/2}(T,\varphi',\mathbf{C}_1(0))
	\leq C \inf_{\varphi' \in \bigcup_{p'=0}^{p-1} \Phi_{\alpha,q,p'}} Q_{\varsigma/2}(T^{(Z)},\varphi',\mathbf{C}_1(0)) , 
\end{equation}
where $C = C(n,m,q,\alpha,\rho) \in (0,\infty)$ is a constant.  Combining \eqref{noncon triangle eqn11} and \eqref{noncon triangle eqn17} gives us \eqref{noncon triangle eqn8}.
}\end{remark}

\begin{proof}[Proof of Theorem~\ref{nonconcentration thm}]  
Suppose that $p \in \{1,2,\ldots,\lfloor q/q_0 \rfloor\}$ such that either $p = 1$ or $p > 1$ and (H4) holds true.  Let $\varphi \in \Phi_{\alpha,q,p}$ and $T$ satisfy Hypothesis~$(\star)$ and Hypothesis~$(\star\star)$.  Let $Z = (\xi,\zeta,\chi) \in \mathcal{B}_q$ such that $T$ is tangent to $P_Z$ at $Z$ and $\mathcal{N}_{T,{\rm Pl}}(Z) \geq \alpha$.  We will show that Theorem~\ref{nonconcentration thm}(A)--(D) hold true, which by mathematical induction will have completed the proof.

\noindent\textit{Proof of Conclusion~(A).}  First we claim that: \textit{there exists a constant $c = c(n,m,q,\alpha) > 0$ such that for every $\rho \in (0,1/4]$ there exists $\varepsilon_0 = \varepsilon_0(n,m,q,\alpha,\rho) \in (0,1)$ and $\beta_0 = \beta_0(n,m,q,\alpha,\rho) \in (0,1)$ such that if $\varphi$ and $T$ satisfy Hypothesis~$(\star)$ and Hypothesis~$(\star\star)$ and if $Z = (\xi,\zeta,\chi) \in \mathcal{B}_q$ with $\mathcal{N}_{T,{\rm Pl}}(Z) \geq \alpha$, then }
\begin{align}\label{nonconA eqn1}
	&c \,(\widehat{E}^2 |\xi|^2 \rho^{n+2\alpha-2} + |\chi|^2 \rho^n + \|A\|^2 \rho^{n+2}) 
	\\ \leq\,& \sum_{j=1}^J \int_{\mathbf{C}_{\rho/4}(\xi,\zeta) \cap \{r > \rho/32\}} \sum_{l=1}^{q_j} 
		|\chi - A(x-\xi,y-\zeta) - D_x \varphi_{j,l}(x,y) \cdot \xi|^2 \,d\|T_{j,l}\|(X) , \nonumber
\end{align}
\textit{where $\widehat{E} = \widehat{E}(T,P_0,\mathbf{C}_1(0))$.}  Here we write $\varphi = \sum_{j=1}^J m_j \varphi_j$ (as in Definition~\ref{Phi-p_defn}) and we let $T_j$ be as in Theorem~\ref{graphical_thm}(B) (with $\gamma = 1/32$ and $\tau = \rho/32$).  
The integrand on the right-hand side of \eqref{nonconA eqn1} is a well-defined Borel measure such that for each open ball $B \subset B_{\rho/4}(\xi,\zeta) \cap \{r > \rho/32\}$, the restriction of the integrand to $B \times \mathbb{R}^m$ is defined in terms of $\varphi_{j,l}$ and $T_{j,l}$ as in \eqref{varphi_localized} and \eqref{T_separated eqn4}.   (Alternatively, since the measures $\|T_{j,l}\|$ are absolutely continuous with respect to $\|T\|$ with Radon-Nikodym derivative equal to the characteristic function $\mathbf{1}_{\op{spt} T_{j,l}}$ of $\op{spt}T_{j,l}$, the right-hand side of \eqref{nonconA eqn1} is the integral of a well-defined single-valued function, expressed locally in terms of $\varphi_{j,l}$ and $\mathbf{1}_{\op{spt} T_{j,l}}$, with respect to the measure $\|T\|$.)
To prove the claim, suppose to the contrary that for every $c > 0$ there exists $\rho \in (0,1/4]$ and for $\nu = 1,2,3,\ldots$ there exists $\varepsilon_{\nu} \rightarrow 0^+$, $\beta_{\nu}\rightarrow 0^+$, $\varphi^{(\nu)} \in \Phi_{\alpha,q,p}$, an $n$-dimensional locally area-minimizing rectifiable current $T^{(\nu)}$ of $\mathbf{C}_1(0)$, $Z_{\nu} = (\xi_{\nu},\zeta_{\nu},\chi_{\nu}) \in \mathcal{B}_q(T^{(\nu)})$, and a linear transformation $A_{\nu} : \mathbb{R}^n \rightarrow \mathbb{R}^m$ such that Hypothesis~$(\star)$ and Hypothesis~$(\star\star)$ hold true with $\varepsilon_{\nu},\beta_{\nu},\varphi^{(\nu)},T^{(\nu)}$ in place of $\varepsilon_0,\beta_0,\varphi,T$, $\mathcal{N}_{T^{(\nu)},{\rm Pl}}(Z_{\nu}) \geq \alpha$, and the tangent plane to $T^{(\nu)}$ at $Z_{\nu}$ is $\{ (x,y) \in \mathbb{R}^n : z = A_{\nu} (x,y) \}$ but 
\begin{align}\label{nonconA eqn2} 
	&\sum_{j=1}^J \int_{\mathbf{C}_{\rho/4}(\xi_{\nu},\zeta_{\nu}) \cap \{r > \rho/32\}} \sum_{l=1}^{q_j} 
		|\chi_{\nu} - A_{\nu}(x-\xi_{\nu},y-\zeta_{\nu}) - D_x \varphi^{(\nu)}_{j,l}(x,y) \cdot \xi_{\nu}|^2 \,d\|T^{(\nu)}_{j,l}\|(X)
	\\<\,& c \,(\widehat{E}_{\nu}^2 |\xi_{\nu}|^2 \rho^{n+2\alpha-2} + |\chi_{\nu}|^2 \rho^n + \|A_{\nu}\|^2 \rho^{n+2}) , \nonumber
\end{align}
where $\widehat{E}_{\nu} = \widehat{E}(T^{(\nu)},P_0,\mathbf{C}_1(0))$.  We write $\varphi^{(\nu)} = \sum_{j=1}^J m_j \varphi^{(\nu)}_j$ with components $\varphi^{(\nu)}_j$ with multiplicity $m_j$ and we let $T^{(\nu)}_j$ be as in Theorem~\ref{graphical_thm}(B) with $T^{(\nu)}$ in place of $T$ (and with $\gamma = 1/32$ and $\tau = \rho/32$).  We let $\varphi^{(\nu)}_{j,l}$ be as in \eqref{varphi_localized} with $\varphi^{(\nu)}_j$ in place of $\varphi_j$ and $T^{(\nu)}_{j,l}$ be as in \eqref{T_separated eqn4} with $T^{(\nu)}_j$ in place of $T_j$.  After passing to a subsequence, we assume that $J$, $m_j$, and $q_j$ are independent of $\nu$.  Set $M_{\nu} = (\widehat{E}_{\nu}^2 |\xi_{\nu}|^2 + |\chi_{\nu}|^2 + \|A_{\nu}\|^2)^{1/2}$.  By Theorem~\ref{graphical_thm}, $\xi_{\nu} \rightarrow 0$ in $\mathbb{R}^2$ as $\nu \rightarrow \infty$.  After passing to a further subsequence, let $\zeta_{\nu} \rightarrow \zeta$ in $\mathbb{R}^{n-2}$, $\widehat{E}_{\nu} \xi_{\nu}/M_{\nu} \rightarrow \kappa$ in $\mathbb{R}^2$, $\chi_{\nu}/M_{\nu} \rightarrow \lambda$ in $\mathbb{R}^m$, and $A_{\nu}/M_{\nu} \rightarrow \Lambda$ as a linear transformation from $\mathbb{R}^n$ to $\mathbb{R}^m$.  Note that $|\kappa|^2 + |\lambda|^2 + \|\Lambda\|^2 = 1$ so that $\kappa,\lambda,\Lambda$ are not all zero.  By Remark~\ref{graphical rmk}(2), after passing to a subsequence, there is a non-zero function $\varphi^{(0)} \in \Phi_{\alpha,q}$ such that $\varphi^{(\nu)}/\widehat{E}_{\nu} \rightarrow \varphi^{(0)}$ uniformly on compact subsets of $\mathbb{R}^n$ as $\nu \rightarrow \infty$ and $C^{-1} \leq \|\varphi^{(0)}\|_{L^2(B_1(0))} \leq C$ for some constant $C = C(n,m,q,\alpha) \in (1,\infty)$.  By Theorem~\ref{graphical_thm}(B) and~\cite[Theorems~32.2 and 34.5]{SimonGMT}, $T^{(\nu)}_{j,l} \rightarrow m_j \llbracket \mathbb{R}^n \times \{0\} \rrbracket$ weakly in each open ball $B \subset \mathbf{C}_{1/32}(0) \cap \{r > \rho/32\}$.  Hence dividing both sides of \eqref{nonconA eqn2} by $M_{\nu}^2$ and letting $\nu \rightarrow \infty$, we deduce that: \textit{for every $c > 0$ there exists $\rho > 0$, $\zeta \in \mathbb{R}^{n-2}$, $\kappa \in \mathbb{R}^2$, $\lambda \in \mathbb{R}^m$, a linear transformation $\Lambda : \mathbb{R}^n \rightarrow \mathbb{R}^m$, and $\varphi^{(0)} \in \Phi_{\alpha,q}$ such that $\kappa,\lambda,\Lambda$ are not all zero, $C^{-1} \leq \|\varphi^{(0)}\|_{L^2(B_1(0))} \leq C$ for some (fixed) constant $C = C(n,m,q,\alpha) \in (1,\infty)$, and }
\begin{align*} 
	&\sum_{j=1}^J m_j \int_{B_{\rho/4}(0,\zeta) \cap \{r > \rho/32\}} \sum_{l=1}^{q_j} 
		|\lambda - \Lambda (x,y-\zeta) - D_x \varphi^{(0)}_{j,l}(x,y) \cdot \kappa|^2 \,dx \,dy 
	\\ \leq\,& c \,(|\kappa|^2 \rho^{n+2\alpha-2} + |\lambda|^2 \rho^n + \|\Lambda\|^2 \rho^{n+2}) . \nonumber 
\end{align*}
Assuming $c > 0$ is sufficiently small, it follows using the mutual orthogonality of $\lambda$, $\Lambda (x,y)$, and $\frac{1}{q_j} \sum_{l=1}^{q_j} D\varphi^{(0)}_{j,l}(x) \cdot \kappa$ in $L^2(B_{\rho/4}(0,\zeta) \cap \{r > \rho/32\})$ that $\kappa \neq 0$ (as otherwise $\kappa,\lambda,\Lambda$ are all zero) and 
\begin{equation*}
	\sum_{j=1}^J m_j \int_{B_{\rho/4}(0,\zeta) \cap \{r > \rho/32\}} \sum_{l=1}^{q_j} |D_x \varphi^{(0)}_{j,l}(x,y) \cdot \kappa|^2 \,dx \,dy
	\leq c \,|\kappa|^2 \rho^{n+2\alpha-2} . 
\end{equation*} 
By translating and rescaling, we may assume that $\zeta = 0$, $\rho = 1$, and $|\kappa| = 1$.  Suppose to the contrary that for $\nu = 1,2,3,\ldots$ there exists $\kappa_{\nu} \in \mathbb{S}^1$ and $\varphi^{(\nu,0)} \in \Phi_{\alpha,q}$ such that $C^{-1} \leq \|\varphi^{(\nu,0)}\|_{L^2(B_1(0))} \leq C$ and 
\begin{equation}\label{nonconA eqn3} 
	\sum_{j=1}^J m_j \int_{B_{1/4}(0) \cap \{r > 1/32\}} \sum_{l=1}^{q_j} |D_x \varphi^{(\nu,0)}_{j,l}(x,y) \cdot \kappa_{\nu}|^2 \,dx \,dy \leq \frac{1}{\nu} .
\end{equation}
After passing to a subsequence, let $\kappa_{\nu} \rightarrow \kappa$ in $\mathbb{S}^1$ and let $\varphi^{(0)} \in \Phi_{\alpha,q}$ such that $\varphi^{(\nu,0)} \rightarrow \varphi^{(0)}$ in $L^2(B_1(0),\mathcal{A}_q(\mathbb{R}^m))$.  Letting $\nu \rightarrow \infty$ in \eqref{nonconA eqn3} 
\begin{equation*}
	D_x \varphi^{(0)}_{j,l}(x,y) \cdot \kappa = 0 \text{ for all } (x,y) \in B_{1/4}(0) \cap \{r > 1/32\} .
\end{equation*}
In other words, $\varphi^{(0)}$ is translation invariant in the $\kappa$ direction.  This contradicts $\varphi^{(0)}$ being a non-zero function in $\Phi_{\alpha,q}$.

Now let $\rho = \rho(n,m,q,p,\alpha) \in (0,1/32]$ and $\tau = \tau(n,m,q,p,\alpha,\rho) \in (0,\rho/32]$ to be later determined.  Assuming $\varepsilon_0 = \varepsilon_0(n,m,q,\alpha,\rho,\tau)$ and $\beta_0 = \beta_0(n,m,q,\alpha,\rho,\tau)$ are sufficiently small, by Theorem~\ref{graphical_thm} we have that $|\xi| \leq \tau \leq \rho/32$.  Using \eqref{noncon triangle2} (with $\tau = \rho/32$) to bound the first term on the right-hand side of \eqref{nonconA eqn1} 
\begin{align}\label{nonconA eqn7}
	&\widehat{E}^2 |\xi|^2 \rho^{n+2\alpha-2} + |\chi|^2 \rho^n + \|A\|^2 \rho^{n+2} 
	\\ \leq\,& C \int_{\mathbf{C}_{\rho}(\xi,\zeta)} \op{dist}^2(X,Z+\Gamma(\op{spt}\op{graph}\varphi)) \,d\|T\|(X) \nonumber
		\\&+ C \int_{\mathbf{C}_{\rho}(\xi,\zeta)} \op{dist}^2(X,\op{spt}\op{graph}\varphi) \,d\|T\|(X) \nonumber 
		\\&+ C \int_{\mathbf{C}_{\rho}(\xi,\zeta) \cap \{|x| \geq 2|\xi|\}} 
			(\widehat{E}^2 \|A\|^2 + \widehat{E}^2 |x|^{2\alpha-4} |\xi|^4) \,d\|T\|(X) \nonumber
		\\&+ C \sum_{j=1}^J m_j \int_{\mathbf{C}_{\rho}(\xi,\zeta) \cap \{|x| \leq 2|\xi|\}} \sum_{l=1}^{q_j} 
			|\chi - A(x-\xi,y-\zeta) - D_x \varphi_{j,l}(x,y) \cdot \xi|^2 \,d\|T\|(X) , \nonumber
\end{align}
where $C = C(n,m,q,\alpha) \in (0,\infty)$ is a constant.  To bound the first term on the right-hand side of \eqref{nonconA eqn7}, by (H4) we can apply Remark~\ref{noncon rescale rmk} to deduce that Hypothesis~$(\star)$ and Hypothesis~$(\star\star)$ hold true with $C \varepsilon_0$, $C (\beta_0+\delta)$, and $\Gamma^{-1}_{\#} \eta_{Z,1/8\#} T$ in place of $\varepsilon_0$, $\beta_0$, and $T$, where $\delta = \delta(n,m,q,\alpha) \in (0,1)$ is a suitably small and $C = C(n,m,q,\alpha) \in (0,\infty)$ is a constant.  Hence we can apply Theorem~\ref{radial noncon thm} with $\Gamma^{-1}_{\#} \eta_{Z,1/8\#} T$ in place of $T$ and with $\sigma = 1/q$ to obtain 
\begin{align}\label{nonconA eqn8}
	&\int_{\mathbf{C}_{1/16}(Z,P_Z)} \frac{\op{dist}^2(X,Z+\Gamma(\op{spt}\op{graph}\varphi))}{|X-Z|^{n+2\alpha-1/q}} \,d\|T\|(X)
	\\ \leq\,& C \int_{\mathbf{C}_{1/8}(Z,P_Z)} \op{dist}^2(X,Z+\Gamma(\op{spt}\op{graph}\varphi)) \,d\|T\|(X) 
		+ C \widehat{E}(T,P_Z,\mathbf{C}_{1/8}(Z,P_Z))^{2+\varsigma} \nonumber 
\end{align}
where $C = C(n,m,q,p,\alpha) \in (0,\infty)$ is a constant.  It follows using \eqref{uniform branch set concl3} from Lemma~\ref{uniform branch set lemma} (with $\mu = 1/q$) that $\op{spt} T \cap \mathbf{C}_{1/32}(\xi,\zeta) \subset \mathbf{C}_{1/16}(Z,P_Z)$.  Moreover, using \eqref{uniform branch set concl5}, $\widehat{E}(T,P_Z,\mathbf{C}_{7\rho'/4}(Z,P_Z)) \leq C \widehat{E} (\rho')^{1/q}$ for all $\rho' \in (0,1/8]$ and thus $\widehat{E}(T,P_Z,\mathbf{C}_{1/8}(Z,P_Z)) \leq C \widehat{E}$, where $C = C(n,m,q) \in (0,\infty)$ are constants.  For each $X = (x,y,z) \in \op{spt} T \cap \mathbf{C}_{1/16}(Z,P_Z)$, by setting $\rho' = 2|X-Z|$ and using~\cite[Lemma~3.5]{KrumWica}, $|z - \zeta| \leq C(n,m,q) \,\widehat{E} |X-Z|^{1+1/q} \leq \frac{1}{2} |X-Z|$ and thus $|X-Z| \leq 2 |(x-\xi,y-\zeta)|$.  That is, $\op{spt} T \cap \mathbf{C}_{\rho}(\xi,\zeta) \subset \mathbf{B}_{2\rho}(Z)$.  Hence also using \eqref{noncon triangle3}, \eqref{nonconA eqn8} gives us 
\begin{align}\label{nonconA eqn9}
	&\rho^{-n-2\alpha+1/q} \int_{\mathbf{C}_{\rho}(\xi,\zeta)} \op{dist}^2(X,Z+\Gamma(\op{spt}\op{graph}\varphi)) \,d\|T\|(X) 
	\\ \leq\,& C \int_{\mathbf{C}_1(0)} \op{dist}^2(X,\op{spt}\op{graph}\varphi) \,d\|T\|(X) + C (\widehat{E}^2 |\xi|^2 + |\chi|^2 + \|A\|^2) 
		+ C \widehat{E}^{2+\varsigma} , \nonumber
\end{align}
where $C = C(n,m,q,p,\alpha) \in (0,\infty)$ are constants.  To bound the third term on the right-hand side of \eqref{nonconA eqn7}, by using $\|T\|(\mathbf{C}_{\rho}(\xi,\zeta)) \leq \|T\|(\mathbf{B}_{2\rho}(Z)) \leq (q+1/2) \,\omega_n(4\rho)^n$ and separately considering the cases $\alpha < 2$ and $\alpha > 2$ 
\begin{align}\label{nonconA eqn10}
	&\int_{\mathbf{C}_{\rho}(\xi,\zeta) \cap \{|x| \geq 2|\xi|\}} (\widehat{E}^2 \|A\|^2 + \widehat{E}^2 |x|^{2\alpha-4} |\xi|^4) \,d\|T\|(X)
	\\ \leq\,& C \widehat{E}^2 \|A\|^2 \rho^n + C\widehat{E}^2 \rho^n (|\xi|^{2\alpha} + \rho^{2\alpha-4} |\xi|^4) , \nonumber 
\end{align}
where $C = C(n,q) \in (0,\infty)$ is a constant.  Similarly, to bound the last term on the right-hand side of \eqref{nonconA eqn7}, notice that assuming $\varepsilon_0$ is sufficiently small, by Lemma~\ref{lip approx thm} and the area formula, $\|T\|(\mathbf{C}_{\rho}(\xi,\zeta) \cap \{|x| \leq 2|\xi|\}) \leq C (|\xi|^2 \rho^{n-2} + \widehat{E}^{2+\varsigma})$, where $C = C(n,m,q) \in (0,\infty)$ is a constant.  This together with Remark~\ref{graphical rmk}(2) gives us 
\begin{align}\label{nonconA eqn11}
	&\sum_{j=1}^J m_j \int_{\mathbf{C}_{\rho}(\xi,\zeta) \cap \{|x| \leq 2|\xi|\}} \sum_{l=1}^{q_j} 
		|\chi - A(x-\xi,y-\zeta) - D_x \varphi_{j,l}(x,y) \cdot \xi|^2 \,d\|T\|(X)
	\\ \leq\,& C (|\xi|^2 \rho^{n-2} + \widehat{E}^{2+\varsigma}) \,(|\chi|^2 + \|A\|^2 \rho^2 + \widehat{E}^2 |\xi|^{2\alpha}) , \nonumber 
\end{align}
where $C = C(n,m,q,\alpha) \in (0,\infty)$ is a constant.  Thus bounding the right-hand side of \eqref{nonconA eqn7} using \eqref{nonconA eqn9}, \eqref{nonconA eqn10}, and \eqref{nonconA eqn11}
\begin{align}\label{nonconA eqn12}
	&\widehat{E}^2 |\xi|^2 \rho^{n+2\alpha-2} + |\chi|^2 \rho^n + \|A\|^2 \rho^{n+2}
	\leq C E_{\varsigma/2}(T,\varphi,\mathbf{C}_1(0))^2 \\&+ C \widehat{E}^2 |\xi|^2 \rho^{n+2\alpha-2} (\rho^{2-1/q} + |\xi|^{2\alpha-2} \rho^{-2\alpha+2} 
		+ |\xi|^2 \rho^{-2} + |\xi|^{2\alpha} \rho^{-2\alpha}) \nonumber 
		\\&+ C |\chi|^2 \rho^n (\rho^{2\alpha-1/q} + |\xi|^2 \rho^{-2}) 
		+ C \|A\|^2 \rho^{n+2} (\rho^{2\alpha-2-1/q} + \widehat{E}^2 \rho^{-2} + |\xi|^2 \rho^{-2}) \nonumber
\end{align}
where $C = C(n,m,q,p,\alpha) \in (0,\infty)$ is a constant.  Assuming $\varepsilon_0 = \varepsilon_0(n,m,q,\alpha,\rho,\tau)$ and $\beta_0 = \beta_0(n,m,q,\alpha,\rho,\tau)$ are sufficiently small, by Theorem~\ref{graphical_thm} we have that $|\xi| \leq \tau$.  Thus choosing $\rho = \rho(n,m,q,p,\alpha)$ and $\tau = \tau(n,m,q,p,\alpha,\rho)$ suitably small, 
 we deduce that conclusion~(A) holds true.

\noindent\textit{Proof of Conclusion~(B).}  Conclusion~(B) follows directly from \eqref{noncon triangle3} and conclusion~(A).

\noindent\textit{Proof of Conclusion~(C).}  By \eqref{noncon triangle3},  
\begin{align}\label{nonconB eqn1}
	&\int_{\mathbf{C}_{\gamma}(0)} \frac{\op{dist}^2(X, \op{spt} \op{graph} \varphi)}{|X-Z|^{n-\sigma}} \,d\|T\|(X) 
	\leq 16^{n-\sigma} \int_{\mathbf{C}_1(0)} \op{dist}^2(X, \op{spt} \op{graph} \varphi) \,d\|T\|(X)
	\\&\hspace{15mm}+ 2 \int_{\mathbf{B}_{1/16}(Z)} \frac{\op{dist}^2(X, Z+\Gamma(\op{spt} \op{graph} \varphi))}{|X-Z|^{n-\sigma}} \,d\|T\|(X) 
		\nonumber
	\\&\hspace{30mm}+ \int_{\mathbf{B}_{1/16}(Z)} \frac{|\chi|^2 + \widehat{E}^2 |\xi|^2 + \|A\|^2}{|X-Z|^{n-\sigma}} \,d\|T\|(X) \nonumber
\end{align}
where $C = C(n,m,q,\alpha,\sigma) \in (0,\infty)$ is a constant.  To bound the second term on the right-hand side of \eqref{nonconB eqn1}, recall that Hypothesis~$(\star)$ and Hypothesis~$(\star\star)$ hold true with $C \varepsilon_0$, $C (\beta_0+\delta)$, and $\Gamma^{-1}_{\#} \eta_{Z,1/8\#} T$, in place of $\varepsilon_0$, $\beta_0$, and $T$, where $\delta = \delta(n,m,q,\alpha) \in (0,1)$ is a suitably small and $C = C(n,m,q,\alpha) \in (0,\infty)$ is a constant.  Hence we can apply Theorem~\ref{radial noncon thm} with $\Gamma^{-1}_{\#} \eta_{Z,1/8\#} T$ in place of $T$ to obtain 
\begin{align*}
	&\int_{\mathbf{C}_{1/16}(Z,P_Z)} \frac{\op{dist}^2(X,Z+\Gamma(\op{spt}\op{graph}\varphi))}{|X-Z|^{n+2\alpha-\sigma}} \,d\|T\|(X)
	\\ \leq\,& C \int_{\mathbf{C}_{1/8}(Z,P_Z)} \op{dist}^2(X,Z+\Gamma(\op{spt}\op{graph}\varphi)) \,d\|T\|(X) 
		+ C \widehat{E}(T,P_Z,\mathbf{C}_{1/8}(Z,P_Z))^{2+\varsigma} 
\end{align*}
where $C = C(n,m,q,p,\alpha,\sigma) \in (0,\infty)$ are constants.  Recalling that $\widehat{E}(T,P_Z,\mathbf{C}_{1/8}(Z,P_Z)) \leq C\widehat{E}$ and using Theorem~\ref{nonconcentration thm}(B) 
\begin{equation}\label{nonconB eqn2}
	\int_{\mathbf{C}_{1/16}(Z,P_Z)} \frac{\op{dist}^2(X,Z+\Gamma(\op{spt}\op{graph}\varphi))}{|X-Z|^{n+2\alpha-\sigma}} \,d\|T\|(X)
	\leq C E_{\varsigma/2}(T,\varphi,\mathbf{C}_1(0)) , 
\end{equation}
where $C = C(n,m,q,p,\alpha,\sigma) \in (0,\infty)$ are constants.  To bound the third term on the right-hand side of \eqref{nonconB eqn1}, by the monotonicity formula $\|T\|(\mathbf{B}_{\rho}(Z)) \leq (q+1/2) \,\omega_n (2\rho)^n$ for all $\rho \in (0,1/16]$.  Hence 
\begin{equation}\label{nonconB eqn3}
	\int_{\mathbf{B}_{1/16}(Z)} \frac{1}{|X-Z|^{n-\sigma}} \,d\|T\|(X)
	\leq C \sum_{k=0}^{\infty} 2^{k(n-\sigma)} \|T\|(\mathbf{B}_{2^{-k-4}}(Z)) 
	\leq C \sum_{k=0}^{\infty} 2^{-k\sigma} \leq C , 
\end{equation}
where $C = C(n,m,q,p,\alpha,\sigma) \in (0,\infty)$ are constants.  Thus by Theorem~\ref{nonconcentration thm}(A) 
\begin{equation}\label{nonconB eqn4}
	\int_{\mathbf{B}_{1/16}(Z)} \frac{|\chi|^2 + \widehat{E}^2 |\xi|^2 + \|A\|^2}{|X-Z|^{n-\sigma}} \,d\|T\|(X)
	\leq C (|\chi|^2 + \widehat{E}^2 |\xi|^2 + \|A\|^2) \leq C E_{\varsigma/2}(T,\varphi,\mathbf{C}_1(0)) , 
\end{equation}
where $C = C(n,m,q,p,\alpha,\sigma) \in (0,\infty)$ are constants.  Combining \eqref{nonconB eqn1}, \eqref{nonconB eqn2}, and \eqref{nonconB eqn4} gives us conclusion~(C).

\noindent\textit{Proof of Conclusion~(D).}  By \eqref{noncon triangle1} 
\begin{align}\label{nonconD eqn1}
	&\int_{\mathbf{C}_{\gamma}(0) \cap \{r > \tau\}} \frac{\op{dist}^2(X,\op{spt}\op{graph}\widehat{\varphi})}{|X-Z|^{n+2\alpha-\sigma}} \,d\|T\|(X)
	\\ \leq\,& \int_{\mathbf{C}_{1/4}(Z)} \frac{\op{dist}^2(X, Z+\Gamma(\op{spt} \op{graph} \varphi))}{|X-Z|^{n+2\alpha-\sigma}} \,d\|T\|(X) \nonumber 
	\\&+ C \int_{\mathbf{C}_{1/4}(Z) \cap \{r > \tau\}} \frac{\widehat{E}^2 \|A\|^2 + \widehat{E}^2 |x|^{2\alpha-4} |\xi|^4}{|X-Z|^{n+2\alpha-\sigma}} \,d\|T\|(X) 
		\nonumber
	\\&+ C \int_{(\mathbf{C}_{\gamma}(0) \setminus \mathbf{C}_{1/4}(Z)) \cap \{r > \tau\}} (\op{dist}^2(X, \op{spt} \op{graph} \varphi) 
		+ |\chi|^2 + \|A\|^2 + \widehat{E}^2 |\xi|^2) \,d\|T\|(X) , \nonumber 
\end{align}
where $C = C(n,m,q,\alpha,\gamma,\sigma) \in (0,\infty)$ are constants.  The first term on the right-hand side of \eqref{nonconD eqn1} is bounded using \eqref{nonconB eqn2}.  To bound the second term on the right-hand side of \eqref{nonconD eqn1}, take $\varepsilon_0$ small enough that by conclusion~(A) we have that $\widehat{E} \tau^{-\alpha} \leq \varepsilon_0 \tau^{-\alpha} < 1$ and $|\xi| \leq \min\{\tau/2,\tau^2\}$.  Hence $|X-Z| \geq \tau/2$ if $X = (x,y,z)$ with $|x| > \tau$ and $|\xi|^2/\tau^4 < 1$.  Thus using \eqref{nonconB eqn3} and conclusion~(A) 
\begin{align}\label{nonconD eqn2}
	&\int_{\mathbf{C}_{1/4}(Z) \cap \{r > \tau\}} \frac{\widehat{E}^2 \|A\|^2 + \widehat{E}^2 |x|^{2\alpha-4} |\xi|^4}{|X-Z|^{n+2\alpha-\sigma}} \,d\|T\|(X)
	\\ \leq\,&16 \widehat{E}^2 (\tau^{-2\alpha} \|A\|^2 + \tau^{-4} |\xi|^4) \int_{\mathbf{C}_{1/4}(Z)} \frac{1}{|X-Z|^{n-\sigma}} \,d\|T\|(X) \nonumber
	\\ \leq\,& C (\|A\|^2 + \widehat{E}^2 |\xi|^2) \leq C E_{\varsigma/2}(T,\varphi,\mathbf{C}_1(0)) , \nonumber
\end{align}
where $C = C(n,m,q,\alpha,\sigma) \in (0,\infty)$ are constants.  The last term on the right-hand side \eqref{nonconD eqn1} is bounded again using conclusion~(A).  Therefore, by combining \eqref{nonconD eqn1}, \eqref{nonconD eqn2}, and conclusion~(A) we obtain conclusion~(D).
\end{proof}

\begin{corollary}\label{nonconcentration cor} 
Let $\alpha = \ell_0/q_0 \neq 2$ where $\ell_0$ and $q_0$ are relatively prime positive integers with $q_0 \leq q$.  Let $\varsigma = \varsigma(n,m,q) \in (0,1)$ be as in Theorem~\ref{relharm separation3 thm}.  For every $\delta \in (0,1/4)$ and $\sigma \in (0,1)$ there exists $\varepsilon_0 = \varepsilon_0(n,m,q,\alpha,\delta,\sigma) \in (0,1)$ and $\beta_0 = \beta_0(n,m,q,\alpha,\delta,\sigma) \in (0,1)$ such that the following holds true.  Let $\varphi \in \Phi_{\alpha,q}$ and $T$ satisfy Hypothesis~$(\star)$ and Hypothesis~$(\star\star)$.  Let $\mathcal{B}_q$ be as in Definition~\ref{introduction B defn} (with $\mu = 1/q$).  Assume that 
\begin{equation*}
	\mathbf{C}_{\delta}(0,y_0) \cap \{ X \in \mathcal{B}_q : \mathcal{N}_{T,{\rm Pl}}(X) \geq \alpha \} \neq \emptyset 
\end{equation*}
for all $y \in B^{n-2}_{1/2}(0)$, then 
\begin{equation}\label{noncon cor concl}
	\int_{\mathbf{C}_{1/2}(0)} \frac{\op{dist}^2(X, \op{spt} \op{graph} \varphi)}{r_{\delta}^{2-\sigma}} \,d\|T\|(X) \leq C E_{\varsigma/2}(T,\varphi,\mathbf{C}_1(0))^2 , 
\end{equation}
where $r_{\delta}(x,y,z) = \max\{|x|,\delta\}$ and $C = C(n,m,q,\alpha,\sigma) \in (0,\infty)$ is a constant (independent of $\delta$).
\end{corollary}

\begin{proof}
This follows from Theorem~\ref{nonconcentration thm}(C) by arguing as in the proof of~\cite[Corollary~3.2]{Sim93}.
\end{proof}

\section{Fine blow-up class}\label{sec:fine blowups} 

\subsection{Definition of the fine blow-up class and notation}  Throughout this section, we let $q \geq 2$ be an integer, $\alpha = \ell_0/q_0 \in (1,2)$ for relatively prime positive integers $\ell_0,q_0$ with $q_0 \leq 2$, and $\varsigma = \varsigma(n,m,q) \in (0,1)$ is as in Theorem~\ref{relharm separation3 thm}.  In Definition~\ref{blowup class defn} below, we define $\mathfrak{B}$ to be the class of all blow-ups $w = (w_j)$ of a sequence $(T^{(\nu)})$ of $n$-dimensional locally area-minimizing rectifiable currents of $\mathbf{C}_1(0)$ relative to a sequence $(\varphi^{(\nu)}) \subset \Phi_{\alpha,q}$ as constructed in Subsection~\ref{sec:blowup subsec}.  $\mathfrak{B}$ is based on a blow-up class introduced in~\cite[Section~9]{KrumWic2} for locally Dirichlet energy minimizing multi-valued functions, which can be regarded ``linear theory'' for locally area-minimizing rectifiable currents at branch points.  For some sequences $\varepsilon_{\nu} \rightarrow 0^+$, $\beta_{\nu} \rightarrow 0^+$, and $\delta_{\nu} \rightarrow 0^+$ and for all $\nu = 1,2,3,\ldots$ we assume that Hypothesis~$(\star)$ and Hypothesis~$(\star\star)$ hold true with $\varepsilon_{\nu},\beta_{\nu},\varphi^{(\nu)},T^{(\nu)}$ in place of $\varepsilon_0,\beta_0,\varphi,T$ and the no $\delta_{\nu}$-gaps condition \eqref{blowups no small gaps} below holds true.  More precisely:

\begin{definition}\label{domain class defn}
$\mathfrak{D}$ is the set of all ordered collections $(m_j,\varphi^{(0)}_j,\varphi^{(\infty)}_j)_{1 \leq j \leq J}$ where: 
\begin{enumerate}[itemsep=3mm,topsep=0mm]
	\item[(I)]  $J$, $m_1,m_2,\ldots,m_J$, and $q_1,q_2,\ldots,q_J$ are positive integers such that $\sum_{j=1}^J m_j q_j = q$; 
	\item[(II)]  for each $j \in \{1,2,\ldots,J\}$, $\varphi^{(0)}_j : \mathbb{R}^n \rightarrow \mathcal{A}_{q_j}(\mathbb{R}^m)$ and $\varphi^{(\infty)}_j : \mathbb{R}^n \rightarrow \mathcal{A}_{q_j}(\mathbb{R}^m)$ are $q_j$-valued functions such that either: 
	\begin{enumerate}[itemsep=3mm,topsep=3mm]
		\item[(a)]  $q_j = q_0$ and $\varphi^{(0)}_j$ and $\varphi^{(\infty)}_j$ are $q_0$-valued functions defined by 
		\begin{equation*}
			\varphi^{(0)}_j(x,y) = \op{Re}(c^{(0)}_j (x_1+ix_2)^{\alpha}), \quad \varphi^{(\infty)}_j(x,y) = \op{Re}(c^{(0)}_j (x_1+ix_2)^{\alpha})
		\end{equation*}
		for all $x = (x_1,x_2) \in \mathbb{R}^2$ and $y \in \mathbb{R}^{n-2}$, where $c^{(0)}_j, c^{(\infty)}_j \in \mathbb{C}^m$ with $|c^{(\infty)}_j| = 1$ and $c^{(0)}_j = \gamma_j c^{(\infty)}_j$ for some $\gamma_j \in [0,\infty)$;
		
		\item[(b)]  $q_j = 1$ and $\varphi^{(0)}_j : \mathbb{R}^n \rightarrow \mathcal{A}_{q_0}(\mathbb{R}^m)$ and $\varphi^{(\infty)}_j : \mathbb{R}^n \rightarrow \mathcal{A}_{q_0}(\mathbb{R}^m)$ are single-valued functions defined by $\varphi^{(0)}(x,y) = \varphi^{(\infty)}(x,y) = 0$ for all $(x,y) \in \mathbb{R}^n$.
	\end{enumerate}
\end{enumerate}
\end{definition}

\begin{definition}\label{blowup class defn} 
$\mathfrak{B}$ is the set of all ordered collections $w = (w_j)_{1 \leq j \leq J}$ of functions 
$$w_j : \op{graph}\varphi^{(\infty)}_j |_{B_1(0) \setminus \{0\} \times \mathbb{R}^{n-2}} \rightarrow \mathcal{A}_{m_j}(\mathbb{R}^m)$$
such that $w = (w_j)$ is the blow-up of a sequence $(T^{(\nu)})_{\nu=1}^{\infty}$ relative to a sequence $(\varphi^{(\nu)})_{\nu=1}^{\infty}$ by the excess $E_{\nu} = E_{\varsigma/2,\varepsilon_{\nu}}(T^{(\nu)},\varphi^{(\nu)},\mathbf{C}_1(0))$ (as in \eqref{graphical E defn}) as constructed in Subsection~\ref{sec:blowup subsec}, where for some sequence $(\varepsilon_{\nu})_{\nu=1}^{\infty}$, $(\beta_{\nu})_{\nu=1}^{\infty}$, and $(\delta_{\nu})_{\nu=1}^{\infty}$ tending to zero:
\begin{enumerate}[itemsep=3mm,topsep=0mm]
	\item[(A)]  $p \in \{1,2,\ldots,\lfloor q/q_0 \rceil\}$ (independent of $\nu$) and $\varphi^{(\nu)} \in \Phi_{\alpha,q,p}$;
	
	\item[(B)]  $T^{(\nu)}$ are $n$-dimensional locally area-minimizing rectifiable currents of $\mathbf{C}_1(0)$ such that 
	\begin{gather*}
		(\partial T^{(\nu)}) \llcorner \mathbf{C}_1(0) = 0, \quad\quad \sup_{X \in \op{spt} T^{(\nu)}} \op{dist}(X,P_0) < \varepsilon_{\nu}, \\
		\|T^{(\nu)}\|(\mathbf{C}_1(0)) < (q+1/2) \,\omega_n, \quad\quad \pi_{\#} T^{(\nu)} = q \llbracket B_1(0) \rrbracket; 
	\end{gather*}
	
	\item[(C)]  
	\begin{equation}\label{blowups hypothesis2}
		Q_{\varsigma/2,\varepsilon_{\nu}}(T^{(\nu)},\varphi^{(\nu)},\mathbf{C}_1(0)) 
			\leq \beta_{\nu} \inf_{\varphi' \in \bigcup_{p'=0}^{p^{(\nu)}-1} \Phi_{\alpha,q,p'}} 
			Q_{\varsigma/2,\varepsilon_{\nu}}(T^{(\nu)},\varphi',\mathbf{C}_1(0))
	\end{equation}
	(where $Q_{\varsigma/2,\varepsilon_{\nu}}(T^{(\nu)},\varphi^{(\nu)},\mathbf{C}_1(0))$ and $Q_{\varsigma/2,\varepsilon_{\nu}}(T^{(\nu)},\varphi',\mathbf{C}_1(0))$ are as in \eqref{graphical E defn});
	
	\item[(D)]  for each $(0,y_0) \in \{0\} \times \mathbb{R}^{n-2} \cap B_{1/2}(0)$
	\begin{equation}\label{blowups no small gaps}
		\mathbf{C}_{\delta_{\nu}}(0,y_0) \cap \{ X \in \mathcal{B}_q(T^{(\nu)}) : \mathcal{N}_{T^{(\nu)},{\rm Pl}}(X) \geq \alpha \} \neq \emptyset , 
	\end{equation} 
	where $\mathcal{B}_q(T^{(\nu)})$ is as in Definition~\ref{introduction B defn} (with $T = T^{(\nu)}$ and $\mu = 1/q$).
\end{enumerate}
\end{definition}

\begin{remark}\label{blowup excess rmk}{\rm 
In Subsection~\ref{sec:blowup subsec}, we constructed blow-ups using the excess $E_{\nu} = E_{\varsigma/2}(T^{(\nu)},\varphi^{(\nu)},\mathbf{C}_1(0))$.  By contrast, here we construct blow-ups with the larger excess $E_{\nu} = E_{\varsigma/2,\varepsilon_{\nu}}(T^{(\nu)},\varphi^{(\nu)},\mathbf{C}_1(0))$ (with $\varepsilon_{\nu}$ in place of $1$).  That way, the error term $\widehat{E}(T^{(\nu)},P_0,\mathbf{C}_1(0))^{2+\varsigma/2}$ is $o(E_{\nu}^2)$ (see for instance Lemma~\ref{blowup norms lemma}).
}\end{remark}

\begin{remark}\label{blowup on domain rmk}{\rm 
Let $w = (w_j) \in \mathfrak{B}$ be the blow-up of a sequence $(T^{(\nu)})$ relative to a sequence $(\varphi^{(\nu)})$ by the excess $E_{\nu} = E_{\varsigma/2,\varepsilon_{\nu}}(T^{(\nu)},\varphi^{(\nu)},\mathbf{C}_1(0))$ (as in Definition~\ref{blowup class defn}).  Express $\varphi^{(\nu)} = \sum_{j=1}^J m_j \varphi^{(\nu)}_j$ for some positive integers $J,m_1,m_2,\ldots,m_J$ and distinct functions $\varphi^{(\nu)}_j : \mathbb{R}^n \rightarrow \mathcal{A}_{q_j}(\mathbb{R}^m)$ (as in Definition~\ref{Phi-p_defn}).  As was discussed in Subsection~\ref{sec:blowup subsec}, after passing to a subsequence for each $w = (w_j) \in \mathfrak{B}$ there exists $(m_j,\varphi^{(0)}_j,\varphi^{(\infty)}_j) \in \mathfrak{D}$ such that:
\begin{enumerate}[itemsep=3mm,topsep=0mm]
	\item[(I)]  for each $j \in \{1,\ldots,J\}$, $m_j$ is the multiplicity of the component $\varphi^{(\nu)}_j$ in accordance with Definition~\ref{Phi-p_defn} (with $\varphi^{(\nu)}$ in place of $\varphi$); 
	
	\item[(II)]  for each $j \in \{1,\ldots,J\}$, precisely one of the following holds true: 
	\begin{enumerate}[itemsep=3mm,topsep=3mm]
		\item[(a)]  $q_j = q_0$, $\varphi^{(\nu)}_j$ are non-zero functions for all $\nu = 1,2,3,\ldots$, and $\varphi^{(\nu)}_j, \varphi^{(0)}_j, \varphi^{(\infty)}_j$ are $q_0$-valued functions given by 
		\begin{gather*}
			\hspace{20mm} \varphi^{(\nu)}_j(x,y) = \op{Re}(c^{(\nu)}_j(x_1+ix_2)^{\alpha}), \\ 
			\hspace{20mm} \varphi^{(0)}_j(x,y) = \op{Re}(c^{(0)}_j(x_1+ix_2)^{\alpha}), \quad 
			\varphi^{(\infty)}_j(x,y) = \op{Re}(c^{(\infty)}_j(x_1+ix_2)^{\alpha})
		\end{gather*}
		for all $x = (x_1,x_2) \in \mathbb{R}^2$ and $y \in \mathbb{R}^{n-2}$, where $c^{(\nu)} \in \mathbb{C}^m \setminus \{0\}$ and 
		\begin{equation*}
			\hspace{20mm} c^{(0)}_j = \lim_{\nu\rightarrow\infty} \frac{c^{(\nu)}_j}{\widehat{E}_{\nu}}, \quad 
			c^{(\infty)}_j = \lim_{\nu\rightarrow\infty} \frac{c^{(\nu)}_j}{|c^{(\nu)}_j|} ,
		\end{equation*}
		where we let $\widehat{E}_{\nu} = \widehat{E}(T^{(\nu)},P_0,\mathbf{C}_1(0))$; 
		
		\item[(b)]  $q_j = 1$, $\varphi^{(\nu)}_j : \mathbb{R}^n \rightarrow \mathbb{R}^m$ are identically zero for all $\nu = 1,2,3,\ldots$, and $\varphi^{(0)}_j : \mathbb{R}^n \rightarrow \mathbb{R}^m$ and $\varphi^{(\infty)}_j : \mathbb{R}^n \rightarrow \mathbb{R}^m$ are also identically zero.
	\end{enumerate}
\end{enumerate}
By Remark~\ref{graphical rmk}(2), 
\begin{equation}\label{blowup on domain rmk eqn1}
	\frac{1}{C} \leq \sum_{j=1}^J \|\varphi^{(0)}_j\|_{L^2(B_1(0))} \leq C 
\end{equation}
for some constant $C = C(n,m,q,\alpha) \in (1,\infty)$.  Note that even in case~(a) where $\varphi^{(\nu)}_j$ are non-zero functions for all $\nu = 1,2,3,\ldots$, we might have that $\varphi^{(0)}_j$ is identically zero.  However, by \eqref{blowup on domain rmk eqn1}, we know that $\varphi^{(0)}_j$ is non-zero for some $j \in \{1,2,\ldots,J\}$.
}\end{remark}

\begin{definition}\label{blowup on domain defn}
Given $(m_j,\varphi^{(0)}_j,\varphi^{(\infty)}_j) \in \mathfrak{D}$, we let $\mathfrak{B}(m_j,\varphi^{(0)}_j,\varphi^{(\infty)}_j)$ be the set of all blow-ups $w = (w_j) \in \mathfrak{B}$ associated with $(m_j,\varphi^{(0)}_j,\varphi^{(\infty)}_j)$ as in Remark~\ref{blowup on domain rmk}.  (Note that here and throughout this section, we write $(m_j,\varphi^{(0)}_j,\varphi^{(\infty)}_j)$ to mean $(m_j,\varphi^{(0)}_j,\varphi^{(\infty)}_j)_{1 \leq j \leq J}$).
\end{definition}

We have, of course, that $\mathfrak{B} = \bigcup_{(m_j,\varphi^{(0)}_j,\varphi^{(\infty)}_j) \in \mathfrak{D}} \mathfrak{B}(m_j,\varphi^{(0)}_j,\varphi^{(\infty)}_j)$.

In the next definition, we identify a general point $x = (x_1,x_2) \in \mathbb{R}^2$ with $re^{i\theta}$, where $(r,\theta)$ ($r \geq 0$, $0 \leq \theta < 2\pi$) are polar coordinates of $(x_1,x_2)$.  Recall that whenever $\varphi^{(\infty)}_j$ is non-zero, $(re^{i\theta},y,\op{Re}(c^{(\infty)}_j r^{\alpha} e^{i\alpha\theta}))$ determines a unique point on $\op{graph} \varphi^{(\infty)}_j$ for each $r \geq 0$ and $\theta \in [0,2\pi q_0)$. 

\begin{definition}\label{homogeneous class defn}  
Given $(m_j,\varphi^{(0)}_j,\varphi^{(\infty)}_j) \in \mathfrak{D}$, we let $\widetilde{\mathfrak L}(m_j,\varphi^{(0)}_j,\varphi^{(\infty)}_j)$ be the set of all ordered collections $\psi = (\psi_j)_{1 \leq j \leq J}$ of functions 
$$\psi_j : \op{graph}\varphi^{(\infty)}_j |_{\mathbb{R}^n \setminus \{0\} \times \mathbb{R}^{n-2}} \rightarrow \mathcal{A}_{m_j}(\mathbb{R}^m)$$
such that: 
\begin{enumerate}[itemsep=3mm,topsep=0mm]
	\item[(a)]  if $q_j = q_0$ and $\varphi^{(\infty)}_j$ is non-zero, then 
	\begin{equation}\label{homogeneous class eqn1}
		\psi_j(re^{i\theta},y,\op{Re}(c^{(\infty)}_j r^{\alpha} e^{i\alpha\theta})) 
			= \sum_{h=1}^{m_j} \llbracket \op{Re}(a_{j,h} r^{\alpha} e^{i\alpha\theta} + \alpha c^{(0)}_j r^{\alpha-1} e^{i(\alpha-1)\theta} (b \cdot y)) \rrbracket 
	\end{equation}
	for every $r > 0$, $\theta \in [0,2\pi q_0)$, and $y \in \mathbb{R}^{n-2}$, where $a_{j,h} \in \mathbb{C}^m$ and $b \in \mathbb{C}^{n-2}$ with $b$ independent of $j$ and $h$; 
	
	\item[(b)]  if $q_j = 1$ and $\varphi^{(0)}_j$ and $\varphi^{(\infty)}_j$ are both identically zero, then 
	\begin{equation}\label{homogeneous class eqn2}
		\psi_j(re^{i\theta},y,0) = (m_j - q_0 N_j) \llbracket 0 \rrbracket 
			+ \sum_{h=1}^{N_j} \sum_{l=0}^{q_0-1} \llbracket \op{Re}(a_{j,h} r^{\alpha} e^{i\alpha(\theta+2\pi l)}) \rrbracket 
	\end{equation}
	for every $r > 0$, $\theta \in [0,2\pi)$, and $y \in \mathbb{R}^{n-2}$, where $N_j$ is an integer with $0 \leq N_j \leq \lfloor m_j/q_0 \rfloor$ and $a_{j,h} \in \mathbb{C}^m \setminus \{0\}$.  (If $N_j = 0$, then $\psi_j(re^{i\theta},y,0) = m_j \llbracket 0 \rrbracket$.  If $N_j = m_j/q_0$, then $m_j - q_0 N_j = 0$ and thus we omit the term $(m_j - q_0 N_j) \llbracket 0 \rrbracket$ on the right-hand side of \eqref{homogeneous class eqn2}.)
\end{enumerate}
We let $\mathfrak{L}(m_j,\varphi^{(0)}_j,\varphi^{(\infty)}_j)$ be the set of all ordered collections of functions $\psi = (\psi_j)_{1 \leq j \leq J}$ (as above) such that $\psi_j$ is given by \eqref{homogeneous class eqn1} and \eqref{homogeneous class eqn2} with $b = 0$.  We let $\widetilde{\mathfrak L} = \bigcup_{(m_j,\varphi^{(0)}_j,\varphi^{(\infty)}_j) \in \mathfrak{D}} \widetilde{\mathfrak L}(m_j,\varphi^{(0)}_j,\varphi^{(\infty)}_j)$ and $\mathfrak{L} = \bigcup_{(m_j,\varphi^{(0)}_j,\varphi^{(\infty)}_j) \in \mathfrak{D}} \mathfrak{L}(m_j,\varphi^{(0)}_j,\varphi^{(\infty)}_j)$. 
\end{definition}

\begin{remark}\label{L scale varphi0 rmk}
If $\psi \in \widetilde{\mathfrak L}(m_j,\varphi^{(0)}_j,\varphi^{(\infty)}_j)$, then $\psi \in \widetilde{\mathfrak L}(m_j,c\varphi^{(0)}_j,\varphi^{(\infty)}_j)$ for all $c \in (0,\infty)$.  In particular, the expressions \eqref{homogeneous class eqn1} and \eqref{homogeneous class eqn2} for $\psi_j$ are unchanged if we replace $\varphi^{(0)}$ and $b$ with $c\varphi^{(0)}$ and $b/c$.
\end{remark}

Let $(m_j,\varphi^{(0)}_j,\varphi^{(\infty)}_j) \in \mathfrak{D}$.  Following the above notation, for each open ball $B \subset\subset B_1(0) \setminus \{0\} \times \mathbb{R}^{n-2}$, there are single-valued harmonic functions $\varphi^{(0)}_{j,l} : B \rightarrow \mathbb{R}^m$ and $\varphi^{(\infty)}_{j,l} : B \rightarrow \mathbb{R}^m$ such that 
\begin{equation}\label{varphi0 localized}
	\varphi^{(0)}_j(x,y) = \sum_{l=1}^{q_j} \llbracket \varphi^{(0)}_{j,l}(x,y) \rrbracket , \quad 
	\varphi^{(\infty)}_j(x,y) = \sum_{l=1}^{q_j} \llbracket \varphi^{(\infty)}_{j,l}(x,y) \rrbracket 
\end{equation}
for all $(x,y) \in B$.  We shall adopt the convention that $\varphi^{(0)}_{j,l}(x,y) = \gamma_j \varphi^{(\infty)}_{j,l}(x,y)$ for all $(x,y) \in B$, where $\gamma_j \in [0,\infty)$ is the constant such that $c^{(0)}_j = \gamma_j c^{(\infty)}_j$ (see Definition~\ref{domain class defn} and note that $\gamma_j = \lim_{\nu\rightarrow\infty} |c^{(\nu)}_j|/\widehat{E}_{\nu}$).  Moreover, if $(\varphi^{(\nu)})$ satisfies the requirements of Remark~\ref{blowup on domain rmk}, we shall adopt the convention that 
\begin{equation}\label{varphi nu localized}
	\varphi^{(\nu)}_j(x,y) = \sum_{l=1}^{q_j} \llbracket \varphi^{(\nu)}_{j,l}(x,y) \rrbracket 
\end{equation}
for all $(x,y) \in B$ (as in \eqref{varphi_localized}), where $\varphi^{(\nu)}_{j,l} : B \rightarrow \mathbb{R}^m$ are single-valued harmonic functions such that if $j \in \{1,2,\ldots,J\}$ such that $\varphi^{(\nu)}_j$ is non-zero for all $\nu$, then $\varphi^{(\nu)}_{j,l}/\widehat{E}_{\nu} \rightarrow \varphi^{(0)}_j$ and $\varphi^{(\nu)}_{j,l}/|c^{(\nu)}_j| \rightarrow \varphi^{(\infty)}_j$ uniformly on $B$ as $\nu\rightarrow\infty$.  For each $w = (w_j) \in \mathfrak{B}(m_j,\varphi^{(0)}_j,\varphi^{(\infty)}_j)$, there exists $m_j$-valued functions $w_{j,l} : B \rightarrow \mathcal{A}_{m_j}(\mathbb{R}^m)$ given by  
\begin{equation}\label{w localized}
	w_{j,l}(x,y) = w_j(x,y,\varphi^{(\infty)}_{j,l}(x,y)) = \sum_{h=1}^{m_j} \llbracket w_{j,l,h}(x,y) \rrbracket 
\end{equation}
for all $(x,y) \in B$.  Since $w = (w_j)$ is component-wise Dirichlet energy minimizing (as in Definition~\ref{v comp minimizing defn}), each function $w_{j,l}$ is locally Dirichlet energy minimizing in $B$.  Moreover, if $B$ does not intersect the $(n-2)$-dimensional singular set of $w_{j,l}$, each $w_{j,l,h}$ can be chosen to be smooth single-valued harmonic functions $w_{j,l,h} : B \rightarrow \mathbb{R}^m$.  For each Lebesgue measurable subset $\Omega \subset B_1(0)$ we let 
\begin{equation*}
	\|w\|_{L^2(\Omega)}^2 = \int_{\Omega} \sum_{j=1}^J \sum_{l=1}^{q_j} \sum_{h=1}^{m_j} |w_{j,l,h}(x,y)|^2 \,dx \,dy , 
\end{equation*}
where $w_{j,l,h}$ are locally defined by \eqref{w localized}.  Similarly, For each $\psi = (\psi_j) \in \widetilde{\mathfrak L}(m_j,\varphi^{(0)}_j,\varphi^{(\infty)}_j)$, there exists $m_j$-valued functions $\psi_{j,l} : B \rightarrow \mathcal{A}_{m_j}(\mathbb{R}^m)$ and single-valued harmonic functions $\psi_{j,l,h} : B \rightarrow \mathbb{R}^m$ such that  
\begin{equation}\label{psi localized}
	\psi_{j,l}(x,y) = \psi_j(x,y,\varphi^{(\infty)}_{j,l}(x,y)) = \sum_{h=1}^{m_j} \llbracket \psi_{j,l,h}(x,y) \rrbracket 
\end{equation}
for all $(x,y) \in B$.  For each Lebesgue measurable subset $\Omega \subset B_1(0)$ we let 
\begin{equation*}
	\|\psi\|_{L^2(\Omega)}^2 = \int_{\Omega} \sum_{j=1}^J \sum_{l=1}^{q_j} \sum_{h=1}^{m_j} |\psi_{j,l,h}(x,y)|^2 \,dx \,dy , 
\end{equation*}
where $\psi_{j,l,h}$ are locally defined by \eqref{psi localized}.

\subsection{Main estimates for fine blow-ups}  

\begin{lemma}\label{integrability lemma} 
Let $(m_j,\varphi^{(0)}_j,\varphi^{(\infty)}_j) \in \mathfrak{D}$ and let $\psi = (\psi_j) \in \widetilde{\mathfrak L}(m_j,\varphi^{(0)}_j,\varphi^{(\infty)}_j)$.  For each $\nu = 1,2,3,\ldots$ let $\varphi^{(\nu)} \in \Phi_{\alpha,q,p}$ satisfies the requirements of Remark~\ref{blowup on domain rmk}, and in particular $\varphi^{(\nu)}/\widehat{E}_{\nu} \rightarrow \varphi^{(0)}$ for some sequence $(\widehat{E}_{\nu})$ of positive numbers with $\widehat{E}_{\nu} \rightarrow 0^+$.  Let $(E_{\nu})$ be a sequence of positive numbers with $E_{\nu} \rightarrow 0^+$ and $E_{\nu} \leq \beta_{\nu} \widehat{E}_{\nu}$, where $\beta_{\nu} \rightarrow 0$.  For each positive integer $\nu$, there exists a rotation $\Gamma_{\nu}$ of $\mathbb{R}^n$ and $\widetilde{\varphi}^{(\nu)} : \mathbb{R}^n \rightarrow \mathcal{A}_q(\mathbb{R}^m)$ such that $\widetilde{\varphi}^{(\nu)} \circ \Gamma_{\nu}^{-1} \in \Phi_{\alpha,q}$, $\|\Gamma_{\nu} - I\| \leq C(n,q,\alpha) \,E_{\nu} \widehat{E}_{\nu}^{-1} \|\psi\|_{L^2(B_1(0))}$ (where $I$ is the identity map on $\mathbb{R}^n$), and:
\begin{enumerate}[itemsep=3mm,topsep=0mm]
	\item[(i)]  for any ball $B \subset \{ (x,y) \in \mathbb{R}^n : |x| \geq 4 E_{\nu} |b|/\widehat{E}_{\nu} \}$ 
	\begin{equation}\label{integrability concl1}
		\widetilde{\varphi}^{(\nu)}(x,y) = \sum_{j=1}^J \sum_{l=1}^{q_j} \sum_{h=1}^{m_j} \llbracket \varphi^{(\nu)}_{j,l}(x,y) 
			+ E_{\nu} \psi_{j,l,h}(x,y) + \mathcal{R}^{(\nu)}_{j,l,h}(x,y) \rrbracket
	\end{equation}
	for all $(x,y) \in B$ with $|x| \geq 4E_{\nu} |b|/\widehat{E}_{\nu}$, where $\mathcal{R}^{(\nu)}_{j,l,h}(x,y) \in \mathbb{R}^m$ such that 
	\begin{equation}\label{integrability concl2}
		|\mathcal{R}^{(\nu)}_{j,l,h}(x,y)| \leq C E_{\nu} |b| \,|c^{(\nu)}_j/\widehat{E}_{\nu}  - c^{(0)}_j| \,|x|^{\alpha-1} 
			+ C \beta_{\nu} E_{\nu} |a_{j,h}| \,|b| \,|x|^{\alpha-1} + C \beta_{\nu} E_{\nu} |b|^2 |x|^{\alpha-2} 
	\end{equation}
	if $\varphi^{(\nu)}_j$ is non-zero and 
	\begin{equation}\label{integrability concl3}
		|\mathcal{R}^{(\nu)}_{j,l,h}(x,y)| \leq C \beta_{\nu} E_{\nu} |a_{j,h}| \,|b| \,|x|^{\alpha-1} + C \beta_{\nu} E_{\nu} |b|^2 |x|^{\alpha-2} 
	\end{equation}
	if $\varphi^{(\nu)}_j$ is identically zero, where $a_{j,h}$ and $b$ are as in Definition~\ref{homogeneous class defn}, $\varphi^{(\nu)}_{j,l}$ as in \eqref{varphi nu localized}, $\psi_{j,l,h}$ is as in \eqref{psi localized}, and $C = C(n,m,q,\alpha) \in (0,\infty)$ is a constant; 
	
	\item[(ii)]  for $\tau \in (0,1/4)$ and $\nu$ sufficiently large 
	\begin{equation}\label{integrability concl4}
		\int_{B_{1/2}(0) \cap \{r < \tau\}} \mathcal{G}(\widetilde{\varphi}^{(\nu)},\varphi^{(\nu)})^2 \leq C \tau^{2\alpha} E_{\nu}^2 \|\psi\|_{L^2(B_1(0))}^2 ,
	\end{equation}
	where $C = C(n,m,q,\alpha) \in (0,\infty)$ is a constant; 	
	
	\item[(iii)]  
	\begin{equation}\label{integrability concl5}
		\int_{B_{1/2}(0)} \sum_{j=1}^J \sum_{l=1}^{q_j} \sum_{h=1}^{m_j} |\psi_{j,l,h}|^2 
			= \lim_{\nu\rightarrow\infty} \frac{1}{E_{\nu}^2} \int_{B_{1/2}(0)} \mathcal{G}(\widetilde{\varphi}^{(\nu)},\varphi^{(\nu)})^2 . 
	\end{equation}
\end{enumerate}
\end{lemma}

\begin{proof} 
Let $\psi$ be as in Definition~\ref{homogeneous class defn}.  Let 
\begin{equation}\label{integrability eqn1}
	B = \left[\begin{array}{ccc} 
		0 & 0 & \op{Re}(b) \\
		0 & 0 & \op{Im}(b) \\
		\op{Re}(b^T) & \op{Im}(b^T) & 0
	\end{array}\right]
\end{equation}
where $b$ is as in Definition~\ref{homogeneous class defn} represented as a row vector and $b^T$ is its transpose.  We define $\Gamma_{\nu} = e^{E_{\nu} B/\widehat{E}_{\nu}}$.  We define $\widetilde{\varphi}^{(\nu)} = \sum_{j=1}^J \widetilde{\varphi}^{(\nu)}_j$ for functions $\widetilde{\varphi}^{(\nu)}_j : \mathbb{R}^n \rightarrow \mathcal{A}_{m_j q_j}(\mathbb{R}^m)$ such that: 
\begin{enumerate}[itemsep=3mm,topsep=0mm]
	\item[(a)]  if $\varphi^{(\nu)}_j$ are non-zero for all $\nu$ and $\varphi^{(\infty)}_j$ is non-zero, then 
	\begin{equation}\label{integrability eqn2}
		\widetilde{\varphi}^{(\nu)}_j(e^{-E_{\nu} B/\widehat{E}_{\nu}} (re^{i\theta},y)) = \sum_{j=1}^J \sum_{h=1}^{m_j} \sum_{l=1}^{q_0} 
			\llbracket \op{Re}((c^{(\nu)}_j + E_{\nu} a_{j,h}) \,r^{\alpha} e^{i\alpha (\theta+2\pi l)} \rrbracket 
	\end{equation}
	for all $r > 0$, $\theta \in [0,2\pi)$, and $y \in \mathbb{R}^{n-2}$; 
	
	\item[(b)]  if $\varphi^{(\nu)}_j$ is identically zero for all $\nu$ and $\varphi^{(0)}_j$ and $\varphi^{(\infty)}_j$ are both identically zero, then $\widetilde{\varphi}^{(\nu)}_j(e^{-E_{\nu} B/\widehat{E}_{\nu}} (x,y)) = E_{\nu} \psi_j(x,y,0)$.
\end{enumerate}
The conclusion of the lemma follows by arguing as in~\cite[Lemma~9.7]{KrumWic2}, with modifications to account for the following.  
In~\cite[Lemma~9.7]{KrumWic2}, $\varphi^{(\nu)},\varphi^{(0)} \in \Phi_{\alpha,q}$ such that $\varphi^{(0)}$ is fixed and $\|\mathcal{G}(\varphi^{(\nu)},\varphi^{(0)})\|_{L^2(B_1(0))} \leq \sqrt{q} \,\varepsilon_{\nu}$ (where $\Phi_{\alpha,q}$ is as in Definition~\ref{Phi-p_defn}).  $\widetilde{\varphi}^{(\nu)}$ is constructed so that $\widetilde{\varphi}^{(\nu)} \circ e^{-E_{\nu} B} \in \Phi_{\alpha,q}$.  By contrast, here $\varphi^{(\nu)},\varphi^{(0)} \in \Phi_{\alpha,q}$ such that $C^{-1} \widehat{E}_{\nu} \leq \|\varphi^{(\nu)}\|_{L^2(B_1(0))} \leq C \widehat{E}_{\nu}$ for some constant $C = C(n,m,q,\alpha) \in (1,\infty)$ and $\varphi^{(\nu)}/\widehat{E}_{\nu} \rightarrow \varphi^{(0)}$ as $\nu\rightarrow\infty$, and $\widetilde{\varphi}^{(\nu)}$ is constructed so that $\widetilde{\varphi}^{(\nu)} \circ e^{-E_{\nu} B/\widehat{E}_{\nu}} \in \Phi_{\alpha,q}$.  
\end{proof}

\begin{lemma}\label{blowup norms lemma}
Let $(m_j,\varphi^{(0)}_j,\varphi^{(\infty)}_j) \in \mathfrak{D}$.  Let $w = (w_j) \in \mathfrak{B}(m_j,\varphi^{(0)}_j,\varphi^{(\infty)}_j)$ be the blow-up of the sequence of currents $(T^{(\nu)})$ by the sequence $(\varphi^{(\nu)})$ in $\Phi_{\alpha,q,p}$ by excess $E_{\nu} = E_{\varsigma/2,\varepsilon_{\nu}}(T^{(\nu)},\varphi^{(\nu)},\mathbf{C}_1(0))$ (as in Definition~\ref{blowup class defn}).  
\begin{enumerate}[itemsep=3mm,topsep=0mm]
	\item[(i)]  For each $\rho \in (0,1/2]$ 
	\begin{equation}\label{blowup norms concl1}
		\lim_{\nu\rightarrow\infty} \frac{E_{\varsigma/2,1}(T^{(\nu)},\varphi^{(\nu)},\mathbf{C}_{\rho}(0))^2}{E_{\nu}^2} 
		= \frac{1}{\omega_n \rho^{n+2}} \int_{B_{\rho}(0)} \sum_{j=1}^J \sum_{l=1}^{q_j} \sum_{h=1}^{m_j} |w_{j,l,h}|^2 .
	\end{equation}
	Hence, in particular $\|w\|_{L^2(B_{1/2}(0))}^2 \leq \omega_n (1/2)^{n+2}$; 
	
	\item[(ii)]  Let $\psi = (\psi_j) \in \widetilde{\mathfrak L}(m_j,\varphi^{(0)}_j,\varphi^{(\infty)}_j)$.  For $\nu = 1,2,3,\ldots$ let $\Gamma_{\nu}$ be a rotation of $\mathbb{R}^n$ and $\widetilde{\varphi}^{(\nu)} : \mathbb{R}^n \rightarrow \mathcal{A}_q(\mathbb{R}^m)$ correspond to $\psi$ as in Lemma~\ref{integrability lemma} with $\widehat{E}_{\nu} = \widehat{E}(T^{(\nu)},P_0,\mathbf{C}_1(0))$ and $E_{\nu} = E_{\varsigma/2,\varepsilon_{\nu}}(T^{(\nu)},\varphi^{(\nu)},\mathbf{C}_1(0))$.  Extend $\Gamma_{\nu}$ to a rotation of $\mathbb{R}^{n+m}$ which fixes vectors in $\{0\} \times \mathbb{R}^m$  Then each $\rho \in (0,1/2]$ 
	\begin{equation}\label{blowup norms concl2}
		\lim_{\nu\rightarrow\infty} \frac{E_{\varsigma/2,1}(\Gamma_{\nu\#} T^{(\nu)},\widetilde{\varphi}^{(\nu)} \circ \Gamma_{\nu}^{-1},\mathbf{C}_{\rho}(0))^2}{E_{\nu}^2} 
			\leq \frac{1}{\omega_n \rho^{n+2}} \int_{B_{\rho}(0)} \sum_{j=1}^J \sum_{l=1}^{q_j} \mathcal{G}(w_{j,l},\psi_{j,l})^2 , 
	\end{equation}
	and 
	\begin{align}\label{blowup norms concl3}
		&\lim_{\nu\rightarrow\infty} \frac{1}{E_{\nu}^2} \sum_{j=1}^J \int_{B_{\rho/2}(0) \cap \{r > \rho/32\}} \sum_{l=1}^{m_j q_j} 
			\op{dist}^2((x,y,\widetilde{\varphi}^{(\nu)}_{j,l}(x,y)), \op{spt} T^{(\nu)}) \,dx \,dy 
			\\&\leq \int_{B_{\rho/2}(0) \cap \{r > \rho/32\}} \sum_{j=1}^J \sum_{l=1}^{q_j} \mathcal{G}(w_{j,l},\psi_{j,l})^2 \nonumber 
	\end{align}
	  On each open ball $B \subset\subset \mathbb{R}^n \setminus \{0\} \times \mathbb{R}^{n-2}$, we let $\widetilde{\varphi}^{(\nu)}(x,y) = \sum_{j=1}^J \widetilde{\varphi}^{(\nu)}_j(x,y)$ and $\widetilde{\varphi}^{(\nu)}_j(x,y) = \sum_{l=1}^{m_j q_j} \llbracket \widetilde{\varphi}^{(\nu)}_{j,l}(x,y) \rrbracket$ for each $(x,y) \in B$, where $\widetilde{\varphi}^{(\nu)}_j : \mathbb{R}^n \rightarrow \mathcal{A}_{m_j q_j}(\mathbb{R}^m)$ is as in the proof of Lemma~\ref{integrability lemma} and $\widetilde{\varphi}^{(\nu)}_{j,l} : B \rightarrow \mathbb{R}^m$ are single-valued harmonic functions. 
\end{enumerate}
\end{lemma}

\begin{proof}
By Corollary~\ref{nonconcentration cor}, for every $\tau \in (0,1/4)$ and $\nu$ sufficiently large (depending on $\tau$) 
\begin{equation}\label{blowup norms eqn1}
	\int_{\mathbf{C}_{1/2}(0) \cap \{r < \tau\}} \op{dist}^2(X, \op{spt} \op{graph} \varphi^{(\nu)}) \,d\|T^{(\nu)}\|(X) \leq C \tau^{2-1/q} E_{\nu}^2 ,
\end{equation}
where $C = C(n,m,q,\alpha) \in (0,\infty)$ is a constant.  The conclusion now follows using \eqref{blowup norms eqn1}, the definition of the blow-up $w$, Lemma~\ref{integrability lemma}, and $\widehat{E}_{\nu}^{2+\varsigma/2}/E_{\nu}^2 \rightarrow 0$.
\end{proof}

\begin{remark}\label{rescale blowup rmk}
Let $(m_j,\varphi^{(0)}_j,\varphi^{(\infty)}_j) \in \mathfrak{D}$ and let $w = (w_j) \in \mathfrak{B}(m_j,\varphi^{(0)}_j,\varphi^{(\infty)}_j)$ be the blow-up of some sequence of currents $(T^{(\nu)})$ relative to a sequence $(\varphi^{(\nu)}) \subset \Phi_{\alpha,q,p}$ as in Definition~\ref{blowup class defn}.  Given $\rho \in (0,1)$ such that $w_j \not\equiv m_j \llbracket 0 \rrbracket$ on $\op{graph}\varphi^{(\infty)}_j |_{B_{\rho}(0)}$ for some $j \in \{1,2,\ldots,J\}$, let $\widetilde{w} = (\widetilde{w}_1,\ldots,\widetilde{w}_J)$ such that $\widetilde{w}_j : \op{graph}\varphi^{(\infty)}_j |_{B_1(0)} \rightarrow \mathcal{A}_{m_j}(\mathbb{R}^m)$ is given by 
\begin{equation*}
	\widetilde{w}_j(x,y,\varphi^{(\infty)}_{j,l}(x,y)) = \frac{w_j(\rho x, \rho y, \varphi^{(\infty)}_{j,l}(\rho x, \rho y))}{\rho^{-n/2} \|w\|_{L^2(B_{\rho}(0))}} 
\end{equation*}
for all $(x,y) \in B_1(0)$, $j \in \{1,\ldots,J\}$, and $l \in \{1,\ldots,q_j\}$, where $\varphi^{(\infty)}_{j,l}$ is as in \eqref{varphi0 localized}.  Then $\widetilde{w} \in \mathfrak{B}$.  In fact, $E_{\varsigma/2,\varepsilon_{\nu}}(T^{(\nu)},\varphi^{(\nu)},\mathbf{C}_{\rho}(0)) \neq 0$ for all sufficiently large $\nu$.  After passing to a subsequence, $\widetilde{w}$ is a fine blow-up of $(\eta_{0,\rho\#} T^{(\nu)})_{\nu=1}^{\infty}$ relative to $(\rho^{\alpha-1} \varphi^{(\nu)})_{\nu=1}^{\infty}$.  $\widetilde{w} \in \mathfrak{B}(m_j,c\varphi^{(0)}_j,\varphi^{(\infty)}_j)$ where $c = \lim_{\nu\rightarrow\infty} \rho^{\alpha-1} \widehat{E}(T^{(\nu)},P_0,\mathbf{C}_1(0))/\widehat{E}(T^{(\nu)},P_0,\mathbf{C}_{\rho}(0)) \in (0,\infty)$.  

{\rm  To see this, note that by \eqref{blowup norms concl1}, $E_{\varsigma/2,\varepsilon_{\nu}}(T^{(\nu)},\varphi^{(\nu)},\mathbf{C}_{\rho}(0)) \neq 0$ for all sufficiently large $\nu$.  Let $v^{(\nu)}_j : \op{graph}\varphi^{(\nu)}_j |_{B_{\rho}(0) \cap \{r > \rho/16\}} \rightarrow \mathbb{R}^m$ and $K^{(\nu)}_j \subset B_{\rho}(0) \cap \{r > \rho/16\}$ be as in Theorem~\ref{graphical_thm}(C) with $\rho,\rho/16,\varphi^{(\nu)},T^{(\nu)},K^{(\nu)}_j,v^{(\nu)}_j$ in place of $\gamma,\tau,\varphi,T,K_j,v_j$.  By Remark~\ref{graphical triangle rmk}(2) and Remark~\ref{graphical rmk}(1)(2) 
\begin{align}\label{rescale blowup rmk eqn1}
	\widehat{E}(T^{(\nu)},P_0,\mathbf{C}_{\rho}(0)) 
	\geq\,& C \rho^{\alpha-1} \|\varphi^{(\nu)}\|_{L^2(B_1(0))} - E_{\varsigma/2,1}(T^{(\nu)},\varphi^{(\nu)},\mathbf{C}_1(0))
	\\ \geq\,& C \rho^{\alpha-1} \widehat{E}(T^{(\nu)},P_0,\mathbf{C}_1(0)) \nonumber 
\end{align}
for all sufficiently large $\nu$, where $\varphi^{(\nu)}_{j,l}$ as in \eqref{varphi nu localized}, $u^{(\nu)}_{j,h,l}$ and $v^{(\nu)}_{j,l}$ are as in \eqref{u defn} and \eqref{v_localized} with $\varphi^{(\nu)},T^{(\nu)},v^{(\nu)}_j,u^{(\nu)}_j$ in place of $\varphi,T,v_j,u_j$, and $C = C(n,m,q,\alpha) > 0$ are constants.  Hence after passing to a subsequence, $c = \lim_{\nu\rightarrow\infty} \rho^{\alpha-1} \widehat{E}(T^{(\nu)},P_0,\mathbf{C}_1(0))/\widehat{E}(T^{(\nu)},P_0,\mathbf{C}_{\rho}(0))$ exists and is in $(0,\infty)$.  By Remark~\ref{noncon rescale rmk} (in the special case $Z = 0$), conditions~(A), (B), and (C) of Definition~\ref{blowup class defn} hold true with $C\varepsilon_{\nu}, C\beta_{\nu}, \eta_{0,\rho\#} T^{(\nu)}, \rho^{1-\alpha} \varphi^{(\nu)}$ in place of $\varepsilon_{\nu}, \beta_{\nu}, T^{(\nu)}, \varphi^{(\nu)}$, where $C = C(n,m,q,\alpha,\rho) \in (0,\infty)$ is a constant.  Using condition (D) of Definition~\ref{blowup class defn}, it is readily verified that condition (D) of Definition~\ref{blowup class defn} holds true with $\delta_{\nu}/\rho, \eta_{0,\rho\#} T^{(\nu)}$ in place of $\delta_{\nu}, T^{(\nu)}$.  By a scaling argument (like in Remark~\ref{T and varphi scaling rmk}) and \eqref{blowup norms concl1}, $\widetilde{w}$ is a fine blow-up of $(\eta_{0,\rho\#} T^{(\nu)})_{\nu=1}^{\infty}$ relative to $(\rho^{\alpha-1} \varphi^{(\nu)})_{\nu=1}^{\infty}$ and $\widetilde{w} \in \mathfrak{B}(m_j,c\varphi^{(0)}_j,\varphi^{(\infty)}_j)$.
}\end{remark}

\begin{lemma}\label{blowup estimates lemma}
Let $(m_j,\varphi^{(0)}_j,\varphi^{(\infty)}_j) \in \mathfrak{D}$ and $w = (w_j) \in \mathfrak{B}(m_j,\varphi^{(0)}_j,\varphi^{(\infty)}_j)$.  The following hold true: 
\begin{enumerate}[itemsep=3mm,topsep=0mm]
	\item[(i)]  For each $\psi \in \widetilde{\mathfrak L}(m_j,\varphi^{(0)}_j,\varphi^{(\infty)}_j)$ and each $\rho \in (0,1/2]$ and $\gamma \in (0,1)$ 
	\begin{gather}
		\label{blowup estimates concl1} \int_{B_{\gamma\rho}(0)} \sum_{j=1}^J \sum_{l=1}^{q_j} 
			R^{2-n} \left|\frac{\partial (w_{j,l}/R^{\alpha})}{\partial R}\right|^2 
			\leq C \rho^{-n-2\alpha} \int_{B_{\rho}(0)} \sum_{j=1}^J \sum_{l=1}^{q_j} \mathcal{G}(w_{j,l},\psi_{j,l})^2  , 
			\end{gather}
	where $w_{j,l}$ is as in \eqref{w localized}, $\psi_{j,l}$ is as in \eqref{psi localized}, $R(x,y) = |(x,y)|$ 
	for each $(x,y) \in \mathbb{R}^n$, and $C = C(n,m,q,\alpha,\gamma) \in (0,\infty)$ is a constant. 
	
	\item[(ii)]  	We have that 
	\begin{equation}\label{blowup estimates concl3} 
		\int_{B_{1/2}(0)} \sum_{j=1}^J \sum_{l=1}^{q_j} \sum_{h=1}^{m_j} |D_y w_{j,l}|^2 \leq C  ,
	\end{equation}
	where 
	$C = C(n,m,q,\alpha,\gamma) \in (0,\infty)$ is a constant. 
		
	\item[(iii)]  	There exists functions $\kappa : B^{n-2}_{1/4}(0) \rightarrow \mathbb{R}^2$, $\lambda : B^{n-2}_{1/4}(0) \rightarrow \mathbb{R}^m$, and $\Lambda : B^{n-2}_{1/4}(0) \rightarrow \mathbb{R}^{m \times n}$ such that for each $\sigma \in (0,2/q)$ and each $\zeta \in B^{n-2}_{1/4}(0)$ 
	\begin{equation}\label{blowup estimates concl4} 
		\int_{B_{1/2}(0)} \sum_{j=1}^J \sum_{l=1}^{q_j} \sum_{h=1}^{m_j} \frac{|w_{j,l,h}(x,y) - \lambda(\zeta) + \Lambda(\zeta) \cdot (x,y-\zeta) 
			+ D_x \varphi^{(0)}_{j,l}(x,y) \cdot \kappa(\zeta)|^2}{|(x,y) - (0,\zeta)|^{n+2\alpha-\sigma}} \,dx \,dy \leq C , 
	\end{equation}
	where $w_{j,l,h}$ is as in \eqref{w localized}, $\varphi^{(0)}_{j,l}$ is as in \eqref{varphi0 localized}, $\mathbb{R}^{m \times n}$ denotes the space of linear maps from $\mathbb{R}^n$ to $\mathbb{R}^m$, and $C = C(n,m,q,\alpha,\sigma) \in (0,\infty)$ is a constant.   If $1 < \alpha < 2$, then $\lambda$ is differentiable and $D\lambda(\zeta) \cdot e = \Lambda(\zeta) \cdot (0,e)$ for all $\zeta \in B^{n-2}_{1/4}(0)$ and $e \in \mathbb{R}^{n-2}$ and 
	\begin{equation}\label{blowup estimates concl5} 
		\|\kappa\|_{C^{0,1-\sigma/2}(B_{1/2}(0))} + \|\lambda\|_{C^{1,\alpha-1-\sigma/2}(B_{1/2}(0))} + \|\Lambda\|_{C^{0,\alpha-1-\sigma/2}(B_{1/2}(0))} \leq C , 
	\end{equation}
	where $C = C(n,m,q,\alpha,\sigma) \in (0,\infty)$ is a constant.  If $\alpha > 2$, then $\Lambda$ is a constant, $\lambda(\zeta) = \lambda(0) + \Lambda \cdot \zeta$ for all $\zeta \in B^{n-2}_{1/4}(0)$, and $\|\kappa\|_{C^{0,1-\sigma/2}(B_{1/2}(0))} + |\lambda(0)| + |\Lambda| \leq C$, where $C = C(n,m,q,\alpha,\sigma) \in (0,\infty)$ is a constant.\vspace{3mm}
	
	Moreover, $\kappa(\zeta)$, $\lambda(\zeta)$, and $\Lambda(\zeta)$ are unique subject to the condition that for some $\sigma \in (0,2/q)$ 
	\begin{equation}\label{blowup estimates concl6} 
		\int_{B_{1/2}(0)} \sum_{j=1}^J \sum_{l=1}^{q_j} \sum_{h=1}^{m_j} \frac{|w_{j,l,h}(x,y) - \lambda(\zeta) + \Lambda(\zeta) \cdot (x,y-\zeta) 
			+ D_x \varphi^{(0)}_{j,l}(x,y) \cdot \kappa(\zeta)|^2}{|(x,y) - (0,\zeta)|^{n+2\alpha-\sigma}} \,dx \,dy < \infty . \hspace{-2mm} 
	\end{equation}
\end{enumerate}
\end{lemma}

\begin{proof}
To see \eqref{blowup estimates concl1}, let $w = (w_j) \in \mathfrak{B}(m_j,\varphi^{(0)}_j,\varphi^{(\infty)}_j)$ be a blow-up of a sequence of locally area-minimizing currents $(T^{(\nu)})$ relative to a sequence $(\varphi^{(\nu)})$ in $\Phi_{\alpha,q,p}$ by excess $E_{\nu} = E_{\varsigma/2,\varepsilon_{\nu}}(T^{(\nu)},\varphi^{(\nu)},\mathbf{C}_1(0))$ (as in Definition~\ref{blowup class defn}).  Let $\psi \in \widetilde{\mathfrak L}(m_j,\varphi^{(0)}_j,\varphi^{(\infty)}_j)$.  For $\nu = 1,2,3,\ldots$, let $\Gamma_{\nu}$ and $\widetilde{\varphi}^{(\nu)}$ be as in Lemma~\ref{integrability lemma} corresponding to $\varphi^{(\nu)}$, $\psi$, $\widehat{E}_{\nu} = \widehat{E}(T^{(\nu)},P_0,\mathbf{C}_1(0))$, and $E_{\nu}$.  Let $\widetilde{p} \in \{p,p+1,\ldots,\lfloor q/q_0 \rfloor\}$ such that $\widetilde{\varphi}^{(\nu)} \circ \Gamma_{\nu}^{-1} \in \Phi_{\alpha,q,\widetilde{p}}$.  Extend $\Gamma_{\nu}$ to a rotation of $\mathbb{R}^{n+m}$ which fixes vectors in $\{0\} \times \mathbb{R}^m$.  Fix $\rho \in (0,1/2]$, and assume that $w |_{B_{\rho}(0)}$ is not identically zero, as otherwise \eqref{blowup estimates concl1} trivially holds true.  Let $\varepsilon_0$ and $\beta_0$ be as in Theorem~\ref{keyest thm}.  By Remark~\ref{noncon rescale rmk} (in the special case $Z = 0$), Hypothesis~$(\star)$ holds true with $\widehat{T}^{(\nu)} = \Gamma_{\nu\#} \eta_{0,\rho\#} T^{(\nu)}$ and $\widehat{\varphi}^{(\nu)} = \rho^{\alpha-1} \widetilde{\varphi}^{(\nu)} \circ \Gamma_{\nu}^{-1}$ in place of $T$ and $\varphi$.  Let $\tau \in (0,\gamma\rho/2)$ be arbitrary and let $u^{(\nu)}_j$ and $v^{(\nu)}_j$ be as in  Theorem~\ref{graphical_thm}(C) with $1/2,T^{(\nu)}, \varphi^{(\nu)},u^{(\nu)}_j,v^{(\nu)}_j$ in place of $\gamma,T,\varphi,u_j,v_j$.  In the case that Hypothesis~$(\star\star)$ holds true with $\widehat{T}^{(\nu)}, \widehat{\varphi}^{(\nu)}$ in place of $T, \varphi$, 
by Theorem~\ref{keyest thm} and Theorem~\ref{graphical_thm}(C) (both with $\widehat{T}^{(\nu)}, \widehat{\varphi}^{(\nu)}$ in place of $T,\varphi$) and by the area formula, for any $\tau \in (0,\gamma\rho/2]$ and all sufficiently large $\nu$ 
\begin{align}\label{blowup estimates eqn1}
	&\int_{B_{\gamma\rho}(0) \cap \{r > \tau\}} \sum_{j=1}^J \sum_{l=1}^{q_j} \sum_{h=1}^{m_j} 
		R^{2-n} \left| \frac{\partial (u^{(\nu)}_{j,l,h}/R^{\alpha})}{\partial R} \right|^2 
	\\ \leq\,& C_1 \rho^{-n-2\alpha} \int_{\mathbf{C}_{\rho}(0)} \op{dist}^2(X,\op{spt}\op{graph}\widetilde{\varphi}^{(\nu)}) \,d\|T^{(\nu)}\|(X) 
		+ C_2 \widehat{E}_{\nu}^{2+\varsigma/2} , \nonumber 
\end{align}
where $u^{(\nu)}_{j,l,h}$ is as in \eqref{u defn} and where $C_1 = C_1(n,m,q,\alpha,\gamma) \in (0,\infty)$ and $C_2 = C_2(n,m,q,\alpha,\gamma,\tau) \in (0,\infty)$ are constants.  By \eqref{u defn} and the homogeneity of $\varphi^{(\nu)}$, $\frac{\partial}{\partial R} (u^{(\nu)}_{j,l,h}/R^{\alpha}) = \frac{\partial}{\partial R} (v^{(\nu)}_{j,l,h}/R^{\alpha})$, where $v^{(\nu)}_{j,l,h}$ is as in \eqref{v_localized} with $\varphi^{(\nu)},u^{(\nu)}_j,v^{(\nu)}_j$ in place of $\varphi,u_j,v_j$.  Let $L = L(m,q) \geq 1$ be an integer and $\boldsymbol{\xi} : \mathcal{A}_q(\mathbb{R}^m) \rightarrow \mathbb{R}^L$ be the Almgren bi-Lipschitz injection (cf.\ ~\cite[Corollary~2.2]{DeLSpaDirMin}).  Letting $w^{(\nu)}_j$ be as in Subsection~\ref{sec:blowup subsec}, it follows from \eqref{blowup eqn3} that after passing to a subsequence $\boldsymbol{\xi} \circ w^{(\nu)}_j \rightarrow \boldsymbol{\xi} \circ w_j$ weakly in $W^{1,2}(\op{graph} \varphi^{(\infty)}_j |_{\Omega'},\mathbb{R}^m)$ for each $\Omega' \subset\subset B_1(0) \setminus \{0\}$.  Using this together with Lemma~\ref{blowup norms lemma} and $\widehat{E}_{\nu}^{2+\varsigma/2}/E_{\nu}^2 \rightarrow 0$ as $\nu\rightarrow\infty$, we can divide both sides of \eqref{blowup estimates eqn1} by $E_{\nu}^2$ and let $\nu\rightarrow\infty$ to obtain \eqref{blowup estimates concl1}.  If on the other hand Hypothesis~$(\star\star)$ does not hold true with $\widehat{T}^{(\nu)}, \widehat{\varphi}^{(\nu)}$ in place of $T, \varphi$, then by Remark~\ref{tilde varphi rmk}(1), there exists $\overline{p} \in \{0,1,2,\ldots,\widetilde{p}-1\}$ and $\overline{\varphi}^{(\nu)} \in \Phi_{\alpha,q,\overline{p}}$ such that 
\begin{equation}\label{blowup estimates eqn2}
	Q_{\varsigma/2,1}(\widehat{T}^{(\nu)},\overline{\varphi}^{(\nu)},\mathbf{C}_1(0)) 
		< 2^{\overline{p}} \beta_0^{1-\overline{p}} Q_{\varsigma/2,1}(\widehat{T}^{(\nu)},\widehat{\varphi}^{(\nu)},\mathbf{C}_1(0))
\end{equation} 
and either $\overline{p} = 0$ or $\overline{p} > 0$ and Hypothesis~$(\star\star)$ holds true with $\widehat{T}^{(\nu)}, \overline{\varphi}^{(\nu)}$ in place of $T^{(\nu)}, \varphi^{(\nu)}$.  By Lemma~\ref{integrability lemma} and Lemma~\ref{blowup norms lemma} 
\begin{equation*}
	Q_{\varsigma/2,1}(\widehat{T}^{(\nu)},\rho^{\alpha-1} \widehat{\varphi}^{(\nu)},\mathbf{C}_1(0)) \leq C E_{\nu} , 
\end{equation*}
where $C = C(n,m,q,\alpha,\rho,\|\psi\|_{L^2(B_1(0))}) \in (0,\infty)$ are constants.  Since $w |_{B_{\rho}(0)}$ is not identically zero, we can use \eqref{rescale blowup rmk eqn1} together with $\hat E_{\nu} \leq \beta_{\nu} E_{\nu}$ to obtain
\begin{equation*}
	Q_{\varsigma/2,1}(\widehat{T}^{(\nu)},\rho^{\alpha-1} \widehat{\varphi}^{(\nu)},\mathbf{C}_1(0))
	\leq C \beta_{\nu} \widehat{E}_{\nu} 
	\leq C \beta_{\nu} \widehat{E}(T^{(\nu)},P_0,\mathbf{C}_{\rho}(0)) 
	= C \beta_{\nu} \widehat{E}(\widehat{T}^{(\nu)},P_0,\mathbf{C}_1(0)) , 
\end{equation*}
where $C = C(n,m,q,\alpha,\rho,\|\psi\|_{L^2(B_1(0))}) \in (0,\infty)$ are constants.  Hence by \eqref{blowup estimates eqn2}, we must have that $\overline{p} \geq 1$.  By applying Theorem~\ref{keyest thm} and Theorem~\ref{graphical_thm}(C) (both with $\widehat{T}^{(\nu)}, \overline{\varphi}^{(\nu)}$ in place of $T,\varphi$), the area formula, and \eqref{blowup estimates eqn2}, we deduce that 
\begin{align}\label{blowup estimates eqn3}
	&\int_{B_{\gamma\rho}(0) \cap \{r > \tau\}} \sum_{j=1}^J \sum_{l=1}^{q_j} \sum_{h=1}^{m_j} 
		R^{2-n} \left| \frac{\partial (u^{(\nu)}_{j,l,h}/R^{\alpha})}{\partial R} \right|^2 
	\\ \leq\,& C_1 \rho^{-n-2\alpha} \int_{\mathbf{C}_{\rho}(0)} \op{dist}^2(X,\op{spt}\op{graph}\widetilde{\varphi}^{(\nu)}) \,d\|T^{(\nu)}\|(X) \nonumber 
		\\&+ C_1 \rho^{-n-2\alpha} \int_{B_{\rho/2}(0) \cap \{r > \rho/32\}} \sum_{j=1}^J \sum_{l=1}^{m_j q_j} 
		\op{dist}^2((x,y,\widetilde{\varphi}^{(\nu)}_{j,l}(x,y)), \op{spt} T^{(\nu)}) \,dx \,dy + C_2 \widehat{E}_{\nu}^{2+\varsigma/2} , \nonumber 
\end{align}
 where $C_1 = C_1(n,m,q,\alpha,\gamma) \in (0,\infty)$ and $C_2 = C_2(n,m,q,\alpha,\gamma,\tau) \in (0,\infty)$.  On each open ball $B \subset\subset \mathbb{R}^n \setminus \{0\} \times \mathbb{R}^{n-2}$, we let $\widetilde{\varphi}^{(\nu)}(x,y) = \sum_{j=1}^J \widetilde{\varphi}^{(\nu)}_j(x,y) = \sum_{j=1}^J \sum_{l=1}^{m_j q_j} \llbracket \widetilde{\varphi}^{(\nu)}_{j,l}(x,y) \rrbracket$ for each $(x,y) \in B$, where $\widetilde{\varphi}^{(\nu)}_j : \mathbb{R}^n \rightarrow \mathcal{A}_{m_j q_j}(\mathbb{R}^m)$ is as in the proof of Lemma~\ref{integrability lemma} and $\widetilde{\varphi}^{(\nu)}_{j,l} : B \rightarrow \mathbb{R}^m$ is single-valued harmonic functions.  Arguing as before, we can divide both sides of \eqref{blowup estimates eqn3} by $E_{\nu}^2$ and let $\nu\rightarrow\infty$ to deduce that \eqref{blowup estimates concl1} holds true.

The estimate \eqref{blowup estimates concl3} follows from Theorem~\ref{keyest thm}, Theorem~\ref{graphical_thm}(C), the area formula, and Lemma~\ref{blowup norms lemma}.
 
To see \eqref{blowup estimates concl4}, let $\zeta \in B^{n-2}_{1/4}(0)$ and observe that by \eqref{blowups no small gaps} there exists $Z_{\nu} = (\xi_{\nu},\zeta_{\nu},\chi_{\nu}) \in \mathbf{C}_{\delta_{\nu}}(0,\zeta) \cap \mathcal{B}_q(T^{(\nu)})$ with $\mathcal{N}_{T^{(\nu)},{\rm Pl}}(Z_{\nu}) \geq \alpha$.  Let $P_{\nu}$ be the unique tangent plane to $T^{(\nu)}$ at $Z_{\nu}$ and $P_{\nu} = \{ (x,y,z) \in \mathbb{R}^{n+m} : z = A_{\nu}(x,y) \}$ for some $A_{\nu} \in \mathbb{R}^{m \times n}$.  Let $v^{(\nu)}_j$ be as in Theorem~\ref{graphical_thm}(C) with $1/2,T^{(\nu)}, \varphi^{(\nu)},v^{(\nu)}_j$ in place of $\gamma,T, \varphi,v_j$.  By Theorem~\ref{nonconcentration thm} and Theorem~\ref{graphical_thm} 
\begin{equation}\label{blowup estimates eqn4}
	\widehat{E}_{\nu} |\xi_{\nu}| + |\chi_{\nu}| + \|A_{\nu}\| \leq C E_{\nu} 
\end{equation}
and for each $\tau \in (0,1/4)$ and $\nu$ sufficiently large (depending on $\tau$) 
\begin{equation}\label{blowup estimates eqn5}
	\int_{B_{1/2}(0) \cap \{r > \tau\}} \sum_{j=1}^J \sum_{l=1}^{q_j} \sum_{h=1}^{m_j} \frac{|v^{(\nu)}_{j,l,h}(x,y) - \chi_{\nu} + A_{\nu}(x-\xi,y-\zeta) 
		+ D_x \varphi^{(\nu)}_{j,l}(x,y) \cdot \xi_{\nu}|^2}{|(x,y) - (\xi_{\nu},\zeta_{\nu})|^{n+2\alpha-\sigma}} \leq C E_{\nu}^2 ,
\end{equation}
here $v^{(\nu)}_{j,l,h}$ is as in \eqref{v_localized} with $\varphi^{(\nu)}_j, v^{(\nu)}_j$ in place of $\varphi_j,v_j$ and $C = C(n,m,q,\alpha,\sigma) \in (0,\infty)$ is a constant.  By \eqref{blowup estimates eqn4}, after passing to a subsequence we can let $\widehat{E}_{\nu} \xi_{\nu}/E_{\nu} \rightarrow \kappa$ in $\mathbb{R}^2$, $\chi_{\nu}/E_{\nu} \rightarrow \lambda$ in $\mathbb{R}^m$, and $A_{\nu}/E_{\nu} \rightarrow \Lambda$ uniformly on compact subsets of $\mathbb{R}^n$ for some $\Lambda \in \mathbb{R}^{m \times n}$.  Recall that $\varphi^{(\nu)}_j/\widehat{E}_{\nu} \rightarrow \varphi^{(0)}_j$ uniformly on compact subsets of $\mathbb{R}^n$.  Thus by also using the definition of the blow-up $w$ and Lemma~\ref{blowup norms lemma}, we can divide both sides of \eqref{blowup estimates eqn5} by $E_{\nu}^2$ and let $\nu \rightarrow \infty$ and $\tau \rightarrow 0^+$ to obtain \eqref{blowup estimates concl4} with $\kappa = \kappa(\zeta)$, $\lambda(\zeta) = \lambda$, and $\Lambda = \Lambda(\zeta)$.  

To see that there are unique $\kappa(\zeta)$, $\lambda(\zeta)$, and $\Lambda(\zeta)$ corresponding to each $\zeta \in B^{n-2}_{1/4}(0)$, suppose to the contrary that for $i = 1,2$ there exists $\kappa_i,\lambda_i,\Lambda_i$ such that \eqref{blowup estimates concl6} holds true with $\kappa_i,\lambda_i,\Lambda_i$ in place of $\kappa(\zeta), \lambda(\zeta), \Lambda(\zeta)$.  Then by the triangle inequality 
\begin{equation*}
	\sum_{j=1}^J m_j \int_{B_{1/2}(0)} \sum_{l=1}^{q_j} \frac{|\lambda_1 - \lambda_2 - (\Lambda_1 - \Lambda_2) \,(x,y-\zeta) 
		+ D_x \varphi^{(0)}_{j,l}(x,y) \cdot (\kappa_1 - \kappa_2)|^2}{|(x,y) - (0,\zeta)|^{n+2\alpha-\sigma}} \,dx \,dy < \infty .
\end{equation*}
Using the mutual $L^2(\mathbb{S}^{n-1})$-orthogonality of $\lambda_1 - \lambda_2$, $(\Lambda_1 - \Lambda_2) \,(x,y)$, and $\frac{1}{q} \sum_{j=1}^J \sum_{l=1}^q m_j D_x \varphi^{(0)}_{j,l}(x,y) \cdot (\kappa_1 - \kappa_2)$ gives us that $\kappa_1 = \kappa_2$, $\lambda_1 = \lambda_2$, and $\Lambda_1 = \Lambda_2$.   

It follows from the construction of $\kappa(\zeta)$, $\lambda(\zeta)$, and $\Lambda(\zeta)$ and \eqref{blowup estimates eqn4} 
\begin{equation}\label{blowup estimates eqn6}
	\sup_{B^{n-2}_{1/4}(0)} |\kappa| + \sup_{B^{n-2}_{1/4}(0)} |\lambda| + \sup_{B^{n-2}_{1/4}(0)} |\Lambda| \leq C 
\end{equation}
for some constant $C = C(n,m,q,\alpha) \in (0,\infty)$.  Let $\zeta_1,\zeta_2 \in B^{n-2}_{1/4}(0)$ and set $\rho = |\zeta_1 - \zeta_2|$.  Assuming $|\zeta_1 - \zeta_2| \leq 1/8$, by the triangle inequality and \eqref{blowup estimates concl4} with $\zeta = \zeta_1,\zeta_2$ 
\begin{align*}
	&\int_{B_{\rho}(0,\zeta_1)} \sum_{j=1}^J \sum_{l=1}^{q_j} |\lambda(\zeta_1) - \lambda(\zeta_2) - \Lambda(\zeta_1) \,(x,y-\zeta_1) 
		+ \Lambda(\zeta_2) \,(x,y-\zeta_2) \\&\hspace{30mm} + D_x \varphi^{(0)}_{j,l}(x,y) \cdot (\kappa(\zeta_1) - \kappa(\zeta_2))|^2 \,dx \,dy 
	\leq C \rho^{n+2\alpha-\sigma} ,
\end{align*}
where $C = C(n,m,q,\alpha,\sigma) \in (0,\infty)$ is a constant.  By \eqref{blowup on domain rmk eqn1} and by again using the mutual $L^2(\mathbb{S}^{n-1})$-orthogonality of constant functions, $(\Lambda(\zeta_1) - \Lambda(\zeta_1)) \,(x,y)$, and $\frac{1}{q} \sum_{j=1}^J \sum_{l=1}^q m_j D_x \varphi^{(0)}_{j,l}(x,y) \cdot (\kappa(\zeta_1) - \kappa(\zeta_1))$ 
\begin{equation}\label{blowup estimates eqn7}
	|\lambda(\zeta_1) - \lambda(\zeta_2) - \Lambda(\zeta_2) \cdot (0,\zeta_1-\zeta_2)| + \rho \,|\Lambda(\zeta_1) - \Lambda(\zeta_2)| 
	+ \rho^{\alpha-1} \,|\kappa(\zeta_1) - \kappa(\zeta_2)| \leq C \rho^{\alpha-\sigma/2} 
\end{equation}
for some constant $C = C(n,m,q,\alpha,\sigma) \in (0,\infty)$.  If $1 < \alpha < 2$, then it follows from \eqref{blowup estimates eqn6} and \eqref{blowup estimates eqn7} that $\lambda$ is differentiable with $D\lambda(\zeta) \cdot e = \Lambda(\zeta) \cdot (0,e)$ for all $\zeta \in B^{n-2}_{1/4}(0)$ and $e \in \mathbb{R}^{n-2}$ and that \eqref{blowup estimates concl4} holds true.  If $\alpha > 2$, then we reach a similar conclusion except \eqref{blowup estimates eqn7} also implies that $\Lambda$ is differentiable with $D\Lambda(\zeta) = 0$ for all $\zeta \in B^{n-2}_{1/4}(0)$, i.e.~$\Lambda$ is constant and consequently $\lambda$ is affine.
\end{proof}

\subsection{Classification of homogeneous degree $\alpha$ blow-ups}  

\begin{theorem}\label{classify in L thm}
Let $w = (w_j) \in \mathfrak{B}$ be homogeneous degree $\alpha$ in the sense that 
\begin{equation*}
	(x,y) \,\mapsto\, \sum_{l=1}^{q_j} \sum_{h=1}^{m_j} \llbracket w_{j,l,h}(x,y) \rrbracket 
\end{equation*}
is a homogeneous degree $\alpha$, $\mathcal{A}_{m_j q_j}(\mathbb{R}^m)$-valued function of $(x,y) \in B_{1/2}(0) \setminus \{0\} \times \mathbb{R}^{n-2}$ for all $j \in \{1,2,\ldots,J\}$, where $w_{j,l,h}$ is as in \eqref{w localized}.  Then $w \in \widetilde{\mathfrak L}$.
\end{theorem}

\begin{proof}
Recall from Lemma~\ref{blowup estimates lemma} that there exists a functions $\kappa : B^{n-2}_{1/4}(0) \rightarrow \mathbb{R}^2$, $\lambda : B^{n-2}_{1/4}(0) \rightarrow \mathbb{R}^m$, and $\Lambda : B^{n-2}_{1/4}(0) \rightarrow \mathbb{R}^{m \times n}$ such that \eqref{blowup estimates concl4} holds true for all $\zeta \in B^{n-2}_{1/4}(0)$ and $\sigma \in (0,1/q)$. 

\noindent\textit{Step~1: Show that $\lambda = 0$ and $\Lambda =0$.}  For each $j \in \{1,2,\ldots,J\}$, define $\overline{w}_j : B_{1/2}(0) \rightarrow \mathcal{A}_{q_j m_j}(\mathbb{R}^m)$ by 
\begin{equation*}
	\overline{w}_j(x,y) = \sum_{l=1}^{q_j} \sum_{h=1}^{m_j} \llbracket w_{j,l,h}(x,y) \rrbracket 
\end{equation*}
for all $(x,y) \in B_{1/2}(0)$.   Let $\overline{w}_{j,a}(x,y) = \frac{1}{q_j m_j} \sum_{l=1}^{q_j} \sum_{h=1}^{m_j} w_{j,l,h}(x,y)$ be the single-valued average of $\overline{w}_j$.  By \eqref{blowup estimates concl4}, for each $\sigma \in (0,2/q)$  
\begin{equation}\label{classify in L eqn1}
	\int_{B_{1/2}(0)} \sum_{j=1}^J \frac{|\overline{w}_{j,a}(x,y) - \lambda(\zeta) + \Lambda(\zeta) \cdot (x,y-\zeta) 
		- D\varphi^{(0)}_{j,a}(x,y) \cdot \kappa(\zeta)|^2}{|(x,y) - (0,\zeta)|^{n+2\alpha-\sigma}} \,dx \,dy \leq C
\end{equation}
for all $\zeta \in B^{n-2}_{1/4}(0)$, where $\varphi^{(0)}_{j,a}(x,y) = \frac{1}{q_j} \sum_{l=1}^{q_j} \varphi^{(0)}_{j,l}(x,y)$ for each $(x,y) \in B_{1/2}(0)$ and where $C = C(n,m,q,\alpha,\sigma) \in (0,\infty)$ is a constant.  
Since $w$ is component-wise Dirichlet energy minimizing on $B_{1/2}(0) \setminus \{0\} \times \mathbb{R}^{n-2}$ (as in Definition~\ref{v comp minimizing defn}), $\overline{w}_{j,a}$ is harmonic on $B_{1/2}(0) \setminus \{0\} \times \mathbb{R}^{n-2}$.  By \eqref{classify in L eqn1} and \eqref{blowup estimates concl5} (with $\sigma = 1/q$), we have that 
\begin{equation*}
	\rho^{-n} \int_{B_{\rho}(0,\zeta)} |\overline{w}_{j,a}|^2 \leq C 
\end{equation*}
for all $\zeta \in B^{n-2}_{1/4}(0)$ and $\rho \in (0,1/4)$, where $C = C(n,m,q,\alpha) \in (0,\infty)$.  Consequently, by the mean value property of harmonic functions, $\overline{w}_{j,a}$ is bounded on $B_{1/4}(0) \setminus \{0\} \times \mathbb{R}^{n-2}$.  It is then standard that $\overline{w}_{j,a}$ extends to a smooth harmonic function on $B_{1/4}(0)$.  
Since $w$ is homogeneous degree $\alpha$ on $B_{1/2}(0) \setminus \{0\} \times \mathbb{R}^{n-2}$, $\overline{w}_{j,a}$ extends to a homogeneous degree $\alpha$ single-valued harmonic function on $\mathbb{R}^n$.

If $\alpha$ is not an integer, then $\overline{w}_{j,a}$ and $\varphi^{(0)}_{j,a}$ are both identically zero on $\mathbb{R}^n$.  Hence using the $L^2(\mathbb{S}^1)$-orthogonality of $\lambda(\zeta)$ and $\Lambda(\zeta) \,(x,y)$, \eqref{classify in L eqn1} gives us gives us that $\lambda(\zeta) = 0$ and $\Lambda(\zeta) = 0$ for all $\zeta \in B^{n-2}_{1/4}(0)$.  If on the other hand $\alpha$ is an integer, then $\overline{w}_{j,a}$ is a homogeneous degree $\alpha$ harmonic polynomial on $\mathbb{R}^n$. Using the series expansion for the harmonic function $\overline{w}_{j,a}$, \eqref{classify in L eqn1} implies that $\lambda(\zeta) = \overline{w}_{j,a}(0,\zeta)$ and $\Lambda(\zeta) = D\overline{w}_{j,a}(0,\zeta)$ for each $\zeta \in B^{n-2}_{1/4}(0)$.  But by Lemma~\ref{blowup estimates lemma}, $\Lambda$ is a constant function on $B^{n-2}_{1/4}(0)$ such that, using the homogeneity of $\overline{w}_{j,a}$, $\Lambda(\zeta) = D\overline{w}_{j,a}(0,0) = 0$ for all $\zeta \in B^{n-2}_{1/4}(0)$.  By similar reasoning, $\lambda$ is a constant function on $B^{n-2}_{1/4}(0)$ such that $\lambda(\zeta) = \overline{w}_{j,a}(0,0) = 0$  for all $\zeta \in B^{n-2}_{1/4}(0)$.

\noindent\textit{Step~2:  Case $\kappa = 0$.}  Now we argue nearly identically to the proof of~\cite[Theorem~10.1]{KrumWic2} in the case that $\lambda = 0$.  (See in particular~\cite[Section~10.1]{KrumWic2}.  Also see~\cite[Theorem~5.8]{KrumWicb}). Fix $j \in \{1,2,\ldots,J\}$ and let $\overline{w} = \overline{w}_j$ (omitting the dependence on $j$).   The only change to the argument of~\cite[Theorem~10.1]{KrumWic2} is that since $\alpha > 1$, the inner variational formula (10.6) of~\cite{KrumWic2} can be proven using the same cut-off function argument that was used to prove the outer variational (10.5) of~\cite{KrumWic2} using the estimates (10.2) and (10.3) of~\cite{KrumWic2}.  Having proved both the outer and inner variational formulas (10.5) and (10.6) of~\cite{KrumWic2}, one argues exactly as in~\cite[Theorem~10.1]{KrumWic2} deducing from (10.2) of~\cite{KrumWic2} that $\overline{w}$ has frequency $= \alpha$ at each point of $\{0\} \times B^{n-2}_{1/16}(0)$.  Hence $\overline{w}(x,y+\zeta) = \overline{w}(x,y)$ for all $(x,y) \in \mathbb{R}^n$ and $\zeta \in \mathbb{R}^{n-2}$.  Since $\overline{w}$ is homogeneous degree $\alpha$ and $w_{j,l}$ is regular on each ball $B \subset\subset B_{1/2}(0) \setminus \{0\} \times \mathbb{R}^{n-2}$, it follows that $w_j$ is of the form of $\psi_j$ in \eqref{homogeneous class eqn1} and \eqref{homogeneous class eqn2} with $b = 0$, i.e.~$w \in \mathfrak{L}$.

\noindent\textit{Step~3:  Reduction to the case $\kappa = 0$.}  Fix $j \in \{1,2,\ldots,J\}$ such that $\varphi^{(0)}_j$ and $\varphi^{(\infty)}_j$ are both non-zero.  Define $w_{j,a} : \op{graph} \varphi^{(\infty)}_j |_{B_{1/2}(0) \setminus \{0\} \times \mathbb{R}^{n-2}} \rightarrow \mathbb{R}^m$ by  
\begin{equation*}
	w_{j,a}(x,y,\varphi_{j,l}(x,y)) = \frac{1}{m_j} \sum_{h=1}^{m_j} w_{j,l,h}(x,y) 
\end{equation*}
for all $(x,y) \in B_{1/2}(0) \setminus \{0\} \times \mathbb{R}^{n-2}$ and $l \in \{1,2,\ldots,q_j\}$.  Notice that we can regard $\varphi^{(\infty)}_j(re^{i\theta},y) = \op{Re}(c^{(\infty)}_j r^{\alpha} e^{i\alpha\theta})$ as a single-valued function which is $2\pi q_0$-periodic in $\theta$, and thus $w_{j,a}(re^{i\theta},y,\op{Re}(c^{(\infty)}_j r^{\alpha} e^{i\alpha\theta}))$ can also be regarded as a single-valued function which is $2\pi q_0$-periodic in $\theta$.  By adapting the proof of~\cite[Lemma~4.10]{Sim93}, we want to consider the Fourier series expansion of $w_{j,a}$ with respect to $\theta$ and show that for each $\iota \in \{1,2\}$ the Fourier coefficient of $D_{x_{\iota}} \varphi^{(0)}_j$ (after suitable normalization) is equal to $r^{\alpha-1} (\widehat{b}_{\iota} \cdot y)$ for some constant $\widehat{b}_{\iota} \in \mathbb{R}^{n-2}$.  It will then follow that $\kappa(\zeta) = (\kappa_1(\zeta),\kappa_2(\zeta))$ is given by $\kappa_{\iota}(\zeta) = b_{\iota} \cdot \zeta$ for all $\zeta \in B^{n-2}_{1/4}(0)$ and all $\iota \in \{1,2\}$, where each $b_{\iota} \in \mathbb{R}^{n-2}$ is a constant. 

Fix $\iota \in \{1,2\}$.  For any arbitrary $\theta_0 \in [0,2\pi)$, let $U_{\theta_0} = \{(re^{i\theta},y) \in B_{1/2}(0) : |\theta-\theta_0| < \pi/2 \}$.  We can express $\varphi^{(0)}_j(x,y) = \sum_{l=1}^{q_j} \llbracket \varphi^{(0)}_{j,l}(x,y) \rrbracket$ and $\varphi^{(\infty)}_j(x,y) = \sum_{l=1}^{q_j} \llbracket \varphi^{(\infty)}_{j,l}(x,y) \rrbracket$ for all $(x,y) \in U_{\theta}$, where $\varphi^{(0)}_{j,l}$ and $\varphi^{(\infty)}_{j,l}$ are single-valued harmonic functions on $U_{\theta_0}$ (like in \eqref{varphi0 localized}).  Let $w_j(x,y,\varphi^{(\infty)}_{j,l}(x,y)) = \sum_{h=1}^{m_j} \llbracket w_{j,l,h}(x,y) \rrbracket$ for each $(x,y) \in U_{\theta_0}$ (as in \eqref{w localized}).  Then for each $\theta_0 \in [0,2\pi)$, $\sum_{l=1}^{q_0} \sum_{h=1}^{m_j} w_{j,l,h}(re^{i\theta},y) \,D_{x_{\iota}} \varphi^{(0)}_{j,l}(e^{i\theta},0)$ is a well-defined single-valued function on $U_{\theta_0}$, which gives us a well-defined single-valued function on $B_{1/2}(0)$ independent of the choice of $\theta_0$.  Let $\widehat{w} : (0,\infty) \times \mathbb{R}^{n-2} \rightarrow \mathbb{R}$ be the homogeneous degree $\alpha$ functions such that   
\begin{equation*}
	\widehat{w}(r,y) = \frac{1}{2\pi q_0 m_j} \int_0^{2\pi} \sum_{l=1}^{q_0} \sum_{h=1}^{m_j} w_{j,l,h}(re^{i\theta},y) 
		\cdot D_{x_{\iota}} \varphi^{(0)}_{j,l}(e^{i\theta},0) \,d\theta 
\end{equation*}
for all $r > 0$ and $y \in \mathbb{R}^{n-2}$ with $r^2 + |y|^2 < 1/4$ (omitting the dependence of $\widehat{w}$ on $j$ and $\iota$).  Since $w$ is component-wise Dirichlet energy minimizing on $B_{1/2}(0) \setminus \{0\} \times \mathbb{R}^{n-2}$, on each ball $B \subset\subset B_{1/2}(0) \setminus \{0\} \times \mathbb{R}^{n-2}$ we can choose $w_{j,l}$ to be locally Dirichlet energy minimizing on $B$ (as in \eqref{w localized}) and thus $\frac{1}{m_j} \sum_{h=1}^{m_j} w_{j,l,h}$ is harmonic on $B$ for each $l \in \{1,\ldots,q_0\}$.  It follows that $\widehat{w}$ is a smooth single-valued function on $(0,\infty) \times \mathbb{R}^{n-2}$ and 
\begin{equation}\label{classify in L eqn2}
	\frac{1}{r} \frac{\partial}{\partial r} \bigg( r \frac{\partial\widehat{w}}{\partial r} \bigg) + \Delta_y \widehat{w} - \frac{(\alpha-1)^2}{r^2} \widehat{w} = 0 
\end{equation}
on $(0,\infty) \times \mathbb{R}^{n-2}$, where $\Delta_y$ denotes the Laplacian with respect to $y \in \mathbb{R}^{n-2}$.  By \eqref{blowup estimates concl4} (recalling that $\lambda = 0$ and $\Lambda = 0$)  
\begin{equation}\label{classify in L eqn3}
	\int_{B^{n-1}_{1/2}(0) \cap \{r > 0\}} \frac{|\widehat{w}(r,y) - r^{\alpha-1} \widehat{\kappa}(\zeta)|^2}{|(r,y-\zeta)|^{n+2\alpha-\sigma}} \,r \,dr \,dy \leq C 
\end{equation}
for all $\zeta \in B^{n-2}_{1/4}(0)$ and $\sigma \in (0,2/q)$, where $\kappa(\zeta) = (\kappa_1(\zeta),\kappa_2(\zeta))$ and $\widehat{\kappa}(\zeta) = \frac{1}{2\pi q_0} \|D\varphi^{(0)}_j\|_{L^2(\mathbb{S}^1 \times \{0\})}^2 \,\kappa_{\iota}(\zeta)$ for each $\zeta \in \mathbb{R}^{n-2}$ and where $C = C(n,m,q,\alpha,\sigma) \in (0,\infty)$ is a constant.  By \eqref{blowup on domain rmk eqn1} and \eqref{blowup estimates concl5}, $|\widehat{\kappa}(y) - \widehat{\kappa}(\zeta)| \leq C |y-\zeta|^{1-\sigma/2}$ for each $\sigma \in (0,2/q)$ and thus by \eqref{classify in L eqn3} 
\begin{equation}\label{classify in L eqn4}
	\int_{B^{n-1}_{\rho}(0,\zeta) \cap \{r > 0\}} |\widehat{w}(r,y) - r^{\alpha-1} \widehat{\kappa}(y)|^2 \,r \,dr \,dy \leq C \rho^{n+2\alpha-\sigma} 
\end{equation}
for all $\zeta \in B^{n-2}_{1/4}(0)$, $\rho \in (0,1/4)$, and $\sigma \in (0,2/q)$, where $C = C(n,m,q,\alpha,\sigma) \in (0,\infty)$ are constants.  For each $\rho \in (0,1/4)$, let $\{B^{n-1}_{\rho}(0,\zeta_k)\}_{1 \leq k \leq N}$ be a collection of balls covering $(0,\rho/2) \times B^{n-2}_{1/4}(0)$ such that $\zeta_k \in B^{n-2}_{1/4}(0)$ and $N \leq C(n) \,\rho^{2-n}$.  Setting $\zeta = \zeta_k$ in \eqref{classify in L eqn4} and summing over $k \in \{1,2,\ldots,N\}$ 
\begin{equation}\label{classify in L eqn5}
	\int_{B^{n-2}_{1/4}(0)} \int_0^{\rho/2} |\widehat{w}(r,y) - r^{\alpha-1} \widehat{\kappa}(y)|^2 \,r \,dr \,dy \leq C \rho^{2+2\alpha-\sigma} 
\end{equation}
for all $\rho \in (0,1/4)$ and $\sigma \in (0,2/q)$, where $C = C(n,m,q,\alpha,\sigma) \in (0,\infty)$ is a constant.  Dividing both sides of \eqref{classify in L eqn5} by $\rho^{2+2\alpha-\sigma/2}$ and integrating over $\rho \in (0,1/4)$ 
\begin{equation}\label{classify in L eqn6}
	\int_{B^{n-2}_{1/4}(0)} \int_0^{1/8} \frac{|\widehat{w}(r,y) - r^{\alpha-1} \widehat{\kappa}(y)|^2}{r^{1+2\alpha-\sigma/2}} \,dr \,dy \leq C 
\end{equation}
for all $\sigma \in (0,2/q)$, where $C = C(n,m,q,\alpha,\sigma) \in (0,\infty)$ is a constant. 

Define $\psi \in C^{\infty}(\mathbb{R}^{n-2},\mathbb{R}^m)$ by $\psi(z) = \widehat{w}(1,z)$ for each $z \in \mathbb{R}^{n-2}$ so that $\widehat{w}(r,y) = r^{\alpha} \psi(y/r)$ for all $r > 0$ and $y \in \mathbb{R}^{n-2}$.  Then \eqref{classify in L eqn2} becomes 
\begin{equation}\label{classify in L eqn7}
	\Delta_z \psi + z_i z_j D_{z_i z_j} \psi - (2\alpha-1) \,z_i D_{z_i} \psi + (2\alpha-1) \,\psi = 0 
\end{equation}
on $\mathbb{R}^{n-2}$.  Let $y = R\omega$ denote polar coordinates on $\mathbb{R}^{n-2}$, where $R > 0$ and $\omega \in \mathbb{S}^{n-3}$.  By \eqref{classify in L eqn6} (with $2\sigma$ in place of $\sigma$), there exists $\kappa \in L^{\infty}(\mathbb{R}^{n-2})$ such that 
\begin{equation}\label{classify in L eqn8}
	\int_{\mathbb{S}^{n-3}} \int_1^{\infty} R^{1-\sigma} \left|\frac{\psi(R\omega)}{R} - \kappa(R\omega)\right|^2 dR \,d\omega \leq C
\end{equation}
for all $\sigma \in (0,1/q)$, where $C = C(n,m,q,\alpha,\sigma) \in (0,\infty)$ is a constant.  For $p \in \{1,2,\ldots,n-2\}$, by \eqref{classify in L eqn7} the function $v = D_p \psi$ satisfies 
\begin{equation}\label{classify in L eqn9}
	\Delta_z v + z_i z_j D_{z_i z_j} v - (2\alpha-3) \,z_i \,D_{z_i} v = 0 
\end{equation}
on $\mathbb{R}^{n-2}$.  Let $\chi_1,\chi_2,\chi_3,\ldots$ be an $L^2(\mathbb{S}^{n-3})$ basis of eigenfunctions for the Laplacian on $\mathbb{S}^{n-3}$ so that $\Delta_{\mathbb{S}^{n-3}} \chi_k + \mu_k (\mu_k+n-4) \,\chi_k = 0$ on $\mathbb{S}^{n-3}$ for some integer $\mu_k \geq 0$.  Write $\psi$ in terms of its Fourier expansion 
\begin{equation*}
	\psi(R\omega) = \sum_{k=0}^{\infty} \gamma_k(R) \,\chi_k(\omega) 
\end{equation*}
for $R > 0$ and $\omega \in \mathbb{S}^{n-3}$.  Replace $\psi(R\omega)$ with $\gamma_k(R) \,\chi_k(\omega)$ and note that \eqref{classify in L eqn7}, \eqref{classify in L eqn8}, and \eqref{classify in L eqn9} continue to hold true for the new functions $\psi$ and $v =D_p \psi$.  Since $\psi$ is bounded near the origin, $\psi$ extends across the origin to a smooth solution to \eqref{classify in L eqn7}.  By \eqref{classify in L eqn7}  
\begin{equation}\label{classify in L eqn10}
	(1+R^2) \,\frac{d^2 \gamma_k}{dR^2} + \left( \frac{n-3}{R} - (2\alpha-1) \,R \right) \frac{d\gamma_k}{dR} 
		+ \left( 2\alpha-1 - \frac{\mu_k (\mu_k + n - 4)}{R^2} \right) \gamma_k = 0
\end{equation}
at each $r \in (0,\infty)$.  
By elementary ODE estimates  
\begin{equation}\label{classify in L eqn11}
	\sup_{\rho/2 < |z| < \rho} |D\psi(z)| \leq \frac{C}{\rho} \sup_{\rho/4 < |z| < 2\rho} |\psi(z)| 
\end{equation}
for some constant $C = C(n,\alpha,k) \in (0,\infty)$.  

Now rewrite \eqref{classify in L eqn9} as 
\begin{equation}\label{classify in L eqn12}
	\frac{\partial}{\partial R} \bigg(g \,\frac{\partial v}{\partial R}\bigg) + \frac{g}{R^2 (1+R^2)} \,\Delta_{\mathbb{S}^{n-3}} v = 0
\end{equation}
on $\mathbb{R}^{n-2}$, where $g = R^{n-3} (1+R^2)^{-(2\alpha+n-6)/2}$.  Since $\frac{\partial v}{\partial R}$ and $\frac{1}{R} \,\nabla_{\mathbb{S}^{n-3}} v$ are bounded as $R \downarrow 0$, we can write \eqref{classify in L eqn12} in the weak form  
\begin{equation}\label{classify in L eqn13}
	\int_{\mathbb{S}^{n-3}} \int_0^{\infty} \left( g \,\frac{\partial v}{\partial R} \,\frac{\partial\phi}{\partial R} 
		+ \frac{g}{R^2 (1+R^2)} \,\nabla^{\mathbb{S}^{n-3}} v \cdot \nabla^{\mathbb{S}^{n-3}} \phi \right) dR \,d\omega = 0
\end{equation}
for all $\phi \in C^{\infty}_c(\mathbb{R}^{n-2})$.  Replace $\phi$ with $v \phi^2$ in \eqref{classify in L eqn13} where $\phi = \phi(R)$ is a smooth radially symmetric cut-off function which is constant near $R = 0$ and vanishes for large $R$ and apply Cauchy's inequality to get  
\begin{align}\label{classify in L eqn14}
	&\int_{\mathbb{S}^{n-3}} \int_0^{\infty} \left( g \left|\frac{\partial v}{\partial R}\right|^2 + \frac{g}{R^2 (1+R^2)} \,|\nabla^{\mathbb{S}^{n-3}} v|^2 \right) \phi(R)^2 
	\,dR \,d\omega \\ \leq\,& 4 \int_{\mathbb{S}^{n-3}} \int_0^{\infty} g \,v^2 \,|\phi'(R)|^2 \,dR \,d\omega . \nonumber
\end{align}
Given $\rho \in (1,\infty)$, let $\phi$ be the logarithmic cut-off function 
\begin{equation*}
	\phi(R) = \begin{cases} 
		1 &\text{if } R \leq \rho \\
		\log(R/\rho^2)/\log(1/\rho) &\text{if } \rho^2 < R < \rho \\
		0 &\text{if } R \geq \rho^2 
	\end{cases} 
\end{equation*}
so that \eqref{classify in L eqn14} together with \eqref{classify in L eqn11} and \eqref{classify in L eqn8} gives us 
\begin{align*}
	&\int_{\mathbb{S}^{n-3}} \int_0^{\rho} \left( g \left|\frac{\partial v}{\partial R}\right|^2 + \frac{g}{R^2 (1+R^2)} \,|\nabla^{\mathbb{S}^{n-3}} v|^2 \right) \,dR \,d\omega 
	\\ \leq\,& \frac{4}{\log^2(\rho)} \int_{\mathbb{S}^{n-3}} \int_{\rho}^{\rho^2} R^{1-2\alpha} \,v^2 \,dR \,d\omega 
	\leq \frac{C}{\log^2(\rho)} \int_{\mathbb{S}^{n-3}} \int_{\rho/2}^{2\rho^2} R^{-1-2\alpha} \,\psi^2 \,dR \,d\omega 
	\\ \leq\,& \frac{C}{\log^2(\rho)} \int_{\mathbb{S}^{n-3}} \int_{\rho/2}^{2\rho^2} R^{1-2\alpha} \left| \frac{\psi(R\omega)}{R} - \kappa(R\omega) \right|^2 
		\,dR \,d\omega + \frac{C\rho^{2-2\alpha} \|\kappa\|_{L^{\infty}(\mathbb{R}^{n-2})}^2}{\log^2(\rho)} 
	\\ \leq\,& \frac{C\rho^{-2\alpha+\sigma}}{\log^2(\rho)} + \frac{C\rho^{2-2\alpha} \|\kappa\|_{L^{\infty}(\mathbb{R}^{n-2})}^2}{\log^2(\rho)} \rightarrow 0 
		\text{ as } \rho \rightarrow\infty ,
\end{align*}
where $\sigma \in (0,1/q)$ and $C = C(n,m,q,\alpha,k,\sigma) \in (0,\infty)$ are constants.  Hence $v$ is a constant function on $\mathbb{R}^{n-2}$; that is, $\psi(z) = \widehat{a} + \widehat{b} \cdot z$ for some $\widehat{a} \in \mathbb{R}$ and $\widehat{b} \in \mathbb{R}^{n-2}$.  Substituting $\psi(z) = \widehat{a} + \widehat{b} \cdot z$ into \eqref{classify in L eqn7} gives us that $\widehat{a} = 0$.  Therefore, $\psi(z) = \widehat{b} \cdot z$ for some $\widehat{b} \in \mathbb{R}^{n-2}$.  In other words, $\widehat{w}(r,y) = r^{\alpha-1} \,(\widehat{b} \cdot y)$ for some $\widehat{b} \in \mathbb{R}^{n-2}$.  

By \eqref{classify in L eqn6}
\begin{equation*}
	\int_{B^{n-2}_{1/2}(0)} \int_0^{1/8} \frac{|\widehat{b} \cdot y - \widehat{\kappa}(y)|^2}{r^{3-\sigma/2}} \,dr \,dy \leq C 
\end{equation*}
for all $\sigma \in (0,2/q)$, where $C = C(n,m,q,\alpha,\sigma) \in (0,\infty)$ is a constant, and thus $\widehat{\kappa}(y) = \widehat{b} \cdot y$ for all $y \in B^{n-2}_{1/4}(0)$.  In other words, there exists $b_1,b_2 \in \mathbb{R}^{n-2}$ such that $\kappa_{\iota}(y) = b_{\iota} \cdot y$ for all $y \in B^{n-2}_{1/4}(0)$ and $\iota \in \{1,2\}$.  Hence we can rewrite \eqref{blowup estimates concl4} (with $\lambda = 0$ and $\Lambda = 0$) as 
\begin{equation}\label{classify in L eqn15}
	\int_{B_{1/2}(0)} \sum_{j=1}^J \sum_{l=1}^{q_j} \sum_{h=1}^{m_j} \frac{|w_{j,l,h}(x,y) - D_{x_1} \varphi^{(0)}_{j,l}(x,y) \,(b_1 \cdot \zeta) 
		- D_{x_2} \varphi^{(0)}_{j,l}(x,y) \,(b_2 \cdot \zeta)|^2}{|(x,y) - (0,\zeta)|^{n+2\alpha-\sigma}} \,dx \,dy \leq C 
\end{equation}
for all $\zeta \in B^{n-2}_{1/4}(0)$ and $\sigma \in (0,2/q)$, where $C = C(n,m,q,\alpha,\sigma) \in (0,\infty)$ is a constant.  Note that by \eqref{blowup estimates concl5}, $|b_1| + |b_2| \leq C(n,m,q,\alpha)$.  

Now define $\Gamma_{\nu} = e^{E_{\nu} B/\widehat{E}_{\nu}}$, where $B$ is as in \eqref{integrability eqn1} with $b = b_1 + ib_2$.  By the argument of~\cite[Theorem~10.1]{KrumWic2} (see the end of~\cite[Section~10.2]{KrumWic2}) and Lemma~\ref{integrability lemma}, we see that the fine blow-up $\widetilde{w} = (\widetilde{w}_j) \in \mathfrak{B}$ of the rotated sequence $\widetilde{T}^{(\nu)} = \Gamma_{\nu\#} T^{(\nu)}$ relative to the sequence $(\varphi^{(\nu)})$ in $\Phi_{\alpha,q,p}$ (by excess $\widetilde{E}_{\nu} = E(\Gamma_{\nu\#} T^{(\nu)},\varphi^{(\nu)},\mathbf{C}_1(0))$) satisfies  
\begin{equation}\label{classify in L eqn16}
	\widetilde{c} \,\widetilde{w}_j(x,y,\varphi^{(\infty)}_{j,l}(x,y)) = w_j(x,y,\varphi^{(\infty)}_{j,l}(x,y)) - D_{x_1} \varphi^{(0)}_{j,l}(x,y) \,(b_1 \cdot y) 
		- D_{x_2} \varphi^{(0)}_{j,l}(x,y) \,(b_2 \cdot y) 
\end{equation}
for all $(x,y) \in B_{1/2}(0)$, $j \in \{1,2,\ldots,J\}$, and $l \in \{1,2,\ldots,q_j\}$, where $\widetilde{c} = \lim_{\nu\rightarrow\infty} \widetilde{E}_{\nu}/E_{\nu}$.  Note that $\widetilde{E}_{\nu} \leq C E_{\nu}$, where $C = C(n,m,q,\alpha) \in (0,\infty)$ is a constant, and thus the limit $\widetilde{c}$ exists and $0 \leq \widetilde{c} \leq C$.  If $\widetilde{c} = 0$, then by \eqref{classify in L eqn16} 
\begin{equation*}
	w_j(x,y,\varphi^{(\infty)}_{j,l}(x,y)) = D_{x_1} \varphi^{(0)}_{j,l}(x,y) \,(b_1 \cdot y) + D_{x_2} \varphi^{(0)}_{j,l}(x,y) \,(b_2 \cdot y) 
\end{equation*}
for all $(x,y) \in B_{1/2}(0)$ and $j \in \{1,2,\ldots,J\}$ so that $w \in \widetilde{\mathfrak L}$ trivially.  On the other hand, if $\widetilde{c} > 0$, then by \eqref{classify in L eqn15} and \eqref{classify in L eqn16}, $\widetilde{w}$ is a homogeneous degree $\alpha$ element of $\mathfrak{B}$ satisfying 
\begin{equation}\label{classify in L eqn17}
	\int_{B_{1/2}(0)} \sum_{j=1}^J \sum_{l=1}^{q_j} \frac{|\widetilde{w}_{j,l}(x,y)|^2}{|(x,y) - (0,\zeta)|^{n+2\alpha-\sigma}} \,dx \,dy \leq C 
\end{equation}
for all $\zeta \in B^{n-2}_{1/4}(0)$ and $\sigma \in (0,2/q)$, where $\widetilde{w}_{j,l}(x,y) = \widetilde{w}_j(x,y,\varphi^{(\infty)}_{j,l}(x,y))$ for each $(x,y) \in B_{1/2}(0)$ (as in \eqref{w localized}).  Applying the special case $\kappa = 0$ from Step~2, we deduce that $\widetilde{w} \in \mathfrak{L}$.  By \eqref{classify in L eqn17}, this means that $w \in \widetilde{\mathfrak L}$.
\end{proof}

\section{Asymptotic decay of fine blow-ups}\label{sec:fine blowups decay}

The main result of this section is 
Theorem~\ref{blowup decay thm}, which asserts that each fine blow-up $w \in \mathfrak{B}$ is decaying to a unique homogeneous degree $\alpha$ blow-up $\psi \in \widetilde{\mathfrak L}$ on approach to the origin.  Our approach utilizes the well-known hole-filling technique in a manner similar to its use in~\cite{Sim93}, and is based on the following lemma.

\begin{lemma}\label{reverseHS lemma}
Let $(m_j,\varphi^{(0)}_j,\varphi^{(\infty)}_j) \in \mathfrak{D}$, $w = (w_j) \in \mathfrak{B}(m_j,\varphi^{(0)}_j,\varphi^{(\infty)}_j)$, and $\psi = (\psi_j) \in \widetilde{\mathfrak L}(m_j,\varphi^{(0)}_j,\varphi^{(\infty)}_j)$ such that 
\begin{equation}\label{reverseHS hyp1}
	\int_{B_{1/2}(0)} \sum_{j=1}^J \sum_{l=1}^{q_j} \mathcal{G}(w_{j,l},\psi_{j,l})^2 
	\leq 2 \inf_{\psi' \in \widetilde{\mathfrak L}(m_j,\varphi^{(0)}_j,\varphi^{(\infty)}_j)} 
		\int_{B_{1/2}(0)} \sum_{j=1}^J \sum_{l=1}^{q_j} \mathcal{G}(w_{j,l},\psi'_{j,l})^2 ,
\end{equation}
where $w_{j,l}$ is as in \eqref{w localized}, $\psi_{j,l}$ is as in \eqref{psi localized}, and $\psi'_{j,l}(x,y) = \psi'_j(x,y,\varphi^{(\infty)}_{j,l}(x,y))$ for each $(x,y) \in B_{1/2}(0)$ where $\varphi^{(\infty)}_{j,l}$ is as in \eqref{varphi0 localized} (cf.~\eqref{psi localized}).  Then 
\begin{equation}\label{reverseHS concl}
	\int_{B_{1/2}(0)} \sum_{j=1}^J \sum_{l=1}^{q_j} \mathcal{G}(w_{j,l},\psi_{j,l})^2 
	\leq C \int_{B_{1/2}(0) \setminus B_{1/8}(0)} \sum_{j=1}^J \sum_{l=1}^{q_j} \left|\frac{\partial (w_{j,l}/R^{\alpha})}{\partial R}\right|^2 ,
\end{equation}
where $R(x,y) = |(x,y)|$ for each $(x,y) \in \mathbb{R}^n$ and $C = C(n,m,q,\alpha) \in (0,\infty)$ is a constant.
\end{lemma}

As in~\cite{KrumWic2} and~\cite{KrumWicb} (and in contrast with the multiplicity one setting of~\cite{Sim93}), the present setting is higher multiplicity and involves functions which are multi-valued.  Thus we proceed by first proving Lemma~\ref{reverseHS lemma0} below, in which we prove the conclusion of Lemma~\ref{reverseHS lemma} assuming a weaker variant \eqref{reverseHS hyp2} of \eqref{reverseHS hyp1} while also assuming that $\psi \in \mathfrak{L}$ and an analog of Hypothesis~$(\star\star)$ for fine blow-ups.

Let $(m_j,\varphi^{(0)}_j,\varphi^{(\infty)}_j) \in \mathfrak{D}$ and $\psi \in \mathfrak{L}(m_j,\varphi^{(0)}_j,\varphi^{(\infty)}_j)$.  Given $j \in \{1,2,\ldots,J\}$ such that $\varphi^{(\infty)}_j$ is non-zero, we know that $\psi_j$ is given by \eqref{homogeneous class eqn1} for some $a_{j,h} \in \mathbb{C}^m$ and with $b = 0$.  Let $d_j(\psi)$ denote the number of distinct values of $a_{j,h}$, i.e.~the number of distinct components $\op{Re}(a_{j,h} r^{\alpha} e^{i\theta\theta})$ as in \eqref{homogeneous class eqn1}.  Similarly, given $j \in \{1,2,\ldots,J\}$ such that $\varphi^{(0)}_j$ and $\varphi^{(\infty)}_j$ are both identically zero, $\psi_j$ is given by \eqref{homogeneous class eqn2} for some positive integer $N$ and some $a_{j,h} \in \mathbb{C}^m \setminus \{0\}$.  Let $d_j(\psi)$ denote the number of distinct values of $a_{j,h}$, i.e.~the number of distinct non-zero components $\op{Re}(a_{j,h} r^{\alpha} e^{i\theta\theta})$ as in \eqref{homogeneous class eqn2}.  Let 
\begin{equation*}
	d(\psi) = \sum_{j=1}^J d_j(\psi). 
\end{equation*}
Note that $1 \leq d_j(\psi) \leq m_j$ if $\varphi^{(\infty)}_j$ is non-zero and $0 \leq d_j(\psi) \leq \lfloor m_j/q_0 \rfloor$ if $\varphi^{(0)}_j$ and $\varphi^{(\infty)}_j$ are both identically zero.  It follows that $p \leq d(\psi) \leq \sum_{j=1}^J m_j q_j/q_0 = q/q_0$; that is, $p \leq d(\psi) \leq \lfloor q/q_0 \rfloor$.  (Here $p$ is the number of non-zero components of $\varphi^{(\nu)}$ as in Definition~\ref{blowup class defn}, and thus $p$ is equal to the number of $j \in \{1,2,\ldots,J\}$ for which $\varphi^{(\infty)}_j$ is non-zero.)  For $s \in \{p,p+1,\ldots,\lfloor q/q_0 \rfloor\}$, let 
\begin{equation*}
	\mathfrak{L}_s(m_j,\varphi^{(0)}_j,\varphi^{(\infty)}_j) = \{ \psi \in \mathfrak{L}(m_j,\varphi^{(0)}_j,\varphi^{(\infty)}_j) : d(\psi) = s \} 
\end{equation*}
and note that $\mathfrak{L}(m_j,\varphi^{(0)}_j,\varphi^{(\infty)}_j) = \bigcup_{s=p}^{\lfloor q/q_0 \rfloor} \mathfrak{L}_s(m_j,\varphi^{(0)}_j,\varphi^{(\infty)}_j)$.

\begin{lemma}\label{reverseHS lemma0}
For every $\overline{M} \in [1,\infty)$ there exists $\overline{\beta} = \overline{\beta}(n,m,q,\alpha,\overline{M}) \in (0,1)$ and $\overline{C} = \overline{C}(n,m,q,\alpha,\overline{M}) \in (0,\infty)$ such that the following holds true.  Let $(m_j,\varphi^{(0)}_j,\varphi^{(\infty)}_j) \in \mathfrak{D}$ and let $p \geq 1$ equal the number of $j \in \{1,2,\ldots,J\}$ for which $\varphi^{(0)}$ and $\varphi^{(\infty)}_j$ are both non-zero.  Let $s \in \{p,p+1,\ldots,\lfloor q/q_0 \rfloor\}$, $w = (w_j) \in \mathfrak{B}(m_j,\varphi^{(0)}_j,\varphi^{(\infty)}_j)$, and $\psi = (\psi_j) \in \mathfrak{L}_s(m_j,\varphi^{(0)}_j,\varphi^{(\infty)}_j)$ such that 
\begin{equation}\label{reverseHS hyp2}
	\int_{B_{1/2}(0)} \sum_{j=1}^J \sum_{l=1}^{q_j} \mathcal{G}(w_{j,l},\psi_{j,l})^2 
	\leq 2 \overline{M}^2 \inf_{\psi' \in \widetilde{\mathfrak L}(m_j,\varphi^{(0)}_j,\varphi^{(\infty)}_j)} 
		\int_{B_{1/2}(0)} \sum_{j=1}^J \sum_{l=1}^{q_j} \mathcal{G}(w_{j,l},\psi'_{j,l})^2 ,
\end{equation}
and either {\rm (i)} $s = p$ or {\rm (ii)} $s > p$ and 
\begin{equation}\label{reverseHS hyp3}
	\int_{B_{1/2}(0)} \sum_{j=1}^J \sum_{l=1}^{q_j} \mathcal{G}(w_{j,l},\psi_{j,l})^2 
	\leq \overline{\beta}^2 \inf_{\psi' \in \bigcup_{s'=p}^{s-1} \mathfrak{L}_{s'}(m_j,\varphi^{(0)}_j,\varphi^{(\infty)}_j)} 
		\int_{B_{1/2}(0)} \sum_{j=1}^J \sum_{l=1}^{q_j} \mathcal{G}(w_{j,l},\psi'_{j,l})^2 
\end{equation}
(where $w_{j,l}$, $\psi_{j,l}$, and $\psi'_{j,l}$ are as in Lemma~\ref{reverseHS lemma}).  Then \eqref{reverseHS concl} holds true with $C = \overline{C}$.
\end{lemma}

\begin{proof}
Fix $\overline{M} \in [1,\infty)$ and integers $1 \leq p \leq s \leq \lfloor q/q_0 \rfloor$.  We proceed by contradiction, supposing that for each $\nu = 1,2,3,\ldots$ there exists $\overline{\beta}_{\nu} \in (0,1)$ with $\overline{\beta}_{\nu} \downarrow 0$, $(m^{(\nu)}_j,\varphi^{(\nu,0)}_j,\varphi^{(\nu,\infty)}_j)_{1 \leq j \leq J^{(\nu)}} \in \mathfrak{D}$, $w^{(\nu)} = (w^{(\nu)}_j) \in \mathfrak{B}(m^{(\nu)}_j,\varphi^{(\nu,0)}_j,\varphi^{(\nu,\infty)}_j)$, and $\psi^{(\nu)} = (\psi^{(\nu)}_j) \in \mathfrak{L}_s(m^{(\nu)}_j,\varphi^{(\nu,0)}_j,\varphi^{(\nu,\infty)}_j)$ such that 
\begin{equation}\label{reverseHS step0 eqn1}
	\int_{B_{1/2}(0)} \sum_{j=1}^{J^{(\nu)}} \sum_{l=1}^{q^{(\nu)}_j} \mathcal{G}(w^{(\nu)}_{j,l},\psi^{(\nu)}_{j,l})^2 
	\leq 2 \overline{M}^2 \inf_{\psi' \in \widetilde{\mathfrak L}(m^{(\nu)}_j,\varphi^{(\nu,0)}_j,\varphi^{(\nu,\infty)}_j)} 
		\int_{B_{1/2}(0)} \sum_{j=1}^J \sum_{l=1}^{q_j} \mathcal{G}(w^{(\nu)}_{j,l},\psi'_{j,l})^2 ,
\end{equation}
and either (i) $s = p$ or (ii) $s > p$ and 
\begin{equation}\label{reverseHS step0 eqn2}
	\int_{B_{1/2}(0)} \sum_{j=1}^{J^{(\nu)}} \sum_{l=1}^{q^{(\nu)}_j} \mathcal{G}(w^{(\nu)}_{j,l},\psi^{(\nu)}_{j,l})^2 
	\leq \overline{\beta}_{\nu}^2 \inf_{\psi' \in \bigcup_{s'=p}^{s-1} \mathfrak{L}_{s'}(m^{(\nu)}_j,\varphi^{(\nu,0)}_j,\varphi^{(\nu,\infty)}_j)} 
		\int_{B_{1/2}(0)} \sum_{j=1}^J \sum_{l=1}^{q_j} \mathcal{G}(w^{(\nu)}_{j,l},\psi'_{j,l})^2 
\end{equation}
but 
\begin{equation}\label{reverseHS step0 eqn3} 
	\int_{B_{1/2}(0) \setminus B_{1/8}(0)} \sum_{j=1}^{J^{(\nu)}} \sum_{l=1}^{q^{(\nu)}_j} \left|\frac{\partial (w^{(\nu)}_{j,l}/R^{\alpha})}{\partial R}\right|^2 
	< \frac{1}{\nu} \int_{B_{1/2}(0)} \sum_{j=1}^{J^{(\nu)}} \sum_{l=1}^{q^{(\nu)}_j} \mathcal{G}(w^{(\nu)}_{j,l},\psi^{(\nu)}_{j,l})^2 ,
\end{equation}
where $w^{(\nu)}_{j,l}$ is as in \eqref{w localized} with $\varphi^{(\nu,\infty)}_j$ and $w^{(\nu)}_j$ in place of $\varphi^{(\infty)}_j$ and $w_j$, $\psi^{(\nu)}_{j,l}$ is as in \eqref{psi localized} with $\varphi^{(\nu,\infty)}_j$ and $\psi^{(\nu)}_j$ in place of $\varphi^{(\infty)}_j$ and $\psi$, and $\psi'_{j,l}$ is as in \eqref{psi localized} with $\varphi^{(\nu,\infty)}_j$ and $\psi'$ in place of $\varphi^{(\infty)}_j$ and $\psi$.  Let 
\begin{equation*}%\label{reverseHS step0 eqn4} 
	F_{\nu} = \left( \int_{B_{1/2}(0)} \sum_{j=1}^J \sum_{l=1}^{q_j} \mathcal{G}(w^{(\nu)}_{j,l},\psi^{(\nu)}_{j,l})^2 \right)^{1/2} 
\end{equation*}
and note that $F_{\nu} > 0$ by \eqref{reverseHS step0 eqn1} and \eqref{reverseHS step0 eqn3}.  To obtain a contradiction, we shall proceed in 5 steps.

\noindent\textit{Step~1.  Selection of currents associated with $w^{(\nu)}$, construction of homogeneous harmonic $q$-valued functions associated with $\psi^{(\nu)}$, and some preliminary bounds.}  Since $w^{(\nu)} \in \mathfrak{B}$, there exists a sequence of locally area-minimizing rectifiable currents $(T^{(\nu,k)})_{k=1}^{\infty}$ in $\mathbf{C}_1(0)$, a sequence $(\varphi^{(\nu,k)})_{k=1}^{\infty}$ in $\Phi_{\alpha,q,p}$, and sequences of positive numbers $(\varepsilon_{\nu,k})_{k=1}^{\infty}$, $(\beta_{\nu,k})_{k=1}^{\infty}$, and $(\delta_{\nu,k})_{k=1}^{\infty}$ tending to zero such that conditions~(A)--(D) of Definition~\ref{blowup class defn} hold true with $T^{(\nu,k)},\varphi^{(\nu,k)},\varepsilon_{\nu,k},\beta_{\nu,k},\delta_{\nu,k}$ in place of $T^{(\nu)},\varphi^{(\nu)},\varepsilon_{\nu},\beta_{\nu},\delta_{\nu}$ and $w^{(\nu)}$ is the blow-up of $(T^{(\nu,k)})$ relative to $(\varphi^{(\nu,k)})$ by excess $E_{\nu,k} = E_{\varsigma/2,\varepsilon_{\nu,k}}(T^{(\nu,k)},\varphi^{(\nu,k)},\mathbf{C}_1(0))$ as $k\rightarrow\infty$.  We require that $\varphi^{(\nu,k)},\varphi^{(\nu,0)}_j,\varphi^{(\nu,\infty)}_j$ satisfy the requirements of Remark~\ref{blowup on domain rmk}.  In particular, 
\begin{equation*}
    \varphi^{(\nu,k)} = \sum_{j=1}^{J^{(\nu)}} m^{(\nu)}_j \varphi^{(\nu,k)}_j 
\end{equation*}
for distinct $\varphi^{(\nu,k)}_j : \mathbb{R}^n \rightarrow \mathcal{A}_{q^{(\nu)}_j}(\mathbb{R}^m)$ (as in Definition~\ref{Phi-p_defn}).  After passing to a subsequence, we shall assume that $J^{(\nu)} = J$, $m^{(\nu)}_j = m_j$, and $q^{(\nu)}_j = q_j$ independent of $\nu$ and for each $j \in \{1,2,\ldots,J\}$ either:
\begin{enumerate}[itemsep=2mm,topsep=0mm]
	\item[(i)]  $q_j = q_0$ and $\varphi^{(\nu,\infty)}_j$ is non-zero,   
	\begin{gather*}
		\varphi^{(\nu,0)}_j(x,y) = \op{Re}(c^{(\nu,0)}_j (x_1+ix_2)^{\alpha}), \quad 
		\varphi^{(\nu,\infty)}_j(x,y) = \op{Re}(c^{(\nu,\infty)}_j (x_1+ix_2)^{\alpha}), \\ 
		\varphi^{(\nu,k)}_j(x,y) = \op{Re}(c^{(\nu,k)}_j (x_1+ix_2)^{\alpha}) \nonumber 
	\end{gather*}
	for all $(x,y) \in \mathbb{R}^n$, where $c^{(\nu,k)}_j \in \mathbb{C}^m \setminus \{0\}$ and setting $\widehat{E}_{\nu,k} = \widehat{E}(T^{(\nu,k)},\varphi^{(\nu,k)},\mathbf{C}_1(0))$ 
	\begin{equation*}
		c^{(\nu,0)}_j = \lim_{k\rightarrow\infty} c^{(\nu,k)}_j/\widehat{E}_{\nu,k}, \quad 
		c^{(\nu,\infty)}_j = \lim_{k\rightarrow\infty} c^{(\nu,k)}_j/|c^{(\nu,k)}_j| ; 
	\end{equation*}
	
	\item[(ii)]  $q_j = 1$ and $\varphi^{(\nu,0)}_j$, $\varphi^{(\nu,\infty)}_j$, and $\varphi^{(\nu,k)}_j$ are all identically zero.
\end{enumerate} 
By the definition of fine blow-up (see Subsection~\ref{sec:blowup subsec}), for each $j \in \{1,2,\ldots,J\}$ and each integer $\nu \geq 1$ there is a sequence of positive numbers $(\tau_{\nu,k})_{k=1}^{\infty}$ tending to zero, Lipschitz functions $v^{(\nu,k)}_j : \op{graph}\varphi^{(\nu,k)}_j |_{B_{1/2}(0) \cap \{r > \tau_{\nu,k}\}} \rightarrow \mathcal{A}_{q_j}(\mathbb{R}^m)$, and closed set $K^{(\nu,k)}_j \subset B_{1/2}(0) \cap \{r > \tau_{\nu,k}\}$ such that conclusion~(C) of Theorem~\ref{graphical_thm} holds true with $T^{(\nu,k)}, \varphi^{(\nu,k)}, v^{(\nu,k)}_j, K^{(\nu,k)}_j$ in place of $T,\varphi,v_j,K_j$ and if $\varphi^{(\nu,k)}_j$ and $\varphi^{(\nu,\infty)}_j$ are non-zero then 
\begin{equation}\label{reverseHSshort step1 eqn2}
    v^{(\nu,k)}_j(re^{i\theta},y,\op{Re}(c^{(\nu,k)}_j r^{\alpha} e^{i\alpha\theta}))/E_{\nu,k} \rightarrow w^{(\nu)}_j(re^{i\theta},y,\op{Re}(c^{(\nu,\infty)}_j r^{\alpha} e^{i\alpha\theta}))
\end{equation}
in the $L^2$-topology on compact subsets of $B_{1/2}(0) \setminus \{0\}\times\mathbb{R}^{n-2}$, whereas if $\varphi^{(\nu,k)}_j$, $\varphi^{(\nu,0)}_j$, and $\varphi^{(\nu,\infty)}_j$ are identically zero then 
\begin{equation}\label{reverseHSshort step1 eqn3}
    v^{(\nu,k)}_j(x,y,0)/E_{\nu,k} \rightarrow w^{(\nu)}_j(x,y,0)
\end{equation}
in the $L^2$-topology on compact subsets of $B_{1/2}(0) \setminus \{0\}\times\mathbb{R}^{n-2}$. 

Construct $\widetilde{\varphi}^{(\nu,k)} \in \Phi_{\alpha,q}$ as in Lemma~\ref{integrability lemma} with $\varphi^{(\nu,k)},\varphi^{(\nu,\infty)}_j,\psi^{(\nu)}_j,E_{\nu,k}$ in place of $\varphi^{(\nu)},\varphi^{(\infty)}_j,\psi_j,E_{\nu}$. In particular, for each $j \in \{1,2,\ldots,J\}$, let 
\begin{equation*}
	\psi^{(\nu)}_j = \sum_{h=1}^{s_j} \widetilde{m}_{j,h} \widetilde{\psi}^{(\nu)}_{j,h} \quad\text{and}\quad \widetilde{\varphi}^{(\nu,k)}_j = \sum_{j=1}^J \sum_{h=1}^{s_j} \widetilde{m}_{j,h} \widetilde{\varphi}^{(\nu,k)}_{j,h}, 
\end{equation*}
where $s_j$, $\widetilde{m}_{j,h}$, and $\widetilde{q}_{j,h}$ be positive integers such that $\sum_{h=1}^{s_j} \widetilde{m}_{j,h} \widetilde{q}_{j,h} = m_j q_j$, $\widetilde{\psi}^{(\nu)}_{j,h} : \op{graph}\varphi^{(\nu)}_j \rightarrow \mathcal{A}_{\widetilde{q}_{j,h}/q_j}(\mathbb{R}^m)$ are distinct components of $\psi^{(\nu)}_j$, and $\widetilde{\varphi}^{(\nu,k)}_{j,h} : \mathbb{R}^n \rightarrow \mathcal{A}_{\widetilde{q}_{j,h}}(\mathbb{R}^m)$ are the corresponding components of $\widetilde{\varphi}^{(\nu,k)}$ (which are described below).  After passing to a subsequence, we shall assume that $s_j$, $\widetilde{m}_{j,h}$, and $\widetilde{q}_{j,h}$ are independent of $\nu$ and for each $j \in \{1,2,\ldots,J\}$ and $h \in \{1,2,\ldots,s_j\}$ one of the following three possibilities holds true: 
\begin{enumerate}[itemsep=2mm,topsep=0mm]
	\item[(a)]  $\varphi^{(\nu,k)}_j$ and $\varphi^{(\nu,\infty)}_j$ are non-zero $q_0$-valued functions for all $\nu,k = 1,2,3,\ldots$;
	\item[(b)]  $\varphi^{(\nu,k)}_j$, $\varphi^{(\nu,0)}_j$, and $\varphi^{(\nu,\infty)}_j$ are identically zero and $\psi^{(\nu)}_{j,h} : \op{graph}\varphi^{(\nu,\infty)} \rightarrow \mathcal{A}_{q_0}(\mathbb{R}^m)$ are non-zero for all $\nu,k = 1,2,3,\ldots$;
	\item[(c)]  $\varphi^{(\nu,k)}_j$, $\varphi^{(\nu,0)}_j$, $\varphi^{(\nu,\infty)}_j$, and $\psi^{(\nu)}_{j,h}$ are identically zero for all $\nu,k = 1,2,3,\ldots$.
\end{enumerate}
In case~(a), $\widetilde{q}_{j,h} = q_0$, $\psi^{(\nu)}_{j,h} : \op{graph}\varphi^{(\nu,\infty)}_j \rightarrow \mathbb{R}^m$ is given by 
\begin{equation}\label{reverseHSshort step1 eqn1}
    \widetilde{\psi}^{(\nu)}_{j,h}(re^{i\theta},y, \op{Re}(c^{(\nu,\infty)}_j r^{\alpha} e^{i\alpha\theta})) = \op{Re}(a^{(\nu)}_{j,h} r^{\alpha} e^{i\alpha\theta})
\end{equation}
for all $r > 0$, $\theta \in [0,2\pi q_0)$, and $y \in \mathbb{R}^{n-2}$, where $a^{(\nu)}_{j,h} \in \mathbb{C}^m$, and $\widetilde{\varphi}^{(\nu,k)}_{j,h} : \mathbb{R}^n \rightarrow \mathcal{A}_{q_0}(\mathbb{R}^m)$ is given by 
\begin{equation}\label{reverseHSshort step1 eqn4}
    \widetilde{\varphi}^{(\nu,k)}(x,y) = \op{Re}((c^{(\nu,k)}_j + E_{\nu,k} a^{(\nu)}_{j,h}) (x_1+ix_2)^{\alpha}) 
\end{equation}
for all $x = (x_1,x_2) \in \mathbb{R}^2$ and $y \in \mathbb{R}^{n-2}$.  In case~(b), $\widetilde{q}_{j,h} = q_0$ and $\psi^{(\nu)}_{j,h} : \op{graph}\varphi^{(\nu,\infty)}_j \rightarrow \mathcal{A}_{q_0}(\mathbb{R}^m)$ and $\widetilde{\varphi}^{(\nu,k)}_{j,h} : \mathbb{R}^n \rightarrow \mathcal{A}_{q_0}(\mathbb{R}^m)$ are given by 
\begin{equation}\label{reverseHSshort step1 eqn5}
    \widetilde{\psi}^{(\nu)}_{j,h}(x,y,0) = \op{Re}(a^{(\nu)}_{j,h} (x_1+ix_2)^{\alpha}) , \quad 
    \widetilde{\varphi}^{(\nu,k)}_{j,h}(x,y,0) = \op{Re}(E_{\nu,k} a^{(\nu)}_{j,h} (x_1+ix_2)^{\alpha})
\end{equation}
for all $x = (x_1,x_2) \in \mathbb{R}^2$ and $y \in \mathbb{R}^{n-2}$, where $a^{(\nu)}_{j,h} \in \mathbb{C}^m \setminus \{0\}$.  In case~(c), $\widetilde{q}_{j,h} = 1$ and $\psi^{(\nu)}_{j,h} : \op{graph}\varphi^{(\nu,\infty)}_j \rightarrow \mathbb{R}^m$ and $\widetilde{\varphi}^{(\nu,k)}_{j,h} : \mathbb{R}^n \rightarrow \mathbb{R}^m$ are both identically zero.

\noindent\textit{Step~2.}  We claim the following: \textit{if $s > p$, then for each sufficiently large $\nu$ and sufficiently large $k$ (depending on $\nu$) }
\begin{align}\label{reverseHSshort step2 eqn1}
	&Q_{\varsigma/2,1}(T^{(\nu,k)}, \widetilde{\varphi}^{(\nu,k)}, \mathbf{C}_{1/2}(0)) 
	\\ \leq\,& \widehat{C} \overline{\beta}_{\nu} \inf_{\varphi' \in \bigcup_{p'=0}^{s-1} \Phi_{\alpha,q,p'}} 
		\op{dist}_{\mathcal H}(\op{spt}\op{graph} \widetilde{\varphi}^{(\nu,k)} \cap \mathbf{C}_{1/2}(0), 
		\op{spt}\op{graph} \varphi' \cap \mathbf{C}_{1/2}(0)) , \nonumber 
\end{align}
\textit{where $\widehat{C} = \widehat{C}(n,m,q,\alpha) \in (0,\infty)$ is a constant.}  

To prove the claim, let $\widehat{C} = \widehat{C}(n,m,q,\alpha) \in (0,\infty)$ be a constant to be later determined and suppose to the contrary that 
\begin{align}\label{reverseHSshort step2 eqn2}
	&Q_{\varsigma/2,1}(T^{(\nu,k)}, \widetilde{\varphi}^{(\nu,k)}, \mathbf{C}_{1/2}(0)) 
	\\ >\,& \widehat{C} \overline{\beta}_{\nu} \inf_{\varphi' \in \bigcup_{p'=0}^{s-1} \Phi_{\alpha,q,p'}} 
		\op{dist}_{\mathcal H}(\op{spt}\op{graph} \widetilde{\varphi}^{(\nu,k)} \cap \mathbf{C}_{1/2}(0), 
		\op{spt}\op{graph} \varphi' \cap \mathbf{C}_{1/2}(0)) . \nonumber
\end{align}
For each $\nu$ and $k$, choose $\widehat{s}^{(\nu,k)} \in \{0,1,\ldots,s-1\}$ and $\widehat{\varphi}^{(\nu,k)} \in \Phi_{\alpha,q,\widehat{s}^{(\nu,k)}}$ such that 
\begin{align}\label{reverseHSshort step2 eqn3}
	&\op{dist}_{\mathcal H}(\op{spt}\op{graph} \widetilde{\varphi}^{(\nu,k)} \cap \mathbf{C}_{1/2}(0), 
		\op{spt}\op{graph} \widehat{\varphi}^{(\nu,k)} \cap \mathbf{C}_{1/2}(0)) 
	\\ <\,& 2 \inf_{\varphi' \in \bigcup_{p'=0}^{s-1} \Phi_{\alpha,q,p'}} 
		\op{dist}_{\mathcal H}(\op{spt}\op{graph} \widetilde{\varphi}^{(\nu,k)} \cap \mathbf{C}_{1/2}(0), 
		\op{spt}\op{graph} \varphi' \cap \mathbf{C}_{1/2}(0)) . \nonumber
\end{align}
By \eqref{reverseHS step0 eqn1} with $\psi' = 0$, $F_{\nu}^2 \leq 2 \overline{M}^2$ and thus $\int_{B_{1/2}(0)} \sum_{j=1}^J |\psi^{(\nu)}_j|^2 \leq 4 \overline{M}^2 + 2$.  Since $s > p$, we can take $\varphi' = \varphi^{(\nu,k)}$ in \eqref{reverseHSshort step2 eqn3} and use the construction of $\widetilde{\varphi}^{(\nu,k)}$ to obtain 
\begin{align*}
	&\op{dist}_{\mathcal H}(\op{spt}\op{graph} \widetilde{\varphi}^{(\nu,k)} \cap \mathbf{C}_{1/2}(0), 
		\op{spt}\op{graph} \widehat{\varphi}^{(\nu,k)} \cap \mathbf{C}_{1/2}(0)) 
	\\ <\,& 2 \op{dist}_{\mathcal H}(\op{spt}\op{graph} \widetilde{\varphi}^{(\nu,k)} \cap \mathbf{C}_{1/2}(0), 
		\op{spt}\op{graph} \varphi^{(\nu,k)} \cap \mathbf{C}_{1/2}(0)) \leq C E_{\nu,k} 
\end{align*}
and thus by the triangle inequality 
\begin{equation}\label{reverseHSshort step2 eqn4}
	\op{dist}_{\mathcal H}(\op{spt}\op{graph} \widehat{\varphi}^{(\nu,k)} \cap \mathbf{C}_{1/2}(0), 
		\op{spt}\op{graph} \varphi^{(\nu,k)} \cap \mathbf{C}_{1/2}(0)) \leq C E_{\nu,k} , 
\end{equation}
where $C = C(n,m,q,\alpha,\overline{M}) \in (0,\infty)$ are constants.  Since by Theorem~\ref{graphical_thm}(A) $\inf_{x \in \mathbb{S}^1} \op{sep}\varphi^{(\nu,k)}(x,0) \geq c \beta_{\nu,k}^{-1} E_{\nu,k}$ where $c = c(n,m,q,\alpha) > 0$ is a constant, we deduce that $\widehat{s}^{(\nu,k)} \geq p$ and 
\begin{equation*}
    \widehat{\varphi}^{(\nu,k)} = \sum_{j=1}^J \sum_{h=1}^{\widehat{s}^{(\nu,k)}_j} \widehat{m}^{(\nu,k)}_{j,h} \widehat{\varphi}^{(\nu,k)}_{j,h}
\end{equation*}
where we let $\widehat{s}^{(\nu,k)}_j$, $\widehat{m}^{(\nu,k)}_{j,h}$, and $\widehat{q}^{(\nu,k)}_{j,h}$ be positive integers such that $\sum_{j=1}^J \widehat{s}^{(\nu,k)}_j = \widehat{s}^{(\nu,k)}$ and $\sum_{j=1}^J \sum_{h=1}^{\widehat{s}^{(\nu,k)}_j} \widehat{m}^{(\nu,k)}_{j,h} \widehat{q}^{(\nu,k)}_{j,h} = q$ and $ $$\widehat{\varphi}^{(\nu,k)}_{j,h} : \mathbb{R}^n \rightarrow \mathcal{A}_{\widehat{q}^{(\nu,k)}_j}(\mathbb{R}^m)$ are distinct $\widehat{q}^{(\nu,k)}_j$-valued functions such that: 
\begin{enumerate}[itemsep=2mm,topsep=0mm]
	\item[(i)] if $\varphi^{(\nu,\infty)}_j$ and $\varphi^{(\nu,k)}_j$ are non-zero, then $\widehat{q}^{(\nu,k)}_j = q_0$ and 
\begin{equation*}
	\widehat{\varphi}^{(\nu,k)}_{j,h}(x,y) = \op{Re}((c^{(\nu,k)}_j + E_{\nu,k} \widehat{a}^{(\nu,k)}_{j,h}) (x_1+ix_2)^{\alpha})
\end{equation*}
for all $x = (x_1,x_2) \in \mathbb{R}^2$ and $y \in \mathbb{R}^{n-2}$, where $\widehat{a}^{(\nu,k)}_{j,h} \in \mathbb{C}^m$ with $|\widehat{a}^{(\nu,k)}_{j,h}| \leq C$ (where $C$ is as in \eqref{reverseHSshort step2 eqn4}); 

    \item[(ii)]  if $\varphi^{(\nu,\infty)}_j$, $\varphi^{(\nu,0)}_j$, and $\varphi^{(\nu,k)}_j$ are identically zero and $\widehat{\varphi}^{(\nu,k)}_{j,h}$ is non-zero, then $\widehat{q}^{(\nu,k)}_j = q_0$ and 
\begin{equation*}
	\widehat{\varphi}^{(\nu,k)}_{j,h}(x,y) = \op{Re}(E_{\nu,k} \widehat{a}^{(\nu,k)}_{j,h} (x_1+ix_2)^{\alpha})
\end{equation*}
for all $x = (x_1,x_2) \in \mathbb{R}^2$ and $y \in \mathbb{R}^{n-2}$, where $\widehat{a}^{(\nu,k)}_{j,h} \in \mathbb{C}^m \setminus \{0\}$ with $|\widehat{a}^{(\nu,k)}_{j,h}| \leq C$ (where $C$ is again as in \eqref{reverseHSshort step2 eqn4}); 

    \item[(iii)]  $\varphi^{(\nu,\infty)}_j$, $\varphi^{(\nu,0)}_j$, $\varphi^{(\nu,k)}_j$, and $\widehat{\varphi}^{(\nu,k)}_{j,h}$ are all identically zero and $\widehat{q}^{(\nu,k)}_j = 1$.
\end{enumerate}
Note that we do not assert that $\widehat{\varphi}^{(\nu,k)}_j$ is associated with an element of $\mathfrak{L}(m_j,\varphi^{(\nu,0)}_j,\varphi^{(\nu,\infty)}_j)$ as in Lemma~\ref{integrability lemma}, as we might have $\sum_{h=1}^{\widehat{s}^{(\nu,k)}_j} \widehat{m}^{(\nu,k)}_{j,h} \widehat{q}^{(\nu,k)}_{j,h} > m_j q_j$ for some $j$.  

Since $\widehat{\varphi}^{(\nu,k)}$ has fewer non-zero components than $\widetilde{\varphi}^{(\nu,k)}$, there exists $j_1 \in \{1,2,\ldots,J\}$ for which there are fewer non-zero components $\widehat{\varphi}^{(\nu,k)}_{j_1,h}$ of $\widehat{\varphi}^{(\nu,k)}$ than there are non-zero components $\widetilde{\varphi}^{(\nu,k)}_{j_1,h}$ of $\widetilde{\varphi}^{(\nu,k)}$.  First let us suppose that $\varphi^{(\nu,k)}_{j_1}$ is non-zero.  Then for each $\theta \in [0,2\pi q_0)$ and $h \in \{1,2,\ldots,s_{j_1}\}$ there exists $h' \in \{1,2,\ldots,\widehat{s}^{(\nu,k)}_{j_1}\}$ such that 
\begin{equation}\label{reverseHSshort step2 eqn5}
	E_{\nu,k} |{\rm Re}((a^{(\nu)}_{j_1,h} - \widehat{a}^{(\nu)}_{j_1,h'}) \,e^{i\alpha\theta})| 
	= \op{dist}(\op{Re}((c^{(\nu,k)}_{j_1} + E_{\nu,k} a^{(\nu)}_{j_1,h}) \,e^{i\alpha\theta}), \op{spt} \widehat{\varphi}^{(\nu,k)}(e^{i\theta},0)) .  
\end{equation}
Since $\widehat{s}^{(\nu,k)}_{j_1} < s^{(\nu)}_{j_1}$, for each $\theta \in [0,2\pi)$ there exists $h_1,h_2 \in \{1,2,\ldots,s^{(\nu)}_{j_1}\}$ with $h_1 \neq h_2$ and $h' \in \{1,2,\ldots,\widehat{s}^{(\nu,k)}_{j_1}\}$ such that \eqref{reverseHSshort step2 eqn5} holds true with $h = h_1,h_2$.  It follows that there exists $h_1,h_2 \in \{1,2,\ldots,s^{(\nu)}_{j_1}\}$ with $h_1 \neq h_2$ and $h' \in \{1,2,\ldots,\widehat{s}^{(\nu,k)}_{j_1}\}$ such that the set $S$ of all $\theta \in [0,2\pi)$ such that \eqref{reverseHSshort step2 eqn5} holds true with $h = h_1,h_2$ satisfies $\mathcal{L}^1(S) \geq \pi/q^3$.  By \eqref{reverseHSshort step2 eqn5} 
\begin{equation*}
	E_{\nu,k} |{\rm Re}((a^{(\nu)}_{j_1,h} - \widehat{a}^{(\nu)}_{j_1,h'}) \,e^{i\alpha\theta})| \leq 2 \op{dist}_{\mathcal H}(\op{spt}\op{graph} \widetilde{\varphi}^{(\nu,k)} \cap \mathbf{C}_{1/2}(0), \op{spt}\op{graph} \widehat{\varphi}^{(\nu,k)} \cap \mathbf{C}_{1/2}(0))
\end{equation*}
for each $\theta \in S$.  Recalling that $\alpha = \ell_0/q_0$, the set $\{\theta \in [0,2\pi) : \sin(\alpha\theta) = 0 \} = (\pi/\alpha)\mathbb{Z} \cap [0,2\pi)$ consists of $2\ell_0$ points.  Since $\mathcal{L}^1(S) \geq \pi/q^3$, there exists $\theta_1,\theta_2 \in S$ such that the angular distance between $\theta_1$ and $\theta_2 + (\pi/\alpha)\mathbb{Z}$ is $\geq \pi/(8\alpha q^4)$ and thus $\sin(\alpha(\theta_1-\theta_2)) \geq \sin(\pi/(8q^4))$.  Hence 
\begin{equation}\label{reverseHSshort step2 eqn6}
	E_{\nu,k} |a^{(\nu)}_{j_1,h_1} - a^{(\nu)}_{j_1,h_2}| 
	\leq  C \op{dist}_{\mathcal H}(\op{spt}\op{graph} \widetilde{\varphi}^{(\nu,k)} \cap \mathbf{C}_{1/2}(0), 
		\op{spt}\op{graph} \widehat{\varphi}^{(\nu,k)} \cap \mathbf{C}_{1/2}(0)) ,  
\end{equation}
where $C = C(n,m,q,\alpha) \in (0,\infty)$ is a constant.  By setting $\psi'_{j_1} = \widetilde{m}_{j_1,h_1} \widetilde{\psi}^{(\nu)}_{j_1,h_2} + \sum_{h \neq h_2} \widetilde{m}_{j_1,h} \widetilde{\psi}^{(\nu)}_{j_1,h}$ and $\psi'_j = \psi^{(\nu)}_j$ if $j \neq j_1$
\begin{equation}\label{reverseHSshort step2 eqn7}
	\inf_{\psi' \in \bigcup_{s'=p}^{s-1} \mathfrak{L}_{s'}(m_j,\varphi^{(\nu,0)}_j, \varphi^{(\nu,\infty)}_j)} \int_{B_{1/2}(0)} \sum_{j=1}^J \sum_{l=1}^{q_j} \mathcal{G}(\psi^{(\nu,k)}_{j,l}, \psi'_{j,l})^2
    \leq C |a^{(\nu)}_{j_1,h_1} - a^{(\nu)}_{j_1,h_2}|^2   
\end{equation}
where $C = C(n,m,q,\alpha) \in (0,\infty)$ is a constant.  By combining \eqref{reverseHSshort step2 eqn2}, \eqref{reverseHSshort step2 eqn3}, \eqref{reverseHSshort step2 eqn6}, and \eqref{reverseHSshort step2 eqn7} 
\begin{align}\label{reverseHSshort step2 eqn8}
	&Q_{\varsigma/2,1}(T^{(\nu,k)}, \widetilde{\varphi}^{(\nu,k)}, \mathbf{C}_{1/2}(0))^2 
		\\ >\,& C \widehat{C}^2 \overline{\beta}_{\nu}^2 E_{\nu,k}^2 \inf_{\psi' \in \bigcup_{s'=p}^{s-1} \mathfrak{L}_{s'}(m_j,\varphi^{(\nu,0)}_j,
		\varphi^{(\nu,\infty)}_j)} \int_{B_{1/2}(0)} \sum_{j=1}^J \sum_{l=1}^{q_j} \mathcal{G}(\psi^{(\nu,k)}_{j,l}, \psi'_{j,l})^2 , \nonumber 
\end{align}
where $C = C(n,m,q,\alpha) > 0$ is a constant.  Dividing both sides of \eqref{reverseHSshort step2 eqn8} by $E_{\nu,k}^2$ and letting $k \rightarrow \infty$ using Lemma~\ref{blowup norms lemma} 
\begin{align}\label{reverseHSshort step2 eqn9}
	&\int_{B_{1/2}(0)} \sum_{j=1}^J \sum_{l=1}^{q_j} \mathcal{G}(w^{(\nu)}_{j,l},\psi^{(\nu)}_{j,l})^2
	\\ \geq\,& \tfrac{1}{2} C \widehat{C}^2 \overline{\beta}_{\nu}^2 \inf_{\psi' \in \bigcup_{s'=p}^{s-1} \mathfrak{L}_{s'}(m_j,\varphi^{(\nu,0)}_j,\varphi^{(\nu,\infty)}_j)}
		\int_{B_{1/2}(0)} \sum_{j=1}^J \sum_{l=1}^{q_j} \mathcal{G}(\psi^{(\nu,k)}_{j,l}, \psi'_{j,l})^2 . \nonumber 
\end{align}
where $C = C(n,m,q,\alpha) > 0$ is a constant.  On the other hand, by the triangle inequality and \eqref{reverseHS step0 eqn2}, for all $\nu$ large enough that $\overline{\beta}_{\nu} \leq 1/2$ 
\begin{equation*}
	\int_{B_{1/2}(0)} \sum_{j=1}^J \sum_{l=1}^{q_j} \mathcal{G}(w^{(\nu)}_{j,l},\psi^{(\nu)}_{j,l})^2
		\leq 4 \overline{\beta}_{\nu}^2 \inf_{\psi' \in \bigcup_{s'=p}^{s-1} \mathfrak{L}_{s'}(m_j,\varphi^{(\nu,0)}_j,\varphi^{(\nu,\infty)}_j)} 
		\int_{B_{1/2}(0)} \sum_{j=1}^J \sum_{l=1}^{q_j} \mathcal{G}(\psi^{(\nu,k)}_{j,l}, \psi'_{j,l})^2 , 
\end{equation*}
which provided $\widehat{C}$ is large enough that $C \widehat{C}^2 > 8$ (where $C$ is as in \eqref{reverseHSshort step2 eqn9}) gives us a contradiction.

Suppose instead that $\varphi^{(\nu,k)}_{j_1}$ is identically zero.  Then for each $\theta \in [0,2\pi)$, $h \in \{1,2,\ldots,s^{(\nu)}_{j_1}\}$, and $l \in \{0,1,\ldots,q_0-1\}$ such that $\widetilde{\psi}^{(\nu)}_{j_1,h}$ is non-zero either 
\begin{equation}\label{reverseHSshort step2 eqn10}
	E_{\nu,k} |{\rm Re}(a^{(\nu)}_{j_1,h} e^{i\alpha(\theta+2\pi l)})| 
	= \op{dist}(\op{Re}(E_{\nu,k} a^{(\nu)}_{j_1,h} e^{i\alpha(\theta+2\pi l)}), \op{spt} \widehat{\varphi}^{(\nu,k)}(e^{i\theta},0))  
\end{equation}
or there exist $l' \in \{0,1,\ldots,q_0-1\}$ and $h' \in \{1,2,\ldots,\widehat{s}^{(\nu,k)}_{j_1}\}$ such that $\widehat{\varphi}^{(\nu)}_{j_1,h'}$ is non-zero and 
\begin{equation}\label{reverseHSshort step2 eqn11}
	E_{\nu,k} |{\rm Re}(a^{(\nu)}_{j_1,h} e^{i\alpha\theta(\theta+2\pi l)} - \widehat{a}^{(\nu)}_{j_1,h'} e^{i\alpha(\theta+2\pi l')})|
	= \op{dist}(\op{Re}(E_{\nu,k} a^{(\nu)}_{j_1,h} \,e^{i\alpha(\theta+2\pi l)}), \op{spt} \widehat{\varphi}^{(\nu,k)}(e^{i\theta},0)) . 
\end{equation}
Since there are fewer non-zero components $\widehat{\varphi}^{(\nu,k)}_{j_1,h}$ than there are non-zero components $\widetilde{\varphi}^{(\nu,k)}_{j_1,h}$, either: 
\begin{enumerate}[itemsep=2mm,topsep=0mm]
	\item[(I)] there exists $h_1 \in \{1,2,\ldots,s^{(\nu)}_{j_1}\}$ and $l \in \{0,1,\ldots,q_0-1\}$ such that $\widetilde{\psi}^{(\nu)}_{j_1,h_1}$ is non-zero and the set $S$ of all $\theta \in [0,2\pi)$ such that \eqref{reverseHSshort step2 eqn10} holds true satisfies $\mathcal{L}^1(S) \geq \pi/(2q)$; 

   \item[(II)]  there exists $h_1,h_2 \in \{1,2,\ldots,s^{(\nu)}_{j_1}\}$, $h' \in \{1,2,\ldots,\widehat{s}^{(\nu,k)}_{j_1}\}$, and $l_1,l_2,l' \in \{0,1,\ldots,q_0-1\}$ such that $(h_1,l_1) \neq (h_2,l_2)$, $\widetilde{\psi}^{(\nu)}_{j_1,h_1}$, $\widetilde{\psi}^{(\nu)}_{j_1,h_2}$, and $\widehat{\varphi}^{(\nu,k)}_{j_1,h'}$ are non-zero and the set $S$ of all $\theta \in [0,2\pi)$ such that \eqref{reverseHSshort step2 eqn11} holds true with $(h,l)=(h_1.l_1),(h_2,l_2)$ satisfies $\mathcal{L}^1(S) \geq \pi/(2q^3)$
\end{enumerate}
In case~(I), arguing as we did to obtain \eqref{reverseHSshort step2 eqn6} above, it follows from \eqref{reverseHSshort step2 eqn10} that  
\begin{equation}\label{reverseHSshort step2 eqn12}
	E_{\nu,k} |a^{(\nu)}_{j_1,h_1}| 
	\leq  C \op{dist}_{\mathcal H}(\op{spt}\op{graph} \widetilde{\varphi}^{(\nu,k)} \cap \mathbf{C}_{1/2}(0), 
		\op{spt}\op{graph} \widehat{\varphi}^{(\nu,k)} \cap \mathbf{C}_{1/2}(0)) ,  
\end{equation}
where $C = C(n,m,q,\alpha) \in (0,\infty)$ is a constant.  By setting $\psi'_{j_1} = \widetilde{m}_{j_1,h_1} \llbracket 0 \rrbracket + \sum_{h \neq h_2} \widetilde{m}_{j_1,h} \widetilde{\psi}^{(\nu)}_{j_1,h}$ and $\psi'_j = \psi^{(\nu)}_j$ if $j \neq j_1$
\begin{equation}\label{reverseHSshort step2 eqn13}
	\inf_{\psi' \in \bigcup_{s'=p}^{s-1} \mathfrak{L}_{s'}(m_j,\varphi^{(\nu,0)}_j, \varphi^{(\nu,\infty)}_j)} \int_{B_{1/2}(0)} \sum_{j=1}^J \sum_{l=1}^{q_j} \mathcal{G}(\psi^{(\nu,k)}_{j,l}, \psi'_{j,l})^2
    \leq C |a^{(\nu)}_{j_1,h_1}|^2   
\end{equation}
where $C = C(n,m,q,\alpha) \in (0,\infty)$ is a constant.  Thus by combining \eqref{reverseHSshort step2 eqn2}, \eqref{reverseHSshort step2 eqn3}, \eqref{reverseHSshort step2 eqn12}, and \eqref{reverseHSshort step2 eqn13} we deduce that \eqref{reverseHSshort step2 eqn8} holds true.  
In case~(II), arguing as we did to obtain \eqref{reverseHSshort step2 eqn6}, it follows using \eqref{reverseHSshort step2 eqn11} that 
\begin{align}\label{reverseHSshort step2 eqn14}
	&E_{\nu,k} |a^{(\nu)}_{j_1,h_1} e^{i2\pi \alpha l_1} - a^{(\nu)}_{j_1,h_2} e^{i 2\pi \alpha l_2}| 
	\\ \leq\,&  C \op{dist}_{\mathcal H}(\op{spt}\op{graph} \widetilde{\varphi}^{(\nu,k)} \cap \mathbf{C}_{1/2}(0), 
		\op{spt}\op{graph} \widehat{\varphi}^{(\nu,k)} \cap \mathbf{C}_{1/2}(0)) ,  \nonumber 
\end{align}
where $C = C(n,m,q,\alpha) \in (0,\infty)$ is a constant.  If $h_1 = h_2$, then \eqref{reverseHSshort step2 eqn14} implies that \eqref{reverseHSshort step2 eqn12} holds true.  Otherwise, $h_1 \neq h_2$ and after replacing $a^{(\nu)}_{j_1,h}$ with $a^{(\nu)}_{j_1,h} e^{i2\pi\alpha l}$ for $(h,l) = (h_1,l_1),(h_2,l_2)$, we may assume that $l_1 = l_2 = 0$ so that \eqref{reverseHSshort step2 eqn14} is the same as \eqref{reverseHSshort step2 eqn6}.  Since either \eqref{reverseHSshort step2 eqn6} or \eqref{reverseHSshort step2 eqn12} holds true, by arguing as above \eqref{reverseHSshort step2 eqn8} holds true.  Having shown that \eqref{reverseHSshort step2 eqn8} holds true in both cases~(I) and (II), we can argue as above to reach a contradiction.

\noindent\textit{Step~3.}  Next we claim that: \textit{for each sufficiently large $\nu$ and sufficiently large $k$ (depending on $\nu$) }
\begin{align}\label{reverseHSshort step3 eqn1}
	&Q_{\varsigma/2,\varepsilon_{\nu,k}}(T^{(\nu,k)}, \widetilde{\varphi}^{(\nu,k)}, \mathbf{C}_{1/2}(0)) 
	\\ \leq\,& C \max\{\overline{M}\beta_{\nu,k},\overline{\beta}_{\nu}\} \inf_{\varphi' \in \bigcup_{p'=0}^{s-1} \Phi_{\alpha,q,p'}} Q_{\varsigma/2,\varepsilon_{\nu,k}}(T^{(\nu,k)}, \varphi', \mathbf{C}_{1/2}(0)) \nonumber 
\end{align}
\textit{for some constant $C = C(n,m,q,\alpha) \in (0,\infty)$.}  To see this, first note that by Lemma~\ref{blowup norms lemma} and \eqref{reverseHS step0 eqn1} with $\psi' = 0$, $Q_{\varsigma/2,\varepsilon_{\nu,k}}(T^{(\nu,k)},\widetilde{\varphi}^{(\nu,k)},\mathbf{C}_{1/2}(0)) \leq C(n) \,E_{\nu,k} F_{\nu} \leq C(n) \,\overline{M} E_{\nu,k}$.  Again using \eqref{reverseHS step0 eqn1} with $\psi' = 0$ together with $F_{\nu} > 0$, $w^{(\nu)}$ is not identically zero.  Hence arguing as in Remark~\ref{rescale blowup rmk}, condition~(C) of Definition~\ref{blowup class defn} holds true with $C\varepsilon_{\nu,k},C\beta_{\nu,k},2\delta_{\nu,k},\eta_{0,1/2\#} T^{(\nu,k)}, \varphi^{(\nu,k)}$ in place of $\varepsilon_{\nu},\beta_{\nu},\delta_{\nu},T^{(\nu)},\varphi^{(\nu)}$, where $C = C(n,m,q,\alpha) \in (0,\infty)$ is a constant.  In particular, 
\begin{equation*}
    Q_{\varsigma/2,\varepsilon_{\nu,k}}(T^{(\nu,k)},\widetilde{\varphi}^{(\nu,k)},\mathbf{C}_{1/2}(0)) \leq C \overline{M} E_{\nu,k} \leq C \overline{M}\beta_{\nu,k} \inf_{\varphi' \in \bigcup_{p'=0}^{p-1} \Phi_{\alpha,q,p'}} Q_{\varsigma/2,\varepsilon_{\nu,k}}(T^{(\nu,k)}, \varphi', \mathbf{C}_{1/2}(0))
\end{equation*}
so that if $p = s$ then \eqref{reverseHSshort step3 eqn1} holds true, and if $p > s$ then 
\begin{equation}\label{reverseHSshort step3 eqn2}
	\varepsilon_{\nu,k}^{-1} \widehat{E}(T^{(\nu,k)}, P_0, \mathbf{C}_{1/2}(0))^{1+\varsigma/2} 
	\leq C \overline{M} \beta_{\nu,k} \inf_{\varphi' \in \bigcup_{p'=0}^{p-1} \Phi_{\alpha,q,p'}} 
		Q_{\varsigma/2,\varepsilon_{\nu,k}}(T^{(\nu,k)},\varphi',\mathbf{C}_{1/2}(0)) 
\end{equation}
for some constant $C = C(n,m,q,\alpha) \in (0,\infty)$.  Suppose that $s > p$ and let $\beta_0 = \beta_0(n,m,q,\alpha,1/2)$ be as in Remark~\ref{graphical triangle rmk}(1).  Using Remark~\ref{tilde varphi rmk}, find $\overline{s}^{(\nu,k)} \in \{p,p+1,\ldots,s\}$ and $\overline{\varphi}^{(\nu,k)} \in \Phi_{\alpha,q,\overline{s}^{(\nu,k)}}$ such that 
\begin{equation*}
	Q_{\varsigma/2,1}(T^{(\nu,k)},\overline{\varphi}^{(\nu,k)},\mathbf{C}_{1/2}(0)) < 2^{q-1} \beta_0^{2-q} \inf_{\varphi' \in \bigcup_{p'=0}^{s-1} \Phi_{\alpha,q,p'}} Q_{\varsigma/2,1}(T^{(\nu,k)}, \varphi', \mathbf{C}_{1/2}(0))
\end{equation*} 
and Hypothesis~$(\star\star)$ holds true with $\eta_{0,1/2\#} T^{(\nu,k)}, 2^{1-\alpha} \overline{\varphi}^{(\nu,k)}$ in place of $T^{(\nu)}, \varphi^{(\nu)}$.  By \eqref{reverseHSshort step2 eqn1} and Remark~\ref{graphical triangle rmk}(1)
\begin{align}\label{reverseHSshort step3 eqn3}
    &Q_{\varsigma/2,1}(T^{(\nu,k)},\widetilde{\varphi}^{(\nu,k)},\mathbf{C}_{1/2}(0)) 
    \\ \leq\,& C \overline{\beta}_{\nu} \op{dist}_{\mathcal H}(\op{spt}\op{graph}\widetilde{\varphi}^{(\nu,k)} \cap \mathbf{C}_1(0),\op{spt}\op{graph}\overline{\varphi}^{(\nu,k)} \cap \mathbf{C}_1(0)) \nonumber 
    \\ \leq\,& C \overline{\beta}_{\nu} Q_{\varsigma/2,1}(T^{(\nu,k)},\widetilde{\varphi}^{(\nu,k)},\mathbf{C}_{1/2}(0)) + C \overline{\beta}_{\nu} Q_{\varsigma/2,1}(T^{(\nu,k)},\overline{\varphi}^{(\nu,k)},\mathbf{C}_{1/2}(0)) , \nonumber
\end{align}
where $C = C(n,m,q,\alpha) \in (0,\infty)$ are constants.  Noting that for $\nu$ sufficiently large we have $C\overline{\beta}_{\nu} < 1/2$ (where $C$ is as in \eqref{reverseHSshort step3 eqn3}), it follows that \eqref{reverseHSshort step3 eqn1} holds true with $1$ in place of $\varepsilon_{\nu,k}$, which together with \eqref{reverseHSshort step3 eqn2} is sufficient to complete the proof of \eqref{reverseHSshort step3 eqn1}.

\noindent\textit{Step~4.  Fine blow-up relative to the components of $\psi^{(\nu)}$.}  In light of conditions (A), (B), and (D) of Definition~\ref{blowup class defn} and \eqref{reverseHSshort step3 eqn1}, we can apply Corollary~\ref{nonconcentration cor} with $\eta_{0,1/2\#} T^{(\nu,k)}$ and $\widetilde{\varphi}^{(\nu,k)}$ in place of $T$ and $\varphi$ to deduce that for every $\tau \in (0,1/4)$ and every sufficiently large $\nu$ and $k$ 
\begin{align}\label{reverseHSshort step4 eqn1}
	&\int_{\mathbf{C}_{1/4}(0) \cap \{r < \tau\}} \op{dist}^2(X, \op{spt} \op{graph} \widetilde{\varphi}^{(\nu,k)}) \,d\|T^{(\nu,k)}\|(X) 
	\\ \leq\,& C \tau^{2-1/q} E_{\varsigma/2,1}(T^{(\nu,k)},\widetilde{\varphi}^{(\nu,k)},\mathbf{C}_{1/2}(0))^2 \nonumber 
\end{align}
where $C = C(n,m,q,\alpha) \in (0,\infty)$ is a constant.  Thus by dividing both sides of \eqref{reverseHSshort step4 eqn1} by $E_{\nu,k}^2$ and letting $k \rightarrow \infty$ using the definition of the blow-up $w^{(\nu)}$ and Lemma~\ref{blowup norms lemma}, for every $\tau \in (0,1/4)$ and every sufficiently large $\nu$ 
\begin{align}\label{reverseHSshort step4 eqn2}
	&\int_{B_{1/4}(0) \cap \{0 < r < \tau\}} \sum_{j=1}^J \sum_{l=1}^{q_j} \mathcal{G}(w^{(\nu)}_{j,l},\psi^{(\nu)}_{j,l})^2 
	\\ \leq\,& C \tau^{2-1/q} \int_{B_{1/2}(0)} \sum_{j=1}^J \sum_{l=1}^{q_j} \mathcal{G}(w^{(\nu)}_{j,l},\psi^{(\nu)}_{j,l})^2 = C \tau^{2-1/q} F_{\nu}^2 . \nonumber 
\end{align}

For each $\nu = 1,2,3,\ldots$, we want to construct blow-ups $\overline{w}^{(\nu)} = (\overline{w}^{(\nu)}_{j,h})$ of $(\eta_{0,1/2\#} T^{(\nu,k)})$ relative to $(2^{1-\alpha} \widetilde{\varphi}^{(\nu,k)})$ by excess $E_{\nu}$ and a corresponding $(\widetilde{m}_{j,h},\widetilde{\varphi}^{(\nu,0)}_{j,h},\widetilde{\varphi}^{(\nu,\infty)}_{j,h}) \in \mathfrak{D}$ as in Subsection~\ref{sec:blowup subsec}.  Then we want to let $\widetilde{\varphi}^{(\nu,0)}_{j,h} \rightarrow \widetilde{\varphi}^{(0)}_{j,h}$ and $\widetilde{\varphi}^{(\nu,\infty)}_{j,h} \rightarrow \widetilde{\varphi}^{(\infty)}_{j,h}$ uniformly on compact subsets of $\mathbb{R}^n$ as $k\rightarrow\infty$ and define $\overline{w}_{j,h} : \op{graph}\widetilde{\varphi}^{(\nu,\infty)}_{j,h} |_{B_{1/2}(0) \setminus \{0\}\times\mathbb{R}^{n-2}} \rightarrow \mathcal{A}_{\widetilde{m}_{j,h}}(\mathbb{R}^m)$ by, roughly speaking, $\overline{w}^{(\nu)}_{j,h}/F_{\nu} \rightarrow \overline{w}_{j,h}$ uniformly over compact subsets of $B_{1/2}(0) \setminus \{0\}\times\mathbb{R}^{n-2}$.  This construction proceeds by separately considering the cases (a), (b), and (c) from Step~1.

First observe that by conditions (A), (B), and (D) of Definition~\ref{blowup class defn} and \eqref{reverseHSshort step3 eqn1}, there exist sequences $(\widetilde{\tau}_{\nu})_{\nu=1}^{\infty}$ and $(\widetilde{\gamma}_{\nu})_{\nu=1}^{\infty}$ such that $0 < \widetilde{\tau}_{\nu} < \widetilde{\gamma}_{\nu} < 1$ and $\widetilde{\tau}_{\nu} \rightarrow 0$ and $\widetilde{\gamma}_{\nu} \rightarrow 1$ as $\nu\rightarrow\infty$ and for each sufficiently large $k$ there exists Lipschitz $\widetilde{m}_{j,h}$-valued functions $\widetilde{v}^{(\nu,k)}_{j,h} : {\rm graph}\big( 2^{1-\alpha} \widetilde{\varphi}^{(\nu,k)}_{j,h} |_{B_{\widetilde{\gamma}_{\nu}}(0) \cap \{r > \widetilde{\tau}_{\nu}\}} \big) \rightarrow \mathcal{A}_{\widetilde{m}_{j,h}}(\mathbb{R}^m)$ and closed sets $\widetilde{K}^{(\nu,k)}_{j,h} \subseteq B_{\widetilde{\gamma}_{\nu}}(0) \cap \{r > \widetilde{\tau}_{\nu}\}$ such that conclusion~(C) of Theorem~\ref{graphical_thm} holds true with $\widetilde{\gamma}_{\nu}, \widetilde{\tau}_{\nu}, \eta_{0,1/2\#} T^{(\nu,k)}, 2^{1-\alpha} \widetilde{\varphi}^{(\nu,k)}, \widetilde{v}^{(\nu,k)}_{j,h}, \widetilde{K}^{(\nu,k)}_{j,h}$ in place of $\gamma, \tau, T, \varphi, v_j,K_j$.

Let's consider case~(a), in which $j \in \{1,2,\ldots,J\}$ and $h \in \{1,2,\ldots,s_j\}$ such that $\varphi^{(\nu,k)}_j$ and $\varphi^{(\nu,\infty)}_j$ are non-zero.  Recall that $\widetilde{\psi}^{(\nu)}_{j,h}$ and $\widetilde{\varphi}^{(\nu,k)}_{j,h}$ are given by \eqref{reverseHSshort step1 eqn1} and \eqref{reverseHSshort step1 eqn4}.  By \eqref{blowups hypothesis2} and Remark~\ref{graphical rmk}(3), $\lim_{k\rightarrow\infty} E_{\nu,k}/\widehat{E}_{\nu,k} = 0$ and $\lim_{k\rightarrow\infty} E_{\nu,k}/|c^{(\nu,k)}_j| = 0$ and thus 
\begin{equation*}
	c^{(\nu,0)}_j = \lim_{k\rightarrow\infty} \frac{c^{(\nu,k)}_j}{\widehat{E}_{\nu,k}} 
		= \lim_{k\rightarrow\infty} \frac{c^{(\nu,k)}_j + E_{\nu,k} a^{(\nu)}_{j,h}}{\widehat{E}_{\nu,k}} , \quad 
    c^{(\nu,\infty)}_j = \lim_{k\rightarrow\infty} \frac{c^{(\nu,k)}_j}{|c^{(\nu,k)}_j|} 
		= \lim_{k\rightarrow\infty} \frac{c^{(\nu,k)}_j + E_{\nu,k} a^{(\nu)}_{j,h}}{|c^{(\nu,k)}_j + E_{\nu,k} a^{(\nu)}_j|} .
\end{equation*}
Hence we shall define  $\widetilde{\varphi}^{(\nu,0)}_{j,h} : \mathbb{R}^n \rightarrow \mathcal{A}_{q_0}(\mathbb{R}^m)$ and $\widetilde{\varphi}^{(\nu,\infty)}_{j,h} : \mathbb{R}^n \rightarrow \mathcal{A}_{q_0}(\mathbb{R}^m)$ by $\widetilde{\varphi}^{(\nu,0)}_{j,h} = c_{\nu} \varphi^{(\nu,0)}_j$ and $\widetilde{\varphi}^{(\nu,\infty)}_{j,h} = \varphi^{(\nu,\infty)}_j$ on $\mathbb{R}^n$, where $c_{\nu} = \lim_{k\rightarrow\infty} 2^{1-\alpha} \widehat{E}_{\nu,k}/\widehat{E}(T^{(\nu,k)},P_0,\mathbf{C}_{1/2}(0)) \in [2^{-n/2-\alpha},\infty)$ (which exists by \eqref{rescale blowup rmk eqn1}).  After passing to a subsequence, let $\widetilde{\varphi}^{(0)}_{j,h} : \mathbb{R}^n \rightarrow \mathcal{A}_{q_0}(\mathbb{R}^m)$ and $\widetilde{\varphi}^{(\infty)}_{j,h} : \mathbb{R}^n \rightarrow \mathcal{A}_{q_0}(\mathbb{R}^m)$ be given by $\widetilde{\varphi}^{(0)}_{j,h}(x,y) = \op{Re}(\widetilde{c}^{(0)}_{j,h} (x_1+ix_2)^{\alpha})$ and $\widetilde{\varphi}^{(\infty)}_{j,h}(x,y) = \op{Re}(\widetilde{c}^{(\infty)}_{j,h} (x_1+ix_2)^{\alpha})$ for each $x = (x_1,x_2) \in \mathbb{R}^2$ and $y \in \mathbb{R}^{n-2}$, where $\widetilde{c}^{(0)}_{j,h} = \lim_{\nu\rightarrow\infty} c_{\nu} c^{(\nu,0)}_j$ and $\widetilde{c}^{(\infty)}_{j,h} = \lim_{\nu\rightarrow\infty} c^{(\nu,\infty)}_j$, so that $\widetilde{\varphi}^{(\nu,0)}_{j,h} \rightarrow \widetilde{\varphi}^{(0)}_{j,h}$ and $\widetilde{\varphi}^{(\nu,\infty)}_{j,h} \rightarrow \widetilde{\varphi}^{(\infty)}_{j,h}$ uniformly on compact subsets of $\mathbb{R}^n$ as $\nu\rightarrow\infty$.  Recalling the definitions of $v^{(\nu,k)}_j$ and $\widetilde{v}^{(\nu,k)}_{j,h}$ via Theorem~\ref{graphical_thm}(C) and the definition of $\widetilde{\varphi}^{(\nu,k)}_{j,h}$ 
\begin{align}\label{reverseHSshort step4 eqn3}
	&v^{(\nu,k)}_j(re^{i\theta},y,\op{Re}(c^{(\nu,k)}_j r^{\alpha} e^{i\alpha\theta})) 
    \\=\,& \sum_{h=1}^{s_j} \sum_{\lambda=1}^{\widetilde{m}_{j,h}} 
		\big\llbracket {\rm Re}(E_{\nu,k} a^{(\nu)}_{j,h} r^{\alpha} e^{i\alpha\theta}) + \tfrac{1}{2} \widetilde{v}^{(\nu,k)}_{j,h,\lambda}(2re^{i\theta},2y,2\op{Re}((c^{(\nu,k)}_j + E_{\nu,k} a^{(\nu)}_{j,h}) \,r^{\alpha} e^{i\alpha\theta})) \big\rrbracket \nonumber 
\end{align}
whenever $(re^{i\theta},y) \in (2K^{(\nu,k)}_j) \cap \bigcap_{h=1}^{s_j} \widetilde{K}^{(\nu,k)}_{j,h}$, where $\widetilde{v}^{(\nu,k)}_{j,h}(X) = \sum_{l=1}^{\widetilde{m}_{j,h}} \llbracket \widetilde{v}^{(\nu,k)}_{j,h,\lambda}(X) \rrbracket$ for each $X \in \op{graph}\widetilde{\varphi}^{(\nu,k)}_{j,h} |_{B_{\widetilde{\gamma}_{\nu}}(0) \cap \{r > \widetilde{\tau}_{\nu}\}}$.  By Theorem~\ref{graphical_thm}(C) and Lemma~\ref{blowup norms lemma}, for each $\widetilde{\tau}_{\nu} \leq \tau < \gamma \leq \widetilde{\gamma}_{\nu}$ and each sufficiently large $k$, 
\begin{equation*}
	\sup_{B_{\gamma}(0) \cap \{r > \tau\}} \sum_{l=1}^{q_0} |\widetilde{v}^{(\nu,k)}_{j,h,l}|^2 
		+ \int_{B_{\gamma}(0) \cap \{r > \tau\}} \sum_{l=1}^{q_0} |D\widetilde{v}^{(\nu,k)}_{j,h,l}|^2 
	\leq C E_{\varsigma/2,1}(T^{(\nu,k)},\widetilde{\varphi}^{(\nu,k)},\mathbf{C}_{1/2}(0))^2 \leq C E_{\nu,k}^2 F_{\nu}^2 ,  
\end{equation*}
where $\widetilde{v}^{(\nu,k)}_{j,h,l}$ is as in \eqref{v_localized} with $\widetilde{\varphi}^{(\nu,k)}_{j,h}$ and $ \widetilde{v}^{(\nu,k)}_{j,h}$ in place of $\varphi_j$ and $v_j$ and $C = C(n,m,q,\alpha,\gamma,\tau) \in (0,\infty)$ are constants.  Hence by the Rellich compactness lemma for multi-valued Sobolev functions~\cite[Proposition~2.11]{DeLSpaDirMin} and~\cite[Theorem~5.2]{DeLSpa1}, after passing to a subsequence there exists $\overline{w}^{(\nu)}_{j,h} : \op{graph} \widetilde{\varphi}^{(\nu,\infty)}_{j,h} |_{B_{\widetilde{\gamma}_{\nu}}(0) \cap \{r > \widetilde{\tau}_{\nu}\}} \rightarrow \mathcal{A}_{\widetilde{m}_{j,h}}(\mathbb{R}^m)$ such that 
\begin{equation*}
	\widetilde{v}^{(\nu,k)}_{j,h}(re^{i\theta},y,2^{1-\alpha} \op{Re}((c^{(\nu,k)}_j + E_{\nu,k} a^{(\nu)}_{j,h}) \,r^{\alpha} e^{i\alpha\theta}))/E_{\nu,k} 
	\rightarrow \overline{w}^{(\nu)}_{j,h}(re^{i\theta},y,\op{Re}(c^{(\nu,\infty)}_j r^{\alpha} e^{i\alpha\theta}))
\end{equation*}
in the $L^2$-topology on $B_{\widetilde{\gamma}_{\nu}}(0) \cap \{r > \widetilde{\tau}_{\nu}\}$ and pointwise $\mathcal{L}^n$-a.e.~on $B_{\widetilde{\gamma}_{\nu}}(0) \cap \{r > \widetilde{\tau}_{\nu}\}$ as $k\rightarrow\infty$.  Moreover, for each $\widetilde{\tau}_{\nu} \leq \tau < \gamma \leq \widetilde{\gamma}_{\nu}$ 
\begin{equation}\label{reverseHSshort step4 eqn4}
	\sup_{B_{\gamma}(0) \cap \{r > \tau\}} \sum_{l=1}^{q_0} |\overline{w}^{(\nu)}_{j,h,l}|^2 \leq C F_{\nu}^2 , 
\end{equation}
where $\varphi^{(\nu,\infty)}_{j,l}$ is as in \eqref{varphi0 localized} with $\varphi^{(\nu,\infty)}_j$ in place of $\varphi^{(\infty)}_j$, $\overline{w}^{(\nu)}_{j,h,l}(x,y) = \overline{w}^{(\nu)}_{j,h}(x,y,\varphi^{(\nu,\infty)}_{j,l}(x,y))$ for each $(x,y) \in B_{\widetilde{\gamma}_{\nu}}(0) \cap \{r > \widetilde{\tau}_{\nu}\}$, and $C = C(n,m,q,\alpha,\gamma,\tau) \in (0,\infty)$ is a constant.  By dividing both sides of \eqref{reverseHSshort step4 eqn3} by $E_{\nu,k}$ and letting $k\rightarrow\infty$ 
\begin{align}\label{reverseHSshort step4 eqn5}
	&w^{(\nu)}_j(re^{i\theta},y,\op{Re}(c^{(\nu,\infty)}_j r^{\alpha} e^{i\alpha\theta})) 
    \\=\,& \sum_{h=1}^{s_j} \sum_{\lambda=1}^{\widetilde{m}_{j,h}} 
		\big\llbracket {\rm Re}(a^{(\nu)}_{j,h} r^{\alpha} e^{i\alpha\theta}) + \tfrac{1}{2} \overline{w}^{(\nu)}_{j,h,\lambda}(2re^{i\theta},2y,2^{\alpha} \op{Re}(c^{(\nu,\infty)}_j r^{\alpha} e^{i\alpha\theta})) \big\rrbracket \nonumber 
\end{align}
whenever $(re^{i\theta},y) \in B_{\widetilde{\gamma}_{\nu}/2}(0) \cap \{r > \widetilde{\tau}_{\nu}/2\}$, where we write $\overline{w}^{(\nu)}_{j,h}(X) = \sum_{\lambda=1}^{\widetilde{m}_{j,h}} \llbracket \overline{w}^{(\nu)}_{j,h,\lambda}(X) \rrbracket$ for each $X \in \op{graph} \widetilde{\varphi}^{(\nu,\infty)}_{j,h} |_{B_{\widetilde{\gamma}_{\nu}}(0) \cap \{r > \widetilde{\tau}_{\nu}\}}$.  By \eqref{reverseHSshort step4 eqn4} and the compactness of locally Dirichlet energy minimizing multi-valued functions, after passing to a subsequence there exists $\overline{w}_{j,h} : \op{graph} \widetilde{\varphi}^{(\infty)}_{j,h} |_{B_1(0) \setminus \{0\} \times \mathbb{R}^{n-2}} \rightarrow \mathcal{A}_{\widetilde{m}_{j,h}}(\mathbb{R}^m)$ such that 
\begin{equation}\label{reverseHSshort step4 eqn6}
	\overline{w}^{(\nu)}_{j,h}(re^{i\theta},y,\op{Re}(c^{(\nu,\infty)}_j r^{\alpha} e^{i\alpha\theta}))/F_{\nu} 
	\rightarrow \overline{w}_{j,h}(re^{i\theta},y,\op{Re}(\widetilde{c}^{(\infty)}_{j,h} r^{\alpha} e^{i\alpha\theta}))
\end{equation}
uniformly on compact subsets of $B_1(0) \setminus \{0\} \times \mathbb{R}^{n-2}$ as $\nu\rightarrow\infty$.  

Let's consider case~(b), in which $j \in \{1,2,\ldots,J\}$ and $h \in \{1,2,\ldots,s_j\}$ such that $\varphi^{(\nu,k)}_j$, $\varphi^{(\nu,0)}_j$, and $\varphi^{(\nu,\infty)}_j$ are identically zero but $\psi^{(\nu)}_{j,h}$ is non-zero.  Recall that $\widetilde{\psi}^{(\nu)}_{j,h}$ and $\widetilde{\varphi}^{(\nu,k)}_{j,h}$ are given by \eqref{reverseHSshort step1 eqn5}.  Since by \eqref{blowups hypothesis2} we have $\lim_{k\rightarrow\infty} E_{\nu,k} a^{(\nu)}_{j,h}/\widehat{E}_{\nu,k} = 0$, we define  $\widetilde{\varphi}^{(\nu,0)}_{j,h} : \mathbb{R}^n \rightarrow \mathcal{A}_{q_0}(\mathbb{R}^m)$ by $\widetilde{\varphi}^{(\nu,0)}_{j,h} = q_0 \llbracket 0 \rrbracket$ on $\mathbb{R}^n$.  Define $\widetilde{\varphi}^{(\nu,\infty)}_{j,h} : \mathbb{R}^n \rightarrow \mathcal{A}_{q_0}(\mathbb{R}^m)$ by $\widetilde{\varphi}^{(\nu,\infty)}_{j,h} = \op{Re}(\widetilde{c}^{(\nu,\infty)}_{j,h} (x_1+ix_2)^{\alpha})$ on $\mathbb{R}^n$, where $\widetilde{c}^{(\nu,\infty)}_{j,h} =  a^{(\nu)}_{j,h}/|a^{(\nu)}_{j,h}|$.  After passing to a subsequence, let $\widetilde{\varphi}^{(0)}_{j,h} : \mathbb{R}^n \rightarrow \mathcal{A}_{q_0}(\mathbb{R}^m)$ and $\widetilde{\varphi}^{(\infty)}_{j,h} : \mathbb{R}^n \rightarrow \mathcal{A}_{q_0}(\mathbb{R}^m)$ be given by $\widetilde{\varphi}^{(0)}_{j,h}(x,y) = q_0 \llbracket 0 \rrbracket$ and $\widetilde{\varphi}^{(\infty)}_{j,h}(x,y) = \op{Re}(\widetilde{c}^{(\infty)}_{j,h} (x_1+ix_2)^{\alpha})$ for each $x = (x_1,x_2) \in \mathbb{R}^2$ and $y \in \mathbb{R}^{n-2}$, where $\widetilde{c}^{(\infty)}_{j,h} = \lim_{\nu\rightarrow\infty} \widetilde{c}^{(\nu,\infty)}_{j,h}$, so that $\widetilde{\varphi}^{(\nu,0)}_{j,h} \rightarrow \widetilde{\varphi}^{(\infty)}_{j,h}$ and $\widetilde{\varphi}^{(\nu,0)}_{j,h} \rightarrow \widetilde{\varphi}^{(\infty)}_{j,h}$ uniformly on compact subsets of $\mathbb{R}^n$ as $\nu\rightarrow\infty$.  Arguing as in case~(a), after passing to a subsequence there exists $\overline{w}^{(\nu)}_{j,h} : \op{graph} \widetilde{\varphi}^{(\nu,\infty)}_{j,h} |_{B_{\widetilde{\gamma}_{\nu}}(0) \cap \{r > \widetilde{\tau}_{\nu}\}} \rightarrow \mathcal{A}_{\widetilde{m}_{j,h}}(\mathbb{R}^m)$ such that 
\begin{equation*}
	\widetilde{v}^{(\nu,k)}_{j,h}(re^{i\theta},y,2^{1-\alpha} \op{Re}(E_{\nu,k} a^{(\nu)}_{j,h} \,r^{\alpha} e^{i\alpha\theta}))/E_{\nu,k} 
	\rightarrow \overline{w}^{(\nu)}_{j,h}(re^{i\theta},y,\op{Re}(\widetilde{c}^{(\nu,\infty)}_{j,h} r^{\alpha} e^{i\alpha\theta}))
\end{equation*}
in the $L^2$-topology on $B_{\widetilde{\gamma}_{\nu}}(0) \cap \{r > \widetilde{\tau}_{\nu}\}$ and pointwise $\mathcal{L}^n$-a.e.~on $B_{\widetilde{\gamma}_{\nu}}(0) \cap \{r > \widetilde{\tau}_{\nu}\}$ as $k\rightarrow\infty$.  After passing to a further subsequence, there exists $\overline{w}_{j,h} : \op{graph} \widetilde{\varphi}^{(\infty)}_{j,h} |_{B_1(0) \setminus \{0\} \times \mathbb{R}^{n-2}} \rightarrow \mathcal{A}_{\widetilde{m}_{j,h}}(\mathbb{R}^m)$ such that 
\begin{equation}\label{reverseHSshort step4 eqn7}
	\overline{w}^{(\nu)}_{j,h}(re^{i\theta},y,\op{Re}(\widetilde{c}^{(\nu,\infty)}_{j,h} r^{\alpha} e^{i\alpha\theta}))/F_{\nu} 
	\rightarrow \overline{w}_{j,h}(re^{i\theta},y,\op{Re}(\widetilde{c}^{(\infty)}_{j,h} r^{\alpha} e^{i\alpha\theta}))
\end{equation}
uniformly on compact subsets of $B_1(0) \setminus \{0\} \times \mathbb{R}^{n-2}$ as $\nu\rightarrow\infty$.  Finally, in case~(c), $j \in \{1,2,\ldots,J\}$ and $h \in \{1,2,\ldots,s_j\}$ such that $\varphi^{(\nu,k)}_j$, $\varphi^{(\nu,0)}_j$, $\varphi^{(\nu,\infty)}_j$, and $\psi^{(\nu)}_{j,h}$ are identically zero.  Recalling that $\widetilde{\varphi}^{(\nu,k)}_{j,h}$ is identically zero, we let $\widetilde{\varphi}^{(\nu,0)}_{j,h} : \mathbb{R}^n \rightarrow \mathbb{R}^m$, $\widetilde{\varphi}^{(\nu,\infty)}_{j,h} : \mathbb{R}^n \rightarrow \mathbb{R}^m$, $\widetilde{\varphi}^{(0)}_{j,h} : \mathbb{R}^n \rightarrow \mathbb{R}^m$, and $\widetilde{\varphi}^{(\infty)}_{j,h} : \mathbb{R}^n \rightarrow \mathbb{R}^m$ all be identically zero.  After passing to a subsequence there exists $\overline{w}^{(\nu)}_{j,h} : \op{graph} \widetilde{\varphi}^{(\nu,\infty)}_{j,h} |_{B_{\widetilde{\gamma}_{\nu}}(0) \cap \{r > \widetilde{\tau}_{\nu}\}} \rightarrow \mathcal{A}_{\widetilde{m}_{j,h}}(\mathbb{R}^m)$ such that 
\begin{equation*}
	\widetilde{v}^{(\nu,k)}_{j,h}(x,y,0)/E_{\nu,k} 
	\rightarrow \overline{w}^{(\nu)}_{j,h}(x,y,0)
\end{equation*}
in the $L^2$-topology on $B_{\widetilde{\gamma}_{\nu}}(0) \cap \{r > \widetilde{\tau}_{\nu}\}$ and pointwise $\mathcal{L}^n$-a.e.~on $B_{\widetilde{\gamma}_{\nu}}(0) \cap \{r > \widetilde{\tau}_{\nu}\}$ as $k\rightarrow\infty$.  After passing to a further subsequence, there exists $\overline{w}_{j,h} : \op{graph} \widetilde{\varphi}^{(\infty)}_{j,h} |_{B_1(0) \setminus \{0\} \times \mathbb{R}^{n-2}} \rightarrow \mathcal{A}_{\widetilde{m}_{j,h}}(\mathbb{R}^m)$ such that 
\begin{equation}\label{reverseHSshort step4 eqn8}
	\overline{w}^{(\nu)}_{j,h}(x,y,0)/F_{\nu} 
	\rightarrow \overline{w}_{j,h}(x,y,0)
\end{equation}
uniformly on compact subsets of $B_1(0) \setminus \{0\} \times \mathbb{R}^{n-2}$ as $\nu\rightarrow\infty$.  Now if $j \in \{1,2,\ldots,J\}$ such that $\varphi^{(\nu,k)}_j$, $\varphi^{(\nu,0)}_j$, and $\varphi^{(\nu,\infty)}_j$ are identically zero, then in place of \eqref{reverseHSshort step4 eqn3} we have 
\begin{align}\label{reverseHSshort step4 eqn9}
	&v^{(\nu,k)}_j(re^{i\theta},y,0) 
    = \sum_{\{ h : \psi^{(\nu)}_{j,h} \equiv 0\}} \sum_{\lambda=1}^{\widetilde{m}_{j,h}} \big\llbracket \tfrac{1}{2} \overline{v}^{(\nu,k)}_{j,h,\lambda}(2re^{i\theta},2y,0) \big\rrbracket 
    \\&+ \sum_{\{ h : \psi^{(\nu)}_{j,h} \not\equiv 0\}} \sum_{l=0}^{q_0-1} \sum_{\lambda=1}^{\widetilde{m}_{j,h}} \big\llbracket {\rm Re}(\widetilde{c}^{(\nu,k)}_{j,h} r^{\alpha} e^{i\alpha(\theta+2\pi l)}) + \tfrac{1}{2} \overline{v}^{(\nu,k)}_{j,h,\lambda}(2re^{i\theta},2y,2\op{Re}(\widetilde{c}^{(\nu,k)}_{j,h} r^{\alpha} e^{i\alpha(\theta+2\pi l)})) \big\rrbracket , \nonumber 
\end{align}
where $\widetilde{c}^{(\nu,k)}_{j,h} = E_{\nu,k} a^{(\nu)}_{j,h}$, and in place of \eqref{reverseHSshort step4 eqn5} we have 
\begin{align}\label{reverseHSshort step4 eqn10}
	&w^{(\nu)}_j(re^{i\theta},y,0) 
    = \sum_{\{ h : \psi^{(\nu)}_{j,h} \equiv 0\}} \sum_{\lambda=1}^{\widetilde{m}_{j,h}} \big\llbracket \tfrac{1}{2} \overline{w}^{(\nu)}_{j,h,\lambda}(2re^{i\theta},2y,0) \big\rrbracket 
    \\&+ \sum_{\{ h : \psi^{(\nu)}_{j,h} \not\equiv 0\}} \sum_{l=0}^{q_0-1} \sum_{\lambda=1}^{\widetilde{m}_{j,h}} \big\llbracket {\rm Re}(a^{(\nu)}_{j,h} r^{\alpha} e^{i\alpha(\theta+2\pi l)}) + \tfrac{1}{2} \overline{w}^{(\nu)}_{j,h,\lambda}(2re^{i\theta},2y,2\op{Re}(\widetilde{c}^{(\nu,\infty)}_{j,h} r^{\alpha} e^{i\alpha(\theta+2\pi l)})) \big\rrbracket \nonumber 
\end{align}
for all $(re^{i\theta},y) \in B_{\widetilde{\gamma}_{\nu}/2}(0) \cap \{r > \widetilde{\tau}_{\nu}/2\}$.

By conditions (A), (B), and (D) of Definition~\ref{blowup class defn} and \eqref{reverseHSshort step3 eqn1}, we can select a suitable diagonal sequence $k = k(\nu)$ such that condition~(C) of Definition~\ref{blowup class defn} holds true with $C\varepsilon_{\nu,k(\nu)},C\overline{\beta}_{\nu}, 2\delta_{\nu,k(\nu)},$ $\eta_{0,1/2\#} T^{(\nu,k(\nu))}, 2^{1-\alpha} \widetilde{\varphi}^{(\nu,k(\nu))}$ in place of $\varepsilon_{\nu},\beta_{\nu},\delta_{\nu},T^{(\nu)},\varphi^{(\nu)}$, where $C = C(n,m,q,\alpha) \in (0,\infty)$ is a constant, and thus after passing to a subsequence there exists a blow-up $\widetilde{w} = (\widetilde{w}_{j,h}) \in \mathfrak{B}$ of $(\eta_{0,1/2\#} T^{(\nu,k(\nu))})$ relative to $(2^{1-\alpha}\widetilde{\varphi}^{(\nu,k(\nu))})$.  We additionally require that $k = k(\nu)$ is chosen so that in case~(a) 
\begin{equation}\label{reverseHSshort step4 eqn11}
    \lim_{\nu\rightarrow\infty} \frac{2^{1-\alpha} c^{(\nu,k(\nu))}_{j,h}}{\widehat{E}(T^{(\nu,k(\nu))},P_0,\mathbf{C}_{1/2}(0))} = \widetilde{c}^{(0)}_{j,h} \quad\text{and}\quad \lim_{\nu\rightarrow\infty} \frac{c^{(\nu,k(\nu))}_{j,h}}{|c^{(\nu,k(\nu))}_{j,h}|} = \widetilde{c}^{(\infty)}_{j,h} 
\end{equation}
and, recalling \eqref{reverseHSshort step1 eqn2}, \eqref{reverseHSshort step1 eqn3}, and Lemma~\ref{blowup norms lemma}, that 
\begin{equation}\label{reverseHSshort step4 eqn12}
    \int_{B_{1/2}(0)\ \cap \{r > 1/\nu\}} \sum_{j=1}^J \sum_{l=1}^{q_j} \mathcal{G}\left(\frac{v^{(\nu,k(\nu))}_{j,l}}{E_{\nu,k(\nu)}}, w^{(\nu)}\right)^2 < \frac{1}{\nu} \,F_{\nu}^2
\end{equation}
and the limit $\widetilde{c} = \lim_{\nu\rightarrow\infty} \widetilde{E}_{\nu}/(E_{\nu} F_{\nu}) \in [0,\omega_n^{-1/2} 2^{n/2+1}]$ exists.  \eqref{reverseHSshort step4 eqn11} guarantees that $\widetilde{w} \in \mathfrak{B}(\widetilde{m}_{j,h},\widetilde{\varphi}^{(0)}_{j,h},\widetilde{\varphi}^{(\infty)}_{j,h})$ where $(\widetilde{m}_{j,h},\widetilde{\varphi}^{(0)}_{j,h},\widetilde{\varphi}^{(\infty)}_{j,h})$ is as constructed above.  By \eqref{reverseHSshort step4 eqn3}, \eqref{reverseHSshort step4 eqn5}, \eqref{reverseHSshort step4 eqn9}, \eqref{reverseHSshort step4 eqn10}, and \eqref{reverseHSshort step4 eqn12} 
\begin{align*}
    &\lim_{\nu\rightarrow\infty} \int_{B_{\gamma}(0)\ \cap \{r > \tau\}} \sum_{j=1}^J \sum_{l=1}^{q_j} \mathcal{G}\left(\frac{\widetilde{v}^{(\nu,k(\nu))}_{j,l}}{E_{\nu,k(\nu)} F_{\nu}}, \frac{\overline{w}^{(\nu)}_{j,l}}{F_{\nu}}\right)^2
    \\=\,& \lim_{\nu\rightarrow\infty} \frac{1}{F_{\nu}^2} \int_{B_{\gamma/2}(0)\ \cap \{r > \tau/2\}} \sum_{j=1}^J \sum_{l=1}^{q_j} \mathcal{G}\left(\frac{v^{(\nu,k(\nu))}_{j,l}}{E_{\nu,k(\nu)}}, w^{(\nu)}_{j,l}\right)^2 = 0 
\end{align*}
so that $\overline{w}_{j,h} = \widetilde{c} \,\widetilde{w}_{j,h}$ on $\op{graph}\widetilde{\varphi}^{(\infty)}_{j,h} |_{B_1(0) \setminus \{0\} \times \mathbb{R}^{n-2}}$.

\noindent\textit{Step~5.  Homogeneity of the fine blow-up relative to the cones associated with $\psi^{(\nu)}$, and the contradiction that completes the proof.}  We claim that $\overline{w} \in \widetilde{\mathfrak L}(\widetilde{m}_{j,h},\widetilde{\varphi}^{(0)}_{j,h},\widetilde{\varphi}^{(\infty)}_{j,h})$.  To see this, without loss of generality assume that $\overline{w}_{j,h}$ is not identically zero for some $j$ and $h$.  By dividing both sides of \eqref{reverseHS step0 eqn3} by $F_{\nu}^2$ and letting $\nu\rightarrow\infty$ using the homogeneity of $\psi^{(\nu)}$, \eqref{reverseHSshort step4 eqn5}, \eqref{reverseHSshort step4 eqn10}, \eqref{reverseHSshort step4 eqn6}, \eqref{reverseHSshort step4 eqn7} and \eqref{reverseHSshort step4 eqn8}, $\overline{w}$ is homogeneous degree $\alpha$ in $(B_{1/2}(0) \setminus B_{1/8}(0)) \setminus \{0\} \times \mathbb{R}^{n-2}$.  Since $\overline{w}_{j,h} = \widetilde{c} \,\widetilde{w}_{j,h}$, $\widetilde{w}$ is homogeneous degree $\alpha$ in $(B_1(0) \setminus B_{1/4}(0)) \setminus \{0\} \times \mathbb{R}^{n-2}$ and component-wise Dirichlet energy minimizing in $B_1(0) \setminus \{0\} \times \mathbb{R}^{n-2}$.  By unique continuation for locally Dirichlet energy minimizing multi-valued functions, $\widetilde{w}$ is homogeneous degree $\alpha$ in $B_1(0) \setminus \{0\} \times \mathbb{R}^{n-2}$.  By Theorem~\ref{classify in L thm}, $\widetilde{w} \in \widetilde{\mathfrak L}(\widetilde{m}_{j,h},\widetilde{\varphi}^{(0)}_{j,h},\widetilde{\varphi}^{(\infty)}_{j,h})$.  Since $\overline{w}_{j,h} = c \,\widetilde{w}_{j,h}$, $\overline{w} \in \widetilde{\mathfrak L}(\widetilde{m}_{j,h},\widetilde{\varphi}^{(0)}_{j,h},\widetilde{\varphi}^{(\infty)}_{j,h})$. 

We define $\overline{\psi}^{(\nu)} = (\overline{\psi}^{(\nu)}_j) \in \widetilde{\mathfrak L}(m_j,\varphi^{(\nu,0)}_j,\varphi^{(\nu,\infty)}_j)$ as follows.  If $\varphi^{(\nu)}_j$ and $\varphi^{(\nu,\infty)}_j$ are not identically zero for all $\nu = 1,2,3,\ldots$, then by Definition~\ref{homogeneous class defn} 
\begin{equation*}
	\overline{w}_{j,h}(re^{i\theta},y,\op{Re}(\widetilde{c}^{(\infty)}_{j,h} r^{\alpha} e^{i\alpha\theta})) 
	= \sum_{\lambda=1}^{\widetilde{m}_{j,h}} \llbracket \op{Re}(\overline{a}_{j,h,\lambda} r^{\alpha} e^{i\alpha\theta} 
		+ \alpha \widetilde{c}^{(0)}_{j,h} r^{\alpha-1} e^{i(\alpha-1)\theta} (\overline{b} \cdot y)) \rrbracket 
\end{equation*}
for all $r > 0$, $\theta \in [0,2\pi q_0)$, and $y \in \mathbb{R}^{n-2}$, where $\overline{a}_{j,h,\lambda} \in \mathbb{C}^m$ and $\overline{b} \in \mathbb{C}^{n-2}$.  Let 
\begin{align*}
    &\overline{\psi}^{(\nu)}_j(re^{i\theta},y,\op{Re}(c^{(\nu,\infty)}_j r^{\alpha} e^{i\alpha\theta})) 
	\\=\,& \sum_{h=1}^{s_j} \sum_{\lambda=1}^{\widetilde{m}_{j,h}} \llbracket \op{Re}((a^{(\nu)}_{j,h} 
		+ 2^{\alpha-1} F_{\nu} \overline{a}_{j,h,\lambda}) \,r^{\alpha} e^{i\alpha\theta} 
		+ \alpha 2^{\alpha-1} F_{\nu} \widetilde{c}^{(0)}_{j,h} r^{\alpha-1} e^{i(\alpha-1)\theta} (\overline{b} \cdot y)) \rrbracket 
\end{align*}
for all $r > 0$, $\theta \in [0,2\pi q_0)$, and $y \in \mathbb{R}^{n-2}$.  If instead $\varphi^{(\nu)}_j$, $\varphi^{(\nu,0)}_j$, and $\varphi^{(\nu,\infty)}_j$ are identically zero for all $\nu = 1,2,3,\ldots$, let 
\begin{align*}
	&\overline{\psi}^{(\nu)}_j(x,y,0) 
    = \sum_{\{ h : \psi^{(\nu)}_{j,h} \equiv 0\}} \sum_{\lambda=1}^{\widetilde{m}_{j,h}} \llbracket 2^{\alpha-1} F_{\nu} \overline{w}_{j,h,\lambda}(x,y,0) \rrbracket    
    \\&+ \sum_{\{ h : \psi^{(\nu)}_{j,h} \not\equiv 0\}} \sum_{l=0}^{q_0-1} \sum_{\lambda=1}^{\widetilde{m}_{j,h}} \llbracket \op{Re}(a^{(\nu)}_{j,h} r^{\alpha} e^{i\alpha(\theta+2\pi l)}) + 2^{\alpha-1} F_{\nu} \overline{w}_{j,h,\lambda}(re^{i\theta},y,\op{Re}(\widetilde{c}^{(\nu,\infty)}_{j,h} r^{\alpha} e^{i\alpha(\theta+2\pi l)})) \rrbracket 
\end{align*} 
for all $(re^{i\theta},y) \in \mathbb{R}^n \setminus \{0\} \times \mathbb{R}^n$.  That way, by \eqref{reverseHSshort step4 eqn5}, \eqref{reverseHSshort step4 eqn10}, \eqref{reverseHSshort step4 eqn6}, \eqref{reverseHSshort step4 eqn7} and \eqref{reverseHSshort step4 eqn8}, for every $\tau \in (0,1/4)$
\begin{align}\label{reverseHSshort step5 eqn1}
	&\lim_{\nu\rightarrow\infty} \frac{1}{F_{\nu}^2} \int_{B_{1/4}(0) \cap \{r > \tau\}} \sum_{j=1}^J \sum_{l=1}^{q_j} \mathcal{G}(w^{(\nu)}_{j,l},\overline{\psi}^{(\nu)}_{j,l})^2 
	\\=\,& \lim_{\nu\rightarrow\infty} \frac{1}{2^{n+2}} \int_{B_{1/2}(0) \cap \{r > 2\tau\}} \sum_{j=1}^J \sum_{h=1}^{s_j} \sum_{l=1}^{\widetilde{q}_j} \mathcal{G}\Bigg(\frac{\overline{w}^{(\nu)}_{j,h,l}}{F_{\nu}},\overline{w}_{j,h,l}\Bigg)^2 
	= 0 . \nonumber 
\end{align}
Arguing as we did to obtain (5.84) of~\cite{KrumWicb}, 
\begin{align}\label{reverseHSshort step5 eqn2}
	\int_{B_{1/2}(0)} \sum_{j=1}^J \sum_{l=1}^{q_j} \mathcal{G}(w^{(\nu)}_{j,l},\overline{\psi}^{(\nu)}_{j,l})^2 
	\leq\,& C \int_{B_{1/4}(0)} \sum_{j=1}^J \sum_{l=1}^{q_j} \mathcal{G}(w^{(\nu)}_{j,l},\overline{\psi}^{(\nu)}_{j,l})^2 
	\\&+ C \int_{B_{1/2}(0) \setminus B_{1/8}(0)} \sum_{j=1}^J \sum_{l=1}^{q_j} \Bigg| \frac{\partial (w^{(\nu)}_{j,l}/R)}{\partial R} \Bigg|^2 , \nonumber
\end{align}
where $\overline{\psi}^{(\nu)}_{j,l}$ is as in \eqref{psi localized} with $\varphi^{(\nu,\infty)}_j$ and $\overline{\psi}^{(\nu)}_j$ in place of $\varphi^{(\infty)}_j$ and $\psi_j$ and $C = C(n,\alpha) \in (0,\infty)$ is a constant.  By the construction of $\overline{\psi}^{(\nu)}$, \eqref{reverseHSshort step4 eqn2}, \eqref{reverseHSshort step5 eqn1}, and \eqref{reverseHS step0 eqn3}, the right-hand side of \eqref{reverseHSshort step5 eqn2} is $o(F_{\nu}^2)$ as $\nu\rightarrow\infty$.  But by \eqref{reverseHS step0 eqn1} with $\psi' = \overline{\psi}^{(\nu)}$, the left-hand side of \eqref{reverseHSshort step5 eqn2} is $\geq c(n,m,q,\alpha)\,\overline{M}^{-2} F_{\nu}^2$ for all sufficiently large $\nu$, giving us a contradiction.
\end{proof}

\begin{proof}[Proof of Lemma~\ref{reverseHS lemma}]
Lemma~\ref{reverseHS lemma} follows from Lemma~\ref{reverseHS lemma0} by arguing as in~\cite[Lemma~11.1]{KrumWic2} (for which Lemma~11.1 follows from Lemma~11.5).
\end{proof}

\begin{lemma}\label{hole filling lemma}
Let $(m_j,\varphi^{(0)}_j,\varphi^{(\infty)}_j) \in \mathfrak{D}$ and $w = (w_j) \in \mathfrak{B}(m_j,\varphi^{(0)}_j,\varphi^{(\infty)}_j)$.  For each $\vartheta \in (0,1/8]$, let $\psi^{(\vartheta)} = (\psi^{(\vartheta)}_j) \in \widetilde{\mathfrak L}(m_j,\varphi^{(0)}_j,\varphi^{(\infty)}_j)$ such that 
\begin{equation}\label{hole filling proj}
	\int_{B_{\vartheta}(0)} \sum_{j=1}^J \sum_{l=1}^{q_j} \mathcal{G}(w_{j,l},\psi^{(\vartheta)}_{j,l})^2 
	\leq 2 \inf_{\psi' \in \widetilde{\mathfrak L}(m_j,\varphi^{(0)}_j,\varphi^{(\infty)}_j)} 
		\int_{B_{\vartheta}(0)} \sum_{j=1}^J \sum_{l=1}^{q_j} \mathcal{G}(w_{j,l},\psi'_{j,l})^2
\end{equation}
where $w_{j,l}$ is as in \eqref{w localized} and $\psi^{(\vartheta)}_{j,l}(x,y) = \psi^{(\vartheta)}_j(x,y,\varphi^{(\infty)}_{j,l}(x,y))$ and $\psi'_{j,l}(x,y) = \psi'_j(x,y,\varphi^{(\infty)}_{j,l}(x,y))$ for each $(x,y) \in B_{1/2}(0)$ where $\varphi^{(\infty)}_{j,l}$ is as in \eqref{varphi0 localized} (as in \eqref{psi localized}).  Then for each $\vartheta_1,\vartheta_2$ with $0 < \vartheta_1 \leq \vartheta_2 \leq 1/8$ 
\begin{equation}\label{hole filling concl}
	\vartheta_1^{-n-2\alpha} \int_{B_{\vartheta_1}(0)} \sum_{j=1}^J \sum_{l=1}^{q_j} \mathcal{G}(w_{j,l},\psi^{(\vartheta_1)}_{j,l})^2 
	\leq C \bigg(\frac{\vartheta_1}{\vartheta_2}\bigg)^{2\mu} \vartheta_2^{-n-2\alpha} \int_{B_{\vartheta_2}(0)} \sum_{j=1}^J \sum_{l=1}^{q_j} 
		\mathcal{G}(w_{j,l},\psi^{(\vartheta_2)}_{j,l})^2 ,
\end{equation}
where $\mu = \mu(n,m,q,\alpha) \in (0,1)$ and $C = C(n,m,q,\alpha) \in (0,\infty)$ are constants (independent of $\vartheta_1$ and $\vartheta_2$).
\end{lemma}

\begin{proof}
Without loss of generality, assume that $\vartheta_2 < 1/2$ and $w$ is not identically zero on $B_{\vartheta_1}(0)$.  Let $\rho \in [\vartheta_1,1/2)$.  By \eqref{blowup estimates concl1} with $\gamma = 1/4$ and $\psi^{(\rho)}$ in place of $\psi$ 
\begin{equation}\label{hole filling eqn1} 
	\int_{B_{\rho/4}(0)} \sum_{j=1}^J \sum_{l=1}^{q_j} R^{2-n} \left|\frac{\partial (w_{j,l}/R^{\alpha})}{\partial R}\right|^2 
	\leq C \rho^{-n-2\alpha} \int_{B_{\rho}(0)} \sum_{j=1}^J \sum_{l=1}^{q_j} \mathcal{G}(w_{j,l},\psi^{(\rho)}_{j,l})^2  , 
\end{equation}
where $C = C(n,m,q,\alpha) \in (0,\infty)$ is a constant.  Let $w^{(\rho)} = (w^{(\rho)}_1,\ldots,w^{(\rho)}_J)$ so that $w^{(\rho)}_j : \op{graph}\varphi^{(\infty)}_j \rightarrow \mathcal{A}_{m_j}(\mathbb{R}^m)$ is given by 
\begin{equation*}
	w^{(\rho)}_j(x,y,\varphi^{(\infty)}_{j,l}(x,y)) = \frac{w_j(2\rho x, 2\rho y, \varphi^{(\infty)}_{j,l}(2\rho x, 2\rho y))}{(2\rho)^{-n/2} \|w\|_{L^2(B_{2\rho}(0))}} 
\end{equation*}
for all $(x,y) \in B_1(0)$, $j \in \{1,\ldots,J\}$, and $l \in \{1,\ldots,q_j\}$, where $\varphi^{(\infty)}_{j,l}$ is as in \eqref{varphi0 localized}.  By Remark~\ref{rescale blowup rmk}, since $w$ is not identically zero in $B_{\rho}(0)$, $w^{(\rho)} \in \mathfrak{B}(m_j,c\varphi^{(0)}_j,\varphi^{(\infty)}_j)$ where $c = \lim_{\nu\rightarrow\infty} (2\rho)^{\alpha-1} \widehat{E}(T^{(\nu)},$ $P_0,\mathbf{C}_1(0))/\widehat{E}(T^{(\nu)},P_0,\mathbf{C}_{2\rho}(0)) \in (0,\infty)$.  By Remark~\ref{L scale varphi0 rmk}, $\psi^{(\rho)}_j \in \widetilde{\mathfrak L}(m_j,c\varphi^{(0)}_j,\varphi^{(\infty)}_j)$.  By Lemma~\ref{reverseHS lemma} with $w^{(\rho)}$ in place of $w$ 
\begin{equation}\label{hole filling eqn2} 
	\rho^{-n-2\alpha} \int_{B_{\rho}(0)} \sum_{j=1}^J \sum_{l=1}^{q_j} \mathcal{G}(w_{j,l},\psi^{(\rho)}_{j,l})^2 
	\leq C \int_{B_{\rho}(0) \setminus B_{\rho/4}(0)} \sum_{j=1}^J \sum_{l=1}^{q_j} R^{2-n} \left|\frac{\partial (w_{j,l}/R^{\alpha})}{\partial R}\right|^2 , 
\end{equation}
where $C = C(n,m,q,\alpha) \in (0,\infty)$ is a constant.  By combining \eqref{hole filling eqn1} and \eqref{hole filling eqn2} 
\begin{equation}\label{hole filling eqn3} 
	\int_{B_{\rho/4}(0)} \sum_{j=1}^J \sum_{l=1}^{q_j} R^{2-n} \left|\frac{\partial (w_{j,l}/R^{\alpha})}{\partial R}\right|^2 
	\leq C_0 \int_{B_{\rho}(0) \setminus B_{\rho/4}(0)} \sum_{j=1}^J \sum_{l=1}^{q_j} R^{2-n} \left|\frac{\partial (w_{j,l}/R^{\alpha})}{\partial R}\right|^2 , 
\end{equation}
where $C_0 = C_0(n,m,q,\alpha) \in (0,\infty)$ is a constant.  Adding $C_0$ times the left-hand side of \eqref{hole filling eqn3} to both sides of \eqref{hole filling eqn3} 
\begin{equation}\label{hole filling eqn4} 
	\int_{B_{\rho/4}(0)} \sum_{j=1}^J \sum_{l=1}^{q_j} R^{2-n} \left|\frac{\partial (w_{j,l}/R^{\alpha})}{\partial R}\right|^2 
	\leq \gamma \int_{B_{\rho}(0)} \sum_{j=1}^J \sum_{l=1}^{q_j} R^{2-n} \left|\frac{\partial (w_{j,l}/R^{\alpha})}{\partial R}\right|^2 
\end{equation}
for all $\rho \in [\vartheta_1,1/2)$, where $\gamma = C_0/(1+C_0) \in (0,1)$.  Noting that \eqref{hole filling concl} obviously holds true if $\vartheta_1 \geq \vartheta_2/4$, we may assume that $\vartheta_1 < \vartheta_2/4$.  Let $N$ be a positive integer such that $4^{-N-1} \vartheta_2 < \vartheta_1 \leq 4^{-N} \vartheta_2$.  By iteratively applying \eqref{hole filling eqn4} with $\rho = 4^{-i} \vartheta_2$ for $i = 1,2,\ldots,N-1$ 
\begin{equation}\label{hole filling eqn5} 
	\int_{B_{\vartheta_1}(0)} \sum_{j=1}^J \sum_{l=1}^{q_j} R^{2-n} \left|\frac{\partial (w_{j,l}/R^{\alpha})}{\partial R}\right|^2 
	\leq C\bigg(\frac{\vartheta_1}{\vartheta_2}\bigg)^{2\mu} \int_{B_{\vartheta_2/4}(0)} \sum_{j=1}^J \sum_{l=1}^{q_j} 
		R^{2-n} \left|\frac{\partial (w_{j,l}/R^{\alpha})}{\partial R}\right|^2 
\end{equation}
where $\mu = -\log\gamma/\log 16$ and $C = C(n,m,q,\alpha) \in (0,\infty)$ is a constant.  Combining \eqref{hole filling eqn2} with $\rho = \vartheta_1$, \eqref{hole filling eqn5}, and \eqref{hole filling eqn1} with $\rho = \vartheta_2$, we obtain \eqref{hole filling concl}.
\end{proof}

\begin{theorem}\label{blowup decay thm} 
Let $(m_j,\varphi^{(0)}_j,\varphi^{(\infty)}_j) \in \mathfrak{D}$ and $w = (w_j) \in \mathfrak{B}(m_j,\varphi^{(0)}_j,\varphi^{(\infty)}_j)$.  There exists a unique $\psi = (\psi_j) \in \widetilde{\mathfrak L}(m_j,\varphi^{(0)}_j,\varphi^{(\infty)}_j)$ such that 
\begin{equation}\label{blowup decay concl} 
	\frac{1}{\rho^{n+2\alpha}} \int_{B_{\rho}(0)} \sum_{j=1}^J \sum_{l=1}^{q_j} \mathcal{G}(w_{j,l},\psi_{j,l})^2 
	\leq C \rho^{2\mu} \int_{B_{1/2}(0)} \sum_{j=1}^J \sum_{l=1}^{q_j} |w_{j,l}|^2 
\end{equation}
where $w_{j,l}$ is as in \eqref{w localized} and $\psi_{j,l}$ is as in \eqref{psi localized}, and $\mu = \mu(n,m,q,\alpha) \in (0,1)$ and $C = C(n,m,q,\alpha) \in (0,\infty)$ are constants.  Additionally, $\|\psi\|_{L^2(B_1(0))} \leq C$ for some constant $C = C(n,m,q,\alpha) \in (0,\infty)$.
\end{theorem}

\begin{proof}
For each $k = 1,2,3,\ldots$, let $\psi^{(8^{-k})} = (\psi^{(8^{-k})}_j) \in \widetilde{\mathfrak L}(m_j,\varphi^{(0)}_j,\varphi^{(\infty)}_j)$ satisfy \eqref{hole filling proj} with $\vartheta = 8^{-k}$.  By Lemma~\ref{hole filling lemma} with $\vartheta_1 \in \{8^{-k},8^{-k-1}\}$ and $\vartheta_2 = 1/8$ together with the homogeneity of $\psi^{(8^{-k})},\psi^{(8^{-k-1})}$ and the triangle inequality 
\begin{equation*}
	\int_{B_1(0)} \sum_{j=1}^J \sum_{l=1}^{q_j} \mathcal{G}(\psi^{(8^{-k})}, \psi^{(8^{-k-1})})^2 
	\leq C 8^{-2k\mu} \int_{B_{1/2}(0)} \sum_{j=1}^J \sum_{l=1}^{q_j} |w_{j,l}|^2 
\end{equation*}
for all $k = 1,2,3,\ldots,$ where $\mu$ is as in Lemma~\ref{hole filling lemma} and $C = C(n,m,q,\alpha) \in (0,\infty)$ is a constant.  Hence there exists $\psi = (\psi_j) \in \widetilde{\mathfrak L}(m_j,\varphi^{(0)}_j,\varphi^{(\infty)}_j)$ such that $\psi^{(8^{-k})}_j \rightarrow \psi_j$ uniformly on $\op{graph}\varphi^{(\infty)}_j |_{B_1(0)}$ and thus \eqref{blowup decay concl} holds true and $\|\psi\|_{L^2(B_1(0))}^2 \leq C \|w\|_{L^2(B_{1/2}(0))}^2 \leq C\omega_n (1/2)^{n+2}$, where $C = C(n,m,q,\alpha) \in (0,\infty)$ is a constant.  Uniqueness of $\psi$ follows by again using Lemma~\ref{hole filling lemma}, the triangle inequality, and the homogeneity of $\psi$.
\end{proof}

\section{Decay of fine excess of area-minimizing currents relative to homogeneous Dirichlet energy minimizing multi-valued functions}\label{sec:fine excess decay}

In this section, we prove our main excess decay lemma for locally area-minimizing rectifiable currents, Theorem~\ref{main excess decay thm}, which plays a key role in the proof of the uniqueness of tangent functions.  First in Lemma~\ref{main excess decay lemma1}, we prove the excess decay lemma assuming that $\varphi$ and $T$ satisfy Hypothesis~$(\star)$ and Hypothesis~$(\star\star)$.  Using an argument of~\cite{Wic14}, in Lemma~\ref{main excess decay lemma2} (and in Theorem~\ref{main excess decay thm}), we replace Hypothesis~$(\star\star)$ with the weaker hypothesis \eqref{main excess decay2 hyp} that the $T$ is much closer to $\varphi$ (in $L^2$) than $T$ is to the plane $P_0 = \mathbb{R}^n \times \{0\}$ provided we allow conclusion~(c) (the improvement of fine excess) and conclusion~(d) to hold at one of a fixed number of smaller scales $\vartheta_1,\vartheta_2,\ldots,\vartheta_{\lfloor q/q_0 \rfloor}$.  For the purposes of obtaining uniform decay estimates for $T$ relative to a tangent function in the next section, having a finite, fixed number of scales to choose from at each step of iterating Theorem~\ref{main excess decay thm} is as good as using a single scale at every step.

\begin{lemma}\label{main excess decay lemma1}
Let $\alpha = \ell_0/q_0 \neq 2$ where $\ell_0$ and $q_0$ are relatively prime positive integers with $q_0 \leq q$.  Let $\varsigma = \varsigma(n,m,q) \in (0,1)$ is as in Theorem~\ref{relharm separation3 thm}.  For each $\vartheta \in (0,1/4)$ there exists $\overline{\varepsilon} = \overline{\varepsilon}(n,m,q,\alpha,\vartheta) \in (0,1)$, $\overline{\beta} = \overline{\beta}(n,m,q,\alpha,\vartheta) \in (0,1)$, and $\overline{\delta} = \overline{\delta}(n,m,q,\alpha,\vartheta) \in (0,1/2)$ such that the following holds true.  Let $p \in \{1,2,\ldots,\lfloor q/q_0 \rfloor\}$, $\varphi \in \Phi_{\alpha,q,p}$ and $T$ satisfy Hypothesis~$(\star)$ and Hypothesis~$(\star\star)$ with $\overline{\varepsilon},\overline{\beta}$ in place of $\varepsilon_0,\beta_0$.  Then either 
\begin{enumerate}[itemsep=2mm,topsep=0mm]
	\item[(A)]  $\mathbf{C}_{\overline{\delta}}(0,y_0) \cap \{ X \in \mathcal{B}_q : \mathcal{N}_{T,{\rm Pl}}(X) \geq \alpha \} = \emptyset$ for some $y_0 \in B^{n-2}_{1/2}(0)$ or 
	
	\item[(B)]  there exists an orthogonal rotation $\Gamma$ of $\mathbb{R}^{n+m}$ which fixes vectors in $\{0\} \times \mathbb{R}^m$ (i.e.~$\Gamma(v) = v$ for all $v \in \{0\} \times \mathbb{R}^m$) and $\widetilde{\varphi} \in \Phi_{\alpha,q}$ such that: 
	\begin{enumerate}[itemsep=2mm,topsep=2mm]
		\item[{\rm (a)}]  $\|\Gamma - I\| \leq \overline{\kappa} E_{\varsigma/2,\overline{\varepsilon}}(T,\varphi,\mathbf{C}_1(0))/\widehat{E}(T,P_0,\mathbf{C}_1(0))$, where $I$ is the identity map on $\mathbb{R}^n$; 	
		\item[{\rm (b)}]  $\op{dist}_{\mathcal H}(\op{spt}\op{graph}\varphi \cap \mathbf{C}_1(0), \op{spt}\op{graph}\widetilde{\varphi} \cap \mathbf{C}_1(0)) 
			\leq \overline{C}_0 E_{\varsigma/2,\overline{\varepsilon}}(T,\varphi,\mathbf{C}_1(0))$; 
		\item[{\rm (c)}]  $Q_{\varsigma/2,\overline{\varepsilon}}(\Gamma_{\#} T,\widetilde{\varphi},\mathbf{C}_{\vartheta}(0)) \leq \overline{\nu} \vartheta^{\alpha-1+\mu} E_{\varsigma/2,\overline{\varepsilon}}(T,\varphi,\mathbf{C}_1(0))$; 
		\item[{\rm (d)}]  $\widehat{E}(T,P_0,\mathbf{C}_{\vartheta}(0)) \geq C_1 \vartheta^{\alpha-1} \|\varphi\|_{L^2(B_1(0))} - \overline{C}_2 Q_{\varsigma/2,\overline{\varepsilon}}(T,\varphi,\mathbf{C}_1(0))$.
	\end{enumerate}
	Here the constants $\overline{\kappa},\overline{\nu},\overline{C}_0,C_1,\overline{C}_2 \in (0,\infty)$ and $\mu \in (0,1)$ each depend only on $n,m,q,\alpha$.
\end{enumerate}
\end{lemma}

\begin{proof}
Fix $\vartheta \in (0,1/4)$.  For $\nu = 1,2,3,\ldots$ let $\varepsilon_{\nu},\beta_{\nu},\delta_{\nu} \in (0,1)$ with $\varepsilon_{\nu}\rightarrow 0$, $\beta_{\nu}\rightarrow 0$, and $\delta_{\nu}\rightarrow 0$ as $\nu\rightarrow\infty$ and let $p \in \{1,2,\ldots,\lfloor q/q_0 \rfloor\}$, $\varphi^{(\nu)} \in \Phi_{\alpha,q,p}$, and $T^{(\nu)}$ be an $n$-dimensional locally area-minimizing rectifiable current of $\mathbf{C}_1(0)$ such that Hypothesis~$(\star)$ and Hypothesis~$(\star\star)$ hold true with $\varepsilon_{\nu},\beta_{\nu},\varphi^{(\nu)},T^{(\nu)}$ in place of $\varepsilon_0,\beta_0,\varphi,T$ and conclusion~(A) of the present lemma does not hold true with $\delta_{\nu},T^{(\nu)}$ in place of $\overline{\delta},T$.  Note that after passing to a subsequence, we assume that $p$ is independent of $\nu$.  To prove the lemma, it suffices to show that for infinitely many $\nu$ there exists an orthogonal rotation $\Gamma_{\nu}$ of $\mathbb{R}^n$ and $\widetilde{\varphi}^{(\nu)} : \mathbb{R}^n \rightarrow \mathcal{A}_q(\mathbb{R}^m)$ such that conclusion~(B) of the present lemma holds true with $\Gamma_{\nu},\varphi^{(\nu)},\widetilde{\varphi}^{(\nu)} \circ \Gamma_{\nu}^{-1},T^{(\nu)}$ in place of $\Gamma,\varphi,\widetilde{\varphi},T$ and with fixed constants $\overline{C}_1 \in (0,\infty)$ depending only on $n$ and $m$ and $\overline{\kappa},\overline{\nu},\overline{C}_0,\overline{C}_2 \in (0,\infty)$ and $\mu \in (0,1)$ depending only on $n,m,q,\alpha$.

After passing to a subsequence, let $w = (w_j) \in \mathfrak{B}$ be the fine blow-up of $(T^{(\nu)})$ relative to $(\varphi^{(\nu)})$ (as in Subsection~\ref{sec:blowup subsec} and Definition~\ref{blowup class defn}).  Let $w \in \mathfrak{B}(m_j,\varphi^{(0)},\varphi^{(\infty)}_j)$, where $(m_j,\varphi^{(0)},\varphi^{(\infty)}_j) \in \mathfrak{D}$ is as in Remark~\ref{blowup on domain rmk}.  By Theorem~\ref{blowup decay thm}, there exists (a unique) $\psi = (\psi_j) \in \widetilde{\mathfrak L}(m_j,\varphi^{(0)},\varphi^{(\infty)}_j)$ such that \eqref{blowup decay concl} holds true and $\|\psi\|_{L^2(B_1(0))} \leq C$, where $\mu = \mu(n,m,q,\alpha) \in (0,1)$ and $C = C(n,m,q,\alpha) \in (0,\infty)$ are constants.  Let $\Gamma_{\nu}$ and $\widetilde{\varphi}^{(\nu)}$ be as in Lemma~\ref{integrability lemma} corresponding to $\psi$, $\varphi^{(\nu)}$, and excesses $E_{\nu} = E_{\varsigma/2,\varepsilon_{\nu}}(T^{(\nu)},\varphi^{(\nu)},\mathbf{C}_1(0))$ and $\widehat{E}_{\nu} = \widehat{E}(T^{(\nu)},P_0,\mathbf{C}_1(0))$.  Extend $\Gamma_{\nu}$ orthogonal rotation of $\mathbb{R}^{n+m}$ which fixes vectors in $\{0\} \times \mathbb{R}^m$.  By Lemma~\ref{integrability lemma} and $\|\psi\|_{L^2(B_1(0))} \leq C$, we immediately have that $\widetilde{\varphi}^{(\nu)} \circ \Gamma_{\nu}^{-1} \in \Phi_{\alpha,q}$ and conclusions~(B)(a) and (B)(b) hold true with $\Gamma_{\nu},\varphi^{(\nu)},\widetilde{\varphi}^{(\nu)} \circ \Gamma_{\nu}^{-1},T^{(\nu)}$ in place of $\Gamma,\varphi,\widetilde{\varphi},T$ and for constants $\overline{\kappa} = \overline{\kappa}(n,m,q,\alpha) \in (0,\infty)$ and $\overline{C}_0 = \overline{C}_0(n,m,q,\alpha) \in (0,\infty)$.  By Lemma~\ref{blowup norms lemma}, \eqref{blowup decay concl} with $\rho = \vartheta$, and $\|w\|_{L^2(B_{1/2}(0))}^2 \leq \omega_n (1/2)^{n+2}$, for all sufficiently large $\nu$ 
\begin{equation}\label{main excess decay1 eqn1}
	\vartheta^{-2\alpha+2} Q_{\varsigma/2,1}(T^{(\nu)},\widetilde{\varphi}^{(\nu)},\mathbf{C}_{\vartheta}(0)) 
	\leq 4C \vartheta^{2\mu} E_{\varsigma/2,\varepsilon_{\nu}}(T^{(\nu)},\varphi^{(\nu)},\mathbf{C}_1(0)) 
\end{equation}
(where $C$ is as in \eqref{blowup decay concl}).  Since conclusion~(A) does not hold true, there exists $Z_{\nu} \in \mathcal{B}_q \cap \mathbf{C}_{\delta_{\nu}}(0)$ with $\mathcal{N}_{T,{\rm Pl}}(Z_{\nu}) \geq \alpha$.  Thus by Remark~\ref{planar freq rmk} and Theorem~\ref{frequency monotonicity thm}(e), for all sufficiently large $\nu$ 
\begin{align}\label{main excess decay1 eqn2}
	\widehat{E}(T,P_0,\mathbf{C}_{\vartheta}(0)) \leq\,& C \widehat{E}(T,P_0,\mathbf{C}_{2\vartheta}(Z_{\nu})) 
	\leq C \vartheta^{\alpha-1} \widehat{E}(T,P_0,\mathbf{C}_{1/2}(Z_{\nu})) \\ \leq\,& C \vartheta^{\alpha-1} \widehat{E}(T,P_0,\mathbf{C}_1(0)) , \nonumber
\end{align}
where $C = C(n,m,q,\alpha) \in (0,\infty)$ are constants.  By \eqref{main excess decay1 eqn1} and \eqref{main excess decay1 eqn2}, conclusion~(B)(c) holds true with $\min\{\mu,\varsigma(\alpha-1)\},\Gamma_{\nu},\varphi^{(\nu)},\widetilde{\varphi}^{(\nu)} \circ \Gamma_{\nu}^{-1},T^{(\nu)}$ in place of $\mu,\Gamma,\varphi,\widetilde{\varphi},T$ and for some constant $\overline{\nu} = \overline{\nu}(n,m,q,\alpha) \in (0,\infty)$.  By Remark~\ref{graphical triangle rmk}(2) with $\rho = \vartheta$, Remark~\ref{graphical rmk}(1), conclusion~(B)(b), and conclusion~(B)(c) 
\begin{align*}
	\widehat{E}(T^{(\nu)},P_0,\mathbf{C}_{\vartheta}(0)) 
	\geq\,& C_1 \vartheta^{\alpha-1} \|\varphi^{(\nu)}\|_{L^2(B_1(0))} - E_{\varsigma/2,1}(T^{(\nu)},\varphi^{(\nu)},\mathbf{C}_{\vartheta}(0))
	\\ \geq\,& C_1 \vartheta^{\alpha-1} \|\varphi^{(\nu)}\|_{L^2(B_1(0))} 
		- \frac{1}{2} E_{\varsigma/2,1}(T^{(\nu)},\widetilde{\varphi}^{(\nu)},\mathbf{C}_{\vartheta}(0)) 
		\\&- C \op{dist}_{\mathcal H}(\op{spt}\op{graph}\varphi \cap \mathbf{C}_1(0), \op{spt}\op{graph}\widetilde{\varphi} \cap \mathbf{C}_1(0)) 
	\\ \geq\,& C_1 \vartheta^{\alpha-1} \|\varphi^{(\nu)}\|_{L^2(B_1(0))} 
		- \overline{C}_2 E_{\varsigma/2,\varepsilon_{\nu}}(T^{(\nu)},\varphi^{(\nu)},\mathbf{C}_1(0)) 
\end{align*}
where $C_1 = C_1(n,m,q,\alpha) > 0$, $C = C(n,m,q,\alpha) \in (0,\infty)$, and $\overline{C}_2 = \overline{C}_2(n,m,q,\alpha) \in (0,\infty)$ are constants.  
\end{proof}

\begin{lemma}\label{main excess decay lemma2}
Let $\alpha = \ell_0/q_0 \neq 2$ where $\ell_0$ and $q_0$ are relatively prime positive integers with $q_0 \leq q$.  Let $\varsigma = \varsigma(n,m,q) \in (0,1)$ is as in Theorem~\ref{relharm separation3 thm}.  For each $p \in \{1,2,\ldots,\lfloor q/q_0 \rfloor\}$ and $\vartheta_1,\vartheta_2,\ldots,\vartheta_p \in (0,1/4)$ there exists $\varepsilon^{(p)} = \varepsilon^{(p)}(n,m,q,\alpha,\vartheta_1,\ldots,\vartheta_p) \in (0,1)$, $\beta^{(p)} = \beta^{(p)}(n,m,q,\alpha,$ $\vartheta_1,\ldots,\vartheta_p) \in (0,1)$, and $\delta^{(p)} = \delta^{(p)}(n,m,q,\alpha,\vartheta_1,\ldots,\vartheta_p) \in (0,1/2)$ such that the following holds true.  Let $\varphi \in \Phi_{\alpha,q,p}$ and $T$ satisfy Hypothesis~$(\star)$ with $\varepsilon^{(p)}$ in place of $\varepsilon_0$ and 
\begin{equation}\label{main excess decay2 hyp}
	Q_{\varsigma/2,\varepsilon^{(p)}}(T, \varphi, \mathbf{C}_1(0)) \leq \beta^{(p)} \widehat{E}(T, P_0, \mathbf{C}_1(0)) .
\end{equation}
Then either 
\begin{enumerate}[itemsep=2mm,topsep=0mm]
	\item[(A)]  $\mathbf{C}_{\delta^{(p)}}(0,y_0) \cap \{ X \in \mathcal{B}_q : \alpha \leq \mathcal{N}_{T,{\rm Pl}}(X) \geq \alpha \} = \emptyset$ for some $y_0 \in B^{n-2}_{1/2}(0)$ or 
	
	\item[(B)]  there exists an orthogonal rotation $\Gamma$ of $\mathbb{R}^{n+m}$ which fixes vectors in $\{0\} \times \mathbb{R}^m$, $\widetilde{\varphi} \in \Phi_{\alpha,q}$, and $j \in \{1,2,\ldots,p\}$ such that: 
	\begin{enumerate}[itemsep=2mm,topsep=2mm]
		\item[(a)]  $\|\Gamma - I\| \leq \kappa^{(p)} Q_{\varsigma/2,\varepsilon^{(p)}}(T,\varphi,\mathbf{C}_1(0))/\widehat{E}(T,P_0,\mathbf{C}_1(0))$; 	
		\item[(b)]  $\op{dist}_{\mathcal H}(\op{spt}\op{graph}\varphi \cap \mathbf{C}_1(0), \op{spt}\op{graph}\widetilde{\varphi} \cap \mathbf{C}_1(0)) 
			\leq C^{(p)}_0 Q_{\varsigma/2,\varepsilon^{(p)}}(T,\varphi,\mathbf{C}_1(0))$; 			
		
			\item[(c)]  $Q_{\varsigma/2,\varepsilon^{(p)}}(\Gamma_{\#} T,\widetilde{\varphi},\mathbf{C}_{\vartheta_j}(0)) \leq \nu^{(p)}_j \vartheta_j^{\alpha-1+\mu} Q_{\varsigma/2,\varepsilon^{(p)}}(T,\varphi,\mathbf{C}_1(0))$. 
		\item[(d)]  $\widehat{E}(T,P_0,\mathbf{C}_{\vartheta_j}(0)) \geq q^{-(p-1)/2} C_1 \vartheta_j^{\alpha-1} \|\varphi\|_{L^2(B_1(0))} - C^{(p)}_2 Q_{\varsigma/2,\varepsilon^{(p)}}(T,\varphi,\mathbf{C}_1(0))$.
	\end{enumerate}
	Here the constants $\mu = \mu(n,m,q,\alpha) \in (0,1)$ and $C_1 = C_1(n,m,q,\alpha) \in (0,\infty)$ are as in Lemma~\ref{main excess decay lemma1}; $\kappa^{(p)},C^{(p)}_0,C^{(p)}_2 \in (0,\infty)$ are constants depending only on $n,m,q,\alpha$ if $p = 1$ and depending only on $n,m,q,\alpha,\vartheta_1,\vartheta_2,\ldots,\vartheta_{p-1}$ if $p > 1$; $\nu^{(p)}_j \in (0,\infty)$ is a constant depending only on $n,m,q,\alpha$ if $j = 1$ and depending only on $n,m,q,\alpha,\vartheta_1,\vartheta_2,\ldots,\vartheta_{j-1}$ if $j \in \{2,3,\ldots,p\}$.
\end{enumerate}
\end{lemma}

\begin{proof} 
We proceed by induction on $p$.  In the base case $p = 1$, set $\varepsilon^{(1)} = \overline{\varepsilon}(n,m,q,\alpha,\vartheta_1)$, $\beta^{(1)} = \overline{\beta}(n,m,q,\alpha,\vartheta_1)$, $\delta^{(1)} = \overline{\delta}(n,m,q,\alpha,\vartheta_1)$, $\kappa^{(1)} = \overline{\kappa}(n,m,q,\alpha)$, $C_0^{(1)} = \overline{C}_0(n,m,q,\alpha)$, $C_2^{(1)} = \overline{C}_2(n,m,q,\alpha)$, and $\nu_1 = \overline{\nu}(n,m,q,\alpha)$, where $\overline{\varepsilon}, \overline{\beta}, \overline{\delta}, \overline{\kappa}, \overline{C}_0, \overline{C}_2,\overline{\nu}$ are as in Lemma~\ref{main excess decay lemma1} (with $\vartheta = \vartheta_1$).  Now the conclusion of the present lemma follows directly from Lemma~\ref{main excess decay lemma1} with $\vartheta = \vartheta_1$.  

Suppose that $p_0 \in \{2,3,\ldots,\lfloor q/q_0 \rfloor\}$ such that Lemma~\ref{main excess decay lemma2} holds true whenever $p \in \{1,2,\ldots,p_0-1\}$.  Let $\vartheta_1,\vartheta_2,\ldots,\vartheta_{p_0} \in (0,1/4)$.  Let $\overline{\varepsilon}, \overline{\beta}, \overline{\delta}, \overline{\kappa}, \overline{C}_0, \overline{C}_2, \overline{\nu}$ be as in Lemma~\ref{main excess decay lemma1} with $\vartheta = \vartheta_1$.  By induction, for each $p \in \{1,2,\ldots,p_0-1\}$ we obtain parameters $\varepsilon^{(p)},\beta^{(p)},\delta^{(p)},\kappa^{(p)},C_0^{(p)},C_2^{(p)}$ be as in Lemma~\ref{main excess decay lemma2} with $\vartheta_2,\vartheta_3,\ldots,\vartheta_{p+1}$ in place of $\vartheta_1,\vartheta_2,\ldots,\vartheta_p$.  We claim that Lemma~\ref{main excess decay lemma2} holds true with $p = p_0$ and with the constants 
\begin{align*}
	\varepsilon^{(p_0)} =\,& \min\{\overline{\varepsilon}\} \cup \{\varepsilon^{(p)} : p = 1,\ldots,p_0-1 \} , \\
	\beta^{(p_0)} =\,& \min\{\beta_0\} \cup \left\{ 2^{-p_0} \overline{\beta}^{p_0} \beta^{(p)} : p = 1,\ldots,p_0-1 \right\} , \\
	\delta^{(p_0)} =\,& \min\{\overline{\delta}\} \cup \{\delta^{(p)} : p = 1,\ldots,p_0-1 \} , \\
	\kappa^{(p_0)} =\,& \max\{\overline{\kappa}\} \cup \{ 2^{p_0} \overline{\beta}^{-p_0} \kappa^{(p)} : p = 1,\ldots,p_0-1 \} , \\
	C_0^{(p_0)} =\,& \max\{\overline{C}_0\} \cup \{ \overline{C} + 2^{p_0} \overline{\beta}^{-p_0} (\overline{C} + C_0^{(p)}) : p = 1,\ldots,p_0-1 \} , \\
	C_2^{(p_0)} =\,& \max\{\overline{C}_2\} \cup \{ 2\overline{C} C_1 \sqrt{\omega_n} 
		+ 2^{p_0} \overline{\beta}^{-p_0} (2\overline{C} C_1 \sqrt{\omega_n} + C_2^{(p)}) : p = 1,\ldots,p_0-1 \} , \\
	\nu_1^{(p_0)} =\,& \overline{\nu} \quad\text{and}\quad 
	\nu_j^{(p_0)} = \max\{ 2^{2p_0} \overline{\beta}^{-2p_0} \nu_{j-1}^{(p)} : p = j-1,\ldots,p_0-1 \} \text{ if } 2 \leq j \leq p_0 , 
\end{align*}
where $\beta_0 = \beta_0(n,m,q,\alpha,p_0,1/2,1/16)$ is as in Theorem~\ref{graphical_thm} (with $\gamma = 1/2$ and $\tau = 1/16$) and $\overline{C} = \overline{C}(n,m,q,\alpha)$ be as in \eqref{main excess decay2 eqn5} below.  To see this, let $\varphi \in \Phi_{\alpha,q,p_0}$ and $T$ be an $n$-dimensional locally area-minimizing rectifiable current of $\mathbf{C}_1(0)$ satisfying Hypothesis~$(\star)$ with $\varepsilon^{(p_0)}$ in place of $\varepsilon_0$ and 
\begin{equation}\label{main excess decay2 eqn1}
	Q_{\varsigma/2,\varepsilon^{(p)}}(T,\varphi,\mathbf{C}_1(0)) \leq \beta^{(p_0)} \widehat{E}(T, P_0, \mathbf{C}_1(0)) .
\end{equation}
If additionally Hypothesis~$(\star\star)$ holds true with $\overline{\beta}$ in place of $\beta_0$, i.e.
\begin{equation*}
	Q_{\varsigma/2,\varepsilon^{(p)}}(T,\varphi,\mathbf{C}_1(0)) 
		< \overline{\beta} \inf_{\varphi' \in \bigcup_{p'=1}^{p-1} \Phi_{\alpha,q,p'}} Q_{\varsigma/2,\varepsilon^{(p)}}(T,\varphi',\mathbf{C}_1(0))
\end{equation*}
then the conclusion of the present lemma follows directly from Lemma~\ref{main excess decay lemma1} with $\vartheta = \vartheta_1$.  If on the other hand 
\begin{equation}\label{main excess decay2 eqn2}
	Q_{\varsigma/2,\varepsilon^{(p)}}(T,\varphi,\mathbf{C}_1(0)) 
		\geq \overline{\beta} \inf_{\varphi' \in \bigcup_{p'=1}^{p-1} \Phi_{\alpha,q,p'}} Q_{\varsigma/2,\varepsilon^{(p)}}(T,\varphi',\mathbf{C}_1(0)) ,
\end{equation}
then arguing as in Remark~\ref{tilde varphi rmk}(1), we can choose an integer $\overline{p} \in \{0,1,2,\ldots,p_0-1\}$ and $\overline{\varphi} \in \Phi_{\alpha,q,\overline{p}}$ such that 
\begin{equation}\label{main excess decay2 eqn3}
	Q_{\varsigma/2,\varepsilon^{(p)}}(T,\overline{\varphi},\mathbf{C}_1(0)) 
		\leq 2^{p_0} \overline{\beta}^{1-p_0} \inf_{\varphi' \in \bigcup_{p'=0}^{p_0-1} \Phi_{\alpha,q,p'}} Q_{\varsigma/2,\varepsilon^{(p)}}(T,\varphi',\mathbf{C}_1(0)) 
\end{equation}
and either $\overline{p} = 0$ or $\overline{p} > 0$ and 
\begin{equation}\label{main excess decay2 eqn4} 
	Q_{\varsigma/2,\varepsilon^{(p)}}(T,\overline{\varphi},\mathbf{C}_1(0)) 
		\leq \overline{\beta} \inf_{\varphi' \in \bigcup_{p'=0}^{\overline{p}-1} \Phi_{\alpha,q,p'}} Q_{\varsigma/2,\varepsilon^{(p)}}(T,\varphi',\mathbf{C}_1(0)) .
\end{equation}
By Remark~\ref{graphical triangle rmk}(1) with $\rho = 1$ and with $\overline{\varphi},\varphi$ in place of $\varphi,\varphi^*$ 
\begin{align}\label{main excess decay2 eqn5} 
	&\op{dist}_{\mathcal H}(\op{spt}\op{graph}\varphi \cap \mathbf{C}_1(0), \op{spt}\op{graph}\overline{\varphi} \cap \mathbf{C}_1(0)) 
	\\ \leq\,& \overline{C} Q_{\varsigma/2,\varepsilon^{(p)}}(T,\varphi,\mathbf{C}_1(0)) + \overline{C} Q_{\varsigma/2,\varepsilon^{(p)}}(T,\overline{\varphi},\mathbf{C}_1(0)) 
		\nonumber 
\end{align}
for some constant $\overline{C} = \overline{C}(n,m,q,\alpha) \in (0,\infty)$.  Let $\varphi = \sum_{j=1}^J m_j \varphi_j$ with distinct components $\varphi_j : \mathbb{R}^n \rightarrow \mathcal{A}_{q_j}(\mathbb{R}^m)$ with integer multiplicity $m_j$ and let $\overline{\varphi} = \sum_{j=1}^{\overline{J}} \overline{m}_j \overline{\varphi}_j$ with distinct components $\overline{\varphi}_j : \mathbb{R}^n \rightarrow \mathcal{A}_{\overline{q}_j}(\mathbb{R}^m)$ with integer multiplicity $\overline{m}_j$ (as in Definition~\ref{Phi-p_defn}).  Let $\varphi(x,y) = \sum_{j=1}^J \sum_{l=1}^{q_j} m_j \llbracket \varphi_{j,l}(x,y) \rrbracket$ for each $(x,y) \in \mathbb{R}^m$ (as in Definition~\ref{Phi-p_defn} and \eqref{varphi_localized}) and similarly let $\overline{\varphi}(x,y) = \sum_{j=1}^{\overline{J}} \sum_{l=1}^{\overline{q}_j} \overline{m}_j \llbracket \overline{\varphi}_{j,l}(x,y) \rrbracket$ for each $(x,y) \in \mathbb{R}^m$, where $\overline{\varphi}_{j,l}(x,y) \in \mathbb{R}^m$.  For each $(x,y) \in \mathbb{R}^m$, $j \in \{1,2,\ldots,J\}$, and $l \in \{1,2,\ldots,q_j\}$ there exists $j' \in \{1,2,\ldots,\overline{J}\}$ and $l' \in \{1,2,\ldots,\overline{q}_{j'}\}$ such that $|\varphi_{j,l}(x,y) - \overline{\varphi}_{j',l'}(x,y)| \leq 2 \op{dist}((x,y,\varphi_{j,l}(x,y)),\op{spt}\op{graph}\overline{\varphi})$.  Thus by the triangle inequality 
\begin{align*}
	|\varphi_{j,l}(x,y)| \leq\,& |\overline{\varphi}_{j',l'}(x,y)| + 2 \op{dist}((x,y,\varphi_{j,l}(x,y)),\op{spt}\op{graph}\overline{\varphi})
	\\ \leq\,& |\overline{\varphi}(x,y)| + 2 \op{dist}((x,y,\varphi_{j,l}(x,y)),\op{spt}\op{graph}\overline{\varphi})
\end{align*}
and thus 
\begin{equation}\label{main excess decay2 eqn6} 
	|\varphi(x,y)| \leq \sqrt{q} \,|\overline{\varphi}(x,y)| 
		+ 2 \sqrt{q} \op{dist}_{\mathcal H}(\op{spt}\op{graph}\varphi \cap \mathbf{C}_1(0), \op{spt}\op{graph}\overline{\varphi} \cap \mathbf{C}_1(0)) .
\end{equation}
Integrating \eqref{main excess decay2 eqn6} over $(x,y) \in B_1(0)$ 
\begin{equation}\label{main excess decay2 eqn7} 
	\|\varphi\|_{L^2(B_1(0))} \leq \sqrt{q} \|\overline{\varphi}\|_{L^2(B_1(0))} + 2 \sqrt{q \omega_n} 
		\op{dist}_{\mathcal H}(\op{spt}\op{graph}\varphi \cap \mathbf{C}_1(0), \op{spt}\op{graph}\overline{\varphi} \cap \mathbf{C}_1(0)) . 
\end{equation}
By \eqref{main excess decay2 eqn2} and \eqref{main excess decay2 eqn3}
\begin{equation}\label{main excess decay2 eqn8} 
	Q_{\varsigma/2,\varepsilon^{(p)}}(T,\overline{\varphi},\mathbf{C}_1(0)) 
	\leq 2^{p_0} \overline{\beta}^{-p_0} Q_{\varsigma/2,\varepsilon^{(p)}}(T,\varphi,\mathbf{C}_1(0)) . 
\end{equation}
In particular, by \eqref{main excess decay2 eqn8} and \eqref{main excess decay2 eqn1} 
\begin{equation}\label{main excess decay2 eqn9} 
	Q_{\varsigma/2,\varepsilon^{(p)}}(T,\overline{\varphi},\mathbf{C}_1(0)) 
		\leq 2^{p_0} \overline{\beta}^{-p_0} \beta^{(p_0)} \widehat{E}(T, P_0, \mathbf{C}_1(0)). 
	 \end{equation}
Since $\sqrt{q} \,2^{p_0} \overline{\beta}^{-p_0} \beta^{(p_0)} < 1$, we must have that $\overline{p} \geq 1$; otherwise, $\overline{p} = 0$ implies that $\overline{\varphi} = q \llbracket 0 \rrbracket$ contradicting \eqref{main excess decay2 eqn9}.  Hence \eqref{main excess decay2 eqn9} gives us that $Q_{\varsigma/2,\varepsilon^{(p)}}(T,\overline{\varphi},\mathbf{C}_1(0)) < \beta^{(\overline{p})} \widehat{E}(T, P_0, \mathbf{C}_1(0))$ and thus we can apply the induction hypothesis to deduce that the conclusion of Lemma~\ref{main excess decay lemma2} holds true with $p = \overline{p}$ and $\overline{\varphi}$ in place of $\varphi$.  It follows using \eqref{main excess decay2 eqn5}, \eqref{main excess decay2 eqn7}, and \eqref{main excess decay2 eqn8} that the conclusion of Lemma~\ref{main excess decay lemma2} holds true.
\end{proof}

\begin{theorem}\label{main excess decay thm}
Let $\alpha = \ell_0/q_0 \neq 2$ where $\ell_0$ and $q_0$ are relatively prime positive integers with $q_0 \leq q$.  Let $\varsigma = \varsigma(n,m,q) \in (0,1)$ is as in Theorem~\ref{relharm separation3 thm}.  For each $\vartheta_1,\vartheta_2,\ldots,\vartheta_{\lfloor q/q_0 \rfloor} \in (0,1/4)$ there exists $\varepsilon = \varepsilon(n,m,q,\alpha,\vartheta_1,\ldots,\vartheta_{\lfloor q/q_0 \rfloor}) \in (0,1)$, $\beta = \beta(n,m,q,\alpha,\vartheta_1,\ldots,\vartheta_{\lfloor q/q_0 \rfloor}) \in (0,1)$, and $\delta = \delta(n,m,q,\alpha,\vartheta_1,\ldots,\vartheta_{\lfloor q/q_0 \rfloor}) \in (0,1/2)$ such that the following holds true.  Let $\varphi \in \Phi_{\alpha,q}$ and $T$ satisfy Hypothesis~$(\star)$ with $\varepsilon$ in place of $\varepsilon_0$ and 
\begin{equation*}
	Q_{\varsigma/2,\varepsilon}(T, \varphi, \mathbf{C}_1(0)) \leq \beta \widehat{E}(T, P_0, \mathbf{C}_1(0)) .
\end{equation*}
Then either 
\begin{enumerate}[itemsep=2mm,topsep=0mm]
	\item[(A)]  $\mathbf{C}_{\delta}(0,y_0) \cap \{ X \in \mathcal{B}_q : \mathcal{N}_{T,{\rm Pl}}(X) \geq \alpha \} = \emptyset$ for some $y_0 \in B^{n-2}_{1/2}(0)$ or 
	
	\item[(B)]  there exists an orthogonal rotation $\Gamma$ of $\mathbb{R}^{n+m}$ which fixes vectors in $\{0\} \times \mathbb{R}^m$, $\widetilde{\varphi} \in \Phi_{\alpha,q}$, and $j \in \{1,2,\ldots,\lfloor q/q_0 \rfloor\}$ such that: 
	\begin{enumerate}[itemsep=2mm,topsep=2mm]
		\item[(a)]  $\|\Gamma - I\| \leq \kappa Q_{\varsigma/2,\varepsilon}(T,\varphi,\mathbf{C}_1(0))/\widehat{E}(T,P_0,\mathbf{C}_1(0))$; 	
		\item[(b)]  $\op{dist}_{\mathcal H}(\op{spt}\op{graph}\varphi \cap \mathbf{C}_1(0), \op{spt}\op{graph}\widetilde{\varphi} \cap \mathbf{C}_1(0)) 
			\leq C_0 Q_{\varsigma/2,\varepsilon}(T,\varphi,\mathbf{C}_1(0))^2$; 
		\item[(c)]  $Q_{\varsigma/2,\varepsilon}(\Gamma_{\#} T,\widetilde{\varphi},\mathbf{C}_{\vartheta_j}(0)) \leq \nu_j \vartheta_j^{\alpha-1+\mu} Q_{\varsigma/2,\varepsilon}(T,\varphi,\mathbf{C}_1(0))$. 
		\item[(d)]  $\widehat{E}(T,P_0,\mathbf{C}_{\vartheta_j}(0)) \geq C_1 \vartheta_j^{\alpha-1} \|\varphi\|_{L^2(B_1(0))} - C_2 Q_{\varsigma/2,\varepsilon}(T,\varphi,\mathbf{C}_1(0))$.
	\end{enumerate}
	Here the constants $\mu = \mu(n,m,q,\alpha) \in (0,1)$ and $C_1 = C_1(n,m) \in (0,\infty)$ are as in Lemma~\ref{main excess decay lemma1}; $\kappa,C_0,C_2 \in (0,\infty)$ are constants depending only on $n,m,q,\alpha$ if $q_0 > q/2$ and depending only on $n,m,q,\alpha,\vartheta_1,\vartheta_2,\ldots,\vartheta_{\lfloor q/q_0 \rfloor - 1}$ if $q_0 \leq q/2$; $\nu_j \in (0,\infty)$ is a constant depending only on $n,m,q,\alpha$ if $j = 1$ and depending only on $n,m,q,\alpha,\vartheta_1,\vartheta_2,\ldots,\vartheta_{j-1}$ if $j \in \{2,3,\ldots,\lfloor q/q_0 \rfloor\}$.
\end{enumerate}
\end{theorem}

\begin{proof}
Set 
\begin{align*}
	\varepsilon =\,& \min\{\varepsilon^{(p)}(n,m,q,\alpha,\vartheta_1,\ldots,\vartheta_p) : p = 1,2,\ldots,\lfloor q/q_0 \rfloor\} , \\
	\beta =\,& \min\{\beta^{(p)}(n,m,q,\alpha,\vartheta_1,\ldots,\vartheta_p) : p = 1,2,\ldots,\lfloor q/q_0 \rfloor\} , \\
	\delta =\,& \min\{\delta^{(p)}(n,m,q,\alpha,\vartheta_1,\ldots,\vartheta_p) : p = 1,2,\ldots,\lfloor q/q_0 \rfloor\} 
\end{align*}
where $\varepsilon^{(p)},\beta^{(p)},\delta^{(p)}$ are as in Lemma~\ref{main excess decay lemma2}.  Set 
\begin{align*}
	\kappa =\,& \max\{\kappa^{(1)}(n,m,q,\alpha)\} \cup \{\kappa^{(p)}(n,m,q,\alpha,\vartheta_1,\ldots,\vartheta_{p-1}) : p = 2,3,\ldots,\lfloor q/q_0 \rfloor\} , \\
	C_0 =\,& \max\{C_0^{(1)}(n,m,q,\alpha)\} \cup \{C_0^{(p)}(n,m,q,\alpha,\vartheta_1,\ldots,\vartheta_{p-1}) : p = 2,3,\ldots,\lfloor q/q_0 \rfloor\} , \\
	C_2 =\,& \max\{C_2^{(1)}(n,m,q,\alpha)\} \cup \{C_2^{(p)}(n,m,q,\alpha,\vartheta_1,\ldots,\vartheta_{p-1}) : p = 2,3,\ldots,\lfloor q/q_0 \rfloor\} 
\end{align*}
if $q_0 < q$ and $\kappa = \kappa^{(1)}(n,m,q,\alpha)$, $C_0 = C_0^{(1)}(n,m,q,\alpha)$, and $C_2 = C_2^{(1)}(n,m,q,\alpha)$ if $q_0 = q$, where $\kappa^{(p)},C_0^{(p)},C_2^{(p)}$ are as in Lemma~\ref{main excess decay lemma2}.  Set 
\begin{equation*}
	\nu_j(n,m,q,\alpha,\vartheta_1,\ldots,\vartheta_{j-1}) = \max\{\nu_j^{(p)}(n,m,q,\alpha,\vartheta_1,\ldots,\vartheta_{j-1}) : p = j,j+1,\ldots,\lfloor q/q_0 \rfloor\} 
\end{equation*}
if $j \in \{2,\ldots,\lfloor q/q_0 \rfloor\}$ and $\nu_1(n,m,q,\alpha) = \max\{\nu_1^{(p)}(n,m,q,\alpha) : p = 1,2,3,\ldots,\lfloor q/q_0 \rfloor\}$, where $\nu_j^{(p)}$ is as in Lemma~\ref{main excess decay lemma2}.  With this choice of parameters, the conclusion of the present theorem follows directly from Lemma~\ref{main excess decay lemma2}.
\end{proof}

\section{$\mathcal{H}^{n-2}$-a.e.~uniqueness of coarse blow-ups and $(n-2)$-rectifiability of the set of branch points of planar frequency $\neq 2$}\label{sec:proof of main results}

\begin{theorem}\label{SandT thm} 
Let $\alpha = \ell_0/q_0 \neq 2$ where $\ell_0$ and $q_0$ are relatively prime positive integers with $q_0 \leq q$.  There exists $\varepsilon = \varepsilon(n,m,q,\alpha) \in (0,1)$, $\beta = \beta(n,m,q,\alpha) \in (0,1)$, $\eta = \eta(n,m,q,\alpha) \in (0,1)$, and $\overline{\mu} = \overline{\mu}(n,m,q,\alpha) \in (0,1)$ such that the following holds true.  Let $\varphi \in \Phi_{\alpha,q}$ be non-zero and locally Dirichlet energy minimizing and let $T$ be an $n$-dimensional locally area-minimizing rectifiable current of $\mathbf{C}_1(0)$ such that 
\begin{gather}
	\label{SandT hyp1} (\partial T) \llcorner \mathbf{C}_1(0) = 0, \quad\sup_{X \in \op{spt} T} \op{dist}(X,P_0) < \varepsilon, \quad
	\pi_{\#} T = q \llbracket B_1(0) \rrbracket,  , \\
	\label{SandT hyp2} Q(T,\varphi,\mathbf{C}_1(0)) \leq \beta \,\widehat{E}(T,P_0,\mathbf{C}_1(0)) ,
\end{gather}
where $Q(T,\varphi,\mathbf{C}_1(0))$ is as in \eqref{Q defn} (with $Z = 0$, $\rho = 1$, and $P = \mathbb{R}^n \times \{0\}$).  Let $\mathcal{B}_q$ be as in Definition~\ref{introduction B defn} (with $\mu = 1/q$).  Then 
\begin{equation*}
	\{Z \in \mathcal{B}_q : \mathcal{N}_{T,{\rm Pl}}(Z) \geq \alpha \} = \Sigma \cup \Gamma
\end{equation*}
where $\Sigma \subseteq L$ for some properly embedded $(n-2)$-dimensional $C^{1,\overline{\mu}}$-submanifold $L$ of $\mathbf{C}_{1/2}(0)$ with $\mathcal{H}^{n-2}(L) \leq 2\omega_{n-2} (1/2)^{n-2}$ and $\Gamma \subseteq \bigcup_{j=1}^{\infty} \mathbf{C}_{\rho_j}(\xi_j,\zeta_j)$ for a countable family of cylinders $\{\mathbf{C}_{\rho_j}(\xi_j,\zeta_j)\}$ with $\rho_j < 1/2$ and $\sum_{j=1}^{\infty} \rho_j^{n-2} \leq 1-\eta$.  Moreover, for $\mathcal{H}^{n-2}$-a.e.~$Z \in \Sigma$, there exists a unique non-zero, homogeneous degree $\alpha$, locally Dirichlet energy minimizing function $\varphi^{(Z)} : P_Z \rightarrow \mathcal{A}_q(P_Z^{\perp})$ such that $\omega_n^{1/2} \|\varphi^{(Z)}\|_{L^2(B_1(0,P_Z))}^{-1} \,\varphi_Z$ is the unique tangent function to $T$ at $Z$, $\dim\op{spine}\varphi^{(Z)} = n-2$, 
\begin{equation}\label{SandT concl}
	\frac{1}{\rho^n} \int_{\mathbf{C}_{\rho}(Z,P_Z)} \op{dist}^2(X,Z+\op{spt}\op{graph}\varphi^{(Z)}) \,d\|T\|(X) \leq C_Z \rho^{2\alpha+2\overline{\mu}} 
\end{equation}
for all $\rho \in (0,1/4]$, and for some $\rho_Z \in (0,1/4]$ and a function $\tau_Z : (0,\rho_Z] \rightarrow (0,1/16]$ with $\lim_{\rho\rightarrow 0^+} \tau_Z(\rho) = 0$ (both depending on $Z$) 
\begin{align}\label{SandT concl2}
	&\op{dist}_{\mathcal H}(\op{spt} T \cap \mathbf{C}_{\rho}(Z,P_Z) \cap \{r_Z > \tau_Z(\rho)\rho\}, \\&\hspace{1in} (Z + \op{spt}\op{graph}\varphi^{(Z)}) \cap \mathbf{C}_{\rho}(Z,P_Z) \cap \{r_Z > \tau_Z(\rho)\rho\}) \leq C_Z \rho^{\alpha+\overline{\mu}} \nonumber 
\end{align}
for all $\rho \in (0,\rho_Z]$, where $P_Z$ is the unique tangent plane to $T$ at $Z$, $r_Z(X) = \op{dist}(\pi_{P_Z}(X),\op{spine}\varphi^{(Z)})$, and $C_Z \in (0,\infty)$ is a constant (depending on $Z$).
\end{theorem}

\begin{proof} 
\textit{Step~1:  Choice of parameters and notation.}  For $j = 1,2,\ldots,\lfloor q/q_0 \rfloor$ choose $\vartheta_j = \vartheta_j(n,m,q,\alpha) \in (0,1/4)$ such that $\vartheta_{j+1} \leq \vartheta_j/8$ if $1 \leq j < \lfloor q/q_0 \rfloor$ and $\nu_j \vartheta_j^{\mu/2} < 1/2$ for all $j$, where $\mu = \mu(n,m,q,\alpha)$, $\nu_1 = \nu_1(n,m,q,\alpha)$, and $\nu_j = \nu_1(n,m,q,\alpha,\vartheta_1,\ldots,\vartheta_{j-1})$ for $2 \leq j \leq \lfloor q/q_0 \rfloor$ are as in Theorem~\ref{main excess decay thm}.  Let $\varepsilon = \varepsilon(n,m,q,\alpha) \in (0,1)$, $\beta = \beta(n,m,q,\alpha) \in (0,1)$, and $\delta = \delta(n,m,q,\alpha) \in (0,1/2)$ be suitable small constants to be later determined.  Set 
\begin{equation*}%\label{SandT eqn1} 
	\mathcal{B}^{\star}_q = \{Z \in \mathcal{B}_q : \mathcal{N}_{T,{\rm Pl}}(X) \geq \alpha \} .
\end{equation*}
For each $Z \in \mathcal{B}^{\star}_q$, let $P_Z$ be the unique tangent plane to $T$ at $Z$. 

Recall that from Remark~\ref{planar freq rmk} that we can apply Lemma~\ref{uniform branch set lemma} with $\mu = 1/q$.  By Theorem~\ref{frequency monotonicity thm}(e), $\widehat{E}(T,P_Z,\mathbf{C}_{\rho}(Z,P_Z)) \leq C \rho^{\alpha-1} \widehat{E}(T,P_0,\mathbf{C}_1(0))$ for all $\rho \in (0,1/4]$.  Given $Y,Z \in \mathcal{B}^{\star}_q$, let $\sigma = |Y-Z|$.  If $\sigma \leq 1/16$, by the triangle inequality and~\cite[Lemma~3.4]{KrumWica}
\begin{align*}
	\|\pi_{P_Y} - \pi_{P_Z}\|^2 
	\leq\,& \frac{2}{\omega_n \sigma^n} \int_{G_n(\mathbf{B}_{\sigma}(Y))} \|\pi_S - \pi_{P_Y}\|^2 \,d|T|(X,S) 
		\\&+ \frac{2}{\omega_n \sigma^n} \int_{G_n(\mathbf{B}_{2\sigma}(Z))} \|\pi_S - \pi_{P_Z}\|^2 \,d|T|(X,S)
	\\ \leq\,& C \widehat{E}(T,P_Y,\mathbf{B}_{2\sigma}(Y))^2 + C \widehat{E}(T,P_Z,\mathbf{B}_{4\sigma}(Z))^2
	\\ \leq\,& C |Y-Z|^{2(\alpha-1)} \widehat{E}(T,P_0,\mathbf{C}_1(0))^2 ,
\end{align*}
where $C = C(n,m,q,\alpha) \in (0,\infty)$ are constants.  By \eqref{uniform branch set concl3}, $\|\pi_{P_Y} - \pi_{P_Z}\| \leq C |X-Y|^{\alpha-1} \widehat{E}(T,P_0,\mathbf{C}_1(0))$ if $|Y-Z| \geq 1/16$, where $C = C(n,m,q,\alpha) \in (0,\infty)$ is a constant.  Let $\widetilde{\Gamma}_Z$ be an orthogonal rotation of $\mathbb{R}^{n+m}$ such that 
\begin{gather}\label{SandT eqn2} 
	\widetilde{\Gamma}_Z(P_Z) = P_0, \quad 
	\|\widetilde{\Gamma}_Z - I\| \leq C \op{dist}_{\mathcal H}(P_Z \cap \mathbf{B}_1(0),P_0 \cap \mathbf{B}_1(0)) , \\
	\|\widetilde{\Gamma}_Y - \widetilde{\Gamma}_Z\| 
    \leq C |Y-Z|^{\alpha-1} \widehat{E}(T,P_0,\mathbf{C}_1(0)) \nonumber 
\end{gather}
for all $Y,Z \in \mathcal{B}^{\star}_q$, where $P_0 = \mathbb{R}^n \times \{0\}$ and $C = C(n,m,q,\alpha) \in (0,\infty)$ is a constant.  

Suppose that $\mathcal{B}^{\star}_q \cap \mathbf{C}_{\vartheta_{\lfloor q/q_0 \rfloor} \delta/8}(0,y_0) = \emptyset$ for some $y_0 \in B^{n-2}_{1/2}(0)$.  Then by Lemma~\ref{nonconcentration3 lemma}, $\mathcal{B}^{\star}_q \subset B^2_{\tau(\varepsilon)}(0) \times \mathbb{R}^{n-2} \times \mathbb{R}^m$ for some $\tau(\varepsilon,\beta) \rightarrow 0$ as $\varepsilon \downarrow 0$ and $\beta \downarrow 0$.  Hence the conclusion of the theorem trivially holds true with $\Sigma = \emptyset$ and $\Gamma = \mathcal{B}^{\star}_q$.  Thus we may assume that $\mathcal{B}^{\star}_q \cap \mathbf{C}_{\vartheta_{\lfloor q/q_0 \rfloor} \delta/8}(0,y_0) \neq \emptyset$ for all $y_0 \in B^{n-2}_{1/2}(0)$.  By Lemma~\ref{nonconcentration3 lemma}, this implies that 
\begin{equation}\label{SandT eqn3} 
	\widetilde{\Gamma}_Z(\mathcal{B}^{\star}_q - Z) \cap \mathbf{C}_{\delta\rho}(0,y_0) \neq \emptyset
\end{equation} 
for all $Z \in \mathcal{B}^{\star}_q \cap \mathbf{C}_{1/4}(0)$, $\rho \in (\vartheta_{\lfloor q/q_0 \rfloor}/4,1/4]$, and $y_0 \in B^{n-2}_{\rho/2}(0)$.

\noindent\textit{Step~2:  Change of coordinates at a branch point $Z$.}  Let $\varsigma = \varsigma(n,m,q) \in (0,1)$ be as in Theorem~\ref{relharm separation3 thm} and let $\beta_0 = \beta_0(n,m,q,\alpha) \in (0,1)$ be a small constant to be later determined.  We claim that there is an integer $\widetilde{p} \in \{1,2,\ldots,p\}$ and $\widetilde{\varphi} \in \Phi_{\alpha,q,\widetilde{p}}$ such that 
\begin{gather}
	\label{SandT eqn4} Q_{\varsigma/2,\varepsilon}(T,\widetilde{\varphi},\mathbf{C}_1(0)) 
		\leq 2^p \beta_0^{-p} Q_{\varsigma/2,\varepsilon}(T,\varphi,\mathbf{C}_1(0)) , \\
	\label{SandT eqn5} Q_{\varsigma/2,\varepsilon}(T,\widetilde{\varphi},\mathbf{C}_1(0)) 
		\leq \beta_0 \inf_{\varphi' \in \bigcup_{p'=0}^{\widetilde{p}-1} \Phi_{\alpha,q,p'}} Q_{\varsigma/2,\varepsilon}(T,\varphi',\mathbf{C}_1(0)) .
\end{gather}
If \eqref{SandT eqn5} holds true with $p,\varphi$ in place of $\widetilde{p},\widetilde{\varphi}$, then we can take $\widetilde{p} = p$ and $\widetilde{\varphi} = \varphi$.  Otherwise, we can argue as in Remark~\ref{tilde varphi rmk}(1) that there exists $\widetilde{p} \in \{0,1,2,\ldots,p-1\}$ and $\widetilde{\varphi} \in \Phi_{\alpha,q,\widetilde{p}}$ such that \eqref{SandT eqn5} holds true and 
\begin{equation}\label{SandT eqn6} 
	Q_{\varsigma/2,\varepsilon}(T,\widetilde{\varphi},\mathbf{C}_1(0)) 
	\leq 2^p \beta_0^{1-p} \inf_{\varphi' \in \bigcup_{p'=0}^{p-1} \Phi_{\alpha,q,p'}} Q_{\varsigma/2,\varepsilon}(T,\varphi',\mathbf{C}_1(0)) .
\end{equation}
It follows from \eqref{SandT eqn6} that \eqref{SandT eqn4} holds true.  By \eqref{SandT eqn4} and \eqref{SandT hyp2} 
\begin{equation*} 
	Q_{\varsigma/2,\varepsilon}(T,\widetilde{\varphi},\mathbf{C}_1(0)) \leq 2^p \beta_0^{-p} \beta \,\widehat{E}(T,P_0,\mathbf{C}_1(0)) 
\end{equation*}
so that provided $2^p \beta_0^{-p} \beta < 1$ we must have that $\widetilde{p} \geq 1$.

Consider an arbitrary point $Z = (\xi,\zeta,\chi) \in \mathcal{B}^{\star}_q \cap \mathbf{C}_{1/4}(0)$.  Set $T^{(Z)} = \widetilde{\Gamma}_{Z\#} \eta_{Z,1\#} T$.  We claim that 
\begin{gather}
	\label{SandT eqn7} \tfrac{1}{C} \,\widehat{E}(T,P_0,\mathbf{C}_1(0)) \leq \widehat{E}(T^{(Z)},P_0,\mathbf{C}_{1/4}(0)) 
		\leq C \widehat{E}(T,P_0,\mathbf{C}_1(0)) , \\
	\label{SandT eqn8} Q_{\varsigma/2,\varepsilon}(T^{(Z)},\varphi,\mathbf{C}_{1/4}(0)) \leq C Q_{\varsigma/2,\varepsilon}(T,\varphi,\mathbf{C}_1(0)) , 
\end{gather}
where $C = C(n,m,q,\alpha) \in (1,\infty)$ is a constant.  Assuming that $\beta_0$ is sufficiently small, in light of \eqref{SandT eqn5} we can argue as we did in Remark~\ref{graphical triangle rmk}(2) (with $\widetilde{\varphi}$ in place of $\varphi$) using Theorem~\ref{graphical_thm}(B)(C) to obtain 
\begin{equation*}
	\|\widetilde{\varphi}\|_{L^2(B_1(0))} \leq \widehat{E}(T,P_0,\mathbf{C}_{1/8}(0,\zeta)) + C E_{\varsigma/2,\varepsilon}(T,\widetilde{\varphi},\mathbf{C}_1(0)) ,
\end{equation*} 
where $C = C(n,m,q,\alpha) \in (0,\infty)$ are constants.  Note that by Theorem~\ref{nonconcentration thm}(A) we have $\op{spt} T \cap \mathbf{C}_{1/8}(0,\zeta) \subset \mathbf{C}_{1/4}(Z,P_Z)$, and that $\op{dist}_{\mathcal H}((Z + P_0) \cap \mathbf{C}_1(0), P_0 \cap \mathbf{C}_1(0)) \leq C |\chi|$.  Hence by triangle inequality, Theorem~\ref{nonconcentration thm}(A), \eqref{SandT eqn4}, and \eqref{SandT hyp2} 
\begin{align}\label{SandT eqn9}
	\|\widetilde{\varphi}\|_{L^2(B_1(0))}
	\leq\,& C \widehat{E}(T,P_Z,\mathbf{C}_{1/4}(Z,P_Z)) + C |\chi| + C \op{dist}_{\mathcal H}(P_Z \cap \mathbf{B}_1(0), P_0 \cap \mathbf{B}_1(0)) 
	\\&+ C E_{\varsigma/2,\varepsilon}(T,\widetilde{\varphi},\mathbf{C}_1(0)) \nonumber
	\\ \leq\,& C \widehat{E}(T^{(Z)},P_0,\mathbf{C}_{1/4}(0)) + C E_{\varsigma/2,\varepsilon}(T,\widetilde{\varphi},\mathbf{C}_1(0)) \nonumber
	\\ \leq\,& C \widehat{E}(T^{(Z)},P_0,\mathbf{C}_{1/4}(0)) + C\beta_0^{-p} \beta \widehat{E}(T,P_0,\mathbf{C}_1(0)) , \nonumber 
\end{align}
where $C = C(n,m,q,\alpha) \in (0,\infty)$ are constants.  By Remark~\ref{graphical rmk}(2), $\widehat{E}(T,P_0,\mathbf{C}_1(0)) \leq C \|\widetilde{\varphi}\|_{L^2(B_1(0))}$.  Thus provided $\beta$ is sufficiently small, \eqref{SandT eqn9} implies the first inequality in \eqref{SandT eqn7}.  The second inequality in \eqref{SandT eqn7} follows by using the triangle inequality, \eqref{uniform branch set concl1}, and \eqref{uniform branch set concl3} to obtain
\begin{equation*}
	\widehat{E}(T,P_Z,\mathbf{C}_{1/4}(Z,P_Z)) \leq C \widehat{E}(T,P_0,\mathbf{C}_1(0)) + C  |\vec P_Z - \vec P_0| 
		\leq C \widehat{E}(T,P_0,\mathbf{C}_1(0)) ,
\end{equation*}
where $C = C(n,m,q) \in (0,\infty)$ are constants.  To obtain \eqref{SandT eqn8}, in light of \eqref{SandT eqn5} we can argue as we did to obtain \eqref{noncon triangle eqn9} and \eqref{noncon triangle eqn10} (with $\widetilde{\varphi}$ in place of $\varphi$) using the triangle inequality, \eqref{noncon triangle3}, Theorem~\ref{graphical_thm}(B), and Theorem~\ref{nonconcentration thm}(A) to get 
\begin{align}\label{SandT eqn10}
	Q_{\varsigma/2,\varepsilon}(T^{(Z)},\widetilde{\varphi},\mathbf{C}_{1/4}(0)) 
	\leq\,& C Q_{\varsigma/2,\varepsilon}(T,\widetilde{\varphi},\mathbf{C}_1(0)) + C |\chi| + C \widehat{E} |\xi| 
		\\&+ C \op{dist}_{\mathcal H}(P_Z \cap \mathbf{B}_1(0), P_0 \cap \mathbf{B}_1(0)) \nonumber
	\\ \leq\,& C Q_{\varsigma/2,\varepsilon}(T,\widetilde{\varphi},\mathbf{C}_1(0)) , \nonumber 
\end{align}
where $\widehat{E} = \widehat{E}(T,P_0,\mathbf{C}_1(0))$ and $C = C(n,m,q,\alpha) \in (0,\infty)$ are constants.  Again noting \eqref{SandT eqn5}, we can argue as we did in Remark~\ref{graphical triangle rmk}(1)  (with $\widetilde{\varphi},\varphi$ in place of $\varphi,\varphi^*$) using Theorem~\ref{graphical_thm}(B)(C) to obtain 
\begin{align}\label{SandT eqn11}
	&\op{dist}_{\mathcal H}(\op{spt} \op{graph} \varphi \cap \mathbf{C}_1(0), \op{spt} \op{graph} \widetilde{\varphi} \cap \mathbf{C}_1(0)) \\ \leq\,& C Q_{\varsigma/2,\varepsilon}(T,\varphi,\mathbf{C}_1(0)) + C E_{\varsigma/2,\varepsilon}(T,\widetilde{\varphi},\mathbf{C}_1(0)) , \nonumber 
\end{align} 
where $C = C(n,m,q,\alpha) \in (0,\infty)$ are constants.  By the triangle inequality, \eqref{SandT eqn10}, \eqref{SandT eqn11}, and \eqref{SandT eqn4} 
\begin{align*}
    Q_{\varsigma/2}(T^{(Z)},\varphi,\mathbf{C}_{1/4}(0)) 
    \leq\,& Q_{\varsigma/2}(T^{(Z)},\widetilde{\varphi},\mathbf{C}_1(0)) \\&+ C \op{dist}_{\mathcal H}(\op{spt} \op{graph} \varphi \cap \mathbf{C}_1(0), \op{spt} \op{graph} \widetilde{\varphi} \cap \mathbf{C}_1(0)) \nonumber 
    \\ \leq\,& C Q_{\varsigma/2,\varepsilon}(T,\varphi,\mathbf{C}_1(0)) + C Q_{\varsigma/2,\varepsilon}(T,\widetilde{\varphi},\mathbf{C}_1(0)) \nonumber
    \\ \leq\,& C \beta_0^{-p} \,Q_{\varsigma/2,\varepsilon}(T,\varphi,\mathbf{C}_1(0)) , \nonumber
\end{align*}
where $C = C(n,m,q,\alpha) \in (0,\infty)$ are constants, and thus \eqref{SandT eqn8} holds true.

\noindent\textit{Step~3:  Iteratively applying the excess decay lemma.}  Let $Z \in \mathcal{B}^{\star}_q \cap \mathbf{C}_{1/4}(0)$.  Let $\Gamma_0 = I$ be the identity map on $\mathbb{R}^{n+m}$, $\varphi_0 = \varphi$, and $s_0 = 1/4$.  We claim that by inductively apply Theorem~\ref{main excess decay thm}, we can show the following: \textit{there exists $k_{\star}(Z) \in \{1,2,3,\ldots\} \cup \{\infty\}$ such that}
\begin{enumerate}[itemsep=2mm,topsep=0mm]
	\item[(i)]  \textit{for each positive integer $k \leq k_{\star}(Z)$ there exists an orthogonal rotation $\Gamma_k$ of $\mathbb{R}^{n+m}$ which fixes vectors in $\{0\} \times \mathbb{R}^m$, $\varphi_k \in \Phi_{\alpha,q}$, $j_k \in \{1,2,\ldots,\lfloor q/q_0 \rfloor\}$, and a radius $s_k$ given by $s_k = \vartheta_{j_k} s_{k-1}$ such that } 
	\begin{align}
		\label{SandT eqn12} &\|\Gamma_k - \Gamma_{k-1}\| 
			\leq \frac{C Q_{\varsigma/2,\varepsilon}(\Gamma_{k-1\#} T^{(Z)},\varphi_{k-1},\mathbf{C}_{s_{k-1}}(0))}{
			\widehat{E}(T^{(Z)},P_0,\mathbf{C}_{s_{k-1}}(0))} , \\
		\label{SandT eqn13} &s_{k-1}^{\alpha-1} \op{dist}_{\mathcal H}(\op{spt}\op{graph}\varphi_{k-1} \cap \mathbf{C}_1(0), 
			\op{spt}\op{graph}\varphi_k \cap \mathbf{C}_1(0)) 
			\\&\hspace{51.2mm} \leq C Q_{\varsigma/2,\varepsilon}(\Gamma_{k-1\#} T^{(Z)},\varphi_{k-1},\mathbf{C}_{s_{k-1}}(0)) , \nonumber \\
		\label{SandT eqn14} &Q_{\varsigma/2,\varepsilon}(\Gamma_{k\#} T^{(Z)},\varphi_k,\mathbf{C}_{s_k}(0))
			\leq \vartheta_{j_k}^{\alpha-1+\mu/2} Q_{\varsigma/2,\varepsilon}(\Gamma_{k-1\#} T^{(Z)},\varphi_{k-1},\mathbf{C}_{s_{k-1}}(0)) , \\
		\label{SandT eqn15} &\widehat{E}(T^{(Z)},P_0,\mathbf{C}_{s_k}(0)) \geq c s_k^{\alpha-1} \|\varphi_{k-1}\|_{L^2(B_1(0))} 
			- C Q_{\varsigma/2,\varepsilon}(\Gamma_{k-1\#} T^{(Z)},\varphi_{k-1},\mathbf{C}_{s_{k-1}}(0)) ,
	\end{align}
	\textit{where $\mu = \mu(n,m,q,\alpha)$ is as in Theorem~\ref{main excess decay thm} and $c = c(n,m) > 0$ and $C = C(n,m,q,\alpha) \in (0,\infty)$ are constants;}
	
	\item[(ii)]  \textit{either $k_{\star}(Z) = \infty$ or $k_{\star}(Z) < \infty$ and for some $y_0 \in B^{n-2}_{s_{k_{\star}(Z)}/2}(0)$.}
	\begin{equation}\label{SandT eqn16}
     \widetilde{\Gamma}_Z(\mathcal{B}^{\star}_q - Z) \cap \mathbf{C}_{\delta\vartheta_{\lfloor q/q_0 \rfloor} s_{k_{\star}(Z)}/2}(0,y_0) = \emptyset .
	\end{equation} 
\end{enumerate}
When $k_{\star}(Z) < \infty$, we shall take $k_{\star}(Z)$ to be the smallest positive integer such that \eqref{SandT eqn16} holds true for some $y_0 \in B^{n-2}_{s_{k_{\star}(Z)}/2}(0)$.

To prove the claim, we shall proceed by induction on $k$.  First let's consider the base case $k = 1$, and verify that the hypotheses of Theorem~\ref{main excess decay thm} hold true with $\eta_{0,1/4\#} T^{(Z)}$ and $4^{1-\alpha} \varphi$ in place of $T$ and $\varphi$.  
  By the triangle inequality, \eqref{SandT hyp1}, and \eqref{uniform branch set concl3}
\begin{equation*}
	\sup_{X \in \op{spt} T^{(Z)} \cap \mathbf{C}_{3/8}(0)} \op{dist}(X,P_0) 
	\leq \sup_{X \in \op{spt} T} \op{dist}(X,P_0) + C |\vec P_Z - \vec P_0| < C \varepsilon ,
\end{equation*}
where $C = C(n,m,q) \in (0,\infty)$ is a constant.  Hence provided $\varepsilon = \varepsilon(n,m,q,\alpha)$ is sufficiently small, by \eqref{SandT hyp1}, the constancy theorem, and \eqref{uniform branch set concl2}
\begin{equation}\label{SandT eqn17}
	(\partial T^{(Z)}) \llcorner \mathbf{C}_{3/8}(0) = 0, \quad (\pi_{\#} T^{(Z)}) \llcorner \mathbf{C}_{3/8}(0) = q \llbracket B_{3/8}(0) \rrbracket ,
\end{equation}
and by~\cite[Lemma~3.7]{KrumWica} and \eqref{oriented tilt and mass}, $\|T^{(Z)}\|(\mathbf{C}_{1/4}(0)) \leq (q+1/2) \,\omega_n (1/4)^n$.  Therefore, Hypothesis~$(\star)$ holds true with $C\varepsilon, \eta_{0,1/4\#}T^{(Z)}$ in place of $\varepsilon_0,T$, where $C = C(n,m,q,\alpha) \in (0,\infty)$ is a constant.  By \eqref{SandT eqn7}, \eqref{SandT eqn8}, and \eqref{SandT hyp2} 
\begin{equation*}%\label{SandT eqn18}
	Q_{\varsigma/2,\varepsilon}(T^{(Z)},\varphi,\mathbf{C}_{1/4}(0)) \leq C\beta \widehat{E}(T^{(Z)},P_0,\mathbf{C}_{1/4}(0))  , 
\end{equation*}
for some constant $C = C(n,m,q,\alpha) \in (0,\infty)$.  Therefore, the hypotheses of Theorem~\ref{main excess decay thm} hold true with $C\varepsilon,C\beta,\eta_{0,1/4\#}T^{(Z)},4^{1-\alpha} \varphi$ in place of $\varepsilon,\beta,T,\varphi$, where $C = C(n,m,q,\alpha) \in (0,\infty)$ is a constant.  The claim in the case $k = 1$ follows by applying Theorem~\ref{main excess decay thm} with $\eta_{0,1/4\#}T^{(Z)}$ and $4^{1-\alpha}\varphi$ in place of $T$ and $\varphi$.  Note that by \eqref{SandT eqn3}, alternative~(A) of Theorem~\ref{main excess decay thm} does not hold true with $\eta_{0,1/4\#}T^{(Z)}$ and $4^{1-\alpha}\varphi$ in place of $T$ and $\varphi$.  Thus by alternative~(B) of Theorem~\ref{main excess decay thm}, conclusion~(i) of the claim holds true with $k = 1$.

Now suppose that for some positive integer $k$, for all $k' \in \{1,2,\ldots,k\}$ conclusion~(i) of the claim holds true with $k'$ in place of $k$.  We want to show that the hypotheses of Theorem~\ref{main excess decay thm} hold true with $\eta_{0,s_k\#} \Gamma_{k\#} T^{(Z)}$ and $s_k^{\alpha-1} \varphi_k$ in place of $T$ and $\varphi$.  Recall that \eqref{SandT eqn17} holds true.  Since $\mathcal{N}_{T^{(Z)},{\rm Pl}}(0) = \mathcal{N}_{T,{\rm Pl}}(Z) \geq \alpha$, by Theorem~\ref{frequency monotonicity thm}(e) and \eqref{SandT eqn7} 
\begin{equation}\label{SandT eqn19}
	\widehat{E}(T^{(Z)},P_0,\mathbf{C}_{4s_k}(0)) \leq C s_k^{\alpha-1} \widehat{E}(T,P_0,\mathbf{C}_1(0)) , 
\end{equation}
where $C = C(n,m,q,\alpha) \in (0,\infty)$ is a constant.  Hence by~\cite[Lemma~3.5]{KrumWica}
\begin{equation*}
	\sup_{X \in \op{spt} T^{(Z)} \cap \mathbf{C}_{2s_k}(0)} \op{dist}(X,P_0) \leq C s_k^{\alpha} \widehat{E}(T,P_0,\mathbf{C}_1(0)) < C \varepsilon ,  
\end{equation*}
where $C = C(n,m,q,\alpha) \in (0,\infty)$ are constants.  Hence by~\cite[Lemma~3.7]{KrumWica} and \eqref{oriented tilt and mass}, $\|T^{(Z)}\|(\mathbf{C}_{s_k}(0)) \leq (q+1/2) \,\omega_n s_k^n$.  Therefore, Hypothesis~$(\star)$ holds true with $C\varepsilon, \eta_{0,s_k\#}T^{(Z)}$ in place of $\varepsilon_0,T$, where $C = C(n,m,q,\alpha) \in (0,\infty)$ is a constant.  By iteratively applying \eqref{SandT eqn14} and also using \eqref{SandT eqn8} 
\begin{align}\label{SandT eqn20}
	Q_{\varsigma/2,\varepsilon}(\Gamma_{i\#}T^{(Z)},\varphi_i, \mathbf{C}_{s_i}(0))
	\leq\,& (4s_i)^{\alpha-1+\mu/2} Q_{\varsigma/2,\varepsilon}(T^{(Z)},\varphi, \mathbf{C}_{1/4}(0)) 
	\\ \leq\,& C s_i^{\alpha-1+\mu/2} Q_{\varsigma/2,\varepsilon}(T,\varphi,\mathbf{C}_1(0)) , \nonumber 
\end{align}
for all integers $0 \leq i \leq k$, where $C = C(n,m,q,\alpha) \in (0,\infty)$ is a constant.  By \eqref{SandT eqn13} and \eqref{SandT eqn20} 
\begin{equation}\label{SandT eqn21}
	\op{dist}_{\mathcal H}(\op{spt}\op{graph}\varphi_i \cap \mathbf{C}_1(0), 
		\op{spt}\op{graph}\varphi_{i+1} \cap \mathbf{C}_1(0)) \leq C s_i^{\mu/2} Q_{\varsigma/2,\varepsilon}(T,\varphi,\mathbf{C}_1(0)) 
	 \end{equation}
for all integers $0 \leq i < k$, where $C = C(n,m,q,\alpha) \in (0,\infty)$ is a constant.  Note that since $\vartheta_l < 1/4$ for all $l \in \{1,2,\ldots,\lfloor q/q_0 \rfloor\}$, $s_i < 4^{-(i-j)} s_j$ for all integers $0 \leq j \leq i \leq k$ and in particular $s_i < 4^{-i-1}$.  Thus by the triangle inequality and \eqref{SandT eqn21}, for all $j \in \{0,1,2,\ldots,k\}$ 
\begin{align}\label{SandT eqn22}
	&\op{dist}_{\mathcal H}(\op{spt}\op{graph}\varphi_j \cap \mathbf{C}_1(0), \op{spt}\op{graph}\varphi_k \cap \mathbf{C}_1(0)) 
	\\ \leq\,& \sum_{i=j}^{k-1} \op{dist}_{\mathcal H}(\op{spt}\op{graph}\varphi_i \cap \mathbf{C}_1(0), 
		\op{spt}\op{graph}\varphi_{i+1} \cap \mathbf{C}_1(0)) \nonumber 
	\\ \leq\,& C \sum_{i=j}^{k-1} 4^{-(i-j)\,\mu/2} s_j^{\mu/2} Q_{\varsigma/2,\varepsilon}(T,\varphi,\mathbf{C}_1(0))  
	\leq C s_j^{\mu/2} Q_{\varsigma/2,\varepsilon}(T,\varphi,\mathbf{C}_1(0)) , \nonumber 
\end{align}
where $C = C(n,m,q,\alpha) \in (0,\infty)$ are constants.  
By \eqref{SandT eqn15}, Remark~\ref{graphical rmk}(2), \eqref{SandT eqn20}, \eqref{SandT eqn22}, and \eqref{SandT hyp2} 
\begin{align}\label{SandT eqn24}
	\widehat{E}(T^{(Z)},P_0,\mathbf{C}_{s_{k-1}}(0)) 
	\geq\,& c s_k^{\alpha-1} \|\varphi_{k-1}\|_{L^2(B_1(0))} 
		- C Q_{\varsigma/2,\varepsilon}(\Gamma_{k-1\#} T^{(Z)},\varphi_{k-1},\mathbf{C}_{s_{k-1}}(0))
	\\ \geq\,& c s_k^{\alpha-1} \|\varphi\|_{L^2(B_1(0))} 
		- C Q_{\varsigma/2,\varepsilon}(\Gamma_{k-1\#} T^{(Z)},\varphi_{k-1},\mathbf{C}_{s_{k-1}}(0)) 
		\nonumber \\&- C s_k^{\alpha-1} \op{dist}_{\mathcal H}(\op{spt}\op{graph}\varphi \cap \mathbf{C}_1(0), 
		\op{spt}\op{graph}\varphi_k \cap \mathbf{C}_1(0)) \nonumber 
	\\ \geq\,& c s_k^{\alpha-1} \widehat{E}(T,P_0,\mathbf{C}_1(0)) - C s_k^{\alpha-1} Q_{\varsigma/2,\varepsilon}(T,\varphi,\mathbf{C}_1(0)) \nonumber 
	\\ \geq\,& c s_k^{\alpha-1} \widehat{E}(T,P_0,\mathbf{C}_1(0)) \nonumber 
\end{align}
where $c = c(n,m,q,\alpha) > 0$ and $C = C(n,m,q,\alpha) \in (0,\infty)$ are constants.  Combining \eqref{SandT eqn19} and \eqref{SandT eqn24}
\begin{equation}\label{SandT eqn25}
	\tfrac{1}{C} \,s_k^{\alpha-1} \widehat{E}(T,P_0,\mathbf{C}_1(0)) \leq \widehat{E}(T^{(Z)},P_0,\mathbf{C}_{s_k}(0)) 
		\leq C s_k^{\alpha-1} \widehat{E}(T,P_0,\mathbf{C}_1(0)) 
\end{equation}
for some constant $C = C(n,m,q,\alpha) \in (1,\infty)$.  By \eqref{SandT eqn20}, \eqref{SandT hyp2}, and \eqref{SandT eqn25} 
\begin{equation}\label{SandT eqn26}
	Q_{\varsigma/2,\varepsilon}(\Gamma_{k\#}T^{(Z)},\varphi_k, \mathbf{C}_{s_k}(0)) \leq C\beta \widehat{E}(T^{(Z)},P_0,\mathbf{C}_{s_k}(0))  , 
\end{equation}
for some constant $C = C(n,m,q,\alpha) \in (0,\infty)$.  In particular, the hypotheses of Theorem~\ref{main excess decay thm} hold true with $C\varepsilon,C\beta,\eta_{0,s_k\#}\Gamma_{k\#}T^{(Z)},s_k^{1-\alpha} \varphi_k$ in place of $\varepsilon,\beta,T,\varphi$, where $C = C(n,m,q,\alpha) \in (0,\infty)$ is a constant.  By \eqref{SandT eqn12}, \eqref{SandT eqn20}, and \eqref{SandT eqn25} 
\begin{equation}\label{SandT eqn27}
	\|\Gamma_i - \Gamma_{i-1}\| \leq \frac{C s_i^{\mu/2} Q_{\varsigma/2,\varepsilon}(T,\varphi,\mathbf{C}_1(0))}{\widehat{E}(T,P_0,\mathbf{C}_1(0))} 
\end{equation}
for all $0 \leq i < k$, where $C = C(n,m,q,\alpha) \in (0,\infty)$ is a constant.  By the triangle inequality, \eqref{SandT eqn27}, and $s_i \leq 4^{-(i-j)} s_j$ for all $i > j$  
\begin{align}\label{SandT eqn28}
	\|\Gamma_j - \Gamma_k\| 
	\leq \sum_{i=j+1}^k \|\Gamma_i - \Gamma_{i-1}\| 
	\leq\,& \sum_{j=i+1}^k \frac{C 4^{-(i-j)\,\mu/2} s_j^{\mu/2} Q_{\varsigma/2,\varepsilon}(T,\varphi,\mathbf{C}_1(0))}{\widehat{E}(T,P_0,\mathbf{C}_1(0))} 
	\\ \leq\,& \frac{C s_j^{\mu/2} Q_{\varsigma/2,\varepsilon}(T,\varphi,\mathbf{C}_1(0))}{\widehat{E}(T,P_0,\mathbf{C}_1(0))} \nonumber 
\end{align}
for all $0 \leq j \leq k$, where $C = C(n,m,q,\alpha) \in (0,\infty)$ are constants.   
The claim now follows by applying Theorem~\ref{main excess decay thm} with $\eta_{0,s_k\#}\Gamma_{k\#}T^{(Z)}$ and $s_k^{\alpha-1} \varphi_k$ in place of $T$ and $\varphi$.  In particular, if alternative~(A) of Theorem~\ref{main excess decay thm} holds true with $\eta_{0,s_k\#}\Gamma_{k\#}T^{(Z)}$ and $s_k^{\alpha-1} \varphi_k$ in place of $T$ and $\varphi$, i.e.~$\Gamma_k \widetilde{\Gamma}_Z(\mathcal{B}^{\star}_q - Z) \cap \mathbf{C}_{\delta s_k/2}(0,y_0) = \emptyset$ for some $(0,y_0) \in B^{n-2}_{s_k/2}(0)$, then since $\|\Gamma_k - I\| \leq C \beta$ (by \eqref{SandT eqn28} with $j = 0$), conclusion~(ii) of the claim holds true with $k_{\star}(Z) = k$ and we stop the iterative procedure.  If instead alternative~(B) of Theorem~\ref{main excess decay thm} holds true with $\eta_{0,s_k\#}\Gamma_{k\#}T^{(Z)}$ and $s_k^{\alpha-1} \varphi_k$ in place of $T$ and $\varphi$, then we deduce that conclusion~(i) of the claim holds true.  To use Theorem~\ref{main excess decay thm}(B)(b) to derive \eqref{SandT eqn13}, we note that since $\varphi_{k-1}$ and $\varphi_k$ are Lipschitz $q$-valued functions with small gradients, $s_{k-1}^{\alpha-1} \op{dist}_{\mathcal H}(\op{spt}\op{graph}\varphi_{k-1} \cap \mathbf{C}_1(0), \op{spt}\op{graph}\varphi_k \cap \mathbf{C}_1(0)) \leq 2 s_{k-1}^{\alpha-1} \|\mathcal{G}(\varphi_{k-1},\varphi_k)\|_{L^{\infty}(B_1(0))} \leq 2 \op{dist}_{\mathcal H}(\op{spt}\op{graph}(s_{k-1}^{\alpha-1} \varphi_{k-1}) \cap \mathbf{C}_1(0), \op{spt}\op{graph}(s_{k-1}^{\alpha-1} \varphi_k) \cap \mathbf{C}_1(0))$.

Having proved the claim, for each positive integer $k$ let 
\begin{equation*}
	\Upsilon_k = \{ Z \in Z \in \mathcal{B}^{\star}_q \cap \mathbf{C}_{1/4}(0) : k_{\star}(Z) = k \} 
\end{equation*}
and let $\Sigma = \Upsilon_{\infty} = \{ Z \in \mathcal{B}^{\star}_q \cap \mathbf{C}_{1/4}(0) : k_{\star}(Z) = \infty \}$.  

\noindent\textit{Step~4:  Uniqueness of tangent functions and submanifold structure of the set $\Sigma$ of branch points with no small gaps.}  Suppose that $Z \in \Sigma$.  For all $\rho \in (0,1/4]$, there exists an integer $k \geq 0$ such that $s_{k+1} < \rho \leq s_k$.  Since $\widetilde{\Gamma}_Z(\mathcal{B}^{\star}_q - Z) \cap \mathbf{C}_{\delta \vartheta_{\lfloor q/q_0 \rfloor} s_k/2}(0,y_0) \neq \emptyset$ for $y_0 \in B^{n-2}_{s_k/2}(0)$
\begin{equation}\label{SandT eqn30}
	\widetilde{\Gamma}_Z(\mathcal{B}^{\star}_q - Z) \cap \mathbf{C}_{\delta \rho/2}(0,y_0) \neq \emptyset 
		\text{ for all } y_0 \in B^{n-2}_{\rho/2}(0) . 
\end{equation}
By \eqref{SandT eqn28}, $(\Gamma_k)_{k=1}^{\infty}$ is a Cauchy sequence as linear maps on $\mathbb{R}^{n+m}$.  Therefore, there exists an orthogonal rotation $\overline{\Gamma}_Z$ of $\mathbb{R}^{n+m}$ which fixes vectors in $\{0\} \times \mathbb{R}^m$ such that $\|\Gamma_k - \overline{\Gamma}_Z\| \rightarrow 0$ as $k \rightarrow \infty$.  By letting $k\rightarrow\infty$ in \eqref{SandT eqn28} 
\begin{equation}\label{SandT eqn31}
	\|\Gamma_j - \overline{\Gamma}_Z\| \leq \frac{C s_j^{\mu/2} Q_{\varsigma/2,\varepsilon}(T,\varphi,\mathbf{C}_1(0))}{\widehat{E}(T,P_0,\mathbf{C}_1(0))}  
\end{equation}
for all integers $j \geq 0$, where $C = C(n,m,q,\alpha) \in (0,\infty)$ is a constant.  In particular, when $j = 0$ 
\begin{equation}\label{SandT eqn32}
	\|\overline{\Gamma}_Z - I\| \leq \frac{C Q_{\varsigma/2,\varepsilon}(T,\varphi,\mathbf{C}_1(0))}{\widehat{E}(T,P_0,\mathbf{C}_1(0))} \leq C\beta
\end{equation}
for some constant $C = C(n,m,q,\alpha) \in (0,\infty)$.  By \eqref{SandT eqn26}, we can use Remark~\ref{graphical rmk}(2) and \eqref{SandT eqn25} to obtain 
\begin{equation}\label{SandT eqn33}
	\tfrac{1}{C} \widehat{E}(T,P_0,\mathbf{C}_1(0)) \leq \|\varphi_k\|_{L^2(B_1(0))} \leq C \widehat{E}(T,P_0,\mathbf{C}_1(0)) 
\end{equation}
for all integers $k \geq 1$, where $C = C(n,m,q,\alpha) \in (1,\infty)$ is a constant.  Hence there exists a subsequence $(\varphi_{k_i})$ of $(\varphi_k)$ and $\overline{\varphi}^{(Z)} \in \Phi_{\alpha,q}$ such that $\varphi_{k_i} \rightarrow \overline{\varphi}^{(Z)}$ uniformly on compact subsets of $\mathbb{R}^n$ as $i \rightarrow\infty$ and 
\begin{equation}\label{SandT eqn34} 
	\tfrac{1}{C} \widehat{E}(T,P_0,\mathbf{C}_1(0)) \leq \|\overline{\varphi}^{(Z)}\|_{L^2(B_1(0))} \leq C \widehat{E}(T,P_0,\mathbf{C}_1(0)) , 
\end{equation}
where $C = C(n,m,q,\alpha) \in (1,\infty)$ is a constant.  By \eqref{SandT eqn22}, $(\op{spt}\varphi_k)_{k=1}^{\infty}$ is a Cauchy sequence of closed subsets of $\mathbf{C}_1(0)$.  Hence $\op{dist}_{\mathcal H}(\op{spt}\op{graph}\varphi_k \cap \mathbf{C}_1(0), \op{spt}\op{graph}\overline{\varphi}^{(Z)} \cap \mathbf{C}_1(0)) \rightarrow 0$ as $k \rightarrow \infty$.  By letting $k\rightarrow\infty$ in \eqref{SandT eqn22}
\begin{equation}\label{SandT eqn35}
	\op{dist}_{\mathcal H}(\op{spt}\op{graph}\varphi_j \cap \mathbf{C}_1(0), \op{spt}\op{graph}\overline{\varphi}^{(Z)} \cap \mathbf{C}_1(0)) 
	\leq C s_j^{\mu/2} Q_{\varsigma/2,\varepsilon}(T,\varphi,\mathbf{C}_1(0))
\end{equation}
for all integers $j \geq 0$, where $C = C(n,m,q,\alpha) \in (0,\infty)$ is a constant.  By \eqref{SandT eqn34}, \eqref{SandT eqn31}, and \eqref{SandT eqn35} 
\begin{align}\label{SandT eqn36}
	&\op{dist}_{\mathcal H}(\op{spt}\op{graph}(\varphi_j \circ \Gamma_j) \cap \mathbf{C}_1(0), 
		\op{spt}\op{graph}(\overline{\varphi}^{(Z)} \circ \overline{\Gamma}_Z) \cap \mathbf{C}_1(0)) 
	\\ \leq\,& \op{dist}_{\mathcal H}(\op{spt}\op{graph}\varphi_j \cap \mathbf{C}_1(0), \op{spt}\op{graph}\overline{\varphi}^{(Z)} \cap \mathbf{C}_1(0)) 
		+ C \widehat{E}(T,P_0,\mathbf{C}_1(0)) \,\|\Gamma_j - \overline{\Gamma}_Z\| \nonumber 
	\\ \leq\,& C s_j^{\mu/2} Q_{\varsigma/2,\varepsilon}(T,\varphi,\mathbf{C}_1(0)) , \nonumber
	\end{align}
for each integer $j \geq 0$, where $C = C(n,m,q,\alpha) \in (0,\infty)$ are constants.  By the triangle inequality, \eqref{SandT eqn20}, and \eqref{SandT eqn36}
\begin{align}\label{SandT eqn37}
	&Q_{\varsigma/2,\varepsilon}(\overline{\Gamma}_{Z\#}T^{(Z)},\overline{\varphi}^{(Z)}, \mathbf{C}_{s_k}(0)) 
	\leq Q_{\varsigma/2,\varepsilon}(\Gamma_{k\#} T^{(Z)},\varphi_k, \mathbf{C}_{s_k}(0)) 
		\\&\hspace{15mm} + C s_k^{\alpha-1} \op{dist}_{\mathcal H}(\op{spt}\op{graph}(\varphi_k \circ \Gamma_k) \cap \mathbf{C}_1(0), 
		\op{spt}\op{graph}(\overline{\varphi}^{(Z)} \circ \overline{\Gamma}_Z) \cap \mathbf{C}_1(0)) \nonumber
	\\ \leq\,& C s_k^{\alpha-1+\mu/2} Q_{\varsigma/2,\varepsilon}(T,\varphi,\mathbf{C}_1(0)) \nonumber 
\end{align}
for each integer $k \geq 0$, where $C = C(n,m,q,\alpha) \in (0,\infty)$ is a constant.  For each $\rho \in (0,1/4]$ there exists an integer $k \geq 0$ such that $s_{k+1} < \rho \leq s_k$ and thus \eqref{SandT eqn25} gives us 
\begin{equation}\label{SandT eqn38}
	\tfrac{1}{C} \,\rho^{\alpha-1} \widehat{E}(T,P_0,\mathbf{C}_1(0)) \leq \widehat{E}(T^{(Z)},P_0,\mathbf{C}_{\rho}(0)) 
		\leq C \rho^{\alpha-1} \widehat{E}(T,P_0,\mathbf{C}_1(0)) 
\end{equation}
and similarly \eqref{SandT eqn37} gives us 
\begin{equation}\label{SandT eqn39}
	E_{\varsigma/2,\varepsilon}(\overline{\Gamma}_{Z\#} T^{(Z)}, \overline{\varphi}^{(Z)},\mathbf{C}_{\rho}(0)) \leq C \rho^{\mu/2} Q_{\varsigma/2,\varepsilon}(T,\varphi,\mathbf{C}_1(0)) ,
\end{equation}
where $C = C(n,m,q,\alpha) \in (1,\infty)$ are constants.  In other words, setting $\varphi^{(Z)} = \widetilde{\Gamma}_Z^{-1} \circ \overline{\varphi}^{(Z)} \circ \overline{\Gamma}_Z \circ \widetilde{\Gamma}_Z$, we deduce that \eqref{SandT concl} holds true. 

It follows from \eqref{SandT concl} that each tangent function of $T$ at $Z$ is equal to a scalar multiple of $\varphi^{(Z)}$, except with possibly different multiplicities on the components.  (Using \eqref{SandT eqn37}, we can show each tangent function to $T$ at $Z$ has non-zero multiplicity on each component of $\varphi^{(Z)}$.)  We claim that there exists a choice of multiplicities for the components of $\varphi^{(Z)}$ such that $\omega_n^{1/2} \|\varphi^{(Z)}\|_{L^2(B_1(0))}^{-1} \,\varphi^{(Z)}$ is the unique tangent function of $T$ at $Z$.  Fix any $\tau \in (0,1/16)$ and let $\widetilde{\varepsilon} = \widetilde{\varepsilon}(n,m,q,\alpha,\tau) \in (0,\varepsilon]$ and $\widetilde{\beta} = \widetilde{\beta}(n,m,q,\alpha,\tau) \in (0,1)$ to be later determined.  By \eqref{SandT eqn17}, \eqref{SandT eqn38},~\cite[Lemma~3.5 and Lemma~3.7]{KrumWica}, and \eqref{oriented tilt and mass}, we know that Hypothesis~$(\star)$ holds true with $\widetilde{\varepsilon},\eta_{0,\rho\#} \overline{\Gamma}_{Z\#}T^{(Z)}$ in place of $\varepsilon_0,T$ for all sufficiently small $\rho > 0$.  Let $p \in \{1,2,\ldots,\lfloor q/q_0 \rfloor\}$ such that $\overline{\varphi}^{(Z)} \in \Phi_{\alpha,q,p}$.  We claim that for all sufficiently large $k$
\begin{equation}\label{SandT eqn40}
	Q_{\varsigma/2,\varepsilon}(\overline{\Gamma}_{Z\#}T^{(Z)},\overline{\varphi}^{(Z)},\mathbf{C}_{s_k}(0)) 
	\leq \widetilde{\beta} \inf_{\varphi' \in \bigcup_{p'=0}^{p-1} \Phi_{\alpha,q,p'}} 
		Q_{\varsigma/2,\varepsilon}(\overline{\Gamma}_{Z\#}T^{(Z)},\varphi',\mathbf{C}_{s_k}(0)) 
\end{equation}
so that Hypothesis~$(\star\star)$ hold true with $\widetilde{\beta},\eta_{0,s_k\#} \overline{\Gamma}_{Z\#}T^{(Z)},s_k^{\alpha-1} \overline{\varphi}^{(Z)}$ in place of $\beta_0,T,\varphi$.  Suppose to the contrary that there exists $k(\nu) \rightarrow \infty$ such that \eqref{SandT eqn40} is false with $k = k(\nu)$.  Arguing as in Remark~\ref{tilde varphi rmk}(1), we can choose $\widetilde{p}_{\nu} \in \{1,2,\ldots,p-1\}$ and $\widetilde{\varphi}^{(\nu)} \in \Phi_{\alpha,q,p_{\nu}}$ such that 
\begin{gather}
	\label{SandT eqn41} Q_{\varsigma/2,\varepsilon}(\overline{\Gamma}_{Z\#}T^{(Z)},\widetilde{\varphi}^{(\nu)},\mathbf{C}_{s_{k(\nu)}}(0)) \leq 2^p \widetilde{\beta}^{-p} Q_{\varsigma/2,\varepsilon}(\overline{\Gamma}_{Z\#}T^{(Z)},\overline{\varphi}^{(Z)},\mathbf{C}_{s_{k(\nu)}}(0))  , \\
	\label{SandT eqn42} Q_{\varsigma/2,\varepsilon}(\overline{\Gamma}_{Z\#}T^{(Z)},\widetilde{\varphi}^{(\nu)},\mathbf{C}_{s_{k(\nu)}}(0)) 
		\leq \widetilde{\beta} \inf_{\varphi' \in \bigcup_{p'=0}^{p_{\nu}-1} \Phi_{\alpha,q,p'}} 
		Q_{\varsigma/2,\varepsilon}(\overline{\Gamma}_{Z\#}T^{(Z)},\varphi',\mathbf{C}_{s_{k(\nu)}}(0)) .
\end{gather}
In light of \eqref{SandT eqn42}, we can apply Remark~\ref{graphical triangle rmk}(1) with $\eta_{0,s_{k(\nu)}\#} \overline{\Gamma}_{Z\#}T^{(Z)},s_{k(\nu)}^{\alpha-1} \overline{\varphi}^{(Z)},s_{k(\nu)}^{\alpha-1} \widetilde{\varphi}^{(\nu)}$ in place of $T,\varphi,\varphi^*$ to obtain 
\begin{align*} 
	&s_{k(\nu)}^{\alpha-1} \op{dist}_{\mathcal H}(\op{spt}\op{graph}\widetilde{\varphi}^{(\nu)} \cap \mathbf{C}_1(0), 
		\op{spt}\op{graph}\overline{\varphi}^{(Z)} \cap \mathbf{C}_1(0)) 
	\\ \leq\,& C Q_{\varsigma/2,\varepsilon}(\overline{\Gamma}_{Z\#}T^{(Z)},\widetilde{\varphi}^{(\nu)},\mathbf{C}_{s_{k(\nu)}}(0)) 
		+ C Q_{\varsigma/2,\varepsilon}(\overline{\Gamma}_{Z\#}T^{(Z)},\overline{\varphi}^{(Z)},\mathbf{C}_{s_{k(\nu)}}(0)) , \nonumber 
\end{align*}
where $C = C(n,m,q,\alpha) \in (0,\infty)$ is a constant.  By \eqref{SandT eqn41} and \eqref{SandT eqn37} 
\begin{align}\label{SandT eqn43} 
	&s_{k(\nu)}^{\alpha-1} \op{dist}_{\mathcal H}(\op{spt}\op{graph}\widetilde{\varphi}^{(\nu)} \cap \mathbf{C}_1(0), 
		\op{spt}\op{graph}\overline{\varphi}^{(Z)} \cap \mathbf{C}_1(0)) 
	\\ \leq\,& C \widetilde{\beta}^{-p} Q_{\varsigma/2,\varepsilon}(\overline{\Gamma}_{Z\#}T^{(Z)},\overline{\varphi}^{(Z)},\mathbf{C}_{s_{k(\nu)}}(0))  
	\leq C \widetilde{\beta}^{-p} s_{k(\nu)}^{\alpha-1+\mu/2} Q_{\varsigma/2,\varepsilon}(T,\varphi,\mathbf{C}_1(0)) , \nonumber 
\end{align}
where $C = C(n,m,q,\alpha) \in (0,\infty)$ are constants.  After passing to a subsequence, there exists $\widetilde{\varphi} \in \Phi_{\alpha,q}$ such that $\widetilde{\varphi}^{(\nu)} \rightarrow \widetilde{\varphi}$ uniformly on compact subsets of $\mathbb{R}^n$.  By dividing \eqref{SandT eqn43} by $s_{k(\nu)}^{\alpha-1}$ and letting $\nu\rightarrow\infty$, we deduce that the supports of the graphs of $\widetilde{\varphi}$ and $\overline{\varphi}^{(Z)}$ are equal, contradicting $\widetilde{\varphi}$ having fewer components than $\overline{\varphi}^{(Z)}$.  Therefore, \eqref{SandT eqn40} must hold true.

Let $\overline{\varphi}^{(Z)} = \sum_{j=1}^{\overline{J}} \overline{m}_j \overline{\varphi}^{(Z)}_j$ with components $\overline{\varphi}^{(Z)}_j : \mathbb{R}^n \rightarrow \mathcal{A}_{\overline{q}_j}(\mathbb{R}^m)$ with multiplicity $\overline{m}_j$ (as in Definition~\ref{Phi-p_defn}).  Noting that in the previous paragraph $\tau$ is arbitrary and can be taken to be $\leq \vartheta_{\lfloor q/q_0 \rfloor}/16$, by \eqref{SandT eqn37} and \eqref{SandT eqn40} we can apply Theorem~\ref{graphical_thm} with $\eta_{0,s_k\#} \overline{\Gamma}_{Z\#}T^{(Z)}$ and $s_k^{\alpha-1} \overline{\varphi}^{(Z)}$ in place of $T$ and $\varphi$, up to changing the multiplicities $\overline{m}_j$ of $\varphi^{(Z)}$, there exist $\rho_Z \in (0,1/4]$, a set 
\begin{equation*}
	\Omega_Z = \mathbf{C}_{\rho_Z}(Z,P_Z) \cap \{Z+X+Y : X \in P_Z,\, Y \in P_Z^{\perp}\text{, and } \op{dist}(X,\op{spine}\varphi^{(Z)}) > \tau_Z(|X|) \,|X| \} 
\end{equation*}
where $\tau_Z : (0,\rho_Z] \rightarrow (0,1/16]$ with $\lim_{\rho\rightarrow 0^+} \tau_Z(\rho) = 0$, and $n$-dimensional locally area-minimizing rectifiable currents $T^{(Z)}_j$ of $\Omega_Z$ such that 
\begin{gather*}
	T \llcorner \Omega_Z = \sum_{j=1}^{\overline{J}} T^{(Z)}_j , \quad\quad
	(\partial T^{(Z)}_j) \llcorner \Omega_Z = 0, \\
	\sup_{X \in \op{spt} T^{(Z)}_j \cap \Omega \cap \mathbf{C}_{\rho}(Z,P_Z)} \op{dist}(X, \op{spt} \op{graph} \varphi^{(Z)}_j) 
		\leq C \rho^{\alpha+\mu/2} Q_{\varsigma/2,\varepsilon}(T,\varphi,\mathbf{C}_1(0)) \text{ for all } \rho \in (0,\rho_Z], \nonumber 
\end{gather*}
where $C = C(n,m,q,\alpha) \in (0,\infty)$ is a constant; in particular, \eqref{SandT concl2} holds true.  Moreover, we can choose the multiplicities $\overline{m}_j$ of $\varphi^{(Z)}$ such that 
\begin{equation*}
	\pi_{\#} T^{(Z)}_j = q_j \overline{m}_j \llbracket \Omega \rrbracket .
\end{equation*}
Let $\rho_{\nu} \rightarrow 0^+$ and $\sigma \in [1,\infty)$.  Following the blow-up procedure of Subsection~\ref{sec:B tanfun sec}, the blow-up of $\eta_{Z,\rho_{\nu}\#} T$ relative to the plane $P_Z$ in $\mathbf{C}_{\sigma}(Z,P_Z)$ is $\omega_n^{1/2} \|\varphi^{(Z)}\|_{L^2(B_1(0,P_Z))}^{-1} \,\varphi^{(Z)}$.  That is, $\omega_n^{1/2} \|\varphi^{(Z)}\|_{L^2(B_1(0,P_Z))}^{-1} \,\varphi^{(Z)}$ is the unique tangent function to $T$ at $Z$.

Let $\beta_0 = \beta_0(n,m,q,\alpha) \in (0,1)$ be a small constant to be later determined.  Fix any $Z \in \Sigma$.  Recall that $p \in \{1,2,\ldots,\lfloor q/q_0 \rfloor\}$ such that $\overline{\varphi}^{(Z)} \in \Phi_{\alpha,q,p}$.  Arguing as in Remark~\ref{tilde varphi rmk}(1), find $p_* \in \{1,2,\ldots,p\}$ and $\varphi_* \in \Phi_{\alpha,q,p}$ such that 
\begin{gather}
	\label{SandT eqn44} Q_{\varsigma/2,\varepsilon}(\overline{\Gamma}_{Z\#}T^{(Z)},\varphi_*,\mathbf{C}_{s_k}(0)) 
	\leq 2^p \beta_0^{-p} Q_{\varsigma/2,\varepsilon}(\overline{\Gamma}_{Z\#}T^{(Z)},\overline{\varphi}^{(Z)}, \mathbf{C}_{s_k}(0)) , \\
	\label{SandT eqn45} Q_{\varsigma/2,\varepsilon}(\overline{\Gamma}_{Z\#}T^{(Z)},\varphi_*,\mathbf{C}_{s_k}(0)) 
	\leq \beta_0 \inf_{\varphi' \in \bigcup_{p'=0}^{p_*-1} \Phi_{\alpha,q,p'}} Q_{\varsigma/2,\varepsilon}(\overline{\Gamma}_{Z\#}T^{(Z)},\varphi',\mathbf{C}_{s_k}(0)) .  
\end{gather}
Assuming $\beta_0$ is sufficiently small, by \eqref{SandT eqn45} we can apply Theorem~\ref{nonconcentration thm}(A) together with \eqref{SandT eqn44}, \eqref{SandT hyp2}, \eqref{SandT eqn25}, and \eqref{SandT eqn37}
\begin{align*}
	&\op{dist}(\overline{\Gamma}_Z \widetilde{\Gamma}_Z (\mathcal{B}^{\star}_q - Z) \cap \mathbf{C}_{s_k/2}(0), \{0\}^2 \times \mathbb{R}^{n-2} \times \{0\}^m) 
	\\ \leq\,& \frac{C s_k \,Q_{\varsigma/2,\varepsilon}(\overline{\Gamma}_{Z\#}T^{(Z)},\varphi_*, \mathbf{C}_{s_k}(0))}{\widehat{E}(T^{(Z)},P_0,\mathbf{C}_{s_k}(0))}
	\leq \frac{C s_k \,Q_{\varsigma/2,\varepsilon}(\overline{\Gamma}_{Z\#}T^{(Z)},\overline{\varphi}^{(Z)}, \mathbf{C}_{s_k}(0))}{
		\widehat{E}(T^{(Z)},P_0,\mathbf{C}_{s_k}(0))}
	\leq C\beta s_k^{1+\mu/2} ,  
\end{align*}
where $C = C(n,m,q,\alpha) \in (0,\infty)$ is a constant.  For each $\rho \in (0,1/8]$ there exists an integer $k \geq 0$ such that $s_{k+1} < 2\rho \leq s_k$ and thus 
\begin{equation}\label{SandT eqn46}
	\op{dist}(\overline{\Gamma}_Z \widetilde{\Gamma}_Z (\mathcal{B}^{\star}_q - Z) \cap \mathbf{C}_{\rho}(0), \{0\}^2 \times \mathbb{R}^{n-2} \times \{0\}^m) 
	\leq C\beta \rho^{1+\mu/2} ,  
\end{equation}
where $C = C(n,m,q,\alpha) \in (0,\infty)$ is a constant.  Next let $Y,Z \in \Sigma$ with $|Y-Z| \leq 1/16$.  Choose an integer $k \geq 0$ such that $s_{k+1}/4 < |Y-Z| \leq s_k/4$.  Set $\widetilde{T} = \overline{\Gamma}_{Y\#} \eta_{0,s_k\#} T^{(Y)}$ and $\widetilde{Z} = s_k^{-1} \overline{\Gamma}_Y \widetilde{\Gamma}_Y(Z - Y)$.  Note that $\widetilde{T}$ is tangent to $\widetilde{P}_{\widetilde Z} = \overline{\Gamma}_Y \widetilde{\Gamma}_Y(P_Z)$ at $\widetilde{Z}$ and let $\widetilde{\Delta}_Z$ be an orthogonal rotation of $\mathbb{R}^{n+m}$ such that 
\begin{equation}\label{SandT eqn47}
    \widetilde{\Delta}_{\widetilde Z}(\widetilde{P}_{\widetilde Z}) = P_0 , \quad \|\widetilde{\Delta}_{\widetilde Z} - I\| \leq C \op{dist}_{\mathcal H}(\widetilde{P}_{\widetilde Z} \cap \mathbf{B}_1(0), P_0 \cap \mathbf{B}_1(0)) \leq C\varepsilon |Y-Z|^{\alpha-1}
\end{equation}
where $P_0 = \mathbb{R}^n \times \{0\}$ and $C = C(n,m,q,\alpha) \in (0,\infty)$ are constants.  Set $\widetilde{T}^{(\widetilde{Z})} = \widetilde{\Delta}_{Z\#} \eta_{\widetilde{Z},1\#} \widetilde{T}$.  By \eqref{SandT eqn17}, \eqref{SandT eqn25},~\cite[Lemma~3.5 and Lemma~3.7]{KrumWica}, and \eqref{oriented tilt and mass}, Hypothesis~$(\star)$ holds true with $C\varepsilon,\widetilde{T}$ in place of $\varepsilon_0,T$ for all $\rho \in (0,1/4]$, where $C = C(n,m,q,\alpha) \in (0,\infty)$ is a constant.  By \eqref{SandT hyp2}, \eqref{SandT eqn25}, and \eqref{SandT eqn37}, 
\eqref{SandT hyp2} holds true with $C\beta,\widetilde{T},$ $s_k^{\alpha-1} \overline{\varphi}^{(Y)}$ in place of $\beta,T,\varphi$.  Hence by arguing as in Step~2, the hypotheses of Theorem~\ref{main excess decay thm} hold true with $C\varepsilon, C\beta, \eta_{0,1/4\#} \widetilde{T}^{(\widetilde{Z})}, (4s_k)^{\alpha-1} \overline{\varphi}^{(Y)}$ in place of $\varepsilon,\beta,T^{(Z)},\varphi$.  By iteratively applying Theorem~\ref{main excess decay thm} as in Step~3 above with $\widetilde{T}^{(\widetilde{Z})},s_k^{\alpha-1} \overline{\varphi}^{(Y)}$ in place of $T^{(Z)},\varphi$ and starting with $\Delta_0 = I$, $\varphi_0 = s_k^{\alpha-1} \overline{\varphi}^{(Y)}$, and $t_0 = 1/4$, for each $i = 1,2,3,\ldots$ there exist an orthogonal rotation $\Delta_i$ of $\mathbb{R}^{n+m}$ which fixes vectors in $\{0\} \times \mathbb{R}^m$, $\psi_i \in \Phi_{\alpha,q}$, $j'_i \in \{1,2,\ldots,\lfloor q/q_0 \rfloor\}$, and a radius $t_i$ given by $t_i = \vartheta_{j'_i} t_{i-1}$ such that \eqref{SandT eqn12}--\eqref{SandT eqn15} hold true with $\widetilde{T}^{(\widetilde{Z})}$, $\Delta_i$, $\psi_i$, $j'_i$, and $t_i$ in place of $T^{(Z)}$, $\Gamma_k$, $\varphi_k$, $j_k$, and $s_k$.  Note that by \eqref{SandT eqn30} with $\rho = t_i s_k$, alternative~(A) of Theorem~\ref{main excess decay thm} does not hold true with $\Delta_{i\#} \eta_{0,t_i\#} \widetilde{T}^{(\widetilde{Z})}$ in place of $T$ and thus alternative~(B) of Theorem~\ref{main excess decay thm} gives us $\Delta_{i+1}$, $\psi_{i+1}$, $j'_{i+1}$, and $t_{i+1}$.  Let $\overline{\Delta}_{\widetilde Z}$ be an orthogonal rotation of $\mathbb{R}^{n+m}$ and $\overline{\psi}^{(\widetilde{Z})} \in \Phi_{\alpha,q}$ such that $\|\Delta_i - \overline{\Delta}_{\widetilde Z}\| \rightarrow 0$ and $\psi_i \rightarrow \overline{\psi}^{(\widetilde{Z})}$ uniformly on compact subsets of $\mathbb{R}^n$ so that $\widetilde{\Delta}_{\widetilde Z}^{-1} \circ \overline{\psi}^{(\widetilde{Z})} \circ \overline{\Delta}_{\widetilde Z} \circ \widetilde{\Delta}_{\widetilde Z}$ is the unique tangent function to $\widetilde{T}$ at $\widetilde{Z}$.  On the other hand $\widetilde{\Gamma}_Z^{-1} \circ \overline{\varphi}^{(Z)} \circ \overline{\Gamma}_Z \circ \widetilde{\Gamma}_Z$ is the unique tangent function tangent to $T$ at $Z$ and thus 
\begin{equation*}
    \widetilde{\Delta}_{\widetilde Z}^{-1} \circ \overline{\psi}^{(\widetilde{Z})} \circ \overline{\Delta}_{\widetilde Z} \circ \widetilde{\Delta}_{\widetilde Z} = \overline{\Gamma}_Y \circ \widetilde{\Gamma}_Y \circ \widetilde{\Gamma}_Z^{-1} \circ \overline{\varphi}^{(Z)} \circ \overline{\Gamma}_Z \circ \widetilde{\Gamma}_Z \circ \widetilde{\Gamma}_Y^{-1} \circ \overline{\Gamma}_Y^{-1}
\end{equation*}
on $\widetilde{P}_{\widetilde Z}$ as the common unique tangent function to $\widetilde{T}$ at $\widetilde{Z}$.  Since the corresponding spines must be equal, 
\begin{equation}\label{SandT eqn49}
    \widetilde{\Delta}_{\widetilde Z}^{-1} \overline{\Delta}_{\widetilde Z}^{-1}(\{0\}^2 \times \mathbb{R}^{n-2} \times \{0\}^m) = \overline{\Gamma}_Y \widetilde{\Gamma}_Y \widetilde{\Gamma}_Z^{-1} \overline{\Gamma}_Z^{-1}(\{0\}^2 \times \mathbb{R}^{n-2} \times \{0\}^m) . 
\end{equation}
In place of \eqref{SandT eqn32} we obtain 
\begin{equation}\label{SandT eqn50}
	\|\overline{\Delta}_{\widetilde Z} - I\| 
	\leq \frac{C Q_{\varsigma/2,\varepsilon}(\overline{\Gamma}_{Y\#} T^{(Y)},\overline{\varphi}^{(Y)}, 
		\mathbf{C}_{s_k}(0))}{\widehat{E}(\overline{\Gamma}_{Y\#} T^{(Y)},P_0,\mathbf{C}_{s_k}(0))} 
	\leq C\beta |Y-Z|^{\mu/2} , 
\end{equation}
where $C = C(n,m,q,\alpha) \in (0,\infty)$ are constants.  By \eqref{SandT eqn49}, \eqref{SandT eqn47}, and \eqref{SandT eqn50} 
\begin{align}\label{SandT eqn51}
    &\op{dist}_{\mathcal H}(\widetilde{\Gamma}_Y^{-1} \overline{\Gamma}_Y^{-1}(\{0\}^2 \times \mathbb{R}^{n-2} \times \{0\}^m) \cap \mathbf{B}_1(0), \widetilde{\Gamma}_Z^{-1} \overline{\Gamma}_Z^{-1}(\{0\}^2 \times \mathbb{R}^{n-2} \times \{0\}^m) \cap \mathbf{B}_1(0)) 
    \\=\,& \op{dist}_{\mathcal H}(\widetilde{\Delta}_{\widetilde Z}^{-1} \overline{\Delta}_{\widetilde Z}^{-1}(\{0\}^2 \times \mathbb{R}^{n-2} \times \{0\}^m) \cap \mathbf{B}_1(0), (\{0\}^2 \times \mathbb{R}^{n-2} \times \{0\}^m) \cap \mathbf{B}_1(0)) \nonumber
    \\ \leq\,& C(\varepsilon+\beta)\, |Y-Z|^{\min\{\alpha-1,\mu/2\}} , \nonumber 
\end{align}
where $C = C(n,m,q,\alpha) \in (0,\infty)$ is a constant.  It follows from \eqref{SandT eqn46} and \eqref{SandT eqn51} that there exist functions $f \in C^{1,\mu/2}(B^{n-2}_{1/2}(0),\mathbb{R}^2)$ and $g \in C^{1,\mu/2}(B^{n-2}_{1/2}(0),\mathbb{R}^m)$ such that $\Sigma \subseteq \{ (f(y),y,g(y)) : y \in B^{n-2}_{1/2}(0) \}$ and $\|f\|_{C^{1,\mu/2}(B^{n-2}_{1/2}(0))} + \|g\|_{C^{1,\mu/2}(B^{n-2}_{1/2}(0))} \leq C \beta$, where $C = C(n,m,q,\alpha) \in (0,\infty)$ is a constant.

\noindent\textit{Step~5:  Covering argument for the sets $\Upsilon_k$ of branch points with small gaps.}  Suppose that $Z \in \Upsilon_{k_0}$ for some positive integer $k_0 < \infty$.  Take $\overline{\Gamma}^{(Z)} = \Gamma_{k_0}$ and $\overline{\varphi}^{(Z)} = \varphi_{k_0}$.  Note that $\overline{\Gamma}^{(Z)}$ and $\varphi^{(Z)}$ are no longer unique (even after changing of the multiplicities of $\varphi^{(Z)}$).  Observe that \eqref{SandT eqn33} holds true for all $k \in \{1,2,\ldots,k_0\}$ and in particular $\|\overline{\varphi}^{(Z)}\|_{L^2(B_1(0))} \leq C \widehat{E}(T,P_0,\mathbf{C}_1(0))$.  Using this together with \eqref{SandT eqn22} and \eqref{SandT eqn28} 
\begin{align}\label{SandT eqn52}
	&\op{dist}_{\mathcal H}(\op{spt}\op{graph}(\varphi_j \circ \Gamma_j) \cap \mathbf{C}_1(0), 
		\op{spt}\op{graph}(\overline{\varphi}^{(Z)} \circ \overline{\Gamma}_Z) \cap \mathbf{C}_1(0)) 
	\\ \leq\,& \op{dist}_{\mathcal H}(\op{spt}\op{graph}\varphi_j \cap \mathbf{C}_1(0), \op{spt}\op{graph}\overline{\varphi}^{(Z)} \cap \mathbf{C}_1(0)) 
		+ C \widehat{E}(T,P_0,\mathbf{C}_1(0)) \,\|\Gamma_j - \overline{\Gamma}_Z\| \nonumber 
	\\ \leq\,& C s_j^{\mu/2} Q_{\varsigma/2,\varepsilon}(T,\varphi,\mathbf{C}_1(0)) , \nonumber
\end{align}
for each integer $j \in \{0,1,2,\ldots,k_0\}$, where $C = C(n,m,q,\alpha) \in (0,\infty)$ are constants.  By the triangle inequality, \eqref{SandT eqn20}, and \eqref{SandT eqn52}
\begin{align}\label{SandT eqn53}
	&Q_{\varsigma/2,\varepsilon}(\overline{\Gamma}_{Z\#}T^{(Z)},\overline{\varphi}^{(Z)}, \mathbf{C}_{s_j}(0)) 
	\leq Q_{\varsigma/2,\varepsilon}(\Gamma_{k\#} T^{(Z)},\varphi_j, \mathbf{C}_{s_j}(0)) 
		\\&\hspace{15mm} + C s_j^{\alpha-1} \op{dist}_{\mathcal H}(\op{spt}\op{graph}(\varphi_j \circ \Gamma_j) \cap \mathbf{C}_1(0), 
		\op{spt}\op{graph}(\overline{\varphi}^{(Z)} \circ \overline{\Gamma}_Z) \cap \mathbf{C}_1(0)) \nonumber
	\\ \leq\,& C s_j^{\alpha-1+\mu/2} Q_{\varsigma/2,\varepsilon}(T,\varphi,\mathbf{C}_1(0)) \nonumber 
\end{align}
for each integer $j \geq 0$, where $C = C(n,m,q,\alpha) \in (0,\infty)$ is a constant.  By the argument leading to \eqref{SandT eqn46}, using \eqref{SandT eqn53} in place of \eqref{SandT eqn37},
\begin{equation*}
	\op{dist}(\overline{\Gamma}_Z^{-1} \widetilde{\Gamma}_Z^{-1} (\mathcal{B}^{\star}_q - Z) \cap \mathbf{C}_{\rho}(0), 
		\{0\}^2 \times \mathbb{R}^{n-2} \times \{0\}^m) \leq C\beta \rho^{1+\mu/2}
\end{equation*}
for all $\rho \in [s_{k_0},1/4]$, where $C = C(n,m,q,\alpha) \in (0,\infty)$ is a constant.  Since $\|\widetilde{\Gamma}_Z - I\| \leq C \widehat{E}(T,P_0,\mathbf{C}_1(0)) < C\varepsilon$ and, by \eqref{SandT eqn28} and \eqref{SandT hyp2}, $\|\overline{\Gamma}_Z - I\| \leq C\beta$, 
\begin{equation}\label{SandT eqn54}
	\op{dist}(\mathcal{B}^{\star}_q \cap \mathbf{C}_{\rho}(Z,P_0), 
		Z + \{0\}^2 \times \mathbb{R}^{n-2} \times \{0\}^m) \leq C(\varepsilon+\beta) \rho 
\end{equation}
for all $\rho \in [s_{k_0},1/4]$, where $C = C(n,m,q,\alpha) \in (0,\infty)$ are constants.  By the definition of $\Upsilon_{k_0}$ 
\begin{equation*}
	\forall Z \in \Upsilon_{k_0} \;\; \exists y_0 \in B^{n-2}_{s_{k_0}/2}(0) \text{ such that } \widetilde{\Gamma}_Z(\mathcal{B}^{\star}_q - Z) 
		\cap \mathbf{C}_{\delta \vartheta_{\lfloor q/q_0 \rfloor} s_{k_0}/2}(0,y_0) = \emptyset
\end{equation*}
Since $\|\widetilde{\Gamma}_Z - I\| \leq C\beta$, provided $\beta$ is sufficiently small 
\begin{equation}\label{SandT eqn55}
	\forall Z \in \Upsilon_{k_0} \;\; \exists Y \in Z + \{0\}^2 \times B^{n-2}_{s_{k_0}/2}(0) \times \{0\}^m \text{ such that } \mathcal{B}^{\star}_q 
		\cap \mathbf{C}_{\delta \vartheta_{\lfloor q/q_0 \rfloor} s_{k_0}/4}(Y) = \emptyset
\end{equation}
Assuming $\delta = \delta(n,m,q,\alpha,\eta)$ is sufficiently small and arguing as in the proof of~\cite[Theorem~1]{Sim93} using \eqref{SandT eqn54} and \eqref{SandT eqn55} in place of (12) and (13) on p.~642 of~\cite{Sim93}, we obtain a covering of $\bigcup_{1 \leq k < \infty} \Upsilon_k$ by a finite collection of balls $\{B^{n+m}_{\rho_j}(X_j)\}_{1 \leq j \leq N}$ such that $\sum_{j=1}^N \rho_j^{n-2} \leq 1-\eta$.
\end{proof}

\begin{theorem}\label{precise rectifiability thm}
Let $q \geq 2$ be an integer and $\mu \in (0,1)$.  There exists $\varepsilon_0 = \varepsilon_0(n,m,q,\mu) \in (0,1)$ such that the following holds true.  Let $T$ be as in Hypothesis~\ref{introduction hyp}, $\mathcal{B}_q$ be as in Definition~\ref{introduction B defn}, and $\mathcal{B}^{(\neq 2)}_q = \{ X \in \mathcal{B}_q : \mathcal{N}_{T,{\rm Pl}}(Z) \neq 2 \}$.  For each $\alpha > 1$ with $\alpha \neq 2$, let $\mathcal{B}_{q,\alpha}$ be the set of all branch points $Z \in \mathcal{B}_q$ such that $\mathcal{N}_{T,{\rm Pl}}(Z) = \alpha$ and $T$ has a tangent function $\varphi$ at $Z$ with $\dim\op{spine}\varphi = n-2$.  Then $\mathcal{B}^{(\neq 2)}_q$ is countably $(n-2)$-rectifiable.  Moreover:
\begin{enumerate}[itemsep=2mm,topsep=0mm]
	\item[(a)]  for any compact set $K \subset \mathbf{C}_{1/2}(0)$, $K \cap \mathcal{B}_{q,\alpha} \neq \emptyset$ for only finitely many $\alpha > 0$; 
	
	\item[(b)]  for every $\alpha > 1$ with $\alpha \neq 2$ and for $\mathcal{H}^{n-2}$-a.e.~$Z \in \mathcal{B}_{q,\alpha}$, there exists a unique tangent function $\varphi^{(Z)} : P_Z \rightarrow \mathcal{A}_q(P_Z^{\perp})$ of $T$ at $Z$ such that $\dim\op{spine}\varphi^{(Z)} = n-2$ and 
	\begin{equation*}
        \frac{1}{\rho^n} \int_{\mathbf{C}_{\rho}(Z,P_Z)} \op{dist}^2(X,Z+\op{spt}\op{graph}\varphi^{(Z)}) \,d\|T\|(X) \leq C_Z \rho^{2\alpha+2\overline{\mu}}
	\end{equation*}
    for all $\rho \in (0,\rho_Z]$ and 
    \begin{align*}
	   &\op{dist}_{\mathcal H}(\op{spt} T \cap \mathbf{C}_{\rho}(Z,P_Z) \cap \{r_Z < \tau_Z(\rho)\rho\}, \\&\hspace{0.85in} (Z+\op{spt}\op{graph}\varphi^{(Z)}) \cap \mathbf{C}_{\rho}(Z,P_Z) \cap \{r_Z < \tau_Z(\rho)\rho\}) \leq C_Z \rho^{\alpha+\overline{\mu}/2} \nonumber 
    \end{align*}
    for all $\rho \in (0,\rho_Z]$, where $\overline{\mu} = \overline{\mu}(n,m,q,\alpha) \in (0,1)$ is a constant, $\rho_Z \in (0,1/4]$ and $C_Z \in (0,\infty)$ are constants (depending on $Z$), and $\tau_Z : (0,\rho_Z] \rightarrow (0,1/16]$ is a function with $\lim_{\rho\rightarrow 0^+} \tau_Z(\rho) = 0$ (also depending on $Z$); 
	
	\item[(c)]  for every $\alpha > 1$ with $\alpha \neq 2$, there is an open set $V_{\alpha} \supset \mathcal{B}_{q,\alpha}$ such that $V_{\alpha} \cap \{X \in \mathcal{B}_q : \mathcal{N}_{T,{\rm Pl}}(X) \geq \alpha \}$ has locally finite $\mathcal{H}^{n-2}$-measure, i.e.~for every $Y \in V_{\alpha} \cap \{X \in \mathcal{B}_q : \mathcal{N}_{T,{\rm Pl}}(X) \geq \alpha \}$ there exists $\rho > 0$ such that $\mathcal{H}^{n-2}(B_{\rho}(Y) \cap \{X \in \mathcal{B}_q : \mathcal{N}_{T,{\rm Pl}}(X) \geq \alpha \}) < \infty$; in particular, $\mathcal{B}_{q,\alpha}$ has locally finite $\mathcal{H}^{n-2}$-measure.
\end{enumerate}
\end{theorem}

\begin{proof} 
The proof is similar to that of~\cite[Theorem~$2'$]{Sim93}, so we will only sketch it here.  Let $T$ be an $n$-dimensional locally area minimizing current of $\mathbf{C}_1(0)$ satisfying Hypothesis~\ref{introduction hyp}.  By Theorem~\ref{B prime stratification thm}, the set $S^{(n-3)}_{q}$ of all points $Z \in \mathcal{B}^{(\neq 2)}_q$ at which every tangent function has spine dimension $\leq n-3$ satisfies $\dim_{\mathcal H} S^{(n-3)}_{q} \leq n-3$.  Recall that the degree of homogeneity $\alpha$ of non-zero locally Dirichlet energy minimizing $q$-valued functions with spine dimension $n-2$ must equal $\ell_0/q_0$ for some relatively prime positive integers $\ell_0,q_0$ with $q_0 \leq q$ and in particular the set of all such $\alpha$ is discrete.  Hence it suffices to fix $\alpha = \ell_0/q_0 \neq 2$, where $\ell_0, q_0$ are relatively prime positive integer with $q_0 \leq q$, and consider the set $\mathcal{B}_{q,\alpha}$ of all points $Z \in \mathcal{B}_q$ at which $\mathcal{N}_{T,{\rm Pl}}(Z) = \alpha$ and least one tangent function has spine dimension $= n-2$. 
 
Let $Z_0 \in \mathcal{B}_{q,\alpha}$ and $\varphi$ be a tangent function of $T$ at $Z_0$ with $\dim \op{spine}\varphi = n-2$.  Recall that from Remark~\ref{planar freq rmk} that we can apply Lemma~\ref{uniform branch set lemma} with $\mu = 1/q$.  By Lemma~\ref{uniform branch set lemma} and the monotonicity formula for planar frequency functions Theorem~\ref{frequency monotonicity thm}, for every $\varepsilon > 0$ and for every suitably small $\sigma_0 \in \big( 0,\frac{1}{16 R(\varepsilon)} \big]$ 
\begin{gather}\label{precise rectifiability eqn1}
	(\partial T) \llcorner \mathbf{C}_{\sigma_0 R(\varepsilon)}(Z_0,P_{Z_0}) = 0, \quad 
	\sup_{X \in T \cap \mathbf{C}_{\sigma_0 R(\varepsilon)}(Z_0,P_{Z_0})} \op{dist}(X,P_0) < \infty, \\ 
	(\pi_{P_{Z_0}\#} T) \llcorner \mathbf{C}_{\sigma_0 R(\varepsilon)}(Z_0,P_{Z_0}) = q \llbracket B_{\sigma_0 R(\varepsilon)}(Z_0,P_{Z_0}) \rrbracket, 
		\nonumber\\
	\widehat{E}(T,P_{Z_0},\mathbf{C}_{7\rho/4}(Z_0,P_{Z_0})) \leq C\varepsilon_0 \Big(\frac{\rho}{\sigma_0 R(\varepsilon)}\Big)^{1/q} 
		\text{ for all } \rho \in (0,\sigma_0 R(\varepsilon)/2] , \nonumber\\
	N_{T,P_{Z_0},Z_0}(\sigma_0 R(\varepsilon)/8) < \alpha + \delta(\varepsilon) , \nonumber
\end{gather}
where $R(\varepsilon)$ and $\delta(\varepsilon)$ are as in Lemma~\ref{lemma2_4}.  By the definition of tangent function in Subsection~\ref{sec:B tanfun sec}, there exists a sequence $\rho_{\nu} \rightarrow 0^+$ such that $\varphi$ is the blow-up of $\eta_{Z_0,\rho_{\nu}\#} \widetilde{T}$ with respect to the plane $P_{Z_0}$ by excess $\widehat{E}(\widetilde{T},P_{Z_0},\mathbf{C}_{\rho_{\nu}}(Z_0,P_{Z_0}))$.  Hence arguing as in the proof of Lemma~\ref{lemma2_4}(ii), we can choose $\sigma_0 = \rho_{\nu}/R(\varepsilon)$ for $\nu$ large enough that 
\begin{equation*}
	Q(T,\widehat{E} (\sigma_0 R(\varepsilon))^{\alpha-1} \varphi,\mathbf{C}_{\sigma_0 R(\varepsilon)}(Z_0,P_{Z_0})) < \varepsilon \,\widehat{E} ,
\end{equation*}
where $\widehat{E} = \widehat{E}(T,P_{Z_0},\mathbf{C}_{\sigma_0 R(\varepsilon)}(Z_0,P_{Z_0}))$ (and \eqref{precise rectifiability eqn1} also continues to hold true).  Set $T_0 = \eta_{Z_0,\sigma_0\#} T$.

For each $\rho_0 \in (0,1/2]$, define the outer measure $\mu_{\rho_0}$ on $\mathbf{C}_1(0)$ by 
\begin{equation*}
	\mu_{\rho_0}(A) = \inf \sum_{i=1}^N \omega_{n-2} \sigma_i^{n-2} 
\end{equation*}
for every set $A \subseteq \mathbf{C}_1(0)$, where the infimum is taken over all finite covers of $A$ by open balls $\mathbf{B}_{\sigma_i}(Y_i)$ for $i = 1,2,\ldots,N$ with $\sigma_i \leq \rho_0$.  Let $\mathcal{B}^+_{\alpha} = \{ X \in \mathcal{B}'_q(T_0) \cap \overline{\mathbf{C}_1(0)} : \Theta(T,X) \geq \alpha \}$ and note that $\eta_{Z_0,\sigma}(\mathcal{B}'_{q,\alpha}) \subseteq \mathcal{B}^+_{\alpha}$.  Choose a cover of $\mathcal{B}^+_{\alpha}$ by a finite collection of open balls $\mathbf{B}_{\sigma_i}(Y_i)$ such that 
\begin{equation*}
	\sum_{i=1}^N \omega_{n-2} \sigma_i^{n-2} \leq \mu_{\rho_0}(\mathcal{B}^+_{\alpha}) + 1
\end{equation*}
Remove the balls $\mathbf{B}_{\sigma_i}(Y_i)$ which do not intersect $\mathcal{B}^+_{\alpha}$ from the collection.  For each $i$, let $Z_i \in \mathbf{B}_{\sigma_i}(Y_i) \cap \mathcal{B}^+_{\alpha}$.  By Lemma~\ref{lemma2_4}, either there exists a non-zero, homogeneous degree $\alpha$, locally Dirichlet energy minimizing $q$-valued function $\varphi_i \in W^{1,2}_{\rm loc}(P_{Z_i},\mathcal{A}_q(P_{Z_i}^{\perp}))$ such that $\dim \op{spine}\varphi_i = n-2$ and 
\begin{equation}\label{precise rectifiability eqn2}
	Q(T_0,\varphi_i,\mathbf{C}_{2\sigma_i}(Z_i,P_{Z_i})) < \varepsilon \,\widehat{E}(T_0,P_{Z_i},\mathbf{C}_{\rho}(Z_i,P_{Z_i})) 
\end{equation}
or there exists an $(n-3)$-dimensional linear subspace $L_i$ of $P_{Z_i}$ such that 
\begin{equation}\label{precise rectifiability eqn3}
	\{ X \in \mathcal{B}'_q(T_0) \cap \overline{\mathbf{C}_{2\sigma_i}(Z_i)} : \mathcal{N}_{T,{\rm Pl}}(X) \geq \alpha \} 
			\subset \{ X \in \mathbb{R}^{n+m} : \op{dist}(X,Z_i+L_i) < \varepsilon\rho \} .
\end{equation}
The conclusion of the theorem now follows by arguing exactly as in~\cite{Sim93} by iteratively applying Theorem~\ref{SandT thm} using the fact that either \eqref{precise rectifiability eqn2} or \eqref{precise rectifiability eqn3} holds true.
\end{proof}

\begin{proof}[Proof of Theorem~\ref{unique tangent function thm} of the Introduction] 
By Theorem~\ref{B prime stratification thm}, for $\mathcal{H}^{n-2}$-a.e.~$Z \in \mathcal{B}^{(\neq 2)}_q$, at least one tangent function $\varphi$ of $T$ at $Z$ has $\dim \op{spine} \varphi = n-2$, and thus $Z \in \mathcal{B}_{q,\alpha}$ where $\alpha = \mathcal{N}_{T,{\rm Pl}}(Z)$.  Hence the conclusion is a direct consequence of Theorem~\ref{precise rectifiability thm}(b).
\end{proof}

\begin{proof}[Proof of Theorem~\ref{B rectifiable thm} of the Introduction]
Let $T$ be an $n$-dimensional locally area minimizing current of $\mathbf{C}_1(0)$ satisfying Hypothesis~\ref{introduction hyp} and let $K \subset \mathbf{C}_{1/2}(0)$ be a compact set.  By Theorem~\ref{precise rectifiability thm}(a), there is a finite set $\{\alpha_1,\alpha_2,\ldots,\alpha_k\}$ such that $\alpha_j < 2$ and $K \cap \mathcal{B}_{q,\alpha_j} \neq \emptyset$ for all $j = 1,2,\ldots,k$.  By Theorem~\ref{precise rectifiability thm}(c), for each $j = 1,2,\ldots,k$ there exists an open set $V_{\alpha_j} \supset \mathcal{B}_{q,\alpha_j}$ such that $V_{\alpha_j} \cap \{X \in \mathcal{B}_q : \mathcal{N}_{T,{\rm Pl}}(Z) \geq \alpha_j \}$ has locally finite $\mathcal{H}^{n-2}$-measure in $V_{\alpha_j}$.  Noting that by Theorem~\ref{frequency monotonicity thm}(c) we have $\mathcal{N}_{T,{\rm Pl}}(Z) \geq 1+\mu$ for all $Z \in \mathcal{B}_q$, set $\alpha_0 = 1+\mu$ and $\alpha_{k+1} = 2$.  For $j = 0,1,2,\ldots,k$ let 
\begin{equation*}
	\Gamma_j = \{ X \in K \cap \mathcal{B}_q : \alpha_j \leq \mathcal{N}_{T,{\rm Pl}}(X) < \alpha_{j+1} \} \cap V_{\alpha_j} 
\end{equation*}
so that $\Gamma_j$ has locally finite $\mathcal{H}^{n-2}$-measure in $V_{\alpha_j}$ and let 
\begin{equation*}
	\widetilde{\Gamma}_j = \{ X \in K \cap \mathcal{B}_q : \alpha_j \leq \mathcal{N}_{T,{\rm Pl}}(X) < \alpha_{j+1} \} \setminus V_{\alpha_j} 
\end{equation*}
so that $\widetilde{\Gamma}_j$ is a subset of the $(n-3)$-stratum $S^{(n-3)}_{q}$ of $\mathcal{B}^{(\neq 2)}_q$ and thus by Theorem~\ref{B prime stratification thm} we have that $\dim_{\mathcal H} \widetilde{\Gamma}_j \leq n-3$.  Recall from Subsection~\ref{sec:stratify sec} that $\mathcal{B}_q$ is a relatively closed subset of $\mathbf{C}_{1/2}(0)$.  By Lemma~\ref{uniform branch set lemma} and~\cite[Corollary~5.5]{KrumWica}, $\mathcal{N}_{T,{\rm Pl}}$ is an upper semi-continuous function on $\mathcal{B}_q$.  It follows that each of $\Gamma_j$ and $\widetilde{\Gamma}_j$ is the intersection of an open set and closed set, and hence is locally compact.  Clearly $\mathcal{B}^{(< 2)}_q \cap K = \bigcup_{j=0}^{\infty} (\Gamma_j \cup \widetilde{\Gamma}_j)$.
\end{proof}

\begin{proof}[Proof of Theorem~\ref{topology thm}] 
\textit{Step~1.}  We claim that: \textit{for every $\beta \in (0,1)$ there exists $\varepsilon = \varepsilon(n,m,q,\alpha,\lambda,\beta) \in (0,1)$ such that assuming the hypotheses of Theorem~\ref{topology thm}, there exists a non-zero, homogeneous degree $\alpha$, locally Dirichlet energy minimizing function $\varphi \in W^{1,2}_{\rm loc}(\mathbb{R}^n,\mathcal{A}_q(\mathbb{R}^m))$ such that $\op{dim}\op{spine}\varphi = n-2$ and 
\begin{gather}
    \label{topology eqn1} Q(T,\varphi,\mathbf{C}_{\lambda}(0)) < \beta \,\widehat{E}(T,P_0,\mathbf{C}_{\lambda}(0)) , \\
    \label{topology eqn2} \{ X \in \mathcal{B}_q \cap \mathbf{B}_{1/4}(0) : \mathcal{N}_{T,{\rm Pl}}(X) \geq \alpha \} \subset \{ X \in \mathbf{B}_{1/4}(0) : \op{dist}(X,\op{spine}\varphi) < \lambda\beta \} , \\
    \label{topology eqn3} \forall Y \in L \;\; \exists ! X \in \op{spine}\varphi \text{ such that } 2\lambda |X| \leq |Y| \text{ and } \pi_L(X) = Y .
\end{gather}}

To see this, fix $\lambda \in (0,1/4)$ and $\beta \in (0,1)$.  For $\nu = 1,2,3,\ldots$ let  $\varepsilon_{\nu}\rightarrow 0^+$, $T^{(\nu)}$ be an $n$-dimensional locally area-minimizing rectifiable current of $\mathbf{C}_1(0)$, and $L_{\nu} \subset \mathbb{R}^{n+m}$ be an $(n-2)$-dimensional linear subspace such that the hypotheses of Theorem~\ref{topology thm} holds true with $\varepsilon_{\nu},T^{(\nu)},L_{\nu}$ in place of $\varepsilon,T,L$.  Set $\widehat{E}_{\nu} = \widehat{E}(T^{(\nu)},P_0,\mathbf{C}_{\lambda}(0))$.  Recall that $(\partial T^{(\nu)})\llcorner\mathbf{C}_1(0) = 0$ and $\Theta(T^{(\nu)},0) = q$.  By \eqref{introduction hyp2}, \eqref{oriented tilt and mass}, and~\cite[Lemma~3.7]{KrumWica}, $\|T^{(\nu)}\|(\mathbf{C}_{7/8}(0)) \leq (q+C(n,m)\,\varepsilon_0^2) \,\omega_n (7/8)^n$.  By Lemma~\ref{uniform branch set lemma}, $\widehat{E}(T^{(\nu)},P_0,\mathbf{C}_{\rho}(0)) \leq C(n,m,q) \,\varepsilon_0 \rho^{1/q}$ for all $\rho \in (0,7/8]$.  Hence by Theorem~\ref{frequency monotonicity thm}(b)(c) and \eqref{topology hyp1}
\begin{gather}
	\label{topology eqn4} N_{T^{(\nu)},P_0,0}(\rho) \leq e^{C\varepsilon_{\nu}^{\varsigma}} (\alpha+\varepsilon_{\nu}) \text{ for all } \rho \in (0,1/2] , \\
	\widehat{E}(T^{(\nu)},P_0,\mathbf{C}_{1/2}(0)) \leq C \widehat{E}_{\nu} , \nonumber 
\end{gather}
where $C= C(n,m,q,\alpha,\lambda) \in (0,\infty)$ are constants.  Let $\varphi \in W^{1,2}_{\rm loc}(B_{1/2}(0),\mathcal{A}_q(\mathbb{R}^n))$ be the blow-up of $T^{(\nu)}$ with respect to the plane $P_0$ by excess $\widehat{E}_{\nu}$.  Since $0 \in \mathcal{B}_q$ and $T$ is tangent to $P_0 = \mathbb{R}^n \times \{0\}$ at the origin, by~\cite[Lemma~6.4 and Corollary~6.6]{KrumWica} we have that $\varphi_a(0) = 0$ and $D\varphi_a(0) = 0$.  Moreover, by~\cite[Lemma~6.5 and Corollary~6.6]{KrumWica} and \eqref{topology eqn4}, $\mathcal{N}_{\varphi}(0) \geq \alpha$ and $N_{\varphi,0}(\rho) \leq \alpha$ for all $\rho \in (0,1/2)$ so that $\varphi$ is homogeneous degree $\alpha$.  Since $\alpha$ is not an integer, $\varphi$ is average-free.  After passing to a subsequence, let $L_{\infty} \subset \mathbb{R}^{n+m}$ be an $(n-2)$-dimensional linear subspace such that $\op{dist}_{\mathcal H}(L_{\nu} \cap \mathbf{B}_1(0),L_{\infty} \cap \mathbf{B}_1(0)) \rightarrow 0$.

We want to show that 
\begin{equation}\label{topology eqn5}
	L_{\infty} \cap \mathbf{B}_{\lambda}(0) \subseteq \pi_{L_{\infty}}(\op{spine}\varphi \cap \overline{\mathbf{B}_{1/4}(0)}) . 
\end{equation}
Let $\zeta \in L_{\infty} \cap \mathbf{B}_{\lambda}(0)$.  Since $L_{\nu}$ converges to $L_{\infty}$, there exists $\zeta_{\nu} \in L_{\nu} \cap \mathbf{B}_{\lambda}(0)$ such that $\zeta_{\nu} \rightarrow \zeta$.  By \eqref{topology hyp2}, there exists $X_{\nu} \in \mathcal{B}_q \cap \mathbf{B}_{1/4}(0)$ such that $\mathcal{N}_{T^{(\nu)},{\rm Pl}}(X_{\nu}) \geq \alpha$ and $\pi_{L_{\nu}}(X_{\nu}) = \zeta_{\nu}$.  After passing to a subsequence, there exists $\xi \in \overline{B_{1/4}(0)}$ such that $X_{\nu} \rightarrow (\xi,0)$.  Since $\varphi$ is average-free and $\mathcal{N}_{T^{(\nu)},{\rm Pl}}(X_{\nu}) \geq \alpha$ for all $\nu$, by~\cite[Corollary~6.6]{KrumWica} we have that $\mathcal{N}_{\varphi}(\xi) \geq \alpha$, i.e.~$\xi \in \op{spine}\varphi$.  Since $\pi_{L_{\nu}}(X_{\nu}) = \zeta_{\nu}$ for all $\nu$, $\pi_{L_{\infty}}(\xi) = \zeta$.  Therefore, for every $\zeta \in L_{\infty} \cap \mathbf{B}_{\lambda}(0)$ there exists $\xi \in \op{spine}\varphi \cap \overline{\mathbf{B}_{1/4}(0)}$ such that $\pi_{L_{\infty}}(\xi) = \zeta$.

By linearity, \eqref{topology eqn5} implies that $L_{\infty} = \pi_{L_{\infty}}(\op{spine}\varphi)$.  Since $\op{spine}\varphi \leq n-2 = \dim L_{\infty}$, $\dim \op{spine}\varphi = n-2$ and $\pi_{L_{\infty}} |_{\op{spine}\varphi} : \op{spine}\varphi \rightarrow L_{\infty}$ is a linear isomorphism.  Arguing as in the proof of Lemma~\ref{lemma2_4}(ii), since $\varphi$ is the blow-up of $T^{(\nu)}$ with respect to the plane $P_0$ by excess $\widehat{E}_{\nu}$, \eqref{topology eqn1} holds true with $T^{(\nu)}$ in place of $T$ for all sufficiently large $\nu$.  \eqref{topology eqn2} follows by arguing as in the proof of Theorem~\ref{B prime stratification thm} and in particular using~\cite[Corollary~6.6]{KrumWica}.  By \eqref{topology eqn5} and since $L_{\nu}$ converges to $L_{\infty}$, $L_{\nu} \cap \mathbf{B}_{\lambda/2}(0) \subset \pi_{L_{\nu}}(\op{spine}\varphi \cap \mathbf{B}_{1/4}(0))$.  Hence $\pi_{L_{\nu}} |_{\op{spine}\varphi} : \op{spine}\varphi \rightarrow L_{\nu}$ is a linear isomorphism and \eqref{topology eqn3} holds true.  In light of the arbitrary choice of sequences $(\varepsilon_{\nu})$, $(T^{(\nu)})$, and $(L_{\nu})$, this completes the proof of the claim.

\noindent\textit{Step~2: Iteratively applying the excess decay lemma and proving the no small gaps condition.}  Let $\varepsilon = \varepsilon(n,m,q,\alpha,\lambda) \in (0,1)$, $\beta = \beta(n,m,q,\alpha,\lambda) \in (0,1)$, and $\delta = \delta(n,m,q,\alpha,\lambda) \in (0,1)$ be suitably small constants.  In particular, assume $\varepsilon = \varepsilon(n,m,q,\alpha,\lambda,\beta) \in (0,1)$ is small enough that that the claim in Step~1 holds true, i.e.~there exists a non-zero, homogeneous degree $\alpha$, locally Dirichlet energy minimizing function $\varphi \in W^{1,2}_{\rm loc}(\mathbb{R}^n,\mathcal{A}_q(\mathbb{R}^m))$ such that $\op{dim}\op{spine}\varphi = n-2$ and \eqref{topology eqn1}--\eqref{topology eqn3} hold true.  After an orthogonal change of coordinates, assume that $\op{spine}\varphi = \{0\} \times \mathbb{R}^{n-2}$ so that $\varphi \in \Phi_{\alpha,q}$.  Replace $T$ and $\varphi$ with $\eta_{0,\lambda\#} T$ and $\lambda^{\alpha-1} \varphi$ so that (after scaling) $T$ is a locally area-minimizing rectifiable current of $\mathbf{C}_1(0)$ satisfying \eqref{SandT hyp1} and \eqref{SandT hyp2} with $\varepsilon_0/\lambda$ in place of $\varepsilon$.  

Let $\vartheta \in (0,1/4)$ such that $\overline{\nu} \vartheta^{\mu/2} < 1/2$, where $\mu = \mu(n,m,q,\alpha)$ and $\overline{\nu} = \overline{\nu}(n,m,q,\alpha)$ are as in Lemma~\ref{main excess decay lemma1}.  Let $Z = (\xi,\zeta,\chi) \in \mathcal{B}^{\star}_q \cap \mathbf{C}_{1/4}(0)$.  Set $T^{(Z)} = \widetilde{\Gamma}_{Z\#} \eta_{Z,1\#} T$, where $\widetilde{\Gamma}_Z$ is as in \eqref{SandT eqn2}.  Set $\Gamma_0 = I$ and $\varphi_0 = \varphi$ and set $s_k = \vartheta^k/4$ for each integer $k \geq 0$.  We claim that by iteratively applying Lemma~\ref{main excess decay lemma1}, we can show that: \textit{for every positive integer $k$ there exists an orthogonal rotation $\Gamma_k$ of $\mathbb{R}^{n+m}$ which fixes vectors in $\{0\} \times \mathbb{R}^m$ and $\varphi_k \in \Phi_{\alpha,q}$ such that } 
\begin{align}
	\label{topology eqn6} &\|\Gamma_k - \Gamma_{k-1}\| 
	\leq \frac{C E_{\varsigma/2,\varepsilon}(\Gamma_{k-1\#} T^{(Z)},\varphi_{k-1},\mathbf{C}_{s_{k-1}}(0))}{\widehat{E}(T^{(Z)},P_0,\mathbf{C}_{s_{k-1}}(0))} , \\
	\label{topology eqn7} &s_{k-1}^{\alpha-1} \op{dist}_{\mathcal H}(\op{spt}\op{graph}\varphi_{k-1} \cap \mathbf{C}_1(0), \op{spt}\op{graph}\varphi_k \cap \mathbf{C}_1(0)) \\&\hspace{51.2mm} \leq C E_{\varsigma/2,\varepsilon}(\Gamma_{k-1\#} T^{(Z)},\varphi_{k-1},\mathbf{C}_{s_{k-1}}(0)) , \nonumber \\
	\label{topology eqn8} &Q_{\varsigma/2,\varepsilon}(\Gamma_{k\#} T^{(Z)},\varphi_k,\mathbf{C}_{s_k}(0)) \leq \vartheta^{\alpha-1+\mu/2} E_{\varsigma/2,\varepsilon}(\Gamma_{k-1\#} T^{(Z)},\varphi_{k-1},\mathbf{C}_{s_{k-1}}(0)) , \\
	\label{topology eqn9} &\widehat{E}(T^{(Z)},P_0,\mathbf{C}_{s_k}(0)) \geq c s_k^{\alpha-1} \|\varphi_{k-1}\|_{L^2(B_1(0))} - C Q_{\varsigma/2,\varepsilon}(\Gamma_{k-1\#} T^{(Z)},\varphi_{k-1},\mathbf{C}_{s_{k-1}}(0)) ,
\end{align}
\textit{where $\mu = \mu(n,m,q,\alpha)$ is as in Lemma~\ref{main excess decay lemma1} and $c = c(n,m) > 0$ and $C = C(n,m,q,\alpha) \in (0,\infty)$ are constants.}

To prove the claim, first observe that $\varphi$ has exactly one non-zero $q_0$-valued component and recall that \eqref{SandT hyp2} hold true.  Thus we can argue as in Step~2 of the proof of Theorem~\ref{SandT thm} with $\widetilde{\varphi} = \varphi$ to show that \eqref{SandT eqn7} and \eqref{SandT eqn8} hold true and that 
\begin{equation*}
    E_{\varsigma/2,\varepsilon}(T^{(Z)},\varphi,\mathbf{C}_{1/4}(0)) \leq C E_{\varsigma/2,\varepsilon}(T,\varphi,\mathbf{C}_1(0)) ,
\end{equation*}
where $C = C(n,m,q,\alpha) \in (0,\infty)$ is a constant.  Now \eqref{topology eqn6}--\eqref{topology eqn9} follow by arguing as in Step~3 of the proof of Theorem~\ref{SandT thm} using Lemma~\ref{main excess decay lemma1} in place of Theorem~\ref{main excess decay thm}.  Let $k \geq 0$ be an integer such that either $k = 0$ or $k > 0$ and for each $i \in \{1,2,\ldots,k\}$ we have found $\Gamma_i$ and $\varphi_i$ satisfying \eqref{topology eqn6}--\eqref{topology eqn9}.  In order to apply Lemma~\ref{main excess decay lemma1} to find $\Gamma_{k+1}$ and $\varphi_{k+1}$, we need to show that for every integer $k \geq 0$ and every $y_0 \in B^{n-2}_{s_k/2}(0)$ 
\begin{equation}\label{topology eqn10}
	\widetilde{\Gamma}_Z(\mathcal{B}^{\star}_q - Z) \cap \mathbf{C}_{\delta\vartheta s_k/2}(0,y_0) \neq \emptyset , 
\end{equation}
Assuming $\beta$ is sufficiently small, by Theorem~\ref{nonconcentration thm}(A) and \eqref{SandT hyp2}, $|\chi| \leq C(n,m,q,\alpha) \,\beta$ and so $|Z + (0,y_0,0)| < 1$.  Hence by \eqref{topology hyp2} there exists $X = (x,y,z) \in \mathcal{B}^{\star}_q \cap \mathbf{B}_{1/(4\lambda)}(0)$ such that 
\begin{equation}\label{topology eqn11}
	\pi_L(X) = \pi_L(Z + (0,y_0,0)) . 
\end{equation}
It follows using Theorem~\ref{nonconcentration thm}(A) that 
\begin{equation}\label{topology eqn12}
	|x - \xi| + |z - \chi| \leq C \beta \max\{|y - \zeta|, s_k\} , 
\end{equation}
where $C = C(n,m,q,\alpha) \in (0,\infty)$ is a constant.  In particular, if $|X-Z| \geq 1/8$, then \eqref{topology eqn12} follows from \eqref{topology eqn2}.  If $s_{i+1}/2 < |X-Z| \leq s_i/2$ for some $i \in \{0,1,2,\ldots,k-1\}$, then by Theorem~\ref{nonconcentration thm}(A) with $\Gamma_{i\#} \eta_{0,s_i\#} T^{(Z)}$ and $s_i^{\alpha-1}\varphi_i$ in place of $T$ and $\varphi$ and by \eqref{SandT eqn26}, $|x-\xi| + |z-\chi| \leq C\beta s_i \leq 2\vartheta^{-1} C\beta |X-Z|$ and thus \eqref{topology eqn12} holds true.  Finally, if $|X-Z| \leq s_k/2$, we again apply Theorem~\ref{nonconcentration thm}(A) with $\Gamma_{k\#} \eta_{0,s_k\#} T^{(Z)}$ and $s_k^{\alpha-1}\varphi$ in place of $T$ and $\varphi$ together with \eqref{SandT eqn26} to similarly deduce that $|x-\xi| + |z-\chi| \leq C\beta s_k$ as in \eqref{topology eqn12}.  By \eqref{topology eqn11} and \eqref{topology eqn12} 
\begin{equation*}
	|\pi_L(0,y-\zeta-y_0,0)| = |\pi_L(x-\xi,0,z-\chi)|  
	\leq C \beta \max\{|y - \zeta|,s_k\} , 
\end{equation*}
where $C = C(n,m,q,\alpha) \in (0,\infty)$ is a constant.  Hence by \eqref{topology eqn3} 
\begin{equation}\label{topology eqn13} 
	|y - \zeta - y_0| \leq C\lambda^{-1} \beta \max\{|y - \zeta|,s_k\} 
\end{equation}
where $C = C(n,m,q,\alpha) \in (0,\infty)$ is a constant.  Note that we must have that $|y - \zeta| \leq s_k$; otherwise, if $|y - \zeta| > s_k$, then by the triangle inequality and \eqref{topology eqn13}, $|y - \zeta| \leq |y - \zeta - y_0| + |y_0| \leq C\lambda^{-1} \beta \,|y - \zeta| + |y_0|$.  Hence provided $\beta$ is small enough that $C\lambda^{-1} \beta < 1/2$ (where $C$ is as in \eqref{topology eqn13}), $|y - \zeta| \leq 2 |y_0| < s_k$, giving us a contradiction.  Now notice that by Theorem~\ref{nonconcentration thm}(A) and \eqref{SandT hyp2}, $\|\widetilde{\Gamma}_Z - I\| \leq C(n,m,q,\alpha) \,\beta$.  Thus assuming $\beta = \beta(n,m,q,\alpha,\delta)$ is sufficiently small, by $\|\widetilde{\Gamma}_Z - I\| \leq C(n,m,q,\alpha) \,\beta$, \eqref{topology eqn12}, and \eqref{topology eqn13} 
\begin{equation*}
	|\widetilde{\Gamma}_Z(X - Z) - (0,y_0,0)| \leq C\beta s_k < \delta \vartheta s_k/2 , 
\end{equation*}
where $C = C(n,m,q,\alpha) \in (0,\infty)$ is a constant.  This completes the proof of both \eqref{topology eqn10} and the claim.

Arguing as in Step~4 of the proof of Theorem~\ref{SandT thm}, we deduce that there exists an orthogonal rotation $\overline{\Gamma}_Z$ of $\mathbb{R}^{n+m}$ which fixes vectors in $\{0\} \times \mathbb{R}^m$ and there exists $\overline{\varphi}^{(Z)} \in \Phi_{\alpha,q}$ such that $\varphi^{(Z)} = \widetilde{\Gamma}_Z^{-1} \circ \overline{\varphi}^{(Z)} \circ \overline{\Gamma}_Z \circ \widetilde{\Gamma}_Z$ is the unique tangent function to $T$ at $Z$, \eqref{SandT eqn37} and \eqref{SandT eqn38} hold true, and in place of \eqref{SandT eqn39} 
\begin{equation}\label{topology eqn14}
    E_{\varsigma/2,\varepsilon}(\overline{\Gamma}_{Z\#} T^{(Z)}, \overline{\varphi}^{(Z)}, \mathbf{C}_{\rho}(0)) \leq C \rho^{\alpha-1+\mu/2} E_{\varsigma/2,\varepsilon}(T, \varphi, \mathbf{C}_1(0))
\end{equation}
for all $\rho \in (0,1/4]$, where $C = C(n,m,q,\alpha) \in (0,\infty)$ is a constant.

By \eqref{SandT hyp1}, Lemma~\ref{uniform branch set lemma}, \eqref{SandT eqn38},~\cite[Lemma~3.5 and Lemma~3.7]{KrumWica}, and \eqref{oriented tilt and mass}, for all $\sigma \in (0,1/4]$ we have that Hypothesis~$(\star)$ holds true with $C(n,m,q,\alpha)\,\varepsilon_0/\lambda, \eta_{0,\sigma\#} T^{(Z)}$ in place of $\varepsilon_0,T$.  By \eqref{SandT hyp2}, \eqref{SandT eqn37}, and \eqref{SandT eqn38}, for every integer $k \geq 0$ we have that Hypothesis~$(\star\star)$ holds true with $C(n,m,q,\alpha)\,\beta, \overline{\Gamma}_{Z\#}\eta_{0,s_k\#} T^{(Z)}, s_k^{\alpha-1} \overline{\varphi}^{(Z)}$ in place of $\beta,T,\varphi$.  Given $\rho \in (0,1/4]$, choose an integer $k \geq 0$ such that $s_{k+1} < \rho \leq s_k$ so that by Theorem~\ref{graphical_thm}(B) with $\gamma = 1/2$ and $ \tau = \vartheta/16$ and with $\overline{\Gamma}_{Z\#}\eta_{0,s_k\#} T^{(Z)}, s_k^{\alpha-1} \overline{\varphi}^{(Z)}$ in place of $T,\varphi$ together with \eqref{topology eqn14}, \eqref{SandT hyp2}, and \eqref{SandT eqn38} 
\begin{align*}
    \sup_{X \in \op{spt} \overline{\Gamma}_{Z\#} T^{(Z)} \cap \mathbf{C}_{\rho/2}(0) \cap \{r > \rho/16\}} \op{dist}(X,\op{spt}\op{graph}\overline{\varphi}^{(Z)}) 
    \leq\,& C s_k^{\alpha+\mu/2} E_{\varsigma/2,\varepsilon}(T, \varphi, \mathbf{C}_1(0))
    \\ \leq\,& C \beta \rho^{\alpha+\mu/2} \widehat{E}(T,P_0,\mathbf{C}_{\rho}(0)),
\end{align*}
where $r(X) = \op{dist}(\pi(X),\{0\} \times \mathbb{R}^{n-2})$ and where $C = C(n,m,q,\alpha) \in (0,\infty)$ are constants.  It follows that Hypothesis~$(\star\star)$ holds true with $C(n,m,q,\alpha)\,\beta, \overline{\Gamma}_{Z\#}\eta_{0,\rho\#} T^{(Z)}, \rho^{\alpha-1} \overline{\varphi}^{(Z)}$ in place of $\beta,T,\varphi$.

\noindent\textit{Step~3:  Parametric representations of the branch set and current near the origin.}  In light of \eqref{topology eqn10}, conclusion~(a) of the theorem holds true as $\mathcal{B}^{\star}_q \cap \mathbf{C}_{1/4}(0) = \{(g(y),y,h(y)) : y \in B^{n-2}_{1/4}(0) \} \cap \mathbf{C}_{1/4}(0)$ for some $g \in C^{1,\mu/2}(B^{n-2}_{1/4}(0),\mathbb{R}^2)$, and $h \in C^{1,\mu/2}(B^{n-2}_{1/4}(0),\mathbb{R}^m)$ such that $g(0) = 0$, $Dg(0) = 0$, $h(0) = 0$, $Dh(0) = 0$, $[Dg]_{\mu/2,B^{n-2}_{1/4}(0)} + [Dh]_{\mu/2,B^{n-2}_{1/4}(0)} \leq C(n,m,q,\alpha) \,\beta$.  (Note that it is clear from Theorem~\ref{graphical_thm} and the discussion below that $\op{sing} T \cap \mathbf{C}_{1/4}(0) = \mathcal{B}^{\star}_q \cap \mathbf{C}_{1/4}(0)$.)

Let $Z = (\xi_0,\zeta_0,\chi_0) \in \mathcal{B}^{\star}_q \cap \mathbf{C}_{1/8}(0)$.  We want to argue as in the proof of Theorem~\ref{graphical_thm} (with $\gamma = 1/2$) to express $T^{(Z)}$ as the graph of a $q$-valued function over domains in $P_Z$ which lie away from the spine of $\varphi^{(Z)}$.  By Theorem~\ref{graphical_thm}(A) with $\overline{\Gamma}_{Z\#}\eta_{0,1/4\#} T^{(Z)}, 4^{1-\alpha} \overline{\varphi}^{(Z)}$ in place of $T,\varphi$ 
\begin{equation}\label{topology eqn15}
	\inf_{x \in \mathbb{S}^1} \op{sep} \overline{\varphi}^{(Z)}(x,0) \geq C \widehat{E}(T,P_0,\mathbf{C}_1(0)) 
\end{equation}
for some constant $C = C(n,m,q,\alpha) > 0$ and \eqref{graphical conclA2} holds true with $\overline{\varphi}^{(Z)}$ in place of $\varphi$.  Let $\sigma \in (0,1/8]$, $\rho \geq \sigma^{1+\mu/(8n)}/4$, and $\zeta \in \mathbb{R}^{n-2}$ with $\rho^2 + |\zeta|^2 \leq (7\sigma/8)^2$.  For each $\kappa \in (0,1]$, let $A_{\rho,\kappa}(\zeta)$ be as in \eqref{annuli_defn} with $\gamma = 1/2$ and let $\mathbf{A}_{\rho,\kappa}(\zeta) = A_{\rho,\kappa}(\zeta) \times \mathbb{R}^m$.  By \eqref{topology eqn14} and $\rho \geq \sigma^{1+\mu/(8n)}/4$ 
\begin{align}\label{topology eqn16}
	\rho \,E_{\varsigma/2,\varepsilon}(\overline{\Gamma}_{Z\#} T^{(Z)}, \overline{\varphi}^{(Z)}, \mathbf{A}_{\rho,1}(\zeta))
	\leq\,& C (\rho/\sigma)^{-n/2} \sigma \,E_{\varsigma/2,\varepsilon}(\overline{\Gamma}_{Z\#} T^{(Z)},\overline{\varphi}^{(Z)},\mathbf{C}_{\sigma}(0)) 
	\\ \leq\,& C \rho^{-n/2} \sigma^{n/2+\alpha+\mu/2} E_{\varsigma/2,\varepsilon}(T,\varphi,\mathbf{C}_1(0)) \nonumber 
	\\ \leq\,& C \sigma^{\alpha+7\mu/16} \nonumber E_{\varsigma/2,\varepsilon}(T,\varphi,\mathbf{C}_1(0)) , \nonumber 
\end{align}
where $C = C(n,m,q,\alpha) \in (0,\infty)$ are constants.  Since Hypothesis~$(\star)$ holds true with $C\varepsilon_0/\lambda$ and $\eta_{0,\sigma\#} T^{(Z)}$ in place of $\varepsilon_0$ and $T$ and by \eqref{topology eqn16}, we can use Lemma~\ref{graphical_lemma0} to obtain 
\begin{align*}
	\sup_{X \in \op{spt} \overline{\Gamma}_{Z\#} T^{(Z)} \cap \mathbf{A}_{\rho,3/4}(\zeta)} \op{dist}(X,\op{spt}\op{graph}\overline{\varphi}^{(Z)}) 
	\leq\,& C \rho \,E_{\varsigma/2,\varepsilon}(T^{(Z)}, \overline{\varphi}^{(Z)}, \mathbf{A}_{\rho,1}(\zeta))
	\\ \leq\,& C \sigma^{\alpha+7\mu/16} E_{\varsigma/2,\varepsilon}(T,\varphi,\mathbf{C}_1(0)) ,
\end{align*}
where $C = C(n,m,q,\alpha) \in (0,\infty)$ are constants.  Since the only component of $\varphi^{(Z)}$ being multiplicity one, we can use the Allard regularity theorem (in place of Almgren's strong Lipschitz approximation theorem) to deduce that there exists $\widetilde{v}^{(Z,\rho,\zeta)} : \op{graph}\varphi^{(Z)} \cap \mathbf{A}_{\rho,1/2}(\zeta) \rightarrow \mathcal{A}_q(\mathbb{R}^m)$ and $\widetilde{u}^{(Z,\rho,\zeta)} : A_{\rho,1/2}(\zeta) \rightarrow \mathcal{A}_q(\mathbb{R}^m)$ such that for each open ball $B \subset A_{\rho,1/2}(\zeta)$ 
\begin{equation}\label{topology eqn17}
	\widetilde{u}^{(Z,\rho,\zeta)}(x,y) = \sum_{l=1}^q \llbracket \overline{\varphi}^{(Z)}_l(x,y) + \widetilde{v}^{(Z,\rho,\zeta)}_l(x,y) \rrbracket 
\end{equation}
for each $(x,y) \in B$, where $\overline{\varphi}^{(Z)}(X) = \sum_{l=1}^q \llbracket \overline{\varphi}^{(Z)}_l(x,y) \rrbracket$ for single-valued harmonic functions $\overline{\varphi}^{(Z)}_l : B \rightarrow \mathbb{R}^m$ and $\widetilde{v}^{(Z,\rho,\zeta)}_l : B \rightarrow \mathbb{R}^m$ is a smooth function given by $\widetilde{v}^{(Z,\rho,\zeta)}_l(x,y) = \widetilde{v}^{(Z,\rho,\zeta)}(x,y,\overline{\varphi}^{(Z)}_l(x,y))$ for each $(x,y) \in B$, and 
\begin{gather}
	\label{topology eqn18} (\overline{\Gamma}_{Z\#} T^{(Z)}) \llcorner \mathbf{A}_{\rho,1/2}(\zeta) = \op{graph} \widetilde{u}^{(Z,\rho,\zeta)} , \\
	\label{topology eqn19} \sup_{(x,y) \in A_{\rho,1/2}(\zeta)} \sum_{l=1}^q |\widetilde{v}^{(Z,\rho,\zeta)}_l(x,y)| 
		\leq C \sigma^{\alpha+7\mu/16} E_{\varsigma/2,\varepsilon}(T,\varphi,\mathbf{C}_1(0)) , 
\end{gather}
where $C = C(n,m,q,\alpha) \in (0,\infty)$ is a constant.  

Let $B \subset A_{\rho,1/2}(\zeta)$ be an arbitrary ball of radius $\rho/128$.  In $B$, $\widetilde{u}^{(Z,\rho,\zeta)}_l = \overline{\varphi}^{(Z)}_l + \widetilde{v}^{(Z,\rho,\zeta)}_l$ solves the minimal surface system 
\begin{equation}\label{topology eqn20}
	D_i \bigg( \sqrt{G(D\widetilde{u}^{(Z,\rho,\zeta)}_l)} \,G^{ij}(D\widetilde{u}^{(Z,\rho,\zeta)}_l) \,D_j \widetilde{u}^{(Z,\rho,\zeta)}_l \bigg) = 0 \text{ in } B , 
\end{equation}
where $G_{ij}(D\widetilde{u}^{(Z,\rho,\zeta)}_l) = \delta_{ij} + D_i \widetilde{u}^{(Z,\rho,\zeta)}_l \cdot D_j \widetilde{u}^{(Z,\rho,\zeta)}_l$ for all $1 \leq i,j \leq n$, $[G^{ij}(D\widetilde{u}^{(Z,\rho,\zeta)}_l)]$ is the inverse matrix of $[G_{ij}(D\widetilde{u}^{(Z,\rho,\zeta)}_l)]$, and $G(D\widetilde{u}^{(Z,\rho,\zeta)}_l)$ is the determinant of $[G_{ij}(D\widetilde{u}^{(Z,\rho,\zeta)}_l)]$.  By \eqref{SandT eqn38}, standard estimates for the minimal surface system, and the arbitrariness of the ball $B$
\begin{equation}\label{topology eqn21}
	\rho^j \sup_{(x,y) \in A_{\rho,3/8}(\zeta)} |D^j \widetilde{u}^{(Z,\rho,\zeta)}(x,y)| 
    \leq C \sigma^{\alpha} \widehat{E}(T,P_0,\mathbf{C}_1(0)) 
\end{equation}
for each integer $j \geq 0$, where $C = C(n,m,q,\alpha,j) \in (0,\infty)$ is a constant.  
Now let $B \subset A_{\rho,3/8}(\zeta)$ be an arbitrary ball of radius $\rho/128$.  Since $\overline{\varphi}^{(Z)}_l$ is harmonic in $B$ and $\widetilde{u}^{(Z,\rho,\zeta)}_l = \overline{\varphi}^{(Z)}_l + \widetilde{v}^{(Z,\rho,\zeta)}_l$ satisfies \eqref{topology eqn20}, 
\begin{equation*}
	\Delta\widetilde{v}^{(Z,\rho,\zeta)}_l = -D_i \bigg( \bigg( \sqrt{G(D\widetilde{u}^{(Z,\rho,\zeta)}_l)} 
		\,G^{ij}(D\widetilde{u}^{(Z,\rho,\zeta)}_l) - \delta_{ij} \bigg) \,D_j \widetilde{u}^{(Z,\rho,\zeta)}_l \bigg) \text{ in } B .
\end{equation*}
Provided $\mu \leq 1/q$ (so that $\alpha \geq 1+\mu$), by \eqref{topology eqn21} and the arbitrariness of the ball $B$ 
\begin{align*}
    \rho^{2+j} \sup_{(x,y) \in A_{\rho,3/8}(\zeta)} \sum_{l=1}^q |D^j \Delta \widetilde{v}^{(Z,\rho,\zeta)}_l(x,y)| \leq\,& C \sigma^{3\alpha} \rho^{-2} \widehat{E}(T,P_0,\mathbf{C}_1(0))^3 
    \\ \leq\,& C \sigma^{\alpha+7\mu/4} \widehat{E}(T,P_0,\mathbf{C}_1(0))^3 
\end{align*}
for each $j \in \{0,1,2\}$, where $C = C(n,m,q,\alpha) \in (0,\infty)$ is a constant.  Hence by the Schauder estimates and \eqref{topology eqn19} 
\begin{align}\label{topology eqn22}
	\sup_{(x,y) \in A_{\rho,1/4}(\zeta)} \sum_{l=1}^q |D^j \widetilde{v}^{(Z,\rho,\zeta)}_l(x,y)| 
    \leq\,& C \sigma^{\alpha+7\mu/16} \rho^{-j} E_{\varsigma/2,\varepsilon}(T,\varphi,\mathbf{C}_1(0)) 
    \\ \leq\,& C \sigma^{\alpha+\mu/8-j} E_{\varsigma/2,\varepsilon}(T,\varphi,\mathbf{C}_1(0)) \nonumber 
\end{align}
for each $j \in \{0,1\}$, where $C = C(n,m,q,\alpha) \in (0,\infty)$ is a constant.  

Next we want to express $T$ as the graph of a $q$-valued function $u$ over $B_{1/16}(0) \subset \mathbb{R}^n$.  By Theorem~\ref{nonconcentration thm}(A) and \eqref{topology eqn14}, for each $Z \in \mathcal{B}^{\star}_q \cap \mathbf{C}_{1/8}(0)$ let $P_Z = \{ (x,y,M^{(Z)}(x,y)) : (x,y) \in \mathbb{R}^n \}$ for some $m \times n$ matrix $M^{(Z)}$ such that $M^{(0)} = 0$ and 
\begin{equation}\label{topology eqn23}
   |M^{(Y)} - M^{(Z)}| \leq C |Y-Z|^{\alpha-1+\mu/2} E_{\varsigma/2,\varepsilon}(T, \varphi, \mathbf{C}_1(0))
\end{equation}
for all $Y,Z \in \mathcal{B}^{\star}_q \cap \mathbf{C}_{1/8}(0)$, where $C = C(n,m,q,\alpha) \in (0,\infty)$ is a constant.  Let $Z = (\xi_0,\zeta_0,\chi_0) \in \mathcal{B}^{\star}_q \cap \mathbf{C}_{1/8}(0)$ and $\rho \in (0,7/64]$ (as above) and set $\zeta = 0$ and $\sigma = 8\rho/7$.  By \eqref{topology eqn18}, \eqref{topology eqn21}, and \eqref{topology eqn23}, there exists $u^{(Z,\rho,0)} : (\xi_0,\zeta_0) + A_{\rho,1/8}(0) \rightarrow \mathcal{A}_q(\mathbb{R}^m)$ such that 
\begin{equation*} 
	T \llcorner (Z + \mathbf{A}_{\rho,1/8}(0)) 
	= \eta^{-1}_{Z,1\#}((\widetilde{\Gamma}_{Z\#}^{-1} \overline{\Gamma}_{Z\#}^{-1}(\op{graph} \widetilde{u}^{(Z,\rho,0)})) \llcorner \mathbf{A}_{\rho,1/8}(0))
	= \op{graph} u^{(Z,\rho,0)} 
\end{equation*}
and 
\begin{align}\label{topology eqn24} 
	&\sup_{(x,y) \in (\xi_0,\zeta_0) + A_{\rho,1/8}(0)} |u^{(Z,\rho,0)}(x,y) - \chi_0 - M^{(Z)} (x-\xi_0,y-\zeta_0)| 
	\\&+ \rho \sup_{(x,y) \in (\xi_0,\zeta_0) + A_{\rho,1/8}(\zeta)} |Du^{(Z,\rho,0)}(x,y) - M^{(Z)}| \nonumber 
	\\&+ \rho^2 \sup_{(x,y) \in (\xi_0,\zeta_0) + A_{\rho,1/8}(\zeta)} |D^2 u^{(Z,\rho,0)}(x,y)| \leq C \rho^{\alpha} \widehat{E}(T,P_0,\mathbf{C}_1(0)) \nonumber
\end{align}
where $u^{(Z,\rho,0)}_a$ denotes the single-valued average of $u^{(Z,\rho,0)}$ and $C = C(n,m,q,\alpha) \in (0,\infty)$ is a constant.  (More generally, if $\sigma \in (0,1/8]$, $\rho \geq \sigma^{1+\mu/(8n)}/4$, and $\zeta \in \mathbb{R}^{n-2}$ with $\rho^2 + |\zeta|^2 \leq (7\sigma/8)^2$, then we can similarly express $\eta^{-1}_{Z,1\#} \widetilde{\Gamma}_{Z\#}^{-1} \overline{\Gamma}_{Z\#}^{-1}(\op{graph} \widetilde{u}^{(Z,\rho,0)})$ as the graph of a $q$-valued function $u^{(Z,\rho,\zeta)}$ over the domain $(\xi_0,\zeta_0) + \pi_{P_0} \widetilde{\Gamma}_Z^{-1} \overline{\Gamma}_Z^{-1}(0,\zeta) + A_{\rho,1/8}(0)$.)  Define $u : B_{1/16}(0) \rightarrow \mathcal{A}_q(\mathbb{R}^m)$ by $u = u^{(Z,\rho,0)}$ on $B_{1/16}(0) \cap ((\xi_0,\zeta_0) + A_{\rho,1/8}(0))$ for all $Z = (\xi_0,\zeta_0,\chi_0) \in \mathcal{B}^{\star}_q \cap \mathbf{C}_{1/8}(0)$ and $\rho \in (0,7/64]$ with $B_{1/16}(0) \cap ((\xi_0,\zeta_0) + A_{\rho,1/8}(0)) \neq \emptyset$ so that 
\begin{equation}\label{topology eqn25} 
	T \llcorner \mathbf{C}_{1/16}(0) = \op{graph} u .
\end{equation}
Given any $Z = (\xi_0,\zeta_0,\chi_0) \in \mathcal{B}^{\star}_q \cap \mathbf{C}_{1/8}(0)$ and $(x,y) \in B_{1/16}(0)$ with $|y - \zeta_0| < |x-\xi_0|/512$ and $|x-\xi_0| < 7/64$, by setting $\rho = |x-\xi_0|$ we have that $(x,y) \in (\xi_0,\zeta_0) + A_{\rho,1/8}(0)$ and thus by \eqref{topology eqn24}  
\begin{align}\label{topology eqn26} 
	&|u(x,y) - \chi_0 - M^{(Z)} (x-\xi_0,y-\zeta_0)| + |x-\xi_0| \,|Du(x,y) - M^{(Z)}|  
	\\&+ |x-\xi_0|^2 |D^2 u(x,y)| \leq C |x-\xi_0|^{\alpha} \widehat{E}(T,P_0,\mathbf{C}_1(0)) , \nonumber
\end{align}
where $C = C(n,m,q,\alpha) \in (0,\infty)$ is a constant.

It follows using \eqref{topology eqn23} and \eqref{topology eqn26} that $u \in C^{1,\min\{1,\alpha-1\}}(B_{1/16}(0),\mathcal{A}_q(\mathbb{R}^m))$.  To see this, let's first show that $u$ is differentiable at each $(g(y_0),y_0) \in B_{1/16}(0)$, where $y_0 \in B^{n-2}_{1/16}(0)$.  Set $Z_0 = (g(y_0),y_0,h(y_0)) \in \mathcal{B}^{\star}_q \cap \mathbf{C}_{1/16}(0)$ and let $(x,y) \in B_{1/16}(0)$ with $|(x,y) - (y_0,g(y_0))| < 1/32$.  If $|y - y_0| < |x - g(y_0)|/512$, then by \eqref{topology eqn26} with $Z = Z_0$ 
\begin{equation*}
    |u(x,y) - h(y_0) - M^{(Z_0)} (x-g(y_0),y-y_0)| \leq C \widehat{E} |x - g(y_0)|^{\alpha} , 
\end{equation*}
where $\widehat{E} = \widehat{E}(T,P_0,\mathbf{C}_1(0))$ and $C = C(n,m,q,\alpha) \in (0,\infty)$ is a constant.  If instead $|y - y_0| \geq |x - g(y_0)|/512$, then by Theorem~\ref{nonconcentration thm}(A) and \eqref{topology eqn14}, the distance of $Z = (g(y),y,h(y))$ to the tangent plane to $T$ at $Z_0 = ((g(y_0),y_0,h(y_0))$ is $\leq C E |y -y_0|^{\alpha+\mu/2}$, i.e.
\begin{equation}\label{topology eqn27}
    |h(y) - h(y_0) - M^{(Z_0)} (g(y)-g(y_0),y-y_0)| 
    \leq C E |y - y_0|^{\alpha+\mu/2} , 
\end{equation}
where $E = E_{\varsigma/2,\varepsilon}(E(T,\varphi,\mathbf{C}_1(0))$ and $C = C(n,m,q,\alpha) \in (0,\infty)$ are constants.  Hence by \eqref{topology eqn23}, \eqref{topology eqn26}, and \eqref{topology eqn27} 
\begin{align}\label{topology eqn28}
    &|u(x,y) - h(y_0) - M^{(Z_0)} (x-g(y_0),y-y_0)| 
    \\ \leq\,& |u(x,y) - h(y) - M^{(Z)} (x-g(y),0)| + |h(y) - h(y_0) + M^{(Z_0)} (g(y)-g(y_0),y-y_0)| \nonumber \\&+ |M^{(Z)} - M^{(Z_0)}| \,|x-g(y)| \nonumber 
    \\ \leq\,& C \widehat{E} |x-g(y)|^{\alpha} + C E |y-y_0|^{\alpha+\mu/2} + C E |y-y_0|^{\alpha-1+\mu/2} |x-g(y)| \nonumber 
    \\ \leq\,& C \widehat{E} |(x-g(y_0),y-y_0)|^{\alpha} , \nonumber
\end{align}
where $C = C(n,m,q,\alpha) \in (0,\infty)$ are constants.  Therefore, $u$ is differentiable at $(g(y_0),y_0)$ with $u(g(y_0),y_0) = q \llbracket h(y_0) \rrbracket$ and $Du(g(y_0),y_0) = q \llbracket M^{(Z_0)} \rrbracket$.  To bound the H\"older coefficient of $Du$, let $(x,y),(x',y') \in B_{1/16}(0)$ with $|(x,y) - (x',y')| < 1/32$.  Assume that $|x' - g(y')| \leq |x - g(y)|$.  If $|(x,y) - (x',y')| < |x - g(y)|/4$, then by \eqref{topology eqn26} with $Z = (g(y),y,h(y))$
\begin{equation*}
	\mathcal{G}(Du(x,y),Du(x',y')) \leq \sup_{B_{2d}(x,y)} |D^2 u| \,|(x,y) - (x',y')| \leq C \widehat{E} \,|(x,y) - (x',y')|^{\min\{1,\alpha-1\}} , 
\end{equation*}
where $d = |(x,y) - (x',y')|$ and $C = C(n,m,q,\alpha) \in (0,\infty)$ is a constant.  If instead $|(x,y) - (x',y')| \geq |x - g(y)|/4$, then by \eqref{topology eqn23} and \eqref{topology eqn26} 
\begin{align*}
	\mathcal{G}(Du(x,y),Du(x',y')) \leq\,& |Du(x,y) - M^{(Z)}| + q |M^{(Z)} - M^{(Z')}| + |M^{(Z')} - Du(x',y')| \\ \leq\,& C \widehat{E} \,|x-g(y)|^{\alpha-1} + C E |y - y'|^{\alpha-1+\mu/2} + C \widehat{E} \,|x'-g(y')|^{\alpha-1} \\ \leq\,& C \widehat{E} \,|(x,y) - (x',y')|^{\alpha-1} , 
\end{align*}
where $Z = (g(y),y,h(y))$, $Z = (g(y'),y,h(y'))$, and $C = C(n,m,q,\alpha) \in (0,\infty)$ are constants.  Therefore, $[Du]_{\min\{1,\alpha-1\},B_{1/16}(0)} \leq C \widehat{E}$, where $C = C(n,m,q,\alpha) \in (0,\infty)$ is a constant.

By the structure of $\widetilde{u}^{(Z,\rho,0)}$ from \eqref{topology eqn17} (with $\zeta = 0$) together with \eqref{SandT hyp2}, Remark~\ref{graphical rmk}(2), \eqref{topology eqn15}, and \eqref{topology eqn22}, there exists $f \in C^{1,\mu/2}(B^2_{1/32^{1/q}}(0) \times B^{n-2}_{1/32}(0),\mathbb{R}^m)$ such that 
\begin{equation}\label{topology eqn29} 
	u(g(y) + x,y) = \sum_{z^q = x} \llbracket f(z,y) \rrbracket 
\end{equation}
for all $x \in B^2_{1/32}(0)$ and all $y \in B^{n-2}_{1/32}(0)$ and thus by \eqref{topology eqn25}, conclusion~(b) holds true.  In particular, $f(\nu z,y) \neq f(z,y)$ for every $q^{\rm th}$ root of unity $\nu \neq 1$, $z \in \mathbf{C}$ with $|z|^q < 1/32$, and $y \in B^{n-2}_{1/32}(0)$, and thus the mapping $(z,y) \mapsto (g(y)+z^q,y,f(z,y))$ is injective.

Let $\varphi^{(0)}$ be the unique tangent function to $T$ at the origin and express $\varphi^{(0)}$ as $\varphi^{(0)}(x_1,x_2,y) = \op{Re}(c (x_1+ix_2)^{p/q})$ for each $x_1,x_2 \in \mathbb{R}$ and $y \in \mathbb{R}^{n-2}$, where $c \in \mathbb{C}^m \setminus \{0\}$.  Recall from~\cite[Theorem~3.2]{MicWhi95} that $c \cdot c = 0$.  In the special case $Z = 0$, $u = \widetilde{u}^{(0,\rho,\zeta)}$ on $B_{1/16}(0) \cap A_{\rho,1/8}(\zeta)$ for all $\sigma \in (0,1/8]$, $\rho \geq \sigma^{1+\mu/(8n)}/4$, and $\zeta \in \mathbb{R}^{n-2}$ with $\rho^2 + |\zeta|^2 \leq (7\sigma/8)^2$ and $B_{1/16}(0) \cap A_{\rho,1/8}(\zeta) \neq \emptyset$.  Hence setting $\Omega = \{ (x,y) \in B_{1/16}(0) : |x| \geq |(x,y)|^{1+\mu/(8n)} \}$, by \eqref{topology eqn17} and \eqref{topology eqn22} (with $\rho = |x|$, $\sigma = 2 |(x,y)|$, and $\zeta = y$) there exists $v : \op{graph}\varphi^{(0)} \cap \Omega \times \mathbb{R}^m \rightarrow \mathcal{A}_q(\mathbb{R}^m)$ such that for each open ball $B \subset \Omega$ 
\begin{equation}\label{topology eqn30} 
	u(x,y) = \sum_{l=1}^q \llbracket \varphi^{(0)}_l(x,y) + v_l(x,y) \rrbracket  
\end{equation}
for each $(x,y) \in B$, where $\varphi^{(0)}(X) = \sum_{l=1}^q \llbracket \varphi^{(0)}_l(x,y) \rrbracket$ for single-valued harmonic functions $\varphi^{(0)}_l : B \rightarrow \mathbb{R}^m$ and $v_l : B \rightarrow \mathbb{R}^m$ is a smooth function given by $v_l(x,y) = v(x,y,\varphi_l(x,y))$ for each $(x,y) \in B$, and 
\begin{equation}\label{topology eqn31} 
	\sum_{l=1}^q |D^j v_l(x,y)| \leq C |(x,y)|^{\alpha+\mu/8-j} E_{\varsigma/2,\varepsilon}(T,\varphi,\mathbf{C}_1(0))  , 
\end{equation}
for all $(x,y) \in \Omega$ and all $j \in \{0,1\}$, where $C = C(n,m,q,\alpha) \in (0,\infty)$ is a constant.  Suppose that $(x,y) \in B_{1/16}(0)$ and $z \in \mathbb{C}$ with $|x| \geq |(x,y)|^{1+\mu/(8n)}$ and $x = g(y) + z^q$.  By \eqref{topology eqn30}, after possibly replacing $f(z,y)$ with $f(\nu z,y)$ for an appropriate $q^{\rm th}$ root of unity $\nu$, \eqref{topology concl3} holds true with $f(z,y) = \op{Re}(cz^p) + v(g(y)+z^q,y,\op{Re}(cz^p))$ and $e(z,y) = v(g(y)+z^q,y,\op{Re}(cz^p))$.  By the chain rule and \eqref{topology eqn31}
\begin{align*}
	|e(z,y)| \leq\,& \sum_{l=1}^q |v_l(x,y)| \leq C E (|z|^q + |y|)^{p/q+\mu/8} , \\
	|D_z e(z,y)| \leq\,& \sum_{l=1}^q |D_x v_l(x,y)| \cdot 2q |z|^{q-1} \leq C E |z|^{q-1} (|z|^q + |y|)^{p/q-1+\mu/8} , \\
	|D_y e(z,y)| \leq\,& \sum_{l=1}^q (|D_x v_l(x,y)| \,|Dg(y)| + |D_y v_l(x,y)|) \leq C E (|z|^q + |y|)^{p/q-1+\mu/8} , 
\end{align*}
where $E = E_{\varsigma/2,\varepsilon}(T,\varphi,\mathbf{C}_1(0))$ and $C = C(n,m,q,\alpha) \in (0,\infty)$ are constants.  Hence \eqref{topology concl4} holds true.  Suppose instead that $(x,y) \in B_{1/16}(0)$ and $z \in \mathbb{C}$ with $|x| < |(x,y)|^{1+\mu/(8n)}$ and $x = g(y) + z^q$.  By Remark~\ref{graphical rmk}(2)
\begin{align*}
	|{\rm Re}(cz^p)| \leq\,& C \widehat{E} |z|^p \leq C \widehat{E} |(x,y)|^{p/q+\mu/(8n)} , \\ 
    |D_z ({\rm Re}(cz^p))| \leq\,& C \widehat{E} |z|^{p-1} \leq C \widehat{E} |(x,y)|^{(p-1)/q+\mu/(8n)} ,
\end{align*}
where $C = C(n,m,q,\alpha) \in (0,\infty)$ are constants.  By Theorem~\ref{nonconcentration thm}(A) and \eqref{topology eqn14}, $|h(y)| \leq C(n,m,q,\alpha) \,E |y|^{\alpha+\mu/2}$.  By the $C^{1,\mu/2}$-estimates for $g$, $|g(y)| \leq C\beta |y|^{1+\mu/2}$ and thus $|x - g(y)| \leq C |(x,y)|^{1+\mu/(8n)}$, where $C = C(n,m,q,\alpha) \in (0,\infty)$ are constants.  Hence setting $Z = (g(y),y,h(y))$ and using $\alpha \geq 1+1/q$, \eqref{topology eqn23}, and \eqref{topology eqn26} 
\begin{align*}
	|u(x,y)| \leq\,& |u(x,y) - h(y) - M^{(Z)} (x-g(y),0)| + |h(y)| + |M^{(Z)}| \,|x - g(y)|
        \\ \leq\,& C \widehat{E} |x - g(y)|^{\alpha} + C E |y|^{\alpha+\mu/2} + C E |y|^{\alpha-1+\mu/2} |x - g(y)|
        \leq C \widehat{E} |(x,y)|^{\alpha+\mu/(8n)} , \\
	|Du(x,y)| \leq\,& |Du(x,y) - M^{(Z)}| +  |M^{(Z)}|
        \\ \leq\,& C \widehat{E} |x-g(y)|^{\alpha-1} + C E |y|^{\alpha-1+\mu/2} 
		\leq C \widehat{E} |(x,y)|^{\alpha-1+\mu/(8nq)} , 
\end{align*}
where $C = C(n,m,q,\alpha) \in (0,\infty)$ are constants.  Again using the chain rule for $x = g(y) + z^q$, it follows that \eqref{topology concl3} and \eqref{topology concl4} hold true, completing the proof of conclusion~(c).  By conclusion~(c), conclusion~(d) holds true with $a = (1 - \nu^p)\,c$.  Note that since $p,q$ are relatively prime, $a \neq 0$.
\end{proof}

\begin{remark}\label{topology pf rmk1}{\rm
It in fact follows from \eqref{topology eqn21} and \eqref{topology eqn23} that $u \in C^{\lfloor\alpha\rfloor,\alpha-\lfloor\alpha\rfloor}(B_{1/16}(0),\mathcal{A}_q(\mathbb{R}^m))$ with $D^j u = q \llbracket 0 \rrbracket$ on $\pi(\mathcal{B}^*_q)$ for all $j \in \{2,\ldots,\lceil \alpha \rceil - 1\}$.  Explicitly, this means the following: 
\begin{enumerate}[itemsep=2mm,topsep=0mm]
    \item[(i)]  $u |_{B_{1/16}(0) \setminus \pi(\mathcal{B}^*_q)}$ locally consists of $q$ single-valued functions which are locally in $C^{\lfloor\alpha\rfloor,\alpha-\lfloor\alpha\rfloor}$ (and in fact are smooth); 
    \item[(ii)]  for each $j \in \{1,2,\ldots,\lfloor\alpha\rfloor-1\}$, $D^j u$ extends to a function on $B_{1/6}(0)$ which is differentiable at each $(x_0,y_0) \in \pi(\mathcal{B}^*_q) \cap B_{1/6}(0)$ with $D^{j+1} u(x_0,y_0) = q \llbracket 0 \rrbracket$ (and thus $D^j u$ is continuous on $B_{1/16}(0)$ respect to the metric $\mathcal{G}$); 
    \item[(iii)]  $D^{\lfloor\alpha\rfloor} u$ is (with respect to the metric $\mathcal{G}$) a $C^{0,\alpha-\lfloor\alpha\rfloor}$-function on $B_{1/16}(0)$.
\end{enumerate}

To see this, we first show that 
\begin{align}\label{topology eqn32} 
	&|u(x,y) - h(y_0) - M^{(Z_0)} (x-g(y_0),y-y_0)| + |(x-g(y_0),y-y_0)| \,|Du(x,y) - M^{(Z_0)}|  
	\\&+ \sum_{j=2}^{\lfloor\alpha\rfloor+1} |(x-g(y_0),y-y_0)|^j \,|D^j u(x,y)| \leq C \widehat{E} \,|(x-g(y_0),y-y_0)|^{\alpha} , \nonumber
\end{align}
for each $Z_0 = (g(y_0),y_0,h(y_0)) \in \mathcal{B}^*_q \cap \mathbf{C}_{1/8}(0)$ and $(x,y) \in B_{1/16}(0)$ with $|x-g(y_0)| < 7/64$, where $\widehat{E} = \widehat{E}(T,P_0,\mathbf{C}_1(0))$ and $C = C(n,m,q,\alpha) \in (0,\infty)$ is a constant.  By the argument leading to \eqref{topology eqn26}, if $|y-y_0| < |x-g(y_0)|/512$, then \eqref{topology eqn32} follows from \eqref{topology eqn21}, which holds true for every integer $j \geq 0$.  If on the other hand $|y-y_0| \geq |x-g(y_0)|/512$, \eqref{topology eqn21} follows from \eqref{topology eqn28} and using \eqref{topology eqn23} and \eqref{topology eqn32} (with $Z = (g(y),y,h(y))$ in place of $Z_0$) to obtain 
\begin{align*}
    |Du(x,y) - M^{(Z_0)}| 
    \leq\,& |Du(x,y) - M^{(Z)}| + |M^{(Z)} - M^{(Z_0)}| \\ \leq\,& C \widehat{E} \,|x - g(y)|^{\alpha-1} + C E |y-y_0|^{\alpha-1+\mu/2}
    \\ \leq\,& C \widehat{E} \,|(x-g(y_0),y-y_0)|^{\alpha-1} 
\end{align*}
and 
\begin{equation*}
    |D^j u(x,y)| \leq C \widehat{E} \,|x-g(y)|^{\alpha-j} \leq C \widehat{E} \,|(x-g(y_0),y-y_0)|^{\alpha-j}
\end{equation*}
for all $j \in \{2,3,\ldots,\lfloor\alpha\rfloor-1\}$, where $E = E_{\varsigma/2,\varepsilon}(T,\varphi,\mathbf{C}_1(0))$ and $C = C(n,m,q,\alpha) \in (0,\infty)$ are constants.  It then follows from \eqref{topology eqn32} that for all $j \in \{1,2,\ldots,\lfloor\alpha\rfloor-1\}$, $D^j u$ is differentiable at each $(g(y_0),y_0) \in \pi(\mathcal{B}^*_q) \cap B_{1/16}(0)$ with $D^j u(g(y_0),y_0) = q \llbracket 0 \rrbracket$, where $y_0 \in B^{n-2}_{1/16}(0)$.  To bound the H\"older coefficient of $D^{\lfloor\alpha\rfloor} u$ in the case that $\alpha > 2$, let $(x,y),(x',y') \in B_{1/16}(0)$ with $|(x,y) - (x',y')| < 1/32$, and assume that $|x'-g(y')| \leq |x-g(y)|$.  If $|(x,y) - (x',y')| < |x-g(y)|/4$, then by \eqref{topology eqn32} with $Z = (g(y),y,h(y))$
\begin{align*}
    \mathcal{G}(D^{\lfloor\alpha\rfloor} u(x,y),D^{\lfloor\alpha\rfloor} u(x',y')) \leq\,& \sup_{B_{2d}(x,y)} |D^{\lfloor\alpha\rfloor+1} u| \,|(x,y) - (x',y')| \\ \leq\,& C \widehat{E} \,|(x,y) - (x',y')|^{\alpha-\lfloor\alpha\rfloor} , 
\end{align*}
where $d = |(x,y) - (x',y')|$ and $C = C(n,m,q,\alpha) \in (0,\infty)$ is a constant.  If instead $|(x,y) - (x',y')| \geq |x-g(y)|/4$, then by \eqref{topology eqn32} 
\begin{align*}
    \mathcal{G}(D^{\lfloor\alpha\rfloor} u(x,y),D^{\lfloor\alpha\rfloor} u(x',y')) \leq\,& |D^{\lfloor\alpha\rfloor} u(x,y)| + |D^{\lfloor\alpha\rfloor} u(x',y')| \\ \leq\,& C \widehat{E} \,|x-g(y)|^{\alpha-\lfloor\alpha\rfloor} + C \widehat{E} \,|x'-g(y')|^{\alpha-\lfloor\alpha\rfloor} \\ \leq\,& C \widehat{E} \,|(x,y) - (x',y')|^{\alpha-\lfloor\alpha\rfloor} , 
\end{align*}
where $C = C(n,m,q,\alpha) \in (0,\infty)$ are constants. 

Note that we do not assert that $f \in C^{\lfloor\alpha\rfloor,\alpha-\lfloor\alpha\rfloor}(B^2_{(1/32)^{1/q}}(0) \times B^{n-2}_{1/32}(0), \mathbb{R}^m)$, as $f$ is given by \eqref{topology eqn29} and we have only shown that $g \in C^{1,\mu/2}(B^{n-2}_{1/16}(0),\mathbb{R}^2)$.  Higher regularity of $f$ is directly related to higher regularity of the singular set of $T$ 
(as discussed in Remark~\ref{topology rmk} of the Introduction).
}\end{remark}

 \bigskip
\hskip-.2in\vbox{\hsize3in\obeylines\parskip -1pt %1\baselinestretch{0.8}
  \small 
Brian Krummel
School of Mathematics \& Statistics 
University of Melbourne
Parkville,VIC  3010, Australia
\vspace{4pt}
{\tt brian.krummel@unimelb.edu.au}} 
\vbox{\hsize3in
\obeylines 
\parskip-1pt 
\small 
Neshan Wickramasekera
DPMMS 
University of Cambridge 
Cambridge CB3 0WB, United Kingdom
\vspace{4pt}
{\tt N.Wickramasekera@dpmms.cam.ac.uk}
%\vspace{4pt}
%\hfill
%\hfill
%\hfill
%\hfill
%\hfill
}

\end{document}